\documentclass[a4paper, 11pt,reqno]{amsart}
\usepackage{amsthm,amssymb,amsmath,graphicx,psfrag,enumitem,mathrsfs}
\usepackage{mathtools}
\usepackage{dsfont}
\usepackage{bm}
\usepackage{standalone}
\usepackage{subcaption}
\usepackage[foot]{amsaddr}
\usepackage[british]{babel}
\usepackage[babel]{microtype}
\usepackage[margin=1in]{geometry}
\usepackage{verbatim}

\usepackage[textsize=footnotesize]{todonotes}

\usepackage{algorithm}
\usepackage{algpseudocode}
\usepackage{tikz}

\usetikzlibrary{backgrounds}

\usepackage[numbers,sort]{natbib}

\makeatletter
\def\NAT@spacechar{~}
\makeatother

\usepackage{enumitem}
\usepackage{hyperref}
\usepackage[capitalise]{cleveref}
\hypersetup{colorlinks=true,
   citecolor=blue,
   filecolor=blue,
   linkcolor=blue,
   urlcolor=blue
}

\crefname{figure}{Figure}{Figures}
\Crefname{figure}{Figure}{Figures}

\allowdisplaybreaks

\newtheorem{definition}{Definition}[section]
\newtheorem{claim}{Claim}[section]

\newtheorem{proposition}[definition]{Proposition}
\newtheorem{theorem}[definition]{Theorem}
\newtheorem{corollary}[definition]{Corollary}
\newtheorem{lemma}[definition]{Lemma}

\newtheorem{conjecture}[definition]{Conjecture}
\newtheorem{question}[definition]{Question}

\newtheorem{remark}[definition]{Remark}
\newenvironment{claimproof}{%
\let\origqed=\qedsymbol%
\renewcommand{\qedsymbol}{$\blacktriangleleft$}%
\begin{proof}}{\end{proof}\let\qedsymbol=\origqed}

\numberwithin{equation}{section}

\newcommand{\bigO}{\ensuremath{\mathcal{O}}}

\newcommand{\cA}{\mathcal{A}}
\newcommand{\cB}{\mathcal{B}}
\newcommand{\cC}{\mathcal{C}}
\newcommand{\cD}{\mathcal{D}}
\newcommand{\cE}{\mathcal{E}}
\newcommand{\cF}{\mathcal{F}}
\newcommand{\cG}{\mathcal{G}}
\newcommand{\cH}{\mathcal{H}}

\newcommand{\cL}{\mathcal{L}}
\newcommand{\cM}{\mathcal{M}}
\newcommand{\cN}{\mathcal{N}}

\newcommand{\cP}{\mathcal{P}}
\newcommand{\cQ}{\mathcal{Q}}

\newcommand{\cS}{\mathcal{S}}
\newcommand{\cT}{\mathcal{T}}
\newcommand{\cU}{\mathcal{U}}
\newcommand{\cV}{\mathcal{V}}
\newcommand{\cW}{\mathcal{W}}
\newcommand{\cX}{\mathcal{X}}

\newcommand{\cZ}{\mathcal{Z}}

\newcommand{\dist}{{\rm dist}}

\newcommand{\Var}{\operatorname{Var}}
\newcommand{\Cov}{\operatorname{Cov}}

\renewcommand{\epsilon}{\varepsilon}
\newcommand{\eps}{\varepsilon}

\newcommand{\eqp}{=_\mathrm{p}}
\newcommand{\neqp}{\neq_\mathrm{p}}
\newcommand{\cupdot}{\mathbin{\dot\cup}}

\newcommand{\COMMENT}[1]{}

\newcommand{\tcr}[1]{}

\usepackage{imakeidx}
\makeindex[columns=3, title=Notation index]

\title{Hamiltonicity of random subgraphs of the hypercube}

\address{School of Mathematics, University of Birmingham, Edgbaston, Birmingham, B15 2TT, United Kingdom.}

\author[P.~Condon]{Padraig Condon}
\email{pxc644@alumni.bham.ac.uk} 
\author[A.~Espuny D\'iaz]{Alberto Espuny D\'iaz}
\email{axe673@alumni.bham.ac.uk}
\author[A.~Gir\~{a}o]{Ant\'onio Gir\~{a}o}
\email{giraoa@bham.ac.uk}
\author[D.~K\"uhn]{Daniela K\"uhn}
\email{d.kuhn@bham.ac.uk}
\author[D.~Osthus]{Deryk Osthus}
\email{d.osthus@bham.ac.uk}

\thanks{This project has received partial funding from the European Research 
Council (ERC) under the European Union's Horizon 2020 research and innovation programme (grant agreement no.~786198, P.~Condon, A.~Espuny D\'iaz, D.~K\"uhn and D.~Osthus).
The research leading to these results was also partially supported by the EPSRC, grant nos.~EP/N019504/1 (A.~Gir\~ao and D.~K\"uhn) and EP/S00100X/1  (D.~Osthus), as well as the Royal Society and the Wolfson Foundation (D.~K\"uhn).\\
An extended abstract of this paper appeared in the \emph{Proceedings of the 2021 ACM-SIAM Symposium on Discrete Algorithms (SODA)}}

\date{\today}

\begin{document}

\begin{abstract}
We study Hamiltonicity in random subgraphs of the hypercube $\cQ^n$.
Our first main theorem is an optimal hitting time result.
Consider the random process which includes the edges of $\mathcal{Q}^n$ according to a uniformly chosen random ordering.
Then, with high probability, as soon as the graph produced by this process has minimum degree $2k$, it contains $k$ edge-disjoint Hamilton cycles, for any fixed $k\in\mathbb{N}$.
Secondly, we obtain a perturbation result:
if $H\subseteq\cQ^n$ satisfies $\delta(H)\geq\alpha n$ with $\alpha>0$ fixed and we consider a random binomial subgraph $\cQ^n_p$ of $\cQ^n$ with $p\in(0,1]$ fixed, then with high probability $H\cup\cQ^n_p$ contains $k$ edge-disjoint Hamilton cycles, for any fixed $k\in\mathbb{N}$.
In particular, both results resolve a long standing conjecture, posed e.g.~by Bollob\'as, that the threshold probability for Hamiltonicity in the random binomial subgraph of the hypercube equals $1/2$.
Our techniques also show that, with high probability, for all fixed $p\in(0,1]$ the graph $\cQ^n_p$ contains an almost spanning cycle.
Our methods involve branching processes, the R\"odl nibble, and absorption.
\end{abstract}

\maketitle

\thispagestyle{empty}

\section{Introduction}\label{introduction}

The $n$-dimensional \emph{hypercube} $\cQ^n$\index{Qn@$\cQ^n$} is the graph whose vertex set consists of all $n$-bit $01$-strings, where two vertices are joined by an edge whenever their corresponding strings differ by a single bit.
The hypercube and its subgraphs have attracted much attention in graph theory and computer science, e.g.~as a sparse network model with strong connectivity properties.
It is well known that hypercubes contain spanning paths (also called \emph{Gray codes} or \emph{Hamilton paths}) and, for all $n\geq2$, they contain spanning cycles (also referred to as \emph{cyclic Gray codes} or \emph{Hamilton cycles}).
Classical applications of Gray codes in computer science are described in the surveys of \citet{Sava97} and \citet{Knuth05}.
Applications of hypercubes to parallel computing are discussed in the monograph of \citet{Leigh92}.

\subsection{Spanning subgraphs in hypercubes}

The systematic study of spanning paths, trees and cycles in hypercubes was initiated in the 1970s.
There is by now an extensive literature about subtrees of the hypercube; see, for instance, results of \citet{BCLR92} about embedding subdivided trees (instigated by processor allocation in distributed computing systems).

As a generalization of Hamilton paths, \citet{CK07} considered the problem of finding a collection of spanning vertex-disjoint paths, given a prescribed set of endpoints.
After several improvements \cite{Chen13,GD08}, this problem was recently resolved by \citet{DGK17}.

The applications of hypercubes as networks in computer science inspired questions about the reliability of its properties.
This led to considering `faulty' hypercubes in which some edges or vertices are missing.
For instance, \citet{CL91} showed that, if $\cQ^n$ has at most $2n-5$ faulty edges and every vertex has (non-faulty) degree at least $2$, then there is a Hamilton cycle in $\cQ^n$ which avoids all faulty edges (and this condition is best possible).
They also showed that the general problem of determining the Hamiltonicity of $\cQ^n$ with a larger number of faulty edges is NP-complete.
More generally, \citet{DG08} studied the existence of spanning collections of vertex-disjoint paths with prescribed endpoints in faulty hypercubes.
(We will apply these results in our proofs, see \cref{connect} for details.)
These can be seen as extremal results about the \emph{robustness} of the hypercube with respect to containing spanning collections of paths and cycles.

\subsection{Hamilton cycles in binomial random graphs}\label{section12}

One of the most studied random graph models is the binomial random graph $G_{n,p}$.
Here we have a (labelled) set of $n$ vertices and we include each edge with probability $p$ independently of all other edges.

Given some monotone increasing graph property $\mathcal{P}$, a function $p^*=p^*(n)$ is said to be a (coarse) \emph{threshold} for $\mathcal{P}$ if $\mathbb{P}[G_{n,p}\in\mathcal{P}]\to1$ whenever $p/p^*\to\infty$ and $\mathbb{P}[G_{n,p}\in\mathcal{P}]\to0$ whenever $p/p^*\to0$.
One can define the stronger notion of a sharp threshold similarly: $p^*=p^*(n)$ is said to be a \emph{sharp threshold} for $\mathcal{P}$ if, for all $\varepsilon>0$, we have that $\mathbb{P}[G_{n,p}\in\mathcal{P}]\to1$ whenever $p\geq(1+\varepsilon)p^*$ and $\mathbb{P}[G_{n,p}\in\mathcal{P}]\to0$ whenever $p\leq(1-\varepsilon)p^*$.
The problem of finding the threshold for the containment of a Hamilton cycle was solved independently by \citet{Posa76} and \citet{Kor77}. 
Furthermore, \citet{Kor77} determined the sharp threshold for Hamiltonicity to be $p^*=\log n/n$. 
These results were later made even more precise by \citet{KS83}.
It is worth noting that $p^*=\log n/n$ is also the sharp threshold for the property of having minimum degree at least $2$.
In this sense, the results about Hamilton cycles in $G_{n,p}$ can be interpreted as saying that the natural obstruction of having sufficiently high minimum degree is also an `almost sufficient' condition.

A property that generalises Hamiltonicity is that of containing $k$ edge-disjoint Hamilton cycles, for some $k\in\mathbb{N}$.
We will present more results in this direction in \cref{introduction3}; for now, let us simply note that the sharp threshold for the containment of $k$ edge-disjoint Hamilton cycles in $G_{n,p}$, for some $k\in\mathbb{N}$ independent of $n$, is $p^*=\log n/n$, i.e.~the same as the threshold for Hamiltonicity.

The study of \emph{robustness} of graph properties has also attracted much attention recently.
For instance, given a graph $G$ which is known to satisfy some property $\mathcal{P}$, consider a random subgraph $G_p$ obtained by deleting each edge of $G$ with probability $1-p$, independently of all other edges.
The problem then is to determine the range of $p$ for which $G_p$ satisfies $\mathcal{P}$ with high probability.
In this setting, a result of \citet{KLS14} asserts that, for any $n$-vertex graph $G$ with minimum degree at least $n/2$, the graph $G_p$ is asymptotically almost surely Hamiltonian whenever $p\gg\log n/n$.
This can be viewed as a robust version of Dirac's theorem on Hamilton cycles.

\subsection{Hamilton cycles in binomial random subgraphs of the hypercube}\label{section13}

Throughout this paper, we will consider random subgraphs of the hypercube and show that the hypercube is robustly Hamiltonian in the above sense.
We will denote by $\cQ^n_p$\index{Qnp@$\cQ^n_p$} the random subgraph of the hypercube obtained by removing each edge of $\cQ^n$ with probability $1-p$ independently of every other edge.

The random graph $\cQ^n_p$ was first studied by \citet{Burt77}, who proved that the sharp threshold for connectivity is $1/2$.
This result was later made more precise by \citet{ES79} and \citet{Bol83}.
As a related problem, \citet{Frieze87} determined the sharp threshold for connectivity in subgraphs of $\cQ^n$ obtained by removing both vertices and edges uniformly at random.
Later, \citet{Bol90} proved that $1/2$ is also the sharp threshold for the containment of a perfect matching in $\cQ^n_p$.
As with the $G_{n,p}$ model, this also coincides with the threshold for having minimum degree at least $1$. 

The main goal of this paper is to study the analogous problem for Hamiltonicity in random subgraphs of the hypercube.
There is a folklore conjecture that the sharp threshold for Hamiltonicity in $\cQ^n_p$ should be $1/2$, i.e.~the same as the threshold for having minimum degree at least $2$.
This question was explicitly asked by \citet{BBPC} at several conferences in the 1980s, in the ICM surveys of \citet{FriezeICM14} and \citet{KOICM14}, as well as in the recent survey of \citet{Frieze19}.
A special case of our first result resolves this problem.

\begin{theorem}\label{thm:thresholdk}
For any $k\in\mathbb{N}$, the sharp threshold for the property of containing $k$ edge-disjoint Hamilton cycles in $\cQ^n_p$ is $p^*=1/2$.
\end{theorem}

For $k=1$, this can be seen as a probabilistic version of the result on faulty hypercubes~\cite{CL91}\COMMENT{Try to relate this more to networks by discussing that this is a good model if, say, edges fail after a certain amount of time, where this time follows some distribution (say, exponential) and is independent for every edge?}, and also as a statement about the robustness of Hamiltonicity in the hypercube.

While, for $p<1/2$, with high probability $\cQ^n_p$ will not contain a Hamilton cycle, it turns out that the reason for this is mostly due to local obstructions (e.g., vertices with degree zero or one).
More precisely, we prove that, for any constant $p\in(0,1/2)$, a.a.s.~the random graph $\cQ^n_p$ contains an almost spanning cycle.

\begin{theorem}\label{thm:almost}
For any $\delta,p\in(0,1]$, a.a.s.~the graph $\cQ^n_p$ contains a cycle of length at least $(1-\delta)2^n$.
\end{theorem}

We believe that the probability bound is far from optimal, in the sense that random subgraphs of the hypercube where edges are picked with vanishing probability should also satisfy this property.

\begin{conjecture}\label{conj1}
Suppose that $p=p(n)$ satisfies that $pn\to\infty$.
Then, a.a.s.~$\cQ^n_p$ contains a cycle of length $(1-o(1))2^{n}$.
\end{conjecture}

\COMMENT{This is a very nice conjecture, I was thinking a little bit and our theorems from the tree section (very slightly tweaked) prove that a.a.s.~there is a bounded degree tree covering all but a vanishing proportion of $\cQ^{n}$. But to get a cycle all of the rest of the paper is useless; even with two layers no tree in the intersection graph appears, $\ell$-cubes roughly never appear since a cube has more edges than the number of cubes containing a vertex. It seems the approach has to be completely different.}
\COMMENT{Should we say something about the $2$-core?}

Similarly, it would be interesting to determine which (long) paths and (almost spanning) trees can be found in $\cQ^n_p$.
Moreover, our methods might also be useful to embed other large subgraphs, such as $F$-factors.

\begin{conjecture}\label{conjFfactors}
Suppose $\eps>0$ and an integer $\ell\geq2$ are fixed and $p \ge 1/2+\eps$.
Then, a.a.s.~$\cQ^n_p$ contains a $C_{2^\ell}$-factor, that is, a set of vertex-disjoint cycles of length $2^\ell$ which together contain all vertices of $\cQ^n$.
\end{conjecture}

\subsection{Hitting time results}\label{introduction3}

Remarkably, the above intuition that having the necessary minimum degree is an `almost sufficient' condition for the containment of edge-disjoint perfect matchings and Hamilton cycles can be strengthened greatly via so-called hitting time results.
These are expressed in terms of random graph processes.
The general setting is as follows.
Let $G$ be an $n$-vertex graph with $m=m(n)$ edges, and consider an arbitrary labelling $E(G)=\{e_1,\ldots,e_m\}$.
The \emph{$G$-process} is defined as a random sequence of nested graphs $\tilde{G}(\sigma)=(G_t(\sigma))_{t=0}^m$\index{Gtilde@$\tilde{G}(\sigma)$}, where $\sigma$ is a permutation of $[m]$ chosen uniformly at random and, for each $i\in[m]_0$, we set $G_i(\sigma)=(V(G),E_i)$, where $E_i\coloneqq\{e_{\sigma(j)}: j\in[i]\}$.
Given any monotone increasing graph property $\mathcal{P}$ such that $G\in\mathcal{P}$, the \emph{hitting time} for $\mathcal{P}$ in the above $G$-process is the random variable $\tau_\mathcal{P}(\tilde{G}(\sigma))\coloneqq\min\{t\in[m]_0: G_t(\sigma)\in\mathcal{P}\}$\index{tau@$\tau_\mathcal{P}(\tilde{G}(\sigma))$}.

Let us denote the properties of containing a perfect matching by $\mathcal{PM}$, Hamiltonicity by $\mathcal{HAM}$, and connectivity by $\mathcal{CON}$, respectively.
For any $k\in\mathbb{N}$, let $\boldsymbol{\delta}k$\index{deltak@$\boldsymbol{\delta}k$} denote the property of having minimum degree at least $k$, and let $\mathcal{HM}k$\index{HMk@$\mathcal{HM}k$} denote the property of containing $\lfloor k/2\rfloor$ edge-disjoint Hamilton cycles and, if $k$ is odd, one matching of size $\lfloor n/2\rfloor$ which is edge-disjoint from these Hamilton cycles\COMMENT{Note that these are all monotone increasing properties.}.
With this notion of hitting times, many of the results about thresholds presented in \cref{section12,section13} can be strengthened significantly.
For instance, \citet{BT85} showed that a.a.s.~$\tau_{\mathcal{CON}}(\tilde{K_n}(\sigma))=\tau_{\boldsymbol{\delta}1}(\tilde{K_n}(\sigma))$ and, if $n$ is even, then a.a.s.~$\tau_{\mathcal{PM}}(\tilde{K_n}(\sigma))=\tau_{\boldsymbol{\delta}1}(\tilde{K_n}(\sigma))$.
\citet{AKS85} and \citet{Bol84} independently proved that a.a.s.~$\tau_{\mathcal{HAM}}(\tilde{K_n}(\sigma))=\tau_{\boldsymbol{\delta}2}(\tilde{K_n}(\sigma))$.
This was later generalised by \citet{BF85}, who proved that, given any $k\in\mathbb{N}$, for $n$ even a.a.s.~$\tau_{\mathcal{HM}k}(\tilde{K_n}(\sigma))=\tau_{\boldsymbol{\delta}k}(\tilde{K_n}(\sigma))$.

A hitting time result for the property of having $k$ edge-disjoint Hamilton cycles when $k$ is allowed to grow with $n$ is still not known, even in $K_n$-processes.
As a slightly weaker notion, consider property $\mathcal{H}$, where we say that a graph $G$ satisfies property $\mathcal{H}$ if it contains $\lfloor\delta(G)/2\rfloor$ edge-disjoint Hamilton cycles, together with an additional edge-disjoint matching of size $\lfloor n/2 \rfloor$ if $\delta(G)$ is odd.
\citet{KKO15}, \citet{KS12} as well as \citet{KO14} proved results for different ranges of $p$ which, together, show that $G_{n,p}$ a.a.s.~satisfies property $\mathcal{H}$.

For graphs other than the complete graph, \citet{Joha18} recently obtained a robustness version of the hitting time results for Hamiltonicity.
In particular, for any $n$-vertex graph $G$ with $\delta(G)\geq(1/2+\varepsilon)n$, he proved that a.a.s.~$\tau_{\mathcal{HAM}}(\tilde{G}(\sigma))=\tau_{\boldsymbol{\delta}2}(\tilde{G}(\sigma))$.
This was later extended to a larger class of graphs $G$ and to hitting times for $\mathcal{HM}2k$, for all $k\in\mathbb{N}$ independent of $n$, by \citet{AK19}.

In the setting of random subgraphs of the hypercube, \citet{Bol90} determined the hitting time for perfect matchings by showing that a.a.s.~$\tau_{\mathcal{PM}}(\tilde{\cQ^n}(\sigma))=\tau_{\mathcal{CON}}(\tilde{\cQ^n}(\sigma))=\tau_{\boldsymbol{\delta}1}(\tilde{\cQ^n}(\sigma))$.
One of our main results (which implies \cref{thm:thresholdk}) is a hitting time result for Hamiltonicity (and, more generally, property $\cH\cM k$) in $\cQ^n$-processes.
Again, this question was raised by \citet{BBPC} at several conferences.

\begin{theorem}\label{thm:hitting}
For all $k\in\mathbb{N}$, a.a.s.~$\tau_{\mathcal{HM}k}(\tilde{\cQ^n}(\sigma))=\tau_{\boldsymbol{\delta}k}(\tilde{\cQ^n}(\sigma))$, that is, the hitting time for the containment of a collection of $\lfloor k/2 \rfloor$ Hamilton cycles and $k-2\lfloor k/2\rfloor$ perfect matchings, all pairwise edge-disjoint, in $\cQ^n$-processes is a.a.s.~equal to the hitting time for the property of having minimum degree at least $k$.
\end{theorem}

We also wonder whether this is true if $k$ is allowed to grow with $n$, and propose the following conjecture which, if true, would be an approximate version of the results of \cite{KS12,KO14,KKO15} in the hypercube.

\begin{conjecture}
For all $p\in (1/2,1]$ and $\eta>0$, a.a.s.~$\cQ^n_p$ contains $(1/2-\eta)\delta(\cQ^n_p)$ edge-disjoint Hamilton cycles.
\end{conjecture}

\subsection{Randomly perturbed graphs}\label{sect:1.5}

A relatively recent area at the interface of extremal combinatorics and random graph theory is the study of \emph{randomly perturbed graphs}.
Generally speaking, the idea is to consider a deterministic dense $n$-vertex graph $H$ (usually satisfying some minimum degree condition) and a random graph $G_{n,p}$ on the same vertex set as $H$.
The question is whether $H$ is close to satisfying some given property $\mathcal{P}$ in the sense that a.a.s.~$H\cup G_{n,p}\in\mathcal{P}$ for some small $p$.
This line of research was sparked off by \citet{BFM03}, who showed that, if $H$ is an $n$-vertex graph with $\delta(H)\geq\alpha n$, for any constant $\alpha>0$, then a.a.s.~$H\cup G_{n,p}$ is Hamiltonian for all $p\geq C(\alpha)/n$.
Other properties that have been studied in this context are e.g.~the existence of powers of Hamilton cycles and general bounded degree spanning graphs~\cite{BMPP20}, $F$-factors~\cite{BTW19} or spanning bounded degree trees~\cite{KKS17,BHKMPP19}.
One common phenomenon in this model is that, by considering the union with a dense graph $H$ (i.e.~a graph $H$ with linear degrees), the probability threshold of different properties is significantly lower than that in the classical $G_{n,p}$ model.
The results for Hamiltonicity~\cite{BFM03} were very recently generalised by \citet{HMMMP20} to allow $\alpha$ to tend to $0$ with $n$ (that is, to allow graphs $H$ which are not dense).

We consider randomly perturbed graphs in the setting of subgraphs of the hypercube. 
To be precise, we take an arbitrary spanning subgraph $H$ of the hypercube, with linear minimum degree, and a random subgraph $\cQ^n_\varepsilon$, and consider $H\cup\cQ^n_\varepsilon$.
(Note here that $\cQ^n_\varepsilon$ is a `dense' subgraph of $\cQ^n$, but for $\varepsilon<1/2$ it will contain both isolated vertices and vertices of very low degrees.)
In this setting, we show the following result.

\begin{theorem}\label{thm:main}
For all $\eps, \alpha\in(0,1]$ and $k\in\mathbb{N}$, the following holds.
Let $H$ be a spanning subgraph of $\cQ^{n}$ such that $\delta(H)\geq\alpha n$.
Then, a.a.s.~$H \cup \cQ^n_\varepsilon$ contains $k$ edge-disjoint Hamilton cycles.
\end{theorem}

We can also allow $H$ to have much smaller degrees, at the cost of requiring a larger probability to find the Hamilton cycles.

\begin{theorem}\label{thm: kedgehit}
For every integer $k\geq 2$, there exists $\eps>0$ such that a.a.s., for every spanning subgraph $H$ of $\cQ^{n}$ with $\delta(H)\geq k$, the graph $H \cup \cQ_{1/2-\eps}^{n}$ contains a collection of $\lfloor k/2 \rfloor$ Hamilton cycles and $k-2\lfloor k/2 \rfloor$ perfect matchings, all pairwise edge-disjoint.
\end{theorem}

Note that \cref{thm: kedgehit} can be viewed as a `universality' result for $H$, meaning that it holds for all choices of $H$ simultaneously.
It would be interesting to know whether such a result can also be obtained for the lower edge probability assumed in \cref{thm:main}, i.e., is it the case that, for all $\eps,\alpha\in(0,1]$, a.a.s.~$G\sim\cQ^n_\eps$ has the property that, for every spanning $H\subseteq \cQ^n$ with $\delta(H)\geq\alpha n$, $G\cup H$ is Hamiltonian?

As we will prove, \cref{thm:thresholdk} follows straightforwardly from \cref{thm:main}, and it follows trivially from \cref{thm:hitting}.
In turn, \cref{thm:hitting} follows from \cref{thm: kedgehit}.
On the other hand, \cref{thm:almost,,thm:main,,thm: kedgehit}, while being proved with similar ideas, are incomparable. 

\subsection{Percolation on the hypercube}

To build Hamilton cycles in random subgraphs of the hypercube, we will consider a random process which can be viewed as a branching process or percolation process on the hypercube.
With high probability, for constant $p>0$, this process results in a bounded degree tree in $\cQ^n_p$ which covers most of the neighbourhood of every vertex in $\cQ^n$, and thus spans almost all vertices of $\cQ^n$.
The version stated below is a special case of \cref{lem: main treereshit}.

\begin{theorem}\label{lem: main tree intro}
For any fixed $\varepsilon,p\in(0,1]$, there exists $D=D(\varepsilon)$ such that a.a.s.~$\cQ^n_p$ contains a tree $T$ with $\Delta(T)\leq D$ and such that $|V(T)\cap N_{\cQ^n}(x)|\geq(1-\varepsilon)n$ for every $x\in V(\cQ^n)$.
\end{theorem}

Further results concerning the local geometry of the giant component in $\cQ^n_p$ for constant $p\in(0,1/2)$ were proved recently by \citet{MSW18}.

The random process we consider in the proof of \cref{lem: main tree intro} can be viewed as a branching random walk (with a bounded number of branchings at each step).
Simpler versions of such processes (with infinite branchings allowed) have been studied by \citet{FP93} and \citet{KKO}, and we will base our analysis on these.
Motivated by our approach, we raise the following question, which seems interesting in its own right.

\begin{question}
Does a non-returning random walk on $\cQ^n$ a.a.s.~visit almost all vertices of~$\cQ^n$?
\end{question}

More generally, there are many results and applications concerning random walks on the hypercube (but allowing for returns). 
For example, motivated by a processor allocation problem, \citet{BC} studied a walk algorithm to embed large (subdivided) trees into the hypercube.
Moreover, the analysis of (branching) random walks is a critical ingredient in the study of percolation thresholds for the existence of a giant component in $\cQ^n_p$.
These have been investigated e.g.~by \citet{BKL, BCVSS} and \citet{HN17}.
Related questions concerning the structure of the giant component are also investigated by \citet{EKK21} and discussed in the survey of \citet{HN14}.

\subsection{Organisation of the paper}

In \cref{sect:outline main} we provide an overview of our ideas and proof methods.
In \cref{notation} we introduce the notation we will use throughout the paper.
In \cref{prelim} we state the different probabilistic tools, as well as some other well-known results, that we will call on, and in \cref{aux} we collect various results on matchings and random subgraphs of the hypercube.
In \cref{nib} we prove \cref{thm: nibble}, our main cube tiling result, and in \cref{section:tree} we prove \cref{lem: main treereshit}, our main near-spanning tree result (see \cref{sect:outline main} for more details on each of these).
Then, in \cref{section8} we prove \cref{thm:main} in the case $k=1$ (see \cref{thm:main1}).
In \cref{sect:thm1} we use this to deduce the general statement of \cref{thm:main}, and also deduce \cref{thm:thresholdk} and explain how to obtain \cref{thm:almost}.
Finally, in \cref{sect:hitting} we show how to modify the proof of \cref{thm:main} to obtain \cref{thm: kedgehit}, and thus our hitting time result (\cref{thm:hitting}).


\section{Outline of the main proofs}\label{sect:outline main}

\subsection{Overall outline}\label{sub:over}

\COMMENT{Some aspects that are currently left out/dishonest in this sketch (let me know if any to be added):\\
The order of the Steps in the main proof is not adhered to (e.g. I mentioned scant molecules after the nibble) for the purpose of clarity.\\
Slices in molecules are left out, and as a result the difference between an external and full skeleton.\\
The concepts of bondlessness and bondlessly surrounded are left out, and therefore their roles in connecting lemmas and rainbow matching lemma.\\
The extension tree is not mentioned.\\
Many of the subtleties of the algorithm are not mentioned/discussed.}
We now sketch the key ideas for the proof of \cref{thm:main}.
We will first prove the case $k=1$, and later use this to deduce the case when $k>1$.
Recall we are given $H\subseteq \cQ^n$ with $\delta(H)\geq\alpha n$, and $G \sim \cQ^n_\eps$, with $\alpha,\varepsilon\in(0,1]$. 
Our aim is to show that a.a.s.~$H\cup G$ is Hamiltonian.

Our approach for finding a Hamilton cycle is to first obtain a spanning tree.
By passing along all the edges of a spanning tree $T$ (with a vertex ordering prescribed by a depth first search), one can create a closed spanning walk $W$ which visits every edge of $T$ twice.
The idea is then to modify such a walk into a Hamilton cycle.
(This approach is inspired by the approximation algorithm for the Travelling Salesman Problem which returns a tour of at most twice the optimal length.)\COMMENT{For the version of the TSP that assumes the triangle inequality. 
You take a minimum spanning tree, and use each edge at most twice to visit every vertex and get back to the start. 
Then use triangle inequality to skip repeated edges.}
More precisely, our approach will be to obtain a near-spanning tree of $\cQ^{n-s}$, for some suitable constant~$s$, and to blow up vertices of this tree into $s$-dimensional cubes (see \cref{newfig}). 
These cubes can then be used to move along the tree without revisiting vertices, which will result in a near-Hamilton cycle $\mathfrak{H}$.
All remaining vertices which are not included in $\mathfrak{H}$ will be absorbed into $\mathfrak{H}$ via absorbing structures that we carefully put in place beforehand.

In \cref{section:outline2,section:outline3,sect:cycle} we outline in more detail how we find a long cycle in $G$ (\cref{thm:almost}).
Note that in \cref{thm:almost} we have $G \sim \cQ^n_\eps$, so a.a.s.~$G$ will have isolated vertices which prevent any Hamilton cycle occurring as a subgraph.
In \cref{section:outline5} we outline how we build on this approach to obtain the case $k=1$ of \cref{thm:main}.
In \cref{section:outline6} we sketch how we obtain \cref{thm:hitting}.

\subsection{Building block I: trees via branching processes.}\label{section:outline2}

We view each vertex in $\cQ^n$ as an $n$-dimensional $01$-coordinate vector.
By fixing the first $s$ coordinates, we fix one of $2^s$ \emph{layers} $L_1, \dots, L_{2^s}$\index{Lilayer@$L_i$ (layer)} of the hypercube, where $s\in \mathbb{N}$ will be constant.
Thus, $L\cong \cQ^{n-s}$ for each layer~$L$.
By considering a Hamilton cycle in $\cQ^s$, we may assume that consecutive layers differ only by a single coordinate on the unique elements of $\cQ^s$ which define them.
Let $G \sim\cQ^n_\eps$. 
For each layer $L$, we let $L(G)\coloneqq G[V(L)]$\index{LG@$L(G)$} and, by momentarily viewing these layers as different subgraphs on the vertex set of $\cQ^{n-s}$, we define the \emph{intersection graph} $I(G)\coloneqq \bigcap_{i=1}^{2^s}L_i(G)$\index{IG@$I(G)$}.
Hence,~$I(G) \sim \cQ^{n-s}_{\eps^{2^s}}$.
We view $I(G)$ as a subgraph of $\cQ^{n-s}$.
We first show that $I(G)$ contains a near-spanning tree $T$ (\cref{lem: main treereshit}).
Thus, a copy of $T$ is present in each of $L_1(G), \dots, L_{2^s}(G)$ simultaneously.

Since the walk $W$ mentioned in \cref{sub:over} passes through each vertex $x$ of $T$ a total of $d_T(x)$ times, it will be important later for $T$ to have bounded degree.\COMMENT{Otherwise a blown-up vertex in this tree could be `overused.'}
In order to guarantee this, we run bounded degree branching processes (see \cref{def: perc}) from several far apart `corners' of the hypercube.
Roughly speaking, $T$ will be formed by taking a union of these processes and removing cycles.
Crucially, the model we introduce for these processes has a joint distribution with $\cQ^{n-s}_{\eps^{2^s}}$, so that $T$ will in fact appear as a subgraph of $I(G)$.
In applying \cref{lem: main treereshit}, we obtain a bounded degree tree $T\subseteq I(G)$ which contains almost all of the neighbours of every vertex of $I(G)$.
We also obtain a `small' \emph{reservoir} set $R\subseteq V(I(G))$, which $T$ avoids and which will play a key role later in the absorption of vertices which do not belong to our initial long cycle.
At this point, both $T$ and $R$ are now present in every layer of the hypercube simultaneously.

\subsection{Building block II: cube tilings via the nibble.}\label{section:outline3}

Let $\ell < s$ be fixed (in our proof, $\ell$ will be a sufficiently large constant and we will take $s=10\ell$ and $n$ sufficiently large compared to these).
In order to gain more local flexibility when traversing the near-spanning tree $T$, we augment $T$ by locally adding a near-spanning $\ell$-cube factor of $I(G)$.
One can use classical results on matchings in almost
regular uniform hypergraphs of small codegree to show that $I(G)$ contains such a collection of $\cQ^\ell$ spanning almost all vertices of $I(G)$.
However, we require the following stronger properties, namely that there exists a collection $\cC$ of vertex disjoint copies of $\cQ^\ell$ in $I(G)$ so that, for each $x \in V(I(G))$,
\begin{enumerate}[label=(\roman*)]
\item\label{item:sketchnibble1} $\cC$ covers almost all vertices in $N_{I(\cQ^n)}(x)$;
\item\label{item:sketchnibble2} the directions spanned by the cubes intersecting $N_{I(\cQ^n)}(x)$ do not correlate too strongly with any given set of directions. (Here, a `direction' refers to the coordinate vector which corresponds to an edge.)
\end{enumerate}
The precise statement is given in \cref{thm: nibble}.
Neither \ref{item:sketchnibble1} nor \ref{item:sketchnibble2} follow from existing results on hypergraph matchings and the proofs strongly rely on geometric properties intrinsic to the hypercube.

To prove \cref{thm: nibble}, we build on the so-called R\"odl nibble.
More precisely, we consider the hypergraph $\cH$, with $V(\cH)= V(\cQ^{n-s})$, where the edge set is given by the copies of $\cQ^{\ell}$ in $I(G)$.
We run a random iterative process where at each stage we add a `small' number of edges from $\cH$ to $\cC$, before removing  all those remaining edges of $\cH$ which `clash' with our selection.
A careful analysis and an application of the Lov\'asz local lemma yield the existence of an instance of this process which terminates in the near-spanning $\ell$-cube factor with the properties required  for \cref{thm: nibble}.

\subsection{Constructing a long cycle.}\label{sect:cycle}

Roughly speaking, we will use $T$ as a backbone to provide `global' connectivity, and will use the near-spanning $\ell$-cube factor $\cC$ and the layer structure to gain high `local' connectivity and flexibility.
We show a representation of this structure in \cref{newfig}.
Let $T \cup \bigcup_{C \in \cC}C \eqqcolon \Gamma' \subseteq I(G)$ and let $\Gamma \subseteq \Gamma'$ be formed by removing all leaves and isolated cubes in $\Gamma'$.
It follows by our tree and nibble results that almost all vertices of $I(G)$ are contained in $\Gamma$.
Note that, for each $v \in V(\cQ^{n-s})=V(I(G))$, there is a unique vertex in each of the $2^s$ layers which corresponds to $v$.
We refer to these $2^s$ vertices as \emph{clones} of $v$ and
 to the collection of these $2^s$ clones as a \emph{vertex molecule}.
Similarly, each $\ell$-cube $C \in \cC$ contained in $\Gamma$ gives rise to a \emph{cube molecule}.
We construct a cycle in $G$ which covers all of the cube molecules (and, therefore, almost all vertices in $\cQ^n$).

\begin{figure}
    \centering
    \includestandalone[width=0.5\textwidth]{figure4}
    \caption{\footnotesize{A representation of the main structure used for the proofs.
    We think of $\cQ^n$ as a `product' of two smaller cubes.
    Each `horizontal' cube represents a copy of $\cQ^{n-s}$, and the red `vertical' cube represents $\cQ^s$.
    All `horizontal' cubes contain a copy of the same tree $T$ and the same cube tiling $\cC$ (which are consistently distributed with respect to the `vertical' cube; this gives rise to `cube molecules').
    When finding a long cycle, cube molecules are highly connected and can be covered by few paths, and the tree is used to join cube molecules to one another.}}
    \label{newfig}
\end{figure}

Let $\Gamma^*$ be the graph obtained from $\Gamma$ by contracting each $\ell$-cube $C \subseteq \Gamma$ into a single vertex.
We refer to such vertices in $\Gamma^*$ as \emph{atomic vertices}, and to all other vertices as \emph{inner tree vertices}.
We run a depth-first search on $\Gamma^*$ to give an order to the vertices.
Next, we construct a \emph{skeleton} which will be the backbone for our long cycle.
The skeleton is an ordered sequence of vertices in $\cQ^n$ which contains the vertices via which our cycle will enter and exit each molecule.  
That is, given an \emph{exit vertex} $v$ for some molecule in the skeleton, the vertex $u$ which succeeds $v$ in the skeleton will be an \emph{entry vertex} for another molecule, and such that $uv \in E(G)$.
Here, a vertex in the skeleton belonging to an inner tree vertex molecule is referred to as both an entry and exit vertex.
(Actually, we will first construct an `external skeleton', which encodes this information.
The skeleton then also prescribes some edges within molecules which go between different layers.)
We use the ordering of the vertices of $\Gamma^*$ to construct the skeleton in a recursive way starting from the lowest ordered vertex.
It is crucial that our tree $T$ has bounded degree (much smaller than $2^s$), so that no molecule is overused in the skeleton.

Once the skeleton is constructed, we apply our `connecting lemmas' (\cref{lem:slicecover,lem:slicecover2}).
These connecting lemmas, applied to a cube molecule with a bounded number of pairs of entry and exit vertices as input (given by the skeleton), provide us with a sequence of vertex-disjoint paths which cover this molecule,  where each path has start and end vertices consisting of an input pair.
The union of all of these paths combined with all edges in $G$ between the successive exit and entry vertices of the skeleton will then form a cycle $\mathfrak{H}\subseteq G$ which covers all vertices lying in the cube molecules (thus proving \cref{thm:almost}).

\subsection{Constructing a Hamilton cycle.}\label{section:outline5}

In order to construct a Hamilton cycle in $H\cup G$, we will absorb the vertices of $V(\cQ^n)\setminus V(\mathfrak{H})$ into $\mathfrak{H}$. 
We achieve this via absorbing structures that we identify for each vertex (see \cref{def:abs}).
To construct these absorbing structures, we will need to use some edges of $H$.
Roughly speaking, to each vertex $v$ we associate a left $\ell$-cube $C^l_v\subseteq \cQ^n$ and a right $\ell$-cube $C^r_v\subseteq \cQ^n$, where $C^l_v,C^r_v$ are both clones of some $\ell$-cubes $C^l, C^r \in \cC$ contained in  $\Gamma$.
We choose these cubes so that $v$ will have a neighbour $u \in V(C^l_v)$ and a neighbour $u' \in V(C^r_v)$, to which we refer as \emph{tips} of the absorbing structure.
Furthermore, $u$ will have a neighbour $w \in V(C^r_v)$, which is also a neighbour of $u'$.
Our near-Hamilton cycle $\mathfrak{H}$ will satisfy the following properties:
\begin{enumerate}[label=(\alph*)]
    \item\label{psa} $\mathfrak{H}$ covers all vertices in $C^l_v \cup C^r_v$ except for $u$, and
    \item\label{psb} $wu' \in E(\mathfrak{H})$.
\end{enumerate}
These additional properties will be guaranteed by our connecting lemmas discussed in \cref{sect:cycle}.
We can then alter $\mathfrak{H}$ to include the segment $wuvu'$ instead of the edge $wu'$, thus absorbing the vertices $u$ and $v$ into $\mathfrak{H}$.
See \cref{fig:abs-struct1} for a representation of the absorbing structure and the absorption process.

\begin{figure}
    \centering
    \begin{tikzpicture}
\draw [gray, thick] (1,1.35) -- (0,2.85) -- (1,4.35) -- (2,2.85) -- cycle;
\node[anchor = south west] at (1.5,3.6){$C^r_v$};
\draw [gray, thick] (-1.5,1.35) -- (-0.5,2.85) -- (-1.5,4.35) -- (-2.5,2.85) -- cycle;
\node[anchor = south east] at (-2,3.6){$C^l_v$};

\draw (-0.25,0) node[circle,fill,inner sep=2pt]{};
\node[anchor = north] at (-0.25,-0.05){$v$};
\draw (1,1.5) node[circle,fill,inner sep=2pt]{};
\node[anchor = north west] at (1,1.6){$u'$};
\draw (-1.5,1.5) node[circle,fill,inner sep=2pt]{};
\node[anchor = north east] at (-1.5,1.5){$u$};
\draw [ultra thick] (1,1.5) -- (-0.25,0) -- (-1.5,1.5);
\draw (0.5,2.25) node[circle,fill,inner sep=2pt]{};
\node[anchor = west] at (0.5,2.25){$w$};
\draw [ultra thick] (1,1.5) -- (0.5,2.25) -- (-1.5,1.5);
\draw [ultra thick, opacity = 0.2] (1,1.5) -- (-0.25,0) -- (-1.5,1.5);
\draw [ultra thick, opacity = 0.2] (1,1.5) -- (0.5,2.25) -- (-1.5,1.5);
\fill [blue, opacity = 0.3] (1,1.35) -- (0,2.85) -- (1,4.35) -- (2,2.85) -- cycle;
\fill [blue, opacity = 0.3] (-1.8,1.8) -- (-1.2,1.8) -- (-0.5,2.85) -- (-1.5,4.35) -- (-2.5,2.85) -- cycle;

\draw [ultra thick, blue] (1.1,1.75) -- (1,1.5) -- (0.5,2.25) -- (0.6,2.5);
\end{tikzpicture}
    \caption{\footnotesize{Representation of the absorbing structure for a vertex $v$.
    The left $\ell$-cube is drawn on the left, and the right $\ell$-cube is drawn on the right.
    When absorbing $v$, we can ensure that all vertices in the shaded part of the cubes are contained in an almost spanning cycle, which in particular contains the edge $wu'$.
    Replacing $wu'$ by $wuvu'$ allows to incorporate both missing edges into the cycle.}}
    \label{fig:abs-struct1}
\end{figure}

The following types of vertices will require absorption.
\begin{enumerate}[label = (\roman*)]
    \item\label{ps1} Every vertex that is not covered by a clone of either some inner tree vertex or of some cube $C \in \cC$ which is contained in $\Gamma$.
    \item\label{ps4} The cycle $\mathfrak{H}$ does not cover all the clones of inner tree vertices and, thus, the uncovered vertices of this type will also have to be absorbed.
\end{enumerate}
However, we will not know precisely which of the vertices described in \ref{ps4} will be covered by $\mathfrak{H}$ and which of these vertices will need to be absorbed until after we have constructed the (external) skeleton.
Moreover,  many potential absorbing structures are later ruled out as candidates (for example, if they themselves contain vertices that will need to be absorbed).
Therefore, it is important that we identify a `robust' collection of many potential absorbing structures for every vertex in $\cQ^n$ at a preliminary stage of the proof.
The precise absorbing structure eventually assigned to each vertex will be chosen via an application of our rainbow matching lemma (\cref{lem: rainbow}) at a late stage in the proof.

We will now highlight the purpose of the reservoir $R$.
Suppose $v \in V(\cQ^n)$ is a vertex which needs to be absorbed via an absorbing structure with left $\ell$-cube $C^l_v$ and left tip $u \in V(C^l_v)$.
Recall that both $u$ and $C^l_v$ are clones of some $u^* \in V(\Gamma)$ and $C^l \in \cC$, where $u^* \in V(C^l)$.
If $u^*$ has a neighbour $w^*$ in $T - V(C^l)$, then it is possible that the skeleton will assign an edge from $u$ to $w$ for the cycle $\mathfrak{H}$ (where $w$ is the clone of $w^*$ in the same layer as $u$).
Given that $u$ is now incident to a vertex outside of $C^l_v$, we can no longer use the absorbing structure with $u$ as a (left) tip (otherwise, we might disconnect $T$).
To avoid this problem, we show that most vertices have many potential absorbing structures whose tips lie in the reservoir $R$ (which $T$ avoids).
Here we make use of vertex degrees of $H$.
A small number of \emph{scant vertices} will not have high enough degree into $R$.
For these vertices we fix an absorbing structure whose tips do not lie in $R$, and then alter $T$ slightly so that these tips are deleted from $T$ and reassigned to $R$.
The fact that scant vertices are few and well spread out from each other will be crucial in being able to achieve this (see \cref{lem:repatch}).

Let us now discuss two problems arising in the construction of the skeleton.
Firstly, let $\cM_C \subseteq \cQ^n$ with $C \in \cC$ be a cube molecule which is to be covered by $\mathfrak{H}$.
Furthermore, suppose one of the clones $C^l_v$ of $C$ belongs to an absorbing structure for some vertex $v$.
Let $u$ be the tip of $C^l_v$ and suppose that $u$ has even parity.
We would like to apply the connecting lemmas to cover $\cM_C - \{u\}$ by paths which avoid $u$. 
But this would now involve covering one fewer vertex of even parity than of odd parity.
This, in turn, has the effect of making the construction of the skeleton considerably more complicated (this construction is simplest when successive entry and exit vertices have opposite parities).\COMMENT{To correct this parity imbalance, a feature would have to be incorporated into the algorithm that identified another molecule with (say) 1 more even than odd vertex to be covered by $\mathfrak{H}$ and give instructions about how $\mathfrak{H}$ should pass between and cover these molecules (correcting the parity imbalance in the process).
}
To avoid this, we assign absorbing structures in pairs, so that, for each $C \in \cC$, either two or no clones of $C$ will be used in absorbing structures.
In the case where two clones are used, we enforce that the tips of these clones will have opposite parities, and therefore each molecule $\cM_C$ will have the same number of even and odd parity vertices to be covered by $\mathfrak{H}$.
We use our robust matching lemma (see \cref{lema:robustmatch}) to pair up the clones of absorbing structures in this way.
To connect up different layers of a cube molecule, we will of course need to have suitable edges between these.
Molecules which do not satisfy this requirement are called `bondless' and are removed from $\Gamma$ before the absorption process (so that their vertices are absorbed).

Secondly, another issue related to vertex parities arises from inner tree vertex molecules. Depending on the degree of an inner tree vertex $v \in V(T)$, the skeleton could contain an odd number of vertices from the molecule $\cM_v$ consisting of all clones of $v$.
All vertices in $\cM_v$ outside the skeleton will need to be absorbed.
But since the number of these vertices is odd,
it would be impossible to pair up (in the way described above) the absorbing structures assigned to these vertices.
To fix this issue, we effectively impose that $\mathfrak{H}$ will `go around $T$ twice'. 
That is, the skeleton will trace through every molecule beginning and finishing at the lowest ordered vertex in $\Gamma^*$.
It will then retrace its steps through these molecules in an almost identical way, effectively doubling the size of the skeleton.
This ensures that the skeleton contains an even number of vertices from each molecule, half of them of each parity.

Finally, once we have obtained an appropriate skeleton, we can construct a long cycle $\mathfrak{H}$ as described in \cref{sect:cycle}. 
For every vertex in $\cQ^n$ which is not covered by $\mathfrak{H}$ we have put in place an absorbing structure,  which is covered  by~$\mathfrak{H}$ as described in \ref{psa} and \ref{psb}.
Thus, as discussed before, we can now use these structures to absorb all remaining vertices into $\mathfrak{H}$ to obtain a Hamilton cycle $\mathfrak{H}' \subseteq H\cup G$, thus proving the case $k=1$ of \cref{thm:main}.

\subsection{Hitting time for the appearance of a Hamilton cycle.}\label{section:outline6}

As mentioned in \cref{sect:1.5}, the hitting time result (\cref{thm:hitting}) follows easily from the perturbation type result \cref{thm: kedgehit}
(see \cref{sect:hit} for the details).
Thus, here we sketch a proof of (the $k=1$ case of) \cref{thm: kedgehit}.
Consider $G\sim \cQ^n_{1/2 - \eps}$.
We show that a.a.s., for any graph $H$ with $\delta(H) \ge 2$, the graph $G \cup H$ is Hamiltonian.
The main additional difficulty faced here is that $G \cup H$ may contain vertices having degree as low as $2$.
For the set $\cU$ of these vertices we cannot hope to use the previous absorption strategy: the neighbours of $v \in \cU$ may not lie in cubes from $\cC$. (In fact, $v$ may not even have a neighbour within its own layer in $G\cup H$.)
To handle such small degree vertices, we first prove that they will be few and well spread out (see \cref{lem: badvertices}).
In \cref{sect:absstruct} we define three types of new `special absorbing structures'.
The type of the special absorbing structure $\mathit{SA}(v)$\index{SAv2@$\mathit{SA}(v)$} for $v$ will depend on whether the neighbours $a,b$ of $v$ in $H$ lie in the same layer as $v$.
In each case, $\mathit{SA}(v)$ will consist of a short path $P_1$ containing the edges $av$ and $bv$, and several other short paths designed to `balance out' $P_1$ in a suitable way. 
(This is further discussed in \cref{sect:absstruct}, see \cref{fig:SAS}.)
These paths will be incorporated into the long cycle $\mathfrak{H}$ described in \cref{sect:cycle}. 
In particular, this allows us to `absorb' the vertices of $\cU$ into $\mathfrak{H}$. 
To incorporate the paths $P_i$ forming $\mathit{SA}(v)$, we will proceed as follows.

Firstly, we make use of the fact that \cref{lem: main treereshit} allows us to choose our near-spanning tree $T$ in such a way that it avoids a small ball around each $v \in \cU$.
Thus, (all clones of) $T$ will avoid $\mathit{SA}(v)$, which has the advantage there will be no interference between $T$ and the special absorbing structures. 
To link up each $\mathit{SA}(v)$ with the long cycle $\mathfrak{H}$, for each endpoint $w$ of a path in $\mathit{SA}(v)$, we will choose an $\ell$-cube in $I(G)$ which suitably intersects $T$ and which contains $w$ (or more precisely, the vertex in $I(G)$ corresponding to $w$). 
Altogether, these $\ell$-cubes allow us to find paths between $\mathit{SA}(v)$ and vertices of $\mathfrak{H}$ which are clones of vertices in $T$. 
The remaining vertices in molecules consisting of clones of these $\ell$-cubes will be covered in a similar way as in \cref{sect:cycle}.
All vertices in these balls around $\cU$ which are not part of the special absorbing structures will be absorbed into $\mathfrak{H}$ via the same absorbing structures used in the proof of \cref{thm:main} to once again obtain a Hamilton cycle $\mathfrak{H}'$.

\subsection{Edge-disjoint Hamilton cycles.}
The results on $k$ edge-disjoint Hamilton cycles can be deduced from suitable versions (\cref{thm:main1,thm: important}) of the case $k=1$.
Those versions are carefully formulated to allow us to repeatedly remove a Hamilton cycle from the original graph.
We deduce \cref{thm:thresholdk} from \cref{thm:main1} in \cref{sect:thm1}, and deduce \cref{thm:hitting} from \cref{thm: important} in \cref{sect:hit}.


\section{Notation}\label{notation}

\textbf{Hierarchies and asymptotics}.
For $n\in\mathbb{Z}$, we denote $[n]\coloneqq\{k\in\mathbb{Z}:1\leq k\leq n\}$\index{[n]@$[n]$} and $[n]_0\coloneqq\{k\in\mathbb{Z}:0\leq k\leq n\}$\index{[n]0@$[n]_0$}.
Whenever we write a hierarchy of parameters, these are chosen from right to left.
That is, whenever we claim that a result holds for $0 < a \ll b \le 1$\index{\ll@$\ll$}, we mean that there exists a non-decreasing function $f\colon [0, 1) \to [0,1)$ such that the result holds for all $a>0$ and all $b \le 1$ with $a \le f(b)$.
We will not compute these functions explicitly.
Hierarchies with more constants are defined in a similar way.
When considering random experiments for a sequence of graphs $(G_n)_{n\in\mathbb{N}}$ with $|V(G_n)|$ tending to infinity with $n$, we say that an event $\mathcal{E}$ holds \emph{asymptotically almost surely} (\emph{a.a.s.}) for $G_n$ if $\mathbb{P}[\mathcal{E}]=1-o(1)$.
When considering asymptotic statements, we will ignore rounding whenever this does not affect the argument.

\textbf{Hypergraphs}.
A \emph{hypergraph} $H$ is an ordered pair $H=(V(H),E(H))$\index{EH1@$E(H)$} where $V(H)$ is called the vertex set and $E(H)\subseteq2^{V(H)}$, the edge set, is a set of subsets of $V(H)$.
If $E(H)$ is a multiset, we refer to $H$ as a \emph{multihypergraph}.
We say that a (multi)hypergraph $H$ is $r$\emph{-uniform} if for every $e\in E(H)$ we have $|e|=r$.
In particular, $2$-uniform hypergraphs are simply called \emph{graphs}.
Given any set of vertices $V'\subseteq V(H)$, we denote the subhypergraph of $H$ \emph{induced} by $V'$ as $H[V']\coloneqq(V',E')$\index{HV'@$H[V']$}, where $E'\coloneqq\{e\in E(H):e\subseteq V'\}$.
We write $H-V'\coloneqq H[V\setminus V']$\index{HV''@$H-V'$}.
Given any set $\hat E\subseteq E(H)$, we will sometimes write $V(\hat E)\coloneqq\{v\in V:\text{there exists } e\in \hat E\text{ such that }v\in e\}$\index{VE@$V(\hat E)$}.

\textbf{Neighbourhoods and degrees}.
Given any (multi)hypergraph $H$ and any vertex $v\in V(H)$, let $E(H,v)\coloneqq\{e\in E(H):v\in e\}$\index{EH2@$E(H,u)$}.
We define the \emph{neighbourhood} of $v$ as $N_H(v)\coloneqq\bigcup_{e\in E(H,v)}e\setminus\{v\}$\index{NG1v@$N_G(v)$}, and we define the \emph{degree} of $v$ by $d_H(v)\coloneqq|E(H,v)|$\index{deg1H1@$d_G(u)$}.
We denote the minimum and maximum degrees of (the vertices in) $H$ by $\delta(H)$\index{deltaH@$\delta(H)$, $\Delta(H)$} and $\Delta(H)$, respectively.
Given any pair of vertices $u,v\in V(H)$, we define $E(H,u,v)\coloneqq\{e\in E(H):\{u,v\}\subseteq e\}$\index{EH3@$E(H,u,v)$}.
The \emph{codegree} of $u$ and $v$ in $H$ is given by $d_H(u,v)\coloneqq|E(H,u,v)|$\index{deg1H2@$d_G(u,v)$}.
Given any set of vertices $W\subseteq V(H)$, we define $N_H(W)\coloneqq \bigcup_{w \in W} N_H(w)$\index{NG1W@$N_G(W)$}. 
We denote $E(H,v,W)\coloneqq\{e\in E(H):v\in e, e\setminus\{v\}\subseteq W\}$\index{EH4@$E(H,u,W)$}, $N_H(v,W)\coloneqq\bigcup_{e\in E(H,v,W)}e\setminus\{v\}$\index{NG1Wv2@$N_G(v,W)$} and $d_H(v,W)\coloneqq|E(H,v,W)|$\index{deg1H3@$d_G(u,W)$}; we refer to the latter two as the neighbourhood and degree of $v$ into $W$, respectively.
Given $A, B \subseteq V(H)$ we denote $E_H(A,B)\coloneqq \{e \in E(H): e\subseteq A\cup B, e \cap A \ne \varnothing, e\cap B \ne \varnothing\}$\index{EHA@$E_H(A,B)$} and $e_H(A,B)\coloneqq|E_H(A,B)|$\index{eHA@$e_H(A,B)$}.
Whenever $A=\{v\}$ is a singleton, we abuse notation and write $E_H(v,B)$ and $e_H(v,B)$.
Thus, $e_H(v,B)$ and $d_H(v,B)$ may be used interchangeably.

\textbf{Distances}.
Given any graph $G$ and two vertices $u,v\in V(G)$, the \emph{distance} $\dist_G(u,v)$\index{dist1@$\dist_G(\cdot,\cdot)$} between $u$ and $v$ in $G$ is defined as the length of the shortest path connecting $u$ and $v$ (and it is said to be infinite if there is no such path).
Similarly, given any sets $A, B \subseteq V(G)$, the \emph{distance} between $A$ and $B$ is given by $\dist_G(A,B)\coloneqq\min_{u\in A, v\in B}{\dist_G(u,v)}$.
For any $r\in\mathbb{N}$, we denote $B_G^r(u)\coloneqq\{v\in V(G):\dist_G(u,v)\leq r\}$\COMMENT{When dealing with cubes, we may equivalently let $B_G^r(u)\coloneqq u+r(\mathcal{D}(\cQ^n)\cup\{\mathbf{0}\})$.}\index{ball@$B_G^r(\cdot)$} and $B_G^r(A)\coloneqq\{v\in V(G):\dist_G(A,v)\leq r\}$; we refer to these sets as the \emph{balls} of radius $r$ around $u$ and $A$, respectively. 

\textbf{Digraphs}.
A \emph{directed graph} (or \emph{digraph}) is a pair $D=(V(D), E(D))$, where $E(D)$ is a set of ordered pairs of elements of $V(D)$.
If no pair of the form $(v,v)$ with $v\in V(D)$ belongs to $E(D)$, we say that $D$ is \emph{loopless}.
Given any $v\in V(D)$, we define its \emph{inneighbourhood} as $N_D^-(v)\coloneqq\{u\in V(D):(u,v)\in E(D)\}$\index{ND@$N_D^-(v)$, $N_D^+(v)$}, and its \emph{outneighbourhood} as $N_D^+(v)\coloneqq\{u\in V(D):(v,u)\in E(D)\}$.
The \emph{indegree} and \emph{outdegree} of $v$ are defined as $d_D^-(v)\coloneqq|N_D^-(v)|$\index{deg1D@$d_D^-(v)$, $d_D^+(v)$} and $d_D^+(v)\coloneqq|N_D^+(v)|$, respectively.
The minimum in- and outdegrees of (the vertices in) $D$ are denoted by $\delta^-(D)$\index{deltaD@$\delta^-(D)$, $\delta^+(D)$} and $\delta^+(D)$, respectively.

\textbf{Matchings}.
Given any multihypergraph or directed graph $(V,E)$, a set $M\subseteq E$ is called a \emph{matching} if its elements are pairwise disjoint.
If the edges of $M$ cover all of $V$, then it is said to be a \emph{perfect matching}.
Given an edge-colouring $c$ of $H$, we say that a matching of $H$ is \emph{rainbow} if each of its edges has a different colour in $c$.

\textbf{Hypercubes}.
We often refer to the $n$-dimensional hypercube $\cQ^n$\index{Qn@$\cQ^n$} as an $n$-\emph{cube} (the $n$ is dropped whenever clear from the context).
Given two vertices $v_1,v_2\in V(\cQ^n)=\{0,1\}^n$, we write $\dist(v_1,v_2)$\index{dist@$\dist(\cdot,\cdot)$} for the Hamming distance between $v_1$ and $v_2$.
Thus, $\{v_1,v_2\}\in E(\cQ^n)$ if and only if $\dist(v_1,v_2)=1$.
Whenever the dimension $n$ is clear from the context, we will use $\mathbf{0}$ to denote the vertex $\{0\}^n$.
Given any $v\in\{0,1\}^n$, we will say that its \emph{parity} is \emph{even} if $\dist(v,\mathbf{0})\equiv0\pmod2$, and we will say that it is \emph{odd} otherwise.
This gives a natural partition of $V(\cQ^n)$ into the sets of vertices with even and odd parities.
Given any two vertices $v_1,v_2\in\{0,1\}^n$, we will write $v_1\eqp v_2$\index{\eq@$\eqp$, $\neqp$} if they have the same parity, and $v_1\neqp v_2$ otherwise.

We will often consider the natural embedding of $V(\cQ^n)$ into $\mathbb{F}_2^n$, which will allow us to use operations on the vertex set: whenever we write $v+u$, for some $u,v\in\{0,1\}^n$, we refer to their sum in $\mathbb{F}_2^n$\index{+@$\cdot+\cdot$}.
Given a vertex $v\in\{0,1\}^n$ and an edge $e=\{x,y\}\in E(\cQ^n)$, we define $v+e$ to be the edge with endvertices $v+x$ and $v+y$.
Given any two sets $A,B\subseteq\{0,1\}^n$, we will use the sumset notation $A+B\coloneqq\{a+b: a\in A,b\in B\}$, and we will abbreviate the $k$-fold sumset $A+\ldots+A$ by $kA$.
Similarly, given any sets $A\subseteq\{0,1\}^n$ and $E\subseteq E(\cQ^n)$, we write $A+E\coloneqq\{a+e: a\in A,e\in E\}$.
Given a graph $G\subseteq\cQ^n$ and a set of vertices $A\subseteq\{0,1\}^n$, $A+G$ will denote the graph with vertex set $A+V(G)$ and edge set $A+E(G)$.
Note that this should never be confused with the notation $G-A$\COMMENT{In sumset notation, since $A\subseteq\mathbb{F}_2^n$, we have that $-A=A$, hence $A+G=G+A=G+(-A)=G-A$, so it should be made clear that we never mean this.}, which will be used exclusively to consider induced subgraphs of $G$.
We will call the unitary vectors in $\mathbb{F}_2^n$ the \emph{directions} of the hypercube.
The set of directions will be denoted by $\mathcal{D}(\cQ^n)$\index{Dir@$\mathcal{D}(\cdot)$}.
Thus, $\mathcal{D}(\cQ^n)=\{\hat e\in\{0,1\}^n:\dist(\hat e,\mathbf{0})=1\}$.
Note that two vertices $v_1,v_2\in\{0,1\}^n$ are adjacent in $\cQ^n$ if and only if there exists $\hat e\in\cD(\cQ^n)$ such that $v_1=v_2+\hat{e}$.
Given any vertex $v\in\{0,1\}^n$ and any set $\mathcal{D}\subseteq\mathcal{D}(\cQ^n)$, we will denote by $\cQ^n(v,\mathcal{D})\coloneqq\cQ^n[v+n(\mathcal{D}\cup\{\mathbf{0}\})]$\index{Qnv@$\cQ^n(v,\mathcal{D})$} the subcube of $\cQ^n$ which contains $v$ and all vertices in $\{0,1\}^n$ which can be reached from $v$ by only adding directions in $\mathcal{D}$.
Given any subcube $Q\subseteq\cQ^n$, we will write $\mathcal{D}(Q)$ to denote the subset of $\mathcal{D}(\cQ^n)$ such that, for any $v\in V(Q)$, we have $Q=\cQ^n(v,\mathcal{D}(Q))$.
Given any direction $\hat e\in\mathcal{D}(Q)$, we will sometimes informally say that $Q$ \emph{uses} $\hat e$\COMMENT{That is, there exist two vertices $x,y\in V(Q)$ such that $x-y=\hat e$.}.
Given two vertices $v_1,v_2\in\{0,1\}^n$, their \emph{differing directions} are all directions in $\mathcal{D}(v_1,v_2)\coloneqq\{\hat{e}\in\mathcal{D}(\cQ^n):\dist(v_1+\hat{e},v_2)<\dist(v_1,v_2)\}$\index{Dirv1v2@$\mathcal{D}(v_1,v_2)$}.
Observe that, if $\dist(v_1,v_2)=d$, then $|\mathcal{D}(v_1,v_2)|=d$ and $\cQ^n(v_1,\mathcal{D}(v_1,v_2))$ is the smallest subcube of $\cQ^n$ which contains both $v_1$ and~$v_2$.


\section{Probabilistic tools}\label{prelim}

Here we list some probabilistic tools that we will use throughout the paper.
The following can be proved easily with the Cauchy-Schwarz inequality.

\begin{proposition}\label{prop: CZ}
Given a non-negative random variable $X$ with finite support, we have that
\begin{equation*}
    \mathbb{P}[X=0] \le 1 - \frac{\mathbb{E}[X]^2}{\mathbb{E}[X^2]}.
\end{equation*}
\end{proposition}
\COMMENT{By rearranging the inequality, we have that
\[\mathbb{P}[X=0] \le 1 - \frac{\mathbb{E}[X]^2}{\mathbb{E}[X^2]}\Longleftrightarrow\mathbb{E}[X]^2\leq\mathbb{E}[X^2](1-\mathbb{P}[X=0]).\]
Assume that the support of $X$ has size $N+1$, with values $\{a_0,a_1,\ldots,a_N\}$, with $a_0 =0$, each happening with probability $p_i$.
Then, the inequality above, by definition, is equivalent to
\[\left(\sum_{i=1}^Na_ip_i\right)^2\leq\left(\sum_{i=1}^Na_i^2p_i\right)\left(\sum_{i=1}^Np_i\right).\]
And this holds by Cauchy-Schwarz.}

Throughout the paper, we will be interested in proving concentration results for different random variables.
We will often need the following Chernoff bound (see e.g.~\cite[Corollary~2.3]{JLR}).

\begin{lemma}\label{lem:Chernoff}
Let $X$ be the sum of\/ $n$ mutually independent Bernoulli random variables and let $\mu \coloneqq \mathbb{E}[X]$.
Then, for all $0<\delta<1$ we have that $\mathbb{P}[X\geq(1+\delta)\mu]\leq e^{-\delta^2\mu/3}$ and $\mathbb{P}[X\leq(1-\delta)\mu]\leq e^{-\delta^2\mu/2}$.
In particular, $\mathbb{P}[|X-\mu|\geq\delta\mu]\leq2e^{-\delta^2\mu/3}$.
\end{lemma}

Similar bounds hold for hypergeometric distributions (see e.g.~\cite[Theorem~2.10]{JLR}).
For $m, n, N \in \mathbb{N}$ with $m, n < N$, a random variable $X$ is said to follow the hypergeometric distribution with parameters $N$, $n$ and $m$ if it can be defined as $X \coloneqq |S \cap [m]|$, where $S$ is a uniformly chosen random subset of $[N]$ of size $n$.

\begin{lemma}\label{lem:ChernoffHyp}
  Suppose $Y$ has a hypergeometric distribution with parameters $N$, $n$ and $m$.
  Then, $\mathbb{P}[|Y - \mathbb{E}[Y ]| \ge t] \le 2e^{-t^2/(3n)}$.\COMMENT{Follows by the remark after \cite[Theorem~2.10]{JLR} that binomial Chernoff implies hypergeometric bounds.
  (The bound in \cref{lem:ChernoffHyp} is weaker than the one in \cref{lem:Chernoff} as $\mu \le n$.)}
\end{lemma}

The following bound will also be used repeatedly (see e.g.~\cite[Theorem A.1.12]{AS16})\COMMENT{The theorem in \emph{The Probabilistic Method} reads as follows: 
If $X=X_i+\ldots+X_n$, where the $X_i$ are Bernoulli random variables and $\mathbb{E}[X_i]=p_i$, and $\beta>1$, then $\mathbb{P}[X\geq\beta np]\leq(e^{\beta-1}\beta^{-\beta})^{np}$.
It is easy to see that \cref{lem:betaChernoff} follows from this:
\[\mathbb{P}[X\geq\beta np]\leq(e^{\beta-1}\beta^{-\beta})^{np}=\left(\frac{e}{\beta}\right)^{\beta np}\cdot\frac{1}{e^{np}}\leq\left(\frac{e}{\beta}\right)^{\beta np}.\]}.

\begin{lemma}\label{lem:betaChernoff}
Let $X$ be the sum of\/ $n$ mutually independent Bernoulli random variables.
Let $\mu \coloneqq \mathbb{E}[X]$, and let $\beta>1$.
Then, $\mathbb{P}[X\geq\beta\mu]\leq\left(e/\beta\right)^{\beta\mu}$.
In particular, we have $\mathbb{P}[X\geq 7 \mu]\leq  e^{-\mu}$.\COMMENT{Alberto: Isn't this already true for $4\mu$? Why did we write $7$?}
\end{lemma}

Given any sequence of random variables $X=(X_1,\ldots,X_n)$ taking values in a set $\Omega$ and a function $f\colon \Omega^n\to\mathbb{R}$, for each $i\in[n]_0$ define $Y_i\coloneqq\mathbb{E}[f(X)\mid X_1,\ldots,X_i]$.
The sequence $Y_0,\ldots,Y_n$ is called the \emph{Doob martingale} for $f$ and $X$.
All the martingales that appear in this paper will be of this form.
To deal with them, we will need the following version of the well-known Azuma-Hoeffding inequality.

\begin{lemma}[Azuma's inequality \cite{Azu67,Hoef63}]\label{lem: Azuma}
Let $X_0,X_1,\ldots$ be a martingale and suppose that $|X_i-X_{i-1}|\leq c_i$ for all $i\in\mathbb{N}$.
Then, for any $n,t\in\mathbb{N}$,
\[\mathbb{P}[|X_n-X_0|\geq t]\leq2\exp\left(\frac{-t^2}{2\sum_{i=1}^nc_i^2}\right).\]
\end{lemma}

The following lemma, which concerns further bounds for martingales, is due to \citet{AKS97} (see also \cite[Theorem~7.4.3]{AS16}).
Here, we describe a version which is tailored to our purposes.
Let $r\in\mathbb{N}$ and let $\mathcal{H}$ be an $r$-uniform hypergraph.
Let $\mathcal{H}'\subseteq\mathcal{H}$ be a random subgraph chosen according to any distribution for which the inclusion of edges are mutually independent.
Let $X$ be a random variable whose value is determined by the presence or absence of the edges of some collection $E'=\{e_1, \dots, e_k\} \subseteq E(\mathcal{H})$ in $\mathcal{H}'$.
Let $p_i$ be the probability that $e_i$ is present in $\mathcal{H}'$.
Let $c_i$ be the maximum value $X$ could change, for some given choice of $\mathcal{H}'$, by changing the presence or absence of $e_i$.
Let $C \coloneqq \max_{i\in[k]}{c_i}$ and $\sigma^2 \coloneqq \sum_{i\in[k]} p_i(1-p_i)c_i^2$.

\begin{lemma}[\citet{AKS97}]\label{lem: AKS}
For all $\alpha>0$ with $\alpha C < 2\sigma$ we have that
\[\mathbb{P}[|X - \mathbb{E}[X]| > \alpha \sigma] \le 2e^{-\alpha^2/4}. \]
\end{lemma}

We will also need the following special case of Talagrand's inequality (see e.g.~\cite[Theorem~7.7.1]{AS16}).
Let $\Omega\coloneqq\prod_{i=1}^{n}\Omega_i$, where each $\Omega_i$ is a probability space.
We say that $f\colon\Omega\to\mathbb{R}$ is \emph{$K$-Lipschitz}, for some $K\in\mathbb{R}$, if for every $x,y\in\Omega$ which differ only on one coordinate we have $|f(x)-f(y)|\leq K$.
We say that $f$ is \emph{$h$-certifiable}, for some $h\colon\mathbb{N}\to\mathbb{N}$, if, for every $x\in\Omega$ and $s\in\mathbb{R}$, whenever $f(x)\geq s$, there exists $I\subseteq[n]$ with $|I|\leq h(s)$ such that every $y\in\Omega$ that agrees with $x$ on the coordinates in $I$ satisfies $f(y)\geq s$.

\begin{lemma}[Talagrand's inequality]\label{lem: Talagrand}
Let $\Omega\coloneqq\prod_{i=1}^{n}\Omega_i$, where each $\Omega_i$ is a probability space. 
Let $X\colon\Omega\to\mathbb{N}$ be $K$-Lipschitz and $h$-certifiable, for some $K\in\mathbb{N}$ and $h\colon\mathbb{N}\to\mathbb{N}$.
Then, for all $b,t\in\mathbb{R}$,
\[\mathbb{P}\left [X\leq b-tK\sqrt{h(b)}\right ]\mathbb{P}[X\geq b] \leq \exp\left ( \frac{-t^2}{4}\right).\]
\end{lemma}

Finally, the Lov\'asz local lemma will come in useful. 
Let $\mathfrak{E} \coloneqq \{\mathcal{E}_1, \mathcal{E}_2, \ldots, \mathcal{E}_m\}$ be a collection of events. 
A \emph{dependency graph} for $\mathfrak{E}$ is a graph $H$ on vertex set $[m]$ such that, for all $i\in[m]$, $\mathcal{E}_i$ is mutually independent of $\{\mathcal{E}_j:j\neq i, j\notin N_H(i)\}$, that is, if $\mathbb{P}[\mathcal{E}_i]=\mathbb{P}[\mathcal{E}_i\mid\bigwedge_{j\in J}\mathcal{E}_j]$ for all $J\subseteq[m]\setminus(N_H(i)\cup\{i\})$.
We will use the following version of the local lemma (it follows e.g.~from~\cite[Lemma 5.1.1]{AS16}).

\begin{lemma}[Lov\'asz local lemma]\label{lem: LLL}
Let $\mathfrak{E} \coloneqq \{\mathcal{E}_1, \mathcal{E}_2, \ldots, \mathcal{E}_m\}$ be a collection of events and let $H$ be a dependency graph for $\mathfrak{E}$.
Suppose that $\Delta(H)\leq d$ and $\mathbb{P}[\mathcal{E}_i]\leq p$ for all $i\in[m]$.
If $ep(d+1)\leq1$, then 
\[\mathbb{P}\left[\bigwedge_{i=1}^m\overline{\mathcal{E}_i}\right]\geq(1-ep)^m.\]
\end{lemma}


\section{Auxiliary results}\label{aux}

\subsection{Results about matchings}\label{sect5matchings}

We will need three auxiliary results to help us find suitable absorbing cube pairs for different vertices.
We will need to preserve the alternating parities of vertices that are absorbed by each molecule.
The first lemma (\cref{lema:robustmatch}) presented in this section, as well as its corollary, will help us to show that all vertices can be paired up in such a way that these parities can be preserved.
The second lemma (\cref{prop: cubematch}) will be used to show that, for each such pair of vertices, there are many possible pairs of absorption cubes.
Finally, the third lemma (\cref{lem: rainbow}) will allow us to assign one of those pairs of absorption cubes to each pair of vertices we need to absorb in such a way that these cube pairs are pairwise vertex disjoint\COMMENT{in the projection.}.

To prove \cref{lema:robustmatch}, as well as \cref{lem: randsethit,thm: maintreeres}, the following consequence of Hall's theorem will be useful.

\begin{lemma}\label{lema:Hall}
  Let $G$ be a bipartite graph with vertex partition $A\cupdot B$.
  Assume that there is some integer $\ell\geq0$ such that, for all $S\subseteq A$, we have $|N(S)|\geq|S|-\ell$.
  Then, $G$ contains a matching which covers all but at most $\ell$ vertices in $A$.
\end{lemma}

Given any graph $G$ and a bipartition $(\mathfrak{A},\mathfrak{B})$ of $V(G)$, we say that $(\mathfrak{A},\mathfrak{B})$ is an \emph{$r$-balanced bipartition} if $||\mathfrak{A}|-|\mathfrak{B}||\leq r$.
Let $G$ be a graph on $n$ vertices, and let $r,d\in \mathbb{N}$ with $r\leq d$.
We say that $G$ is \emph{$d$-robust-parity-matchable} with respect to an $r$-balanced bipartition $(\mathfrak{A},\mathfrak{B})$ if, for every $S\subseteq V(G)$ such that $|S|\leq d$ and $|\mathfrak{A}\setminus S|=|\mathfrak{B}\setminus S|$, the graph $G-S$ contains a perfect matching~$M$ with the property that every edge $e\in M$ has one endpoint in $\mathfrak{A}\setminus S$ and one endpoint in $\mathfrak{B}\setminus S$.

Given two disjoint sets of vertices $A$ and $B$, the binomial random bipartite graph $G(A,B,p)$ is obtained by adding each possible edge with one endpoint in $A$ and the other in $B$ with probability $p$ independently of every other edge.
Given any two bipartite graphs on the same vertex set, $G_1=(A,B,E_1)$ and $G_2=(A,B,E_2)$, and any $\alpha\in\mathbb{R}$, we define $\Gamma_{G_1,G_2}^\alpha(A)$\index{Gamma@$\Gamma_{G_1,G_2}^\alpha(A)$} as the graph with vertex set $A$ where any two vertices $x,y\in A$ are joined by an edge whenever $|N_{G_1}(x)\cap N_{G_2}(y)|\geq\alpha|B|$ or $|N_{G_1}(y)\cap N_{G_2}(x)|\geq\alpha|B|$.

\begin{lemma}\label{lema:robustmatch}
  Let $d, k, r \in\mathbb{N}$ and $\alpha,\varepsilon,\beta >0$ be such that $r\leq d$, $1/k\ll 1/d,\varepsilon,\alpha$ and $\beta\ll \varepsilon,\alpha$. 
  Then, any bipartite graph $G=G(A,B,E)$ with $|B|=n\geq|A|\geq k$ such that $d_G(x)\geq\alpha n$ for every $x\in A$ satisfies the following with probability at least $1-2^{-10n}$:
  for any $r$-balanced bipartition of $A$ into $(\mathfrak{A},\mathfrak{B})$, the graph $\Gamma_{G,G(A,B,\varepsilon)}^\beta(A)$ is $d$-robust-parity-matchable with respect to $(\mathfrak{A},\mathfrak{B})$.
\end{lemma}

\begin{proof}
Let $\Gamma\coloneqq\Gamma_{G, G(A,B,\varepsilon)}^\beta(A)$.
Let $\Gamma'$ be the auxiliary digraph with vertex set $A$ where, for any pair of vertices $x,y\in A$, there is a directed edge from $x$ to $y$ if $|N_G(x)\cap N_{G(A,B,\varepsilon)}(y)|\geq\beta n$.
Observe that the graph obtained from $\Gamma'$ by ignoring the directions of its edges and identifying the possible multiple edges is exactly $\Gamma$, which means that
$\delta(\Gamma)\geq\delta^+(\Gamma')$.

Given any two vertices $x,y\in A$, by \cref{lem:Chernoff} we have that\COMMENT{We have that
\[\mathbb{P}[(x,y)\notin E(\Gamma')]=\mathbb{P}[|N_G(x)\cap N_{G(A,B,\varepsilon)}(y)|<\beta n].\]
Since $N_G(x)$ is a fixed set of size $\geq\alpha n$, the random variable $X\coloneqq|N_G(x)\cap N_{G(A,B,\varepsilon)}(y)|$ follows a binomial distribution with parameters $|N_G(x)|$ and $\varepsilon$.
We then have that $\mathbb{E}[X]=\varepsilon|N_G(x)|\geq\varepsilon\alpha n$.
Then, applying \cref{lem:Chernoff}, we conclude that
\begin{align*}
    \mathbb{P}[X<\beta n]&=\mathbb{P}[X<(\beta/(\varepsilon\alpha))\varepsilon\alpha n]\leq\mathbb{P}[X<(\beta/(\varepsilon\alpha))\mathbb{E}[X]]\\
    &\leq e^{-(1-\beta/(\varepsilon\alpha))^2\mathbb{E}[X]/2}\leq e^{-(1-\beta/(\varepsilon\alpha))^2\varepsilon\alpha n/2}\leq e^{-\varepsilon\alpha n/3}.
\end{align*}
}
\[\mathbb{P}[(x,y)\notin E(\Gamma')]=\mathbb{P}[|N_G(x)\cap N_{G(A,B,\varepsilon)}(y)|<\beta n]\leq e^{-\varepsilon\alpha n/3}.\]
Furthermore, for a fixed $x\in A$, observe that the events that $(x,y)\notin E(\Gamma')$, for all $y\in A\setminus\{x\}$, are mutually independent.
Therefore, $d^+_{\Gamma'}(x)$ is a sum of independent Bernoulli random variables.
Let $m\coloneqq|A|$.
If $d^+_{\Gamma'}(x)<4m/5$, that means that there is a set of $m/5$ vertices $Y\subseteq A\setminus\{x\}$ such that $(x,y)\notin E(\Gamma')$ for all $y\in Y$.
We then conclude that\COMMENT{We have that $\binom{m}{m/5}\leq2^m$.
Therefore,
\begin{align*}
    \mathbb{P}[d^+_{\Gamma'}(x)<4m/5]&\leq\sum_{Y\in\binom{A\setminus\{x\}}{m/5}}\mathbb{P}[(x,y)\notin E(\Gamma')\text{ for all }y\in Y]\leq\binom{m}{m/5}\left(e^{-\varepsilon\alpha n/3}\right)^{m/5}\\
    &\leq 2^{-\varepsilon\alpha nm/(15\log2)+m}\leq 2^{-\varepsilon\alpha nm/11}\leq 2^{-20n},
\end{align*}
where in the last and second-to-last inequalities we use the fact that $n\geq m\geq k\gg (\varepsilon\alpha)^{-1}$.
}
\[\mathbb{P}[d^+_{\Gamma'}(x)<4m/5]\leq\sum_{Y\in\binom{A\setminus\{x\}}{m/5}}\mathbb{P}[(x,y)\notin E(\Gamma')\text{ for all }y\in Y]\leq\binom{m}{m/5}e^{-\varepsilon\alpha nm/15}\leq 2^{-20n}.\]
By a union bound over the choice of $x$, we conclude that 
\[\mathbb{P}[\delta(\Gamma)<4m/5]\leq\mathbb{P}[\delta^+(\Gamma')<4m/5]\leq m2^{-20n}\leq 2^{-10n}.\]

Now, condition on the event that the previous holds.
Fix any $r$-balanced bipartition $(\mathfrak{A},\mathfrak{B})$ of $A$ and let $\Gamma_{(\mathfrak{A},\mathfrak{B})}$ be the bipartite subgraph of $\Gamma$ induced by this bipartition.
Fix any set $S\subseteq A$ with $|S|\leq d$ and $|\mathfrak{A}\setminus S|=|\mathfrak{B}\setminus S|$.
We have that $\delta(\Gamma_{(\mathfrak{A},\mathfrak{B})} - S)\geq4m/5-m/2-d-r>m/4$\COMMENT{Here we use that $m\gg d\geq r$.}.
Therefore, by \cref{lema:Hall}, $\Gamma_{(\mathfrak{A},\mathfrak{B})}-S$ contains a perfect matching\COMMENT{Indeed, for any set $C\subseteq\mathfrak{A}$ with $|C|\leq m/4$, the degree condition trivially implies that $|N_{\Gamma_{(\mathfrak{A},\mathfrak{B})}-S}(C)|\geq|C|$, and for any $C\subseteq\mathfrak{A}$ with $|C|>m/4$ we have that $|\mathfrak{A}\setminus C|<m/4$ and, therefore, $N_{\Gamma_{(\mathfrak{A},\mathfrak{B})}-S}(C)=\mathfrak{B}$.}.
\end{proof}

While \cref{lema:robustmatch} will be used in the proof of \cref{thm:main1} in \cref{section8}, we will instead need to use the following \cref{cor: robustmatch} in the proof of \cref{thm: important} in \cref{sect:hitting}.
Let $G$ be a graph on $2n$ vertices.
Let $(\mathfrak{A},\mathfrak{B})$ be a balanced bipartition of $V(G)$ and let  $(V_1,\ldots,V_k)$, for some $k\in \mathbb{N}$, be a partition of $V(G)$.
Given any $d\in\mathbb{N}$, we say that $G$ is $d$-robust-parity-matchable with respect to $(\mathfrak{A},\mathfrak{B})$ \emph{clustered in} $(V_1,\ldots,V_k)$ if, for every $S\subseteq V(G)$ with $|S|\leq d$ and $|S\cap\mathfrak{A}|=|S\cap\mathfrak{B}|$, the graph $G-S$ contains a perfect matching $M$ such that every edge $e\in M$ has one endpoint in $\mathfrak{A}\setminus S$ and one endpoint in $\mathfrak{B}\setminus S$ and, for every $e=\{x,y\}\in M$, if $x \in V_j$ then $y\in V_{j-1}\cup V_{j}\cup V_{j+1}$ (where we take indices cyclically).

\begin{corollary}\label{cor: robustmatch}
Let $d, k, t \in\mathbb{N}$ and $\alpha,\varepsilon, \beta >0$ be such that $1/k\ll 1/d,\varepsilon,\alpha$ and $\beta\ll \varepsilon,\alpha$.
Let $G=G(A,B,E)$ be a bipartite graph and $(A_1,\ldots,A_t)$ be a partition of $A$ such that
\begin{itemize}
    \item $|B|=n\geq|A|$,
    \item for every $i\in [t]$, we have that $|A_i|\geq k$ is even,
    \item for every $x\in A$, we have $d_G(x)\geq\alpha n$.
\end{itemize} 
Then, the following holds with probability at least $1-2^{-9n}$:
for each $i\in [t]$ and for any balanced bipartition of $A_i$ into $(\mathfrak{A}_i,\mathfrak{B}_i)$, the graph $\Gamma_{G, G(A,B,\varepsilon)}^\beta(A)$ is $d$-robust-parity-matchable with respect to $(\bigcup_{i=1}^{t}\mathfrak{A}_i,\bigcup_{i=1}^{t}\mathfrak{B}_i)$ clustered in $(A_1,\ldots,A_t)$.
\end{corollary}

\begin{proof}
Given any set $C\subseteq A$, for each $i\in[t]$, let $C_i\coloneqq C\cap A_i$.
Given any bipartition $(\mathfrak{A},\mathfrak{B})$ of $A$, we write $C^{\mathfrak{A}}\coloneqq C\cap\mathfrak{A}$ and $C^{\mathfrak{B}}\coloneqq C\cap\mathfrak{B}$.
Throughout this proof, we consider the indices in $[t]$ to be taken cyclically.

For each set $D\subseteq A$ with $|D|\leq d$, and for each $i\in [t]$, we apply \cref{lema:robustmatch} to the graph $G[D_i\cup A_{i-1},B]$, with $2d$, $k$, $d$, $\alpha$, $\varepsilon$ and $\beta$ playing the roles of $d$, $k$, $r$, $\alpha$, $\varepsilon$ and $\beta$, respectively. 
Then, by a union bound over all choices of $D$ and all choices of $i\in[t]$\COMMENT{We are multiplying by at most $t\binom{n}{d}$, and observe that $t\leq n/k$ and $d\ll k< n$, so $t\binom{n}{d}\ll 2^n$.}, the following holds with probability at least $1-2^{-9n}$.
For each $i\in[t]$, consider any balanced bipartition $(\mathfrak{A}_i,\mathfrak{B}_i)$ of $A_i$.
Consider any $D\subseteq A$ with $|D_i|\leq d$ for each $i\in [t]$.
Then, for each $i\in [t]$, the graph $\Gamma^{\beta}_{G[D_i\cup A_{i-1},B], G(D_i\cup A_{i-1},B,\varepsilon)}(D_i\cup A_{i-1})$ is $2d$-robust-parity-matchable with respect to $((\mathfrak{A}_i\cap D)\cup \mathfrak{A}_{i-1},(\mathfrak{B}_i\cap D) \cup \mathfrak{B}_{i-1})$.
Condition on the event that the above holds.

Now, for each $i\in[t]$, fix a balanced bipartition $(\mathfrak{A}_i,\mathfrak{B}_i)$ of $A_i$.
Let $\mathfrak{A}\coloneqq\bigcup_{i=1}^{t}\mathfrak{A}_i$ and $\mathfrak{B}\coloneqq\bigcup_{i=1}^{t}\mathfrak{B}_i$.
Let $S\subseteq A$ be a subset of size $|S|\leq d$ such that $|S^\mathfrak{A}|=|S^\mathfrak{B}|$.
We want to show that $\Gamma_{G, G(A,B,\varepsilon)}^\beta(A)-S$ contains a perfect matching $M$ such that every edge $e\in M$ has one endpoint in $\mathfrak{A}\setminus S$ and one endpoint in $\mathfrak{B}\setminus S$ and, for every $e=\{x,y\}\in M$, if $x \in A_j$ then $y\in A_{j-1}\cup A_{j}\cup A_{j+1}$.
We begin by proving the following claim.

\begin{claim}\label{robmatchclaim}
There exists a set $D\subseteq A\setminus S$ satisfying the following properties:
\begin{enumerate}[label=$(\mathrm{RM}\arabic*)$]
    \item\label{robmatchclaim1} for every $i\in [t]$ we have $|D_i|\leq d$, and
    \item\label{robmatchclaim2} for every $i\in [t]$ we have $|D_{i+1}^\mathfrak{A}\cup D_i^\mathfrak{B}\cup S_i^\mathfrak{B}|=|D_{i+1}^\mathfrak{B}\cup D_i^\mathfrak{A}\cup S_i^\mathfrak{A}|$.
\end{enumerate}
\end{claim} 

\begin{claimproof}
We will construct one such set $D$ by constructing the sets $D_i\subseteq A_i$ inductively.
We will argue by induction on $i\in[t]$ in decreasing order.
Let $D_t\coloneqq\varnothing$.
Now, suppose that, for some $i\in[t-1]$, we have already constructed the sets $D_j\subseteq A_j$ for all $j\in[t]\setminus[i]$. 
Then, let $D_i\subseteq A_i\setminus S_i$ be a smallest set such that\COMMENT{Think of it this way: we will construct a matching such that the vertices of $D_{i+1}$ are connected to vertices in $A_i$, and the vertices in $D_i$ are connected to vertices in $A_{i-1}$, which avoids $S_i$.
For this to be a perfect matching, we must have a balanced bipartition.
That is, we need that $|\mathfrak{A}_i|+|D_{i+1}^\mathfrak{A}|-|D_i^\mathfrak{A}|-|S_i^\mathfrak{A}|=|\mathfrak{B}_i|+|D_{i+1}^\mathfrak{B}|-|D_i^\mathfrak{B}|-|S_i^\mathfrak{B}|$ (because all the sets are disjoint).
Since $|\mathfrak{A}_i|=|\mathfrak{B}_i|$, this is equivalent to the statement here.}
\begin{equation}\label{equa:robmatchclaim}
    |D_{i+1}^\mathfrak{A}\cup D_i^\mathfrak{B}\cup S_i^\mathfrak{B}|=|D_{i+1}^\mathfrak{B}\cup D_i^\mathfrak{A}\cup S_i^\mathfrak{A}|.
\end{equation}
Observe that either $D_i=D_i^\mathfrak{A}$ or $D_i=D_i^\mathfrak{B}$.
Furthermore, observe that $||D_{i+1}^\mathfrak{A}\cup S_i^\mathfrak{B}|-|D_{i+1}^\mathfrak{B}\cup S_i^\mathfrak{A}||\leq|D_{i+1}\cup S_i|$.
Therefore, there exists a set $D_i\subseteq A_i$ as required with $|D_i|\leq|D_{i+1}\cup S_i|$\COMMENT{Because all the sets in the expression above are disjoint and we are taking a smallest set $D$ which satisfies it.
If, say, we have that $|D_{i+1}^\mathfrak{A}\cup S_i^\mathfrak{B}|\geq|D_{i+1}^\mathfrak{B}\cup S_i^\mathfrak{A}|$, then we must have $|D_i|=|D_i^\mathfrak{A}|=|D_{i+1}^\mathfrak{A}\cup S_i^\mathfrak{B}|-|D_{i+1}^\mathfrak{B}\cup S_i^\mathfrak{A}|$.
Furthermore, such a set must exist since $k\gg d$ and, as we will see, we require $|D_i|\leq d$.}.

In order to prove that this results in a set $D$ which satisfies the required properties, consider the following.
First, by following the induction above, we have that $|D_t|=0$ and $|D_i|\leq|D_{i+1}|+|S_i|$, hence $|D_i|\leq\sum_{j=1}^t|S_j|=|S|\leq d$ for all $i\in[t]$, thus \ref{robmatchclaim1} holds.
On the other hand, \ref{robmatchclaim2} holds by \eqref{equa:robmatchclaim} for all $i\in[t-1]$, so we must prove that it also holds for $i=t$.
But this follows by summing \eqref{equa:robmatchclaim} over all $i\in[t-1]$, and using the fact that $|S^\mathfrak{A}|=|S^\mathfrak{B}|$\COMMENT{We have that 
\begin{align*}
    \sum_{i=1}^{t-1}|D_{i+1}^\mathfrak{A}\cup D_i^\mathfrak{B}\cup S_i^\mathfrak{B}|&=\sum_{i=1}^{t-1}|D_{i+1}^\mathfrak{B}\cup D_i^\mathfrak{A}\cup S_i^\mathfrak{A}|\\
    \iff \sum_{i=1}^{t-1}|D_{i+1}^\mathfrak{A}|+\sum_{i=1}^{t-1}|D_i^\mathfrak{B}|+\sum_{i=1}^{t-1}|S_i^\mathfrak{B}|&=\sum_{i=1}^{t-1}|D_{i+1}^\mathfrak{B}|+\sum_{i=1}^{t-1}|D_i^\mathfrak{A}|+\sum_{i=1}^{t-1}|S_i^\mathfrak{A}|\\
    \iff \sum_{i=1}^{t-1}|D_{i+1}^\mathfrak{A}|-\sum_{i=1}^{t-1}|D_i^\mathfrak{A}|-\sum_{i=1}^{t-1}|S_i^\mathfrak{A}|&=\sum_{i=1}^{t-1}|D_{i+1}^\mathfrak{B}|-\sum_{i=1}^{t-1}|D_i^\mathfrak{B}|-\sum_{i=1}^{t-1}|S_i^\mathfrak{B}|\\
    \iff |D_t^\mathfrak{A}|-|D_1^\mathfrak{A}|-(|S^\mathfrak{A}|-|S_t^\mathfrak{A}|)&=|D_t^\mathfrak{B}|-|D_1^\mathfrak{B}|-(|S^\mathfrak{B}|-|S_t^\mathfrak{B}|)\\
    \iff |D_t^\mathfrak{A}|-|D_1^\mathfrak{A}|+|S_t^\mathfrak{A}|&=|D_t^\mathfrak{B}|-|D_1^\mathfrak{B}|+|S_t^\mathfrak{B}|,
\end{align*}
where in the last inequality we use the fact that $|S^\mathfrak{A}|=|S^\mathfrak{B}|$ by assumption.
Finally, the last expression is equivalent to $|D_1^\mathfrak{B}|+|D_t^\mathfrak{A}|+|S_t^\mathfrak{A}|=|D_1^\mathfrak{A}|+|D_t^\mathfrak{B}|+|S_t^\mathfrak{B}|$, which is what we wanted to prove.}.
\end{claimproof}

Let $D$ be the set given by \cref{robmatchclaim}.
Now, for each $i\in[t]$, let $J_i\coloneqq D_i\cup S_i$.
By \cref{robmatchclaim} \ref{robmatchclaim1} we have that $|J_i|\leq 2d$.
Furthermore, by \cref{robmatchclaim} \ref{robmatchclaim2} it follows that $|(\mathfrak{A}_i\cup D_{i+1}^\mathfrak{A})\setminus J_i|=|(\mathfrak{B}_i\cup D_{i+1}^\mathfrak{B})\setminus J_i|$.
By the conditioning above, this means that $\Gamma^{\beta}_{G[D_{i+1}\cup A_i,B], G(D_{i+1}\cup A_i,B,\varepsilon)}(D_{i+1}\cup A_i)-J_i$ contains a perfect matching $M_i$ such that every edge of $M_i$ has one endpoint in $(\mathfrak{A}_i\cup D_{i+1}^\mathfrak{A})\setminus J_i$ and one endpoint in $(\mathfrak{B}_i\cup D_{i+1}^\mathfrak{B})\setminus J_i$.
Finally, let $M\coloneqq\bigcup_{i=1}^tM_i$.
It is clear that $M$ satisfies the required conditions.
The statement follows\COMMENT{Everything after the conditioning works for any arbitrary bipartition and set $S$ as described in the statement.}.
\end{proof}

The second lemma will be stated in terms of directed graphs.

\begin{lemma}\label{prop: cubematch}
  Let $c,C>0$ and let $\alpha\in(0,1/(1+c/C))$.
  Let $D$ be a loopless $n$-vertex digraph such that
  \begin{enumerate}[label=$(\mathrm{\roman*})$]
      \item\label{item:digraphindeg} for every $A\subseteq V(D)$ with $|A| \ge \alpha n$ we have $\sum_{v\in A}d^-(v)\geq c\alpha n$, and
      \item\label{item:digraphoutdeg} for every $B\subseteq V(D)$ with $|B|\le c\alpha n/C$ we have $\sum_{v\in B}d^+(v)\leq c\alpha n$.
  \end{enumerate}
  Then, $D$ contains a matching $M$ with $|M|>c\alpha n/(2C)$.
\end{lemma}

\begin{proof}
Assume for a contradiction that the largest matching $M$ in $D$ has size $|M|\leq c\alpha n/(2C)$.
Since $\alpha<1/(1+c/C)$, there exists a set $A\subseteq V(D)\setminus V(M)$ with $|A|\ge \alpha n$, and thus, by \ref{item:digraphindeg}, $\sum_{v\in A}d^-(v)\geq c\alpha n$.
Since $M$ is the largest matching, all edges that enter $A$ must come from vertices of $M$ (otherwise, we could add one such edge to $M$, finding a larger matching).
However, by \ref{item:digraphoutdeg} and since $M$ contains at least one edge, the number of edges going out of $V(M)$ is less than $c\alpha n$\COMMENT{Because $|M|$ edges belong to $M$ and, thus, do not leave $M$.}, a contradiction.
\end{proof}

For convenience, we state the third lemma in terms of rainbow matchings in hypergraphs.

\begin{lemma}\label{lem: rainbow}
Let $n,r\in \mathbb{N}$ and let $\mathcal{H}$ be an $n$-edge-coloured $r$-uniform multihypergraph.
Then, for any $m\geq 10$, the following holds.
Suppose ${\mathcal{H}}$ satisfies the following two properties:
\begin{enumerate}[label=$(\mathrm{\roman*})$]
    \item For every $i\in [n]$, there are at least $m$ edges of colour $i$. 
    \item $\Delta({\mathcal{H}})\leq m/(6r)$.
\end{enumerate}
Then, there exists a rainbow matching of size $n$. 
\end{lemma}

\begin{proof}
The idea is to pick a random edge from each colour class and prove that with non-zero probability this results in a rainbow matching. 
First, for each $i\in [n]$, let $M_i$ be a set of $m$ edges of colour $i$. 
We choose an edge from each $M_i$ uniformly at random, independently of the other choices. 
For any $i,j\in[n]$ with $i\neq j$ and for any two edges $e \in M_i$ and $e' \in M_j$ for which $e\cap e'\neq \varnothing$, we denote by $A_{e,e'}$ the event that both $e$ and $e'$ are picked.
We observe that
\[ \mathbb{P}[A_{e,e'}]= \left(\frac{1}{m}\right)^2.\]
Moreover, note that every event $A_{e,e'}$ is independent of all other events $A_{f,f'}$ but at most $2 m\cdot r\cdot \Delta(\cH)\leq m^2/3$. 
Indeed, this holds because $A_{e,e'}$ can only depend on those events which involve at least one edge from either colour $i$ or colour $j$. 
Applying now \cref{lem: LLL}, we deduce that with non-zero probability no event $A_{e,e'}$ occurs, as required. 
\end{proof}


\subsection{Properties of random subgraphs of the hypercube}

In this section we state and prove some basic properties of random subgraphs of the hypercube.
The first one guarantees that the degrees of all vertices are linear in the dimension.

\begin{lemma}\label{lem:mindegQeps}
Let $0<\delta \ll\varepsilon\leq1/2$.
Then, we a.a.s.~have that $\delta(\mathcal{Q}_{1/2+\varepsilon}^n)\geq\delta n$.
\end{lemma}

\begin{proof}
Let $p\coloneqq1/2+\varepsilon$.
Fix any $v\in\{0,1\}^n$.
Throughout this proof, we write $d(v)$ to refer to the degree of $v$ in $\mathcal{Q}_p^n$.

Note that $d(v)$ follows a binomial distribution with parameters $n$ and $p$.
Since $\delta<1/2$, it follows that\COMMENT{We use here the fact that $\binom{n}{i}p^i(1-p)^{n-i}$ is increasing in $i$ for $i\leq n/2$.
Indeed, $p^i(1-p)^{n-i}$ is increasing in $i$ since $p>1-p$, and $\binom{n}{i}$ in increasing for $i\leq n/2$, so their product is increasing in this range.}
\[\mathbb{P}[d(v)\leq\delta n]\leq \delta n\binom{n}{\delta n}p^{\delta n}(1-p)^{(1-\delta)n}.\]
Using the Stirling formula, we conclude that\COMMENT{Because of the well-known bounds 
\[\sqrt{2\pi n}\left(\frac ne\right)^ne^{\frac{1}{12n+1}}\leq n!\leq\sqrt{2\pi n}\left(\frac ne\right)^ne^{\frac{1}{12n}}\]
and the Taylor expansion 
\[e^{-\frac1x}=1-\frac1x+\frac{1}{2x^2}-\frac{1}{6x^3}+\ldots=1-\bigO(x^{-1}),\]
we have that
\[\binom{n}{\delta n}=(1\pm\bigO(n^{-1}))\frac{1}{\sqrt{2\pi\delta(1-\delta)n}}\left(\frac{1}{\delta^\delta(1-\delta)^{1-\delta}}\right)^n.\]
The claim follows by multiplying the terms in the previous expression to this one.}
\[\mathbb{P}[d(v)\leq\delta n]\leq (1+\bigO(n^{-1}))\sqrt{\frac{\delta n}{2\pi(1-\delta)}}\left(\left(\frac{p}{\delta}\right)^\delta\left(\frac{1-p}{1-\delta}\right)^{1-\delta}\right)^n.\]
By the union bound, it now suffices to show that 
\[\left(\frac{p}{\delta}\right)^\delta\left(\frac{1-p}{1-\delta}\right)^{1-\delta}=\left(\frac{1+2\varepsilon}{2\delta}\right)^\delta\left(\frac{1-2\varepsilon}{2(1-\delta)}\right)^{1-\delta}<\frac12,\]
but this follows since $\delta\ll\varepsilon$\COMMENT{We know that $(1/x)^x$ tends to $1$ both as $x\to0^+$ and as $x\to1^-$.
Therefore, \[\lim_{\delta\to0^+}\left(\frac{1+2\varepsilon}{2\delta}\right)^\delta\left(\frac{1-2\varepsilon}{2(1-\delta)}\right)^{1-\delta}=\frac{1-2\varepsilon}{2}<\frac12,\]
and hence there exists some $\delta_0=\delta_0(\varepsilon)$ such that
\[\left(\frac{1+2\varepsilon}{2\delta}\right)^\delta\left(\frac{1-2\varepsilon}{2(1-\delta)}\right)^{1-\delta}<\frac12\]
for all $\delta\leq\delta_0$.}.
\end{proof}

Furthermore, it will be important to show that the number of vertices whose degree deviates from the expected degree is small.

\begin{lemma}\label{lem:degconcentrates}
Let $\varepsilon\in(0,1)$ and $a\in(1/2,1)$.
Let $X$ be the number of vertices $v\in\{0,1\}^n$ for which $d_{\mathcal{Q}_\varepsilon^n}(v)\neq\varepsilon n\pm n^a$.
Then\COMMENT{We can change the $6$ or the exponent of $n$ for anything we need later on.}, \[\mathbb{P}[X\geq e^{-n^{2a-1}/(6\varepsilon)}2^n]\leq2e^{-n^{2a-1}/(6\varepsilon)}.\]
\end{lemma}

\begin{proof}
Throughout this proof, we use $d(v)$ for $d_{\mathcal{Q}_\varepsilon^n}(v)$.
Note that $d(v)$ follows a binomial distribution with parameters $n$ and $\varepsilon$, so $\mathbb{E}[d(v)]=\varepsilon n$ and, by \cref{lem:Chernoff},
\begin{equation*}\label{equa:Pbadvertex}
    \mathbb{P}[d(v)\neq\varepsilon n\pm n^a]\leq2e^{-n^{2a-1}/(3\varepsilon)}.
\end{equation*}
We then have that $\mathbb{E}[X]\leq 2^{n+1}e^{-n^{2a-1}/(3\varepsilon)}$, and the statement follows by Markov's inequality.
\end{proof}

\begin{remark}\label{rmk:degdev}
In particular, for any $a\in(1/2,1)$ we have that a.a.s.~the number of vertices whose degree deviates from the expectation by more than $n^a$ is at most $e^{-n^{2a-1}/(6\varepsilon)}2^n$.
\end{remark}

It will be important to show that a.a.s.~there are not too many of the above vertex type `close' to any given vertex.

\begin{lemma}\label{lem:badverticesdonotclump}
Let $\varepsilon\in(0,1)$, $a\in(2/3,1)$ and $\ell\in\mathbb{N}$.
Then, for any $b>2-2a$\COMMENT{We can possibly get $\geq$ here, but I think that will require being a bit more careful and giving a stronger upper bound on $\varepsilon$.}, a.a.s.~there are no vertices $v\in\{0,1\}^n$ for which $|\{u\in B^\ell(v): d_{\mathcal{Q}_\varepsilon^n}(u)\neq\varepsilon n\pm n^a\}|\geq n^{b}$.
\end{lemma}

\begin{proof}
First, note that we may assume that $b<a$ (otherwise, choose a value $b'\in(2-2a,a)$\COMMENT{It is for this to be nonempty that we need $a>2/3$.} and prove the statement for this value, which in turn implies the result for $b$).
Throughout this proof, we write $d(v)$ for $d_{\mathcal{Q}_\varepsilon^n}(v)$.
Fix any vertex $v\in\{0,1\}^n$. 
Let $X(v)\coloneqq|\{u\in B^\ell(v): d(u)\neq\varepsilon n\pm n^a\}|$.

If $X(v)\geq n^{b}$, there exists a set $A\subseteq B^\ell(v)$ of size $|A|=n^{b}$ such that $d(u)\neq\varepsilon n\pm n^a$ for all $u\in A$.
We call such a set $A$ \emph{bad}.
We now bound the probability that such a bad set exists.
Given any set $A\in\binom{B^\ell(v)}{n^{b}}$, for each $u\in A$ let $d^A(u)\coloneqq|N_{\cQ_\varepsilon^n}(u)\setminus A|$.
Observe that $d^A(u)=d(u)\pm n^{b}$ and, since $b<a$, for any $u\in A$ we have that, if $d(u)\neq\varepsilon n\pm n^a$, then $d^A(u)\neq\varepsilon n\pm n^a/2$.

Fix a set $A\in\binom{B^\ell(v)}{n^{b}}$.
Observe that $\mathbb{E}[d^A(u)]\in[\varepsilon n(1-n^{b-1}),\varepsilon n]$ for all $u\in A$.
Furthermore, the variables $\{d^A(u): u\in A\}$ are mutually independent, and each of them follows a binomial distribution.
By \cref{lem:Chernoff}, for each $u\in A$ we have that\COMMENT{We have that $d^A(u)$ is stochastically dominated by a random variable $Y\sim\text{Bin}(n,\varepsilon)$ and stochastically dominates a random variable $Z\sim\text{Bin}(n-n^{b},\varepsilon)$.
Using \cref{lem:Chernoff}, it follows that 
\[
    \mathbb{P}[d^A(u)\geq\varepsilon n+n^a/2]\leq\mathbb{P}[Y\geq\varepsilon n+n^a/2]=\mathbb{P}\left[Y\geq\left(1+\frac{n^{a-1}}{2\varepsilon}\right)\varepsilon n\right]\leq e^{-n^{2a-1}/(12\varepsilon)}.
\]
For $Z$, observe that $\mathbb{E}[Z]=\varepsilon n(1-1/n^{1-b})$.
Furthermore, observe that $\frac{(1-1/n^{1-b})(1 - n^{a-1}/(4\eps))}{1 - n^{a-1}/(2\eps)} > 1$ since the numerator is at least
$1 - \frac{n^{a-1}}{4\eps} - \frac{1}{n^{1-b}} > 1 - \frac{n^{a-1}}{2\eps}$ since $ \frac{n^{a-1}}{4\eps} > \frac{1}{n^{1-b}}$ as $b<a$.
By \cref{lem:Chernoff},
\begin{align*}
    \mathbb{P}[d^A(u)\leq\varepsilon n-n^a/2]&\leq\mathbb{P}[Z\leq\varepsilon n-n^a/2]\leq\mathbb{P}\left[Z\leq\varepsilon n\left(1-\frac{n^{a-1}}{2\varepsilon}\right)\frac{(1 - n^{a-1}/(4\eps))(1-\frac{1}{n^{1-b}})}{1 - n^{a-1}/(2\eps)}\right]\\
    & = \mathbb{P}\left[Z\leq\left(1-\frac{n^{a-1}}{4\varepsilon}\right)\mathbb{E}[Z]\right] \leq e^{-n^{2a-1}(1-1/n^{1-b})/(32\varepsilon)}\leq e^{-n^{2a-1}/(40\varepsilon)}.
\end{align*}
Adding both of these yields the result.}
\[\mathbb{P}[d(u)\neq\varepsilon n\pm n^a]\leq\mathbb{P}[d^A(u)\neq\varepsilon n\pm n^a/2]\leq2e^{-n^{2a-1}/(40\varepsilon)}.\]
Since the variables $d^A(u)$ are mutually independent, it follows that
\[\mathbb{P}[A\text{ is bad}]\leq\left(2e^{-n^{2a-1}/(40\varepsilon)}\right)^{n^{b}}.\]

Now consider a union bound over all possible choices of $A$ and all choices of $v$.
It suffices to prove that
\[\binom{\ell n^\ell}{n^{b}}\left(2e^{-n^{2a-1}/(40\varepsilon)}\right)^{n^{b}}<2.1^{-n}.\]
Since $\binom{\ell n^\ell}{n^{b}}\leq(e\ell n^{\ell-b})^{n^{b}}$, it suffices to show that
\[n^{b}(1+\ln 2+ \ln \ell + (\ell-b)\ln n-n^{2a-1}/(40\varepsilon))<-n\ln 2.1.\]
But this follows for $n$ sufficiently large, from the fact that $b>2-2a$.
\end{proof}

Next we show that, in any ball of radius $\ell$, the number of vertices whose degree is far from the expected is much smaller (at most a constant) if we allow larger deviations for the degrees.
Even more, we can prove a similar statement if we restrict the degrees to some linear subsets of the total neighbourhood in $\cQ^n$.
Recall that, for any vertex $v\in\{0,1\}^n$, any graph $G\subseteq\cQ^n$ and a set $S\subseteq N_{\cQ^n}(v)$, we denote $d_G(v,S) = |N_G(v)\cap S|$\index{deg1H3@$d_G(u,W)$}.

\begin{lemma}\label{lem:vbadverticesdonotclump}
Let $\varepsilon,\delta,\gamma\in(0,1)$ and $\ell\in\mathbb{N}$.
For each $v\in\{0,1\}^n$, let $S(v)\subseteq N_{\cQ^n}(v)$ satisfy $|S(v)|\geq\gamma n$.
Let $\mathcal{E}$ be the event that there are no vertices $v\in\{0,1\}^n$ for which $|\{u\in B^\ell(v):d_{\mathcal{Q}_\varepsilon^n}(u, S(u))\neq(1\pm\delta)\varepsilon |S(u)|\}|\geq100/(\delta^2\varepsilon\gamma)$.
Then, for $n$ sufficiently large, $\mathbb{P}[\mathcal{E}]\geq1-e^{-4n}$.
\end{lemma}

\begin{proof}
Throughout this proof, we write $d(v)$ for $d_{\mathcal{Q}_\varepsilon^n}(v)$ and $d(v, S)$ for $d_{\mathcal{Q}_\varepsilon^n}(v, S)$, for any set $S$.

Let $C\coloneqq\lceil100/(\delta^2\varepsilon\gamma)\rceil$.
Fix any vertex $v\in\{0,1\}^n$ and $A\in\binom{B^\ell(v)}{C}$.
Observe that for any $u\in A$, if $d(u, S(u))\neq(1\pm\delta)\varepsilon|S(u)|$, then $d(u, S(u)\setminus A)\neq(1\pm\delta/2)\varepsilon |S(u)|$.
Observe that $\mathbb{E}[d(u, S(u)\setminus A)]\in[\varepsilon(|S(u)|-C),\varepsilon|S(u)|]$ for all $u\in A$.
Furthermore, the variables $\{d(u, S(u)\setminus A): u\in A\}$ are mutually independent, and each of them follows a binomial distribution.
By \cref{lem:Chernoff}, for each $u\in A$ we have that, for $n$ sufficiently large,\COMMENT{We have that $d(u, S(u)\setminus A)$ is stochastically dominated by a random variable $Y\sim\text{Bin}(|S(u)|,\varepsilon)$ and stochastically dominates a random variable $Z\sim\text{Bin}(|S(u)|-C,\varepsilon)$.
Using \cref{lem:Chernoff}, it follows that 
\[
    \mathbb{P}[d(u, S(u)\setminus A)\geq(1+\delta/2)\varepsilon|S(u)|]\leq\mathbb{P}[Y\geq(1+\delta/2)\varepsilon|S(u)|]\leq e^{-\delta^2\varepsilon|S(u)|/12}\leq e^{-\delta^2\varepsilon\gamma n/12}.
\]
For $Z$, observe that $\mathbb{E}[Z]=\varepsilon(|S(u)|-C)$.
By \cref{lem:Chernoff},
\[\mathbb{P}[d(u, S(u)\setminus A)\leq(1-\delta/2)\varepsilon |S(u)|]\leq\mathbb{P}[Z\leq(1-\delta/2)\varepsilon |S(u)|]\leq\mathbb{P}\left[Z\leq\left(1-\delta/3\right)\mathbb{E}[Z]\right]\leq e^{-\delta^2\varepsilon (|S(u)|-C)/18}\leq e^{-\delta^2\varepsilon\gamma n/19}.\]
Adding both of these yields the result.}
\[\mathbb{P}[d(u)\neq(1\pm\delta)\varepsilon|S(u)|]\leq\mathbb{P}[d(u, S(u)\setminus A)\neq(1\pm\delta/2)\varepsilon|S(u)|]\leq2e^{-\delta^2\varepsilon\gamma n/19}\leq e^{-\delta^2\varepsilon\gamma n/20}.\]
We say that $A$ is bad if $d(u, S(u))\neq(1\pm\delta)\varepsilon|S(u)|$ for all $u\in A$.
Since the variables $d(u, S(u)\setminus A)$ are mutually independent, it follows that
\[\mathbb{P}[A\text{ is bad}]\leq\left(e^{-\delta^2\varepsilon\gamma n/20}\right)^C\leq e^{-5n}.\]
Observe that $\mathcal{E}$ holds if there are no bad sets $A$.
By a union bound over all choices of $v$ and all choices of $A$, it follows that
\[\mathbb{P}[\overline{\mathcal{E}}]\leq2^n\binom{\ell n^\ell}{C}e^{-5n}\leq e^{-4n}.\qedhere\]
\end{proof}

In more generality than \cref{lem:degconcentrates,lem:vbadverticesdonotclump}, we will need to use the fact that in $\cQ^n_\eps$ the directions in which the neighbours of a vertex lie are not correlated too much between vertices.
Given a graph $G\subseteq\cQ^n$, for any set $S\subseteq\mathcal{D}(\cQ^n)$ and any vertex $x\in\{0,1\}^n$, we denote $N_{G,S}(x)\coloneqq\{x+\hat e:\hat e\in S,\{x,x+\hat e\}\in E(G)\}$\index{NGS@$N_{G,S}(x)$} and $d_{G,S}(x)\coloneqq|N_{G,S}(x)|$\index{degGeS@$d_{G,S}(x)$}.
Similarly, for any $y\in\{0,1\}^n$ such that $\dist(x,y)=1$, we denote $N_{G,S,x}(y)\coloneqq N_{G,S}(y)\setminus \{x\}$\index{NGS2@$N_{G,S,x}(y)$} and $d_{G,S,x}(y)\coloneqq|N_{G,S,x}(y)|$\index{degGeS2@$d_{G,S,x}(y)$}.

\begin{lemma}\label{lem:directionconcentrates}
  For every $\varepsilon,\delta\in(0,1)$, a.a.s.~the following holds for every $x\in\{0,1\}^n$: for any set $S\subseteq\mathcal{D}(\cQ^n)$ with $|S|\geq\delta n$, all but at most $100/(\varepsilon\delta)$ vertices $y\in N_{\cQ^n}(x)$ satisfy $d_{\cQ^n_\varepsilon,S,x}(y)\geq 2\varepsilon|S|/3$.
\end{lemma}

\begin{proof}
Fix a vertex $x\in\{0,1\}^n$ and a set $S\subseteq\mathcal{D}(\cQ^n)$ with $|S|\geq\delta n$.
Choose any vertex $y\in N_{\cQ^n}(x)$, and consider the variable $X(y)\coloneqq d_{\cQ^n_\varepsilon,S}(y)$.
It suffices to prove that a.a.s.~$X(y)>2\varepsilon|S|/3$ for all but $100/(\varepsilon\delta)$ vertices $y\in N_{\cQ^n}(x)$.
Observe that $X(y)\sim\mathrm{Bin}(|S|,\varepsilon)$, so $\mathbb{E}[X(y)]=\varepsilon|S|$ and, by \cref{lem:Chernoff},
\[\mathbb{P}[X(y)\leq2\varepsilon|S|/3]\leq e^{-\varepsilon|S|/18}\leq e^{-\varepsilon\delta n/18}.\]
Observe, furthermore, that $N_{\cQ^n}(x)$ is an independent set in $\cQ^n$, hence the variables $\{X(y): y\in N_{\cQ^n}(x)\}$ are mutually independent.
It follows that the probability that at least $100/(\varepsilon\delta)$ vertices $y\in N_{\cQ^n}(x)$ do not satisfy the bound is at most $\binom{n}{100/(\varepsilon\delta)}e^{-5n}$\COMMENT{If the statement is not satisfied (for the fixed choices of $x$ and $S$), that means that there exists a set of $100/(\varepsilon\delta)$ vertices $y\in N_{\cQ^n}(x)$ which do not satisfy the required bound.
The number of such sets is $\binom{n}{100/(\varepsilon\delta)}$, and the probability that all vertices in any such set fail the bound simultaneously is at most $(e^{-\varepsilon\delta n/18})^{100/(\varepsilon\delta)}<e^{-5n}$.}.
Finally, by a union bound over all choices of $S$ and $x$, we conclude that the statement fails with probability at most $2^{3n}e^{-5n}=o(1)$\COMMENT{We have $2^n$ choices for $x$, at most $2^n$ choices for $S$, and $2^n$ is also a bound for $\binom{n}{100/(\varepsilon\delta)}$.}.
\end{proof}

We are also interested in the number of subcubes in which each vertex lies.
Given a graph $G$, a vertex $v$ and any $\ell\in\mathbb{N}$, we denote the number of copies of $\cQ^\ell$ in $G$ which contain $v$ by $d_G^\ell(v)$\index{degGl@$d_G^\ell(u)$}.
It is easy to give trivial upper bounds for this number by considering its value in $\cQ^n$.
Indeed, for all $v\in\{0,1\}^n$ we have that
\begin{equation}\label{equa:l-degbound}
    d^\ell_{\cQ^n_\varepsilon}(v)\leq\binom{n}{\ell}.
\end{equation}

\begin{lemma}\label{lem:RNdegree}
  Let $\varepsilon\in (0,1)$, $a\in(1/2,1)$ and $\ell\in\mathbb{N}$.
  Then, a.a.s.~all but at most $2^{n}e^{-n^{2a-1}/(6\varepsilon)}$ vertices $v\in V(\cQ^n_\varepsilon)$ satisfy 
  \begin{equation}\label{equa:RNdegree}
      d^\ell_{\cQ^n_\varepsilon}(v)=(1\pm\bigO(n^{a-1}))\frac{\varepsilon^{2^{\ell-1}\ell}}{\ell!}n^\ell.
  \end{equation}
\end{lemma}

\begin{proof}
Throughout this proof, we will write $d(v)$ for $d_{\cQ_\varepsilon^n}(v)$, $N(v)$ for $N_{\cQ_\varepsilon^n}(v)$ and $d^\ell(v)$ for $d^\ell_{\cQ_\varepsilon^n}(v)$.

Fix a vertex $v\in\{0,1\}^n$, reveal all the edges incident to $v$, and condition on the event that $d(v)=\varepsilon n\pm n^a$.
We then have that\COMMENT{Without conditioning, the expectation is $\binom{n}{\ell}\varepsilon^{2^{\ell-1}\ell}$, which is asymptotically the same as we get here.}\COMMENT{This can be computed as follows.
The number of possible cubes $\cQ^\ell$ is given by $\binom{\varepsilon n\pm n^a}{\ell}$, as any cube that contains $v$ must use edges which are incident to $v$ and we have conditioned on the event of having $\varepsilon n\pm n^a$ such edges.
The probability that any such cube lies in our graph is the probability of each edge being there to the power of the number of edges (because edges are chosen independently).
The number of edges in a cube $\cQ^\ell$ is $2^\ell\ell/2$, but $\ell$ such edges are guaranteed to be there (again, because of the condition).}
\begin{equation}\label{eqn:cubedeg}
\mathbb{E}[d^\ell(v)]=\binom{\varepsilon n\pm n^a}{\ell}\varepsilon^{(2^{\ell-1}-1)\ell}=(1\pm\bigO(n^{a-1}))\frac{n^\ell}{\ell!}\varepsilon^{2^{\ell-1}\ell}.
\end{equation}

Let $\mathcal{D}(v)\subseteq\mathcal{D}(\cQ^n)$ be the set of directions such that $N(v)=v+\mathcal{D}(v)$.
Consider the graph $\Gamma_\ell(v)\coloneqq\cQ^n[v+\ell(\mathcal{D}(v)\cup\{\mathbf{0}\})]$.
For each $i\in[\ell]$, let $L_i\subseteq E(\Gamma_\ell(v))$ be the set of edges which are at distance $i-1$ from $v$ in $\cQ^n$.
Note that these sets partition $E(\Gamma_\ell(v))$ and that $|L_i|=i\binom{d(v)}{i}$\COMMENT{the number of edges at distance $i-1$ is equal to the number of vertices at distance $i$ in $\Gamma_\ell(v)$ (which is $\binom{d(v)}{i}$) times the degree of these vertices into vertices which are at distance $i-1$ from $v$ ($i$).}.
Let $m\coloneqq|E(\Gamma_\ell(v))|$ and $m_j\coloneqq\sum_{i=1}^j|L_i|$ for all $j\in[\ell]$.
Label the edges of $\Gamma_\ell(v)$ as $e_1,\ldots,e_m$ in such a way that all the edges in $L_1$ come first, then the edges in $L_2$, and so on, until covering all the edges in $L_\ell$.
For each $j\in[m]$, let $X_j$ be the indicator random variable that $e_j\in E(\cQ_\varepsilon^n)$ (recall that we condition on the neighbourhood of $v$ being revealed and $v$ being good).
We now consider an edge-exposure martingale given by the variables $Y_{j}\coloneqq\mathbb{E}[d^\ell(v)\mid X_1,\ldots,X_j]$, for $j\in[m]_0$.
This is a Doob martingale with $Y_{d(v)}=\mathbb{E}[d^\ell(v)]$\COMMENT{after the conditioning} and $Y_m=d^\ell(v)$.

We must now bound the differences $|Y_j-Y_{j-1}|$, for all $j\in[m]$.
Observe that the maximum change in the expected number of $\ell$-dimensional cubes in $\cQ^n_\eps$ containing $v$ when a new edge $e_i$ is revealed is bounded from above by the number of such cubes in $\Gamma_\ell(v)$ containing $e_i$.
Given any $k\in[\ell]\setminus\{1\}$ and any $i\in[m_k]\setminus[m_{k-1}]$, we claim that the number of copies of $\cQ^\ell$ in $\Gamma_\ell(v)$ containing $e_i$ is bounded by $\binom{d(v)-k}{\ell-k}$ (recall that for all $i\in[m_1]$ we have that $Y_i=Y_{i-1}$).
Indeed, let $e_i=\{x,y\}$ with $\dist(x,v)=k-1$, and let $\mathcal{D}_k(x)\subseteq\mathcal{D}(v)$ be the set of $k-1$ directions such that $\cQ^n(v,\mathcal{D}_k(x))$ contains $x$.
Then, any copy $\cQ$ of $\cQ^\ell$ in $\Gamma_\ell(v)$ containing $e_i$ must satisfy that $\mathcal{D}(\cQ)\subseteq\mathcal{D}(v)$ contains $\mathcal{D}_k(x)$, the direction given by $y-x$, and any other $\ell-k$ of the directions in $\mathcal{D}(v)$, for which there is the claimed number of choices.
Therefore, we conclude that
\begin{equation}\label{eqn:martdif}
\sum_{i=1}^m|Y_i-Y_{i-1}|^2\leq\sum_{k=2}^\ell k\binom{d(v)}{k}\binom{d(v)-k}{\ell-k}^2= \frac{1}{((\ell-2)!)^2}d(v)^{2\ell-2}\left(1\pm\bigO\left(\frac{1}{d(v)}\right)\right).
\end{equation}
Hence, by \cref{lem: Azuma}, for $n$ sufficiently large we have that\COMMENT{We have that
\begin{align*}
    \mathbb{P}\left[|d^\ell(v)-\mathbb{E}[d^\ell(v)]|\geq\sqrt{\frac{2}{\varepsilon}}\frac{1}{(\ell-2)!}d(v)^{\ell-1/2}\right] &\leq2\exp\left(\frac{-\frac{2}{\varepsilon((\ell-2)!)^2}d(v)^{2\ell-1}}{\frac{2}{((\ell-2)!)^2}d(v)^{2\ell-2}\left(1\pm\bigO\left(\frac{1}{d(v)}\right)\right)}\right)\\
    &=2\exp\left(\left(1\pm\bigO\left(\frac{1}{d(v)}\right)\right)\frac{-d(v)}{\varepsilon}\right)\\
    &=2\exp\left(-\left(1\pm\bigO\left(n^{-1}\right)\pm\bigO\left(n^{a-1}\right)\right)n\right)\\
    &< 2.1^{-n}.
\end{align*}}
\[\mathbb{P}\left[|d^\ell(v)-\mathbb{E}[d^\ell(v)]|\geq \sqrt{\frac{2}{\varepsilon}}\frac{1}{(\ell-2)!}d(v)^{\ell-1/2}\right]\leq 2.1^{-n}.\]

Finally, by \cref{rmk:degdev} combined with a union bound on all vertices $v$ such that $d(v)=\varepsilon n\pm n^a$, we conclude that a.a.s.~all such vertices satisfy that $d^\ell(v)=(1\pm\bigO(n^{a-1}))n^\ell\varepsilon^{2^{\ell-1}\ell}/\ell!$.\COMMENT{
In the proof we conditioned on $v$ being good, that is $d(v) = \eps n \pm n^a$.
So we showed that \[
\mathbb{P}\left[|d^\ell(v)-\mathbb{E}[d^\ell(v)]|\geq \sqrt{\frac{2}{\varepsilon}}\frac{1}{(\ell-2)!}d(v)^{\ell-1/2} \wedge v \text{ good}\right] \le \mathbb{P}\left[|d^\ell(v)-\mathbb{E}[d^\ell(v)]|\geq \sqrt{\frac{2}{\varepsilon}}\frac{1}{(\ell-2)!}d(v)^{\ell-1/2}\mid v \text{ good}\right] \le 2.1^{-n}.\]
It follows by a union bound over all vertices that a.a.s.~if  any vertex is good it can't be that $|d^\ell(v)-\mathbb{E}[d^\ell(v)]|\geq \sqrt{\frac{2}{\varepsilon}}\frac{1}{(\ell-2)!}d(v)^{\ell-1/2}$.
Now by \cref{rmk:degdev} we have also that a.a.s.~all but at most $2^ne^{-n^{2a -1}/(6\eps)}$ are not good and the result follows by a union of these two a.a.s.~statements.
}
\end{proof}

\begin{remark}\label{remark:RNdegree}
In particular, the proof of \cref{lem:RNdegree} shows that a.a.s.~all vertices $v\in\{0,1\}^n$ which satisfy $d_{\cQ^n_\varepsilon}(v)=\varepsilon n\pm n^a$ also satisfy \eqref{equa:RNdegree}.
Therefore, by \cref{lem:badverticesdonotclump}, for any $r\in\mathbb{N}$ and $a \in (2/3, 1)$, a.a.s.~in any ball of radius $r$, all but at most $n^{2-2a+\eta}$ vertices satisfy \eqref{equa:RNdegree}, where $\eta>0$ is an arbitrarily small constant.
\end{remark}

In more generality, we will need to bound the number of subcubes which contain a given pair of vertices.
Given a graph $G$, two vertices $u$ and $v$, and any $\ell\in\mathbb{N}$, we denote the number of copies of $\cQ^\ell$ in $G$ which contain both $u$ and $v$ by $d_G^\ell(u,v)$\index{degGl2@$d_G^\ell(u,v)$}.
Again, we can easily give upper bounds for this number in $\cQ^n_\varepsilon$ by considering its value in $\cQ^n$.
Indeed, for all $u,v\in\{0,1\}^n$ we have that
\begin{equation}\label{equa:l-codegbound}
    d^\ell_{\cQ^n}(u,v)\leq\binom{n}{\ell-\dist(u,v)}\leq n^{\ell-\dist(u,v)}
\end{equation}
(here, we understand that $\binom{n}{a}=0$ for all $a<0$).\COMMENT{Indeed, observe that any cube containing both $u$ and $v$ must use all the directions given by $\mathcal{D}(u,v)$, of which there are $\dist(u,v)$, and the remaining directions can be chosen in any way.}

We will also need the property that the cubes containing a given vertex use different directions quite evenly.
More precisely, given any graph $G\subseteq\cQ^n$, any set $S\subseteq\mathcal{D}(\cQ^n)$, two vertices $x,y\in\{0,1\}^n$, an integer $\ell\in\mathbb{N}$ and a real $t\in\mathbb{R}$, we denote by $d^\ell_{G,S,t,x}(y)$\index{degGl3@$d_{G,S,t,x}^\ell(y)$} the number of copies $C$ of $\cQ^\ell$ which contain $y$, do not contain $x$, and satisfy $|\mathcal{D}(C)\cap S|\geq t$.

\begin{lemma}\label{lem:RNdirections}
  Let $0< 1/\ell \ll \delta <1$, with $\ell\in\mathbb{N}$.
  Let $\varepsilon,\eta\in(0,1)$ and $a\in(2/3,1)$.
  Then, a.a.s.~the following holds for every $x\in\{0,1\}^n$: for any set $S\subseteq\mathcal{D}(\cQ^n)$ with $|S|\geq\delta n$, all but at most $n^{2(1+\eta-a)}$ vertices $y\in N_{\cQ^n}(x)$ satisfy
  \[d^\ell_{\cQ^n_\varepsilon,S,\ell^{1/2},x}(y)\geq\frac{\varepsilon^{2^{\ell-1}\ell}}{2\ell!}n^\ell.\]
\end{lemma}

\begin{proof}
Throughout this proof, for any $x\in\{0,1\}^n$ and any $y\in N_{\cQ^n}(x)$, we will write $N(y)$ for $N_{\cQ^n_\varepsilon}(y)$, $d(y)$ for $d_{\cQ^n_\varepsilon}(y)$, $d_{S,x}(y)$ for $d_{\cQ^n_\varepsilon,S,x}(y)$, $d^\ell(y)$ for $d^\ell_{\cQ^n_\varepsilon}(y)$ and $d^\ell_{S,\ell^{1/2},x}(y)$ for $d^\ell_{\cQ^n_\varepsilon,S,\ell^{1/2},x}(y)$.

Fix a set $S\subseteq\mathcal{D}(\cQ^n)$ with $|S|\geq\delta n$.
Let $D\coloneqq\varepsilon^{2^{\ell-1}\ell}n^\ell/\ell!$ and fix a vertex $x\in\{0,1\}^n$.
Consider any $y\in N_{\cQ^n}(x)$; reveal all edges incident to $y$ and condition on the event that $d(y)=\varepsilon n\pm n^{a}$ and $d_{S,x}(y)\geq\varepsilon|S|/2$.
By \cref{lem:badverticesdonotclump,lem:directionconcentrates}, we have that a.a.s.~all but at most $n^{2(1+\eta-a)}$ vertices $y\in N_{\cQ^n}(x)$ satisfy this event.\COMMENT{Apply \cref{lem:badverticesdonotclump} with $b=2-2a+\eta$. Since the number of bad vertices given by \cref{lem:directionconcentrates} is constant, the result follows.}

Let $\mathcal{D}(y)$ be the set of directions such that $N(y)\setminus\{x\}=y+\mathcal{D}(y)$.
Thus, $|\cD(y)| = d(y) \pm1$.
Let $\alpha\coloneqq|S\cap\mathcal{D}(y)|/n$, and note that $\varepsilon\delta/2\leq\alpha\leq\varepsilon+{n^{a-1}}$.
Similar to the proof of \eqref{eqn:cubedeg}, we have that
\[\mathbb{E}[d^\ell_{S,\ell^{1/2},x}(y)]=\varepsilon^{(2^{\ell-1}-1)\ell}\sum_{i=\lceil\ell^{1/2}\rceil}^\ell\binom{\alpha n}{i}\binom{\varepsilon n-\alpha n\pm (n^a+1)}{\ell-i}\geq 3D/4.\]

Consider the graph $\Gamma_\ell(y)\coloneqq\cQ^n[y+\ell(\mathcal{D}(y)\cup\{\mathbf{0}\})]$.
For each $i\in[\ell]$, let $L_i\subseteq E(\Gamma_\ell(y))$ be the set of edges which are at distance $i-1$ from $y$.
Note that these sets partition $E(\Gamma_\ell(y))$ and that $|L_i|=i\binom{d(y)\pm 1}{i}$.
Let $m\coloneqq|E(\Gamma_\ell(y))|$ and $m_j\coloneqq\sum_{i=1}^j|L_i|$ for all $j\in[\ell]$.
Label the edges of $\Gamma_\ell(y)$ as $e_1,\ldots,e_m$ in such a way that all the edges in $L_1$ come first, then the edges in $L_2$, and so on, until covering all the edges in $L_\ell$.
For each $j\in[m]$, let $X_j$ be the indicator random variable that $e_j\in E(\cQ_\varepsilon^n)$.
We now consider an edge-exposure martingale given by the variables $Y_{j}\coloneqq\mathbb{E}[d^\ell_{S,\ell^{1/2},x}(y)\mid X_1,\ldots,X_j]$, for $j\in[m]_0$.
This is a Doob martingale with $Y_{d(y)}=\mathbb{E}[d^\ell_{S,\ell^{1/2},x}(y)]$\COMMENT{after the conditioning} and $Y_m=d^\ell_{S,\ell^{1/2},x}(y)$.

In order to bound the differences $|Y_j-Y_{j-1}|$, for all $j\in[m]$, observe that the maximum change in the expected number of $\ell$-dimensional cubes in $\cQ^n_\eps$ containing $y$ when a new edge $e_i$ is revealed is bounded from above by the number of such cubes in $\Gamma_\ell(y)$ containing $e_i$.
In particular, this is an upper bound for the maximum change in the expected number of $\ell$-dimensional cubes in $\cQ^n_\eps$ containing $y$, not containing $x$, and whose directions intersect $S$ in a set of size at least $\ell^{1/2}$, when a new edge $e_i$ is revealed.
Thus, similarly as in \eqref{eqn:martdif}, it follows that\COMMENT{Given any $k\in[\ell]\setminus\{1\}$ and any $i\in[m_k]\setminus[m_{k-1}]$, we claim that the number of copies of $\cQ^\ell$ containing $e_i$ is bounded by $\binom{d(y)-k}{\ell-k}$(The proof is the same as in the proof of \cref{lem:RNdegree}.) (recall that for all $i\in[m_1]$ we have that $Y_i=Y_{i-1}$).
\[\sum_{i=1}^m|Y_i-Y_{i-1}|^2\leq\sum_{k=2}^\ell k\binom{d(y)}{k}\binom{d(y)-k}{\ell-k}^2= \frac{1}{((\ell-2)!)^2}d(y)^{2\ell-2}\left(1\pm\bigO\left(\frac{1}{d(y)}\right)\right).\]}
\[\sum_{i=1}^m|Y_i-Y_{i-1}|^2\leq \frac{1}{((\ell-2)!)^2}d(y)^{2\ell-2}\left(1\pm\bigO\left(\frac{1}{d(y)}\right)\right).\]
Hence, by \cref{lem: Azuma}, we have that\COMMENT{We have that
\begin{align*}
    \mathbb{P}\left[d^\ell_{S,t,x}(y)\leq D/2\right] &\leq\mathbb{P}\left[d^\ell_{S,t,x}(y)\leq\mathbb{E}[d^\ell_{S,t,x}(y)]-D/4\right] \leq\mathbb{P}\left[|d^\ell_{S,t,x}(y)-\mathbb{E}[d^\ell_{S,t,x}(y)]|\geq D/4\right]\\
    &\leq2\exp\left(-\frac{D^2}{32}\frac{((\ell-2)!)^2}{d(y)^{2\ell-2}(1\pm\bigO(1/d(y)))}\right)=2\exp\left(-\frac{\varepsilon^{2^\ell \ell} n^{2\ell}}{32(\ell!)^2}\frac{((\ell-2)!)^2}{(\varepsilon n\pm n^a)^{2\ell-2}(1\pm\bigO(1/n))}\right)\\
    &=\exp\left(-\Theta(n^2)\right).
\end{align*}}
\[\mathbb{P}\left[d^\ell_{S,\ell^{1/2},x}(y)\leq D/2\right]\leq e^{-\Theta(n^2)}.\]

Finally, the statement follows by a union bound on all sets $S\subseteq\mathcal{D}(\cQ^n)$ with $|S|\geq\delta n$, on all vertices $x\in\{0,1\}^n$, and on all vertices $y\in N_{\cQ^n}(x)$ such that $d(y)=\varepsilon n\pm n^a$ and $d_{S,x}(y)\geq\varepsilon|S|/2$.\COMMENT{
Let $y$ be good if $d(y) = \eps n \pm n^a$ and $d_{S,x}(y)\ge \eps|S|/2$.
We've shown 
\[\mathbb{P}\left[d^\ell_{S,\ell^{1/2},x}(y)\leq D/2 \wedge y \text{ good}\right]\leq  \mathbb{P}\left[d^\ell_{S,\ell^{1/2},x}(y)\leq D/2 \mid y \text{ good}\right] \le e^{-\Theta(n^2)}.\]
So it follows by a union bound over all vertices that each good $y$ must be such that $d^\ell_{S,\ell^{1/2},x}(y)> D/2$.
We have that a.a.s.~there are at most $n^{2(1+\eta - \alpha)}$ vertices which are not good,
so the proof follows by a union bound of the two a.a.s.~statements.}
\end{proof}


\section{Tiling random subgraphs of the hypercube with small cubes}\label{nib}

Throughout this section, we will consider auxiliary hypergraphs to obtain information about subgraphs of the $n$-dimensional hypercube.
The general idea will be to apply the so-called `R\"odl nibble' to achieve this.
Roughly speaking, the R\"odl nibble is a randomised iterative process which, given an almost regular uniform hypergraph with small codegrees, finds a matching covering all but a small proportion of the vertices.
The basic idea is the following.
Let $H$ be an almost regular uniform hypergraph with small codegrees.
Consider a random subset of the edges $E\subseteq E(H)$, where each edge is taken independently with the same probability.
If this probability is chosen carefully, then one can show that, with high probability, $E$ is `almost' a matching and that the hypergraph resulting after the deletion of all vertices covered by $E$ is still almost regular and has small codegrees.
This allows one to iterate the process until all but a small fraction of the vertices have been covered.
This approach is the basis for the proof of our main result in this section, \cref{thm: nibble}.
The main auxiliary result is \cref{lem:nibble}, which shows that in each iteration of the process we have the properties we require.
In particular, we require our matching to satisfy several additional `local' properties.
This means our application of the nibble will require strong concentration results, as well as the use of the Lov\'asz local lemma.
It is also worth noting that our result relies strongly on the geometry of the hypercube, and cannot be stated for general hypergraphs.

\subsection{The R\"odl nibble}\label{subsec: nibble}

Given $\ell\in\mathbb{N}$ and any graph $G\subseteq\cQ^n$, we will denote by $H_\ell(G)$\index{HlG@$H_\ell(G)$} the $2^\ell$-uniform hypergraph with vertex set $V(G)$ where a set of vertices $W\subseteq\{0,1\}^n$ with $|W|=2^\ell$ forms a hyperedge if and only if $G[W]\cong\cQ^\ell$.
Observe that the vertex set of $H_\ell(G)$ is (a subset of) $\{0,1\}^n$.
Hence, we can use the underlying notation of directions we have considered for hypercubes so far.
In particular, given any pair of vertices $x,y\in V(H_\ell(G))$, any set $S\subseteq\mathcal{D}(\cQ^n)$ and a real $t\in\mathbb{R}$, we denote by $d_{H_\ell(G),S,t,x}(y)$\index{degGl4@$d_{H_\ell(G),S,t,x}(y)$} the number of hyperedges $e\in E(H_\ell(G))$ which contain $y$, do not contain $x$, and satisfy $|\mathcal{D}(G[e])\cap S|\geq t$.
Note that, with the notation from \cref{lem:RNdirections}, $d_{H_\ell(G),S,t,x}(y)=d^\ell_{G,S,t,x}(y)$.
In order to simplify notation, for any vertex $x\in\{0,1\}^n$ and any sets $Y\subseteq N_{\cQ^n}(x)$ and $E\subseteq E(H_\ell(G))$, we let $E_x(Y)\coloneqq\{e\in E:\dist(x,e)=1,e\cap Y\neq\varnothing\}$\index{ExY@$E_x(Y)$}.
If $E$ is the set of all edges of a given hypergraph $H \subseteq H_\ell(G)$, we may sometimes denote this by $E_x(H,Y)$\index{ExY2@$E_x(H,Y)$}.
Furthermore, it is worth noting that, for hypergraphs $H_\ell(G)$ defined as above, the inequality $d_H(x)\leq\sum_{y\in V(H)\setminus\{x\}}d_H(x,y)$, which holds for all hypergraphs $H$\COMMENT{assuming they have no edges of size $1$...} and all $x\in V(H)$, can be improved to the following: for every $\ell\geq2$ and every $x\in V(H_\ell(G))$,
\begin{equation}\label{equa:codegsum}
    d_{H_\ell(G)}(x)\leq\sum_{y\in V(H_\ell(G))\cap N_{\cQ^n}(x)}d_{H_\ell(G)}(x,y).
\end{equation}

The following observations will also come in useful.

\begin{remark}\label{rem:RNcodeg}
Let $\ell,t\in\mathbb{N}$ and $G\subseteq\cQ^n$, and let $H\coloneqq H_\ell(G)$.
Let $x\in V(H)$ and $e\in E(H)$ be such that $\dist(x,e)=t$.
Then, there is a unique vertex $y\in e$ such that $\dist(x,y)=t$\COMMENT{
\begin{proof}
We will think about this using levels.
We may assume without loss of generality that $x=\mathbf{0}$ is the unique vertex in level $0$.
The edge $e$ corresponds to a copy $C$ of $\cQ^\ell$ in $\cQ^n$.
It is clear that any cube $C$ contains exactly one vertex in its lowest level (otherwise, one can show that there is a vertex in a lower level), and this is the unique vertex $y$ such that $\dist(x,y)=\dist(x,e)$ (because the distance between any vertex $y$ and the empty vertex is the level to which $y$ belongs).
\end{proof}}.
Furthermore, for every $e'\in E(H)$ such that $x\in e'$ we have that $e\cap e'\neq\varnothing$ if and only if $y\in e'$\COMMENT{
\begin{proof}
Assume without loss of generality that $x=\mathbf{0}$.
By the previous remark, there is a unique vertex $y$ such that $\dist(x,y)=t$.
In particular, $y$ has exactly $t$ ones among its coordinates.
Now assume that there is an edge $e'$ containing $x$ and intersecting $e$.
Assume $e'$ contains some vertex $z\in e$.
As $z\in e$, we know that $z$ has a $1$ in every coordinate where $y$ has a $1$ (by the minimality of the distance of $y$).
Therefore, as $\mathcal{D}(e')$ (where we see $e'$ as a copy of $\cQ^\ell$) must contain all direction needed to reach $z$ (otherwise $z\notin e'$), in particular it contains all directions needed to reach $y$.
But then $y\in e'$.
\end{proof}}.
In particular, $|\{e'\in E(H):x\in e',e\cap e'\neq\varnothing\}|=d_{H}(x,y)$.
\end{remark}

\begin{remark}\label{rem:RNcodeg2}
Let $\ell,t\in\mathbb{N}$ and $G\subseteq\cQ^n$, and let $H\coloneqq H_\ell(G)$.
Let $x\in V(H)$ and $Y\subseteq N_{\cQ^n}(x)$.
Let $e\in E(H)$ be such that $\dist(x,e)=\dist(Y,e)=t$.
Let $Y'\coloneqq\{y\in Y:\dist(y,e)=t\}$.
Then, $|Y'|\leq\ell$ and none of the edges in $E_x(H,Y\setminus Y')$ intersects $e$.
\end{remark}

\begin{remark}\label{rem:RNExCAP}
Let $\ell\in\mathbb{N}$ and $G\subseteq\cQ^n$, and let $H\coloneqq H_\ell(G)$.
Let $x\in V(H)$ and $Y\subseteq N_{\cQ^n}(x)$.
Then, for any $e\in E_x(H,Y)$, we have $|\{e'\in E_x(H,Y\setminus e):e\cap e'\neq\varnothing\}|=\bigO(n^{\ell-1})$.
\end{remark}

\begin{proof}
We may assume that $G=\cQ^n$.
Given any vertex $v\in e$ with $\dist(x,v)=i$, $1\le i\le \ell+1$, we have that $|\{y\in Y\setminus e:\dist(y,v)<i\}|\leq i-1$.\COMMENT{This is the set of all vertices $y \in Y \setminus e$ such that there is an edge $e'\in E_x(H,Y)$ containing both $y$ and $v$, as for all other vertices $z\in Y\setminus e$, any edge containing both $z$ and $v$ must also contain $x$; the bound follows by considering any of the possible $i$ directions from $x$ to $v$ minus the direction used to reach $e\cap Y$.}
The conclusion follows by summing over all vertices $v\in e$ and using \eqref{equa:l-codegbound}.\COMMENT{We get
\begin{align*}
    |\{e'\in E_x(H,Y\setminus e):e\cap e'\neq\varnothing\}|&\leq\sum_{v\in e}\sum_{\substack{y\in Y\setminus e\\\dist(y,v)=\dist(x,v)-1}}d^\ell_{H}(y,v)\\
    &\leq\sum_{i=1}^{\ell+1}\sum_{\substack{v\in e\\\dist(x,v)=i}}(i-1)n^{\ell-i+1}\\
    &=\sum_{i=2}^{\ell+1}\binom{\ell}{i-1}(i-1)n^{\ell-i+1}=\bigO(n^{\ell-1}).
\end{align*}}
\end{proof}

\begin{remark}\label{rem:basicequationIneedhere}
Let $k,n\in\mathbb{N}$ and $A\subseteq\{0,1\}^n$.
Then,
\[|\{v\in\{0,1\}^n:\dist(v,A)=k\}|\leq|A|\binom{n}{k}.\]
\end{remark}

Consider $\ell\in\mathbb{N}$, $G\subseteq\cQ^n$ and $H\coloneqq H_\ell(G)$.
Recall that each edge of $H$ corresponds to an $\ell$-dimensional subcube of $G$.
Let $e\in E(H)$, $E\subseteq E(H)$ and $S\subseteq\mathcal{D}(\cQ^n)$.
We define the \emph{significance of $e$ in $S$} as $\sigma(e,S)\coloneqq|\mathcal{D}(e)\cap S|$\index{sigmae1@$\sigma(e,S)$}.
Given any $t\in\mathbb{R}$, we say that $e$ is \emph{$t$-significant in $S$} if $\sigma(e,S)\geq t$.
We define the \emph{significance of $E$ in $S$} as $\sigma(E,S)\coloneqq\sum_{e\in E}\sigma(e,S)$\index{sigmae2@$\sigma(E,S)$}.
We denote $\Sigma(E,S,t)\coloneqq\{e\in E:\sigma(e,S)\geq t\}$.\index{sigmaEST@$\Sigma(E,S,t)$}
In particular, $\sigma(E,S)\geq t|\Sigma(E,S,t)|$.

With this, we are now ready to state the main auxiliary result in this section.
This shows that, given $H=H_\ell(G)$, under suitable conditions about the degrees, the codegrees and the local distribution of the edges of $H$ along the directions of the cube (namely, that the edges are significant in every large set of directions), one iteration of our nibble process will yield a subset of edges which is locally close to a matching, satisfies several local properties that we require of our matching (namely, the edges given by the nibble are sufficiently significant in large sets of directions, and not too significant in any given direction), and its deletion yields a hypergraph which still satisfies almost the same suitable conditions for further iterations.

\begin{lemma}\label{lem:nibble}
Let $\ell,k,K\in\mathbb{N}$ with $k>\ell\geq2$ and let $\beta\in(0,1]$.
Let $G \subseteq \cQ^n$ and let $H\coloneqq H_\ell(G)$.
Fix $x\in\{0,1\}^n$.
Let $A_0\coloneqq N_{\cQ^n}(x)$ and, for each $i\in[K]$, let $A_i\subseteq A_0$ be a set of size $|A_i|\geq\beta n$\COMMENT{These are not necessarily distinct.}.
Assume that there exist two constants $a\in(3/4,1)$ and $\gamma\in(0,1]$, and $D=\Theta(n^\ell)$, such that
\begin{enumerate}[label=$(\mathrm{P}\arabic*)$]
    \item\label{RNPdegree} for every $y\in V(H)\cap B^k_{\cQ^n}(x)$ we have $d_H(y) = (1 \pm \bigO(n^{a-1}))D$;
    \item\label{RNPvertexdistribution} for every $i\in[K]_0$ we have $|V(H)\cap A_i|=(1\pm\bigO(n^{a-1}))\gamma|A_i|$.
\end{enumerate}
Then, for all $\varepsilon\ll1/\ell$, the following holds.

Let $E'\subseteq E(H)$ be a random subset of $E(H)$ obtained by adding each edge with probability $\varepsilon/D$, independently of every other edge.
Let $E''\subseteq E'$ be the set of all edges not intersecting any other edge of $E'$.
Then, $E'$, $E''$, $V'\coloneqq V(H)\setminus V(E')$ and $H'\coloneqq H[V']$ satisfy the following:
\begin{enumerate}[label=$(\mathrm{N}\arabic*)$]
    \item\label{RNCremainingvertices} with probability at least $1-e^{-\Theta(n^{1/2})}$, for every $i\in[K]_0$ we have 
    \[|V'\cap A_i|=(1\pm\bigO(n^{a-1}))e^{-\varepsilon}\gamma|A_i|;\]
    \item\label{RNCmatchingsize} with probability at least $1-e^{-\Theta(n^{1/2})}$, for every $i\in[K]_0$ we have 
    \[|V(E'')\cap A_i|\geq\varepsilon(1-2^{\ell+1}\varepsilon)\gamma|A_i|;\]
    \item\label{RNCsingledirection} with probability at least $1-e^{-\Theta(n^{1/2})}$, for every $\hat e\in\mathcal{D}(\cQ^n)$ we have
    \[|\Sigma(E'_x(A_0),\{\hat e\},1)|=o(n^{1/2});\]
    \item\label{RNCdegree} with probability at least $1-e^{-\Theta(n^{2a-1})}$, for every $y\in V(H')\cap B^{k-\ell}_{\cQ^n}(x)$ we have
    \[d_{H'}(y)=(1\pm\bigO(n^{a-1}))e^{-(2^\ell-1)\varepsilon}D.\]
\end{enumerate}

If, in addition to \ref{RNPdegree} and \ref{RNPvertexdistribution}, there exist $c,\delta\in(0,1]$ such that
\begin{enumerate}[label=$(\mathrm{P}\arabic*)$,start=3]
    \item\label{RNPdirections} for every $i\in[K]_0$ and every $S\subseteq\mathcal{D}(\cQ^n)$ with $|S|\geq\delta n$ we have 
    \[|\Sigma(E_x(H,A_i),S,\ell^{1/2})|\geq(1-\bigO(n^{a-1}))c\gamma|A_i|D,\]
\end{enumerate}
then $E'$ and $H'$ also satisfy the following:
\begin{enumerate}[label=$(\mathrm{N}\arabic*)$,start=5]
    \item\label{RNCdirectionsleft} with probability at least $1-e^{-\Theta(n^{1/2})}$, for every $i\in[K]_0$ and every $S\subseteq\mathcal{D}(\cQ^n)$ with $|S|\geq\delta n$ we have 
    \[|\Sigma(E_x(H',A_i),S,\ell^{1/2})|\geq(c-\varepsilon)e^{-(2^\ell-1)\varepsilon}\gamma|A_i|D;\]
    \item\label{RNCdirectionscovered} for any fixed $S\subseteq\mathcal{D}(\cQ^n)$ with $|S|\geq\delta n$, with probability at least $1-e^{-\varepsilon c^2\gamma\beta n/100}$, for every $i\in[K]_0$ and $n$ sufficiently large we have
    \[|V(\Sigma(E'_x(A_i),S,\ell^{1/2}))\cap A_i|\geq\varepsilon c^2\gamma|A_i|/8.\]
\end{enumerate}
\end{lemma}

\begin{proof}
We begin by noting that, since $H=H_\ell(G)$, the following two properties hold:
\begin{enumerate}[label=$(\mathrm{P}\arabic*)$,start=4]
    \item\label{RNPcodegreefar} for all $y,z \in V(H)$ such that $\dist(y,z)>\ell$, we have that $d_H(y,z)=0$;
    \item\label{RNPcodegreeclose} for each $i\in[\ell]$, for all $y,z \in V(H)$ with $\dist(y,z)=i$, we have that $d_H(y,z)=\bigO(D/n^i)$.
\end{enumerate}
Indeed, by the definition of $H$, both follow from \eqref{equa:l-codegbound}.
These will be used repeatedly throughout the proof.

We next observe another simple property which will be useful later in the proof.
Fix any $i\in[K]_0$.
Note that, by \ref{RNPdegree} and \ref{RNPvertexdistribution}, $|E_x(H,A_i)|=\sum_{y\in A_i\cap V(H)}d_H(y)\pm \ell n^\ell=(1\pm\bigO(n^{a-1}))\gamma|A_i|D$\COMMENT{$|E_x(H,A_i)|$ is the number of edges of $H$ which contain any vertex in $A_i$ and do not contain $x$. 
It is clear that the number is at most the sum of the degrees of each vertex in $A_i\cap V(H)$.
However, we must subtract the number of edges which contain $x$ and intersect $A_i$, for which an upper bound would be the degree of $x$ in $H$.
But $x$ might not be a vertex of $H$, in which case the degree is not well defined.
So we want to give a general upper bound for this quantity.
For instance, we can consider the degree of $x$ in $H_\ell(\cQ^n)$, which by \eqref{equa:l-degbound} is $\binom{n}{\ell}\leq n^\ell$. 
(We multiply by $\ell$ because the same edge with $x$ could be counted up to $\ell$ times.)
This proves the first inequality, and the second follows from \ref{RNPdegree} and \ref{RNPvertexdistribution}, and the fact that $n^{a-1}D|A_i| > n^\ell$.}.
Therefore, $\mathbb{E}[|E_x'(A_i)|]=(1\pm\bigO(n^{a-1}))\varepsilon\gamma|A_i|$ and, by \cref{lem:Chernoff}\COMMENT{$\mathbb{P}[|E_x'(A_i)|\neq(1\pm\bigO(n^{a-1}))\varepsilon\gamma|A_i|]\leq2e^{-\Theta(n^{2a-2})\varepsilon\gamma|A_i|/3}=e^{-\Theta(n^{2a-1})}$.},
\begin{equation}\label{equa:RNchosenedges}
    \mathbb{P}[|E_x'(A_i)|\neq(1\pm\bigO(n^{a-1}))\varepsilon\gamma|A_i|]=e^{-\Theta(n^{2a-1})}.
\end{equation}
By a union bound over all $i\in[K]_0$, we conclude that, with probability $1-e^{-\Theta(n^{2a-1})}$, we have $|E_x'(A_i)|=(1\pm\bigO(n^{a-1}))\varepsilon\gamma|A_i|$ for all $i\in[K]_0$.

\textbf{\ref{RNCremainingvertices}: On the number of vertices in $\boldsymbol{A_i}$ remaining in $\boldsymbol{H'}$.}

In order to prove that \ref{RNCremainingvertices} holds, fix $i\in[K]_0$.
Let $Y_i\coloneqq V(H)\cap A_i$ and fix any vertex $y\in Y_i$.
By \ref{RNPdegree} we have that $d_{H}(y)=(1\pm\bigO(n^{a-1}))D$.
Therefore,\COMMENT{The probability that $y\in V'$ is the probability that none of the edges containing $y$ is picked.
Since the edges are chosen independently,
\begin{align*}
    \mathbb{P}[y\notin V(E')]&=(1-\varepsilon/D)^{(1\pm\bigO(n^{a-1}))D}=(e^{-\varepsilon/D}\pm\bigO(1/D^2))^{(1\pm\bigO(n^{a-1}))D}=(1\pm\bigO(1/D))e^{-(1\pm\bigO(n^{a-1}))\varepsilon}\\
    &=(1\pm\bigO(1/D))(1\pm\bigO(n^{a-1}))^\varepsilon e^{-\varepsilon}=(1\pm\bigO(n^{a-1})) e^{-\varepsilon},
\end{align*}
where the second inequality is by considering a Taylor expansion and the third inequality can be shown by considering the binomial expansion.}
\[\mathbb{P}[y\in V']=(1-\varepsilon/D)^{(1\pm\bigO(n^{a-1}))D}=(1\pm\bigO(n^{a-1}))e^{-\varepsilon}.\]
Thus, by \ref{RNPvertexdistribution}, $\mathbb{E}[|V'\cap A_i|]=(1\pm\bigO(n^{a-1}))e^{-\varepsilon}\gamma|A_i|$\COMMENT{$\mathbb{E}[|V'\cap A_i|]=(1\pm\bigO(n^{a-1}))e^{-\varepsilon}|Y_i|=(1\pm\bigO(n^{a-1}))e^{-\varepsilon}(1\pm\bigO(n^{a-1}))\gamma|A_i|=(1\pm\bigO(n^{a-1}))e^{-\varepsilon}\gamma|A_i|$.}.
We must now prove that $|V'\cap A_i|$ concentrates with high probability.
However, the events $\{y\notin V(E')\}_{y\in Y_i}$ are not necessarily independent\COMMENT{Actually, they are independent if $x\notin V(H)$ (see \cref{rem:RNcodeg}), but we have no guarantee of this.
So they are not necessarily independent.}.

In order to consider independent events, let $E^*\coloneqq\{e\in E(H):x\notin e,e\cap A_i\neq\varnothing\}$ and, for each $y\in Y_i$, let $d^*_{H}(y)\coloneqq|\{e\in E^*:y\in e\}|$.
By \ref{RNPcodegreeclose}, we have $d^*_{H}(y)=d_{H}(y)\pm\bigO(D/n)=(1\pm\bigO(n^{a-1}))D$\COMMENT{Note that, if $x\notin V(H)$, then $d_H^*(y)=d_H(y)$ and we are done.
Otherwise, we have that $d_H^*(y)=d_H(y)-d_H(x,y)$, and the conclusion follows by \ref{RNPcodegreeclose}.}.
Let $V^*\coloneqq Y_i\setminus V(E'\cap E^*)$.
For every $y\in Y_i$ we have that 
\[\mathbb{P}[y\notin V(E'\cap E^*)]=(1-\varepsilon/D)^{(1\pm\bigO(n^{a-1}))D}=(1\pm\bigO(n^{a-1}))e^{-\varepsilon},\]
hence $\mathbb{E}[|V^*|]=(1\pm\bigO(n^{a-1}))e^{-\varepsilon}\gamma|A_i|$\COMMENT{$\mathbb{E}[|V^*|]=(1\pm\bigO(n^{a-1}))e^{-\varepsilon}|Y_i|=(1\pm\bigO(n^{a-1}))e^{-\varepsilon}(1\pm\bigO(n^{a-1}))\gamma|A_i|=(1\pm\bigO(n^{a-1}))e^{-\varepsilon}\gamma|A_i|$.}.
Furthermore, the events $\{y\notin V(E'\cap E^*)\}_{y\in Y_i}$ are mutually independent, so by \cref{lem:Chernoff} we have that\COMMENT{Recall that, $|A_i|=\Theta(n)$.
We have
\[\mathbb{P}[|V^*|\neq(1\pm\bigO(n^{a-1}))e^{-\varepsilon}\gamma|A_i|]\leq2e^{-\Theta(n^{2a-2})(1\pm\bigO(n^{a-1}))e^{-\varepsilon}\gamma|A_i|/3}=e^{-\Theta(n^{2a-1})}.\]}
\begin{equation}\label{equa:RNdegbound0.5}
    \mathbb{P}[|V^*|\neq(1\pm\bigO(n^{a-1}))e^{-\varepsilon}\gamma|A_i|]\leq e^{-\Theta(n^{2a-1})}.
\end{equation}
As $V'\cap A_i\subseteq V^*$, we conclude that\COMMENT{We have
\[\mathbb{P}[|V'\cap A_i|\geq(1+\bigO(n^{a-1}))e^{-\varepsilon}\gamma|A_i|]\leq\mathbb{P}[|V^*|\geq(1+\bigO(n^{a-1}))e^{-\varepsilon}\gamma|A_i|]\leq e^{-\Theta(n^{2a-1})}.\]}
\begin{equation}\label{equa:RNdegbound1}
    \mathbb{P}[|V'\cap A_i|\geq(1+\bigO(n^{a-1}))e^{-\varepsilon}\gamma|A_i|]\leq e^{-\Theta(n^{2a-1})}.
\end{equation}

In order to obtain the lower tail concentration, observe that 
\[|V'\cap A_i|=|V^*|-|V^*\cap V(E'\setminus E^*)|,\]
so it will suffice to show that the last term in the previous expression is small with high probability.
Let $\hat E\coloneqq\{e\in E(H):|e\cap A_i|>1\}$.
Note that $|V^*\cap V(E'\setminus E^*)|\leq2^\ell|\hat E\cap E'|$\COMMENT{We can get the stronger bound that $|V^*\cap V(E'\setminus E^*)|\leq\ell|\hat E\cap E'|$, but it is not necessary.
To see this stronger bound, observe the following.
By \cref{rem:RNcodeg}, given an edge $e$ which contains two neighbours $y,y'$ of $x$, there is a unique vertex $z$ such that $\dist(x,z)=\dist(x,e)$.
Since $\dist(x,y)=\dist(x,y')=1$, this implies that $\dist(x,e)=0$ and, thus, $x\in e$.
Now, clearly, each edge which contains $x$ cannot contain more than $\ell$ of its neighbours (because edges correspond to $\ell$-cubes by the construction of the hypergraph).}, so it suffices to bound this quantity.
By \ref{RNPcodegreeclose}, $|\hat E|=\bigO(D)$\COMMENT{For any pair $y_1,y_2\in A_i$ we have $d_{H}(y_1,y_2)=\bigO(D/n^2)$, and there are $\Theta(n^2)$ such pairs of vertices.}.
Since edges are picked independently, we have that $\mathbb{E}[|\hat E\cap E'|]=\bigO(1)$ and, by \cref{lem:betaChernoff}\COMMENT{We have
\[\mathbb{P}[|\hat E\cap E'|>n^{1/2}]\leq(\bigO(1/n^{1/2}))^{\Theta(n^{1/2})}=e^{-\Theta(n^{1/2}\log n)}\leq e^{-\Theta(n^{1/2})}.\]}, $\mathbb{P}[|\hat E\cap E'|>n^{1/2}]\leq e^{-\Theta(n^{1/2})}$.
Combining this with \eqref{equa:RNdegbound0.5}, we conclude that $\mathbb{P}[|V'\cap A_i|\leq(1-\bigO(n^{a-1}))e^{-\varepsilon}\gamma|A_i|]\leq e^{-\Theta(n^{1/2})}$\COMMENT{We have that
\begin{align*}
    \mathbb{P}[|V'\cap A_i|\leq(1-\bigO(n^{a-1}))e^{-\varepsilon}\gamma|A_i|] &=\mathbb{P}[|V^*|-|V^*\cap V(E'\setminus E^*)|\leq(1-\bigO(n^{a-1}))e^{-\varepsilon}\gamma|A_i|]\\
    &\leq\mathbb{P}[|V^*|-2^\ell|\hat E\cap E'|\leq(1-\bigO(n^{a-1}))e^{-\varepsilon}\gamma|A_i|]\\
    &\leq\mathbb{P}[|V^*|\leq(1-\bigO(n^{a-1}))e^{-\varepsilon}\gamma|A_i|] + \mathbb{P}[|\hat E\cap E'|\geq n^{1/2}]\\
    &\leq e^{-\Theta(n^{2a-1})} + e^{-\Theta(n^{1/2})}=e^{-\Theta(n^{1/2})}.
\end{align*}
where in the third line we use the bound $a>3/4$, and in the fourth we combine \eqref{equa:RNdegbound0.5} and the last claim before this comment with the bound $a>3/4$.}.
Together with \eqref{equa:RNdegbound1}, the previous yields
\begin{equation}\label{equa:RNdegbound2}
    \mathbb{P}[|V'\cap A_i|\neq(1\pm\bigO(n^{a-1}))e^{-\varepsilon}\gamma|A_i|]\leq e^{-\Theta(n^{1/2})}.
\end{equation}
The statement of \ref{RNCremainingvertices} follows by a union bound over all $i\in[K]_0$.

\textbf{\ref{RNCmatchingsize}: On the number of vertices in $\boldsymbol{A_i}$ covered by the matching.}

We now prove \ref{RNCmatchingsize}.
Fix $i\in[K]_0$ and let $Y_i\coloneqq V(H)\cap A_i$.
Observe that 
\begin{equation}\label{equa:nibbleextraD}
    |V(E'')\cap A_i|=|(V(H)\setminus V')\cap A_i|-|V(E'\setminus E'')\cap A_i|.
\end{equation}
By \eqref{equa:RNdegbound2} and \ref{RNPvertexdistribution} we have that $|(V(H)\setminus V')\cap A_i|=(1\pm\bigO(n^{a-1}))(1-e^{-\varepsilon})\gamma|A_i|$ with probability at least $1-e^{-\Theta(n^{1/2})}$, so let us consider the last term in \eqref{equa:nibbleextraD}.

Given any vertex $y\in Y_i$, by abusing notation, let $d_{E'}(y)\coloneqq|\{e\in E':y\in e\}|$.
Observe that $y\in V(E'\setminus E'')$ if and only if $d_{E'}(y)\geq2$ or $d_{E'}(y)=1$ and, for the edge $e\in E'$ such that $y\in e$, there exists $z\in e\setminus\{y\}$ such that $d_{E'\setminus\{e\}}(z)\geq1$.
Let $\mathcal{B}(y)$ be the event that, conditioned on $d_{E'}(y)=1$, there exists such a vertex $z$.
Then, for any $y\in Y_i$ we have
\begin{equation}\label{equa:comment741}
\mathbb{P}[y\in V(E'\setminus E'')]\leq\mathbb{P}[d_{E'}(y)\geq2]+\mathbb{P}[d_{E'}(y)=1]\mathbb{P}[\mathcal{B}(y)].
\end{equation}
Observe that, by \ref{RNPdegree}, $d_{E'}(y)\sim\mathrm{Bin}((1\pm\bigO(n^{a-1}))D,\varepsilon/D)$.
Then, it is easy to check that
\begin{equation}\label{equa:comment742}
\mathbb{P}[d_{E'}(y)=1]=(1\pm\bigO(n^{a-1}))\varepsilon e^{-\varepsilon}\quad\text{ and }\quad\mathbb{P}[d_{E'}(y)\geq2]=1-(1\pm\bigO(n^{a-1}))(1+\varepsilon)e^{-\varepsilon}.
\end{equation}
By a union bound and the fact that $\mathbb{P}[d_{E'\setminus e}(z)\geq1]\leq\mathbb{P}[d_{E'}(z)\geq1]$ for every $e\in E(H)$ and $z\in e\setminus\{y\}$, we also have that
\begin{equation}\label{equa:comment743}
    \mathbb{P}[\mathcal{B}(y)]\leq(1\pm\bigO(n^{a-1}))(2^\ell-1)(1-e^{-\varepsilon}).
\end{equation}
Combining \eqref{equa:comment741}--\eqref{equa:comment743}, for $n$ sufficiently large we have that $\mathbb{P}[y\in V(E'\setminus E'')]\leq2^\ell\varepsilon^2$.
Hence, by considering all $y\in Y_i$ and \ref{RNPvertexdistribution}, we conclude that
\begin{equation}\label{equa:RNdegexpect}
    \mathbb{E}[|V(E'\setminus E'')\cap A_i|]\leq(1+\bigO(n^{a-1}))2^\ell\varepsilon^2\gamma|A_i|.
\end{equation}

In order to prove concentration, we will resort to Talagrand's inequality.
Consider $X\coloneqq|V(E'\setminus E'')\cap A_i|$.
This is a random variable on the probability space given by the product of the probability spaces associated with each edge of $H$ being present in $E'$.
In this setting, it is easy to see that $X$ is a $\ell2^\ell$-Lipschitz function\COMMENT{Let us consider first the variable $Y=|E'\setminus E''|$.
We must consider two cases.\\
First, consider the case where, given the collection $E'$, we add a new edge $e$.
The value of $Y$ will change by $1$ for each edge $e'\in E'$ which intersects $e$ and does not intersect any other edge in $E'$.
The set of such edges $e'$ must form a matching, so the number of such edges is bounded from above by the size of $e$, that is, $2^\ell$.\\
Similarly, consider now the case where, given the collection $E'$, we remove an edge $e$.
We want to know how many edges $e'$ intersect $e$ and no other edge in $E'$ (as only these will change the value of $Y$, by one unit each).
The same bound as above holds trivially.\\
Thus, $Y$ is $2^\ell$-Lipschitz.
Since no edge can intersect the neighbourhood of $x$ in more than $\ell$ vertices, it follows that the value of $X$ cannot change by more than $\ell\cdot 2^\ell$.}.
Furthermore, $X$ is $h$-certifiable for $h\colon\mathbb{N}\to\mathbb{N}$ given by $h(s)=2s$.
(Indeed, consider any possible outcome $E'$ and assume $X(E')\geq s$ for some $s\in\mathbb{N}$.
This means there is a set of at least $s$ vertices in $A_i$, each of which belongs to an edge of $E'$ which, in turn, intersects some other edge of $E'$.
Fix any $s$ such vertices and, for each of these, fix one edge of $E'$ containing the vertex and one other edge of $E'$ intersecting the previous edge.
This gives a set $E^\diamond$ of at most $2s$ edges such that any other possible outcome $E^*$ with $E^\diamond\subseteq E^*$ satisfies $X(E^*)\geq s$.)
Thus, by \cref{lem: Talagrand}, for any real values $b$ and $t$ we have that
\[\mathbb{P}\left[X\leq b-t\ell2^\ell\sqrt{2b}\right]\mathbb{P}[X\geq b]\leq e^{-t^2/4}.\]
By considering the change of variables $c=b-t\ell2^\ell\sqrt{2b}$, we conclude that, for any reals $c$ and $t$\COMMENT{This follows by solving a polynomial of degree $2$. 
We want to rewrite $b$ as a function of $c$ and $t$, so we consider the polynomial $x^2-t\ell2^{\ell+1/2}x-c=0$, where $x=b^{1/2}$, and solve it.
As we know $x$ must be positive, there is a unique solution.},
\begin{equation}\label{equa:taldeg2}
    \mathbb{P}\left[X\leq c\right]\mathbb{P}\left[X\geq\left(t\ell2^{\ell+1/2}+\left(t^2\ell^22^{2\ell+1}+4c\right)^{1/2}\right)^2/4\right]\leq e^{-t^2/4}.
\end{equation}
Let $c\coloneqq3\cdot2^{\ell-1}\varepsilon^2\gamma|A_i|$ and $t\coloneqq\Theta(n^{a-1/2})$.
By Markov's inequality, we have that $\mathbb{P}[X\leq c]\geq1/4$\COMMENT{By Markov's inequality, we have that $\mathbb{P}[X\geq c]\leq\mathbb{E}[X]/c\leq2/3+\bigO(n^{a-1})$, and the claim follows by taking the complement.} for $n$ sufficiently large.
By substituting these into \eqref{equa:taldeg2}, we conclude that\COMMENT{We have that
\begin{align*}
    \left(\frac{t\ell2^{\ell+1/2}+\left(t^2\ell^22^{2\ell+1}+4c\right)^{1/2}}{2}\right)^2&=\left(\frac{\left(4c+\bigO(n^{2a-1})\right)^{1/2}+\bigO(n^{a-1/2})}{2}\right)^2\\
    &=\left(\frac{\left(4c(1+\bigO(n^{2a-2}))\right)^{1/2}+\bigO(n^{a-1/2})}{2}\right)^2\\
    &=\left(\frac{2c^{1/2}(1+\bigO(n^{2a-2}))+\bigO(n^{a-1/2})}{2}\right)^2\\
    &=\left(c^{1/2}(1+\bigO(n^{2a-2}))(1+\bigO(n^{a-1}))\right)^2\\
    &=(1+\bigO(n^{a-1}))c.
\end{align*}
}
\begin{equation}\label{equa:RNmatchingsizepre}
    \mathbb{P}\left[X\geq(1+\bigO(n^{a-1}))c\right]\leq e^{-\Theta(n^{2a-1})}.
\end{equation}
From this and \eqref{equa:RNdegbound2}, it follows that\COMMENT{Using the fact that $e^{-\varepsilon}\leq1-\varepsilon+\varepsilon^2/2$, we have
\begin{align*}
    &\mathbb{P}[|V(E'')\cap A_i|\leq\varepsilon(1-2^{\ell+1}\varepsilon)\gamma|A_i|] =\mathbb{P}[|V(E'')\cap A_i|\leq(\varepsilon-2^{\ell+1}\varepsilon^2)\gamma|A_i|]\\
    \leq\,&\mathbb{P}[|(V(H)\setminus V')\cap A_i|\neq(1\pm\bigO(n^{a-1}))(1-e^{-\varepsilon})\gamma|A_i|]+\mathbb{P}\left[X\geq(1+\bigO(n^{a-1}))c\right],
\end{align*}
where the second line follows from the observation that $(1\pm\bigO(n^{a-1}))(1-e^{-\varepsilon}-3\cdot2^{\ell-1}\varepsilon^2)\geq (\varepsilon-2^{\ell+1}\varepsilon^2)$.
}
\[\mathbb{P}[|V(E'')\cap A_i|\leq\varepsilon(1-2^{\ell+1}\varepsilon)\gamma|A_i|]\leq e^{-\Theta(n^{1/2})}.\]
The statement of \ref{RNCmatchingsize} follows by a union bound over all $i\in[K]_0$.

\textbf{\ref{RNCsingledirection}: On the significance of $\boldsymbol{E'}$ in any direction.}

In order to prove \ref{RNCsingledirection}, we first observe that there are `few' edges in $E_x(H,A_0)$ which use any given direction.
Indeed, given any vertex $y\in V(H)\cap A_0$ and any direction $\hat e\in\mathcal{D}(\cQ^n)$, the number of edges $e\in E(H)$ containing $y$ and such that $\hat e\in\mathcal{D}(e)$ equals the codegree of $y$ and $y+\hat e$.
Therefore, by \ref{RNPcodegreeclose}, there are $\bigO(D/n)$ such edges and, adding over all vertices $y\in V(H)\cap A_0$, we conclude that $|\Sigma(E_x(H,A_0),\{\hat e\},1)|=\bigO(D)$.
Since $|\Sigma(E_x'(A_0),\{\hat e\},1)|\sim\mathrm{Bin}(|\Sigma(E_x(H,A_0),\{\hat e\},1)|,\varepsilon/D)$, it immediately follows that $\mathbb{E}[|\Sigma(E_x'(A_0),\{\hat e\},1)|]=\bigO(1)$ and, by \cref{lem:betaChernoff}\COMMENT{By \cref{lem:betaChernoff}, we have that
\begin{align*}
    \mathbb{P}[|\Sigma(E_x'(A_0),\{\hat e\},1)|=\Omega(n^{1/2})]&=\mathbb{P}[|\Sigma(E_x'(A_0),\{\hat e\},1)|=\Omega(n^{1/2})\mathbb{E}[|\Sigma(E_x'(A_0),\{\hat e\},1)|]]\\
    &\leq\left(\frac{e}{\Omega(n^{1/2})}\right)^{\Theta(n^{1/2})}=\left(\bigO(n^{-1/2})\right)^{\Theta(n^{1/2})}\leq e^{-\Theta(n^{1/2}\log n)}\leq e^{-\Theta(n^{1/2})}.
\end{align*}
},
\[\mathbb{P}[|\Sigma(E_x'(A_0),\{\hat e\},1)|=\Omega(n^{1/2})]\leq e^{-\Theta(n^{1/2})}.\]
The statement of \ref{RNCsingledirection} follows by a union bound over all directions $\hat e\in\mathcal{D}(\cQ^n)$.

\textbf{\ref{RNCdegree}: On the degrees in $\boldsymbol{H'}$.}

We now want to bound the degrees of vertices in $H'$ in order to prove \ref{RNCdegree}.
Consider any vertex $y\in V(H)$ such that $\dist(x,y)\leq k-\ell$.
Condition on the event that $y\in V'$.
First, observe that, by \ref{RNPdegree} and \ref{RNPcodegreeclose}\COMMENT{Since $y\in V'$, we are implicitly conditioning on the event that none of the edges containing $y$ belong to $E'$.
Each edge containing $y$ (of which there are $(1\pm\bigO(n^{a-1}))D$ by \ref{RNPdegree}, since $\dist(x,y)\leq k$) will belong to $H'$ if and only if no edge containing each of its other vertices is chosen by $E'$.
As all vertices in these edges have (roughly) the same degree (here we use again \ref{RNPdegree} together with the fact that $\dist(x,y)\leq k-\ell$, since all vertices $v$ in edges which contain $y$ therefore satisfy $d(x,v)\leq k$) and the codegrees are small (by \ref{RNPcodegreeclose}), this occurs with probability $(1-\varepsilon/D)^{(2^\ell-1)(1\pm\bigO(n^{a-1}))D}$.
Therefore,
\begin{align*}
    \mathbb{E}[d_{H'}(y)]&=(1\pm\bigO(n^{a-1}))D(1-\varepsilon/D)^{(2^\ell-1)(1\pm\bigO(n^{a-1}))D}\\
    &=(1\pm\bigO(n^{a-1}))D(e^{-\varepsilon/D}\pm\bigO(1/D^2))^{(2^\ell-1)(1\pm\bigO(n^{a-1}))D}\\
    &=(1\pm\bigO(n^{a-1}))D(1\pm\bigO(1/D))e^{-\varepsilon(2^\ell-1)(1\pm\bigO(n^{a-1}))}\\
    &=(1\pm\bigO(n^{a-1}))D(1\pm\bigO(1/D))(1\pm\bigO(n^{a-1}))e^{-\varepsilon(2^\ell-1)}=(1\pm\bigO(n^{a-1}))e^{-(2^\ell-1)\varepsilon}D.
\end{align*}
},
\begin{equation}\label{equa:RNexpdegree}
    \mathbb{E}[d_{H'}(y)]=(1\pm\bigO(n^{a-1}))D(1-\varepsilon/D)^{(2^\ell-1)(1\pm\bigO(n^{a-1}))D}=(1\pm\bigO(n^{a-1}))e^{-(2^\ell-1)\varepsilon}D.
\end{equation}

In order to bound the probability that $d_{H'}(y)$ deviates from its expectation, we will apply \cref{lem: AKS}.
Observe that the value of $d_{H'}(y)$ is determined by the presence or absence of the edges of $E^\bullet\coloneqq\{e\in E(H):\text{there exists }e'\in E(H)\text{ such that }y\in e'\setminus e,e\cap e'\neq\varnothing\}$ in $E'$.
Note that, for each $e\in E^\bullet$, the maximum possible change in the value of $d_{H'}(y)$ due to the presence or absence of $e$ is $c_e\coloneqq|\{e'\in E(H):y\in e',e\cap e'\neq\varnothing\}|$.
Let $C\coloneqq\max_{e\in E^\bullet}{c_e}$ and $\sigma^2\coloneqq\sum_{e\in E^\bullet}(\varepsilon/D)(1-\varepsilon/D)c_e^2$.
We must now estimate the value of $\sigma$.

Partition $E^\bullet$ into sets $E_i$, $i\in[\ell]$, given by $E_i\coloneqq\{e\in E^\bullet:\dist(y,e)=i\}$\COMMENT{Observe that, for every edge $e$ such that $\dist(A,e)>\ell$ we have that $e\notin E^\bullet$ by definition.}.
Observe that, by \ref{RNPdegree} and \cref{rem:basicequationIneedhere}, for all $i\in[\ell]$ we have\COMMENT{Maybe here, being strict, we have to add more reasons why this holds.}
\begin{equation}\label{equa:sizeEi}
    |E_i|=\bigO(n^iD).
\end{equation}
Furthermore, for each $i\in[\ell]$ and each $e\in E_i$, it follows from \cref{rem:RNcodeg} and \ref{RNPcodegreeclose} that\COMMENT{By \cref{rem:RNcodeg}, for each $i\in[\ell]$ and each $e\in E_i$ there exists a unique vertex $z\in e$ such that $\dist(y,e)=\dist(y,z)=i$.
Also by \cref{rem:RNcodeg}, we then have that $c_e=d_H(y,z)$, and the claim follows by \ref{RNPcodegreeclose}.}
\begin{equation}\label{equa:boundce}
    c_e=\bigO(D/n^i).
\end{equation}

In order to apply \cref{lem: AKS}, we will need to show that $\sigma$ is not too small.
For this, we claim that
\begin{equation}\label{equa:c1}
\begin{minipage}[c]{0.65\textwidth}
there exist $\Theta(nD)$ edges $e\in E_1$ such that $c_e=\Theta(D/n)$.
\end{minipage}\ignorespacesafterend 
\end{equation} 
Indeed, an averaging argument using \eqref{equa:codegsum} together with \ref{RNPdegree} shows that there are $\Theta(n)$ vertices $z\in V(H)\cap N_{\cQ^n}(y)$ such that $d_H(y,z)=\Theta(D/n)$\COMMENT{Assume that for every $z\in V(H)\cap N_{\cQ^n}(y)$ but at most $o(n)$ of them we had that $d_H(y,z)=o(D/n)$.
Then, by \eqref{equa:codegsum} we have that $d_H(y)\leq |V(H)\cap N_{\cQ^n}(y)|\cdot o(D/n)+o(n)\cdot\Theta(D/n)=\bigO(n)\cdot o(D/n)+o(n)\cdot\Theta(D/n)=o(D)$, a contradiction to \ref{RNPdegree}.
So we must have that at least $\Theta(n)$ vertices $z\in V(H)\cap N_{\cQ^n}(y)$ satisfy that $d_H(y,z)=\Theta(D/n)$.
(In particular, this implies that in $H$ all neighbourhoods must be linear.)}.
Let $Z$ be the set of all those vertices $z$.
For each $z\in V(H)\cap N_{\cQ^n}(y)$, let $E_1(z)\coloneqq\{e\in E_1:z\in e\}$.
By \cref{rem:RNcodeg}, for every $e\in E_1(z)$, $z$ is the unique vertex in $e$ such that $\dist(y,z)=1$, so this gives a partition of $E_1$.
By \ref{RNPdegree} and \ref{RNPcodegreeclose}, we have that $|E_1(z)|=\Theta(D)$ for every $z\in Z$.
Again by \cref{rem:RNcodeg}, for every $z\in Z$ and every $e\in E_1(z)$ we have that $c_e=d_H(y,z)$.
\eqref{equa:c1} now follows.

In particular, \eqref{equa:c1} combined with \eqref{equa:boundce} shows that $C=\Theta(D/n)$.
Combining \eqref{equa:sizeEi}--\eqref{equa:c1}, it follows that\COMMENT{We have that 
\[\sigma^2=\sum_{e\in E^\bullet}\frac\varepsilon D\left(1-\frac\varepsilon D\right)c_e^2=\sum_{i=1}^\ell\sum_{e\in E_i}\frac\varepsilon D\left(1-\frac\varepsilon D\right)c_e^2.\]
Now, by \eqref{equa:sizeEi} and \eqref{equa:boundce}, we have that
\[\sigma^2=\sum_{i=1}^\ell|E_i|\frac{\varepsilon}{D}\bigO(D^2/n^{2i})=\sum_{i=1}^\ell\bigO(D^2/n^i)=\bigO(D^2/n).\]
The lower bound follows similarly by noting that 
\[\sum_{e\in E_1}\frac\varepsilon D\left(1-\frac\varepsilon D\right)c_e^2=\Theta(D^2/n)\]
(this follows from \eqref{equa:c1}).
}
\[\sigma^2=\Theta(D^2/n).\]
Now, by setting $\alpha\coloneqq\Theta(n^{a-1/2})$, we observe that $\alpha=o(\sigma/C)$ and, thus, by \cref{lem: AKS} and \eqref{equa:RNexpdegree}\COMMENT{We have that
\[\mathbb{P}[d_{H'}(y)\neq(1\pm\bigO(n^{a-1}))e^{-(2^\ell-1)\varepsilon}D]=\mathbb{P}[d_{H'}(y)\neq\mathbb{E}[d_{H'}(y)]\pm\alpha\sigma]\leq 2e^{-\alpha^2/4}=e^{-\Theta(n^{2a-1})}.\]
},
\[\mathbb{P}[d_{H'}(y)\neq(1\pm\bigO(n^{a-1}))e^{-(2^\ell-1)\varepsilon}D]\leq e^{-\Theta(n^{2a-1})}.\]
The statement of \ref{RNCdegree} follows by a union bound over all vertices $y\in V(H)\cap B_{\cQ^n}^{k-\ell}(x)$\COMMENT{Formally, we have proved that, for any $y\in V(H)\cap B_{\cQ^n}^{k-\ell}(x)$, 
\[\mathbb{P}[d_{H'}(y)\neq(1\pm\bigO(n^{a-1}))e^{-(2^\ell-1)\varepsilon}D\mid y\in V']\leq e^{-\Theta(n^{2a-1})}.\]
It follows that
\[\mathbb{P}[y\in V'\wedge d_{H'}(y)\neq(1\pm\bigO(n^{a-1}))e^{-(2^\ell-1)\varepsilon}D]\leq e^{-\Theta(n^{2a-1})}.\]
The statement of \ref{RNCdegree} fails if there is some $y\in V'$ such that $d_{H'}(y)\neq(1\pm\bigO(n^{a-1}))e^{-(2^\ell-1)\varepsilon}D$.
Thus, by a union bound, the probability that the statement fails is at most $(n+1)^{k-\ell}e^{-\Theta(n^{2a-1})}=e^{-\Theta(n^{2a-1})}$.}.

\textbf{\ref{RNCdirectionsleft}: On the significance of $\boldsymbol{H'}$ in any large set of directions.}

We now turn our attention to \ref{RNCdirectionsleft}.
Fix $i\in[K]_0$ and $S\subseteq\mathcal{D}(\cQ^n)$ with $|S|\geq\delta n$.
By \eqref{equa:RNchosenedges} we have that 
\begin{equation}\label{equa:RNchosenedgesbis}
    |E_x'(A_i)|=(1\pm\bigO(n^{a-1}))\varepsilon\gamma|A_i|
\end{equation} 
with probability at least $1-e^{-\Theta(n^{2a-1})}$.
Furthermore, by \eqref{equa:RNdegbound2}, we have that 
\begin{equation}\label{equa:RNdegbound2bis}
    |V'\cap A_i|=(1\pm\bigO(n^{a-1}))e^{-\varepsilon}\gamma|A_i|
\end{equation}
with probability at least $1-e^{-\Theta(n^{1/2})}$.
Reveal $E_x'(A_i)$, as well as all edges in $E'$ which contain $x$ and intersect $A_i$, and condition on the event that \eqref{equa:RNchosenedgesbis} and \eqref{equa:RNdegbound2bis} hold (note that this event is determined by the edges we have revealed).
For the remainder of the proof of \ref{RNCdirectionsleft}, all probabilistic statements refer to probabilities when revealing all other edges in $E'$.

Let $X'\coloneqq|\Sigma(E_x(H',A_i),S,\ell^{1/2})|$.
Note that $X'$ is a sum of indicator random variables, one for each edge in $\Sigma(E_x(H,A_i),S,\ell^{1/2})$; we will refer to those edges $e\in\Sigma(E_x(H,A_i),S,\ell^{1/2})$ for which we have $\mathbb{P}[e\in\Sigma(E_x(H',A_i),S,\ell^{1/2})]\neq0$ as \emph{potential edges}.
The set of potential edges is denoted by $E_P$.
We now want to prove a lower bound on $|E_P|$.
We know that $\Sigma(E_x(H',A_i),S,\ell^{1/2})\subseteq\Sigma(E_x(H,A_i),S,\ell^{1/2})$.
By \ref{RNPdirections} we have that $|\Sigma(E_x(H,A_i),S,\ell^{1/2})|\geq(1-\bigO(n^{a-1}))c\gamma|A_i|D$.
Any edge of $\Sigma(E_x(H,A_i),S,\ell^{1/2})$ whose endpoint in $A_i$ does not lie in $V'$ is not a potential edge.
By \eqref{equa:RNdegbound2bis}, \ref{RNPdegree} and \ref{RNPvertexdistribution}, the number of such edges is at most $(1+\bigO(n^{a-1}))(1-e^{-\varepsilon})\gamma|A_i|D$.
Furthermore, some of the edges in $E_x'(A_i)$ may intersect other edges in $\Sigma(E_x(H,A_i),S,\ell^{1/2})$ (and, if this happens, the latter are not potential edges)\COMMENT{We only have to consider the revealed edges from $E_x'(A_i)$ here, and not any other revealed edge $e$ which contains $x$ and intersects $A_i$ (since if $e\cap e'\neq\varnothing$ for some $e'\in\Sigma(E_x(H,A_i),S,\ell^{1/2})$ then the unique vertex in $e'\cap A_i$ must lie in $e$ and so $e'$ has already been counted above).}.
By \cref{rem:RNExCAP} and \eqref{equa:RNchosenedgesbis}, the number of such non-potential edges is $\bigO(D)$\COMMENT{We will obtain an upper bound by looking at the supersets (ignoring the directions given by $S$).
By \eqref{equa:RNchosenedgesbis}, there are $\Theta(n)$ edges in $E'_x(A_i)$ already revealed.
For any of these edges $e\in E'_x(A_i)$, apply \cref{rem:RNExCAP} to obtain that there are $\bigO(n^{\ell-1})$ edges in $E_x(H,A_i)$ that intersect $e$. 
By adding over all choices for the edge $e\in E'_x(A_i)$, we conclude that the number of edges `forbidden' by $E'_x(A_i)$ is $\Theta(n)\cdot\bigO(n^{\ell-1})=\bigO(n^\ell)=\bigO(D)$, which is the conclusion we wanted.
Note that in this development we have ignored the edges that we have revealed and contain $x$; this is because there can be no $e\in E_x(H,A_i)$ with an endpoint in $V'$ which intersects any such edge.}.
Combining these bounds, we conclude that $|E_P|\geq(1-\bigO(n^{a-1}))(c-(1-e^{-\varepsilon}))\gamma|A_i|D$.
Now, each of these potential edges contributes to $X'$ if and only if none of its vertices lie in any edge in $E'$.
By \ref{RNPdegree} and \ref{RNPcodegreeclose}, it follows that, for each $e\in E_P$, $\mathbb{P}[e\in\Sigma(E_x(H',A_i),S,\ell^{1/2})]=(1\pm\bigO(n^{a-1}))e^{-(2^\ell-1)\varepsilon}$\COMMENT{This follows from an argument in the proof of \ref{RNCdegree}. Note that here is where we need that $k>\ell$.} and, therefore,
\begin{equation}\label{equa:RNexpdegree2}
    \mathbb{E}[X']\geq(1-\bigO(n^{a-1}))(c-(1-e^{-\varepsilon}))e^{-(2^\ell-1)\varepsilon}\gamma|A_i|D.
\end{equation}

In order to prove concentration we will resort once more to \cref{lem: AKS}.
Let $E^\bullet_i\coloneqq\{e\in E(H):e\cap A_i=\varnothing\text{ and there exists }e'\in\Sigma(E_x(H,A_i),S,\ell^{1/2})\text{ such that }e\cap e'\neq\varnothing\}$\COMMENT{Note that $E^\bullet_i$ does not contain any edge containing $x$, even if this is not obvious from the definition, as there may be edges containing $x$ which do not intersect $A_i$.
However, all such edges cannot intersect any edges intersecting $A_i$ and not containing $x$ (see \cref{rem:RNcodeg}).
Therefore, by the condition $e\cap e'\neq\varnothing$, all edges containing $x$ are excluded.}.
The value of $X'$ is determined uniquely by the presence or absence of the edges of $E^\bullet_i$ in $E'$.
For each $e\in E^\bullet_i$, the maximum change in the value of $X'$ due to the presence or absence of $e$ can be bounded by $c_e\coloneqq|\{e'\in \Sigma(E_x(H,A_i),S,\ell^{1/2}):e'\cap e\neq\varnothing\}|$.
Let $C\coloneqq\max_{e\in E^\bullet_i}{c_e}$ and $\sigma^2\coloneqq\sum_{e\in E^\bullet_i}(\varepsilon/D)(1-\varepsilon/D)c_e^2$.
We must now estimate the value of $\sigma$.

Partition $E^\bullet_i$ into sets $E_i^j$, $j\in[\ell]$, given by $E_i^j\coloneqq\{e\in E^\bullet_i:\dist(e,A_i)=j\}$\COMMENT{Observe that, for every edge $e$ such that $\dist(e,A_i)>\ell$ we have that $e\notin E^\bullet_i$ by definition.}.
Observe that, by \ref{RNPdegree} and \cref{rem:basicequationIneedhere}, for all $j\in[\ell]$ we have
\begin{equation}\label{equa:sizeEi2}
    |E_i^j|=\bigO(n^{j+1}D).
\end{equation}
Furthermore, for each $j\in[\ell]$ and each $e\in E_i^j$, we claim that
\begin{equation}\label{equa:boundce2}
    c_e=\bigO(D/n^j).
\end{equation}
This follows by a case analysis combining \cref{rem:RNcodeg,rem:RNcodeg2} and \ref{RNPcodegreeclose}.
Indeed, we consider three cases for each $e\in E_i^j$, depending on its distance to $x$.
For each $j\in[\ell]$, let $E_i^{j+}\coloneqq\{e\in E_i^j:\dist(x,e)=j+1\}$, $E_i^{j\bullet}\coloneqq\{e\in E_i^j:\dist(x,e)=j\}$ and $E_i^{j-}\coloneqq\{e\in E_i^j:\dist(x,e)=j-1\}$.
For each $e\in E_i^{j+}$, the claimed bound on $c_e$ follows from \cref{rem:RNcodeg} and \ref{RNPcodegreeclose}.
For each $e\in E_i^{j-}$, \cref{rem:RNcodeg} can be used to show that $c_e=0$.
Finally, consider any $e\in E_i^{j\bullet}$.
By \cref{rem:RNcodeg2}, there are at most $\ell$ vertices $z\in A_i$ such that $\dist(z,e)=j$.
By \cref{rem:RNcodeg}, for each such $z$ there is a unique $y(z)\in e$ such that $\dist(z,y(z))=\dist(z,e)=j$.
Again by \cref{rem:RNcodeg}, we have $c_e\leq \sum_{z\in A_i:\dist(z,e)=j}d_H(z,y(z))$ (by \cref{rem:RNcodeg2}, any edge in $E_x(H,A_i)$ intersecting $e$ must contain some $x\in A_i$ with $\dist(z,e)=j$), and the claim follows by \ref{RNPcodegreeclose}.

In particular, we claim that
\begin{equation}\label{equa:c2}
\begin{minipage}[c]{0.65\textwidth}
there exist $\Theta(n^2D)$ edges $e\in E_i^1$ such that $c_e=\Theta(D/n)$.
\end{minipage}\ignorespacesafterend 
\end{equation} 
Indeed, by \ref{RNPdegree}, \ref{RNPvertexdistribution} and \ref{RNPdirections}, there are at least $c\gamma|A_i|/2$ vertices $y\in A_i\cap V(H)$ such that $d_{H,S,\ell^{1/2},x}(y)\geq cD/2$\COMMENT{Indeed, assume this is not true.
Then, 
\begin{align*}
    |\Sigma(E_x(H,A_i),S,\ell^{1/2})|&=\sum_{y\in A_i\cap V(H)}d_{H,S,\ell^{1/2},x}(y)<\frac{c}{2}\gamma|A_i|(1+\bigO(n^{a-1}))D+(1+\bigO(n^{a-1}))\left(1-\frac{c}{2}\right)\gamma|A_i|\frac{c}{2}D\\
    &=(1\pm\bigO(n^{a-1}))(1/2+(1-c/2)/2)\gamma|A_i|cD<(1-\bigO(n^{a-1}))\gamma|A_i|cD,
\end{align*}
where in the first inequality we combine \ref{RNPdegree} and \ref{RNPvertexdistribution} with the assumption and in the last inequality we observe that, since $c>0$, we have $1/2+(1-c/2)/2<1$.
This is a contradiction on \ref{RNPdirections}.}.
Let $U_i$ denote the set of these vertices.
Then, an averaging argument using \eqref{equa:codegsum} together with \ref{RNPdegree} shows that, for each $y\in U_i$, there are $\Theta(n)$ vertices $z\in V(H)\cap (N_{\cQ^n}(y)\setminus\{x\})$ such that $d_{H,S,\ell^{1/2},x}(y,z)=\Theta(D/n)$\COMMENT{Assume that for every $z\in V(H)\cap N_{\cQ^n}(y)$ but at most $o(n)$ of them we had that $d_{H,S,\ell^{1/2},x}(y,z)=o(D/n)$ (note that, by \ref{RNPdegree}, we know that there are $\Theta(n)$ vertices $z\in V(H)\cap N_{\cQ^n}(y)$).
Then, by an equation analogous to \eqref{equa:codegsum} when restricting to the subset $\Sigma(E(H,A_i),S,\ell^{1/2})$, we have that $d_{H,S,\ell^{1/2},x}(y)\leq \Theta(n)\cdot o(D/n)+o(n)\cdot\Theta(D/n)=o(D)$, a contradiction.
So we must have that at least $\Theta(n)$ vertices $z\in V(H)\cap N_{\cQ^n}(y)$ satisfy that $d_{H,S,\ell^{1/2},x}(y,z)=\Theta(D/n)$.}.
For each $y\in U_i$, let $Z_i(y)$ be the set of such vertices.
Now, fix any $y\in U_i$ and, for each $z\in Z_i(y)$, let $E_i^1(z)\coloneqq\{e\in E_i^1:z\in e\}$.
By \ref{RNPdegree} and \ref{RNPcodegreeclose}, we have that $|E^1_i(z)|=\Theta(D)$ for every $z\in Z_i(y)$\COMMENT{Note that here we are using the fact that $z$ has at most two neighbours in $A_i$.}; furthermore, by \cref{rem:RNcodeg}, for every $e\in E^1_i(z)$, $z$ is the unique vertex in $e$ such that $\dist(y,z)=\dist(y,e)=1$.
Then, for every $z\in Z_i(y)$ and every $e\in E_1(z)$ we have that
\[c_e\geq d_{H,S,\ell^{1/2},x}(y,z)=\Theta(D/n).\]
\eqref{equa:c2} now follows by considering all vertices $y\in U_i$\COMMENT{Here, we use again the fact that $z$ has at most two neighbours in $A_i$.}.

In particular, \eqref{equa:c2} combined with \eqref{equa:boundce2} shows that $C=\Theta(D/n)$.
Combining \eqref{equa:sizeEi2}--\eqref{equa:c2}, it follows that\COMMENT{We have that 
\[\sigma^2=\sum_{e\in E_i^\bullet}\frac\varepsilon D\left(1-\frac\varepsilon D\right)c_e^2=\sum_{j=1}^\ell\sum_{e\in E_i^j}\frac\varepsilon D\left(1-\frac\varepsilon D\right)c_e^2.\]
Now, by \eqref{equa:sizeEi2} and \eqref{equa:boundce2}, we have that
\[\sigma^2=\sum_{j=1}^\ell|E_i^j|\frac{\varepsilon}{D}\bigO(D^2/n^{2j})=\sum_{j=1}^\ell\bigO(D^2/n^{j-1})=\bigO(D^2).\]
The lower bound follows similarly by noting that 
\[\sum_{e\in E_i^1}\frac\varepsilon D\left(1-\frac\varepsilon D\right)c_e^2=\Theta(D^2)\]
(this follows from \eqref{equa:c2}).
}
\[\sigma^2=\Theta(D^2).\]
Now, by setting $\alpha\coloneqq n^{a-1}\mathbb{E}[X']/\sigma=\Theta(n^a)$, we observe that $\alpha=o(\sigma/C)$ and, thus, by \cref{lem: AKS} and \eqref{equa:RNexpdegree2}\COMMENT{We have that
\begin{align*}
    \mathbb{P}[X'<(c-\varepsilon)e^{-(2^\ell-1)\varepsilon}\gamma|A_i|D]&\leq\mathbb{P}[X'<(1-\alpha\sigma/\mathbb{E}[X'])(1-\bigO(n^{a-1}))(c-(1-e^{-\varepsilon}))e^{-(2^\ell-1)\varepsilon}\gamma|A_i|D]\\
    &\leq\mathbb{P}[X'<(1-\alpha\sigma/\mathbb{E}[X'])\mathbb{E}[X']]\\
    &\leq\mathbb{P}[|X'-\mathbb{E}[X']|>\alpha\sigma]\leq2e^{-\alpha^2/4}= e^{-\Theta(n^{2a})}.
\end{align*}
},
\[\mathbb{P}[X'<(c-\varepsilon)e^{-(2^\ell-1)\varepsilon}\gamma|A_i|D]\leq e^{-\Theta(n^{2a})}.\]
Since this holds for every $S\subseteq\mathcal{D}(\cQ^n)$ with $|S|\geq\delta n$, by a union bound we conclude that the same holds simultaneously for every such set $S$\COMMENT{Here we use the fact that $a>1/2$.}.
Recall, however, that this holds after conditioning on the event that \eqref{equa:RNchosenedgesbis} and \eqref{equa:RNdegbound2bis} hold, which happens with probability $1-e^{-\Theta(n^{1/2})}$.
Taking this into account and using a union bound over all choices of $i\in[K]_0$, the statement of \ref{RNCdirectionsleft} follows.

\textbf{\ref{RNCdirectionscovered}: On the significance of $\boldsymbol{E'}$ in a large fixed set of directions.}

We finally turn our attention to \ref{RNCdirectionscovered}.
Fix $i\in[K]_0$.
Let $U_i$ denote the set of vertices $y\in A_i\cap V(H)$ such that $d_{H,S,\ell^{1/2},x}(y)\geq cD/2$.
By \ref{RNPdegree}, \ref{RNPvertexdistribution} and \ref{RNPdirections} we have $|U_i| \ge c\gamma|A_i|/2$.
\COMMENT{Indeed, assume this is not true.
Then, 
\begin{align*}
    |\Sigma(E_x(H,A_i),S,\ell^{1/2})|&=\sum_{y\in A_i\cap V(H)}d_{H,S,\ell^{1/2},x}(y)<\frac{c}{2}\gamma|A_i|(1+\bigO(n^{a-1}))D+(1+\bigO(n^{a-1}))\left(1-\frac{c}{2}\right)\gamma|A_i|\frac{c}{2}D\\
    &=(1\pm\bigO(n^{a-1}))(1/2+(1-c/2)/2)\gamma|A_i|cD<(1-\bigO(n^{a-1}))\gamma|A_i|cD,
\end{align*}
where in the first inequality we combine \ref{RNPdegree} and \ref{RNPvertexdistribution} with the assumption and in the last inequality we observe that, since $c>0$, we have $1/2+(1-c/2)/2<1$.
This is a contradiction on \ref{RNPdirections}.}
Let $V_i\coloneqq V(\Sigma(E'_x(A_i),S,\ell^{1/2}))\cap U_i$.
For each $y\in U_i$ we have that
\[\mathbb{P}[y\notin V_i]\leq(1-\varepsilon/D)^{cD/2}=(1\pm\bigO(1/D))e^{-\varepsilon c/2}.\]
Thus, we conclude that
\[\mathbb{E}[|V_i|]\geq(1-\bigO(1/D))(1-e^{-\varepsilon c/2})c\gamma|A_i|/2.\]
Note that the events $\{y\notin V_i\}_{y\in U_i}$ are mutually independent.
Hence, by \cref{lem:Chernoff}\COMMENT{We have
\begin{align*}
    \mathbb{P}[|V_i|\leq(1-e^{-\varepsilon c/2})c\gamma|A_i|/3] &=\mathbb{P}\left[|V_i|\leq\frac23(1-e^{-\varepsilon c/2})c\gamma|A_i|/2\right]\leq\mathbb{P}\left[|V_i|\leq(1\pm\bigO(1/D))\frac23\mathbb{E}[|V_i|]\right]\\
    &\leq e^{-(1\pm\bigO(1/D))\mathbb{E}[|V_i|]/18}\leq e^{-(1\pm\bigO(1/D))(1-e^{-\varepsilon c/2})c\gamma|A_i|/36}\leq e^{-\varepsilon c^2\gamma\beta n/99}.
\end{align*}
In the last inequality we use the fact that $1-e^{-\varepsilon c/2}\geq\varepsilon c/2-\varepsilon^2c^2/8\geq9\varepsilon c/24$.},
\[\mathbb{P}[|V_i|\leq(1-e^{-\varepsilon c/2})c\gamma|A_i|/3]\leq e^{-\varepsilon c^2\gamma\beta n/99}.\]
Finally, note that $(1-e^{-\varepsilon c/2})c\gamma|A_i|/3\geq(\varepsilon c/2-\varepsilon^2c^2/8)c\gamma|A_i|/3\geq\varepsilon c^2\gamma|A_i|/8$.
The statement of \ref{RNCdirectionscovered} follows by a union bound over all $i\in[K]_0$.
\end{proof}

\subsection{Iterating the nibble}

By making use of \cref{lem:nibble}, we can now prove the main result of this section.
Roughly speaking, \cref{thm: nibble} states that, for any constant $\varepsilon>0$ and $\ell\in\mathbb{N}$, with high probability the random graph $\cQ^n_\varepsilon$ contains a set of $\ell$-dimensional cubes which are vertex-disjoint, cover all but a small proportion of the vertices of $\cQ^n_\varepsilon$, and are `sufficiently significant' with respect to every large set of directions, while not being `too significant' with respect to any given direction.

By analogy with the notation introduced before \cref{lem:nibble}, given any $\ell\in\mathbb{N}$, any $S\subseteq\mathcal{D}(\cQ^n)$ and any copy $C$ of $\cQ^\ell$ with $C\subseteq\cQ^n$, we define the \emph{significance of\/ $C$ in $S$} as $\sigma(C,S)\coloneqq|\mathcal{D}(C)\cap S|$\index{sigmaC@$\sigma(C,S)$}.\COMMENT{We will not actually use the name `significance', but we use the notation.}
Similarly, given any set $\mathcal{C}$ of $\ell$-dimensional cubes in $\cQ^n$, we define the \emph{significance of\/ $\mathcal{C}$ in $S$} as $\sigma(\mathcal{C},S)\coloneqq\sum_{C\in\mathcal{C}}\sigma(C,S)$\index{sigmaC2@$\sigma(\mathcal{C},S)$}.
We also denote $\Sigma(\mathcal{C},S,t)\coloneqq\{C\in\mathcal{C}:\sigma(C,S)\geq t\}$\index{sigmaCST@$\Sigma(\mathcal{C},S,t)$}.
Given any $x\in\{0,1\}^n$ and any $Y\subseteq N_{\cQ^n}(x)$, we denote $\mathcal{C}_x(Y)\coloneqq\{C\in\mathcal{C}:\dist(x,C)=1,V(C)\cap Y\neq\varnothing\}$\index{CxY@$\mathcal{C}_x(Y)$}.
In particular, we will write $\mathcal{C}_x\coloneqq\mathcal{C}_x(N_{\cQ^n}(x))$\index{Cx@$\mathcal{C}_x$}.

\begin{theorem}\label{thm: nibble}
Let $\varepsilon,\delta,\alpha,\beta\in(0,1)$ and $K,\ell \in \mathbb{N}$ be such that $1/\ell\ll\alpha \ll \beta$.
For each $x\in\{0,1\}^n$, let $A_0(x)\coloneqq N_{\cQ^n}(x)$ and, for each $i\in[K]$, let $A_i(x)\subseteq A_0(x)$ be a set of size $|A_i(x)|\geq\beta n$.
Then, the graph $\cQ^n_\varepsilon$ a.a.s.~contains a collection $\mathcal{C}$ of vertex-disjoint copies of $\cQ^\ell$ such that the following properties are satisfied for every $x\in\{0,1\}^n$:
\begin{enumerate}[label=$(\mathrm{M}\arabic*)$]
    \item\label{RNFmatchingsize} $|A_0(x) \cap V(\mathcal{C})| \geq (1-\delta)n$;
    \item\label{RNFsingledirection} for every $\hat e\in\mathcal{D}(\cQ^n)$ we have $|\Sigma(\mathcal{C}_x,\{\hat e\},1)|=o(n^{1/2})$;
    \item\label{RNFdirections} for every $i\in[K]_0$ and every $S\subseteq\mathcal{D}(\cQ^n)$ with $\alpha n/2\leq|S|\leq\alpha n$ we have
    \[|\Sigma(\mathcal{C}_x(A_i(x)),S,\ell^{1/2})|\geq|A_i(x)|/3000.\]
\end{enumerate}
\end{theorem}

\begin{proof}
Let $n_0,k\in\mathbb{N}$ and $\varepsilon'>0$ be such that $1/n_0\ll1/k\ll\varepsilon'\ll1/\ell,\delta$ and $1/n_0\ll\varepsilon$, and let $n\geq n_0$.
Let $H\coloneqq H_\ell(\cQ^n_\varepsilon)$. 
Observe that, with the notation from \cref{lem:RNdegree,lem:RNdirections}, for any $x\in\{0,1\}^n$, $y\in A_0(x)$, $S\subseteq\mathcal{D}(\cQ^n)$ and $t\in\mathbb{R}$ we have that $d_H(x)=d^\ell_{\cQ^n_\varepsilon}(x)$ and $d_{H,S,t,x}(y)=d^\ell_{\cQ^n_\varepsilon,S,t,x}(y)$.
Let $D_1\coloneqq\varepsilon^{2^{\ell-1}\ell}n^\ell/\ell!$.

\begin{claim}
We a.a.s.~have that, for every $x\in\{0,1\}^n$ and every $S\subseteq\mathcal{D}(\cQ^n)$ with $|S|\geq\alpha n/2$,
\begin{enumerate}[label=$(\mathrm{C}\arabic*)$]
    \item \label{equa:RNitCond1}  $|\{y\in B^{2\ell}_{\cQ^n}(x):d_H(y)\neq(1\pm\bigO(n^{-1/8}))D_1\}|\leq n^{7/8}$, and
    \item \label{equa:RNitCond2}  $|\{y\in A_0(x):d_{H,S,\ell^{1/2},x}(y)<D_1/2\}|\leq n^{3/4}$.
\end{enumerate}
\end{claim}

\begin{claimproof}
\ref{equa:RNitCond1} holds a.a.s.~by \cref{remark:RNdegree} applied with $a=7/8$ and $2\ell$ playing the role of $r$.
\ref{equa:RNitCond2} holds a.a.s.~by applying \cref{lem:RNdirections} with $a=3/4$ and $\alpha/2$ playing the role of $\delta$.
\end{claimproof}

Now, we condition on \ref{equa:RNitCond1} and \ref{equa:RNitCond2} and will show that there exists a collection $\mathcal{C}$ of vertex-disjoint copies of $\cQ^\ell$ in $\cQ^n_\varepsilon$ satisfying \ref{RNFmatchingsize}--\ref{RNFdirections}, as desired.
In order to do this, we would like to apply \cref{lem:nibble} to $H$ with $a=7/8$ and $D=D_1$.
However, $H$ does not satisfy all the required properties.
It is worth noting that it does satisfy \ref{RNPcodegreefar} and \ref{RNPcodegreeclose}, which follow immediately from \eqref{equa:l-codegbound}.
The argument now will be to modify $H$ slightly so that \cref{lem:nibble} can be applied, independently of the choice of $x\in\{0,1\}^n$, and then iterate.

\begin{claim}\label{claimnibble}
There exists $H_1\subseteq H$ which satisfies \ref{RNPdegree}--\ref{RNPdirections} with $7/8$, $D_1$, $1$, $1/2$ and $\alpha/2$ playing the roles of $a$, $D$, $\gamma$, $c$ and $\delta$, respectively, for every $x\in\{0,1\}^n$ and every value of $k>\ell$.
\end{claim}

\begin{claimproof}[Proof of \cref{claimnibble}]
We construct $H_1$ by removing from $H$ all vertices $y\in\{0,1\}^n$ such that $d_{H}(y)\neq(1\pm\bigO(n^{-1/8}))D_1$.
We first need to show that this deletion does not substantially decrease the degrees of other vertices.
In fact, we claim that, for any $y\in\{0,1\}^n$,
\begin{equation}\label{equa:degH1}
\begin{minipage}[c]{0.75\textwidth}
if $d_{H}(y)=(1\pm\bigO(n^{-1/8}))D_1$, then $d_{H_1}(y)=(1\pm\bigO(n^{-1/8}))D_1$.
\end{minipage}\ignorespacesafterend 
\end{equation} 
Indeed, consider any vertex $y\in\{0,1\}^n$ which satisfies $d_{H}(y)=(1\pm\bigO(n^{-1/8}))D_1$.
By \ref{RNPcodegreefar}, the removal of any vertex $z$ such that $\dist(y,z)>\ell$ does not affect the degree of $y$.
Furthermore, by \ref{equa:RNitCond1}, the number of vertices $z\in B^\ell_{\cQ^n}(y)$ such that $d_{H}(z)\neq(1\pm\bigO(n^{-1/8}))D_1$ is at most $n^{7/8}$.
Let $Z$ be the set of all such vertices.
By \ref{RNPcodegreeclose}, the number of edges incident to $y$ which are removed because of some fixed $z\in Z$ is $\bigO(D_1/n)$.
By adding over all $z\in Z$, we have that the number of edges incident to $y$ which have been removed is $\bigO(D_1n^{-1/8})$, and the claim follows.
One can similarly prove that, for any $y\in V(H_1)$ and any $S\subseteq\mathcal{D}(\cQ^n)$\COMMENT{Note that the edges considered here are subsets of the edges considered above, hence the number of edges deleted can be bounded by the number of edges deleted above.},
\begin{equation}\label{equa:degH1directions}
\begin{minipage}[c]{0.75\textwidth}
if $d_{H,S,\ell^{1/2},x}(y)\geq D_1/2$, then $d_{H_1,S,\ell^{1/2},x}(y)\geq(1-\bigO(n^{-1/8}))D_1/2$.
\end{minipage}\ignorespacesafterend 
\end{equation} 
By \eqref{equa:degH1}, we now have that $H_1$ satisfies \ref{RNPdegree} with the parameters stated in \cref{claimnibble} for every $x\in\{0,1\}^n$.
It also follows that \ref{RNPvertexdistribution} holds for every $x\in\{0,1\}^n$ since, by \ref{equa:RNitCond1}, at most $n^{7/8}$ vertices are removed from the neighbourhood of any vertex in $H$.
Finally, by \eqref{equa:degH1directions} and \ref{equa:RNitCond2}, we also have that $H_1$ satisfies \ref{RNPdirections} for every $x\in\{0,1\}^n$.
\end{claimproof}

Therefore, we are in a position to apply \cref{lem:nibble}.
The argument now will be as follows.
In order to prove that a collection of cubes as described in the statement exists, we will take a random collection of cubes in an iterative manner.
We will prove that such a collection satisfies the desired properties locally with high probability.
Then, we will apply the local lemma to extend the properties to the whole hypergraph.
For this, it is important to define the probability space we are working with.

Fix a vertex $x\in\{0,1\}^n$.
We now proceed iteratively.
Let $G_1$ be the graph obtained by deleting from $\cQ^n_\varepsilon$ all vertices $y\in\{0,1\}^n$ which do not satisfy $d^\ell_{\cQ^n_\varepsilon}(y)=(1\pm\bigO(n^{-1/8}))D_1$.
Note that $H_1=H_\ell(G_1)$.
Let $i\in[k]$ and suppose that we have already defined $G_i$ and $H_i$, where $H_i=H_\ell(G_i)$.
Choose a random set of edges $E_i\subseteq E(H_i)$ by adding each edge in $E(H_i)$ to $E_i$ independently with probability $\varepsilon'/D_i$.
Then, define $H_{i+1}\coloneqq H_i-V(E_i)$.
(Observe that $H_{i+1}=H_\ell(G_{i+1})$, where $G_{i+1}\coloneqq G_i-V(E_i)$.)
Finally, let $D_{i+1}\coloneqq e^{-(2^\ell-1)\varepsilon'}D_i$, and iterate for $k$ steps.

The randomized process above defines a probability space on the sequences of outcomes of each iteration of the process.
Formally, the process, when iterated, results in a random sequence $E^k\coloneqq(E_1,\ldots,E_k)$ of sets of edges of $H_1$.
Note that, for each $i\in[k]$, the hypergraph $H_{i+1}$ is uniquely determined by $(E_1,\ldots,E_i)$, and $H_1$ does not depend on any of these sets; thus, the sequence $E^k$ encodes all the information about the outcome of the iterative process.
For any $i\in[k]$, we will write $E^i\coloneqq(E_1,\ldots,E_i)$. 
We will write $\mathbb{P}_{E^0}[E_1]\coloneqq\mathbb{P}[\text{process outputs }E_1\text{ on input }H_1]$ and, for each $i\in[k]\setminus\{1\}$, we will write $\mathbb{P}_{E^{i-1}}[E_i]\coloneqq\mathbb{P}[\text{process outputs }E_i\text{ on input }H_i]$ (where $H_i$ is determined by $E^{i-1}$ for all $i\geq2$).
Whenever needed, we will treat $E^0$ as an empty sequence.
Note that the choice of the process in any iteration affects the probability distribution on all subsequent iterations.
For each $i\in[k]_0$, let $\Omega^i$ be the set of all sequences $E^i=(E_1,\ldots,E_i)$ such that, for all $j\in[i]$, $\mathbb{P}_{E^{j-1}}[E_j]>0$, and let $\Omega\coloneqq\Omega^k$.
Given any $\boldsymbol{\omega}=E^k=(E_1,\ldots,E_k)\in\Omega$, we write $\boldsymbol{\omega}^i\coloneqq E^i$.
Consider any $\boldsymbol{\omega}=E^k=(E_1,\ldots,E_k)\in\Omega$.
The probability distribution on the outputs of the iterative process is given by $\mathbb{P}_\Omega[\boldsymbol{\omega}]\coloneqq\prod_{j=1}^k\mathbb{P}_{E^{j-1}}[E_j]$.
Similarly, the distribution on the outputs after $i$ iterations of the process is given by $\mathbb{P}_{\Omega^i}[\boldsymbol{\omega}^i]\coloneqq\prod_{j=1}^i\mathbb{P}_{E^{j-1}}[E_j]$.
Observe that, given $i\in[k]$ and $\boldsymbol{\omega}'\in\Omega^i$, we have that\COMMENT{We prove the more general statement that, for any $j\in[k]\setminus[i]$, we have that $\mathbb{P}_{\Omega^j}[\boldsymbol{\omega}^i=\boldsymbol{\omega}']=\mathbb{P}_{\Omega^i}[\boldsymbol{\omega}']$.
Indeed,
\begin{align*}
    \mathbb{P}_{\Omega^j}[\boldsymbol{\omega}^i=\boldsymbol{\omega}']&=\sum_{\substack{(E_{i+1},\ldots,E_j):\\\boldsymbol{\omega}'(E^{i+1},\ldots,E^j)\in\Omega^j}}\mathbb{P}_{\Omega^{j}}[\boldsymbol{\omega}'(E_{i+1},\ldots,E_j)]=\sum_{\substack{(E_{i+1},\ldots,E_j):\\\boldsymbol{\omega}'(E^{i+1},\ldots,E^j)\in\Omega^j}}\prod_{\kappa=1}^j\mathbb{P}_{E^{\kappa-1}}[E_\kappa]\\
    &=\mathbb{P}_{\Omega^i}[\boldsymbol{\omega}']\sum_{\substack{(E_{i+1},\ldots,E_j):\\\boldsymbol{\omega}'(E^{i+1},\ldots,E^j)\in\Omega^j}}\prod_{\kappa=i+1}^j\mathbb{P}_{E^{\kappa-1}}[E_\kappa]=\mathbb{P}_{\Omega^i}[\boldsymbol{\omega}'].
\end{align*}
Note that here we also use that, for each $\boldsymbol{\omega}'\in\Omega^i$, there exists $\boldsymbol{\omega}=(E_1,\ldots,E_k)\in\Omega$ such that $E^i=\boldsymbol{\omega}'$ (this holds since $\boldsymbol{\omega}\coloneqq\boldsymbol{\omega}'(\varnothing,\ldots,\varnothing)$ is always a valid choice, even if $i=1$ and $\boldsymbol{\omega}'=E(H_1)$).}
\begin{equation}\label{equa:RNnewdisteq}
    \mathbb{P}_\Omega[\boldsymbol{\omega}^i=\boldsymbol{\omega}']=\mathbb{P}_{\Omega^i}[\boldsymbol{\omega}'].
\end{equation}

In particular, we wish to apply \cref{lem:nibble} in each iteration of the process.
In order to do so, we will restrict ourselves to a suitable subspace of $\Omega$ by conditioning (again, in an iterative way).
Let $k_1\coloneqq\lfloor1/(3\varepsilon')\rfloor$\COMMENT{This value is chosen so that $c_i>1/6$ for all $i\in[k_1]$.
See the definition of $c_i$ below.}, $\gamma_1\coloneqq1$ and $c_1\coloneqq1/2$.
For each $i\in[k_1]$, we proceed as follows.
Given $E^{i-1}\in\Omega^{i-1}$ (and thus the hypergraph $H_i$), let $\mathcal{A}_i(E^{i-1})$ be the event that $E_i$, $E_{i}'\coloneqq\{e\in E_i:e\cap V(E_i\setminus\{e\})=\varnothing\}$, $V_i\coloneqq V(H_i)\setminus V(E_i)$ and $H_{i+1}=H_i[V_i]$ satisfy \ref{RNCremainingvertices}--\ref{RNCdirectionsleft} with $D_i$, $(k-i+1)\ell+1$, $\gamma_i$, $c_i$, $\varepsilon'$ and $\alpha/2$ playing the roles of $D$, $k$, $\gamma$, $c$, $\varepsilon$ and $\delta$, respectively (in all iterations we will use $a=7/8$). 
Then, let $\gamma_{i+1}\coloneqq e^{-\varepsilon'}\gamma_i$ and $c_{i+1}\coloneqq c_i-\varepsilon'$, and iterate.
Note that the definition of $k_1$ guarantees that, for all $i\in[k_1]$, we have $\gamma_i\geq2/3$\COMMENT{$\gamma_i=e^{-(i-1)\varepsilon'}\geq e^{-k_1\varepsilon'}= e^{-\lfloor1/(3\varepsilon')\rfloor\varepsilon'}\geq e^{-1/3}\geq2/3$} and $c_i\geq1/6$.

\begin{claim}\label{RNclaim3}
For any $i\in[k_1]$, let $E^{i-1}\in\Omega^{i-1}$ be such that $H_i$ satisfies \ref{RNPdegree}--\ref{RNPdirections} with $D_i$, $(k-i+1)\ell+1$, $\gamma_i$, $c_i$ and $\alpha/2$ playing the roles of $D$, $k$, $\gamma$, $c$ and $\delta$, respectively.
Then,
\[\mathbb{P}_{E^{i-1}}[\mathcal{A}_i(E^{i-1})]\geq1-e^{-\Theta(n^{1/2})}.\]
\end{claim}

\begin{claimproof}
This follows immediately from \cref{lem:nibble}\COMMENT{This follows by a union bound on each of the properties \ref{RNCremainingvertices}--\ref{RNCdirectionsleft} given by the nibble.}.
\end{claimproof}

This will naturally lead us into applying \cref{lem:nibble} iteratively.
Indeed, in any given iteration, assume that $H_i$ satisfies \ref{RNPdegree}--\ref{RNPdirections} with $D_i$, $(k-i+1)\ell+1$, $\gamma_i$, $c_i$ and $\alpha/2$ playing the roles of $D$, $k$, $\gamma$, $c$ and $\delta$, respectively.
Note that, for $i=1$, by \cref{claimnibble}, these properties hold for every choice of $x\in\{0,1\}^n$ (but recall that we have now fixed $x$).
If $\mathcal{A}_i(E^{i-1})$ holds, then, because of \ref{RNCremainingvertices}, \ref{RNCdegree} and \ref{RNCdirectionsleft}, the next hypergraph $H_{i+1}$ satisfies \ref{RNPdegree}--\ref{RNPdirections} with $D_{i+1}$, $(k-i)\ell+1$, $\gamma_{i+1}$, $c_{i+1}$ and $\alpha/2$ playing the roles of $D$, $k$, $\gamma$, $c$ and $\delta$, respectively, so \cref{lem:nibble} can be applied again\COMMENT{The idea here is that, with these new constants, $H_{i+1}$ satisfies the conditions for \cref{lem:nibble} and we can iterate.
Indeed, recalling that $D_{i+1}=e^{-(2^\ell-1)\varepsilon'}$, the definition for $\gamma_{i+1}$ follows from \ref{RNCremainingvertices}, and the definition of $c_{i+1}$ follows from \ref{RNCdirectionsleft}.
To see this last one, in particular, observe that setting $c_{i+1}\coloneqq(c_i-\varepsilon')e^{\varepsilon'}$ would be enough to satisfy the condition. 
The definition we chose is simpler and sufficient for the calculations.}.

As discussed above, in order to apply \cref{lem:nibble} fully in each iteration, we must condition on the event that certain properties are satisfied after the previous iteration (namely, the corresponding event $\mathcal{A}_i$ holds).
For each $j\in[k_1]_0$, let $\Omega_*^j\coloneqq\{(E_1,\ldots,E_j)\in\Omega^j:E_i\in\mathcal{A}_i(E^{i-1})\text{ for all }i\in[j]\}$.
We denote $\Omega_*\coloneqq\Omega_*^{k_1}$. 
Using \cref{RNclaim3}, it now easily follows by induction that, for any $i\in[k_1]$,\COMMENT{The base case $i=1$ follows directly from \cref{RNclaim3}, given that $H_1$ satisfies \ref{RNPdegree}--\ref{RNPdirections} by \cref{claimnibble}.
Now consider the general inductive step.
By \eqref{equa:RNnewdisteq} and \cref{RNclaim3}, we have that
\begin{align*}
    \mathbb{P}_{\Omega}\left[\boldsymbol{\omega}^i\in\Omega^i_*\right]&=\mathbb{P}_{\Omega^i}\left[\Omega^i_*\right]=\sum_{E_1\in\mathcal{A}_1(E^0)}\ldots\sum_{E_i\in\mathcal{A}_i(E^{i-1})}\mathbb{P}_{\Omega^{i}}[E^{i}]=\sum_{E_1\in\mathcal{A}_1(E^0)}\ldots\sum_{E_i\in\mathcal{A}_i(E^{i-1})}\prod_{j=1}^i\mathbb{P}_{E^{j-1}}[E_{j}]\\
    &=\sum_{E_1\in\mathcal{A}_1(E^0)}\ldots\sum_{E_{i-1}\in\mathcal{A}_{i-1}(E^{i-2})}\left(\prod_{j=1}^{i-1}\mathbb{P}_{E^{j-1}}[E_{j}]\left(\sum_{E_{i}\in\mathcal{A}_{i}(E^{i-1})}\mathbb{P}_{E^{i-1}}[E_{i}]\right)\right)\\
    &=\sum_{E_1\in\mathcal{A}_1(E^0)}\ldots\sum_{E_{i-1}\in\mathcal{A}_{i-1}(E^{i-2})}\left(\prod_{j=1}^{i-1}\mathbb{P}_{E^{j-1}}[E_{j}]\right)\mathbb{P}_{E^{i-1}}[\mathcal{A}_i(E^{i-1})]\\
    &=\mathbb{P}_{E^{i-1}}[\mathcal{A}_i(E^{i-1})]\sum_{E_1\in\mathcal{A}_1(E^0)}\ldots\sum_{E_{i-1}\in\mathcal{A}_{i-1}(E^{i-2})}\mathbb{P}_{\Omega^{i-1}}[E^{i-1}]\\
    &=\mathbb{P}_{E^{i-1}}[\mathcal{A}_i(E^{i-1})]\,\mathbb{P}_{\Omega^{i-1}}\left[\Omega^{i-1}_*\right]=\mathbb{P}_{E^{i-1}}[\mathcal{A}_i(E^{i-1})]\,\mathbb{P}_{\Omega}\left[\boldsymbol{\omega}^i\in\Omega^{i-1}_*\right]\\
    &\geq(1-e^{-\Theta(n^{1/2})})(1-e^{-\Theta(n^{1/2})})=1-e^{-\Theta(n^{1/2})}.
\end{align*}}
\begin{equation}\label{equa:RNconditioningspaces}
    \mathbb{P}_{\Omega}[\boldsymbol{\omega}^i\in\Omega^i_*]\geq1-e^{-\Theta(n^{1/2})}.
\end{equation}

Now fix a set of directions $S\subseteq\mathcal{D}(\cQ^n)$ with $\alpha n/2\leq|S|\leq\alpha n$.
For each $i\in[k_1]$, let $\mathcal{B}_i(S)\subseteq\Omega^i$ be the event that $|V(\Sigma({E_{i}}_x(A_j(x)),S,\ell^{1/2}))\cap A_j(x)|<\varepsilon'c_i^2\gamma_i|A_j(x)|/8$ for some $j\in[K]_0$.
In other words, $\mathcal{B}_i(S)$ is the event that, in the $i$-th iteration, \ref{RNCdirectionscovered} fails for $S$.
Let 
\begin{align*}
    a_i&\coloneqq\max_{E^{i-1}\in\Omega^{i-1}_*}\mathbb{P}_{\Omega^i}[\boldsymbol{\omega}\in\mathcal{B}_i(S)\mid\boldsymbol{\omega}^{i-1}=E^{i-1}]=\max_{E^{i-1}\in\Omega^{i-1}_*}\mathbb{P}_{E^{i-1}}[\mathcal{B}_i(S)],\\
    b_i&\coloneqq\max_{E^{i-1}\in\Omega^{i-1}_*}\mathbb{P}_{\Omega^i}[\boldsymbol{\omega}\in\mathcal{B}_i(S)\cap\mathcal{A}_i(E^{i-1})\mid\boldsymbol{\omega}^{i-1}=E^{i-1}]=\max_{E^{i-1}\in\Omega^{i-1}_*}\mathbb{P}_{E^{i-1}}[\mathcal{B}_i(S)\cap\mathcal{A}_i(E^{i-1})],\\
    f_i&\coloneqq\max_{E^{i-1}\in\Omega^{i-1}_*}\mathbb{P}_{\Omega^i}[\mathcal{A}_i(E^{i-1})\mid\boldsymbol{\omega}^{i-1}=E^{i-1}]=\max_{E^{i-1}\in\Omega^{i-1}_*}\mathbb{P}_{E^{i-1}}[\mathcal{A}_i(E^{i-1})].
\end{align*}
By \cref{lem:nibble}\ref{RNCdirectionscovered}, for each $i\in[k_1]$ we have that\COMMENT{By the choice of $k_1$, we have that $\gamma_i=e^{-(i-1)\varepsilon'}\geq e^{-k_1\varepsilon'}= e^{-\lfloor1/(3\varepsilon')\rfloor\varepsilon'}\geq e^{-1/3}\geq2/3$ and $c_i\geq1/6$ for all $i\in[k_1]$.
Then, by \cref{lem:nibble}\ref{RNCdirectionscovered}, we have that
\[a_i\leq e^{-\varepsilon'c_i^2\gamma_i\beta n/100}\leq e^{-\varepsilon'\beta n/5400}.\]}
\[b_i\leq a_i\leq e^{-\varepsilon'\beta n/5400}\eqqcolon d.\]
Let $\mathcal{I}(S)$ be the set of indices $i\in[k_1]$ in which $\mathcal{B}_i(S)$ holds. 
Note that, for any set $\mathcal{I}\subseteq[k_1]$, by \eqref{equa:RNnewdisteq} we have that $\mathbb{P}_\Omega[\mathcal{I}(S)=\mathcal{I}]=\mathbb{P}_{\Omega^{k_1}}[\mathcal{I}(S)=\mathcal{I}]$.
Using the definitions above and induction on $k_1$\COMMENT{The proof is similar to that of \eqref{equa:RNconditioningspaces}.}, it follows that 
\[f\coloneqq\mathbb{P}_{\Omega^{k_1}}[(\mathcal{I}(S)=\mathcal{I})\wedge\Omega_*]\leq\prod_{i\in\mathcal{I}}b_i\prod_{i\in[k_1]\setminus\mathcal{I}}f_i\leq d^{|\mathcal{I}|}.\]
Let $X=X(S)\coloneqq|\mathcal{I}(S)|$.
By adding over all sets $\mathcal{I}\subseteq[k_1]$ with $|\mathcal{I}|\geq k_1/2$, we conclude that\COMMENT{For the second inequality, we use the definition of $k_1$ and the fact that $1/n\ll\varepsilon'\ll1$ to show that $2^{k_1}d^{k_1/2}\leq2^{1/(3\varepsilon')}e^{-\varepsilon'\beta k_1 n/10800}\leq2^{1/(3\varepsilon')}e^{-\beta n/40000}\leq e^{-\beta n/50000}$.}
\begin{equation}\label{equa:nibthmnewfeb}
    \mathbb{P}_{\Omega^{k_1}}[(X\geq k_1/2)\wedge\Omega_*]\leq 2^{k_1}d^{k_1/2}\leq e^{-\beta n/50000}.
\end{equation}

Let $\mathcal{B}\subseteq\Omega^{k_1}$ be the event that there exists a set $S\subseteq\mathcal{D}(\cQ^n)$ with $\alpha n/2\leq|S|\leq\alpha n$ such that the event $\mathcal{B}_i(S)$ holds in at least $k_1/2$ iterations.
A union bound on \eqref{equa:nibthmnewfeb} over all choices of $S\subseteq\mathcal{D}(\cQ^n)$ with $\alpha n/2\leq|S|\leq\alpha n$ allows us to conclude that $\mathbb{P}_{\Omega^{k_1}}[\mathcal{B}\wedge\Omega_*]\leq e^{-\Theta(n)}$\COMMENT{Here is where we use the choice of $\alpha$.
Indeed, we want that 
\[\sum_{i=\alpha n/2}^{\alpha n}\binom{n}{i}e^{-\beta n/50000}\leq e^{-\Theta(n)}.\]
We can bound $\sum_{i=\alpha n/2}^{\alpha n}\binom{n}{i}\leq\alpha n\binom{n}{\alpha n}\leq\alpha n\left(\frac{en}{\alpha n}\right)^{\alpha n}\leq e^{\sqrt{\alpha}n}$.
The statement follows by the choice of $\alpha\ll\beta$.}.
Finally, combining this with \eqref{equa:RNnewdisteq} and \eqref{equa:RNconditioningspaces}, we have that\COMMENT{We have that
\begin{align*}
    \mathbb{P}_\Omega[\boldsymbol{\omega}^{k_1}\in\mathcal{B}]&=\mathbb{P}_\Omega[\boldsymbol{\omega}^{k_1}\in\mathcal{B}\mid\boldsymbol{\omega}^{k_1}\in\Omega_*]\,\mathbb{P}_\Omega[\boldsymbol{\omega}^{k_1}\in\Omega_*]+\mathbb{P}_\Omega[\boldsymbol{\omega}^{k_1}\in\mathcal{B}\mid\boldsymbol{\omega}^{k_1}\in\overline{\Omega_*}]\,\mathbb{P}_\Omega[\boldsymbol{\omega}^{k_1}\in\overline{\Omega_*}]\\
    &\leq\mathbb{P}_\Omega[\boldsymbol{\omega}^{k_1}\in\mathcal{B}\mid\boldsymbol{\omega}^{k_1}\in\Omega_*]+\mathbb{P}_\Omega[\boldsymbol{\omega}^{k_1}\in\overline{\Omega_*}]=\frac{\mathbb{P}_\Omega[\boldsymbol{\omega}^{k_1}\in\mathcal{B}\cap\Omega_*]}{\mathbb{P}_\Omega[\boldsymbol{\omega}^{k_1}\in\Omega_*]}+\mathbb{P}_\Omega[\boldsymbol{\omega}^{k_1}\in\overline{\Omega_*}]\\
    &=\frac{\mathbb{P}_{\Omega^{k_1}}[\mathcal{B}\wedge\Omega_*]}{\mathbb{P}_\Omega[\boldsymbol{\omega}^{k_1}\in\Omega_*]}+\mathbb{P}_\Omega[\boldsymbol{\omega}^{k_1}\in\overline{\Omega_*}]\leq e^{-\Theta(n)}+ e^{-\Theta(n^{1/2})}.
\end{align*}
}
\begin{equation}\label{equa:RNBadprob}
    \mathbb{P}_\Omega[\boldsymbol{\omega}^{k_1}\in\mathcal{B}]\leq\mathbb{P}_\Omega[\boldsymbol{\omega}^{k_1}\in\mathcal{B}\mid\boldsymbol{\omega}^{k_1}\in\Omega_*]+\mathbb{P}_\Omega[\boldsymbol{\omega}^{k_1}\in\overline{\Omega_*}]\leq e^{-\Theta(n^{1/2})}.
\end{equation}

Now, for all $i\in[k]\setminus[k_1]$, we iterate as above with the difference that we no longer require \ref{RNPdirections}.
Thus, we can no longer guarantee that \ref{RNCdirectionsleft} or \ref{RNCdirectionscovered} hold with high probability, but we still have \ref{RNCremainingvertices}--\ref{RNCdegree} as above.
To be more precise, given any $i\in[k]\setminus[k_1]$ and $E^{i-1}\in\Omega^{i-1}$ (and thus the hypergraph $H_i$), let $\mathcal{A}_i(E^{i-1})$ be the event that $E_i$, $E_i'\coloneqq\{e\in E_i:e\cap V(E_i\setminus\{e\})=\varnothing\}$, $V_i\coloneqq V(H_i)\setminus V(E_i)$ and $H_{i+1}=H_i[V_i]$ satisfy \ref{RNCremainingvertices}--\ref{RNCdegree} with $D_i$, $(k-i+1)\ell+1$, $\varepsilon'$ and $\gamma_i$ playing the roles of $D$, $k$, $\varepsilon$ and $\gamma$, respectively.
Then, let $\gamma_{i+1}\coloneqq e^{-\varepsilon'}\gamma_i$ and iterate.
Similarly to \cref{RNclaim3}, we can now show the following.

\begin{claim}\label{RNclaim4}
For any $i\in[k]\setminus[k_1]$, let $E^{i-1}$ be such that $H_i$ satisfies \ref{RNPdegree} and \ref{RNPvertexdistribution} with $D_i$, $(k-i+1)\ell+1$ and $\gamma_i$ playing the roles of $D$, $k$ and $\gamma$, respectively.
Then,
\[\mathbb{P}_{E^{i-1}}[\mathcal{A}_i(E^{i-1})]\geq1-e^{-\Theta(n^{1/2})}.\]
\end{claim}

\begin{claimproof}
This follows immediately from \cref{lem:nibble}\COMMENT{This follows by a union bound on each of the properties \ref{RNCremainingvertices}--\ref{RNCdegree} given by the nibble.}.
\end{claimproof}

Now, assume that $H_i$ satisfies \ref{RNPdegree} and \ref{RNPvertexdistribution} with $D_i$, $(k-i+1)\ell+1$ and $\gamma_i$ playing the roles of $D$, $k$ and $\gamma$, respectively.
Note that, for $i=k_1+1$, we have these properties by conditioning on $\mathcal{A}_{k_1}(E^{k_1-1})$.
Then, if $\mathcal{A}_i(E^{i-1})$ holds, because of \ref{RNCremainingvertices} and \ref{RNCdegree}, the hypergraph $H_{i+1}$ satisfies \ref{RNPdegree} and \ref{RNPvertexdistribution} with $D_{i+1}$, $(k-i)\ell+1$ and $\gamma_{i+1}$ playing the roles of $D$, $k$ and $\gamma$, respectively, so we may apply \cref{lem:nibble} again.

Let $\mathcal{A}\coloneqq\{(E_1,\ldots,E_k)\in\Omega:E^{k_1}\in\Omega_*\cap\overline{\mathcal{B}},E_i\in\mathcal{A}_i(E^{i-1})\text{ for all }i\in[k]\setminus[k_1]\}$.
By combining \eqref{equa:RNconditioningspaces}, \eqref{equa:RNBadprob} and \cref{RNclaim4}, observe that $\mathbb{P}_\Omega[\mathcal{A}]\geq1-e^{-\Theta(n^{1/2})}$\COMMENT{This can be proved as follows.
For each $j\in[k]\setminus[k_1]$, let $\Omega^j_*\coloneqq\{(E_1,\ldots,E_j)\in\Omega^j: E_i\in\mathcal{A}_i(E^{i-1})\text{ for all }i\in[j]\}$.
Then, observe that
\begin{align*}
    \mathbb{P}_\Omega[\mathcal{A}]&\leq\mathbb{P}_\Omega[\boldsymbol{\omega}^{k_1}\in\overline{\mathcal{B}}\cap\Omega_*]\prod_{i=k_1+1}^k\max_{E^{i-1}\in\Omega^{i-1}_*}\mathbb{P}_{\Omega}[\boldsymbol{\omega}^i\in\mathcal{A}_i(E^{i-1})\mid\boldsymbol{\omega}^{i-1}=E^{i-1}]\\
    &\leq (1-e^{-\Theta(n^{1/2})}-e^{-\Theta(n^{1/2})})\prod_{i=k_1+1}^k(1-e^{-\Theta(n^{1/2})})=1-e^{-\Theta(n^{1/2})},
\end{align*}
where in the second inequality we use \eqref{equa:RNconditioningspaces}, \eqref{equa:RNBadprob} and \cref{RNclaim4}.}.
For any $(E_1,\ldots,E_k)\in\mathcal{A}$, let $E\coloneqq\bigcup_{i=1}^kE_{i}$ and $E'\coloneqq\bigcup_{i=1}^kE_i'$.
Note that $E'$ is a matching by construction, that is, it corresponds to a collection $\cC'$ of vertex-disjoint copies of $\cQ^\ell$ in $\cQ^n_\varepsilon$.
We will now show that $\cC'$ satisfies \ref{RNFmatchingsize}--\ref{RNFdirections} for our fixed vertex $x$\COMMENT{The first property follows by adding over all $i\in[k]$ the conclusion given by \ref{RNCmatchingsize} in each iteration.
Here is where the choice of $k$ is needed.
Indeed, for any $j\in[K]_0$ we have that 
\begin{align*}
    |A_j(x)\cap V(E')|&=\sum_{i=1}^k|A_j(x)\cap V(E_i')|\geq\sum_{i=1}^k\varepsilon'(1-2^{\ell+1}\varepsilon')\gamma_i|A_j(x)|=\varepsilon'(1-2^{\ell+1}\varepsilon')|A_j(x)|\sum_{i=1}^ke^{-(i-1)\varepsilon'}\\
    &=\varepsilon'(1-2^{\ell+1}\varepsilon')\frac{1-e^{-k\varepsilon'}}{1-e^{-\varepsilon'}}|A_j(x)|.
\end{align*}
We want this to be at least $(1-\delta)|A_j(x)|$, so it suffices to have 
\[1-e^{-k\varepsilon'}\geq\frac{(1-\delta)(1-e^{-\varepsilon'})}{\varepsilon'(1-2^{\ell+1}\varepsilon')}.\]
Since $1-e^{-\varepsilon'}\leq\varepsilon'$, having
\[1-e^{-k\varepsilon'}\geq\frac{1-\delta}{1-2^{\ell+1}\varepsilon'}.\]
guarantees that what we want holds.
And since $1/k\ll\varepsilon'\ll1/\ell,\delta$, this holds.\\
The second property follows trivially by adding the conclusion of \ref{RNCsingledirection} over all $i\in[k]$.}.
Indeed, \ref{RNFmatchingsize} and \ref{RNFsingledirection} hold for $x$ since $(E_1,\ldots,E_k)\in\mathcal{A}$ implies that \ref{RNCmatchingsize} and \ref{RNCsingledirection} hold in each iteration (note that \ref{RNFmatchingsize} follows from the case $i=0$ of \ref{RNCmatchingsize}).
In order to prove \ref{RNFdirections} for $x$, consider the following.
For each $S\subseteq\mathcal{D}(\cQ^n)$ with $\alpha n/2\leq|S|\leq\alpha n$, there are at least $k_1/2$ iterations $i\in[k_1]$ in which \ref{RNCdirectionscovered} holds (for all $j\in[K]_0$).
For each such set $S$, let $\mathcal{I}^*(S)\subseteq[k_1]$ be the set of indices of all such iterations.
In particular, for each $i\in\mathcal{I}^*(S)$ we have that $|V(\Sigma({E_i}_x(A_j(x)),S,\ell^{1/2}))\cap A_j(x)|\geq\varepsilon'c_i^2\gamma_i|A_j(x)|/8\geq\varepsilon'|A_j(x)|/432$\COMMENT{For all $i\in[k_1]$ we have that $\gamma_i\geq2/3$ and $c_i\geq1/6$, as was already discussed.} for all $j\in[K]_0$.
Furthermore, by \ref{RNCremainingvertices} and \ref{RNCmatchingsize}, we know that the number of vertices of $A_j(x)$ covered by $E_i\setminus E_i'$ satisfies $|V(E_i\setminus E_i')\cap A_j(x)|\leq2^{\ell+2}(\varepsilon')^2|A_j(x)|$\COMMENT{Note that \ref{RNCremainingvertices} implies that $|V(H_{i-1})\cap A_j(x)|=(1\pm\bigO(n^{-1/8}))\gamma_i|A_j(x)|$ and, thus, $|V(E_i)\cap A_j(x)|=(1\pm\bigO(n^{-1/8}))(1-e^{-\varepsilon'})\gamma_i|A_j(x)|$.
Thus, we have that
\begin{align*}
    |V(E_i\setminus E_i')\cap A_j(x)|&=|V(E_i)\cap A_j(x)|-|V(E_i')\cap A_j(x)|\\
    &\leq(1+\bigO(n^{-1/8}))(1-e^{-\varepsilon'})\gamma_i|A_j(x)|-\varepsilon'(1-2^{\ell+1}\varepsilon')\gamma_i|A_j(x)|\\
    &=(1+\bigO(n^{-1/8}))(1-e^{-\varepsilon'}-\varepsilon'+2^{\ell+1}(\varepsilon')^2)\gamma_i|A_j(x)|\\
    &\leq2^{\ell+2}(\varepsilon')^2\gamma_i|A_j(x)|\leq2^{\ell+2}(\varepsilon')^2|A_j(x)|.
\end{align*}
}, so $|V(\Sigma({E_i'}_x(A_j(x)),S,\ell^{1/2}))\cap A_j(x)|\geq(1/432-2^{\ell+2}\varepsilon')\varepsilon'|A_j(x)|$ for all $j\in[K]_0$.
By adding this over all $i\in\mathcal{I}^*(S)$, we conclude that $|V(\Sigma(E'_x(A_j(x)),S,\ell^{1/2}))\cap A_j(x)|\geq|A_j(x)|/3000$\COMMENT{We have that
\begin{align*}
    |V(\Sigma(E'_x(A_j(x)),S,\ell^{1/2}))\cap A_j(x)|&=\sum_{i\in\mathcal{I}(S)}|V(\Sigma({E_i'}_x(A_j(x)),S,\ell^{1/2}))\cap A_j(x)|\\
    &\geq(1/432-2^{\ell+2}\varepsilon')\varepsilon'|A_j(x)|k_1/2\\
    &\geq(1/432-2^{\ell+2}\varepsilon')2|A_j(x)|/13\geq|A_j(x)|/3000.
\end{align*}
} for all $j\in[K]_0$, as we wanted to see.

In order to show that there exists a matching which satisfies \ref{RNFmatchingsize}--\ref{RNFdirections} simultaneously for all $x\in\{0,1\}^n$, let $\mathcal{E}(x)$ be the event that they hold for $x$.
By the above discussion, we have that $\mathbb{P}_\Omega[\overline{\mathcal{E}(x)}]\leq \mathbb{P}_\Omega[\overline{\cA}]\leq e^{-\Theta(n^{1/2})}$ for each $x\in\{0,1\}^n$.
Furthermore, throughout the iterative process, the presence or absence in $E$ of any edges $e$ such that $\dist(x,e)>k\ell+1$ does not have any effect on $\mathcal{E}(x)$, so $\mathcal{E}(x)$ is mutually independent of all events $\mathcal{E}(y)$ with $\dist(x,y)\geq3k\ell$\COMMENT{Need to have at least $2(k\ell+2)+\ell-2$ ($k\ell+2$ for the distance to $x$ and $y$ and the extra $\ell-2$ for the diameter of the edge). What we write is just simpler and also works.}.
Thus, by \cref{lem: LLL}, we conclude that there is a choice of $E$ which satisfies \ref{RNFmatchingsize}--\ref{RNFdirections} for every $x\in\{0,1\}^n$\COMMENT{$n^{3k\ell}$ is a bound on the degree of the dependency graph (follows from \cref{rem:basicequationIneedhere}, say), and the bound needed for \cref{lem: LLL} holds trivially for $n$ sufficiently large.}.
\end{proof}


\section{Near-spanning trees in random subgraphs of the hypercube}\label{section:tree}

In this section we present our results on  bounded degree near-spanning trees in $\cQ^n_\varepsilon$.
In \cref{section:tree1} we prove the main result of this section (\cref{lem: main treereshit}).
This implies that with high probability there exists a near-spanning bounded degree tree in $\cQ^n_\eps$, which covers most of the neighbourhood of every vertex whilst avoiding a small random set of vertices, to which we refer as a reservoir.
In \cref{section:tree2} we prove \cref{thm: maintreeres}, which allows us to extend the tree using vertices of the reservoir such that (amongst others) the proportion of uncovered vertices is even smaller.\COMMENT{We will apply the following theorem with $\eps' = 1/4$.
It follows that each vertex not in the tree will have $3/4$ of its neighbours in the tree.
By revealing more probability and using a matching Hall-type argument we can absorb almost all vertices into the tree except for an exceptional set which is exponentially small in size.
Furthermore, these exceptional vertices do not clump together and therefore the neighbourhood property will be satisfied for any constant $\eps'$.}
Finally, in \cref{lem:repatch} we show that, if some number of small local obstructions is prescribed, the tree given by \cref{thm: maintreeres} can be slightly modified to avoid these obstructions.
For convenience, throughout this section, we move away from the algebraic notation for the hypercube to a more combinatorial notation.

We (re-)define the hypercube by setting $V(\cQ^n)\coloneqq\mathcal{P}([n])$\index{Qn@$\cQ^n$} and joining two vertices $u,v\in\mathcal{P}([n])$ by an edge if and only if $||u|-|v||=1$ and $u\subseteq v$ or $v\subseteq u$.
In this setting, directions correspond to the elements in $[n]$, and following a direction $i\in[n]$ from a vertex $v\in\mathcal{P}([n])$ means adding $i$ to $v$ if $i\notin v$, or deleting it from $v$ if $i\in v$.
Note that there is a natural partition of $V(\cQ^n)$ into sets such that every vertex of a set has the same size.
Given any set $S\subseteq [n]$, we denote $S^{(t)}\coloneqq \{X\subseteq S: |X|=t\}$\index{St@$S^{(t)}$}. 
We will denote by $L_i$\index{Lilevel@$L_i$ (level)}, for $i\in[n]_0$, the set of all vertices $v\in V(\cQ^n) = \mathcal{P}([n])$ with $|v|=i$ (that is, $L_i=[n]^{(i)}$), and we will refer to these sets as \emph{levels}.
This notation is especially useful because of the natural notion of containment of vertices, which provides a partial order on the vertices of $\cQ^n$.
Given any graph $G\subseteq\cQ^n$, for a vertex $x \in L_i$, we refer to the neighbours of $x$ in $G$ lying in $L_{i+1}$ as \emph{up-neighbours}, and to the neighbours of $x$ in $L_{i-1}$ as \emph{down-neighbours}, and denote these sets by $N_G^{\uparrow}(x)$\index{NGup@$N_G^{\uparrow}(x)$, $N_G^{\downarrow}(x)$} and $N_G^{\downarrow}(x)$, respectively.
We write $d_G^{\uparrow}(x)\coloneqq|N_G^{\uparrow}(x)|$\index{degGup@$d_G^{\uparrow}(x)$, $d_G^{\downarrow}(x)$} and $d_G^{\downarrow}(x)\coloneqq|N_G^{\downarrow}(x)|$.
Whenever the subscript is omitted, we mean that $G=\cQ^n$.
We will say that a path $P=v_1\ldots v_k$ in $\cQ^n$ is a \emph{chain} if its vertices satisfy the relation $v_1\subseteq\ldots\subseteq v_k$, and refer to it as a $v_1$-$v_k$ chain.

In more generality, because of the symmetries of the hypercube, this notation can be extended with respect to any vertex $v\in V(\cQ^n)$ by defining, for each $i\in[n]_0$, $L_i(v)\coloneqq\{u\in V(\cQ^n):\dist(u,v)=i\}$\index{Lilevel2@$L_i(v)$}.
One can then define up-neighbours and down-neighbours with respect to $v$, and use the notations $N_G^{\uparrow}(x,v)$\index{NGup2@$N_G^{\uparrow}(x,v)$, $N_G^{\downarrow}(x,v)$}, $N_G^{\downarrow}(x,v)$, $d_G^{\uparrow}(x,v)$\index{degGup2@$d_G^{\uparrow}(x,v)$, $d_G^{\downarrow}(x,v)$} and $d_G^{\downarrow}(x,v)$, for $G\subseteq\cQ^n$.
We say that a path $P=v_1\ldots v_k$ in $\cQ^n$ is a \emph{chain with respect to $v$} if its vertices satisfy that, if $v_1\in L_j(v)$ for some $j\in[n]_0$, then for all $\ell\in[k]\setminus\{1\}$ we have $v_\ell\in L_{j+\ell-1}(v)$, and refer to it as a $v_1$-$v_k$ chain\COMMENT{A chain is simply a shortest path between two vertices, lying in different levels (with respect to $v$).}.
Given any graph $G\subseteq\cQ^n$, for any $i\in[n]$ and $v\in V(\mathcal{Q}^n)$, we will write $E_G(L_{i-1}(v),L_{i}(v))$ for the set of edges of $G$ whose endpoints lie in the levels $L_{i-1}(v)$ and $L_{i}(v)$, respectively.
We will drop the subscript whenever $G=\cQ^n$.

\subsection{Constructing a bounded degree near-spanning tree}\label{section:tree1}

Our goal in this subsection is to prove \cref{lem: main treereshit} below.
Given a graph $G$ and $\delta\in[0,1]$, let $\mathit{Res}(G,\delta)$\index{Res@$\mathit{Res}(G,\delta)$} be a probability distribution on subsets of $V(G)$, where $R \sim \mathit{Res}(G,\delta)$ is obtained by adding each vertex $v\in V(G)$ to $R$ with probability $\delta$, independently of every other vertex.
We will refer to this set $R$ as a \emph{reservoir}.

\begin{theorem}\label{lem: main treereshit}
Let $0 < 1/D, \delta \ll \eps' \le 1/2$, and let $\eps, \gamma\in(0,1]$ and $k\in \mathbb{N}$.
Then, the following holds a.a.s.
Let $\cA\subseteq V(\cQ^{n})$ with the following two properties:
\begin{enumerate}[label=$(\mathrm{P}\arabic*)$]
    \item for any distinct $x,y \in \cA$ we have $\dist(x,y)\geq \gamma n$, and
    \item $B^{k+2}_{\cQ^n}(\cA) \cap \{ \varnothing, [n],[\lceil n/2 \rceil], [n]\setminus [\lceil n/2 \rceil]\}=\varnothing$.\COMMENT{Note, it follows that you can `take out' at most a $1/n$ proportion of any level, which is fine.
    The $+2$ will be to ensure we don't take out any of the levels $L_1$ or $L_2$  for when we need to show there's a Hamilton cycle between both.}
\end{enumerate} 
Let $R \sim \mathit{Res}(\cQ^n, \delta)$. 
Then, there exists a tree $T \subseteq \cQ^n_\eps-(R \cup B^{k}_{\cQ^n}(\cA))$ such that
\begin{enumerate}[label=$(\mathrm{T}\arabic*)$]
\item $\Delta(T) < D$, 
\item\label{prop22hit} for all $x \in V(\cQ^n) \setminus B^{k}_{\cQ^n}(\cA)$, we have that $|N_{\cQ^n}(x) \cap V(T)| \ge (1-\varepsilon') n$.
\end{enumerate}
\end{theorem}

The set $\cA$ will be important in the proof of \cref{thm:hitting,thm: kedgehit}, where it will play the role of the set $\cU$ of vertices of small degree.
In the proof of \cref{thm:thresholdk,thm:main} we can take $\cA = \varnothing$.

To prove \cref{lem: main treereshit}, we will consider suitable `branching-like' processes which start at the `bottom' of the hypercube, and grow `upwards'.
The tree will be formed by considering unions of such processes.
The precise definition of the model we use is given in \cref{def: perc}. 
Crucially, there is a joint distribution of this branching-like process model and the binomial model  $\cQ^n_\eps$.\COMMENT{This guarantees these graphs can be found as subgraphs of $\cQ^n_\eps$ for any $\eps>0$.}
These processes are analysed and constructed in the results leading up to \cref{lem: chain count}.
Subgraphs obtained from the processes are then connected into a tree in \cref{lem: treeres2hit}.

We begin with a formal description of our model.
We denote by $\mathbf{p}= (p_0, \dots, p_{n-1})\in[0,1]^n$\index{p@$\mathbf{p}$} an $n$-component vector of probabilities.
We now describe a distribution on subgraphs of $\cQ^n$ which is biased with respect to the number of edges between different levels of the hypercube.

\begin{definition}[Level-biased subgraphs of $\cQ^n$]\label{def:biased}
Given $n \in \mathbb{N}$ and $\mathbf{p} = (p_0, \dots, p_{n-1})\in[0,1]^n$, let $\cW^n_{\mathbf{p}}$\index{Wnp@$\cW^n_{\mathbf{p}}$} be a distribution on subgraphs of $\cQ^n$ where $W \sim \cW^n_{\mathbf{p}}$ is generated as follows: we set $V(W)\coloneqq V(\cQ^n)$ and, for each $i \in[n-1]_0$, each $e\in E(L_i,L_{i+1})$ is included in $W$ with probability $p_i$, independently of all other edges. 
\end{definition}

Roughly speaking, the above model has the advantage that, by choosing our probabilities $p_i$ appropriately, it will allow us to generate subgraphs of $\cQ^n$ where each vertex has the same number of up-neighbours in expectation.
Moreover, note that there is a joint distribution of $\cW^n_\mathbf{p}$ and $\cQ^n_{p}$ such that we have $\cW^n_\mathbf{p} \subseteq \cQ^n_{p}$, where $p$ is the maximum component of $\mathbf{p}$.\COMMENT{
We have for the $p$ level of the hypercube that the subgraphs are distributed the same here. 
For levels with smaller probabilities $p'$, say, we have that these are equal in between these levels in distribution to $\cQ^n_{p'} \subseteq \cQ^n_{p}$ and the claim follows.}

We are now in a position to define one further distribution on subgraphs of $\cQ^n$.
We will search for a near-spanning tree for $\cQ^n_\varepsilon$ in the graphs generated according to this distribution.

\begin{definition}[Percolation graph $\cP(n,\mathbf{p},M)$]\label{def: perc}
Given $n, M \in \mathbb{N}$ and $\mathbf{p} = (p_0, \dots, p_{n-1})\in[0,1]^n$,
we define $\cP(n,\mathbf{p}, M)$\index{PnpM1@$\cP(n,\mathbf{p},M)$} to be a distribution on subgraphs of $\cQ^n$ where $P \sim \cP(n,\mathbf{p},M)$ is generated as follows.
Let $R\sim \mathit{Res}(\cQ^n, 1/100)$ and $W \sim \cW^n_{\mathbf{p}}$.
For each $x\in V(\cQ^n)$, if $d_W^\uparrow(x)\geq M$, let $B(x)\subseteq N^\uparrow_W(x)$ be a uniformly random set of size $M$ (otherwise, let $B(x)\coloneqq\varnothing$), and let $E(x)$ be the set of edges joining $x$ to each $y\in B(x)$.
Let $W'$ be the spanning subgraph of $W$ with edge set $\bigcup_{x\in V(\cQ^n)}E(x)$.
The graph $P\subseteq\cQ^n$ is then given by setting $P\coloneqq W'-R$.
\end{definition}

\begin{remark}\label{rmk:perc-indep}
Observe that, given any two distinct edges $e,e'\in E(\cQ^n)$, the events $e\in E(W')$ and $e'\in E(W')$ are mutually dependent if and only if for some $i \in [n]$ we have $e, e' \in E(L_{i-1},L_i)$ with $e \cap e' = \{v\}$ for some $v \in L_{i-1}$.
Otherwise, these events are independent.
In particular, if $e\in E(L_{i-1},L_i)$ and $e'\in E(L_{j-1},L_j)$ with $i\neq j$, then these events are always independent.
\end{remark}

Note that $\cP(n,\mathbf{p}, M) \subseteq \cW^n_\mathbf{p}$ by definition, and therefore we have a joint distribution of $\cP(n,\mathbf{p}, M)$ and $\cQ^n_{p}$ such that $\cP(n,\mathbf{p}, M) \subseteq \cQ^n_{p}$, where $p$ is the maximum component of $\mathbf{p}$.

\begin{definition}[Feasible $(n,\mathbf{p}, M)$]\label{feasible}
We say that the tuple $(n,\mathbf{p}, M)$ is \emph{feasible} if 
\begin{enumerate}[label=$(\mathrm{\roman*})$]
    \item\label{feasible1} $p_i=0$ for all $9n/10<i<n$,\COMMENT{The reason this is being enforced here is because it means the $p_i$ for $i> 9n/10$ have to be 0 and can be forgotten about later. Otherwise we'd have to keep saying they were small enough that we could contain these $P$ graphs in a $Q^n_\eps$ where we know what the $\eps$ is.}
    \item\label{feasible2} $max_{i \in [n-1]_0} p_i < 1/10$ and $M>1600$,
    \item\label{feasible3} there exists $t \in \mathbb{R}$ with $600<t<100M$ such that $P \sim \cP(n,\mathbf{p}, M)$ satisfies $\mathbb{P}[e \in E(P)]=t/n$ for all $e \in \bigcup_{i=0}^{\lfloor9n/10\rfloor}E(L_i, L_{i+1})$\COMMENT{
We are stopping at edges up to level $9n/10$ here because above say level $n-M$, there are no edges by definition. 
So it becomes convenient for what we want later to do it this way.}.
\end{enumerate}
\end{definition}

\begin{remark}\label{rmk: W'}
Let $(n,\mathbf{p}, M)$ be feasible, where $\mathbf{p}=(p_0, \dots, p_{n-1})$.
Note that $p_0$ determines the value of $p_i$ for all $i\in[\lfloor9n/10\rfloor]$.\COMMENT{If $e \in E(L_i, L_{i+1})$ and $e' \in E(L_{j}, L_{j+1})$, then the only variables (given that $M$ is fixed for both) which determine whether or not these edges are present are $p_i$ and $p_j$.
By extension, $p_0$ determines all other values.}
Furthermore, let $P \sim \cP(n,\mathbf{p}, M)$. 
We can generate $P$ by first sampling $W \sim \cW^n_{\mathbf{p}}$ and $R \sim \mathit{Res}(\cQ^n, 1/100)$, and then defining the graph $W'$ as described in \cref{def: perc}.
Let $t' \coloneqq t/(\frac{99}{100})^2$, where $t$ is as in \cref{feasible}\ref{feasible3}.
Since $(n, \mathbf{p}, M)$ is feasible, for all $e\in \bigcup_{i=0}^{\lfloor9n/10\rfloor}E(L_i, L_{i+1})$ we have \COMMENT{Because to get from $W'$ to $P$, you just remove $R$. 
Any edge has a $(99/100)^2$ probability of surviving this.}
\[\mathbb{P}[e \in E(W')]= t'/n.\]
Furthermore, for all $i\in[\lfloor9n/10\rfloor]_0$, given $e, e' \in E(L_i, L_{i+1})$ with $e \ne e'$, we have that\COMMENT{
By \cref{rmk:perc-indep}, we have equality except in the case where $e, e'$ share a vertex in $L_i$.
In such a case, let $e = \{x,y\}$ and $e' =\{x, z\}$.
To get $W'$ first we generate $W$.
If $\{y,z\}\nsubseteq N^\uparrow_W(x)$, then the probability that both edges are in $W'$ is $0$, so we may condition on the event that $\{y,z\}\subseteq N^\uparrow_W(x)$.
Furthermore, condition on the event that $|N^\uparrow_W(x)| = D$.
If $D < M$, all probabilities are $0$.
If $D\geq M$, we have that
\[\mathbb{P}[e,e' \in E(W')] = (1/100)^2 M(M-1)/(D(D-1)) \leq (1/100)^2 M^2/D^2 = \mathbb{P}[e \in E(W')]^2.\]
By adding over all values of $D$, the remark follows.
}
\[\mathbb{P}[e,e' \in E(W')] \le \mathbb{P}[e \in E(W')]^2 = (t'/n)^2.\]
\end{remark}

Let us justify the last inequality in \cref{rmk: W'}.
By \cref{rmk:perc-indep}, the events that $e\in E(W')$ and $e'\in E(W')$ are independent (in which case we have equality) except in the case where $e$ and $e'$ share a vertex in $L_i$.
In such a case, let $e=\{x,y\}$ and $e'=\{x,z\}$.
In order to generate $W'$, we first generate $W\sim\mathcal{W}^n_{\mathbf{p}}$.
If $\{y,z\}\nsubseteq N^\uparrow_W(x)$, then the probability that both edges are in $W'$ is $0$, so we may condition on the event that $\{y,z\}\subseteq N^\uparrow_W(x)$.
Furthermore, condition on the event that $|N^\uparrow_W(x)| = D$, for some $2\leq D\leq n-i$.
Let us denote the probabilities in this conditional space by $\mathbb{P}^*$.
If $D < M$, then $\mathbb{P}^*[e \in E(W')]=\mathbb{P}^*[e'\in E(W')]=0$.
On the other hand, if $D\geq M$, we have that
\[\mathbb{P}^*[e,e' \in E(W')] = \frac{M}{D}\frac{M-1}{D-1} \leq \frac{M^2}{D^2} = \mathbb{P}^*[e \in E(W')]^2.\]
The claimed inequality follows since it holds uniformly over all values of $D$.

From here on, where it is clear from the context, we will use $p_0, \dots, p_{n-1}$ to denote the components of each probability vector $\mathbf{p}$, and will use $t$ to denote the value $t$ in \cref{feasible} and $t'$ to denote the value $t'$ in \cref{rmk: W'}.  

\begin{proposition}\label{prop: feasible}
For all $\eps \in(0,1/10)$\COMMENT{This bound can be improved.}, $M>1600$, and  $n \in \mathbb{N}$ such that $0< 1/n\ll1/M,\varepsilon$ there exists a tuple $(n,\mathbf{p}, M)$ which is feasible and such that $p_{i} \le \eps$ for all $i\in [n-1]_0$.
\end{proposition}

\begin{proof}
Let $P \sim \cP(n, \mathbf{p}, M)$, for some  $\mathbf{p}$ which will be determined later.
We generate $P$ by first sampling $W \sim \cW^n_\mathbf{p}$ and $R\sim \mathit{Res}(\cQ^n, 1/100)$.
Let $j\in[\lfloor9n/10\rfloor]_0$ be fixed and let $e \in E(L_j, L_{j+1})$.
Let $x \in L_j$ be incident with $e$.
Let $\mathcal{A}$ be the event that $e\in E(W)$.
For each $k \in [n-j]_0$, let $\mathcal{B}_k$ be the event that $d_{W}^\uparrow(x) = k$.
Let $\mathcal{C}$ be the event that $e\in E(P)$.
For each $i\in[n-M]_0$, let 
\[f_i(y)\coloneqq\Big(\frac{99}{100}\Big)^2\frac{M}{n-i}\sum_{k = M}^{n-i}\binom{n-i}{k}y^k(1-y)^{n-i-k}.\]
Then, we have that\COMMENT{The third equality follows from \cref{def: perc}.}
\begin{align*}
\mathbb{P}[\mathcal{C}] & = \sum_{k = M}^{n-j}\mathbb{P}[\mathcal{C}\mid \mathcal{A}\wedge \mathcal{B}_k]\,\mathbb{P}[\mathcal{A}\mid \mathcal{B}_k]\,\mathbb{P}[\mathcal{B}_k]\\
&= \sum_{k = M}^{n-j}\Big(\frac{99}{100}\Big)^2\frac{M}{k}\frac{k}{n-j}\binom{n-j}{k}p_j^k(1-p_j)^{n-j-k} =f_j(p_j).
\end{align*}

Let $m \coloneqq \min_{i \in[\lfloor9n/10\rfloor]_0} f_i(\eps)$.
\begin{claim}\label{claim: constant bound}
We have $\frac{600}{n} < m < \frac{100M}{n}$.
\end{claim}

\begin{claimproof}[Proof of \cref{claim: constant bound}]
Let $i \in[\lfloor9n/10\rfloor]_0$ be such that $f_i(\eps) = m$.
Clearly\COMMENT{Here we use that $i \le 9n/10$ and that the summation is part of a binomial sum that sums to 1.}, 
\[f_i(\eps) \le \frac{100M}{n}\sum_{k = M}^{n-i}\binom{n-i}{k}\eps^k(1-\eps)^{n-i-k} < \frac{100M}{n}.\]
Moreover, we have that 
\[f_i(\eps) = \Big(\frac{99}{100}\Big)^2\frac{M}{n-i}\sum_{k = M}^{n-i}\binom{n-i}{k}\eps^k(1-\eps)^{n-i-k} > \frac{1}{2}\Big(\frac{99}{100}\Big)^2\frac{M}{n},\] 
as $M \le \eps(n-i)/2$.
\end{claimproof}

For each $i\in[\lfloor9n/10\rfloor]_0$ such that $f_i(\eps) > m$, since $f_i(x)$ is continuous, by the intermediate value theorem we can choose some $p_i \in (0, \eps)$ such that $f_i(p_i) = m$.
This determines the probability vector $\mathbf{p} = (p_0, \dots, p_{n-1})$.
By \cref{claim: constant bound}, the tuple $(n,\mathbf{p}, M)$ is feasible, hence the statement is satisfied.
\end{proof}

In order to construct the near-spanning tree, we will generate a graph $P \sim \cP(n,\mathbf{p},M)$, for some feasible $(n,\mathbf{p},M)$, and will be interested in whether or not there exists a chain in $P$ from some vertex $x\in L_m$ to some vertex $y \in L_{m'}$, for $m' > m$ and $x \subseteq y$.
Note that the presence and absence of such chains in $P$ are highly dependent.
Thus, in order to show that such chains exist with high probability, we will consider the number of $x$-$y$ chains and bound its variance.
We do so in the following lemma.
In order to state it, we first need to set up some notation.

Given $x \in L_m$ and $y \in L_{m'}$ with $m' \ge m$, we denote by $\mathcal{X}_{x,y}$\index{Xxy@$\mathcal{X}_{x,y}$} the collection of $x$-$y$ chains in $\cQ^n$. 
For each $X \in \mathcal{X}_{x,y}$ and any graph $G \subseteq \cQ^n$, let $Y_X(G)$ be the corresponding indicator variable which takes value $1$ if $X\subseteq G$ and $0$ otherwise.
Let $Y_{x,y}(G) \coloneqq \sum_{X \in \mathcal{X}_{x,y}} Y_X(G)$.
Whenever $G$ is clear from the context, we will simply write $Y_{x,y}$\index{Yxy@$Y_{x,y}$}.
We define 
\[\Delta(Y_{x,y}) \coloneqq \sum_{\substack{(X,X')\in\mathcal{X}_{x,y}^2\\X\neq X'}} \Cov[Y_X Y_{X'}],\]
so $\Var[Y_{x,y}]=\Delta(Y_{x,y})+\sum_{X\in\mathcal{X}_{x,y}}\Var[Y_X]$.\index{Delta1Y@$\Delta(Y_{x,y})$}

\begin{lemma}\label{lem: cov}
Let $P \sim \cP(n, \mathbf{p}, M)$, where $(n, \mathbf{p}, M)$ is feasible with $0< 1/n \ll 1/M$.
Let $1\leq m<m'\leq 9n/10$ with $m'-m+1\geq n/4 -1$.
Let $x \in L_m$ and $y \in L_{m'}$ with $x \subseteq y$. 
Then,
\[\Delta(Y_{x,y})\leq2\mathbb{E}[Y_{x,y}]^2.\]
\end{lemma}

The proof of \cref{lem: cov} makes use of the analysis in the proof of a similar lemma of \citet[Lemma 7]{KKO}.
In order to shorten our analysis here, we first state a partial result which follows from the analysis of \cite{KKO}.
For this, we first need to give some more definitions.

Fix $x \in L_m$ and $y \in L_{m'}$ with $m' \ge m$ and $x\subseteq y$.
Observe that $|\mathcal{X}_{x,y}|=\dist(x,y)!=(m'-m)!$ depends only on the distance between $x$ and $y$.
For each $k\in[n]$, let $R_k\coloneqq(k-1)!$\index{Rk@$R_k$}.
Given any $X, X' \in \mathcal{X}_{x,y}$ with $X\neq X'$, let $i(X,X')\coloneqq |V(X) \cap V(X')|-2$\index{iX@$i(X,X')$}, let $s(X,X')$\index{sX@$s(X,X')$} be the number of connected components of $X-V(X')$, and let $\ell(X,X')$\index{lX@$\ell(X,X')$} be the largest order over these components.\COMMENT{That is, the size of the component of $X - V(X')$ with the most number of vertices.}

Next, we define the set of possible intersection patterns for two chains.
Let $k\coloneqq m'-m+1$.
Given any chains $X, X' \in \mathcal{X}_{x,y}$, let $A(X,X')$\index{AX@$A(X,X')$} be the collection of indices $a \in [k-2]$ for which $X$ and $X'$ agree on their $(a+1)$-th elements (where we consider $x$ to be the first element of $X$ and $X'$).
An \emph{admissible $(i,\ell,s)$-pattern}\index{(i@$(i,\ell,s)$-pattern} is a set $A\subseteq [k-2]$ with $|A|=i$ such that the longest interval of consecutive elements in $[k-2] \setminus A$ contains exactly $\ell$ elements and such that the number of maximal intervals of consecutive elements in $[k-2] \setminus A$ is exactly $s$.
We denote by $\mathcal{A}_{i,\ell,s}$\index{Ai@$\mathcal{A}_{i,\ell,s}$} the set of all admissible $(i,\ell,s)$-patterns. 
Furthermore, we define $C_{i,\ell,s} \coloneqq |\mathcal{A}_{i,\ell,s}|$\index{Ci@$C_{i,\ell,s}$}.
Note that any pair of chains $X, X' \in \mathcal{X}_{x,y}$ with $i(X,X') = i$, $\ell(X,X') = \ell$ and $s(X,X') = s$ define an admissible $(i,\ell,s)$-pattern $A(X,X') \in \cA_{i,\ell,s}$.

Given a chain $X \in \cX_{x,y}$ and a pattern $A \in \mathcal{A}_{i,\ell,s}$, let $F(A)$ be the number of chains $X' \in \cX_{x,y}$ such that $A(X,X') = A$.
(Note that the definition of $F(A)$ is independent of $X$.)
Let $F_{i,\ell,s} \coloneqq \max_{A \in \mathcal{A}_{i,\ell,s}} F(A)$\index{Fi@$F_{i,\ell,s}$}.
Observe that $F_{i,\ell,s}$ is an upper bound on the number of chains $X'$ with $A(X,X') = A$.

Finally, for each triple $(i,\ell,s)\in[k-3]_0\times[k-2]^2$, let 
\[\Delta_{i,\ell,s}\coloneqq \sum_{\substack{(X,X')\in\mathcal{X}_{x,y}^2,\, X \ne X'\\i(X,X')=i,\, \ell(X,X')=\ell,\, s(X,X')=s}}\mathbb{E}[Y_X Y_{X'}].\]
Furthermore, let\index{Delta1i@$\Delta_{i,\ell,s}$} 
\[\Delta_0(Y_{x,y}) \coloneqq \sum_{\substack{(X,X')\in\mathcal{X}_{x,y}^2\\i(X,X')=0}} \Cov[Y_X Y_{X'}]   \qquad \text{ and } \qquad \Delta_1(Y_{x,y}) \coloneqq \sum_{\substack{(X,X')\in\mathcal{X}_{x,y}^2\\i(X,X')\in[k-3]}} \Cov[Y_X Y_{X'}].\]
Thus, $\Delta(Y_{x,y})=\Delta_0(Y_{x,y})+\Delta_1(Y_{x,y})$\index{Delta0@$\Delta_0(Y_{x,y})$, $\Delta_1(Y_{x,y})$}.
Note that, by summing $\Delta_{i,\ell,s}$ over all triples $(i,\ell,s)\in[k-3]\times[k-2]^2$, we obtain an upper bound for $\Delta_1(Y_{x,y})$.\COMMENT{In the following lemma it might be nicer to rewrite as $\Delta(Y) \le \mathbb{E}[Y]^2$ and the same for \cref{lem: cov} as that's all we'll ultimately use.}

\begin{lemma}[\cite{KKO}]\label{lem: KKO}
 For all $M>100$ there exists $n_0$ such that, for all $n\ge n_0$, the following holds.
 Let $x \in L_1$ and $y \in L_{n-1}$ with $x \subseteq y$. 
 Let $p \ge M/(2n)$.
 Let $Q \subseteq \cQ^n$ be a random subgraph chosen according to any distribution such that
 \[\frac{\Delta_{i,\ell,s}}{\mathbb{E}[Y_{x,y}]^2} \le  \frac{C_{i,\ell,s}F_{i,\ell,s}}{R_{n-1}p^{i}},\]
 for each possible choice of $(i,\ell,s)\in[k-3]\times[k-2]^2$\COMMENT{Need to avoid $i=0$ and $X=X'$.}.
 Then,
 \[\Delta_1(Y_{x,y}) \le \frac{100}{M} \mathbb{E}[Y_{x,y}]^2.\]
\end{lemma}

With this, we are finally ready to prove \cref{lem: cov}.

\begin{proof}[Proof of \cref{lem: cov}]
Let $P \sim \cP(n, \mathbf{p}, M)$, where $(n, \mathbf{p}, M)$ is feasible.
Recall, from \cref{def: perc}, that $P$ is generated by first sampling a set $R \sim \mathit{Res}(\cQ^n, 1/100)$ and a graph $W \sim \cW^n_{\mathbf{p}}$.
We then generate the graph $W'$ by choosing, for each $v \in \bigcup_{i=0}^{\lfloor9n/10\rfloor}L_i$, a set of $M$ up-neighbours uniformly at random from the set of up-neighbours $v$ has in $W$, provided $d^\uparrow_{W'}(v)\ge M$ (and by setting $d^\uparrow_{W'}(v) \coloneqq 0$ otherwise).
Let $t'\coloneqq t/(\frac{99}{100})^2$.
Thus, for all $e\in \bigcup_{i=0}^{\lfloor9n/10\rfloor}E(L_i, L_{i+1})$ we have by \cref{rmk: W'} that
\[\mathbb{P}[e \in W']=t'/n.\]

Let $k\coloneqq m'-m+1$, and let $X$ be a fixed $x$-$y$ chain in $\cQ^n$.
By \cref{rmk:perc-indep} it follows that\COMMENT{There are $k-1$ edges, each appearing in $W'$ with probability $t'/n$ independently of each other.  
There is then a $(99/100)^k$ probability that none of the $k$ vertices in the path are in $R$. (Again, each vertex being present in $R$ is independent of everything else.)}
\begin{equation}\label{eq:Yxyexp}
    \mathbb{E}[Y_{x,y}] = R_k \mathbb{P}[X\subseteq P] = R_k(t'/n)^{k-1}\Big(\frac{99}{100}\Big)^k.
\end{equation} 
Furthermore, for all $(i,\ell,s)\in[k-3]_0\times[k-2]^2$, we have that
\begin{equation}\label{eq:Delta1}
    \Delta_{i,\ell,s} \leq R_kC_{i,\ell,s}F_{i,\ell,s}\left(\frac{99t'}{100n}\right)^{2k-i-2}.
\end{equation}
To see this, note that we may first choose an $x$-$y$ chain $X$, for which there are $R_k$ choices.
Next, we choose an admissible $(i,\ell,s)$-pattern $A \in \cA_{i,\ell,s}$, of which there are $C_{i,\ell,s}$.
We then have at most $F_{i,\ell,s}$ choices for $x$-$y$ chains $X'$ with $A(X,X') = A$.
Next, we bound the number of vertices and edges of $X\cup X'$.
It is clear that $X$ has $k$ vertices and $k-1$ edges, and $|V(X')\setminus V(X)|=k-i-2$.
Moreover, observe that $|E(X')\setminus E(X)| \ge k-i-1$.\COMMENT{To see that there are at least $k-i-1$ edges in $E(X')\setminus E(X)$, consider the following argument. 
There are $k-i-2$ non-intersection points. 
For each of these vertices in $X'$ we have an edge incident to it (from below, say) which is not in $X$.
Furthermore, the last non-intersection vertex must also have another edge (going above), since both $X$ and $X'$ finish at the same vertex.
Therefore, there must be at least $k-i-1$ new edges.}
The bound finally follows by considering the probability that all these vertices and edges are present in $P$ and by \cref{rmk: W'}.\COMMENT{The final part of the remark tells us that when we consider two distinct up-edges coming from the same vertex $x$, the probability that both are present simultaneously in $W'$ is less that the probability of one of them being present squared (this is what we use).}

We are going to compute bounds for $\Delta_0(Y_{x,y})$ and $\Delta_1(Y_{x,y})$ separately, and then combine them to obtain the result.
We begin with a bound for $\Delta_1(Y_{x,y})$.
Combining \eqref{eq:Yxyexp} and \eqref{eq:Delta1}, it follows that, for all $(i,\ell,s)\in[k-3]\times[k-2]^2$,\COMMENT{We have that
\begin{align*}
    \frac{\Delta_{i,\ell,s}}{\mathbb{E}[Y_{x,y}]^2} &\le \frac{R_kC_{i,\ell,s}F_{i,\ell,s}(t'\frac{99}{100})^{2k-i-2}n^{2k-2}}{n^{2k-i-2}R_k^2(t')^{2k-2}(\frac{99}{100})^{2k}} = \frac{C_{i,\ell,s}F_{i,\ell,s}n^{i}}{R_k(t')^{i}(\frac{99}{100})^{i+2}} = \frac{C_{i,\ell,s}F_{i,\ell,s}}{R_k}\left(\frac{n}{t'}\right)^i\left(\frac1{\frac{99}{100}}\right)^{i+2}\\
    &=\frac{C_{i,\ell,s}F_{i,\ell,s}}{R_k}\left(\frac{n}{t'(\frac{99}{100})^{\frac{i+2}{i}}}\right)^i \leq\frac{C_{i,\ell,s}F_{i,\ell,s}}{R_k}\left(\frac{n}{t'(\frac{99}{100})^3}\right)^i,
\end{align*}
where the last inequality follows since $i\geq1$.} 
\begin{equation*}\label{eqn:1}
\frac{\Delta_{i,\ell,s}}{\mathbb{E}[Y_{x,y}]^2} \leq \frac{C_{i,\ell,s}F_{i,\ell,s}}{R_k}\left(\frac{n}{t'}\right)^i\left(\frac1{\frac{99}{100}}\right)^{i+2} \leq \frac{C_{i,\ell,s}F_{i,\ell,s}}{R_k}\left(\frac{n}{t'(\frac{99}{100})^3}\right)^i.
\end{equation*}
Note that $(\frac{m'-m+2}{n})t'(\frac{99}{100})^3 > 100$\COMMENT{We have $m' -m +1 \ge n/4 - 1$ so $(\frac{m'-m+2}{n}) \ge 1/4$.
Furthermore, we have $t'(\frac{99}{100})^3 = t(\frac{99}{100})> 600(\frac{99}{100})>400$ since $t$ comes from a feasible set of values $(n, \mathbf{p}, M)$.}.
It follows that we can apply \cref{lem: KKO} with $(\frac{m'-m+2}{n})t'(\frac{99}{100})^3$ and $m'-m+2$  playing the roles of $M$ and $n$ and $p = t'(\frac{99}{100})^3/n$\COMMENT{Note that $\left(\frac{m' -m +2}{(\frac{m'-m+2}{n})t'(\frac{99}{100})^3}\right)^i$ = $\left(\frac{n}{t'(\frac{99}{100})^3}\right)^i$.} to obtain that 
\begin{equation}\label{eq:Delta1bound}
    \Delta_1(Y_{x,y}) \le \frac{100}{(\frac{m'-m+2}{n})t'(\frac{99}{100})^3} \mathbb{E}[Y_{x,y}]^2\leq \mathbb{E}[Y_{x,y}]^2.
\end{equation}

We now turn our attention to $\Delta_0(Y_{x,y})$.
For any two chains $X,X'\in\mathcal{X}_{x,y}$ such that $i(X,X')=0$, we have that $X\cup X'$ has $2k-2$ vertices and the same number of edges.
Therefore, by \cref{rmk:perc-indep,rmk: W'} we have $\mathbb{E}[Y_XY_{X'}]\le(\frac{99t'}{100n})^{2k-2}$, and by \eqref{eq:Yxyexp} we have that\COMMENT{We may bound the number of pairs of chains by the total number of pairs of chains, which is $R_k^2$.
Then, using \eqref{eq:Yxyexp}, we have that
\begin{align*}
    \Delta_0(Y_{x,y})&\leq R_k^2\left(\left(\frac{t'}{n}(\frac{99}{100})\right)^{2k-2}-\left(\frac{t'}{n}\right)^{2k-2}(\frac{99}{100})^{2k}\right)=R_k^2\left(\frac{t'}{n}(\frac{99}{100})\right)^{2k-2}(1-(99/100)^2)\\
    &\leq2\frac{1}{100} R_k^2\left(\frac{99t'}{100n}\right)^{2k-2}=\frac{\frac{2}{100}}{(\frac{99}{100})^2}\mathbb{E}[Y_{x,y}]^2.
\end{align*}
}
\begin{equation}\label{eq:Delta0bound}
    \Delta_0(Y_{x,y})\leq R_k^2\Big(\frac{99t'}{100n}\Big)^{2k-2}\Big(1-\Big(\frac{99}{100}\Big)^2\Big)\leq\mathbb{E}[Y_{x,y}]^2.
\end{equation}

The conclusion follows immediately by combining \eqref{eq:Delta1bound} and \eqref{eq:Delta0bound}.
\end{proof}

In order to proceed further, we will consider unions of independent graphs $P\sim\cP(n,\mathbf{p},M)$.

\begin{definition}\label{def:percC}
Let $n,M,C \in \mathbb{N}$ and $\mathbf{p}\in[0,1]^n$.
We define $\cP^C(n,\mathbf{p}, M)$\index{PnpM2@$\cP^C(n,\mathbf{p}, M)$} to be a distribution on subgraphs of $\cQ^n$ such that $P \sim \cP^C(n,\mathbf{p}, M)$ is generated by taking $C$ independently generated graphs $P_i \sim \cP(n,\mathbf{p}, M)$ and setting $P \coloneqq \bigcup_{i=1}^C P_i$.
For each $i \in [C]$, there is a set $R_i \sim \mathit{Res}(\cQ^n, 1/100)$ associated with $P_i$.
Let $R \coloneqq \bigcap_{i=1}^C R_i$.
We say that $R$ is the \emph{reservoir associated with} $P$.
\end{definition}

It follows from \cref{def: perc,def:percC} that there is a joint distribution of $\cP^C(n,\mathbf{p}, M)$ and $\cQ^n_{\min\{1,Cp\}}$ such that $\cP^C(n,\mathbf{p}, M) \subseteq \cQ^n_{\min\{1,Cp\}}$, where $p = \max_{i\in[n-1]_0} p_i$.
Note that for all $x \in V(\cQ^n)$ we have that $\mathbb{P}[x \in R] = (1/100)^C$.

Our next goal is to prove that, by choosing constants appropriately, there is a high probability that there exists an $x$-$y$ chain in $P \sim \cP^C(n,\mathbf{p},M)$, even if we restrict the set of `valid' chains to a significant subset of the total.
For this, we will make use of \cref{lem: cov}. 
Given any vertices $x \in L_m$ and $y \in L_{m'}$ with $x \subseteq y$, any set $\mathcal{Z}\subseteq\mathcal{X}_{x,y}$, and any graph $G\subseteq \cQ^n$, we denote the number of $x$-$y$ chains $X\in\mathcal{Z}$ such that $X\subseteq G$ by $Y(\mathcal{Z},G)$.

\begin{corollary}\label{cor: boost}
For $n, C\in \mathbb{N}$ and $\eta, \alpha>0$ such that $0 < 1/n \ll 1/C \ll \eta, \alpha$ and any feasible $(n, \mathbf{p}, M)$  with $0< 1/n \ll 1/M$, the following holds.
Let $1\leq m<m'\leq 9n/10$ with $m'-m+1\geq n/4-1$.
Let $x \in L_m$ and $y \in L_{m'}$ with $x \subseteq y$. 
Let $\mathcal{Z}_{x,y}\subseteq \mathcal{X}_{x,y}$ be such that $|\mathcal{Z}_{x,y}|\geq\alpha|\mathcal{X}_{x,y}|$.
Let $P \sim \cP^C(n,\mathbf{p}, M)$.
Then, 
\[\mathbb{P}[Y(\mathcal{Z}_{x,y},P)>0] \ge 1 - \eta.\]
\end{corollary}

\begin{proof}
For each $i \in [C]$, let $P_i \sim \cP(n, \mathbf{p}, M)$, and let $P \coloneqq \bigcup_{i=1}^{C}P_i$.
Let $Y_i\coloneqq Y_{x,y}(P_i)$ and $Z_i\coloneqq Y(\mathcal{Z}_{x,y},P_i)$, and let
\[\Delta(Z_i)\coloneqq\sum_{\substack{(X,X')\in\mathcal{Z}_{x,y}^2\\X\neq X'}}\Cov[Y_X(P_i)Y_{X'}(P_i)].\]
Note that\COMMENT{Let $X$ be a sum of $n$ indicator random variables $X_i$.
Then, we have 
\begin{align*}
    \mathbb{E}[X^2]-\mathbb{E}[X]^2=\Var(X)&=\sum_{i=1}^n\Var(X_i)+\sum_{i\neq j}\Cov(X_i,X_j)=\sum_{i=1}^n(\mathbb{E}[X_i^2]-\mathbb{E}[X_i]^2)+\sum_{i\neq j}\Cov(X_i,X_j)\\
    &\leq\mathbb{E}[X]+\Delta(X).
\end{align*}
\eqref{eq:chaincovariance} follows by reordering.}
\begin{equation}\label{eq:chaincovariance}
    \mathbb{E}[Y_i^2] \leq \Delta(Y_i) + \mathbb{E}[Y_i] + \mathbb{E}[Y_i]^2.
\end{equation}
We also have 
\begin{equation}\label{eq:chaintrivial}
    \mathbb{E}[Z_i^2]\leq\mathbb{E}[Y_i^2]
\end{equation}
and, since all $x$-$y$ chains are equiprobable,
\begin{equation}\label{eq:chaintrivial2}
    \mathbb{E}[Z_i]^2\geq\alpha^2\mathbb{E}[Y_i]^2.
\end{equation}

Let $k\coloneqq m'-m+1$.
By \eqref{eq:Yxyexp}, we have that $\mathbb{E}[Y_i] = R_k(t'/n)^{k-1}(99/100)^{k}$, where $t'$ is the value given in \cref{rmk: W'}.
Recall that $R_k=|\mathcal{X}_{x,y}| = (k-1)!$.
We have by Stirling's formula that $\mathbb{E}[Y_i] >1$.\COMMENT{Recall that $t'=t/(1-\delta)^2$, where here $\delta=1/100$. We have by Stirling's formula that
\[
    \mathbb{E}[Y_i] > (k-1)!\left(\frac{t'}{2n}\right)^{k-1} > \left(\frac{k}{2e}\right)^{k-1}\left(\frac{t'}{2n}\right)^{k-1} > \left(\frac{t}{99}\right)^{k-1}>1.
\]
}
Therefore, $\mathbb{E}[Y_i] \le \mathbb{E}[Y_i]^2$.
Moreover, it follows by \cref{lem: cov} that $\Delta(Y_i) \le 2\mathbb{E}[Y_i]^2$.
So $\mathbb{E}[Y_i^2] \le 4\mathbb{E}[Y_i]^2$ by \eqref{eq:chaincovariance}.
Combining this with  \eqref{eq:chaintrivial}, \eqref{eq:chaintrivial2} and \cref{prop: CZ} we obtain
\[\mathbb{P}[Z_i = 0] \le 1 - \frac{\mathbb{E}[Z_i]^2}{\mathbb{E}[Z_i^2]} \le 1 - \frac{\alpha^2\mathbb{E}[Y_i]^2}{\mathbb{E}[Z_i^2]} \le 1 - \frac{\alpha^2\mathbb{E}[Y_i]^2}{\mathbb{E}[Y_i^2]} \le 1 - \alpha^2/4.\]
It follows that 
\[\mathbb{P}[Y(\mathcal{Z}_{x,y},P)=0] = \prod_{i \in [C]} \mathbb{P}[Z_i=0] \le (1-\alpha^2/4)^C \le \eta.\qedhere\]
\end{proof}

When performing our analysis on the structure of $P$, the dependence of chains on each other becomes difficult to take into account.
In order to deal with this issue, we will show that, with high probability, it suffices to consider only chains which lie in some large subsets of the total sets of chains, with the property that the presence or absence of a chain in one of these large subsets is independent from chains of all other subsets.
(Note that \cref{cor: boost} works for these sets of chains as long as they are not too small.)
The next two lemmas guarantee the existence of such sets.
In \cref{lem: chain count} we prove that, assuming  $x, x' \in L_{m}$, and $y, y' \in L_{m'}$, where $y,y'$ are far apart,  one can construct very large sets of chains between the pairs $x,y$ and $x',y'$, which are independent in the sense described above.
Then, in \cref{lem: randsethit} we will prove that we can pick many endpoints $y \in L_{m'}$ in such a way that they are suitably far apart.

Given $0\leq m<m'\leq n$, let $x,x' \in L_m$ and $y,y' \in L_{m'}$ with $x\subseteq y$ and $x'\subseteq y'$. 
We denote by $\mathcal{X}_{x,y}^{\neg x',y'}$\index{Xxy2@$\mathcal{X}_{x,y}^{\neg x',y'}$} the collection of chains $X \in \mathcal{X}_{x,y}$ for which there is no $X' \in \mathcal{X}_{x',y'}$ with $V(X)\cap V(X')\neq\varnothing$.

\begin{lemma}\label{lem: chain count}
For all $n \ge 100$, the following holds.
Let $1\leq m<m'\leq n-1$ be such that $n/4 -1\leq k\coloneqq m'-m+1\leq n/2$.
Let $x, x' \in L_m$ and $y, y' \in L_{m'}$ with $x \subseteq y$ and $x' \subseteq y'$ be such that  $\dist(x,x')=2$ and $\dist(y,y')\geq9k^2/(10n)$\COMMENT{If there are no vertices satisfying this, the statement is vacuously true.}.
Then, 
\[|\mathcal{X}_{x,y}^{\neg x',y'}| \ge \left(1-\frac{60000}{n}\right)|\mathcal{X}_{x,y}|.\]
\end{lemma}

\begin{proof}
We may assume that $x\cup x' \subseteq y\cap y'$, since otherwise $\mathcal{X}_{x,y}^{\neg x',y'} = \cX_{x,y}$.
Let $b \coloneqq |y\cap y'|$.
We have that $b\leq m'-9k^2/(20n)$.\COMMENT{Think of the vertices as $0$-$1$ coordinate vectors.
If they are distance $s$ apart, then there are $s$ coordinates on which they differ.
None of these $s$ coordinates can appear in $y \cap y'$.
Then, $y$ (say) has $s/2$ of these $s$ on which they differ. 
The rest of $y$ can still appear in $y \cap y'$, giving $m'-s/2$ $1$ coordinates.}
Let $H$ denote the smallest subcube of $\cQ^n$ which contains both $x \cup x'$ and $y \cap y'$\COMMENT{Note that since $x\cup x' \subseteq y\cap y'$ we have that $H$ may be defined as the subcube spanned by the vertices of the chains of $\mathcal{X}_{x \cup x',y \cap y'}$.}.
For each $i \in [b]\setminus[m]$, let $\mathcal{X}^i_{x,y}\subseteq\mathcal{X}_{x,y}$ be the set of chains $X\in\mathcal{X}_{x,y}$ such that $V(X)\cap L_i \cap V(H)\neq\varnothing$.
Note that $\mathcal{X}_{x,y}^{\neg x',y'} \supseteq \cX_{x,y} \setminus \bigcup_{i \in [b]\setminus [m]}\cX_{x,y}^i$ and
\begin{equation}\label{equa:chainsubset2}
    |\mathcal{X}^i_{x,y}|=\binom{b-m-1}{i-m-1}(m'-i)!(i-m)!.
\end{equation}
Indeed,\COMMENT{Note that the set of all $x$-$y$ chains can be produced by choosing a maximal set of incomparable sets $x\subseteq z\subseteq y$ and then considering the concatenation of all possible $x$-$z$ chains with all possible $z$-$y$ chains.} there are $\binom{b-m-1}{i-m-1}$ choices to fix an element $z \in V(H) \cap L_i$.
(To see this, consider that $H$ is itself a cube of dimension $b-m-1$, and we are choosing a vertex $z$ from the $(i-m-1)$-th level of this cube.)
Then, there are $(i-m)!$ $x$-$z$ chains, and $(m'-i)!$ $z$-$y$ chains.

Recall that $|\mathcal{X}_{x,y}| = (k-1)!$.
By comparing this with \eqref{equa:chainsubset2} and simplifying, for all $i\in[b]\setminus[m]$ we obtain \COMMENT{We have that
\[\frac{|\mathcal{X}^i_{x,y}|}{|\mathcal{X}_{x,y}|}=\frac{\binom{b-m-1}{i-m-1}(m'-i)!(i-m)!}{(k-1)!}=\frac{(b-m-1)!(m'-i)!(i-m)!}{(i-m-1)!(b-i)!(k-1)!}=\frac{(b-m-1)!(m'-i)!(i-m)}{(b-i)!(m'-m)!}.\]
Then, the upper bound follows by observing that $b\leq m'$.}
\begin{equation}\label{equa:chainsubset3}
    \frac{|\mathcal{X}^i_{x,y}|}{|\mathcal{X}_{x,y}|}=\frac{i-m}{m'-m}\prod_{j=1}^{i-m-1}\frac{b-m-j}{m'-m-j}\leq\frac{i-m}{k-1}.
\end{equation}

We now split the analysis into two cases.
First, when $i$ is small, we bound \eqref{equa:chainsubset3} directly.
For all $i\in[b]\setminus[m]$ with $i \le m + 64$, it follows from \eqref{equa:chainsubset3} that\COMMENT{We have that
\[\frac{|\mathcal{X}^i_{x,y}|}{|\mathcal{X}_{x,y}|}\leq\frac{i-m}{k-1}\leq\frac{64}{n/4-2}=\frac{256}{n-8}\leq\frac{300}{n},\] 
where the last inequality holds for $n\geq100$.}
\begin{equation}\label{equa:chainsubset4}
    \frac{|\mathcal{X}^i_{x,y}|}{|\mathcal{X}_{x,y}|}\leq\frac{64}{n/4-1}\leq\frac{300}{n}.
\end{equation}
On the other hand, for each $i\in[b-1]\setminus[m]$, by \eqref{equa:chainsubset3} we have that\COMMENT{We have that 
\[\frac{|\mathcal{X}^i_{x,y}|}{|\mathcal{X}_{x,y}|}\frac{i+1-m}{i-m}\frac{b-i}{m'-i}=\frac{i+1-m}{i-m}\frac{b-i}{m'-i}\frac{i-m}{m'-m}\prod_{j=1}^{i-m-1}\frac{b-m-j}{m'-m-j}=\frac{i+1-m}{m'-m}\prod_{j=1}^{i-m}\frac{b-m-j}{m'-m-j}=\frac{|\mathcal{X}^{i+1}_{x,y}|}{|\mathcal{X}_{x,y}|}.\]
Then, the upper bound follows by observing that (recall $n \ge 100$)
\[\frac{b-i}{m'-i}\leq\frac{b}{m'}\leq1-\frac{9k^2}{20nm'}\leq1-\frac{9k^2}{20n^2} = 1 - \frac{9(n/4-1)^2}{20n^2}\leq1-\frac{9n^2}{20\cdot16n^2} + \frac{18}{20\cdot16\cdot 100}<\frac{312}{320}=\frac{39}{40}.\]}
\[\frac{|\mathcal{X}^{i+1}_{x,y}|}{|\mathcal{X}_{x,y}|}=\frac{|\mathcal{X}^i_{x,y}|}{|\mathcal{X}_{x,y}|}\frac{i+1-m}{i-m}\cdot\frac{b-i}{m'-i}\leq\frac{39}{40}\frac{i+1-m}{i-m}\frac{|\mathcal{X}^i_{x,y}|}{|\mathcal{X}_{x,y}|}.\]
For all $i\in[b-1]\setminus[m+64]$\COMMENT{That is, $i\geq m+65$}, this yields\COMMENT{We have that $(i-m+1)/(i-m)$ is maximised when $i$ is minimised, and the case $i=m+65$ yields exactly $(66/65)(39/40)=99/100$.}
\begin{equation}\label{equa:chainsubset5}
    \frac{|\mathcal{X}^{i+1}_{x,y}|}{|\mathcal{X}_{x,y}|}\leq \frac{99}{100}\frac{|\mathcal{X}^i_{x,y}|}{|\mathcal{X}_{x,y}|}.
\end{equation}

Finally, by combining \eqref{equa:chainsubset4} and \eqref{equa:chainsubset5}, and considering a geometric series, we conclude that\COMMENT{We have shown that passing from $i$ to $i+1$ above has the effect of passing along a geometric series.
That is, we have 
\[\frac{\sum_{i=m+1}^{b}|\mathcal{X}^{i}_{x,y}|}{|\mathcal{X}_{x,y}|}\le \frac{300}{n}(64 + 99/100 + (99/100)^2 +\ldots) \le \frac{300}{n}\left(100 + \frac{1}{1-(99/100)}\right)=\frac{60000}{n}.\]}
\[\frac{|\mathcal{X}_{x,y}^{\neg x',y'}|}{|\mathcal{X}_{x,y}|}\geq\frac{|\mathcal{X}_{x,y}| - \sum_{i=m+1}^{b}|\mathcal{X}^{i}_{x,y}|}{|\mathcal{X}_{x,y}|} \ge 1 - \frac{60000}{n}.\qedhere\]
\end{proof}

\begin{remark}\label{rmk:2}
\cref{lem: chain count} holds similarly if $\dist(x,x')\geq9k^2/(10n)$ and $\dist(y,y')=2$.\COMMENT{This follows by considering the cube from `top to bottom'. 
All calculations in the proof can be done in the exact same way.}
\end{remark}

\begin{proposition}\label{prop: ballradiushit}
Let $0< 1/n\ll \gamma, 1/k\leq 1$, where $n, k \in \mathbb{N}$, and let $S \subseteq V(\cQ^{n})$ be such that, for all distinct $x, x' \in S$, we have $\dist(x,x') \geq \gamma n$.
Then, for any $y \in L_m$ such that $m \geq n/8$, and for every $\gamma m/2\leq t\leq (1-\gamma/2)m$, we have $|y^{(t)}\cap B^{k}_{\cQ^n}(S)| \le |y^{(t)}|2^{-\gamma n/200}$.
\end{proposition}

\begin{proof}
Let $m \ge n/8$ and $y \in L_m$.
Let $\gamma m/2\leq t\leq (1-\gamma/2)m$ and let $S' \subseteq S$ be the set of all those $x \in S$ for which $B^k_{\cQ^n}(x) \cap y^{(t)} \ne \varnothing$.
We have that\COMMENT{Where the final inequality uses that $|B^k_{\cQ^n}(x)| \le 2n^k$.
The 2 is because for $k=1$ we have $|B^k_{\cQ^n}(x)| = n+1$.} \begin{equation}\label{equa:100}
|B^k_{\cQ^n}(S)\cap y^{(t)}| \le \sum_{x \in S'}|B^k_{\cQ^n}(x)\cap y^{(t)}| \le 2n^k|S'|.
\end{equation}
Moreover, for every $x, x' \in S'$ we have that $B^{\gamma n/3}_{\cQ^n}(x) \cap B^{\gamma n/3}_{\cQ^n}(x')=\varnothing$, and, therefore,
\begin{equation}\label{equa:101}
|S'| ({\min_{x \in S'}}{|B^{\gamma n/3}_{\cQ^n}(x)\cap y^{(t)}|})\le |y^{(t)}|.
\end{equation}

\begin{claim}
For every $x \in S'$ we have $|B^{\gamma n/3}_{\cQ^n}(x)\cap y^{(t)}| \ge 2^{\gamma m/20}$.
\end{claim}

\begin{claimproof}
Let $x' \in B^k_{\cQ^n}(x) \cap y^{(t)}$.
Let $z \subseteq x'$ be such that $z \in L_{\gamma m/7}$ (recall that $t \ge \gamma m/2$).
Since $y\in L_m$ we have that $|y \setminus x'| = m-t$.
Let $z' \subseteq y\setminus x'$ be such that $z' \in L_{\gamma m/7}$ (recall that $t \le (1-\gamma/2)m$).
It follows that $(x'\setminus z)\cup z' \in B^{\gamma n/3}_{\cQ^n}(x)\cap y^{(t)}$.
Note that there are $\binom{t}{\gamma m/7}$ choices for $z$ and $\binom{m-t}{\gamma m/7}$ choices for $z'$.
It follows that\COMMENT{For the final inequality, note that $\binom{n}{m} \ge (n/m)^{m}$ and, therefore, \[\binom{m-t}{\gamma m/7} \ge ((m-t)/(\gamma m/7))^{\gamma m/7} \ge ((\gamma m/2)/(\gamma m/ 7))^{\gamma m/7} \ge 3.5^{\gamma m/7}.\]
Something similar holds for $\binom{t}{\gamma m/7}$.}
\[\left |B^{\gamma n/3}_{\cQ^n}(x) \cap y^{(t)}\right| \geq \binom{m-t}{\gamma m/7}\binom{t}{\gamma m/7}\geq 2^{\gamma m/20}.\qedhere\]
\end{claimproof}

Combining \eqref{equa:100}, \eqref{equa:101} and the above claim we have 
\[|B^k_{\cQ^n}(S)\cap y^{(t)}| \le  \frac{2n^k|y^{(t)}|}{{\min_{x \in S'}}{|B^{\gamma n/3}_{\cQ^n}(x)\cap y^{(t)}|}} \le 2n^k |y^{(t)}|2^{-\gamma m/20} \le |y^{(t)}|2^{-\gamma n/200}.\qedhere\]
\end{proof}

Given $x,y \in V(\cQ^n)$ with $x \subseteq y$ and $S\subseteq V(\cQ^n)$, we denote by $\mathcal{X}_{x,y}^{\neg S}$\index{Xxy3@$\mathcal{X}_{x,y}^{\neg S}$} the collection of chains $X \in \mathcal{X}_{x,y}$ for which $V(X) \cap S = \varnothing$.
\begin{lemma}\label{lem: pathball}\label{lem:pathsinthit}
Let $0 <1/n\ll \gamma, 1/k\leq 1$ where $n, k \in \mathbb{N}$, and let $S \subseteq V(\cQ^{n})$ be such that for all $x, x' \in S$ we have $\dist(x,x') \ge \gamma n$. 
Let $x, y \in V(\cQ^{n})\setminus B^k_{\cQ^n}(S)$ with $x \subseteq y$ and $m\coloneqq \dist(x,y) \geq n/8$. 
Then,  $|\mathcal{X}_{x,y}^{\neg B^k_{\cQ^n}(S)}| \ge 3m!/4$.  
\end{lemma}

\begin{proof}
We may assume that $x=\varnothing$ and $y=[m]$, where $m\geq n/8$.
Let $\cX_{x,y}^i$ denote the collection of chains $X \in \cX_{x,y}$ for which $V(X) \cap L_i \cap B^k_{\cQ^n}(S) \ne \varnothing$.
We have
\begin{equation}\label{eqm:imp}
|\cX_{x,y} \setminus \mathcal{X}_{x,y}^{\neg B^k_{\cQ^n}(S)}| \le \sum_{i=1}^{m-1}|\cX^i_{x,y}|.
\end{equation}
Furthermore, by \cref{prop: ballradiushit} (with $\gamma/2$ playing the role of $\gamma$), for all $\gamma m/4\le i \le  (1-\gamma/4)m$ we have that\COMMENT{There are at most $\binom{m}{i}2^{-\gamma m/400}$ vertices in $y^{(t)} \cap B^k_{\cQ^n}(S)$ by \cref{prop: ballradiushit}. 
Each of these can be part of $i!(m-i)!$ many $x$-$y$ chains.} 
\begin{equation}\label{eqn:chainhitting}
    |\cX^i_{x,y}| \le \binom{m}{i}2^{-\gamma m/400}i!(m-i)!=2^{-\gamma m/400}m!.
\end{equation}

Next, we consider the case $i \in [\gamma m/4]$, where first we prove the following claim.

\begin{claim}\label{claim:200}
For all $i \in [4k]$ we have $|\cX^i_{x,y}|\le (2n)^{i-1}(k+1)i!(m-i)!$.
\end{claim}

\begin{claimproof}
Observe that $|S \cap \bigcup_{i=1}^{\gamma n/2-1} L_i| \le 1$.
If $S \cap \bigcup_{i=1}^{\gamma n/2-1} L_i=\varnothing$, then $\cX^i_{x,y}=\varnothing$ for all $i \in [4k]$, so assume $|S \cap \bigcup_{i=1}^{\gamma n/2-1} L_i| = 1$.
Let $v$ be the unique vertex in $S \cap \bigcup_{i=1}^{\gamma n/2-1} L_i$. 
Then, $B^k_{\cQ^n}(v)\cap L_i = B^k_{\cQ^n}(S)\cap L_i$ for each $i \in [4k]$.
Thus, in order to prove \cref{claim:200}, it suffices to show that $|B^k_{\cQ^n}(v)\cap L_i| \le (2n)^{i-1}(k+1)$ for each $i \in [4k]$.

We will proceed by induction on $i$. 
Since $\varnothing=x\notin B^k_{\cQ^n}(v)$, it follows that $|B^k_{\cQ^n}(v) \cap L_1| \le k+1$\COMMENT{This is the value we get if $v \in L_{k+1}$.
If $v$ was in a lower level then $\varnothing=x\notin B^k_{\cQ^n}(v)$ and we get a contradiction.
If $v$ was in a higher level then $B^k_{\cQ^n}(v)\cap L_1 = \varnothing$.}, so the base case holds.

Now, suppose that $|B^k_{\cQ^n}(v) \cap L_{i-1}| \le (2n)^{i-2}(k+1)$ for some $2 \le i \le 4k$.
Consider first the case where $v \in L_j$ for some $i \le j \le i+k$.
In this case, any $u \in L_i \cap B^k_{\cQ^n}(v)$ satisfies either
\begin{enumerate}[label=$(\mathrm{\roman*})$]
    \item\label{item:treenewthing1} $u \subseteq v$ or
    \item\label{item:treenewthing2} there is a $v$-$u$ path of length at most $k$ whose penultimate vertex lies in $L_{i-1}$.\COMMENT{To see the second we first `lose' all the coordinate $1$'s from $v$ which $u$ does not have, taking us to a level below $L_i$.
    We then take the $1$'s we need to get back `up' to $u$.} 
\end{enumerate} 
There are $\binom{j}{i} \le \binom{k+i}{i}$ choices for $u$ satisfying \ref{item:treenewthing1},
whereas by applying induction to the penultimate vertex in such paths it follows that there are at most $n(2n)^{i-2}(k+1)$ choices for $u$ satisfying \ref{item:treenewthing2}.\COMMENT{By the induction hypothesis we can reach at most $(2n)^{i-2}(k+1)$ vertices in $L_{i-1}$ in $k-1$ steps or fewer from $v$.
Each of these has $n$ (at most) choices to go up to $L_i$ for the final vertex in these paths.}
Altogether, we have 
\[|B^k_{\cQ^n}(v)\cap L_i| \le \binom{k+i}{i} + n(2n)^{i-2}(k+1) \le (2n)^{i-1}(k+1).\]

The case where $v \in L_j$ for some $i-k \le j< i$ is handled similarly.
This completes the induction step and the proof of the claim.
\end{claimproof}

Recall that $|S\cap B^{\gamma m/2-1}_{\cQ^n}(x)| \le 1$.
It follows that for all $i \in [\gamma m/3]$ we have that $|B^k_{\cQ^n}(S)\cap L_i|  \le n^k$ and, therefore, $|\cX_{x,y}^i| \le n^ki!(m-i)!$.
Suppose $|S\cap B^{\gamma m/2-1}_{\cQ^n}(x)|=1$, and let $v$ be the unique vertex in $S \cap B^{\gamma m/2-1}_{\cQ^n}(x)$. 
Let $j\in[\gamma m/2-1]$ be such that $v \in L_j$.
It follows by \cref{claim:200} that\COMMENT{To get the final inequality note that $i!(m-i)!$ is decreasing from 1 up to $\gamma m/2$.
For the first term, note that it is increasing on the range $i=1$ up to $i=4k$.
To see this, note that as we pass from $i$ to $i+1$ we lose a power of $\le m$ but gain a power of $2n(i+1)$. 
Finally note that the first term in the last line, that is $(2n)^{2k-1}(k+1)(2k)!(m-2k)! > (2k)(n^k)(2k)!(m-2k)!$ (so every term in the final summation is greater than the maximum summand in the second term of the second last line). } 
\begin{align}\label{eqn:lastlast}
\sum_{i=1}^{\gamma m/3}|\cX^i_{x,y}| &\le \sum_{i=j-k}^{j+k}|\cX^i_{x,y}|\le \begin{cases}
  \sum_{i=j-k}^{j+k}(2n)^{i-1}(k+1)i!(m-i)! & \text{if } j\le3k,\\    
  \sum_{i=j-k}^{j+k}n^ki!(m-i)! & \text{if } 3k < j < \gamma m/2
  \end{cases}\nonumber \\
& \le \sum_{i=2k}^{4k}(2n)^{i-1}(k+1)i!(m-i)!.
\end{align}
If $S \cap B^{\gamma m/2}_{\cQ^n}(x) = \varnothing$, then this trivially holds too.
By the symmetry of the hypercube, we also have that 
\begin{equation}\label{eqn:201}
\sum_{i=m-\gamma m/3}^{m}|\cX^i_{x,y}| \le \sum_{i=2k}^{4k}(2n)^{i-1}(k+1)i!(m-i)!.
\end{equation}
Therefore, by \eqref{eqm:imp}--\eqref{eqn:201} we have\COMMENT{The terms in the second sum are of the form $cm!/n$ for some constant $c$ which depends on $k$ (and there are at most $2k$ of these terms).}
\[|\cX_{x,y} \setminus \mathcal{X}_{x,y}^{\neg B^k_{\cQ^n}(S)}| \le \sum_{i=\gamma m/3}^{m - \gamma m/3}2^{-\gamma m/400}m! + 2\sum_{i=2k}^{4k}(2n)^{i-1}(k+1)i!(m-i)! \le m!/4.\qedhere\]
\end{proof}

\begin{lemma}\label{lem: randsethit}
Let $0<1/n\ll\eta, 1/k',\gamma \leq 1$ and $n/2\leq k < n$  with $n, k', k \in \mathbb{N}$.
Let $\cA\subseteq V(\cQ^{n})$ be such that, for all $x\in V(\cQ^{n})$, we have that $|B_{\cQ^{n}}^{\gamma n}(x)\cap \cA|\leq 1$. 
Let $y \in L_k$ and let $s \coloneqq \lfloor (k+1)/2\rfloor$.
Then, there exist three sets of vertices $A=\{a_1, \ldots, a_{(1-\eta)n}\}\subseteq L_1$, $B=\{b_1, \ldots, b_{(1-\eta)n}\}\subseteq N_{Q^n}(y)$ and $C=\{c_1, \ldots, c_{(1-\eta)n}\}\subseteq L_s$ such that
\begin{enumerate}[label=$(\mathrm{\roman*})$]
    \item\label{lem: randset1hit} for each pair $i,j\in[(1-\eta)n]$ with $i\neq j$ we have $\dist(c_i,c_j)\geq9s^2/(10n)$, 
    \item \label{lem: randset3hit} $B^{k'}_{\cQ^n}(\cA)\cap C=\varnothing$, and
    \item\label{lem: randset2hit} for each $i \in [(1-\eta)n]$ we have $a_i \subseteq c_i \subseteq b_i$.
\end{enumerate}
\end{lemma}

\begin{proof}
Choose $k$ vertices $c_1,\dots,c_k\in y^{(s)}$ independently and uniformly at random\COMMENT{Repetitions are allowed.}.
Then, choose $n-k$ vertices $c'_{k+1}, \dots, c'_n\in y^{(s-1)}$ independently and uniformly at random.
For each $i\in[n]\setminus[k]$, choose an element $a_i\in[n]\setminus y$ such that all the $a_i$ are distinct, and let $c_i\coloneqq c_i'\cup\{a_i\}\in L_s$\COMMENT{Note that $|[n] \setminus y|=n-k$, so we use each element here exactly once. They can be assigned to each $c_i'$ arbitrarily.}.
For each $i\in[n]\setminus[k]$, let $b_i\in N^\uparrow(y)$ be the unique vertex  such that $a_i\in b_i$, so that when viewing each $a_i$ now as a 1-element set, we have $a_i\subseteq c_i\subseteq b_i$ for all $i\in[n]\setminus[k]$.

Note that, for each pair $i,j \in [n]$ with $i\neq j$, we have that\COMMENT{
Indeed, assume first that $i,j\in[k]$.
Think of $0$-$1$ coordinate vectors.
Each vector has $s$ $1$s, so for each single coordinate $\ell\in y$ we have that $\mathbb{P}[c_i(\ell)=1]=s/k$.
Furthermore, since the choices of $c_i$ and $c_j$ are independent, for each $\ell\in y$ we have that $\mathbb{P}[c_i(\ell)=c_j(\ell)=1]=s^2/k^2$.
Since $|y|=k$, we expect $c_i$ and $c_j$ to overlap in $s^2/k$ positions.\\
Assume now that $i,j\in[n]\setminus[k]$.
As above, we expect $c_i'$ and $c_j'$ to overlap in $(s-1)^2/k$ positions.
Since $c_i$ and $c_j$ are obtained by adding two distinct elements to each of $c_i$ and $c_j$, we conclude that $\mathbb{E}[|c_i\cap c_j|]=(s-1)^2/k$. \\
Finally, assume that $i\in k$ and $j\in[n]\setminus[k]$.
By the same argument, we have that $\mathbb{E}[|c_i\cap c_j'|]=s(s-1)/k$, so $\mathbb{E}[|c_i\cap c_j|]=s(s-1)/k$.}
\begin{equation}\label{equa:chain_inters_exp2}
    \mathbb{E}[|c_i\cap c_j|]\leq s^2/k.
\end{equation}
Assume that we reveal each $c_i$ in turn.
We then have that, for each $i\in[n]\setminus\{1\}$, the variables $|c_i\cap c_j|$ with $j\in[i-1]$ are hypergeometric.
Thus, by \cref{lem:ChernoffHyp} and \eqref{equa:chain_inters_exp2}, for each pair $i,j \in [n]$ with $i\neq j$ we have that\COMMENT{We have that
\begin{align*}
    \mathbb{P}[|c_i\cap c_j|\geq 21s^2/(20k)]\leq\mathbb{P}[||c_i\cap c_j|-\mathbb{E}[|c_i\cap c_j|]|\geq s^2/(20k)]\leq2e^{-s^4/(1200k^2s)}=2e^{-s^3/(1200k^2)}.
\end{align*}
Now note that $s^3/k^2\geq k/10$ and thus, since $k\geq n/2$, $s^3/k^2\geq n/20$.
Since $n$ is sufficiently large, we conclude that
\[\mathbb{P}[|c_i\cap c_j|\geq 21s^2/(20k)]\leq e^{-n/25000}.\]
}
\[\mathbb{P}[|c_i\cap c_j|\geq 21s^2/(20k)]\leq e^{-n/25000}.\]
By a union bound, it follows that a.a.s.~for all pairs $i, j \in [n]$ with $i\neq j$ we have $|c_i\cap c_j|<21s^2/(20k)$.
Therefore, since $\dist(c_i,c_j)=2(|c_i|-|c_i\cap c_j|)>2s(1-21s/(20k))$ and using the definition of $s$, we conclude that $\dist(c_i,c_j) \geq 9s/10$\COMMENT{We have that $\dist(c_i,c_j)=2(|c_i|-|c_i\cap c_j|)>2s(1-21s/(20k))$.
By the definition of $s$, we have that $s/k\leq100/198$ (for $n$ sufficiently large).
Thus, $\dist(c_i,c_j)>2s(1-21s/(20k))\geq2s(1-2100/(20\cdot198))\geq2s(1-11/20)=9s/10$.}.
In particular, \ref{lem: randset1hit} holds a.a.s.

Next, let $S_1\coloneqq y^{(s)}\cap B^{k'}_{\cQ^n}(\cA)$ and $S_2\coloneqq y^{(s-1)}\cap B^{k'+1}_{\cQ^n}(\cA)$.
By applying \cref{prop: ballradiushit} first with $k'$, $\cA$ and $s$ playing the roles of $k$, $S$ and $t$, respectively, and then with $k'+1$, $\cA$ and $s-1$ playing the roles of $k$, $S$ and $t$, respectively, we obtain that $|S_1| \le \binom{k}{s}2^{-\gamma n/200}$ and $|S_2| \le \binom{k}{s-1}2^{-\gamma n/200}$.\COMMENT{For $S_2$ we have applied \cref{prop: ballradiushit} with $k'+1$ playing the role of $k$ (because if $i \in [n]\setminus[k]$, and $c_i \in B^{k'}_{\cQ^n}(\cA)$ then $c_i' \in S_2$).}
Therefore, for all $i \in [k]$ we have 
\[\mathbb{P}[c_i \in S_1] \le \frac{\binom{k}{s}2^{-\gamma n/200}}{\binom{k}{s}}  = 2^{-\gamma n/200}.\]
For all  $i \in [n]\setminus [k]$ we have $\mathbb{P}[c_i \in S_1] \le \mathbb{P}[c_i' \in S_2]$\COMMENT{Since $z\in S_1$ implies $N^\downarrow(z) \subseteq S_2$.} and similarly we have $ \mathbb{P}[c_i' \in S_2] \le 2^{-\gamma n/200}$.
It now follows by a union bound that \ref{lem: randset3hit} holds a.a.s.\COMMENT{$\mathbb{P}[$\ref{lem: randset3hit} not holding$] \le n2^{-\gamma n/200}$.}

Next, consider an auxiliary bipartite graph $H$ with parts $y^{(1)}$ and $\{c_1, \dots, c_{k}\}$ and the following edge set.
For each $i \in [k]$ and $a \in y^{(1)}$, let $\{a,c_i\}$ be an edge whenever $a \in c_i$.
Thus, for each $i \in [k]$ we have that $d_H(c_i) = s$.
Furthermore, it follows by \cref{lem:Chernoff} that a.a.s.~for all $a \in y^{(1)}$ we have $d_H(a) = (1 \pm \eta/2)s$\COMMENT{We have that $\mathbb{E}[d_H(a)] = (s/k)k=s$.
Therefore, by \cref{lem:Chernoff} we have 
\[\mathbb{P}[d_H(a) \neq (1 \pm \eta/2)s]=\mathbb{P}[|d_H(a)- \mathbb{E}[d_H(a)]| \geq \eta\mathbb{E}[d_H(a)]/2] \leq2e^{-\eta^2s/12}\leq e^{-\eta^2 n/50}.\]
In the last calculation we use that $s$ is linear in $n$.
A union bound over all vertices gives the result.}.
Condition on this.
Then, for each $X \subseteq y^{(1)}$, since we have $e_H(N_H(X), y^{(1)}) \ge e_H(N_H(X), X)$, it follows that $|N(X)| \geq |X| - \eta k/2$\COMMENT{We have that $e_H(N_H(X), L_1\cap y) \ge e_H(X, N_H(X))$.
Therefore, by the bounds we have derived on the degrees, we have $|N(X)|s \ge |X|s(1 - \eta/2)$, so $|N(X)| \ge |X| - \eta|X|/2)$.
The claim follows since $|X| \le k$.}.
Therefore, by \cref{lema:Hall} we have a matching of size $(1-\eta/2)k$ in $H$.
 
Similarly, a.a.s.~we have a matching of size $(1- \eta/2)k$ in the analogous bipartite graph $H'$ with parts $N^\downarrow(y)$ and $\{c_1, \dots, c_{k}\}$, where for each $i \in [k]$ and $b \in N^\downarrow(y)$ we have that $\{b,c_i\}$ is an edge whenever $c_i\subseteq b$.\COMMENT{We set up the corresponding auxiliary graph $H'$ with parts $\{c_1, \dots, c_{k}\}$ and $N^\downarrow(y)$, where for each $i \in [k]$ and $b \in N^\downarrow(y)$ we add the edge $\{b,c_i\}$ if and only if $c_i\subseteq b$.
Note that $d_{H'}(c_i) = k-s = s\pm 1$ for all $i \in [k]$.
Furthermore, for $b\in N^\downarrow(y)$ we have that
\[\mathbb{E}[d_{H'}(b)] = k\frac{\binom{k-1}{s}}{\binom{k}{s}} = k(k-s)/k = k-s = s\pm 1.\]
Thus, $H'$ is a random graph with the same distribution as $H$ so the result follows. (The $s\pm 1$ can be absorbed into $\eta/3$, say.)}
By concatenating these matchings (and relabelling the indices if necessary), it follows that a.a.s.~there is an ordering $\{a_1, \dots, a_k\}$ of the elements of $y$ and an ordering $\{b_1, \dots, b_k\}$ of the vertices of $N^\downarrow(y)$ such that, for all $i\in [(1- \eta)k]$, we have $a_i \subseteq c_i \subseteq b_i$.
Furthermore, as explained before, by construction, for all $i\in[n]\setminus[k]$, we have $a_i \subseteq c_i \subseteq b_i$.
Thus, \ref{lem: randset2hit} holds a.a.s.

Finally, given that each of \ref{lem: randset1hit}, \ref{lem: randset3hit}, \ref{lem: randset2hit} holds a.a.s., there must exist a choice of $c_1,\ldots,c_n$ such that \ref{lem: randset1hit}-\ref{lem: randset2hit} hold simultaneously.
\end{proof}

We are now in a position to combine the results we have shown so far to prove the following key lemma, which is used to provide a base structure for the near-spanning tree which we seek. 

\begin{lemma}\label{lem: treeres2hit}
Let $0< 1/n \ll 1/C \ll \eps' \leq 1/2$, and $0 <1/n \ll 1/k',\gamma \leq 1/2$, where $n, k', C \in \mathbb{N}$, and let $(n, \mathbf{p}, M)$ be feasible with $0 < 1/n \ll 1/M$. 
Moreover, let $\cA\subseteq V(\cQ^{n})$ be such that, for all $x\in V(\cQ^{n})$, we have $|B^{\gamma n}_{\cQ^n}(x)\cap \cA|\leq 1$ and $\varnothing \notin B^{k'}_{\cQ^n}(\cA)$.
Then, with probability at least $1-e^{-50n}$ we have that $P \sim \cP^C(n, \mathbf{p}, M)$ satisfies the following:
for all $y \in \bigcup_{i=\lceil n/2 \rceil}^{\lfloor9n/10\rfloor}L_i\setminus B^{k'}_{\cQ^n}(\cA)$, there exists a collection of chains $\mathcal{X}_y$ such that, for all $X\in\mathcal{X}_y$, we have $X\subseteq P - B^{k'}_{\cQ^n}(\cA)$, one of the endpoints of $X$ belongs to $L_1$, and 
\[\Big|N_{\cQ^n}(y)\cap\bigcup_{X\in\mathcal{X}_y}V(X)\Big|\geq(1-\eps')n.\]
\end{lemma}

\begin{proof}
Fix $\eta>0$ such that $0 < 1/n \ll \eta \ll \eps'$, and let $m\coloneqq480000$.
Fix a vertex $y \in L_k\setminus B^{k'}_{\cQ^n}(\cA)$ for some $n/2\le k \le 9n/10$.
Let $s \coloneqq \lfloor (k+1)/2\rfloor$.
By \cref{lem: randsethit} with $\eta/2$ playing the role of $\eta$, there exists a collection of vertices $\{c_1, \dots, c_{(1-\eta/2)n}\} \subseteq  L_s$ such that $B^{k'}_{\cQ^n}(c_i)\cap \cA = \varnothing$ and $\dist(c_i,c_j)\geq9s^2/(10n)$ for all pairs $i,j \in [(1-\eta/2)n]$ with $i \ne j$; an ordering $b_1, \dots, b_n$ of $N_{\cQ^n}(y)$, and an ordering $a_1, \dots, a_n$ of $L_1$, such that for all $i\in [(1-\eta/2)n]$ we have $a_i\subseteq c_i \subseteq b_i$.
For each $i \in [(1-\eta/2)n]$, we call $(a_i, b_i, c_i)$ a \emph{triple}.
Note that $|B^{k'}_{\cQ^n}(\cA)\cap (L_1 \cup N_{\cQ^n}(y))| \le 2(k'+1)$,\COMMENT{There's at most 1 vertex $x \in \cA$ in $B^{\gamma n}_{\cQ^n}(\varnothing)$ and furthermore, this vertex is not in $B^{k'}_{\cQ^n}(\varnothing)$.
It follows that this vertex $x$ is at least distance $k'+1$ from $\varnothing$ and so $B^{k'}_{\cQ^n}(x)$ can have at most $k'+1$ elements in $L_1$. 
A similar argument work for $y$ instead of $\varnothing$.}
and hence we may assume for each $i \in [(1-\eta)n]$ that $(a_i, b_i, c_i)$ forms a triple where $a_i, b_i \notin B^{k'}_{\cQ^n}(\cA)$.
We denote by $\cT$ the collection of all such triples.
Partition $[(1-\eta)n]$ into two sets $\mathcal{I}_1\coloneqq\{i\in[(1-\eta)n]:b_i\in N^\downarrow(y)\}$ and $\mathcal{I}_2\coloneqq[(1-\eta)n]\setminus\mathcal{I}_1$.
Let $\mathcal{A}_1\coloneqq\{a_i:i\in\mathcal{I}_1\}$, $\mathcal{A}_2\coloneqq\{a_i:i\in\mathcal{I}_2\}$, $\mathcal{B}_1\coloneqq\{b_i:i\in\mathcal{I}_1\}$, $\mathcal{B}_2\coloneqq\{b_i:i\in\mathcal{I}_2\}$, $\mathcal{C}_1\coloneqq\{c_i:i\in\mathcal{I}_1\}$ and $\mathcal{C}_2\coloneqq\{c_i:i\in\mathcal{I}_2\}$.
Note that $k-\eta n \le |\cC_1| \le k$.\COMMENT{Given that $k \ge n/2$.}

We first turn our attention to $\mathcal{A}_1$, $\mathcal{B}_1$ and $\mathcal{C}_1$.
Partition $\mathcal{A}_1$, $\mathcal{B}_1$ and $\mathcal{C}_1$ into sets $\mathcal{A}^1, \dots, \mathcal{A}^{m}$, $\mathcal{B}^1, \dots, \mathcal{B}^{m}$ and $\mathcal{C}^1, \dots, \mathcal{C}^{m}$, respectively, each of size at least $\lfloor (k - \eta n)/m\rfloor$ and at most 2$\lfloor (k - \eta n)/m\rfloor$, and such that, for every triple $(a,b,c)\in \cT$ there exists $j\in[m]$ such that $a \in \mathcal{A}^j$, $b \in \mathcal{B}^j$ and $c \in \mathcal{C}^j$.\COMMENT{It will help with the analysis later to have such a partition, when considering the set of chains between vertices. 
We will throw away chains that overlap (that is, are contained in both $\mathcal{X}_{a^i_j, b^i_j}$ and $\mathcal{X}_{a^i_{j'}, b^i_{j'}}$, say).
By \cref{lem: chain count} we won't have to throw away too many chains, but without this further partition of the vertices we would still be throwing away too much. 
This fixes that problem.}
For each $i \in [m]$, write $\mathcal{A}^i=\{a^i_1,\ldots,a^i_{|\cA^i|}\}$, $\mathcal{B}^i=\{b^i_1,\ldots,b^i_{|\cA^i|}\}$ and $\mathcal{C}^i=\{c^i_1,\ldots,c^i_{|\cA^i|}\}$, where the labeling is such that $(a^i_j, b^i_j, c^i_j)\in\mathcal{T}$ for each $j\in[|\mathcal{A}^i|]$.
For each $i \in [m]$ and $j \in [|\mathcal{A}^i|]$, we define the set $\mathcal{Z}_{a^i_j, c^i_j}\subseteq\mathcal{X}_{a^i_j, c^i_j}$ as the set of all chains $X\in\mathcal{X}_{a^i_j, c^i_j}$ which, for all $j'\in[|\mathcal{A}^i|]\setminus\{j\}$, neither intersect any chain $X'\in\mathcal{X}_{a^i_{j'},c^i_{j'}}$, nor $B^{k'}_{\cQ^n}(\cA)$.
By \cref{lem: chain count,lem: pathball} and the definition of $m$, we have that\COMMENT{Note that by \cref{lem: chain count} we have that, for each $j'\in[|\mathcal{A}^i|]\setminus\{j\}$, at most a $60000/n$ proportion of $\mathcal{X}_{a^i_j,c^i_j}$ is not disjoint from chains in $\mathcal{X}_{a^i_{j'},c^i_{j'}}$.
We have at most $2\lfloor k(1-\eta)/m\rfloor\leq2k(1-\eta)/m\leq9n/(5m)\leq2n/m$ choices to range over $j'$, which gives at most a $120000/m$ proportion of $\mathcal{X}_{a^i_j, c^i_j}$ overlapping with chains in any other $X_{a^i_{j'},c^i_{j'}}$.
And $120000/m = 1/4$.
On the other hand, by \cref{lem: pathball} we also loose at most another $1/4$ proportion.}
\begin{equation}\label{eq: half size hit}
|\mathcal{Z}_{a^i_j, c^i_j}| \ge \frac{1}{2}|\mathcal{X}_{a^i_j, c^i_j}|.
\end{equation}
For each triple $(a,b,c)\in \cT$ and any graph $G\subseteq\cQ^n$, let $I_{a,c}(G)$ take value $1$ if $Y(\mathcal{Z}_{a,c},G)>0$, and $0$ otherwise.
(Recall that $Y(\mathcal{Z}_{a,c},G)$ denotes the number of chains $X \in \cZ_{a,c}$ with $X \subseteq G$.)
For each $i \in [m]$, let $I_i(G) \coloneqq \sum_{j \in [|\mathcal{A}^i|]}I_{a^i_j, c^i_j}(G) =  \sum_{j \in [|\mathcal{A}^i|]}I_{a^i_j, c^i_j}(G - B^{k'}_{\cQ^n}(\cA))$.

We are now in a position to consider $P\sim\cP^C(n,\mathbf{p},M)$.
Recall that $P$ is generated by sampling $C$ independent graphs $P_i$, where $P_i\sim\cP(n,\mathbf{p},M)$.
In each $P_i$ we can give bounds on the probability that certain chains appear. 
Note that, for each $i\in[C]$ and each fixed $i'\in[m]$ we have that, for every pair $j,j'\in[|\mathcal{A}^{i'}|]$ with $j\neq j'$, the variables $Y(\mathcal{Z}_{a^{i'}_j,c^{i'}_j},P_i)$ and $Y(\mathcal{Z}_{a^{i'}_{j'},c^{i'}_{j'}},P_i)$ are independent (and, therefore, $I_{a^{i'}_{j},c^{i'}_{j}}(P_i)$ and $I_{a^{i'}_{j'},c^{i'}_{j'}}(P_i)$ are independent too).
Since $C$ is a large constant, this independence will allow us to boost the probability that these chains appear in $P - B^{k'}_{\cQ^n}(\cA)$.
The analysis is broken into two steps.

\begin{claim}\label{claim:probboost}
With probability at least $1-2e^{-75n}$, the graph $P \sim P^{C}(n, \mathbf{p}, M)$ satisfies the following.
\begin{enumerate}[label=$(\arabic*)$]
    \item\label{itm:cp1} $P - B^{k'}_{\cQ^n}(\cA)$ contains an $a$-$c$ chain for at least $(1-\eps'/2)k$ of the triples $(a,b,c) \in \cT$ with $c \in \cC_1$.
    \item\label{itm:cp2} $P - B^{k'}_{\cQ^n}(\cA)$ contains a $c$-$b$ chain for at least $(1-\eps'/2)k$ of the triples $(a,b,c) \in \cT$ with $c \in \cC_1$.
\end{enumerate}
\end{claim}

\begin{claimproof}
We show that \ref{itm:cp1} and \ref{itm:cp2} each hold with probability $1-e^{-75n}$. 
The result then follows by a union bound.

For \ref{itm:cp1}, let $C' \coloneqq \sqrt{C}$.
By \eqref{eq: half size hit}, we can apply \cref{cor: boost} with $(\eps')^2$, $1/2$ and $C'$ playing the roles of $\eta$, $\alpha$ and $C$, respectively.
Thus, for $P'\sim\cP^{C'}(n,\mathbf{p},M)$, for all $i\in[m]$ and $j\in[|\mathcal{A}^i|]$ we have that
\begin{equation*}\label{equa:boosting1}
    \mathbb{P}[I_{a^i_j,c^i_j}(P')=1] = \mathbb{P}[Y(\mathcal{Z}_{a^i_j,c^i_j},P')>0] \ge 1 - (\eps')^2.
\end{equation*}
It follows that for all $i \in [m]$ we have $\mathbb{E}[I_i(P')]\geq(1-(\eps')^2)|\mathcal{A}^i|$ and, therefore, by \cref{lem:Chernoff},
\begin{equation}\label{eqn: ind rand}
\mathbb{P}[I_i(P') > |\mathcal{A}^i|(1-(\eps')^{3/2})] > 1 - e^{-(\eps')^3n/(25\cdot10^6)}.
\end{equation}

Let $P \sim \cP^{C}(n,\mathbf{p}, M)$, and note that $P$ can be generated by sampling $C'$ independent graphs $P_j'\sim \cP^{C'}(n,\mathbf{p}, M)$ and considering their union.
For each $i \in [m]$, let $\cE^i$ be the event that $I_i(P)>|\mathcal{A}^i|(1-(\eps')^{3/2})$.
It follows from \eqref{eqn: ind rand} that, for each $i\in[m]$, we have\COMMENT{For each $P_{j}'$, by \eqref{eqn: ind rand} we know that the probability that $I_i(P_j')\leq|\mathcal{A}^i|(1-(\eps')^{3/2})$ is at most $e^{-(\eps')^3n/(2\cdot10^7)}$.
Since the $P_j'$ are independent, it follows that the probability that $I_i(P)\leq|\mathcal{A}^i|(1-(\eps')^{3/2})$ is at most $(e^{-(\eps')^3n/(2\cdot10^7)})^{C'}\leq e^{-100n}$. 
The claim follows by taking the complement.} 
 $\mathbb{P}[\cE^i] > 1-e^{-100n}$.
Now let $\cE$ be the event that, for all $i\in[m]$, $\mathcal{E}_i$ holds. 
It follows by a union bound that
\[\mathbb{P}[\cE] \ge 1- e^{-75n}.\]
Thus, with probability at least $1-e^{-75n}$ the graph $P - B^{k'}_{\cQ^n}(\cA)$ contains an $a$-$c$ chain for at least $(1-(\eps')^{3/2})|\mathcal{C}_1|$ of the triples $(a,b,c) \in \cT$ with $c \in \cC_1$.
Since $|\mathcal{C}_1|\geq(1-2\eta)k$\COMMENT{$|\mathcal{C}_1|\geq k-\eta n\geq(1-2\eta)k$ since $n\leq 2k$.}, $P - B^{k'}_{\cQ^n}(\cA)$ contains an $a$-$c$ chain for at least $(1-\eps'/2)k$ of the triples $(a,b,c) \in \cT$ with $c \in \cC_1$.\COMMENT{Here we are using that $\eta \ll \eps'$.}

To show \ref{itm:cp2}, for each triple $(a,b,c) \in \cT$ with $c\in\mathcal{C}_1$, one can consider the set $\mathcal{X}_{c,b}$ and define sets $\mathcal{Z}_{c,b}$ and variables $I_{c,b}(G)$ analogously to the proof of \ref{itm:cp1}.
Then, by \cref{cor: boost}, \cref{lem: chain count} together with \cref{rmk:2}, and \cref{lem: pathball}\COMMENT{Note that we can apply these results since $k-s=k-\lfloor(k+1)/2\rfloor\geq k-(k+1)/2=(k-1)/2\geq n/4-1$, which is the condition we need.}, the same argument as above shows that, with probability at least $1-e^{-75n}$, the graph $P - B^{k'}_{\cQ^n}(\cA)$ contains a $c$-$b$ chain for at least $(1-\eps'/2)k$ of the triples $(a,b,c) \in \cT$ with $c \in \mathcal{C}_1$\COMMENT{This is a  `symmetric' proof given how we've set things up using the geometry of the hypercube.}.
\end{claimproof}

It follows by \cref{claim:probboost} that with probability at least $1-2e^{-75n}$ we have that $P - B^{k'}_{\cQ^n}(\cA)$ contains an $a$-$b$ chain for at least $(1-\eps')k$ of the triples $(a,b,c) \in \cT$ with $c\in\mathcal{C}_1$.
We can prove an analogous result for the sets $\mathcal{A}_2$, $\mathcal{B}_2$ and $\mathcal{C}_2$.
More specifically, we can show that with probability at least $1-2e^{-75n}$, for $P\sim P^C(n,\mathbf{p},M)$, the graph $P - B^{k'}_{\cQ^n}(\cA)$ contains an $a$-$b$ chain for at least $(1-\eps')(n-k)$ of the triples $(a,b,c)\in \cT$ with $c\in\mathcal{C}_2$.\COMMENT{The only thing that we used the first time was that $k$ was linear in $n$.
Since $k\le 9n/10$, the same result follows.}
Combining this with the previous, it follows that, with probability at least $1-4e^{-75n}$, $P - B^{k'}_{\cQ^n}(\cA)$ contains an $a$-$b$ chain for at least $(1-\eps')n$ of the triples $(a,b,c)\in \cT$.
Finally, the result follows by a union bound over all $y\in\bigcup_{i=\lceil n/2\rceil}^{\lfloor9n/10\rfloor}L_i\setminus B^{k'}_{\cQ^n}(\cA)$\COMMENT{The result holds with probability at least $1-4|\bigcup_{i=n/2}^{9n/10}L_i|e^{-75n}\geq1-2^{n+2}e^{-75n}\geq1-e^{-50n}$.}.
\end{proof}

Let $F$ be the union of all chains given by \cref{lem: treeres2hit} (applied with $k' \coloneqq k$).
Then, $F$ satisfies \ref{prop22hit} in \cref{lem: main treereshit} for all vertices $x \in \bigcup_{i=\lceil n/2\rceil}^{\lfloor9n/10\rfloor}L_i\setminus B^{k}_{\cQ^n}(\cA)$.
However, we need this property to hold for every $x \in V(\cQ^n) \setminus B^k_{\cQ^n}(\cA)$.
Recall the discussion in the beginning of this section where, due to the symmetries in the hypercube, we can `redefine' any vertex $v \in V(\cQ^n)$ to be the empty set $\varnothing$.
As discussed, this leads to a redefined notion of levels in the hypercube where, for each $i\in[n]_0$, we let $L_i(v)\coloneqq\{u\in V(\cQ^n):\dist(u,v)=i\}$.
The notion of a chain in this setting was also discussed.

When we consider this generalised setting, by replacing $L_i$ with $L_i(v)$ in \cref{def:biased,def: perc,def:percC}, we obtain a distribution on subgraphs of $\cQ^n$ which we denote by $\cP^C_v(n,\mathbf{p},M)$\index{PnpM3@$\cP^C_v(n,\mathbf{p},M)$}.
(Note, again, that there is a joint distribution of $\cP^C_v(n,\mathbf{p}, M)$ and $\cQ^n_{\min\{1,Cp\}}$ such that $\cP^C_v(n,\mathbf{p}, M) \subseteq Q^n_{\min\{1,Cp\}}$, where $p = \max_{i\in[n-1]_0}p_i$.)
Then, for any fixed $v\in V(\cQ^n)$, \cref{lem: treeres2hit} holds in this setting by replacing chains by chains with respect to $v$.
Intuitively, we may think of this simply as growing several branching processes rooted at different vertices of the hypercube.
This will be crucial in proving \ref{prop22hit}.

Note that $F$ may have unbounded degrees and also may be disconnected.
To turn $F$ into a bounded degree forest we will later delete suitable edges.
To make it connected without significantly raising any vertex degrees we will apply the following lemma.

\begin{lemma}\label{lem: hamconnect}
For $n \in \mathbb{N}$ such that $0 < 1/n \ll  \delta \leq 1/50$ and $0<\eps < 1/2$, the following holds a.a.s.
Let $R \sim \mathit{Res}(\cQ^n, \delta)$.
Then, there exists a cycle in $\cQ^n_\eps[(L_1 \cup L_2)\setminus R]$ which covers $L_1 \setminus R$. 
\end{lemma}

\begin{proof}
Let $R \sim \mathit{Res}(\cQ^n, \delta)$.
Let $\mathcal{A}$ be the event that $|R \cap L_1| \ge n/4$.
By \cref{lem:betaChernoff} we have that $\mathbb{P}[\mathcal{A}] \le e^{-\Theta(n)}$\COMMENT{For each $i\in[n]$, let $X_i$ be an indicator variable which takes value $1$ if $\{i\}\in R$, and $0$ otherwise. 
Let $X \coloneqq |R\cap L_1|= \sum_{i=1}^n X_i$.
We have $\mathbb{E}[X] = \delta n$.
Let $\beta \coloneqq \delta^{-1}/4$.
By \cref{lem:betaChernoff} we have that 
\[\mathbb{P}[X\geq n/4]=\mathbb{P}[X\geq\beta\delta n]\leq(e/\beta)^{\beta\delta n}=(4\delta e)^{n/4}\leq e^{-n/4}.\]}.
Expose $R\cap L_1$ and condition on the event that $\mathcal{A}$ does not occur.

Note that for each pair of vertices $x,y \in L_1$ there exists a unique vertex $z \in L_2 \cap N_{\cQ^n}(x) \cap N_{\cQ^n}(y)$ (in particular, $z=x\cup y$).
Let $H$ be an auxiliary graph with vertex set $L_1\setminus R$, where we include an edge between $x$ and $y$ if $x\cup y\notin R$ and $\{x,x\cup y\}, \{y,x\cup y\}\in E(\cQ^n_\eps)$.
By definition, a Hamilton cycle in $H$ would correspond uniquely to a cycle in $\cQ^n_\eps[(L_1 \cup L_2)\setminus R]$ covering $L_1 \setminus R$.
Note that $H$ has the same distribution as a binomial random graph $G \sim G_{n-|R\cap L_1|,p}$, where $p = (1-\delta)\eps^2$.
Let $\mathcal{B}$ be the event that there exists a Hamilton cycle in $H$.
As, after conditioning on $\cA$ not holding, $G_{n-|R\cap L_1|,p}$ is a.a.s.~Hamiltonian (see e.g.~\cite{Kor77,Posa76}), it follows that
\[\mathbb{P}[\mathcal{B}] \geq \mathbb{P}[\mathcal{B}\mid\overline{\mathcal{A}}]\mathbb{P}[\overline{\mathcal{A}}] \geq (1-o(1))(1- e^{-\Theta(n)}) = 1 - o(1).\qedhere\]
\end{proof}

We are now in a position to prove \cref{lem: main treereshit}.

\begin{proof}[Proof of \cref{lem: main treereshit}]
Choose constants $M, C \in \mathbb{N}$ such that $1/D,\delta \ll 1/C, 1/M \ll \eps'$.
By \cref{prop: feasible}, there exists a tuple $(n,\mathbf{p},M)$ which is feasible and such that $\max_{i\in[n-1]_0} p_i \le \eps/(5C)$.
Let $x_1\coloneqq\varnothing$, $x_2\coloneqq[\lceil n/2 \rceil]$, $x_3\coloneqq[n]\setminus x_2$ and $x_4\coloneqq[n]$.
For each $j \in [4]$, let $P_j \sim \cP_{x_j}^C(n, \mathbf{p}, M)$ be sampled independently, and let $R_j$ be the reservoir associated with $P_j$.
Let $R\coloneqq\bigcap_{j\in[4]}R_j$, and note that $R\sim\mathit{Res}(\cQ^n,1/10^{8C})$.
Finally, let $Q \sim \cQ^n_{\eps/5}$ be independent of all other previous choices.
Recall that, for each $j \in [4]$, there is a joint distribution of $\cP_{x_j}^C(n, \textbf{p}, M)$ and $\cQ^n_{\eps/5}$ such that $\cP_{x_j}^C(n, \textbf{p}, M)\subseteq\cQ^n_{\eps/5}$ (see the discussion after \cref{def:percC}). 
It follows that there is a joint distribution of $\bigtimes_{j=1}^{4}\cP_{x_j}^C(n, \mathbf{p}, M)\times\cQ^n_{\varepsilon/5}$ and $\cQ^n_{\varepsilon}$ such that $P_1\cup P_2\cup P_3\cup P_4\cup Q \subseteq \cQ^n_\eps$.\COMMENT{See the discussion after \cref{def:percC} which explains why $P_1 \subseteq Q^n_{\eps/5}$, say.}
Therefore, it suffices to show that we can find the desired tree $T$ in $(P_1\cup P_2\cup P_3\cup P_4\cup (Q - R))- B^k_{\cQ^n}(\cA)$.

For each $j \in [4]$, let $A_j \coloneqq \bigcup_{i=\lceil n/2\rceil}^{\lfloor9n/10\rfloor}L_i(x_j)\setminus B^k_{\cQ^n}(\cA)$, and let $\mathcal{E}_j$ be the event that, for all $y \in A_j$, the graph $P_j - B^k_{\cQ^n}(\cA)$ contains a collection $\mathcal{X}^j_y$ of chains with respect to $x_j$, where each chain $X\in\mathcal{X}^j_y$ has an endpoint in $L_1(x_j)$ (and thus in $L_1(x_j)\setminus (R_j\cup B^k_{\cQ^n}(\cA))$), and such that at least $(1-\eps')n$ of the neighbours of $y$ in $\cQ^n$ are covered by the union of the chains in $\mathcal{X}^j_y$.
Note that $\mathcal{E}_j$ is equivalent to saying that the union of the chains in $\cX^j_y$ satisfies \ref{prop22hit} for all $y \in A_j$.
For each $j \in [4]$ we have by \cref{lem: treeres2hit} that $\mathbb{P}[\mathcal{E}_j] \ge 1-e^{-50n}$.
Condition on the event that $\mathcal{E}_j$ holds for all $j\in[4]$.

For each $j\in[4]$, let $F_j\subseteq\cQ^n$ be given by $F_j\coloneqq\bigcup_{y\in A_j}\bigcup_{X\in\mathcal{X}_y^j}X$.
For each $j\in[4]$, let $G_j\subseteq F_j$ be defined by removing, for each $y\in V(F_j)\setminus\{x_j\}$, all edges of $F_j$ joining $y$ to its down-neighbours with respect to $x_j$ except for one (if $y$ has one such down-neighbour in $F_j$).\COMMENT{Note here that the notion of down-neighbour changes for each $j\in[4]$.
In general, the set of down-neighbours of $y$ with respect to a vertex $v$ is $N^\downarrow_v(y)\coloneqq\{x\in N(y):\dist(x,v)<\dist(y,v)\}$, as was already discussed at the beginning of the section.
This is only ever used here.}
In particular, it follows that each connected component of $G_j$ is a tree and contains one vertex in $L_1(x_j)$, and that $\Delta(G_j)\leq CM+1$.
Since $G_j$ has the same vertex set as $F_j$, we have that $G_j$ satisfies \ref{prop22hit} for all $y\in A_j$.
Furthermore, note that $V(\cQ^n)\setminus B^k_{\cQ^n}(\cA) = \bigcup_{j=1}^{4}A_j$.
Therefore, the graph $G\coloneqq\bigcup_{j\in[4]}G_j$ satisfies \ref{prop22hit} and $\Delta(G)\leq4CM+4$.

Since $B^{k+2}_{\cQ^n}(\cA) \cap \{ \varnothing, [n],[\lceil n/2\rceil], [n]\setminus [\lceil n/2\rceil]\}=\varnothing$ it follows that $B^{k}_{\cQ^n}(\cA) \cap (L_1(x_j) \cup L_2(x_j)) = \varnothing$ for each $j \in [4]$.
Let $\mathcal{E}_5$ be the event that, for each $j \in [4]$, $Q[L_1(x_j)\cup L_2(x_j)]-R$ contains a cycle $C_j$ which covers $L_1(x_j)\setminus R$.\COMMENT{
The edge sets are disjoint between the different cycles that we consider, so we don't need to worry about dependencies.}
By four applications of \cref{lem: hamconnect} (applied with $x_j$ playing the role of $\varnothing$) we have that $\mathbb{P}[\mathcal{E}_5] = 1 - o(1)$.
Condition on the event that this occurs.

Let $H \coloneqq G \cup \bigcup_{j \in [4]}C_j$.
It follows that $H$ is connected\COMMENT{For each $j\in[4]$ we have that $G_j\cup C_j$ is connected. 
Furthermore we have that the chains from $P_1$ and $P_4$ both have some intersection with the chains from $P_2$ as well as those from $P_3$, and therefore everything is connected.} and $\Delta(H) \le 4CM+6$. 
In order to complete the proof, let $T \subseteq H$ be a spanning tree of $H$. 
\end{proof}

\subsection{Extending the tree}\label{section:tree2}

Roughly speaking, in \cref{lem: main treereshit} we showed that, for any $\varepsilon>0$, given a reservoir chosen at random, the random graph $\cQ^n_\varepsilon$ a.a.s.~contains a bounded-degree tree $T'$ which avoids the reservoir and satisfies the local property that, for every vertex $x\in V(\cQ^n)$, all but a fixed small proportion of its neighbours are covered by $T'$.
Our goal in this section is to show that $T'$ can be extended into a tree $T$  where the proportion of uncovered vertices (in each neighbourhood) is even smaller, while still retaining the bounded degree property.
The precise statement is the following.

\begin{theorem}\label{thm: maintreeres}
For all $0<1/n\ll1/\ell,\eps\leq1$, where $n,\ell \in \mathbb{N}$, the following holds.
Let $R, W \subseteq V(\cQ^n)$ and let $T' \subseteq \cQ^n-(R\cup W)$ be a tree.
For each $x \in V(\cQ^n)\setminus W$, let $Z(x) \subseteq  N_{\cQ^n}(x)\cap V(T')$ be such that $|Z(x)|\geq3n/4$.\COMMENT{In application, these sets $Z(x)$ will allow us to extend the tree we get from \cref{lem: main treereshit} such that new leaves will be connected to tree vertices which are also contained in cubes later on.}
Then, a.a.s.~there exists a tree $T$ with $T' \subseteq T \subseteq (\cQ^n_\eps \cup T') - W$ such that
\begin{enumerate}[label=$(\mathrm{TC}\arabic*)$]
\item $\Delta(T)\leq\Delta(T')+1$;
\item for all $x \in V(\cQ^n)$, we have that $|B^\ell_{\cQ^n}(x) \setminus (V(T) \cup W)| \leq  n^{3/4}$\COMMENT{Alberto: Following the same proof exactly, this can be improved to a constant. We don't need it, but I don't see any reason not to write the stronger result, if it requires no extra effort.}, and
\item for each $x \in V(T)\cap R$, we have that $d_T(x)=1$ and the unique neighbour $x'$ of $x$ in $T$ is such that $x' \in Z(x)$.
\end{enumerate}
\end{theorem}

\begin{proof}
Let $Q \sim \cQ^n_\eps$.
For each $x \in V(\cQ^n)\setminus W$ we have $3\eps n/4 \le\mathbb{E}[e_{Q}(x, Z(x))] \le \eps n$.
Let $S_1\coloneqq\{x\in V(\cQ^n) : d_{Q}(x) > 11\eps n/10\}$, $S_2\coloneqq\{x\in V(\cQ^n)\setminus W : e_{Q}(x, Z(x)) < 2\eps n/3\}$ and $S\coloneqq S_1 \cup S_2$. 
Let $\mathcal{E}_1$ be the event that there exists no vertex $x \in V(\cQ^n)$ such that $|B^\ell_{\cQ^n}(x)\cap S_1|\geq n^{1/2}$.
By \cref{lem:vbadverticesdonotclump} we have that $\mathbb{P}[\mathcal{E}_1] \ge 1 - e^{-4n}$.\COMMENT{We let the sets $S(u) = N_{Q^n}(u)$ when applying \cref{lem:vbadverticesdonotclump}. 
We're also subbing in way worse bounds here---note that lemma gives constant number of bad vertices!}
Similarly, let $\mathcal{E}_2$ be the event that there exists no vertex $x \in V(\cQ^n)$ such that $|B^\ell_{\cQ^n}(x)\cap S_2|\geq n^{1/2}$.
By \cref{lem:vbadverticesdonotclump} we have that $\mathbb{P}[\mathcal{E}_2] \ge 1 - e^{-4n}$.\COMMENT{Note in the application of of \cref{lem:vbadverticesdonotclump} we need a set $S(v)$  for each vertex $v$.
We let the set $S(v) = Z(v)$ and for vertices in $W$, we let $Z(v)$ be any arbitrary set of size $3n/4$ in $N_{\cQ^n}(v)$.
It follows that there's at most $n^{1/2}$ bad vertices of this type, so even less when removing the vertices in $W$ (i.e. leaving us with the vertices in $S_2$).}
Condition on $\mathcal{E}_1 \wedge \mathcal{E}_2$ holding, that is, that there is no vertex $x \in V(\cQ^n)$ such that $|B^\ell_{Q}(x)\cap S|\geq n^{3/4}$.

Given $\mathcal{E}_1 \wedge \mathcal{E}_2$, let $H$ be an auxiliary bipartite graph with parts $A\coloneqq V(T')\setminus S$ and $B\coloneqq V(\cQ^n) \setminus (V(T')\cup W \cup S)$, where we include an edge between $a \in A$ and $b \in B$ whenever $\{a,b\} \in E(Q)$ and $a \in Z(b)$.
By definition of $S$ we have for all $a \in A$ that\COMMENT{Since $a$ is not in $S$ (and also $a \notin W $ since $V(T') \cap W = \varnothing$) it has degree at most $11\eps n/10$ in $Q$, of which at least $2\eps n/3$ of such neighbours are in $Z(a)\subseteq V(T') \subseteq V(\cQ^n)\setminus B$.} 
\[d_H(a) \le 11\eps n/10 -2\eps n/3 < \eps n/2.\]
Furthermore, we have for all $b \in B$ that
\[d_H(b) \ge 2\eps n/3 - n^{3/4} > \eps n/2.\]
Since for all $X \subseteq B$ we have $e_H(N_H(X), B) \ge e_H(X, N_H(X))$, it follows that $|N_H(X)| \ge |X|$\COMMENT{Using the previous two inequalities, we have that
\[|N_H(X)|\eps n/2\geq e_H(N_H(X), B) \ge e_H(X, N_H(X))\geq|X|\eps n/2.\]}.
Thus, by \cref{lema:Hall}, $H$ contains a matching covering all of $B$.
This corresponds to a matching in $\cQ \sim \cQ^n_\varepsilon$.
The statement follows by setting $T$ to be the union of $T'$ and this matching.
\end{proof}


\subsection{The repatching lemma}\label{section:tree3}

Later we will apply \cref{lem: main treereshit} to obtain a tree $T$ and a reservoir $R$ in $\cQ^n_\eps$ which is disjoint from $V(T)$. 
To carry out the absorption step later on, it will be important that for each vertex some proportion of its neighbourhood consists of vertices in $R$.\COMMENT{
This will help avoid the following problem. 
Suppose $x$ is to be absorbed via two cubes with entry vertices $y$ and $z$. 
If $y$ and $z$ are both in the tree it's possible that the algorithm will enter/leave these cubes via $y$ and $z$.
In such a case these cubes cannot be used to absorb $x$.
Having neighbours $y$ and $z$ which are not in $T$ will allow us to avoid this problem by using such vertices instead.}
However, the tree produced by \cref{lem: main treereshit} (and the subsequent application of \cref{thm: maintreeres}) will result in a small number of vertices with few or no neighbours in $R$.
The following repatching lemma will be called on to deal with such vertices, by slightly altering $T$.
\COMMENT{We will write the following lemma in the form where it is to be applied to a single vertex. 
In practice, we will apply it iteratively to a set of bad vertices.
These bad vertices will be well spread out and therefore we can incorporate them into the $F$ sets below.}
\COMMENT{The idea will be to apply this lemma to a bad molecule after already deciding which vertices will be absorbed together as pairs, to keep parities intact (the robust matching stuff is already done at this stage).
Given $x$ and $y$ from different layers of the bad molecule which are to be absorbed `as a pair' we know that their neighbourhoods have linear overlap in the projection of the cubes. 
The set $e$ corresponds to some edge $z_1z_2$ in this overlapping neighbourhood where $z_1$ and $z_2$ could absorb $x$ and $y$ (in their respective layers).
The set $B(e)$ corresponds to the neighbours of $z_1z_2$ in the tree. 
These need to be reconnected in the tree if we are to remove $z_1$ and $z_2$ and place them in the reservoir.
Finally the sets $F$ correspond to vertices used in previous patchings when we iterate this lemma.
If we remove some vertices and place them into the reservoir, we don't want to incorporate them back into the tree. 
The lack of clumping later on will and the fact that these connecting sets which we find here are `small' will keep things ok here.}

Given a graph $P$ and $S \subseteq V(P)$ we say that \emph{$S$ is connected in $P$} if the vertices of $S$ lie in the same component of $P$.

\begin{lemma}\label{lem:repatch}
Let $0 < 1/n \ll c,\eps, 1/f, 1/D$, where $f, D \in \mathbb{N}$.
Given a fixed $x \in V(\cQ^n)$, let $C(x) \subseteq N_{\cQ^n}(x) \times N_{\cQ^n}(x)$ be such that $|C(x)| \geq cn$\COMMENT{Just a subset of matched edges in the neighbourhood of $x$.} and such that, for all distinct $(y_1,z_1), (y_2, z_2) \in C(x)$, we have  $\{y_1,z_1\} \cap \{y_2, z_2\} = \varnothing$.
Furthermore, for each $(y,z) \in C(x)$, let $B(y,z) \subseteq (N_{\cQ^n}(y) \cup N_{\cQ^n}(z)) \setminus\{x\}$ with $|B(y,z)| < D$.\COMMENT{Just the tree neighbours of the vertices of $y$ and $z$. 
In the applications we will be able to assume that $y$ and $z$ do not have $x$ as a neighbour in the tree, because $x$ will have at most $M$ neighbours in the tree so we can discard a constant number of $(y,z)$ pairs which cover these cases and still then have a linear number of $(y,z)$ pairs of this form.}
Then, with probability at least $1-e^{-5n}$, for every $F \subseteq V(\cQ^n)$ with $|F|\le f$, there exist a pair $(y,z) \in C(x)$ with $y,z\notin F$ and a graph $P \subseteq \cQ^n_{\eps} - \{y,z\}$ with $|V(P)| < 5 D$\COMMENT{The reason for 5 here is as follows. Notice the paths we construct by \ref{cl3} are length 4, so they have at most 5 vertices. We create at most $D$ of these paths, so 5D is an upper bound here.} such that
\begin{enumerate}[label=$(\mathrm{R}\arabic*)$]
    \item \label{lab2} $B(y,z) \cap N_{\cQ^n}(y)$ is connected in $P$, and so is $B(y,z) \cap N_{\cQ^n}(z)$.
    \item \label{lab4} $V(P) \cap F = \varnothing$.
\end{enumerate}
\end{lemma}

\begin{proof}
We provide a counting argument to show there exist edge-disjoint graphs $P_1, \dots, P_{\eps' n} \subseteq \cQ^n$ such that, if any is present in $\cQ^n_{\eps}$, then it would satisfy \ref{lab2} and \ref{lab4} for some $(y,z) \in C(x)$. 
We will then prove that, with high probability, one of the $P_i$ must be present in $\cQ^n_\eps$.
Note that we may assume $x=\varnothing$.
By passing to a subset of $C(x)$ and replacing $c$ with $c/(30D)$ if necessary, we may also assume that $|C(x)|=cn$ and $2Dc < 1/10$\COMMENT{otherwise, we may take some $C'(x) \subseteq C(x)$ with $|C'(x)| = c'n$ such that $2Dc' < 1/10$ and prove the statement for this set, which is sufficient for showing the result.}.
Similarly, by passing to a suitable subset of $C(x)$, we may assume that, for any distinct $(y,z), (y',z') \in C(x)$, we have that $B(y,z) \cap B(y',z') = \varnothing$\COMMENT{
Note that, by considering $x=\varnothing$, the vertices in $B(y,z)$ each have size $2$.
This means they each have at most one neighbour of size $1$ other than $y$ or $z$, and each such neighbour could rule out at most one other $(y',z')$ pair from our considerations.
Since $|B(y,z)| <D$, we remove at most a $1-1/D$ fraction of the pairs $(y',z')$ to get this vertex-disjoint property, that is, we can guarantee that we keep at least a $1/D$ proportion of the pairs.}.

Fix any $F \subseteq V(\cQ^n)$ with $|F|\le f$.
We update $C(x)$ by removing any pair $(y,z)\in C(x)$ for which $(\{y,z\}\cup B(y,z))\cap F\neq\varnothing$.
It follows that $|C(x)| \geq cn - 2f$.\COMMENT{Since $B(y,z) \cap B(y',z') = \varnothing$, a forbidden vertex from $F$ can appear in at most $1$ set $B(y,z)$, as well as in at most 1 set $\{y',z'\}$, and we have that $|F| \leq f$.}
Now, for each $(y,z) \in C(x)$ and for each $w\in\{y,z\}$, let $A_w \coloneqq N_{\cQ^n}(w) \cap B(y,z)$,\COMMENT{
Note that as we range over $(y,z)$ we have that we will never define another $A_y$ by definition of the $(y,z)$ pairs being vertex disjoint.
We could have that the case of $(y,y)$ in our definition but this doesn't matter here.}
and let $x^w_1, \dots, x^w_{|A_{w}|}$ be the vertices of $A_{w}$.

\begin{claim}\label{claim: counting paths2}
For each $e = (y,z) \in C(x)$,  $w \in \{y,z\}$ and $i \in [|A_w|-1]$, there exists a collection $\mathcal{P}^{w}_i$ of subgraphs of $\cQ^n$ such that the following hold:
\begin{enumerate}[label=$(\mathrm{RC}\arabic*)$]
    \item\label{cl1} $|\mathcal{P}^{w}_i| \ge n/2$ and for each $P \in \mathcal{P}^{w}_i$ we have $V(P) \cap (F\cup \{y,z\}) = \varnothing$.
    \item  \label{cl3} Every $P \in \mathcal{P}^{w}_i$ is an $(x_i^w,x_{i+1}^w)$-path of length $4$. 
    \item \label{cl6} The graphs in $\mathcal{P}^w_i$ are pairwise edge-disjoint.
     \item\label{cl2} For every $e' =(y',z') \in C(x)$ with $e' \neq e$, every $w' \in \{y',z'\}$ and every $j \in [|A_{w'}|-1]$, the graphs in $\mathcal{P}^{w}_i$ are edge-disjoint from those in $\mathcal{P}^{w'}_j$.
\end{enumerate}
\end{claim}

(Note that we do not require the paths in $\cP_i^w$ to be edge-disjoint from those in $\cP_{i'}^{w'}$ when $w, w' \in \{y,z\}$ are distinct and $i \in [|A_w|-1]$, $i' \in [|A_{w'}|-1]$.)

\begin{claimproof}[Proof of \cref{claim: counting paths2}]
Let $e_1, \dots, e_{cn}$ be an ordering of the elements of $C(x)$, where for each $k \in [cn]$ we have that $e_k = (y_k, z_k)$.
Note that, for each $k \in [cn]$, each $w\in\{y_k,z_k\}$ and all $i\in[|A_{w}|]$, we have that $|x^w_i|=2$, and for each $i,j \in [|A_{w}|]$ with $i\neq j$ we have that $\dist(x^w_i, x^w_j) = 2$, with $x^w_i\cap x^w_j=w$. 

Suppose that, for some $1< k \le cn$, every $j\in[k-1]$, every $w\in\{y_j,z_j\}$ and every $i\in[|A_w|-1]$, we have found a collection $\mathcal{P}^{w}_i$ which satisfies \ref{cl1}--\ref{cl2}.
We now show that, for each $w \in \{y_k, z_k\}$ and each $i \in [|A_w| -1]$, a suitable choice for $\mathcal{P}^{w}_i$ exists.
We construct the set $\cP^w_i$ as follows.
Let $v_1\coloneqq x^w_i\setminus w$ and $v_2\coloneqq x^w_{i+1}\setminus w$.
For each $d\in[n]\setminus(x^w_i\cup x^w_{i+1})$, let $P_d \subseteq \cQ^n$ be the path which passes through the following vertices in successive order:
\[x^w_i, x^w_i\cup\{d\}, x^w_i\cup\{d\}\cup v_2, (x^w_i\cup\{d\}\cup v_2)\setminus v_1=x^{w}_{i+1}\cup\{d\}, x^w_{i+1}.\]
Note that each path $P_d$ has length $4$ and that $V(P_d)\cap \{y_k,z_k\} = \varnothing$.
Furthermore, for any distinct $d,d'\in[n]\setminus(x^w_i\cup x^w_{i+1})$, it is clear that $P_d$ and $P_{d'}$ are internally disjoint, and hence, are edge-disjoint.
To avoid $F$ as well as the edges of any previously chosen paths we set 
\[\mathcal{P}_i^w\coloneqq\left\{P_d: d\in[n]\setminus(x^w_i\cup x^w_{i+1}); x^{w}_i\cup\{d\},x^{w}_{i+1}\cup\{d\}\notin N\!\left(\bigcup_{j=1}^{k-1}B(y_j, z_j)\right)\!\!; V(P_d)\cap F=\varnothing\right\}.\]

It follows that $\cP^w_i$ satisfies \ref{cl3} and \ref{cl6}.
Recall that $V(P_d)\cap \{y_k,z_k\} = \varnothing$.
Therefore, to see that \ref{cl1} holds, note that, for all distinct $(y',z'), (y'',z'') \in C(x)$ and all $x' \in B(y',z'), x'' \in B(y'',z'')$, since $B(y',z')\cap B(y'',z'')=\varnothing$, we have that $x'$ and $x''$ contain at most one common neighbour in the third level $L_3$ of $\cQ^n$.\COMMENT{We are only interested in intersections with vertices of size $3$ because the paths we consider are of the form where we always `go up', that is, the second vertex in each path has size $3$.}
Since $|\bigcup_{j=1}^{k-1}B(y_j,z_j)|<Dcn$, there are at most $2Dcn < n/10$ choices for $d$ such that $x^{w}_i\cup\{d\} \in N(\bigcup_{j=1}^{k-1}B(y_j, z_j))$, or $x^{w}_{i+1}\cup\{d\} \in N(\bigcup_{j=1}^{k-1}B(y_j, z_j))$.\COMMENT{
Again, think of $x$ as the empty set. 
There are at most $cn$ previous iterations, and each $B(y,z)$ has size at most $D$ in each of these iterations. 
Now a vertex in $B(y,z)$ can have at most one neighbour as a $x^{w}_j \cup\{d\}$ (say), as two vertices in the same level can't share two up-neighbours. 
The same is true for the other endpoint of the path.
Hence, we remove at most $2Dcn$ paths from consideration.}
Furthermore, since $|F|\leq f$, it follows that there are still at least $n/2$ suitable choices for $d$, that is, \ref{cl1} holds as desired.
Additionally, \ref{cl2} holds by construction; indeed, since neither the second nor the fourth vertex of each path in $\cP_i^w$ lies in some path in $\bigcup_{j\in[k-1]}\bigcup_{w'\in\{y_j,z_j\}}\bigcup_{i'\in[|A_{w'}|-1]}\mathcal{P}_{i'}^{w'}$, the paths in $\mathcal{P}_i^w$ must be edge-disjoint from all the paths in $\bigcup_{j\in[k-1]}\bigcup_{w'\in\{y_j,z_j\}}\bigcup_{i'\in[|A_{w'}|-1]}\mathcal{P}_{i'}^{w'}$. 
Thus, we can proceed by induction and create a suitable collection $\mathcal{P}^w_i$ for each $k\in[cn]$, $w \in \{y_k, z_k\}$ and $i \in [|A_w|-1]$.
\end{claimproof}

For each $e = (y,z) \in C(x)$, $w \in \{y,z\}$ and $i \in [|A_w|-1]$, let $\mathcal{P}^w_i$ be the collection of subgraphs given by \cref{claim: counting paths2}.
Note that, for any choice of $P_1 \in \mathcal{P}^w_1, \dots, P_{|A_w|-1} \in \mathcal{P}^w_{|A_w|-1}$, we have that $A_w$ is connected in $P_w \coloneqq \bigcup_{j=1}^{|A_w|-1} P_j$.
To complete the proof, we now show that, on passing to $\cQ^n_\eps$, with high probability there will exist some $e = (y,z) \in C(x)$ and some $P_y$ and $P_z$ of the above form such that $P_y \cup P_z \subseteq \cQ^n_\eps$.
Moreover, note that each such choice of $P_y \cup P_z$ satisfies \ref{lab2} and \ref{lab4} for our fixed $F$ and $|P_y \cup P_z| \le 5D$.
Since $P_y \cup P_z\subseteq B^{4}_{\cQ^n}(x)$, \cref{lem:repatch} will then follow by a union bound over all choices of $F \subseteq B^{4}_{\cQ^n}(x)$ with $|F| \le f$.

Let $Q \sim \cQ^n_\eps$.
Consider $e = (y,z) \in C(x), w \in \{y,z\}$ and $i \in [|A_w|-1]$.
Let $P \in \mathcal{P}^w_i$ and recall that $P$ has length $4$.
It follows that $\mathbb{P}[P \nsubseteq Q] = 1-\eps^4$. 
Let $\mathcal{E}^w_i$ be the event that there exists some $P \in \mathcal{P}^w_i$ such that $P \subseteq Q$.
Since $|\mathcal{P}^w_i| \ge n/2$ and paths in $\mathcal{P}^w_i$ are edge-disjoint by \ref{cl6}, we have that $\mathbb{P}[\mathcal{E}^w_i] \ge 1 - (1-\eps^4)^{n/2}$.
Let $\mathcal{E}_{e} \coloneqq \bigwedge_{w \in \{y,z\}}\bigwedge_{i \in [|A_w|-1]}\mathcal{E}^w_i$.
Since $|A_y|+|A_z|\leq 2D$, we have that\COMMENT{The first inequality follows by a union bound. 
For the second, since $1-x\leq e^{-x}$, we have that $1 - 2D(1-\eps^4)^{n/2} \geq 1 - 2De^{-\eps^4n/2} \geq 1 -e^{-\eps^4 n/4}$, where the last inequality holds for $n$ sufficiently large.}
\[\mathbb{P}[\mathcal{E}_{e}] \ge 1 - 2D(1-\eps^4)^{n/2} \geq 1 -e^{-\eps^4 n/4}.\]
Finally, let $\mathcal{E}$ be the event that there exists some $e \in C(x)$ such that the event $\mathcal{E}_{e}$ occurs.
It follows by \ref{cl2} that, for $e, e' \in C(x)$ with $e \ne e'$, the event $\mathcal{E}_{e}$ is independent of $\mathcal{E}_{e'}$.
Therefore, since $|C(x)| \geq cn$, we have that\COMMENT{We have
\[\mathbb{P}[\mathcal{E}] \geq 1 -(e^{-\eps^4 n/4})^{cn} = 1-e^{-\eps^4cn^2/4}.\]}
\[\mathbb{P}[\mathcal{E}] \geq 1 - e^{-\eps^4 c n^2/4}.\]

Recall that by \ref{cl3} it now suffices to consider a union bound over all choices of $F \subseteq B^{4}_{\cQ^n}(x)$ with $|F|\leq f$.\COMMENT{That is, by restricting any choice of $F$ to $B^{4}_{\cQ^n}(x)$, we have at most this many cases to consider.}
The result follows since 
\[1 - f\binom{n^{4}}{f} e^{-\eps^4 c n^2/4}> 1-e^{-5n}.\qedhere\]
\end{proof}


\section{Hamilton cycles in randomly perturbed dense subgraphs of the hypercube}\label{section8}

In this section, we introduce a few more auxiliary lemmas and combine them with the tools we have developed so far to prove the following result.

\begin{theorem}\label{thm:main1}
For every $\eps, \alpha\in(0,1]$ and $c>0$, there exists $\Phi\in \mathbb{N}$ such that the following holds.
Let $H \subseteq \cQ^{n}$ be a spanning subgraph with $\delta(H)\geq\alpha n$ and let $G\sim\cQ^n_\varepsilon$.
Then, a.a.s.~there is a subgraph $G'\subseteq G$ with $\Delta(G')\leq\Phi$ such that, for every $F\subseteq \cQ^{n}$ with $\Delta(F)\leq c\Phi$, the graph $((H\cup G)\setminus F)\cup G'$ is Hamiltonian.
\end{theorem}

Note that \cref{thm:main1} trivially implies the case $k=1$ of \cref{thm:main}.
In fact, in \cref{sect:thm1} we will use \cref{thm:main1} to prove \cref{thm:main} in full generality.
For this derivation, we will need the stronger conditions imposed in the statement of \cref{thm:main1}.
More precisely, the formulation of \cref{thm:main1} involving a `forbidden' graph $F$ and a `protected' graph $G'$ is designed to make repeated applications of \cref{thm:main1} possible in order to take out $k$ edge-disjoint Hamilton cycles.
When finding the $i$-th Hamilton cycle, the protected graph will contain all the essential ingredients for this, while the forbidden graph will contain all previously chosen Hamilton cycles as well as the protected graphs for the entire set of Hamilton cycles (see \cref{sect:thm1} for details).

The first step of the proof of \cref{thm:main1} will be to consider a particular partition of the hypercube into subcubes.
The structure of this partition will be used extensively throughout the rest of the paper, so we first introduce the necessary notation in the next subsection. 
Then, in \cref{sect:bondless} we prove several results regarding this structure, concerning its properties in $\cQ^n_\varepsilon$ and with respect to a reservoir $R\sim Res(\cQ^n, \delta)$.
In \cref{connect}, we will prove our \emph{connecting lemmas}, which provide sets of paths in (sub)cubes which (roughly speaking) link up pairs of vertices and, together, span all vertices of these (sub)cubes.
We prove \cref{thm:main1} in \cref{sect:mainpf}.
Finally, we deduce \cref{thm:thresholdk,thm:almost,thm:main} from \cref{thm:main1} in \cref{sect:thm1}.

\subsection{Layers, molecules, atoms and absorbing structures}\label{sect8notation}

Throughout this section, given any two vectors $u$ and $v$, we will write $uv$ for their concatenation.
Consider $\cQ^n$ and some $s\in\mathbb{N}$, with $2\leq s<n$.
We divide $\cQ^n$ into $2^s$ vertex-disjoint copies of $\cQ^{n-s}$ as follows:
for each $a\in\{0,1\}^s$, we consider the set of vertices $V_a\coloneqq\{av: v\in\{0,1\}^{n-s}\}$, and consider the graph $\cQ(a)\coloneqq\cQ^n[V_a]$.\COMMENT{Alternatively, we can think of this partition of the hypercube by seeing it as a `blow-up' of $\cQ^s$ in which each vertex is replaced by a copy of $\cQ^{n-s}$ and each edge $e$ is replaced by a perfect matching of size $2^{n-s}$ in the direction given by $e$.}
We will refer to each $\cQ(a)$ as an $s$-\emph{layer} of $\cQ^n$ ($s$ will be dropped whenever clear from the context).
Given $0\leq\ell\leq n-s$, we will refer to any copy of a cube $\cQ^\ell$ in one of the $s$-layers as an $\ell$-\emph{atom} (again, $\ell$ will be dropped whenever clear from the context).

Fix a Hamilton cycle $\mathcal{C}$ of $\cQ^s$.
By abusing notation, whenever necessary, we assume that the coordinate vector of each vertex of $\mathcal{C}$ is concatenated with $n-s$ $0$'s.
$\cC$ induces a cyclical ordering on $\{0,1\}^s$, which we will label as $a_1,\ldots,a_{2^s}$.
In turn, this gives a cyclical ordering on the set of layers. 
In this section, for each $i\in[2^s]$, we denote $L_i\coloneqq\cQ(a_i)$\index{Lilayer@$L_i$ (layer)} (as opposed to \cref{section:tree}, where $L_i$ denoted the $i$-th level of the hypercube).
Given an $\ell$-atom $\mathcal{A}$ in an $s$-layer $\cQ(a)$, we refer to $\mathcal{M}(\mathcal{A})\coloneqq\mathcal{A}+V(\mathcal{C})$ as an $(s,\ell)$-\emph{molecule} (again, the parameters will be dropped when clear from the context).
Thus $\cM(\cA)$ is the vertex-disjoint union of $2^s$ copies of $\cQ^\ell$.
We refer to an $(s,0)$-molecule as a \emph{vertex molecule} and an $(s,\ell)$-molecule for $\ell\geq1$ as a \emph{cube molecule}.
Observe that, if we label the atoms in a molecule cyclically following the labelling of the layers, then $\cQ^n$ contains a perfect matching between any two consecutive atoms where all edges are in the same direction as the corresponding edge in $\mathcal{C}$\COMMENT{(and no other edges exist between the two atoms)}.
Whenever we work with molecules, we consider this cyclical order implicitly.
In particular, whenever we refer to a molecule $\cM = \cM(\cA) = \cA_1\cup \dots \cup \cA_{2^s}$, the cyclical order $\cA_1\cup \dots \cup \cA_{2^s}$ of the $\cA_i$ is that induced by $\cC$. 
Given a molecule $\cM(\cA) = \cA_1\cup \dots \cup \cA_{2^s}$, a \emph{slice} $\cM^* \subseteq \cM(\cA)$ will consist of the subgraph of $\cM(\cA)$ induced by its intersection with some number of consecutive layers, i.e.~$\cM^* = \cA_{a+1}\cup \dots \cup \cA_{a+t}$ for some $a,t \in [2^s]$.
Alternatively, given any $a\in V(\mathcal{C})$, any path $P\subseteq\mathcal{C}$ and any atom $\mathcal{A}\subseteq\cQ(a)$, $P$ determines a slice of $\mathcal{M}(\mathcal{A})$ by setting $\mathcal{M}^*\coloneqq\mathcal{A}+V(P)$.

Consider $i\in[2^s]$ and the cyclical ordering of the layers given by $\mathcal{C}$.
Given any subgraph $G \subseteq \cQ^n$, we will often denote the restriction of $G$ to the $i$-th layer by $L_i(G)$\index{LG@$L(G)$}, that is, $L_i(G) \coloneqq G[V(L_i)]$.
Given any $v \in \{0,1\}^{n-s}$, we will refer to the vertex $a_iv$ as the $i$-th \emph{clone} of $v$.
In general, when it is clear from the context, we will also refer to the $i$-th clone of a cube $C \subseteq \cQ^{n-s}$ (as well as other subgraphs), which, analogously, will be the corresponding copy in $L_i$ of $C$.
In particular, the $i$-th layer $L_i$ is the $i$-th clone of $\cQ^{n-s}$.

As we already discussed in \cref{sect:outline main}, in order to prove our results we will first construct a near-spanning cycle and then absorb the remaining vertices into this cycle.
We will achieve this by using the following absorbing structure.

\begin{definition}[Absorbing $\ell$-cube pair]\label{def:abs}
Let $\ell,n\in\mathbb{N}$, and let $G\subseteq \cQ^n$.
Given a vertex $x \in V(\cQ^n)$, an \emph{absorbing $\ell$-cube pair for $x$} in $G$, which we denote by $(C^l,C^r)$\index{(Cl@$(C^l,C^r)$}, is a subgraph of $G$ which consists of two vertex-disjoint $\ell$-dimensional cubes $C^l,C^r\subseteq G$ and three edges $e, e^l, e^r\in E(G)$ satisfying the following properties:
\begin{enumerate}[label=$(\mathrm{AP}\arabic*)$]
    \item\label{ACPprop1} $|V(C^l) \cap N_{\cQ^n}(x)| = |V(C^r) \cap N_{\cQ^n}(x)| = 1$;
    \item\label{ACPprop4} $e^l$ and $e^r$ are the unique edges from $x$ to $C^l$ and $C^r$, respectively\COMMENT{Note that the fact that there is a unique edge follows from \cref{rem:RNcodeg}.};
    \item\label{ACPprop2} the unique vertex $y\in V(C^l) \cap N_{\cQ^n}(x)$ satisfies $\dist(y, C^r) = 1$, and
    \item\label{ACPprop3} $e$ is the unique edge from $y$ to $C^r$.
\end{enumerate}

We will refer to $C^l$ as the \emph{left absorption cube} and to $C^r$ as the \emph{right absorption cube}.
Given an absorbing $\ell$-cube pair $(C^l, C^r)$ we refer to $y$ as the \emph{left absorber tip}, and to the unique vertex $z \in V(C^r) \cap N_{\cQ^n}(x)$ as the \emph{right absorber tip}.
We refer to the unique vertex $z' \in e\setminus\{y\}$ as the \emph{third absorber vertex}.
\end{definition}

We remark here that, given an absorbing $\ell$-cube pair for $x$, the third absorber vertex $z'$ and the right absorber tip $z$ are adjacent.
Moreover, together with $x$ and the left absorber tip $y$, they form a cycle of length $4$.
See \cref{fig:abs-struct1} for a representation of the absorbing $\ell$-cube pairs.

\subsection{Bondless and bondlessly surrounded molecules}\label{sect:bondless}

Given any graph $G\subseteq\cQ^n$, we will say that an $(s,\ell)$-molecule $\mathcal{M}=\mathcal{A}_1\cup\dots\cup\mathcal{A}_{2^s}\subseteq\cQ^n$, where $\mathcal{A}_i$ is the $i$-th clone of some $\ell$-cube $\mathcal{A}\subseteq \cQ^{n-s}$, is \emph{bonded} in $G$ if, for all $i\in[2^s]$\COMMENT{(recall that we consider some fixed cyclical ordering between the different layers of $G$)}, $G$ contains at least $100$ edges  between $\mathcal{A}_i$ and $\mathcal{A}_{i+1}$ whose endpoint in $\cA_i$ has even parity and at least 100 such edges whose endpoint in $\cA_i$ has odd parity. \COMMENT{These can only be matching edges.}
Otherwise, we call it \emph{bondless} in $G$.
Furthermore, given a collection $\cU$ of $(s,\ell)$-molecules in $G$, we say that $\mathcal{M} \in \cU$ is \emph{bondlessly surrounded} in $G$ (with respect to $\mathcal{U}$) if there exists some vertex $v \in V(\mathcal{M})$ which has at least $n/{2^{\ell+5s}}$\COMMENT{This number has to be sufficient so that, when we multiply by all of the vertices in a molecule (that is, $2^{\ell +s}$), it is still small.} neighbours in $\cQ^n$ which are part of $(s,\ell)$-molecules of $\mathcal{U}$ which are bondless in $G$. 
Both bondless and bondlessly surrounded molecules create difficulties in applying the rainbow matching lemma (\cref{lem: rainbow}), which in turn is used to assign absorption structures to vertices.\COMMENT{Loosely, they correspond to some molecule potentially being overused when it comes to finding absorbing cube structures.}
Therefore, it will become important that we bound the number of each, and show that they are well spread out.

\begin{lemma}\label{lem:moleculegood}
  Let $\varepsilon>0$ and $\ell,s,n\in\mathbb{N}$ be such that $s<n$, $\ell\leq n-s$ and $1/\ell\ll\varepsilon$.
  Then, for any $(s,\ell)$-molecule $\mathcal{M}\subseteq\cQ^n$, the probability that it is bondless in $\cQ^n_{\varepsilon}$ is at most $2^{s+1-\varepsilon 2^{\ell}/4}$.
\end{lemma}

\begin{proof}
Fix an $(s,\ell)$-molecule $\mathcal{M} = \cA_1 \cup \dots \cup \cA_{2^s} \subseteq\cQ^n$.
Consider a pair of consecutive atoms $\mathcal{A}_i,\mathcal{A}_{i+1}\subseteq\mathcal{M}$, for some $i\in[2^s]$.
Let $X_i$ be the number of edges between $\mathcal{A}_i$ and $\mathcal{A}_{i+1}$ in $\cQ^n_{\varepsilon}$ whose endpoint in $\mathcal{A}_i$ is odd, and let $Y_i$ be the number of such edges whose endpoint in $\mathcal{A}_i$ is even.
We have that $X_i,Y_i\sim\text{Bin}(2^{\ell-1},\varepsilon)$.
By \cref{lem:Chernoff}, it follows that\COMMENT{Observe that $\mathbb{E}[X_i]=\varepsilon 2^{\ell-1}$.
By \cref{lem:Chernoff}, we have that
\begin{align*}
    \mathbb{P}[X_i<100]&=\mathbb{P}[X_i\leq99]=\mathbb{P}\left[X_i\leq\left(1-\left(1-\frac{99}{\varepsilon 2^{\ell-1}}\right)\right)\varepsilon 2^{\ell-1}\right]\\
    &\leq e^{-\left(1-\frac{99}{\varepsilon 2^{\ell-1}}\right)^2\varepsilon 2^{\ell-1}/2} =2^{-\log_2e\left(1-\frac{99}{\varepsilon 2^{\ell-1}}\right)^2\varepsilon 2^{\ell-1}/2}\leq 2^{-\varepsilon 2^{\ell}/4},
\end{align*}
where in the last inequality we use the fact that $1/\ell\ll\varepsilon$ to get the bound $\log_2e\left(1-\frac{99}{\varepsilon 2^{\ell-1}}\right)^2\geq1$.}
\[\mathbb{P}[X_i<100]\leq2^{-\varepsilon 2^{\ell}/4},\]
and the same bound holds for $\mathbb{P}[Y_i<100]$.
By a union bound over all $i\in[2^s]$, we conclude that
\[\mathbb{P}[\mathcal{M}\text{ is bondless in }\cQ_\varepsilon^n]\leq2^{s+1-\varepsilon 2^{\ell}/4}.\qedhere\]
\end{proof}

\begin{lemma}\label{lem: T2B don't clump}
Let $\varepsilon\in(0,1)$ and $\ell,n\in\mathbb{N}$ with $0 < 1/n\ll1/\ell\ll\varepsilon$, and let $s \coloneqq 10\ell$.
Let $\mathfrak{M}$ be a collection of vertex-disjoint $(s,\ell)$-molecules $\mathcal{M}\subseteq\cQ^n$.
For each $x\in V(\cQ^n)$, let $N^\mathfrak{M}(x) \coloneqq \{\mathcal{M}\in\mathfrak{M}: \dist(x, \cM) =1\}$.
Assume that the following holds for every $x \in V(\cQ^n)$:
\begin{enumerate}[label=$(\mathrm{BS})$]
    \item\label{fact2} for any direction $\hat e\in\mathcal{D}(\cQ^n)$, there are at most $\sqrt{n}$ molecules $\mathcal{M} \in N^\mathfrak{M}(x)$ such that $\hat{e} \in \cD(\cA)$ for all atoms $\cA \in \cM$.
\end{enumerate}
Then, with probability at least $1-2^{-n^{9/8}}$, for every $x\in V(\cQ^n)$ we have that $B^{\ell^2}_{\cQ^n}(x)$ intersects at most $n^{1/3}$ molecules from $\mathfrak{M}$ which are bondlessly surrounded in $\cQ^n_\varepsilon$.
\end{lemma}

\begin{proof}
We begin by fixing an arbitrary vertex $x\in V(\cQ^n)$ and an arbitrary set $B\subseteq\mathfrak{M}$ of $n^{1/3}$ molecules which intersect $B^{\ell^2}_{\cQ^n}(x)$.
We will estimate the probability that all of the molecules in $B$ are bondlessly surrounded in $\cQ^n_\varepsilon$, by considering the neighbourhoods of the different vertices which make up these molecules.
If the probability of being bondlessly surrounded was independent over different molecules and vertices, then this would be a straightforward calculation.
However, there are dependencies which we must consider: namely, when two different molecules have edges to the same third molecule.
We will first bound the number of such configurations in $\cQ^n$.
Since the molecules in $\mathfrak{M}\supseteq B$ are vertex-disjoint, it follows that, if two of these molecules are adjacent in $\cQ^n$, then all of their atoms are pairwise adjacent in each of the layers, via clones of the same edges.
Thus, we can restrict the analysis to a single layer.

Fix a layer $L$ and let $\mathfrak{A}$ be the collection of atoms obtained by intersecting each molecule $\mathcal{M}\in\mathfrak{M}$ with $L$.
Let $\mathfrak{A}_B\subseteq\mathfrak{A}$ be the set of such atoms whose molecules lie in $B$.
Fix an atom $\mathcal{A}\in \mathfrak{A}_B$, and let $y \in V(\mathcal{A})$ be a fixed vertex.
We say an atom $\mathcal{A}'\in\mathfrak{A}$ is $y$\emph{-dependent} if there exists $\mathcal{A}''\in \mathfrak{A}_B$, $\mathcal{A}''\neq\mathcal{A}$, such that $\dist(y,\mathcal{A'})=\dist(\mathcal{A}',\mathcal{A}'')=1$. 
The following claim will allow us to bound the number of $y$-dependent atoms.

\begin{claim}\label{claim:cubeneighbours}
Fix $\mathcal{A}''\in \mathfrak{A}_B$ with $\mathcal{A}''\neq\mathcal{A}$.
Then, the number of $\cA' \in \mathfrak{A}$ for which $\dist(y,\mathcal{A'})=\dist(\mathcal{A}',\mathcal{A}'')=1$ is at most $(\ell+2)2^\ell(2 + \sqrt{n})$.
\end{claim}

\begin{claimproof}
Let $z \in V(\cA'')$; we are going to bound the number of choices of atoms $\cA' \in \mathfrak{A}$ for which $\dist(y,\mathcal{A'})=\dist(\mathcal{A}',z)=1$, and later we will add over all possible choices of $z$.
Note that we may assume $|\cD(y,z)|\leq\ell+2$, as otherwise there cannot be any atom $\mathcal{A}'\in\mathfrak{A}$ satisfying our condition.
Let $\hat{e} \in \cD(y,z)$.
We are going to bound the number of atoms $\cA' \in \mathfrak{A}$ such that $\dist(y,\mathcal{A'})=\dist(\mathcal{A}',z)=1$ for which $\hat e \in \cD(\cA')$ and those for which $\hat e \notin \cD(\cA')$ separately.

Consider first the number of options for an atom $\cA' \in \mathfrak{A}$ such that $\dist(y,\mathcal{A'})=\dist(\mathcal{A}',z)=1$ and $\hat{e} \notin \cD(\cA')$.
Since any path between $y$ and $z$ must contain an edge in the direction of $\hat e$, we must have either $y+\hat{e} \in V(\cA')$ or $z+\hat{e} \in V(\cA')$.
Since all the atoms in $\mathfrak{A}$ are vertex-disjoint, this leaves at most two options for $\cA'\in \mathfrak{A}$.\COMMENT{The unique (if it even exists) atom such that $y+\hat{e} \in V(\cA')$, and such that $z+\hat{e} \in V(\cA')$.}
Alternatively, suppose $\hat{e} \in \cD(\cA')$.
Then, by \ref{fact2}  applied with $y$ playing the role of $x$, we have at most $\sqrt{n}$ possibilities for $\cA'\in \mathfrak{A}$.
Thus, in total there are at most $2+\sqrt{n}$ possible choices for $\mathcal{A}'$.
Finally, by considering all $\hat e\in\cD(y,z)$ and all $z \in V(\cA'')$, we prove the claim.
\end{claimproof}

By considering all possibilities for $\mathcal{A''}\in \mathfrak{A}_B$, since $|\mathfrak{A}_B| = n^{1/3}$, it follows by \cref{claim:cubeneighbours} that the number of $y$-dependent atoms is  at most $n^{6/7}$.\COMMENT{We have $n^{1/3} 2^\ell(2 + \sqrt{n}) < n^{6/7}$.}
For each $y\in V(\mathcal{A})$, let $N'(y)\subseteq N^\mathfrak{M}(y)$ be given by removing from $N^{\mathfrak{M}}(y)$ all molecules which contain a $y$-dependent atom.
It follows that $|N'(y)|=|N^\mathfrak{M}(y)| - o(n)$\COMMENT{We have that $|N'(y)|\geq |N^\mathfrak{M}(y)|- n^{6/7}$.} for every $y\in V(\mathcal{A})$.

Let $\cM_\cA \in \mathfrak{M}$ be the molecule containing $\cA$.
For each vertex $y \in V(\mathcal{A})$, let $\mathcal{E}_y$ be the event that $N'(y)$ contains at least $n/{2^{\ell+5s+1}}$ molecules $\mathcal{M}\in\mathfrak{M}$ which are bondless in $\cQ^n_\varepsilon$.\COMMENT{This is half the number of bondless molecules required for the definition of bondlessly surrounded. 
That is, this is a necessary condition for being bondless, even if we assume that all molecules in $N^\mathfrak{M}(y)\setminus N'(y)$ are bondless.
So showing that this happens with very small probability suffices.}\COMMENT{Note if this is not true for all $x$ then a molecule will have degree at most $n/\ell^3$ in $\mathcal{H}$ (in the application of the rainbow lemma).}
Then, $|N'(y)| \ge n/{2^{\ell+5s+2}}$.
Moreover, we only consider here those vertices $y \in V(\cA)$ for which $|N^\mathfrak{M}(y)| \ge n/2^{\ell + 5s +1}$, since otherwise $y$ cannot contribute towards $\cM_\cA$ being bondlessly surrounded. 
Fix such a vertex $y$.
Let $Y$ be the number of molecules $\mathcal{M}\in N'(y)$ which are bondless in $\cQ^n_\varepsilon$.
Note that  $Y$ is a sum of independent indicator variables. 
By \cref{lem:moleculegood}, we have that $\mathbb{E}[Y] \leq 2^{s+1 -\eps2^\ell/4}n$.
In order to derive a lower bound for $\mathbb{E}[Y]$, note that the probability that an $(s,\ell)$-molecule $\cM$ is bondless can be bounded from below by the probability that there are no edges between two fixed consecutive atoms $\cA_1, \cA_2 \subseteq \cM$, whose endpoints in $\cA_1$ are even.
This occurs with probability $(1-\varepsilon)^{2^{\ell-1}}$.
Thus, 
\[\mathbb{E}[Y]\geq (1-\varepsilon)^{2^{\ell-1}}|N'(y)| \ge (1-\varepsilon)^{2^{\ell-1}}(n/{2^{\ell+5s+2}}).\]
By \cref{lem:betaChernoff}, we have that $\mathbb{P}[\mathcal{E}_y] \leq 2^{-cn}$,
for some constant $c>0$ which depends on $\ell$ and $\varepsilon$.\COMMENT{
We have  \[\mathbb{E}[Y]\geq(1-\varepsilon)^{2^{\ell-1}}(n/{2^{\ell+5s+2}})=2^{-c'2^\ell}n,\] for some constant $c'$ which depends on $\eps$ and $\ell$.
Then, by \cref{lem:betaChernoff} and since $s = 10\ell$ we have that 
\begin{align*}
    \mathbb{P}[Y\geq n/2^{\ell+5s+1}]& \le \mathbb{P}[Y\geq 7\mathbb{E}[Y]]\leq e^{-\mathbb{E}[Y]}\leq 2^{-c'n},
\end{align*}
for some $c'$ which depends on $\ell$ and $\varepsilon$.}
For each atom $\mathcal{A}\in \mathfrak{A}_B$, let $\cB_\cA$ be the event that there exists a vertex $y\in V(\mathcal{A})$ such that $\mathcal{E}_y$ holds.
Let $\cB \coloneqq \bigwedge _{\cA \in \mathfrak{A}_B}\cB_\cA$.
Note that the definition of $N'(y)$ ensures that the events $\cB_\cA$ with $\cA \in \mathfrak{A}_B$ are pairwise independent.
Thus, \COMMENT{The first inequality follows by considering a union bound for all $y\in V(\mathcal{A})$, and then observing that the events $\cB_\cA$ for different atoms $\cA$ are independent, by the definition of $N'(y)$.}
\[\mathbb{P}[\mathcal{B}] \le (2^{\ell-cn})^{n^{1/3}} < 2^{-n^{5/4}}.\]
In turn, this means that the probability that all molecules $\mathcal{M}\in B$ are bondlessly surrounded is bounded from above by $2^{-n^{5/4}}$.
\cref{lem: T2B don't clump} now follows by a union bound over the $2^n$ choices for $x$ and the at most $\binom{n^{\ell^2}}{n^{1/3}}$ choices for $B$.\COMMENT{Note that $\binom{n^{\ell^2}}{n^{1/3}}\leq n^{\ell^2n^{1/3}}=e^{\ell^2n^{1/3}\log n}$.
We have that $2^n\cdot e^{\ell^2n^{1/3}\log n}\cdot2^{-n^{-5/4}} < 2^{- n^{9/8}}$.}
\end{proof}

Finally, we will show that `scant' molecules are not too clustered. (We will later define a vertex molecule as `scant' --with respect to a graph $H$ and a reservoir $R$-- if one of its vertices $v_i$ has the property that few of its neighbours lie in the $i$-th clone of $R$.)

\begin{lemma}\label{lem: scant don't clump}
Let $C,s,n\in\mathbb{N}$ such that $0<1/n\ll1/C\ll\alpha,\delta\leq1$ and $1/n\ll1/s$.
Let $H \subseteq \cQ^n$ be such that $\delta(H) \geq \alpha n$.
For each $v \in V(\cQ^{n-s})$ and each $i\in[2^s]$, let $v_i$ be the $i$-th clone of $v$, and let $\mathcal{M}_v\coloneqq\{v_i: i\in[2^s]\}$. 
Let $R \sim \mathit{Res}(\cQ^{n-s}, \delta)$ and, for each $i \in [2^s]$, let $R_i$ be the $i$-th clone of $R$.
Let 
\[B \coloneqq \{\mathcal{M}_v\mid v \in  V(\cQ^{n-s}), \text{ there exists } i \in [2^s]: e_{H}(v_i, R_i) < \alpha\delta n/4\}.\]
Let $\mathcal{E}$ be the event that there exists some $u \in V(\cQ^{n-s})$ such that $B_{\cQ^{n-s}}^{10\ell}(u)$ contains more than $C$ vertices $v\in V(\cQ^{n-s})$ with $\mathcal{M}_v \in B$.
Then, $\mathbb{P}[\mathcal{E}]< e^{-n}$.
\end{lemma}

\begin{proof}
Let $u \in V(\cQ^{n-s})$ and let $D\subseteq B_{\cQ^{n-s}}^{10\ell}(u)$ be a set of $C$ vertices. 
Let $D'\coloneqq\bigcup_{x,y\in D:x\neq y}N_{\cQ^{n-s}}(x)\cap N_{\cQ^{n-s}}(y)$.
Since any pair of distinct vertices share at most two neighbours, we have that $|D'| \leq 2\binom{C}{2}$.
For each $i \in [2^s]$, we denote the $i$-th clone of $D'$ by $D'_i$, and let $R'_{i}\coloneqq R_{i}\setminus D'_{i}$.

For each $x\in V(\cQ^{n})$, let $i(x)$ be the unique index $i\in[2^s]$ such that $x\in V(L_i)$.
Observe that $e_H(x, V(L_{i(x)})) >2\alpha n/3$\COMMENT{Since $s$ is constant.} for every $x\in V(\cQ^{n})$.
For each $x\in V(\cQ^{n})$, let $\mathcal{E}_x$ be the event that $e_H(x, R_{i(x)})\leq\alpha\delta n/4$, and let $\mathcal{E}'_x$ be the event that $e_H(x, R'_{i(x)})\leq\alpha\delta n/4$\COMMENT{Note that if $\mathcal{E}_x$ occurs, then $\mathcal{E}_x'$ occurs, that is, $\mathbb{P}[\mathcal{E}_x]\leq\mathbb{P}[\mathcal{E}'_x]$.}.
It follows by \cref{lem:Chernoff} that $\mathbb{P}[\mathcal{E}'_x]\leq e^{-\alpha\delta n/16}$\COMMENT{Let $Y\coloneqq e_H(x, R'_{i(x)})$.
We have that $\mathbb{E}[Y]\geq\delta(\alpha n-s-2\binom{C}{2})\geq\alpha\delta n/2$.
It follows by \cref{lem:Chernoff} that
\[\mathbb{P}[Y\leq\alpha\delta n/4]\leq\mathbb{P}[Y\leq\mathbb{E}[Y]/2]\leq e^{-\mathbb{E}[Y]/8}\leq e^{-\alpha\delta n/16}.\]
}
for all $x\in V(\cQ^n)$.
For each $v\in V(\cQ^{n-s})$, let $\mathcal{E}_v$ and $\mathcal{E}_v'$ be the events that there exists $i\in[2^s]$ such that $\mathcal{E}_{v_i}$ and $\mathcal{E}_{v_i}'$ hold, respectively.
By a union bound, it follows that $\mathbb{P}[\mathcal{E}_v']\leq 2^se^{-\alpha\delta n/16}$ for all $v\in V(\cQ^{n-s})$.
Finally, let $\mathcal{E}_D$ and $\mathcal{E}'_D$ be the events that $\mathcal{E}_v$ and $\mathcal{E}_v'$, respectively, hold for every $v\in D$.
Note that the events in the collection $\{\mathcal{E}_v': v\in V(\cQ^{n-s})\}$ are mutually independent.
Furthermore, since the event $\mathcal{E}_x$ implies $\mathcal{E}_x'$ for all $x\in V(\cQ^n)$, we have that
\[\mathbb{P}[\mathcal{E}_D] \leq \mathbb{P}[\mathcal{E}'_D] \leq (2^se^{-\alpha\delta n/16})^{C} < e^{-5n}.\]
Taking a union bound over all vertices $u$ and over all choices of $D$ we obtain the result.\COMMENT{$2^{n}n^{10C\ell}e^{-5n}<e^{-n}$.}
\end{proof}


\subsection{Connecting cubes}\label{connect}

The hypercube satisfies some robust connectivity properties.
The problem of (almost) covering $\cQ^n$ with disjoint paths has been extensively studied.

In order to create a long cycle, which can be used to absorb all remaining vertices, while preserving the absorbing structure, we will make use of the robust connectivity  properties of the hypercube.
In particular, we will need several results which guarantee that, given any prescribed pairs of vertices in a slice, there is a spanning collection of vertex-disjoint paths, each of which uses the vertices of one of the given pairs as endpoints.
We will also need similar results for almost spanning collections of paths, where these paths avoid a given prescribed vertex.
Throughout this subsection we denote by $uv$ the edge between two given adjacent vertices $u$ and $v$ (instead of $\{u,v\}$).

The following lemma will be essential for us.
It follows from some results of \citet[Corollary~5.2]{DG08}\COMMENT{There is a huge amount of literature about this type of results. For the case i) this is not the original reference, although it seems to be the first reference for case ii). Should we add a different reference for i)?}.

\begin{lemma}\label{lem:connectcubes}
  For all $n\geq100$, the graph $\cQ^n$ satisfies the following.
  \begin{enumerate}[label=$(\mathrm{\roman*})$]
      \item\label{lem:connectcubesnormal} Let $m\in[25]$ and let $\{u_i,v_i\}_{i\in[m]}$ be disjoint pairs of vertices with $u_i\neqp v_i$ for all $i\in[m]$.
      Then, there exist $m$ vertex-disjoint paths $\mathcal{P}_1,\ldots,\mathcal{P}_m\subseteq\cQ^n$ such that, for each $i\in[m]$, $\mathcal{P}_i$ is a $(u_i,v_i)$-path, and $\bigcup_{i\in[m]}V(\mathcal{P}_i)=V(\cQ^n)$.
      \item\label{lem:connectcubesavoid} Let $x\in V(\cQ^n)$.
      Let $m\in[25]$ and let $\{u_i,v_i\}_{i\in[m]}$ be disjoint pairs of vertices of $\cQ^n-\{x\}$ such that $u_1,v_1\neqp x$ and $u_i \neqp v_i$ for all $i \in [m] \setminus\{1\}$.
      Then, there exist $m$ vertex-disjoint paths $\mathcal{P}_1,\ldots,\mathcal{P}_m\subseteq\cQ^n$ such that, for each $i\in[m]$, $\mathcal{P}_i$ is a $(u_i,v_i)$-path, and $\bigcup_{i\in[m]}V(\mathcal{P}_i)=V(\cQ^n)\setminus\{x\}$.
      \item\label{lem:connectcubesRareParities} Let $\{u_i,v_i\}_{i\in[2]}$ be disjoint pairs of vertices with $u_i\eqp v_i$ for all $i\in[2]$ and $u_1\neqp u_2$.
      Then, there exist two vertex-disjoint paths $\mathcal{P}_1,\mathcal{P}_2\subseteq\cQ^n$ such that, for each $i\in[2]$, $\mathcal{P}_i$ is a $(u_i,v_i)$-path, and $V(\mathcal{P}_1)\cup V(\mathcal{P}_2)=V(\cQ^n)$.
  \end{enumerate}
\end{lemma}

We now motivate the statement (as well as the proof) of \cref{lem:slicecover}, which is the main result of this subsection.
We are given a slice $\mathcal{M}^*$ of a molecule $\mathcal{M}\subseteq\cQ^n$ which is bonded in a graph $G\subseteq\cQ^n$.
Furthermore, we are given collections of vertices $L$, $R$ (which are part of absorbing cube structures), and $S$  (which, when constructing a long cycle, will be used to enter and leave $\mathcal{M}^*$).
More specifically, we have that
\begin{itemize}
    \item $L$ will have size $0$ or $2$, and will consist of left absorber tips.
    If it has size $2$, the vertices will have opposite parities.
    These must be avoided by our connecting paths, so that we can make use of the absorbing structures we have put in place  (see the discussion in \cref{sect:outline main}).
    \item $R$ will consist of the pairs of right absorber tip and third absorber vertex.
    These must be connected via an edge with the paths we find.
    \item $S$ will consist of a set of pairs of vertices $\{u,v\}$ with $u\neqp v$.
    Later, when creating a long cycle, $u$ will be a vertex through which we enter $\mathcal{M}^*$ from a different molecule, and $v$ will be the next vertex from which we leave $\mathcal{M}^*$ (with respect to some ordering).
    Each of our paths will be a $(u,v)$-path, for some such pair $\{u,v\}$.
\end{itemize}

In order to find our paths, we will call on \cref{lem:connectcubes}.
To illustrate this, suppose $\mathcal{M}^*$ consists of the atoms $\mathcal{A}_1, \dots, \mathcal{A}_t$, for some $t \in \mathbb{N}$. 
Suppose that $S = \{u,v\}$ with $u \in V(\mathcal{A}_1)$ and $v \in V(\mathcal{A}_t)$.
Furthermore, suppose that $L,R = \varnothing$.
To construct a path from $u$ to $v$, we will first specify the edges used to pass between different atoms.
For all $k \in [t-1]$, we choose an edge $v_k^{\uparrow}u_{k+1}^{\uparrow}$ from $\mathcal{A}_k$ to $\mathcal{A}_{k+1}$, thus $v_k^{\uparrow}\neqp u_{k+1}^{\uparrow}$.
For technical reasons, we aim to have all the vertices $u_{k+1}^{\uparrow}$ of the same parity as $u$.
We can then apply \cref{lem:connectcubes} to find a path from $u_{k+1}^{\uparrow}$ to $v_{k+1}^{\uparrow}$ which covers all of $V(\cA_{k+1})$.
Together with the edges $v_k^\uparrow u_{k+1}^\uparrow$, all these paths will form a single path from $u$ to $v$ which spans $V(\cM^*)$.
In the more general setting where $u \in V(\cA_i)$ and $v \in V(\cA_j)$ with $1<i<j<t$, the $(u,v)$-path we construct would first pass down to $\cA_1$, then up to $\cA_t$ and, finally, back down to $\cA_j$.

When $L\neq\varnothing$, due to vertex parities, the following issue can arise.
Suppose $L = \{x,y\}$ with $x \in V(\cA_1)$, $u \in V(\cA_2)$, $y \in V(A_3)$ and $v \in V(\cA_j)$ for some $j>3$ (and $R = \varnothing$).
Furthermore, suppose that both $u$ and $x$ have odd parity.
In line with the above description, the vertex $u^{\downarrow}_1$, through which we enter $\cA_1$, would have odd parity.
It follows that, since $x$ also has odd parity, we cannot hope to construct a path which starts at $u^{\downarrow}_1$ and covers all of $V(\cA_1) \setminus\{x\}$.
The solution will be instead to pass up to $\cA_3$ first (and, in general, to whichever atom contains~$y$).
Recall that, since $x$ has odd parity, $y$ must have even parity.
We specify a vertex $u^{\uparrow}_3$ of odd parity, through which we enter $\cA_3$, but then also specify a vertex $v^{\downarrow}_3$ of odd parity from which we will leave $\cA_3$ to reenter $\cA_2$.
We now arrive back in $\cA_2$ with a vertex $u^{\downarrow}_2$ of even parity.
We will specify another vertex $v^{\downarrow}_2$ of odd parity from which we leave $\cA_2$ and a vertex $u^{\downarrow}_1$ of even parity through which we enter $\cA_1$.
In this way, we can now apply \cref{lem:connectcubes} to find a path which starts at $u^{\downarrow}_1$ and covers all of $V(\cA_1) \setminus\{x\}$, and which can be extended into a path from $u$ to $v$ covering all of $V(\cM^*)\setminus L$. 

There are several other instances which must be dealt with in a similar way.
This is formalised by \cref{lem:slicecover}.
Before proving this lemma, however, we need the following definition.

\begin{definition}[$(u,j,F,R)$-alternating parity sequence]\label{def:seqtotal}
Let $\ell, s, t, n \in \mathbb{N}$ with $t \leq 2^s$\index{(uj@$(u,j,F,R)$-alternating parity sequence} and $2\leq\ell\leq n-s$.
Let $G \subseteq \cQ^n$.
Let $\cM= \cA_1 \cup \dots \cup \cA_{2^s}\subseteq\cQ^n$ be an $(s,\ell)$-molecule and let $\cM^* = \cA_{a+1} \cup \dots \cup \cA_{a+t}$, for some $a\in[2^s]$, be a slice of $\cM$.  
Let $u \in V(\cA_i)$, for some $i \in [a+t]\setminus[a]$.
Let $j\in[a+t]\setminus[a]$, and let $F,R \subseteq V(\cM^*)$.
Suppose $i \le j$.
Let $I_R \coloneqq \{k \in [j-i]_0: |R \cap V(\cA_{i+k})| \ge 1\}$.
Assume that the following properties hold:
\begin{itemize}
    \item For all $k \in [j-i]_0$ we have that $|R \cap V(\cA_{i+k})| \in \{0,2\}$.
    \item For each $k \in I_R$, the vertices in $R\cap V(\cA_{i+k})$ are adjacent in $\cQ^n$, and we write $R\cap V(\cA_{i+k})=\{w_k, z_k\}$ so that $w_k\neqp u$.
\end{itemize}
Let $\cS' = (u_0,v_1,u_1, \dots, v_{j-i}, u_{j-i})$ be a sequence of vertices satisfying the following properties:
\begin{enumerate}[label=$(\mathrm{P}\arabic*)$,start=0]
    \item\label{def:aps1} If $u\in R$, then $u_0\coloneqq w_{0}$; otherwise, $u_0\coloneqq u$.
    \item\label{def:aps2} For each $k \in [j-i]$ we have that $u_{k}\eqp u$.
    \item\label{def:aps3} For each $k \in [j-i]$ we have that $v_k \in V(\cA_{i+k-1})$, $u_k \in V(\cA_{i+k})$ and $v_ku_k\in E(G)$.
    \item\label{def:aps4} The vertices of\/ $\cS'$ other than $u_0$ avoid $F\cup R$.
\end{enumerate}

A \emph{$(u,j,F,R)$-alternating parity sequence $\cS$ in $G$} is a sequence obtained from any sequence $\cS'$ which satisfies \ref{def:aps1}--\ref{def:aps4} as follows.
For each $k\in I_R\cap[j-i]$, replace each segment $(v_{k}, u_{k})$ of $\cS'$ by $(v_{k}, u_{k}, w_{k}, z_{k})$.

The case $i>j$ is defined similarly by replacing each occurrence of $[j-i]$ and $[j-i]_0$ in the above by $[i-j]$ and $[i-j]_0$, and each occurrence of  $\cA_{i+k}$ and $\cA_{i+k-1}$ by $\cA_{i-k}$ and $\cA_{i-k+1}$.
\end{definition}

Given an alternating parity sequence $\cS$, we will denote by $\cS^{-}$ the sequence obtained from $\cS$ by deleting its initial element.

\begin{lemma}\label{lem:slicecover}
Let $n,s,\ell \in \mathbb{N}$ be such that $s\geq4$ and $100\leq\ell\leq n-s$.
Let $G\subseteq\cQ^n$ and consider any $(s,\ell)$-molecule $\mathcal{M}=\mathcal{A}_1\cup\dots\cup\mathcal{A}_{2^s}\subseteq\cQ^n$ which is bonded in $G$.
Let $\mathcal{M}^*=\mathcal{A}_{a+1}\cup\dots\cup\mathcal{A}_{a+t}$, for some $a \in [2^s]$ and $t \ge 10$, be a slice of $\cM$.
Moreover, consider the following sets.
\begin{enumerate}[label=$(\mathrm{C}\arabic*)$]
    \item\label{itm:conn1.1} Let $L\subseteq V(\mathcal{M}^*)$ be a set of size $|L|\in\{0,2\}$ such that, if $L=\{x,y\}$, then $x\in V(\mathcal{A}_i)$ and $y\in V (\mathcal{A}_j)$ with $i \neq j$ and $x \neqp y$.
    \item\label{itm:conn1.2} Let $R\subseteq V(\mathcal{M}^*)\setminus L$ be a (possibly empty) set of vertices with $|R|\leq10$ such that, for all $k\in[a+t]\setminus[a]$, we have $|R\cap V(\mathcal{A}_k)|\in\{0,2\}$ and, if $|R\cap V(\mathcal{A}_k)|=2$, then $R\cap V(\mathcal{A}_k)=\{w_k,z_k\}$ satisfies that $w_kz_k \in E(\cM^*)$ and, if $|L|=2$, then $k\notin\{i,j\}$.
    \item\label{itm:conn1.3} Let $m\in[14]$ and consider $m$ vertex-disjoint pairs $\{u_r,v_r\}_{r\in[m]}$, where $u_r,v_r\in V(\mathcal{M}^*)\setminus L$ and $u_r\neqp v_r$ for all $r\in[m]$, such that, for each $r\in[m]$, we have $u_r\in V(\mathcal{A}_{i_r})$ and $v_r\in V(\mathcal{A}_{j_r})$\COMMENT{Possibly $i_r=j_r$.}.
    Assume, furthermore, that for each $t' \in [t]$ we have that $|\bigcup_{r \in [m]}\{u_r, v_r\} \cap V(\mathcal{A}_{a+ t'})\cap R|\leq1$.
\end{enumerate}
Then, there exist vertex-disjoint paths $\mathcal{P}_1,\ldots,\mathcal{P}_m\subseteq\mathcal{M}^*\cup G$ such that, for each $r\in[m]$, $\mathcal{P}_r$ is a $(u_r,v_r)$-path, $\bigcup_{r\in[m]}V(\mathcal{P}_r)=V(\mathcal{M}^*)\setminus L$, and every pair $\{w_{k},z_{k}\}$ with $k \in [a+t]\setminus[a]$ is an edge of some $\mathcal{P}_{r}$.
\end{lemma}

\begin{proof}
By relabelling the atoms, we may assume that $\mathcal{M^*}=\mathcal{A}_1\cup\dots\cup\mathcal{A}_t$.\COMMENT{Here we implicitly assume that the indices of vertices in the lemma statement ($i_1, j_1$ etc.) have all been relabelled as well.}
Let $S\coloneqq\{u_r,v_r: r\in[m]\}$.
By relabelling the vertices, we may assume that $i_r\leq j_r$ for all $r\in[m]$ and (if $L\neq\varnothing$) $i<j$.
Let $I_L\coloneqq\{k\in[t]: L\cap V(\mathcal{A}_k) \neq\varnothing\}$, $I_R\coloneqq\{k\in[t]: R\cap V(\mathcal{A}_k)\cap S\neq\varnothing\}$
and $R^*\coloneqq R \setminus \bigcup_{k \in I_R}V(\cA_{k})$.
Note that $I_L = \varnothing$ or $I_L = \{i,j\}$ and $I_L \cap I_R = \varnothing$. 
For each $r\in[m]$, let $I_R^r\coloneqq\{k \in \{i_r, j_r\} : R \cap V(\cA_k)\cap\{u_r, v_r\}\neq\varnothing\}$, so that $I_R = \bigcup_{r=1}^mI_R^r$.
Without loss of generality, we may also assume that, for each $r \in [m]$, if $u_r \in R$, then $u_r = z_{i_r}$, and if $v_r\in R$, then $v_r = w_{j_r}$.
Similarly, for each $k\in[t]\setminus I_R$, if $R\cap V(\mathcal{A}_k)=\{w_k,z_k\}$, we may assume that $w_k\neqp u_1$.

For each $r\in[m]$, we will create a list $\mathcal{L}_r$ of vertices.
We will refer to $\mathcal{L}_r$ as the \emph{skeleton} for $\mathcal{P}_r$.
We will later use these skeletons to construct the vertex-disjoint paths via \cref{lem:connectcubes}. 
For each $r\in[m]$, we will write $L^*_r$ for the (unordered) set of vertices in $\mathcal{L}_r$.
In order to construct each $\mathcal{L}_r$, we will start with an empty list and update it in (possibly) several steps, by concatenating alternating parity sequences.
Whenever $\mathcal{L}_r$ is updated, we implicitly update $L^*_r$.
In the end, for each $r\in[m]$ we will have a list of vertices $\mathcal{L}_r = (x_1^r, \dots, x_{\ell_r}^r)$.
For each $r\in[m]$ and $k \in [t]$, let $I_r(k) \coloneqq \{h\in[\ell_r-1]:2\nmid h\text{ and }x^r_h, x^r_{h+1} \in V(\cA_{k})\}$.
We will require the $\cL_r$ to be pairwise vertex-disjoint.
Furthermore, we will require that they satisfy the following properties:
\begin{enumerate}[label=$(\mathcal{L}\arabic*)$]
    \item\label{itm:conn0} For all $r\in[m]$ we have that $\ell_r$ is even.
    \item\label{itm:conn1} For all $r\in[m]$ and $h \in[\ell_r-1]$, if $h$ is odd, then $x^r_h,x^r_{h+1} \in V(\cA_{k})$, for some $k \in [t]$; if $h$ is even, then $x^r_hx^r_{h+1} \in E(G\cup\mathcal{M}^*)$.
    \item\label{itm:conn2} For all $k\in[t]$ we have that $1 \le|I_1(k)|\le 6$ and $|I_r(k)|\le 1$ for all $r\in[m]\setminus\{1\}$.
\end{enumerate}
\begin{enumerate}[label=$(\mathcal{L}4)_1$]
   \item\label{itm:conn31} For each $k \in [t] \setminus(I_L \cup I^1_R)$ and each $h \in I_1(k)$, we have $x_h^1 \neqp x_{h+1}^1$.
   For each $k \in I_L \cup I^1_R$, for all but one $h \in I_1(k)$ we have $x^1_h \neqp x_{h+1}^1$, while for the remaining index $h \in I_1(k)$ we have that $x^1_h \eqp x_{h+1}^1$ and their parity is opposite to that of the unique vertex in $L \cap V(\cA_k)$ if $k \in I_L$ and to that of the unique vertex in $\{w_k, z_k\} \cap \{u_1, v_1\}$ if $k \in I_R^1$.
\end{enumerate}
\begin{enumerate}[label=$(\mathcal{L}4)_r$]
    \item\label{itm:conn32} For each $r \in [m]\setminus\{1\}$, the following holds.
    For each $k \in [t]\setminus I^r_R$ and each $h \in I_r(k)$, we have $x_h^r \neqp x_{h+1}^r$.
    For each $k \in I_R^r$, for all but one $h \in I_r(k)$ we have $x^r_h \neqp x_{h+1}^r$, while for the remaining index $h \in I_r(k)$ we have that $x^r_h \eqp x_{h+1}^r$ and their parity is opposite to that of the unique vertex in $\{w_k, z_k\} \cap \{u_r, v_r\}$. 
\end{enumerate}
\begin{enumerate}[label=$(\mathcal{L}\arabic*)$,start=5]
    \item\label{itm:conn4} For each $r\in[m]$, we have the following. 
    If $u_r \notin R$, then $u_r = x_1^r$.
    If $v_r \notin R$, then $v_r = x^r_{\ell_r}$.
    If $u_r \in R$ (and thus $u_r = z_{i_r}$), then $w_{i_r} = x_1^r$ and $u_r \notin L^*_1 \cup \dots \cup L^*_m$.
    If $v_r \in R$ (and thus $v_r = w_{j_r}$), then $z_{j_r} = x_{\ell_r}^r$ and $v_r \notin L^*_1 \cup \dots \cup L^*_m$.
    \item\label{itm:conn5} Every pair $(w_k, z_k)$ with $\{w_k, z_k\} \subseteq R^*$ is contained in $\cL_1$ and $z_k$ directly succeeds $w_k$.
\end{enumerate}

We begin by constructing $\mathcal{L}_1$.
Let $\cL_1\coloneqq \varnothing$ and let  $F\coloneqq L \cup R \cup S$.
If $i_1 =1$ and $R^*\cap V(\cA_1) = \{w_1, z_1\}$, then let $\cS_1 \coloneqq (u_1, w_1, z_1)$.
If $i_1=1$ and $u_1\in R$, then let $\cS_1\coloneqq(u_1)$.
Otherwise, let $\cS_1$ be a $(u_1, 1, F, (R \cap V(\cA_{i_1}))\cup (R^*\cap V(\cA_{1})))$-alternating parity sequence.
Let $\mathcal{L}_1\coloneqq\cS_1$.
Note that the existence of such a sequence $\cS_1$ is guaranteed by our assumption that $\mathcal{M}$ is bonded in $G$.
To see this, note that all edges of $G$ required by $\cS_1$ (that is, the pairs $\{v_k, u_k\}$ in \cref{def:seqtotal}) need to be chosen so that they do not have an endpoint in $F$; given any particular pair of consecutive atoms, this forbids at most $32$ edges between these two atoms ($28$ because of $S$ and $4$ because of $L\cup R$\COMMENT{Note that an atom cannot contain vertices from $L$ and $R$ simultaneously, so the worst case here is when the two consecutive atoms contains 2 vertices from $R$ each--giving 4 forbidden edges.}).

We will now update $\cL_1$. 
While doing so, we will update $F$ and consider several alternating parity sequences.
The existence of each of these follows a similar argument to the above.
For any given pair of consecutive atoms, every time we update $F$, the set of forbidden edges will increase its size by at most $3$.
We will update $F$ at most four times, so $F$ will forbid at most $44$ edges between any pair of consecutive atoms.
Thus, by the definition of bondedness, each of the alternating parity sequences required below actually exists.

Let $u^\downarrow_1$ be the last vertex in $\cL_1$.
Note that $u_1^\downarrow \eqp u_1$ by \cref{def:seqtotal}\ref{def:aps2}.
We update $F$ as $F\coloneqq F \cup L^*_1$.
For the next step in the construction of $\mathcal{L}_1$, there are three cases to consider, depending on the size of $L$ and, if $|L|=2$\COMMENT{We need this here so that $x$ is defined.}, the relative parities of $x$ and $u_1$.
If $i_1=1$ and $u_1\in R$, let $R^\diamond\coloneqq R^*\cup\{w_1,z_1\}$; otherwise, let $R^\diamond\coloneqq R^*$.

\textbf{Case 1}: $L=\varnothing$.\\
Let $\cS_2$ be a $(u^\downarrow_1, t, F, R^\diamond)$-alternating parity sequence.
If $i_1=1$ and $u_1\in R$, update $\mathcal{L}_1$ as $\mathcal{L}_1\coloneqq\mathcal{S}_2$.
Otherwise, update $\mathcal{L}_1$ as $\mathcal{L}_1\coloneqq\mathcal{L}_1\cS_2^{-}$.
Update $F\coloneqq F \cup L_1^*$.

\textbf{Case 2}: $|L|=2$ and $x\neqp u_1$.\\
Let $\cS_2$ be a $(u^\downarrow_1, i, F, R^\diamond)$-alternating parity sequence.
If $i_1=1$ and $u_1\in R$, update $\mathcal{L}_1$ as $\mathcal{L}_1\coloneqq\mathcal{S}_2$.
Otherwise, update $\mathcal{L}_1\coloneqq\mathcal{L}_1\cS_2^{-}$.
Update $F\coloneqq F \cup L_1^*$.
Choose any vertex $u_i^*\in V(\mathcal{A}_i)$ with $u_i^*\neqp u_1$\COMMENT{
Note. We are going to discard the special vertex $u^*_i$ here in the next line.
This is also why we can pick it from anywhere in $V(\mathcal{A}_i)$.
Note also that the last vertex in a alternating parity sequence by definition has the same parity as $u_1$. 
So picking $u_i^*$ of parity opposite to that of $u_1$ is the same as picking it of parity opposite to that of the final vertex of $S_2$, which is really what we want here.
This means we will be left with a sequence which begins with a vertex of the same parity as $u_j^{\uparrow1}$ (the last vertex we added to our list).
This `change in parity' is what will allow us to correct the parity imbalance caused by the left absorber vertices in this case, as discussed in the bla bla at the start.}, and let $\cS_3$ be a $(u_i^*, j, F, R^\diamond)$-alternating parity sequence.
Update $\mathcal{L}_1\coloneqq\mathcal{L}_1\cS_3^{-}$ and $F\coloneqq F \cup L_1^*$.
Let $v^-$ be the final vertex of $\cS_2$, and let $v^+$ be the second vertex of $\cS_3$.
Note that $v^-$ and $v^+$ appear consecutively in $\cL_1$ and that $v^- \eqp v^+ \eqp u_1 \neqp x$.
Finally, choose any vertex $u_j^*\in V(\mathcal{A}_j)$ with $u_j^*\eqp u_1$, let $\cS_4$ be a $(u_j^*, t, F, R^\diamond)$-alternating parity sequence, and update $\mathcal{L}_1\coloneqq\mathcal{L}_1\cS_4^{-}$ and $F\coloneqq F \cup L_1^*$.
Let $w^-$ be the final vertex of $\cS_3$, and let $w^+$ be the second vertex of $\cS_4$.
We then have that $w^-$ and $w^+$ appear consecutively in $\cL_1$, and $w^-\eqp w^+ \neqp u_1\eqp y$ (recall that $x\neqp y$).
Moreover, the final vertex $u_t^\uparrow$ of $\cL_1$ satisfies $u_t^\uparrow \eqp u_1$. 

\textbf{Case 3}: $|L|=2$ and $x\eqp u_1$.\\
Let $\cS_2$ be a $(u^\downarrow_1, j, F, R^\diamond)$-alternating parity sequence.
If $i_1=1$ and $u_1\in R$, update $\mathcal{L}_1$ as $\mathcal{L}_1\coloneqq\mathcal{S}_2$; otherwise, update $\mathcal{L}_1\coloneqq\mathcal{L}_1\cS_2^{-}$.
Update $F\coloneqq F \cup L_1^*$.
Next, let $u^*_j \in V(\mathcal{A}_j)$ be a vertex with $u^*_j\neqp u_1$ and let $\cS_3$ be a 
$(u^*_j, i, F, \varnothing)$-alternating parity sequence.
 Update $\mathcal{L}_1\coloneqq\mathcal{L}_1\cS_3^{-}$ and $F\coloneqq F \cup L_1^*$.
Finally, let $u^*_i \in  V(\mathcal{A}_i)$ be a vertex with $u^*_i\eqp u_1$ and let $\cS_4$ be a 
$(u^*_i, t, F, R^*\cap\bigcup_{k=j+1}^tV(\mathcal{A}_k))$-alternating parity sequence.
Update $\mathcal{L}_1\coloneqq\mathcal{L}_1\cS_4^{-}$ and $F\coloneqq F \cup L_1^*$.

In each of the three cases, let $u_t^\uparrow$\COMMENT{We use $t$ here as the index because the last vertex is in $\cA_t$.} denote the last vertex in $\mathcal{L}_1$.
Note that, by \cref{def:seqtotal}\ref{def:aps2}, we have $u_t^\uparrow \eqp u_1$, and recall that $v_1 \neqp u_1$.
Let $\cS_5$ be a $(u^{\uparrow}_t, j_1, F, \varnothing)$-alternating parity sequence.
Update $\mathcal{L}_1\coloneqq\mathcal{L}_1\cS_5^{-}$.
Again by \cref{def:seqtotal}\ref{def:aps2}, we have that the final vertex $u^*$ of $\cL_1$ is such that $u^* \eqp u_t^\uparrow \eqp u_1 \neqp v_1$.
Finally, if $v_1\in R$,  update $\mathcal{L}_1\coloneqq\mathcal{L}_1(z_{j_1})$; otherwise, update it as $\mathcal{L}_1\coloneqq\mathcal{L}_1(v_1)$.
Observe that $\mathcal{L}_1$ satisfies \ref{itm:conn0}--\ref{itm:conn2}, \ref{itm:conn31}, \ref{itm:conn4} and \ref{itm:conn5} for the case $r=1$ by construction.

We now construct $\mathcal{L}_r$ for all $r\in[m] \setminus \{1\}$.
For each $r\in[m]\setminus\{1\}$, we proceed iteratively as follows.
Let $\mathcal{L}_r\coloneqq\varnothing$ and $F_r\coloneqq L \cup R \cup S \cup\bigcup_{r'\in[r-1]}L_{r'}^*$.
Let $\cS^r$ be a $(u_r, j_r, F_r, R\cap V(\cA_{i_r}))$-alternating parity sequence and update  $\mathcal{L}_r$ as $\mathcal{L}_r\coloneqq\cS^r$.
If $v_r\in R$, update $\mathcal{L}_r\coloneqq\mathcal{L}_r(z_{j_r})$; otherwise, update $\mathcal{L}_r\coloneqq\mathcal{L}_r(v_r)$.
Note that each sequence $\cS^r$ requires the existence of at most one edge of $G$, which has to avoid $F_r$, between any pair of consecutive atoms of $\cM^*$.
In a similar way to what was discussed above, at most three choices of such edges can be forbidden every time we add a new alternating parity sequence to $F$.
Since for each $r\in[m]\setminus\{1\}$ we consider one new sequence, by the time we consider $F_m$ we have increased the number of forbidden edges by at most $3(m-1)\leq 39$.
This gives a total of at most $83$ forbidden edges and, thus, the existence of the sequences $\cS^r$ is guaranteed by the assumption that $\mathcal{M}$ is bonded in $G$.
Moreover, the lists $\cL_1, \dots, \cL_r$ now satisfy \ref{itm:conn0}--\ref{itm:conn5}.\COMMENT{
\ref{itm:conn0}, \ref{itm:conn1} and \ref{itm:conn32} like before, and to see \ref{itm:conn2} for $r>1$ note in these cases the list just consists of a sequence in one direction.
Finally, we enforced both \ref{itm:conn4} explicitly in our constructions above.} 

We are now in a position to apply \cref{lem:connectcubes}.
For each $k \in [t]$, let $t_k\coloneqq \sum_{r \in [m]}|I_r(k)|$.
Furthermore, for any $r\in[m]$ and $k\in[t]$, for each $h\in I_r(k)$, we refer to the pair $x^r_h, x^r_{h+1}$ as a $\emph{matchable pair}$.
By \ref{itm:conn2}, \ref{itm:conn31}, \ref{itm:conn32} and \cref{lem:connectcubes}\ref{lem:connectcubesnormal}, each atom $\mathcal{A}_k$ with $k\in[t]\setminus(I_L\cup I_{R})$ can be covered by $t_k$ vertex-disjoint paths, each of whose endpoints are a matchable pair contained in $\mathcal{A}_k$.
Similarly, by \ref{itm:conn2}, \ref{itm:conn31}, \ref{itm:conn32} and \cref{lem:connectcubes}\ref{lem:connectcubesavoid}, each atom $\mathcal{A}_k$ with $k\in I_L\cup I_R$ contains $t_k$ vertex-disjoint paths, each of whose endpoints are a matchable pair in $\cA_k$ such that the union of these $t_k$ paths covers precisely $V(\mathcal{A}_k)\setminus (L\cup (S\cap R))$.
(Recall that by \ref{itm:conn1.2} and \ref{itm:conn1.3} the set $V(\cA_k) \cap (L\cup (S\cap R))$ consists of a single vertex if $k \in I_L \cup I_R$.)
For each matchable pair $x^r_h, x^r_{h+1}$ in $\cA_k$, let us denote the corresponding path by $\mathcal{P}_{x^r_h, x^r_{h+1}}$.

The paths $\cP_1, \dots, \cP_m$ required for \cref{lem:slicecover} can now be constructed as follows.
For each $r\in[m]$, let $\mathcal{P}_r$ be the path obtained from the concatenation of the paths $\mathcal{P}_{x^r_{h},x^r_{h+1}}$, for each odd $h\in[\ell_r]$, via the edges $x^r_{h}x^r_{h+1}$ for $h\in[\ell_r-1]$ even.
By \ref{itm:conn4}, if $\mathcal{P}_r$ does not contain $u_r$, then $\mathcal{P}_r$ starts in $w_{i_r}$, and $u_r$ does not lie in any other path; therefore, we can update $\mathcal{P}_r$ as $\mathcal{P}_r\coloneqq u_r\mathcal{P}_r$.
Similarly, if $\mathcal{P}_r$ does not contain $v_r$, then $\mathcal{P}_r$ ends in $z_{j_r}$ and $v_r$ does not lie in any other path, and thus we can update $\mathcal{P}_r$ as $\mathcal{P}_r\coloneqq\mathcal{P}_rv_r$.
It follows that $\bigcup_{r\in[m]}V(\mathcal{P}_r) = V(\cM^*)\setminus L$, and thus the paths $\mathcal{P}_r$ are as required in \cref{lem:slicecover}.
\end{proof}

We also need the following simpler result.
Its proof follows similar ideas as those present in the proof of \cref{lem:slicecover}.
For the sake of completeness, we include the proof of \cref{lem:slicecover2} in \cref{app:connect}.
We point out here that \cref{lem:connectcubes}\ref{lem:connectcubesRareParities} is only needed for this proof.

\begin{lemma}\label{lem:slicecover2}
Let $n,s,\ell\in\mathbb{N}$ be such that $4\leq s$ and $100\leq\ell\leq n-s$.
Let $G\subseteq\cQ^n$ and consider any $(s,\ell)$-molecule $\mathcal{M}=\mathcal{A}_1\cup\dots\cup\mathcal{A}_{2^s}\subseteq\cQ^n$ which is bonded in $G$.
Let $\cM^* = \cA_{a+1} \cup \dots \cup \cA_{a+t}$, for some $a \in [2^s]$ and $t\geq10$, be a slice of $\cM$.
Moreover, consider the following sets.

\begin{enumerate}[label=$(\mathrm{C}'\arabic*)$]
\item Let $L\subseteq V(\mathcal{M}^*)$ be a set of size $|L|\in\{0,2\}$ such that, if $L=\{x,y\}$, then $x\in V(\mathcal{A}_i)$ and $y\in V(\mathcal{A}_j)$, with $i\neq j$ and $x \neqp y$.
\item Let $R\subseteq V(\mathcal{M}^*)\setminus L$ be a (possibly empty) set of vertices with $|R|\leq10$ such that, for all $k\in[a+ t]\setminus[a]$, we have $|R\cap V(\mathcal{A}_k)|\in\{0,2\}$ and, if $|R\cap V(\mathcal{A}_k)|=2$, then $R\cap V(\mathcal{A}_k)=\{w_k,z_k\}$ satisfies that $w_kz_k \in E(\cM^*)$ and, if $|L|=2$, then $k\notin\{i,j\}$.
\item Consider two vertex-disjoint pairs $\{u_r,v_r\}_{r\in[2]}$ with $u_1,u_2\in V(\mathcal{A}_{a+1})\setminus L$ and $v_1,v_2\in V(\mathcal{A}_{a+t})\setminus L$ such that $u_1\neqp u_2$, $v_1\neqp v_2$, $u_1\eqp v_1$, and $|\{u_1,u_2\}\cap R|,|\{v_1,v_2\}\cap R|\leq1$.
\end{enumerate}

Then, there exist two vertex-disjoint paths $\mathcal{P}_1,\mathcal{P}_2\subseteq\mathcal{M}^*\cup G$ such that, for each $r\in[2]$, $\mathcal{P}_r$ is a $(u_r,v_r)$-path, $V(\mathcal{P}_1)\cup V(\mathcal{P}_2)=V(\mathcal{M}^*)\setminus L$, and every pair of the form $\{w_{k},z_{k}\}\subseteq R$ with $k \in [a+t]\setminus[a]$ is an edge of either $\mathcal{P}_{1}$ or $\mathcal{P}_{2}$.
\end{lemma}


\subsection{Proof of \texorpdfstring{\cref{thm:main1}}{Theorem 8.1}}\label{sect:mainpf}

\begin{proof}[Proof of \cref{thm:main1}]
Let $1/D, \delta' \ll 1$, and let\COMMENT{I am going to list here the relationships we need. \\
$\varepsilon,\alpha\leq1$ are given by the statement, so we should not need to say they are $\leq1$ again.\\
$1/D,\delta'\ll1/5$ are given by \cref{lem: main treereshit}.\\
$1/S'\ll\alpha,\delta'$. This is needed in Step~2, where a bad event is when there's more than $S'$ scant molecules in a ball.\\
$1/k^*\ll1/D,\delta',\varepsilon,\alpha$ and $\beta\ll\delta',\varepsilon,\alpha$ come from \cref{lema:robustmatch}.\\
$\alpha'\ll\beta$ from the nibble.\\
$1/\ell\ll\alpha'$ used in \ref{itm:di1} comment.\\
$\delta\ll1/\ell$ is also needed for the nibble, and used in other places.\\
Note that $s$ (and similarly $q$) should not be part of the hierarchy, as the statement is not true if $s$ is far larger than $\ell$?} 
\begin{equation}\label{equa:hierarchy1}
    0 < 1/n_0 \ll \delta \ll  1/\ell \ll 1/k^*, \alpha' \ll \beta, 1/S' \ll 1/c, 1/D, \delta', \varepsilon, \alpha,
\end{equation}
where $n_0,\ell,k^*,S',D\in\mathbb{N}$.
Our proof assumes that $n$ tends to infinity; in particular, $n\geq n_0$.
Let $s\coloneqq 10\ell$, \[\Phi\coloneqq12\ell\] and \[\Psi\coloneqq c\Phi.\]

Observe that $\cQ^n[\{0,1\}^s\times\{0\}^{n-s}]\cong\cQ^s$ contains a Hamilton cycle. 
We fix an ordering of the layers $L_1,\ldots,L_{2^s}$ of $\cQ^n$ induced by this Hamilton cycle (as defined in \cref{sect8notation}).
If we view these layers as different subgraphs on the vertex set of the same copy of $\cQ^{n-s}$, then for any $G\subseteq\cQ^n$ we can define $I(G)\coloneqq \bigcap_{i=1}^{2^s} L_i(G)$\index{IG@$I(G)$}.\COMMENT{Intersection over edges.}\COMMENT{We will treat these intersection graphs as auxiliary graphs, not as subgraphs of $\cQ^n$.}
In particular, the \emph{intersection graph} $I$ of the layers is defined as $I\coloneqq \bigcap_{i=1}^{2^s} L_i$.
So actually $I\cong\cQ^{n-s}$.
Recall that $L_i(G)=G[V(L_i)]$.
Note that, if $\cG\subseteq I(G)$, then there is a clone of $\cG$ in $L_i(G)$, for each $i\in[2^s]$.
For each layer $L$, we denote by $\cG_L$ the clone of $\cG$ in $L(G)$.
Observe that, for any $\eta\in[0,1]$, we have $I(\cQ^n_\eta)\sim\cQ^{n-s}_{\eta^{2^s}}$.
We will sometimes write $G_I$\index{GI@$G_I$} for the subgraph of $I$ where, for each $e\in E(I)$, we have $e\in E(G_I)$ whenever $G$ contains some clone of $e$ (thus, $G_I$ is the `union' of the subgraphs that $G$ induces on each layer).

For each $i \in [7]$, let $\eps_i \coloneqq \eps/7$ and let $G_i \sim \cQ^n_{\eps_i}$, where these graphs are taken independently. 
It is easy to see that $\bigcup_{i=1}^7 G_i \sim \cQ^n_{\eps'}$ for some $\eps' < \eps$.
Thus, it suffices to show that a.a.s.~there is a graph $G'\subseteq\bigcup_{i=1}^7 G_i$ with $\Delta(G')\leq\Phi$ such that, for every $F\subseteq \cQ^{n}$ with $\Delta(F)\leq \Psi$, the graph $((H\cup\bigcup_{i=1}^7 G_i)\setminus F)\cup G'$ is Hamiltonian.
We now split our proof into several steps.\\

\textbf{Step~1: Finding a tree and a reservoir.}
Consider the probability space $\Omega \coloneqq\cQ^{n-s}_{\eps_1^{2^s}} \times \mathit{Res}(\cQ^{n-s},\delta')$ (with the latter defined as in \cref{section:tree1}), so that, given $R\sim\mathit{Res}(I,\delta')$, we have that $(I(G_1),R)\sim\Omega$.

Let $\mathcal{E}_1$ be the event that there exists a tree $T \subseteq I(G_1)-R$ such that the following hold:
\begin{enumerate}[label=$(\mathrm{TR}\arabic*)$]
\item $\Delta(T) < D$, and
\item for all $x \in V(I)$, we have that $|N_{I}(x)\cap V(T)|\geq4(n-s)/5$.
\end{enumerate}
It follows from \cref{lem: main treereshit}, with $n-s$, $\eps_1^{2^{s}}$, $\delta', \varnothing$ and $1/5$ playing the roles of $n$, $\eps$, $\delta, \cA$ and $\eps'$, respectively, that $\mathbb{P}_\Omega[\mathcal{E}_1]=1-o(1)$.\\

\textbf{Step~2: Identifying scant molecules.}
For each $v \in V(I)$, let $\cM_v$ denote the vertex molecule $\cM_v \coloneqq\{av: a \in \{0,1\}^s\}$\COMMENT{$\cM_v$ gives us the molecule consisting of the different clones of $v$ in $\cQ^n$.}.
We say a vertex molecule $\cM_v$ is \emph{scant} if there exist some layer $L$ and some vertex $x \in V(\cM_v \cap L)$ such that $d_H(x, R_L) < \alpha\delta' n/10$, where $R_L$ is the clone of $R$ in $L$.\COMMENT{Note here that scant is defined with respect to the deterministic graph $H$.}
Let $\mathcal{E}_2$ be the event that there exists some $x\in V(I)$ such that there are more than $S'$ vertices $v \in B_{I}^{10\ell}(x)$ satisfying that $\mathcal{M}_v$ is scant (where $S'$ satisfies \eqref{equa:hierarchy1}).
It follows from \cref{lem: scant don't clump} with $S'$ and $\delta'$ playing the roles of $C$ and $\delta$ that $\mathbb{P}_\Omega[\mathcal{E}_2] < e^{-n}$.
Let $\mathcal{E}_1^*\coloneqq\mathcal{E}_1\wedge\overline{\mathcal{E}_2}$.
Therefore, $\mathbb{P}_{\Omega}[\mathcal{E}_1^*]=1-o(1)$.

Condition on $\mathcal{E}_1^*$ holding.
Then, $G_1$ satisfies the following: there exist a set $R\subseteq V(I)$\index{R1@$R$} and a tree $T\subseteq I(G_1)-R$\index{T1@$T$} such that the following hold:
\begin{enumerate}[label=$(\mathrm{T}\arabic*)$]
\item\label{itm:tree1} $\Delta(T) < D$;
\item\label{itm:tree2} for all $x \in V(I)$, we have that $|N_{I}(x) \cap V(T)| \ge 4(n-s)/5$, and
\item\label{itm:scant} for every $x\in V(I)$, $B_{I}^{10\ell}(x)$ contains at most $S'$ vertices $v$ such that $\cM_v$ is scant.
\end{enumerate}
Recall this implies clones of $T$ and $R$ satisfying \ref{itm:tree1}--\ref{itm:scant} exist simultaneously in each layer of $G_1$.\\

\textbf{Step~3: Finding robust matchings for each slice.}
Recall from \cref{section:outline5} that we will absorb vertices in pairs, where each pair consists of two clones $x'$, $x''$ of the same vertex $x \in V(I)$.
In this step, for each $x \in V(I)$ and for each set of clones of $x$ that may need to be absorbed, we find a pairing of these clones so that we can later build suitable absorbing $\ell$-cube pairs for each such pair of clones.
We will find this pairing separately for each slice of the vertex molecule $\cM_x$.
Considering each slice separately has the advantage that the chosen pairs are `localised'.
This will be convenient later when linking up the paths used to absorb these vertices.
Accordingly, we now partition the set of layers into sets of consecutive layers as follows.
Let 
\begin{equation}\label{qandt}
q\coloneqq2^{10Dk^*} \qquad \text{ and let } \qquad t\coloneqq2^s/q,
\end{equation}
where $k^*$ satisfies \eqref{equa:hierarchy1}.
For each $j\in[t]$, let $S_j\coloneqq\bigcup_{i=(j-1)q+1}^{jq}L_i$.
Given any molecule $\mathcal{M}$, we consider the slices $\mathcal{S}_j(\mathcal{M})\coloneqq S_j\cap\mathcal{M}$.
We denote by $\cS(\cM)$ the collection of all these slices of $\cM$.

Let $V_{\mathrm{sc}} \subseteq V(I)$\index{Vsc@$V_{\mathrm{sc}}$} be the set of all vertices $x\in V(I)$ such that $\cM_x$ is scant.
Recall $G_2\sim \cQ^{n}_{\eps_2}$.
For each $v\in V(I)\setminus V_{\mathrm{sc}}$ and each $\cS\in\mathcal{S}(\cM_v)$, we define the following auxiliary bipartite graphs.
Let $H(\cS)\coloneqq (V(\cS),N_{I}(v),E_{H})$, where $E_{H}$ is defined as follows.
Consider $v'\in V(\mathcal{S})$ and let $L^{v'}$ be the layer which contains $v'$.
Let $w\in N_{I}(v)$, and let $w_{L^{v'}}$ be the clone of $w$ in $L^{v'}$.
Then, $\{v',w\}\in E_{H}$ if and only if $w\in R$ and $\{v',w_{L^{v'}}\}\in E(H)$.
Note that $d_{H(\cS)}(v') \ge \alpha \delta' n/10$ for all $v'\in V(\mathcal{S})$ since $\cS$ is a slice of a vertex molecule which is not scant.
Similarly, we define $G_2(\cS)\coloneqq (V(\cS),N_{I}(v),E_{G_2})$, where $\{v',w\}\in E_{G_2}$ if and only if $\{v',w_{L^{v'}}\}\in E(G_2)$.

Note that the partition of $V(\mathcal{S})$ into vertices of even and odd parity is a balanced bipartition.
Define the graph $\Gamma^{\beta}_{H(\mathcal{S}),G_2(\cS)}(V(\cS))$ as in \cref{sect5matchings} (where $\beta$ satisfies \eqref{equa:hierarchy1}).
Note that, by definition, we have that $V(\Gamma^{\beta}_{H(\mathcal{S}),G_2(\cS)}(V(\cS)))=V(\cS)$.
Furthermore, by definition, 
\begin{enumerate}[label=$(\mathrm{RM})$]
    \item\label{item:RM} given any $w_1,w_2\in V(\mathcal{S})$, we have that $\{w_1,w_2\}\in E(\Gamma^{\beta}_{H(\mathcal{S}),G_2(\cS)}(V(\cS)))$ if and only if $|N_{H(\cS)}(w_1)\cap N_{G_2(\cS)}(w_2)|\geq \beta (n-s)$ or $|N_{G_2(\cS)}(w_1)\cap N_{H(\cS)}(w_2)|\geq \beta (n-s)$.\COMMENT{Here we are using that $|B| = n-s$ in \cref{lema:robustmatch}.}
\end{enumerate}
By applying \cref{lema:robustmatch} with $d=24D$, $r=0$, $\alpha=\alpha\delta'/10$\COMMENT{This gives us the condition in the statement because $\mathcal{M}_v$ is not scant.}, $\eps=\eps_2$, $n= n-s$, $k=q=2^{10Dk^*}$, $\beta = \beta$ and $G=H(\cS)$, we obtain that, with probability at least $1-2^{-10(n-s)}\geq 1-2^{-8n}$, the graph $\Gamma^{\beta}_{H(\mathcal{S}),G_2(\cS)}(V(\cS))$ is $24D$-robust-parity-matchable with respect to the partition of $V(\mathcal{S})$ into vertices of even and odd parity\COMMENT{Note that here we are using the fact that $|V(\cS)|=q>k^*$, where $k^*$ is given in the starting hierarchy.}.

We would like to proceed as above for slices in scant molecules; however, recall that scant molecules contain vertices with few or no neighbours in the reservoir, and therefore we must adapt our approach.
For each $v\in V_{\mathrm{sc}}$ and each $\cS\in\cS(\cM_v)$, we define an auxiliary bipartite graph $H(\mathcal{S})$ and $G_2(\cS)$
as above, except that we omit the condition that $w\in R$ for the existence of an edge in $H(\mathcal{S})$. 
By applying \cref{lema:robustmatch} again, we obtain that, with probability at least $ 1-2^{-8n}$, the graph $\Gamma^{\beta}_{H(\mathcal{S}),G_2(\cS)}(V(\cS))$ is $24D$-robust-parity-matchable with respect to the partition of $V(\cS)$ into vertices of even and odd parity.

By a union bound over all $v\in V(I)$ and all slices $\mathcal{S}\in\mathcal{S}(\mathcal{M}_v)$, we have that a.a.s.~the graph $\Gamma^{\beta}_{H(\mathcal{S}),G_2(\cS)}(V(\cS))$ is $24D$-robust-parity-matchable (with respect to the partition of $V(\mathcal{S})$ into vertices of even and odd parity) for every slice $\mathcal{S}$\COMMENT{This happens with probability at least $1-2^n 2^{-8n}$.}, where $H(\mathcal{S})$ is as defined above in each case.
We condition on this event holding and call it $\cE_2^*$.
Thus, for each slice $\mathcal{S}$ and each set $X\subseteq V(\mathcal{S})$ with $|X|\leq 24D$ which contains as many odd vertices as even vertices, there exists a perfect matching $\mathfrak{M}(\mathcal{S},X)$ in the bipartite graph with parts consisting of the even and odd vertices of $V(\mathcal{S})\setminus X$, respectively, and edges given by $\Gamma^{\beta}_{H(\mathcal{S}),G_2(\cS)}(V(\cS))$.
For each slice $\cS$, we denote by $\mathfrak{M}(\cS)$\index{MS@$\mathfrak{M}(\cS)$, $\mathfrak{M}(v)$} the set of edges contained in the union (over all $X$) of the matchings $\mathfrak{M}(\mathcal{S},X)$ (without multiplicity).
Furthermore, for each $e=\{w_\mathrm{e},w_\mathrm{o}\} \in \mathfrak{M}(\cS)$, we let $N(e)\coloneqq(N_{H(\cS)}(w_\mathrm{e})\cap N_{G_2(\cS)}(w_\mathrm{o}))\cup(N_{G_2(\cS)}(w_\mathrm{e})\cap N_{H(\cS)}(w_\mathrm{o}))$\index{Ne@$N(e)$}\COMMENT{Note that $N(e)\subseteq V(I)$, by definition.}.
By \ref{item:RM}, we have $|N(e)|\geq\beta(n-s)\geq\beta n/2$.
For each $v\in V(I)$, let $\mathfrak{M}(v)\coloneqq\bigcup_{\cS \in \cS(\cM_v)} \mathfrak{M}(\mathcal{S})$.
Let $K \coloneqq \max_{v \in V(I)}|\mathfrak{M}(v)|$.\COMMENT{In other words, $K$ is the maximum number of edges in the union of all matchings over all the robust-parity-matchable graphs given by each slice of each molecule.}
In particular, we have that $K \leq \binom{2^s}{2}$\COMMENT{For all $v\in V(I)$, $\mathcal{M}_v$ contains $2^s$ vertices.
Since $\bigcup_{\cS \in \cS(\cM_v)} \mathfrak{M}(\mathcal{S})$ can be seen as a graph on vertex set $V(\mathcal{M}_v)$, we get the claim as a trivial upper bound.}.\\

\textbf{Step~4: Obtaining an appropriate cube factor via the nibble.}
For each $x\in V(I)$, consider the multiset  $\mathfrak{A}(x)\coloneqq\{N(e): e \in \mathfrak{M}(x)\}$.
If  $|\mathfrak{A}(x)| < K$, we artificially increase its size to $K$ by repeating any of its elements.
Label the sets in $\mathfrak{A}(x)$ arbitrarily as $\mathfrak{A}(x)=\{A_1(x),\ldots,A_K(x)\}$.
Thus, if $x\in V(I)\setminus V_\mathrm{sc}$, then $A_i(x)\subseteq R$ for all $i\in[K]$.

Let $\mathcal{C}$ be any collection of subgraphs $C$ of $I$ such that $C\cong\cQ^\ell$ for all $C\in\mathcal{C}$.
For any vertex $x\in V(I)$ and any set $Y\subseteq N_I(x)$, let $\mathcal{C}_x(Y)\subseteq\mathcal{C}$\index{CxY@$\mathcal{C}_x(Y)$} be the set of all $C\in\mathcal{C}$ such that $x\notin V(C)$ and $Y\cap V(C)\neq\varnothing$, and let $\mathcal{C}_x\coloneqq\mathcal{C}_x(N_I(x))$\index{Cx@$\mathcal{C}_x$}.

Recall $G_3\sim \cQ^{n}_{\eps_3}$ and $I(G_3) \sim \cQ^{n-s}_{\eps^{2^s}_3}$. 
We now apply \cref{thm: nibble} to the graph $I(G_3)$, with $\varepsilon_3^{2^s}, \alpha', \delta/2$, $\beta/2$, $K$ and $\ell$ playing the roles of $\varepsilon, \alpha$, $\delta$, $\beta$, $K$ and $\ell$, respectively, and using the sets $A_i(x)$ given above, for each $x\in V(I)$ and $i\in[K]$.
(Recall $\alpha'$ satisfies \eqref{equa:hierarchy1}.)
Thus, a.a.s.~we obtain a collection $\mathcal{C}$\index{Cw1@$\mathcal{C}$} of vertex-disjoint copies of $\cQ^\ell$ in  $I(G_3)$, such that the following properties hold for every $x\in V(I)$: 
\begin{enumerate}[label=$(\mathrm{N}\arabic*)$]
    \item\label{itm:nib1} $|\cC_x| \geq (1-\delta/2)(n-s) \geq (1-\delta)n$.
    \item\label{itm:nib3} For every direction $\hat{e}\in\mathcal{D}(I)$ we have that $|\Sigma(\mathcal{C}_x,\{\hat e\},1)|=o(n^{1/2})$.
    \item\label{itm:nib2} For every $i\in[K]$ and every $S\subseteq\mathcal{D}(I)$ with $\alpha' (n-s)/2\leq|S|\leq \alpha' (n-s)$ we have
    \[|\Sigma(\mathcal{C}_x(A_i(x)),S,\ell^{1/2})|\geq|A_i(x)|/3000\geq\beta n/6000.\]
\end{enumerate}
Condition on the above event holding and call it $\cE_3^*$.\\

\textbf{Step~5: Absorption cubes.}
For each $x \in V(I)$ and $i \in [K]$, we define an auxiliary digraph $\mathfrak{D} = \mathfrak{D}(A_i(x))$ on vertex set $A_i(x)-\{x\}$ (seen as a set of directions of $\mathcal{D}(I)$)\COMMENT{Clearly, this is a subset of $\mathcal{D}(I)$.
However, we may have a consistency problem, because we sometimes also regard directions here as directions of $\mathcal{D}(\cQ^n)$, that is, vectors of length $n-s$ and $n$ `at the same time'. 
This is important because later we add directions in $A_i(x)-\{x\}$ to vertices of $\cQ^n$, as well as to vertices in $I$.
We will try to make sure that this does not happen, but if we cannot avoid it, we will need to come up with some way to fix this.} by adding a directed edge from $\hat{e}$ to $\hat{e}'$ if there is a cube $C^r\in\cC_x(A_i(x))$ such that $x+\hat{e}\in V(C^r)$ and $\hat{e}'\in\mathcal{D}(C^r)$. 
In this way, an edge from $\hat{e}$ to $\hat{e}'$ in $\cD$ indicates that the cube $C^r$ could be used as a right absorber cube for $x$, if combined with a vertex-disjoint left absorber cube with tip $x + \hat{e}'$.
Observe that, for all $\hat e\in A_i(x)-\{x\}$,
\begin{equation}\label{equa:step4.1}
    d^+_\mathfrak{D}(\hat e)\in[\ell]_0.
\end{equation}
Furthermore, it follows by \ref{itm:nib2} that any set $S \subseteq V(\mathfrak{D})$ with $|S| = \alpha'n/2$ satisfies\COMMENT{directed edges, that is, edges that go into $S$.}\COMMENT{We have that $|V(\mathfrak{D})|=|A_i(x)|\geq\beta n/2$, so by \ref{itm:nib2} we have at least $\ell^{1/2}\beta n/6000$ edges into $S$. (To see this, recall the definition of significance and a short remark given in the corresponding paragraph. This is done after \cref{rem:basicequationIneedhere}.)} 
\begin{equation}\label{equa:step4.2}
    e_\mathfrak{D}(V(\mathfrak{D}), S) \ge \ell^{1/2}\beta n/6000 >\ell^{1/2}\beta^2 n.
\end{equation}

Recall that $A_i(x)=N(\{x_1,x_2\})$ for some $\{x_1,x_2\} \in \mathfrak{M}(\cS)$, where $\cS \in \cS(\cM_x)$ is some slice of $\cM_x$.
Note that $x_1,x_2\in\mathcal{M}_x$, and let $L^j$ be the layer containing $x_j$ for each $j\in[2]$.
We say that $x_1$ and $x_2$ are the vertices (or clones of $x$) which \emph{correspond} to the pair $(x,i)$.
Let $(\hat{e},\hat{e}')\in E(\mathfrak{D})$ and, for each $j\in[2]$, let $e_j$ be the clone of $\{x+\hat{e}',x+\hat{e}'+\hat e\}$ in $L^j$.
It follows that there is a cube $C^r \in \cC_x(A_i(x))$ such that $e_j$ connects the clone $C_j$ of $C^r$ to the clone of $x +\hat{e}'$ in $L^j$.

Recall $G_4 \sim \cQ^n_{\eps_4}$.
Let $\mathfrak{D}' \subseteq \mathfrak{D}$ be the subdigraph which retains each edge $(\hat{e},\hat{e}')\in E(\mathfrak{D})$ if and only if the edges $e_1, e_2$ described above are both present in $G_4$.
Note that each edge of $\mathfrak{D}$ is therefore retained independently of every other edge with probability $\eps_4^2$.
By \cref{lem:Chernoff}, \eqref{equa:step4.1} and \eqref{equa:step4.2}, it follows that $\mathfrak{D}'$ satisfies the following with probability at least $1-e^{-10n}$:\COMMENT{We only need to justify \ref{itm:di1}, as \ref{itm:di2} follows immediately from \eqref{equa:step4.1}.
By \eqref{equa:step4.2}, for each $A\subseteq V(\mathfrak{D})$ with $|A|=\alpha' n/2$ we have that $A$ has at least $\beta^2 \ell^{1/2}n$ in-edges in $\mathfrak{D}$.
Let $Y$ be the number of these edges retained in $\mathfrak{D}'$.
Then, $\mathbb{E}[Y] \ge \eps_4^2 \beta^2 \ell^{1/2}n$.
By \cref{lem:Chernoff}, we have that 
\[\mathbb{P}[Y\leq \eps_4^3 \beta^2 \ell^{1/2}n]\leq\mathbb{P}[Y\leq\eps_4^2 \beta^2 \ell^{1/2}n/2]\leq\mathbb{P}[Y\leq\mathbb{E}[Y]/2]\leq e^{-\mathbb{E}[Y]/8}\leq e^{-\eps_4^2 \beta^2 \ell^{1/2}n/8}.\]
In order to prove that the claim holds for every set $A$ as described, consider a union bound.
Over all sets $A$, we have a union bound of $\binom{n}{\alpha'n/2}\leq\left(\frac{2e}{\alpha'}\right)^{\alpha'n/2}$.
By the choice that $1/\ell\ll\alpha'$, we have that $\ell^{1/2}$ kills this exponent.}
\begin{enumerate}[label=$(\mathrm{DG}\arabic*)$]
    \item\label{itm:di1} for every $A\subseteq V(\mathfrak{D})$ with $|A|=\alpha'n/2$ we have $\sum_{v\in A}d_{\mathfrak{D}'}^-(v)\geq \eps_4^3 \beta^2 \ell^{1/2}n$, and
    \item\label{itm:di2} for every $B\subseteq V(\mathfrak{D})$ we have that $\sum_{v\in B}d_{\mathfrak{D}'}^+(v)\leq \ell |B|$.
\end{enumerate}
Recall that $\mathfrak{D} = \mathfrak{D}(A_i(x))$.
By a union bound, \ref{itm:di1} and \ref{itm:di2} hold a.a.s.~for all $x\in V(I)$ and $i\in[K]$.
We condition on this event and call it $\cE^*_4$.

For each $x \in V(I)$ and $i \in [K]$, recall that \ref{item:RM} and the definition of $A_i(x)$ imply that $|A_i(x)|\ge \beta(n-s)$.
Thus, it follows by \cref{prop: cubematch} with $|A_i(x)|$, $2\alpha'/\beta$, $\eps_4^3\beta^3\ell^{1/2}/(2\alpha')$ and $\ell$ playing the roles of $n$, $\alpha$, $c$ and $C$, respectively\COMMENT{With these values, by \ref{itm:di1} we have, for all $A\subseteq V(\mathfrak{D}(A_i(x)))$ with $|A|\geq2\alpha'|A_i(x)|/\beta\geq\alpha'n$, that 
\[\sum_{v\in A}d_{\mathfrak{D}'}^-(v)\geq \eps_4^3 \beta^2 \ell^{1/2}n\geq\frac{\eps_4^3\beta^3\ell^{1/2}}{2\alpha'}\frac{2\alpha'}{\beta}|A_i(x)|,\]
so the first condition for \cref{prop: cubematch} holds.
The second condition then holds easily by \ref{itm:di2}.
Finally, the condition in the statement of \cref{prop: cubematch} about the range of $\alpha$ holds by the choice of $1/\ell\ll\alpha'$.}, that there exists a matching $M''(A_i(x))$ of size at least $\frac{\eps^3_4\beta^2}{2\ell^{1/2}}|A_i(x)| \ge  \eps_4^3 \beta^3n/(3\ell^{1/2})$ in $\mathfrak{D}'(A_i(x))$.

Next, for each $x \in V(I)$ and $i \in [K]$, we remove from $M''(A_i(x))$ all edges $(\hat e,\hat e')\in M''(A_i(x))$ such that $x+\hat{e}'$ does not lie in any cube of $\cC_x(A_i(x))$.
We denote the resulting matching by $M'(A_i(x))$.
Note that, by \ref{itm:nib1}, we have
\begin{equation}\label{equa:step4.3}
    |M'(A_i(x))|\geq\eps_4^3 \beta^3n/(3\ell^{1/2}) - \delta n\geq n/\ell.
\end{equation}

Consider $A_i(x)$, for some $x\in V(I)$ and $i\in[K]$, and let $x_1$, $x_2$ be the clones of $x$ which correspond to $(x,i)$.
As before, for each $j \in [2]$, let $L^j$ be the layer containing $x_j$.
Recall \cref{def:abs} and note that, by construction, we have the following. 
\begin{enumerate}[label=$(\mathrm{AB}\arabic*)$]
    \item \label{itm:abss2} For each edge $(\hat e, \hat e')\in M'(A_i(x))$, there is an absorbing $\ell$-cube pair $(C^l,C^r)$ for $x$ in $I$ such that, for each $j \in [2]$, the clone $(C^l_j, C^r_j)$ of $(C^l,C^r)$ in $L^j$ is an absorbing $\ell$-cube pair for $x_j$ in $H \cup G_2 \cup G_3 \cup G_4$.
    In particular, the edge joining the left absorber tip to the third absorber vertex lies in $G_4$.
    Moreover, $C^l, C^r \in \cC_x(A_i(x))\subseteq \cC$ and $(C^l,C^r)$ has left and right absorber tips $x+\hat e'$ and $x+\hat e$, respectively.
    Furthermore, for each $x\in V(I)\setminus V_\mathrm{sc}$, these tips lie in $R$.
    We refer to $(C_1^l,C_1^r)$ and $(C_2^l,C_2^r)$ as the absorbing $\ell$-cube pairs for $x_1$ and $x_2$ \emph{associated} with $(\hat e,\hat e')$.
\end{enumerate}

Thus, the graph $H\cup G_2\cup G_3\cup G_4$, contains at least $n/\ell$ absorbing $\ell$-cube pairs for each of the clones $x_1$ and $x_2$ of $x$ associated with edges in $M'(A_i(x)) \subseteq \mathfrak{D}(A_i(x))$.
Moreover, since $M'(A_i(x))$ is a matching, for each $j \in [2]$ these absorbing $\ell$-cube pairs for $x_j$ are pairwise vertex-disjoint apart from $x_j$.

For ease of notation, we will often consider the absorbing $\ell$-cube pair $(C^l,C^r)$ for $x$ in $I$ which $(C_1^l,C_1^r)$ and $(C_2^l,C_2^r)$ are clones of, and use it as a placeholder for both of its clones.
By slightly abusing notation, we will refer to $(C^l, C^r)$ as the \emph{absorbing $\ell$-cube pair associated} with $(\hat e,\hat e')$\COMMENT{We can do this because the vertices $x_1,x_2$ are determined by the choice of $x\in V(I)$ and $i\in[K]$.}.
Note, however, that $(C^l,C^r)$ is not necessarily an absorbing $\ell$-cube pair for $x$ in $I(H\cup G_2\cup G_3\cup G_4)$.\COMMENT{When we take the intersection graph $I(H\cup G_2\cup G_3\cup G_4)$ we might miss out on the edges from $x$ to the tips of the cubes as well as the edge from $x+\hat{e}'$ to $C^r$.}\\

\textbf{Step~6: Removing bondless molecules.}
Recall  $G_5\sim \cQ^{n}_{\eps_5}$.
In this step, we consider the edges between the different layers.\COMMENT{In fact this is the only place where we will be interested in these edges between layers, so we could in theory spend all of our probability here if we liked, as these edges are independent from all the edges we consider in the other steps.}

For each $C\in \cC$, let $\cM_C$ denote the cube molecule consisting of the clones of $C$.
Let $\cC' \subseteq \cC$\index{Cw2@$\cC'$} be the set of cubes $C \in \cC$ for which $\cM_C$ is bonded in $G_5$.
By an application of \cref{lem:moleculegood}, for each $C\in \cC$ we have that
\[\mathbb{P}[C\notin\cC']=\mathbb{P}[\cM_C\text{ is bondless in }G_5]\leq2^{s+1-\eps_52^{\ell}/4}\leq 2^{-\eps2^{\ell}/30}.\]

For each $x\in V(I)$, let $A_0(x)\coloneqq N_I(x)$.\COMMENT{We need the $0$ case here and in what follows to know there are at most $n/\ell^4$ bondless cubes sitting on the neighbourhood of $x$.
It will be used in the first display of the next step in giving a lower bound for the sets $Z(x)$.}
For each $i \in [K]_0$, let $\mathcal{E}(x,i)$ be the event that $|\mathcal{C}_x(A_i(x))\setminus\mathcal{C}'|>n/\ell^4$.
Since the cubes $C\in\mathcal{C}$ are vertex-disjoint, the events that the molecules $\mathcal{M}_C$ are bondless in $G_5$ are independent.
Therefore, we have that 
\begin{equation*}
\mathbb{P}[\mathcal{E}(x,i)] \leq \binom{n}{n/\ell^4}(2^{-\eps2^{\ell}/30})^{n/\ell^4}\leq 2^{-10n}.
\end{equation*}
Let $\mathcal{E}_4\coloneqq\bigvee_{x\in V(I)}\bigvee_{i\in[K]_0}\mathcal{E}(x,i)$.
By a union bound over all $x \in V(I)$ and $i \in [K]_0$\COMMENT{This is a union bound over at most $2^{n-s}\binom{2^{s}}{2}$ terms.
Actually, note that proving it for $i=0$ implies the event holds for all other $i$.}, it follows that
\begin{equation}\label{equa:step5.1}
    \mathbb{P}[\mathcal{E}_4]\leq2^{-8n}.
\end{equation}

Let $\cC_{\mathrm{bs}}\subseteq\cC$ be the set of all $C\in\cC$ such that $\cM_C$ is bondlessly surrounded in $G_5$ (with respect to $\{\cM_{C'} : C' \in \cC\}$). 
For each $x\in V(I)$, let $\mathcal{E}(x)$ be the event that there are more than $n^{1/3}$ cubes $C\in\mathcal{C}_{\mathrm{bs}}$ which intersect $B_I^{\ell^2}(x)$.
Let $\mathcal{E}_5\coloneqq\bigvee_{x\in V(I)}\mathcal{E}(x)$.
By  \ref{itm:nib3}, we may apply \cref{lem: T2B don't clump} with $\varepsilon_5$ playing the role of $\varepsilon$ to conclude that 
\begin{equation}\label{equa:step5.2}
    \mathbb{P}[\mathcal{E}_5]\leq 2^{-n^{9/8}}.
\end{equation}

Now let $\mathcal{E}_5^*\coloneqq\overline{\mathcal{E}_4}\wedge\overline{\mathcal{E}_5}$.
It follows from \eqref{equa:step5.1} and \eqref{equa:step5.2} that $\mathcal{E}_5^*$ occurs a.a.s.
Condition on this event.

Let $\mathcal{C}''\coloneqq\mathcal{C}'\setminus\cC_{\mathrm{bs}}$\index{Cw3@$\cC''$}.
For each $x\in V(I)$ and each $i\in[K]$, let
\begin{enumerate}[label=$(\mathrm{AB}\arabic*)$]\setcounter{enumi}{1}
    \item\label{itm:abss1}  $M(A_i(x)) \subseteq M'(A_i(x))$ consist of all edges $(\hat e,\hat e')\in M'(A_i(x))$ whose associated absorbing $\ell$-cube pair $(C^l,C^r)$ satisfies that $C^r,C^l\in\mathcal{C}''$.
\end{enumerate} 
By combining \eqref{equa:step4.3} with the further conditioning, it follows that, for each $x\in V(I)$ and each $i\in[K]$,\COMMENT{We have that $|M(A_i(x))|\geq n/\ell-n/\ell^4-n^{1/3}\geq n/\ell^{2}$.
Recall that we are considering a matching, hence each cube whose cube molecule is bondless in $G_5$ affects at most one edge of the matching.}
\begin{equation}\label{equa:step5.3}
    |M(A_i(x))|\geq n/\ell-n/\ell^4-n^{1/3}\geq n/\ell^{2}.
\end{equation}

Consider any $x \in V(I)$ and $i \in [K]$, and let $x_1, x_2$ be the two clones of $x$ corresponding to $(x,i)$.
Then, at this point, for each $j \in [2]$, $H\cup G_2\cup G_3\cup G_4$ contains at least $n/\ell^2$ vertex-disjoint (apart from  $x_j$) absorbing $\ell$-cube pairs for $x_j$ such that each of these absorbing $\ell$-cube pairs $(C^l, C^r)$ is associated with an edge of $M(A_i(x))$, and for each $C \in \{C^l, C^r\}$ the corresponding cube molecule $\cM_C$ is bonded in $G_5$ and (within the collection $\{\cM_{C'}: C' \in \cC\}$ of all cube molecules) $\cM_C$ is not bondlessly surrounded in $G_5$.\\

\textbf{Step~7: Extending the tree $\boldsymbol{T}$.}
For each $x \in V(I)$, let $Z(x)\coloneqq N_{I}(x)\cap V(T) \cap \big(\bigcup_{C\in\cC''}V(C)\big)$.
It follows by \ref{itm:tree2}, \ref{itm:nib1} and our conditioning on the event $\mathcal{E}_5^*$ that, for each $x \in V(I)$, we have that
\[|Z(x)| \ge 4(n-s)/5 - \delta n - n/\ell^4 - n^{1/3} \ge 3n/4.\]

Recall $G_6\sim \cQ^{n}_{\eps_6}$. 
We apply \cref{thm: maintreeres} with $\varepsilon_6^{2^s}$, $2$, $T$, $R, \varnothing$ and the sets $Z(x)$ playing the roles of $\varepsilon$, $\ell$\COMMENT{We might need $\ell>2$ for hitting time result.}, $T'$, $R, W$ and $Z(x)$, respectively.
Combining this with \ref{itm:tree1}, we conclude that\COMMENT{with probability at least $1 - 2e^{-\eps^{6}n}$} a.a.s.~there exists a tree $T'$\index{T2@$T'$} such that $T\subseteq T'\subseteq I(G_6)\cup T$ and the following hold:
\begin{enumerate}[label=$(\mathrm{ET}\arabic*)$]
\item\label{itm: ET1} $\Delta(T') < D+1$;
\item\label{itm: ET2} for all $x \in V(I)$, we have that $|B^{2}_{I}(x)\setminus V(T')| \le n^{3/4}$;
\item\label{itm:ET3} for each $x \in V(T')\cap R$, we have that $d_{T'}(x)=1$ and the unique neighbour $x'$ of $x$ in $T'$ is such that $x' \in Z(x)$. 
\end{enumerate}
We condition on the above event holding and call it $\cE^*_6$.

At this point, for each $x \in V(I)$ and each $i \in [K]$, we redefine the set $M(A_i(x))$.
\begin{enumerate}[label=$(\mathrm{AB}\arabic*)$]\setcounter{enumi}{2}
\item\label{itm:AB3} Let $M(A_i(x))$ retain only those edges whose associated absorbing $\ell$-cube pair $(C^l, C^r)$ satisfies that both $C^l$ and $C^r$ intersect $T'$ in at least $2$ vertices.\COMMENT{Otherwise, we might use them as an absorbing pair but the tree might not pass into them at all. 
The $2$ is because one of these vertices could be the absorbing tip for a scant molecule.
Then after repatching, the tree wouldn't enter this cube unless there's also another vertex.}
\end{enumerate}
It follows from \eqref{equa:step5.3} and \ref{itm: ET2} that 
\begin{equation}\label{eqn:maix1}
    |M(A_i(x))| \ge  n/\ell^2 - n^{3/4} > 4n/\ell^3.
\end{equation}

\textbf{Step~8: Fixing a collection of absorbing $\boldsymbol{\ell}$-cube pairs for the vertices in scant molecules.}
Recall $G_7 \sim \cQ^n_{\eps_7}$.
Consider any  $x \in V_{\mathrm{sc}}$ and $j \in [K]$.
Recall from Step~3 that the tips of the cubes of the absorbing $\ell$-cube pair associated with a given edge in $M(A_j(x))$ may not lie in the reservoir $R$.
Roughly speaking, we will alter $T'$ so that the tips are relocated from the tree $T'$ to the reservoir $R$. 

We start by redefining the matchings $M(A_j(x))$ as follows: 
for each $x\in V_{\mathrm{sc}}$ and each $j \in [K]$, remove from $M(A_j(x))$ all edges $(\hat e,\hat e')$ such that $N_{T'}(x)\cap\{x+\hat e,x+\hat e'\}\neq\varnothing$.\COMMENT{The reason for this is that we can only apply \cref{lem:repatch} to pairs $(y,z)$ for which neither have $x$ as a neighbour.}
It follows from \eqref{eqn:maix1} and \ref{itm: ET1} that, for all $x\in V(I)$ and $j \in [K]$,
\begin{equation}\label{equa:step7.1}
    |M(A_j(x))| \ge  4n/\ell^3 - D > 2n/\ell^3.
\end{equation}

For each $x\in V_{\mathrm{sc}}$, each $j \in [K]$ and each matching $M'\subseteq M(A_j(x))$ with $|M'|\geq n/\ell^3$, let $\mathcal{E}'(x,j,M')$ be the following event:\COMMENT{
The reason for these submatchings is the following. 
Suppose we take a molecule $\cM_x$ and we apply the patching lemma to connect $x_1, x_2 \in \cM_x$.
Next we want to repatch $x_1, x_3 \in \cM_x$ say.
But suppose the lemma repatches $x_1$ and $x_3$ by using a tip from the repatched $x_1$ and $x_2$.
We can't avoid this by saying put this tip in the bad sets $B$ to avoid, because the path already avoids these vertices (it avoids tips).
So it could be used twice and submatchings help avoid this. (It's a problem for different molecules $\cM_x$ and $\cM_{x'}$.)}
\begingroup
\addtolength\leftmargini{-0.2in}
\begin{quote}
    for every set $B\subseteq V(I)$ with $|B|<2^{\ell+s+3}\Psi KS'$\COMMENT{Taking $f=2^{\ell+s+3}\Psi KS'$ is enough.}, there exists an edge $\vec{e} \in M'$, whose associated absorbing $\ell$-cube pair $(C^l, C^r)$ has tips $x^l$ and $x^r$, for which there exists a subgraph $P(\vec{e},B)\subseteq I(G_7) - \{x^l,x^r\}$ such that $|V(P(\vec{e},B))|<21D/2$, $V(P(\vec{e},B))\cap B=\varnothing$, and both $N_{T'}(x^l)$ and $N_{T'}(x^r)$ are connected in $P(\vec{e},B)$.
\end{quote}
\endgroup
\noindent For a graph $P(\vec{e},B)$ as above, we will refer to $x^l$ and $x^r$ as the tips \emph{associated} with $P(\vec{e},B)$, and refer to $(C^l,C^r)$ as the absorbing $\ell$-cube pair \emph{associated} with $P(\vec{e},B)$.
(Recall that, if $\vec{e} = (\hat{e}, \hat{e}')$, then $x^l = x+ \hat{e}$ and $x^r = x + \hat{e}'$.)

By invoking \cref{lem:repatch} with $n-s$, $\eps_7^{2^s}$, $1/\ell^3$, $2^{\ell+s+3}\Psi KS'$, $2D+2$ and the sets $\{(x+ \hat{e}, x+ \hat{e}'): (\hat{e}, \hat{e}') \in M'\}$
and $(N_{T'}(x+\hat e)\cup N_{T'}(x+\hat e'))_{(\hat e,\hat e')\in M'}$ playing the roles of $n$, $\eps$, $c$, $f$, $D$, $C(x)$ and $(B(y,z))_{(y,z)\in C(x)}$, respectively (note that, in order for this to satisfy the condition $x \notin B(y,z)$ from \cref{lem:repatch}, we must have $\{x,x+\hat e\},\{x,x+\hat e'\}\notin E(T')$; this is guaranteed by the deletion of matching edges before \eqref{equa:step7.1}), we have that $\mathcal{E}'(x,j,M')$ holds with probability at least $1-2^{-5(n-s)}$.
Let $\mathcal{E}_7^*\coloneqq\bigwedge_{x \in V_{\mathrm{sc}}}\bigwedge_{j \in [K]}\bigwedge_{M'\subseteq M(A_j(x)):|M'|\geq n/\ell^3} \mathcal{E}'(x,j,M')$.
By a union bound over all $x\in V_{\mathrm{sc}}$, $j\in[K]$ and $M'\subseteq M(A_j(x))$ such that $|M'|\geq n/\ell^3$, it follows that $\mathbb{P}[\mathcal{E}_7^*]\geq1-2^{-2n}$\COMMENT{There are at most $2^{n-s}$ choices for $x$ and at most $\binom{2^s}{2}$ choices for $j$. Then, $2^n$ is a generous upper bound for the number of choices for $M'$.}.

Condition on the event that $\mathcal{E}_7^*$ holds.
It follows that, for each $x \in V_{\mathrm{sc}}$, $j \in [K]$, $M'\subseteq M(A_j(x))$ with $|M'|\geq n/\ell^3$ and any $B \subseteq V(I)$ with $|B|<2^{\ell+s+3}\Psi KS'$, there exists a subgraph $P(x,j,M',B)\subseteq I(G_7)$ with $|V(P(x,j,M',B))|<21D/2$\COMMENT{The application of \cref{lem:repatch} says that $|V(P(x,j,M',B))|<5(2D+1)<21D/2$.} which avoids $B\cup \{x^l, x^r\}$,  where $x^l$ and $x^r$ are the tips associated with $P(x,j,M',B)$, and such that both $N_{T'}(x^l)$ and $N_{T'}(x^r)$ are connected in $P(x,j,M',B)$. 
Moreover, by choosing $P(x,j, M', B)$ minimal, we may assume that it consists of at most two components, and each such component contains either $N_{T'}(x^l)$ or $N_{T'}(x^r)$.\COMMENT{The purpose for this assumption is so that $P$ doesn't consist of the path structure we construct in the repatching lemma in addition to other (far away) vertices that contribute nothing. Such other vertices would have to be taken into account for no clumping arguments to come.}

Let $\iota\coloneqq|V_{\mathrm{sc}}|$ and let $x_1, \dots, x_\iota$ be an ordering of $V_{\mathrm{sc}}$.
For each $i \in [\iota]$, $j \in [K]$ and $k\in [2^{s+1}\Psi]$, by ranging over $i$ first, then $j$, and then $k$\COMMENT{for i do: {for j do: for k {do:{}}}}, we will iteratively fix a graph $P(x_i,j,k,M'_{i,j,k},B_{i,j,k})$ as above.
In particular, this graph will have an absorbing $\ell$-cube pair with tips $x^l_{i,j,k}$ and $x^r_{i,j,k}$ associated with it. 
After the graph $P(x_i,j,k,M'_{i,j,k},B_{i,j,k})$ is fixed, so are these tips.
Let $\mathcal{J}_{i,j,k}\coloneqq([i-1]\times[K]\times [2^{s+1}\Psi])\cup\{(i,j',k'):(j',k')\in[j-1]\times[2^{s+1}\Psi]\}\cup \{(i,j,k''):k''\in [k-1]\}$ and suppose that we have already fixed $P(x_{i'},j',k',M'_{i',j',k'},B_{i',j',k'})$ for all $(i',j',k')\in \mathcal{J}_{i,j,k}$ such that these $P(x_{i'},j',k',M'_{i',j',k'},B_{i',j',k'})$ are vertex-disjoint from each other and from the set $\{x^l_{i',j',k'}, x^r_{i',j',k'}: (i',j',k')\in \mathcal{J}_{i,j,k}\}$ of tips associated with all these $P(x_{i'},j',k',M'_{i',j',k'},B_{i',j',k'})$.
In order to fix $P(x_i,j,k,M'_{i,j,k},B_{i,j,k})$, we first define the sets $B_{i,j,k}$ and $M'_{i,j,k}$.
Let $M'_{i,j,k}$ be obtained from $M(A_j(x_i))$ as follows.
Remove all edges whose associated absorbing $\ell$-cube pair $(C^l,C^r)$ satisfies $(V(C^l)\cup V(C^r))\cap\{x^l_{i',j',k'},x^r_{i',j',k'}:(i',j',k')\in\mathcal{J}_{i,j,k}\}\neq\varnothing$\COMMENT{This says that one of the cubes has already been used to potentially absorb some other vertex, hence we cannot use it again.}.
Remove all edges $(\hat e,\hat e')\in M(A_j(x_i))$ such that $\{x_i+\hat e,x_i+\hat e'\}\cap\bigcup_{(i',j',k')\in\mathcal{J}_{i,j,k}} V(P(x_{i'},j',k',M'_{i',j',k'},B_{i',j',k'}))\neq\varnothing$ too\COMMENT{By removing these pairs, we are guaranteed that the tree after each iteration does not change the sets $B(y,z)$ that we want to repatch.
Note that this is not needed, but otherwise we would need to explicitly update the tree and these sets in each iteration, whereas now we do it all in one step at the end.}.
Note that, by \eqref{equa:step7.1} and \ref{itm:scant}, it follows that $|M'_{i,j,k}|\geq n/\ell^3$.\COMMENT{Consider first the first removal of edges.
For each triple $(i',j',k')\in\mathcal{J}_{i,j,k}$, we have already fixed the tips of the absorbing $\ell$-cube pairs.
The number of triples that may affect the matching $M_{i,j,k}'$ is bounded by those pairs such that $x_{i'}$ is at distance at most $\ell+2\leq10\ell$ from $x_i$, of which there are at most $S'K2^{s+1}\Psi$ by \ref{itm:scant}.
Furthermore, each triple may affect at most two distinct edges of $M(A_j(x_i))$.
Thus, we have that $|M'_{i,j,k}|\geq 2n/\ell^3-2^{s+2}K\Psi S'$.\\
Now consider the second removal of edges.
Since each graph $P(x_{i'},j',k',M'_{i',j',k'},B_{i',j',k'})$ contains at most $21D/2$ vertices and consists of at most $2$ components, each of which contains either $N_{T'}(x^l_{i',j',k'})$ or $N_{T'}(x^r_{i',j',k'})$, all vertices $x_i$ affected by it are at distance at most $21D/2 + 2\leq10\ell$ from $x_{i'}$.
That is, for a fixed triple $(i,j,k)$, each path $P(x_{i'},j',k',M'_{i',j',k'},B_{i',j',k'})$ can delete at most $21D/2$ edges from $M(A_j(x_i))$, and there are at most $S'K2^{s+1}\Psi$ choices for $(i',j',k')$ such that the vertex $x_{i'}$ is `close' to $x_i$ (again by \ref{itm:scant}).\\
Thus, overall, $|M'_{i,j,k}|\geq 2n/\ell^3-2^{s+2}K\Psi S'-21\cdot2^s\Psi DKS'\geq n/\ell^3$.}
Let $B_{i,j,k}$ be the set of vertices $y \in B^{\ell/2}_I(x_i)\subseteq B^{10\ell}_I(x_i)$ such that at least one of the following holds:
\begin{enumerate}[label=$(\mathrm{P}\arabic*)$]
    \item\label{itm:bad7.1} there exists $(i',j',k')\in\mathcal{J}_{i,j,k}$ such that $y\in V(P(x_{i'},j',k',M'_{i',j',k'},B_{i',j',k'}))$;
    \item\label{itm:bad7.2} there exists $(i',j',k')\in\mathcal{J}_{i,j,k}$ such that $y$ lies in the absorbing $\ell$-cube pair associated with $P(x_{i'},j',k',M'_{i',j',k'},B_{i',j',k'})$.
    \COMMENT{Really all we need is that it avoids the tip and at least one other vertex from each cube associated with a previous iteration which also intersects the tree. This one extra vertex will ensure that the tree still intersects the cube after repatching, so that the absorbing structure we're putting in place will in fact be taken into the large cycle we construct.}
\end{enumerate}
Note that $|B_{i,j,k}|<2^{s+\ell+3}\Psi KS'$ by \ref{itm:scant}\COMMENT{$2^{s+1}\Psi KS'$ is a strict upper bound on the number of choices of $(i',j',k')$ such that $\dist(x_i,x_{i'})\leq2\ell$ (note that any vertex $x_{i'}$ which is further away cannot have any effect over $B_{i,j,k}$).
$21D/2$ is an upper bound for the order of each of the graphs $P(x_i',j',k',M'_{i',j',k'},B_{i',j',k'})$.
Finally, $2\cdot2^\ell$ is an upper bound for the order of each pair of cubes.
In conclusion, $|B_{i,j,k}|<2^{s+1}\Psi KS'(21D/2+2\cdot2^\ell)\leq \Psi KS'2^{\ell+s+3}$.}.
We then fix $P(x_i,j,k,M'_{i,j,k},B_{i,j,k})$ to be the graph guaranteed by our conditioning on $\mathcal{E}_7^*$.
Observe that, by the choice of $B_{i,j,k}$, we have that $P(x_i,j,k,M'_{i,j,k},B_{i,j,k})$ is vertex-disjoint from $\bigcup_{(i',j',k') \in \mathcal{J}_{i,j,k}}P(x_{i'},j',k',M'_{i',j',k'},B_{i',j',k'})$.
We denote by $(C^l(x_i,j,k),C^r(x_i,j,k))$ the absorbing $\ell$-cube pair for $x_i$ associated with $P(x_i,j,k,M'_{i,j,k},B_{i,j,k})$.
By the choice of $M'_{i,j,k}$, we have that 
\begin{enumerate}[label=$(\mathrm{CD})$]
    \item\label{itm:DisjAbsCubes} for all $(i',j',k')\in \mathcal{J}_{i,j,k}$, $C^l(x_i,j,k)$ and $C^r(x_i,j,k)$ are both vertex-disjoint from $C^l(x_{i'},j',k')$ and $C^r(x_{i'},j',k')$. 
\end{enumerate}

Let $\cC^\mathrm{sc}_1\coloneqq\{(C^l(x_i,j,k),C^r(x_i,j,k)):(i,j,k)\in[\iota]\times[K]\times [2^{s+1}\Psi]\}$.
Let $P'\coloneqq\{x^l_{i,j,k},x^r_{i,j,k}:(i,j,k)\in[\iota]\times[K]\times[2^{s+1}\Psi]\}$ and $P\coloneqq\bigcup_{i\in[\iota],j\in[K],k\in[2^{s+1}\Psi]}P(x_i,j,k,M'_{i,j,k},B_{i,j,k})$.
Recall that $P(x_i,j,k,M'_{i,j,k},B_{i,j,k})$ avoids the tips $x^l_{i,j,k}$ and $x^r_{i,j,k}$ associated with it.
It follows from this, \ref{itm:bad7.2}, and the definition of $M'_{i,j,k}$ that $P' \cap V(P) = \varnothing$.
Let $T'''\coloneqq T'[V(T')\setminus P'] \cup P$.  
Note that $T'''$ is connected by the definition of $\cE'(x,j, M')$.
Let $T''$\index{T3@$T''$} be a spanning tree of $T'''$.
By \ref{itm: ET1} and the fact that the graphs $P(x_i,j,k,M'_{i,j,k},B_{i,j,k})$ are vertex-disjoint and satisfy $|V(P(x_i,j,k,M'_{i,j,k},B_{i,j,k}))| < 21D/2$, it follows that
\begin{equation}\label{equa:T''bound}
    \Delta(T'')\leq 12D.
\end{equation}
Define the (new) reservoir $R' \coloneqq (R \cup P')\setminus V(P)$\index{R2@$R'$}.

At this point, for each $x \in V(I)\setminus V_\mathrm{sc}$ and each $i \in [K]$, we redefine the set $M(A_i(x))$ as follows.
\begin{enumerate}[label=$(\mathrm{AB}\arabic*)$]\setcounter{enumi}{3}
\item\label{itm:AB4} Let $M(A_i(x))$ retain only those edges whose associated absorbing $\ell$-cube pair $(C^l, C^r)$ satisfies that both $C^l$ and $C^r$ are vertex-disjoint from both cubes of all absorbing $\ell$-cube pairs of $\cC^\mathrm{sc}_1$ and both tips $x^l$ and $x^r$ satisfy that $x^l,x^r\in R\setminus V(P)\subseteq R'$.
\end{enumerate}
Note that, by \ref{itm:scant}, we have $|B^{\ell+1}_I(x)\cap V(P)|\leq21\cdot2^s\Psi DKS'$ and $|B^{\ell+1}_I(x)\cap V(\bigcup_{(C^l,C^r)\in\cC^\mathrm{sc}_1}(C^l\cup C^r))|\leq4\cdot2^{\ell+s}\Psi KS'$.
Combining this with \eqref{eqn:maix1} and \ref{itm:abss2}, it follows that 
\begin{equation}\label{eqn:maix2}
    |M(A_i(x))| \geq 4n/\ell^3-(21D+4\cdot2^\ell)2^s\Psi KS' > n/\ell^3.
\end{equation}

\textbf{Step~9: Fixing a collection of absorbing $\boldsymbol{\ell}$-cube pairs for the vertices in non-scant molecules.}
At this point, we still do not know which vertices will need to be absorbed eventually into an almost spanning cycle, but we can already determine the vertices in $I$ whose clones the vertices to be absorbed will be (the reason for this will be apparent later, see Step~13).
Recall that $\cC'$ and $\cC''$ were defined in Step~6.
Let $\cC'''\coloneqq \{C\in \cC': V(C)\cap V(T'')\neq \varnothing\}$\index{Cw4@$\cC'''$} and let $V_\mathrm{abs}\coloneqq V(I)\setminus \bigcup_{C\in \cC'''} V(C)$\index{Vabs@$V_\mathrm{abs}$}.
We will now fix a collection of absorbing $\ell$-cube pairs for all vertices in each vertex molecule $\mathcal{M}_x$ with $x\in V_\mathrm{abs}\setminus V_\mathrm{sc}$.

First, recall from \ref{itm:scant} that, for all $x \in V(I)$, we have that $|B^{10\ell}_I(x)\cap V_{\mathrm{sc}}| \leq S'$.
Thus, in constructing $T''$, we removed at most $2^{s+2}\Psi KS'$ vertices in $B^\ell_I(x)$ from $T'$.
Therefore, it follows from \ref{itm: ET2} that, for all $x \in V(I)$, we have 
\begin{equation}\label{eqn:moveC}
|B^{2}_{I}(x)\setminus V(T'')| \leq 2n^{3/4}.
\end{equation}
For all $x \in \bigcup_{C\in\cC''}V(C)$, we claim that 
\begin{eqnarray}\label{eqn:movedCi}
|N_I(x)\cap V(T'')\cap\bigcup_{C\in\cC'}V(C)|\geq(1-2^{1-\ell-5s})n.
\end{eqnarray}
To see that this holds, combine \ref{itm:nib1}, \eqref{eqn:moveC} and the definition of bondlessly surrounded molecules. \COMMENT{We combine \ref{itm:nib1} and \eqref{eqn:moveC}, but now with the definition of bondlessly surrounded.
It follows that at least $(1-\delta)n-n/2^{\ell+5s}-2n^{3/4}\geq(1-2/2^{\ell+5s})n$ of the neighbours of $x$ in $I$ lie simultaneously in bonded cubes and the tree $T''$.}

Recall also the definition of $\mathfrak{M}(x)$ from Step~3.

\begin{claim}\label{claim:rainbowapplnew}
For each $x\in V_\mathrm{abs}\setminus V_\mathrm{sc}$ and each $e\in \mathfrak{M}(x)$, there exists a set $\mathcal{C}_1^\mathrm{abs}(e)$ of $2^{s+1}\Psi$ absorbing $\ell$-cube pairs $(C_{k}^l(e),C_{k}^r(e))\subseteq I$, for $k\in[2^{s+1}\Psi]$, which satisfies the following:
\begin{enumerate}[label=$(\mathrm{\roman*})$]
    \item for all $x\in V_\mathrm{abs}\setminus V_\mathrm{sc}$, $e\in \mathfrak{M}(x)$ and $k\in[2^{s+1}\Psi]$, the absorbing $\ell$-cube pair $(C_{k}^l(e),C_{k}^r(e))$ is associated with some edge in $M(A_j(x))$, for some $j \in [K]$, and
    \item for all $x,x'\in V_\mathrm{abs}\setminus V_\mathrm{sc}$, all $e\in \mathfrak{M}(x)$ and $e'\in \mathfrak{M}(x')$, and all $k,k'\in [2^{s+1}\Psi]$ with $(x,e,k)\neq (x',e',k')$, the absorbing $\ell$-cube pairs $(C_{k}^l(e),C_{k}^r(e))$ are vertex-disjoint (except for $x$ in the case when $x=x'$).
\end{enumerate}
\end{claim}

\begin{claimproof}
Let $\cV \coloneqq\bigcup_{x\in V_\mathrm{abs}\setminus V_\mathrm{sc}}\mathfrak{M}(x)$.
Let $K'\coloneqq|\cV|$, and let $f_1, \dots, f_{K'}$ be an ordering of the edges in $\cV$.
Given any $i \in [K']$, the edge $f_{i}$ corresponds to a pair $(x,j(i))$ (in the sense that $A_{j(i)}(x)=N(f_i)$, see Step~4), where $x\in V_\mathrm{abs}\setminus V_\mathrm{sc}$ and $j(i)\in [K]$.
Let $\mathfrak{C}_i$ be the collection of at least $n/\ell^3$ absorbing $\ell$-cube pairs for $x$ in $I$ guaranteed by \eqref{eqn:maix2}. 
In particular, each of these absorbing $\ell$-cube pairs $(C^l, C^r)$ is associated with an edge of $M(A_{j(i)}(x))$ and, by \ref{itm:abss1}, satisfies $C^l, C^r \in \cC''$.

Let $\cH$ be the $2^{s+1}\Psi K'$-edge-coloured auxiliary multigraph with $V(\mathcal{H})\coloneqq\mathcal{C}''$, which contains an edge between $C$ and $C'$ of colour $(i,k)\in [K']\times [2^{s+1}\Psi]$ whenever $(C,C') \in \mathfrak{C}_i$ or $(C',C)\in \mathfrak{C}_i$.
In particular, $\cH$ contains at least $n/\ell^3$ edges of each colour.
We now bound $\Delta(\cH)$.
Consider any $C \in V(\cH)$.
Note that, for each edge $e$ of $\cH$ incident to $C$, there exists some $x=x(e) \in V_\mathrm{abs}\setminus V_\mathrm{sc}$ such that $C$ together with some other cube $C' \in V(\cH)$ forms an absorbing $\ell$-cube pair for $x$.
In particular, $x$ must be adjacent to $C$ in $I$.
Moreover, if $e$ has colour $(i,k)$, then $f_i\in \mathfrak{M}(x)$ (and it has corresponding pair $(x,j(i))$ for some $j(i)\in[K]$).
Since $f_i \in \mathfrak{M}(x)$ and $|\mathfrak{M}(x)|\leq \binom{2^{s}}{2}$, it follows that each vertex $y$ which is adjacent to $C$ in $I$ can play the role of $x$ for at most $2^{s+1}\Psi\cdot 2^{2s}$ edges of $\cH$ incident to $C$.
Thus, $d_{\cH}(C)$ is at most $2^{s+1}\Psi\cdot 2^{2s}$ times the number of vertices $y\in V_\mathrm{abs}\setminus V_\mathrm{sc}$ which are adjacent to $C$ in $I$.
Recall that $V_{\mathrm{abs}} = V(I) \setminus \bigcup_{C\in \cC'''} V(C)$.
Together with \eqref{eqn:movedCi}, this implies that the number of vertices in $V_{\mathrm{abs}}$ which are adjacent to $C$ is at most $|C|n/{2^{\ell +5s -1}}$.
Thus, $d_{\cH}(C) \le 2^{s+1}\Psi 2^{2s}|C|n/{2^{\ell +5s -1}}\le n/\ell^4$. \COMMENT{Here we used that $\Phi= 12\ell$ and $\Psi= c\Phi\leq \ell \Phi$. }

Since each colour class has size at least $n/\ell^3$ and $\Delta(\mathcal{H})\leq n/\ell^4$, by \cref{lem: rainbow}, $\mathcal{H}$ contains a rainbow matching of size $2^{s+1}\Psi K'$.\COMMENT{We apply \cref{lem: rainbow} with $r= 2$, $m=n/\ell^3$, and note that $\Delta(\cH) < n/\ell^{4} < m/(6r)$.}
For each $(i,k)\in [K']\times [2^{s+1}\Psi]$, let $(C_{k}^l(f_i),C_{k}^r(f_i))\in\mathfrak{C}_i$ be the absorbing $\ell$-cube pair of colour $(i,k)$ in this rainbow matching.
\end{claimproof}

Recall that, for any $x\in V(I)$, each index $i\in[K]$ is given by a unique edge $e \in \mathfrak{M}(x)$ via the relation $N(e) = A_i(x)$.
For each $x\in V_\mathrm{abs}\setminus V_\mathrm{sc}$ and each $i\in[K]$, let $\cC^\mathrm{abs}_1(x,i)\coloneqq\cC^\mathrm{abs}_1(e)$, where $e$ is the unique edge given by the relation above, be the set of absorbing $\ell$-cube pairs guaranteed by \cref{claim:rainbowapplnew}.
Similarly, for each $k\in[2^{s+1}\Psi]$, let $(C_{k}^l(x,i),C_{k}^r(x,i))\coloneqq(C_{k}^l(e),C_{k}^r(e))$.\COMMENT{This change of notation is not needed here or in the next paragraph, as taking union over all $i$ is the same as taking unions over all $e$. However, this notation is more convenient later for the list of properties \ref{itm:C6} and their proof, as we want to have consistent notation there.}\\

Let $G\coloneqq\bigcup_{i=1}^{7} G_i$.
For each $x\in V_\mathrm{abs}\setminus V_\mathrm{sc}$ and each $i\in [K]$, let $G^{*}(x,i)\subseteq I$ be the graph consisting of all edges between the left absorber tip and third absorber vertex of every absorbing $\ell$-cube pair in $\cC_1^\mathrm{abs}(x,i)$.
Let $G^{\bullet}\subseteq I$ be the graph consisting of all edges between the left absorber tip and third absorber vertex of every absorbing $\ell$-cube pair in $\cC_1^\mathrm{sc}$.
Let $G^{*}\coloneqq G^\bullet\cup\bigcup_{x\in V_\mathrm{abs}\setminus V_\mathrm{sc}}\bigcup_{i\in[K]}G^*(x,i)\subseteq I$.
Recall that, given any graph $\cG\subseteq I$, for each layer $L$, we denote by $\cG_L$ the clone of $\cG$ in $L$.
Let $G_4^*\coloneqq G_4\cap\bigcup_{i=1}^{2^s}G^{*}_{L_i}$.
Furthermore, let $G_5^*\subseteq G_5$ consist of all edges of $G_5$ which have endpoints in different layers.
We let $G'\subseteq G$ be the spanning subgraph with edge set 
\[E(G')\coloneqq E(G_4^{*})\cup E(G_5^*)\cup\bigcup_{C\in \cC'} E(\cM_{C})\cup\bigcup_{i=1}^{2^s} E(T''_{L_i}).\]
Note that, using \eqref{equa:T''bound}, we have that $\Delta(G')\leq\Phi$. \COMMENT{every vertex in $G_4^{*}$ has degree at most $1$, each vertex in $G_5^*$ has degree at most $s=10\ell$, and since the molecules are vertex disjoint, each vertex in a molecule has degree at most $\ell$.
Putting everything together, we conclude that $\Delta(G')\leq 1+s+\ell+12D\leq12\ell=\Phi$.}

Now, let $F\subseteq \cQ^{n}$ be any graph with $\Delta(F)\leq \Psi $.
Recall that we denote by $F_{I}\subseteq I$ the graph which contains every edge $\{x,y\}\in E(I)$ such that there exists an edge $e=\{x',y'\}\in E(F)$ with $x'\in \cM_{x}$ and $y'\in \cM_{y}$.

Note that $T'' \subseteq I(G')$, $R'\subseteq V(I)$, and $C\subseteq I(G')$ for every $C\in \cC'$.
Recall the definitions of $\cC''$ from Step~6 and $\cC'''$ from Step~9.
Combining all the previous steps, we claim that the following hold (conditioned on the events $\cE_1^*,\ldots,\cE_7^*$, which occur a.a.s.).
\begin{enumerate}[label=$(\mathrm{C}\arabic*)$]
    \item\label{itm:C2} $\Delta(T'')\leq 12D$.
    \item\label{itm:C1} Any vertex $x\in R'\cap V(T'')$ is a leaf of $T''$.
    Furthermore, if $x\in R'\cap V(T'')$, then its unique neighbour $x'$ in $T''$ satisfies that $x'\in Z(x)$ (where $Z(x)$ is as defined in Step~7).
    \item\label{itm:C5bis} For all $x \in V(I)$, we have that $|N_I(x)\cap V(T'')\cap\bigcup_{C\in\cC''}V(C)| \geq (1-2/\ell^4)n$.
    \item\label{itm:C3} For each $x \in V_{\mathrm{sc}}$ and $i \in [K]$, there is an absorbing $\ell$-cube pair $(C^l(x,i),C^r(x,i))$ for $x$ in $I$, which is associated with some edge $e \in M(A_i(x))$.
    In particular, $(C^l(x,i),C^r(x,i))$ is as described in \ref{itm:abss2} (recall also \ref{itm:abss1}), that is, there are two absorbing $\ell$-cube pairs $(C^l_1(x,i),C^r_1(x,i))$ and $(C^l_2(x,i),C^r_2(x,i))$ in $H\cup G$, associated with $e\in M(A_i(x))$, for the clones $x_1$ and $x_2$ of $x$ which correspond to $(x,i)$.
    Additionally, each of these absorbing $\ell$-cube pairs $(C^l(x,i),C^r(x,i))$ satisfies the following:
    \begin{enumerate}[label=$(\mathrm{C}4.\arabic*)$]
        \item\label{itm:C30} $(C^l_1(x,i), C^r_1(x,i)) \cup (C^l_2(x,i), C^r_2(x,i)) - V(\cM_x) \subseteq G'$;
        \item\label{itm:C31} the tips $x^{l}$ of $C^l(x,i)$ and $x^{r}$ of $C^r(x,i)$ lie in $R'\setminus V(T'')$, and $\{x,x^{l}\},\{x,x^{r}\}\notin E(F_I)$;
        in particular, the tips $x_1^l,x_1^r$ of $(C^l_1(x,i),C^r_1(x,i))$ and $x_2^l,x_2^r$ of $(C^l_2(x,i),C^r_2(x,i))$ satisfy that $\{x_1,x_1^{l}\},\{x_1,x_1^{r}\},\{x_2,x_2^{l}\},\{x_2,x_2^{r}\}\in E((H\cup G)\setminus F)$;
        \item\label{itm:C325} $C^l(x,i),C^r(x,i)\in\mathcal{C}''\cap\mathcal{C}'''$, and
        \item\label{itm:C33} for any $x' \in V_{\mathrm{sc}}$ and $i' \in [K]$ with $(x',i')\neq(x,i)$ we have that $C^l(x,i)$, $C^r(x,i)$, $C^l(x',i')$ and $C^r(x',i')$ are vertex-disjoint.
    \end{enumerate}
    Let $\cC^\mathrm{sc}$ denote the collection of these absorbing $\ell$-cube pairs.
    \item\label{itm:C6} For each $x \in V_{\mathrm{abs}}\setminus V_{\mathrm{sc}}$ and $i \in [K]$, there is an absorbing $\ell$-cube pair $(C^l(x,i),C^r(x,i))$ for $x$ in $I$, which is associated with some edge in $M(A_i(x))$.
    In particular, $(C^l(x,i),C^r(x,i))$ is as described in \ref{itm:abss2}, that is, there are two absorbing $\ell$-cube pairs $(C^l_1(x,i),C^r_1(x,i))$ and $(C^l_2(x,i),C^r_2(x,i))$ in $H\cup G$, associated with $e\in M(A_i(x))$, for the clones $x_1$ and $x_2$ of $x$ which correspond to $(x,i)$.
    Moreover, each of these absorbing $\ell$-cube pairs $(C^l(x,i),C^r(x,i))$ satisfies the following:
    \begin{enumerate}[label=$(\mathrm{C}5.\arabic*)$]
        \item\label{itm:C50} $(C^l_1(x,i), C^r_1(x,i)) \cup (C^l_2(x,i), C^r_2(x,i)) - V(\cM_x) \subseteq G'$;
        \item\label{itm:C52} the tips $x_i^{l}$ of $C^l(x,i)$ and $x_i^{r}$ of $C^r(x,i)$ lie in $R'$, and $\{x,x_i^{l}\},\{x,x_i^{r}\}\notin E(F_I)$; in particular, the tips $x_1^l,x_1^r$ of $(C^l_1(x,i),C^r_1(x,i))$ and $x_2^l,x_2^r$ of $(C^l_2(x,i),C^r_2(x,i))$ satisfy that $\{x_1,x_1^{l}\},\{x_1,x_1^{r}\},\{x_2,x_2^{l}\},\{x_2,x_2^{r}\}\in E((H\cup G)\setminus F)$;
        \item\label{itm:C54} $C^l(x,i),C^r(x,i)\in\mathcal{C}''\cap\mathcal{C}'''$;
        \item\label{itm:C55} for any $x' \in V_\mathrm{abs}\setminus V_{\mathrm{sc}}$ and $i' \in [K]$ with $(x',i')\neq(x,i)$ we have that $C^l(x,i)$, $C^r(x,i)$, $C^l(x',i')$ and $C^r(x',i')$ are vertex-disjoint, and
        \item\label{itm:C51} both $C^l(x,i)$ and $C^r(x,i)$ are vertex-disjoint from all cubes of absorbing $\ell$-cube pairs in $\cC^\mathrm{sc}$.
    \end{enumerate}
    Let $\cC^\mathrm{\neg sc}$ denote the collection of these absorbing $\ell$-cube pairs.
\end{enumerate}
Indeed, \ref{itm:C2} is given in \eqref{equa:T''bound}.
\ref{itm:C1} holds by \ref{itm:ET3} and the fact that $P'\cap V(T'')=\varnothing$.
\ref{itm:C5bis} follows by combining \ref{itm:nib1}, the conditioning on $\mathcal{E}_5^*$, and \eqref{eqn:moveC}\COMMENT{To see this, combine \ref{itm:nib1}, the conditioning on $\mathcal{E}_5^*$, and \eqref{eqn:moveC} to see that at least $(1-\delta)n-n/\ell^4-n^{1/3}-2n^{3/4}\geq(1-2/\ell^4)n$ of the neighbours of $x$ in $I$ lie simultaneously in bonded cubes which are not bondlessly surrounded and the tree $T''$.}.
\ref{itm:C3} follows from the construction of $P$ and $T''$ in Step~8.
Indeed, for each $x\in V_\mathrm{sc}$ and $i\in[K]$, consider the collection of absorbing $\ell$-cube pairs $\{(C^l(x,i,k),C^r(x,i,k))\}_{k\in[2^{s+1}\Psi]}$ defined in Step~8.
Since $\Delta(F)\leq\Psi$, it follows that $d_{F_I}(x)\leq2^s\Psi$, and thus there must exist some absorbing $\ell$-cube pair in this collection such that the edges joining its tips to $x$ do not belong to $F_I$.
Fix one such absorbing $\ell$-cube pair and call it $(C^l(x,i),C^r(x,i))$.
Then, \ref{itm:C30} holds by the definition of $G'$ combined with \ref{itm:abss2}, and \ref{itm:C31} holds by the definition of $R'$ and $T''$ combined with \ref{itm:abss2}, while \ref{itm:C33} holds by \ref{itm:DisjAbsCubes}.
On the other hand, \ref{itm:C325} follows because of the definition of the set $M(A_i(x))$ in \ref{itm:abss1} and \ref{itm:AB3}.
Finally, consider \ref{itm:C6}.
For each $x \in V_\mathrm{abs}\setminus V_{\mathrm{sc}}$ and $i \in [K]$, consider the collection $\cC^\mathrm{abs}_1(x,i)$ of $2^{s+1}\Psi$ absorbing $\ell$-cube pairs for $x$ in $I$ guaranteed by \cref{claim:rainbowapplnew}.
For each of these absorbing $\ell$-cube pairs we have that \ref{itm:C54} holds by \ref{itm:abss1}, \ref{itm:AB3} and the fact that, by \ref{itm:AB4}, their intersection with $T''$ contains their intersection with $T'$\COMMENT{The reason for this is as follows. We did not repatch any absorbing $\ell$-cube pairs for non scant vertices (and we had already removed all clashes with potential absorbing $\ell$ cube pairs between scant and non scant vertices in \ref{eqn:maix1}). So all of these cubes still have their original intersections with $T'$.}.
Similarly, \ref{itm:C55} holds by \cref{claim:rainbowapplnew}, and \ref{itm:C51} holds because of \ref{itm:AB4}.
Finally, note that $\Delta(F_I) \le 2^s\Psi$.
It follows that there exists a choice of $(C^l(x,i),C^r(x,i))\in\cC^\mathrm{abs}_1(x,i)$ such that $\{x,x_i^{l}\},\{x,x_i^{r}\}\notin E(F_I)$.
Then, \ref{itm:C50} and \ref{itm:C52} hold by the definition of $G'$, \ref{itm:abss2} and \ref{itm:AB4}.\\

\textbf{Step~10: Constructing auxiliary trees $\boldsymbol{T^*}$ and $\boldsymbol{\tau_0}$.}
From this point on, every step will be deterministic. 
Let $T^*$ be obtained from $T''$ by removing all leaves of $T''$ which lie in $R'$.

We will now construct an auxiliary tree $\tau_0$\index{tau0@$\tau_0$}, which will be used in the construction of an almost spanning cycle.
We start by defining an auxiliary multigraph $\Gamma'$ as follows.
First, let $\Gamma_1\coloneqq T^*\cup\bigcup_{C\in\mathcal{C}'}C$.\COMMENT{That is, no bondless cubes.} 
(Recall that $\cC'$ is the collection of all $C\in\cC$ for which $\cM_C$ is bonded in $G_5$, see Step~6.)
Let $\Gamma_2$ be the graph obtained by iteratively removing all leaves from $\Gamma_1$ until all vertices have degree at least $2$.
Observe that, after this is achieved, the resulting graph still contains all cubes $C \in \cC'$.
Let $\Gamma_3$ be obtained from $\Gamma_2$ by removing all connected components which consist of a single cube $C \in \cC'$.
Now, let $\Gamma'$ be the multigraph obtained by contracting each cube $C \in \cC'$ such that $C\subseteq\Gamma_3$ into a single vertex.
We refer to the vertices resulting from contracting such cubes as \emph{atomic vertices}, and to the remaining vertices in $\Gamma'$ as \emph{inner tree vertices}.
Given $C \in \cC$ and $j \in [2^s]$, we call $\cA = \cM_C \cap L_j$ an \emph{atom}.
We continue to identify each inner tree vertex $v$ with the vertex $v \in V(I)$ from which it originated in $\Gamma_1$. 
Observe that $\Gamma'$ is connected\COMMENT{If there is an edge joining two atoms, this edge belongs to the tree, so it belongs to the unique connected component which contains the tree. Therefore, in $\Gamma_2$, all connected components must either contain the tree or be `isolated' atoms.}, and \ref{itm:C2} implies that
\begin{equation}\label{equa:step9.1}
\begin{minipage}[c]{0.7\textwidth}
$d_{\Gamma'}(v)\leq12D$ for all inner tree vertices, and $\Delta(\Gamma')\leq 12\cdot2^\ell D$.
\end{minipage}\ignorespacesafterend 
\end{equation} 
Given an atomic vertex $v\in V(\Gamma')$, let $C(v) \in \cC$ be the cube which was contracted to $v$ in the construction of $\Gamma'$, and let \index{Mofv@$\mathcal{M}(v)$}$\mathcal{M}(v)\coloneqq \cM_{C(v)}$.
Furthermore, for each $j\in[2^s]$, let\COMMENT{so $\cA_j(v)$ is the atom obtained by intersecting $\mathcal{M}(v)$ with the $j$-th layer.}
\[\index{Ajofv@$\mathcal{A}_j(v)$}\mathcal{A}_j(v) \coloneqq \mathcal{M}(v) \cap L_j.\]
Similarly, for any $v\in V(\Gamma')$ which is an inner tree vertex, we define $\mathcal{M}(v)\coloneqq \cM_v$. 
Observe that every edge $e\in E(\Gamma')$ corresponds to a unique edge $e' \in I(G')$.
We say that $e$ \emph{originates} from $e'$.
We denote by $D(e)\in\mathcal{D}(I)$ the direction of $e'$ in $I$\COMMENT{Because this is a direction in $I$, it is a vector with $n-s$ coordinates.
We will later add this direction to vertices on $V(\cQ^n)$, which have $n$ coordinates.
So here we have a small problem with notation.
Somehow, we need to say that we look at this direction as the direction in the $n$-cube.}.
By abusing notation, we will sometimes also view $D(e)$ as a direction in $\cQ^n$.

Next, we fix any atomic vertex $v_0\in V(\Gamma')$.
We define a (simple) auxiliary labelled rooted tree $\tau_0=\tau_0(v_0)$ by performing a depth-first search on $\Gamma'$ rooted at $v_0$ and then iteratively removing all leaves which are inner tree vertices.
This results in a tree $\tau_0$ rooted at an atomic vertex $v_0$ and all whose leaves are atomic vertices.
Let $m\coloneqq|V(\tau_0)|-1$, and let the vertices of $\tau_0$ be labelled as $v_0,v_1,\ldots,v_m$, with the labelling given by the order in which each vertex is explored by the depth-first search performed on $\Gamma'$.
For each $i\in[m]$, we define $\tau_i$ as the maximal subtree of $\tau_0$ which contains $v_i$ and all whose vertices have labels which are at least as large as $i$.
Given any vertex $x\in V(I)$, we say that $x$ is \emph{represented} in $\tau_0$ if $x\in V(\tau_0)$ or there exists some atomic vertex $v\in V(\tau_0)$ such that $x\in V(C(v))$.
Similarly, we say that a cube $C\in\mathcal{C}$ is \emph{represented} in $\tau_0$ if there exists an atomic vertex $v\in V(\tau_0)$ such that $C=C(v)$.
We will sometimes also say that $\mathcal{M}_x$ or $\mathcal{M}_C$ are represented in $\tau_0$, respectively.

The tree $\tau_0$ will be the backbone upon which we construct our long cycle.
First, we need to set up some more notation.
For each $i\in[m]_0$, let $p_i\coloneqq d_{\tau_i}(v_i)$ and let $N_{\tau_i}(v_i)=\{u_1^i,\ldots,u_{p_i}^i\}$\COMMENT{Maintaining the order from the labelling of $\tau_0$, say, but I think this is not important.}.
It follows from \eqref{equa:step9.1} that
\begin{equation}\label{equa:step9.2}
\begin{minipage}[c]{0.7\textwidth}
$p_i \leq 12D-1$ if $v_i$ is an inner tree vertex\COMMENT{We need the number of layers in each slice to be larger than this (say, at least $10$ times).}, and $\Delta(\tau_0)\leq12\cdot2^\ell D$\COMMENT{We need the number of slices in each molecule to be larger than this (say, at least $10$ times).}.
\end{minipage}\ignorespacesafterend 
\end{equation} 
For each $i\in[m]_0$ and $k\in[p_i]$, let $e_k^i\coloneqq\{v_i,u_k^i\}$, let $f_k^i\coloneqq D(e_k^i)$, and let $j_k^i$ be the label of $u_k^i$ in $\tau_0$, that is, $u_k^i=v_{j_k^i}$.
For any $k\in[p_i]$, we will sometimes refer to $i$ as the \emph{parent index} of $j_k^i$.
Furthermore, for each $i\in[m]_0$ such that $v_i$ is an atomic vertex, and for each $k\in[p_i]$, consider the unique edge $e_{i,k}' \in E(I(G'))$ so that $e_k^i$ originates from $e_{i,k}'$, and let \index{nuik@$\nu_k^i$}$\nu_k^i$ be the endpoint of $e_{i,k}'$ in $C(v_i)$.
Finally, for each $i\in[m]_0$, we define a parameter $\Delta(v_i)$ recursively by setting
\begin{equation}\label{equa:Delta}
\Delta(v_i)\coloneqq
\begin{cases}
0&\text{ if }v_i\text{ is an atomic vertex which is a leaf of }\tau_0,\\
\sum_{k=1}^{p_i}\Delta(u_k^i)&\text{ if }v_i\text{ is an atomic vertex which is not a leaf of }\tau_0,\\
p_i+1+\sum_{k=1}^{p_i}\Delta(u_k^i)&\text{ if }v_i\text{ is an inner tree vertex.}
\end{cases}
\end{equation}
This parameter $\Delta(v_i)$ will be used to keep track of parities throughout the following steps.
Note that $\Delta(v_i)$ counts the number of times a depth first search of $\tau_i$ (starting and ending at $v_i$) traverses an inner tree vertex.

Consider the partition of all molecules into slices of size $q$ introduced at the beginning of Step~3, where $q$ is as defined in \eqref{qandt}.
Given any $v\in V(\tau_0)$, we denote the slices of its molecule by $\mathcal{M}_1(v),\ldots,\mathcal{M}_{t}(v)$, where $t$ is as defined in \eqref{qandt}.
Thus, for each $i \in [t]$ we have that \index{Mofv2@$\mathcal{M}_i(v)$}$\cM_i(v) = \bigcup_{j=(i-1)q+1}^{iq}\cA_j(v)$.
For each $i\in [m]_0$, we are going to assign an \emph{input slice} $\mathcal{M}_{b(i)}(v_i)$ to each vertex $v_i$.
We do so by recursively assigning an \emph{input index} $b(i)\in[t]$ to each $i\in[m]_0$.
We begin by letting $b(0)\coloneqq1$.
Then, for each $i\in[m]_0$ and each $k\in[p_i]$, we set
\[
b(j_k^i)\coloneqq
\begin{cases}
b(i)&\text{ if }v_i\text{ is an inner tree vertex,}\\
b(i)+k-1\pmod{t}&\text{ if }v_i\text{ is an atomic vertex.}
\end{cases}
\]
Note that the bound on $\Delta(\tau_0)$ in \eqref{equa:step9.2} and the definition of $t$ in \eqref{qandt} imply that $b(j^i_k)\neq b(j^i_{k'})$ whenever $v_i$ is an atomic vertex and $k\neq k'$.\\

\textbf{Step~11: Finding an external skeleton for $\boldsymbol{T^*}$.}
Our next goal is to find an almost spanning cycle in $G'$ by using $\tau_0$ to explore different molecules in a given order.
For this, we are going to generate a \emph{skeleton}; this will be an ordered list of vertices which we will denote by $\mathcal{L}$.
In order to construct $\mathcal{L}$, we will construct disjoint \emph{partial skeletons} $\mathcal{L}_i$ and $\hat{\mathcal{L}}_i$ for all $i\in[m]$ in an inductive way.
Each of these skeletons will start and end in the input slice for the vertex $v_i$ which is being considered.
These partial skeletons will depend on the starting and ending vertices of $\mathcal{M}_{b(i)}(v_i)$ which are provided for each of them.
Therefore, given two distinct starting vertices $x,\hat{x}\in V(\mathcal{M}_{b(i)}(v_i))$ and two distinct ending vertices $y,\hat y\in V(\mathcal{M}_{b(i)}(v_i))$, we will denote the partial skeletons by $\mathcal{L}_i(x,y)$ and $\hat{\mathcal{L}}_i(\hat x,\hat y)$, respectively.

The first step in the construction of $\mathcal{L}$ is to construct a set of vertices $L^\bullet$, to which we will refer as an \emph{external skeleton}, and for which we will in turn construct \emph{partial external skeletons} in an inductive way.
The external skeleton will be essential in determining which vertices will not be covered by the almost spanning cycle, and hence need to be absorbed.
Roughly speaking, the external skeleton will contain
\begin{enumerate}[label=$(\mathrm{\roman*})$]
    \item all vertices where the almost spanning cycle enters and leaves each cube molecule represented in $\tau_0$, and
    \item all vertices which are not in cube molecules and are needed to connect cube molecules to each other (that is, some clones of inner tree vertices).
\end{enumerate} 
On the other hand, all vertices in a vertex molecule represented in $\tau_0$ by an inner tree vertex which do not belong to the external skeleton will have to be absorbed.

For each  $i\in[m]$, given the starting and ending vertices $x,y,\hat x,\hat y\in V(\mathcal{M}_{b(i)}(v_i))$ for $\mathcal{L}_i(x,y)$ and $\hat{\mathcal{L}}_i(\hat x,\hat y)$, we will denote the corresponding partial external skeleton by $L_i^\bullet(x,y,\hat x,\hat y)$.

The external skeleton is constructed recursively.
The partial external skeletons are the result of each recursive step, assuming that the starting and ending points have been defined.
Roughly speaking, for each $i\in[m]$, we will define partial external skeletons for any possible starting and ending vertices.
The starting and ending vertices which we actually use are then fixed by the partial external skeleton whose index is the parent of $i$.
Ultimately, all of them will be fixed when defining the external skeleton $L^\bullet$.

Let $\mathcal{M}_{\mathrm{Res}}\subseteq V(\cQ^n)$\index{MRes@$\mathcal{M}_{\mathrm{Res}}$} be the union of all the clones of $R'$.
We will construct an external skeleton $L^\bullet$ which satisfies the following properties:
\begin{enumerate}[label=$(\mathrm{ES}\arabic*)$]
    \item\label{item:preskprop1} For each $i\in[m]$ such that $v_i$ is an inner tree vertex, $L^\bullet\cap V(\mathcal{M}_{b(i)}(v_i))$ contains exactly $2p_i+2$ vertices, half of them of each parity, and $L^\bullet\cap(V(\mathcal{M}(v_i))\setminus V(\mathcal{M}_{b(i)}(v_i)))=\varnothing$.
    \item\label{item:preskprop2} For each $i\in[m]$ such that $v_i$ is an atomic vertex, $L^\bullet\cap V(\mathcal{M}(v_i))$ contains exactly $4p_i+4$ vertices.
    If $v_i$ is not a leaf of $\tau_0$, eight of these vertices (four of each parity) lie in $V(\mathcal{M}_{b(i)}(v_i))$, and four (two of each parity) lie in each $V(\mathcal{M}_{b(i)+k}(v_i))$ with $k\in[p_i-1]$.
    If $v_i$ is a leaf, then all four of these vertices lie in $V(\mathcal{M}_{b(i)}(v_i))$.
    \item\label{item:preskprop3} $L^\bullet\cap V(\mathcal{M}(v_0))$ contains exactly $4p_0$ vertices, four of them (two of each parity) lying in each $V(\mathcal{M}_k(v_0))$ with $k\in[p_0]$.
    \item\label{item:preskprop4} The sets described in \ref{item:preskprop1}--\ref{item:preskprop3} partition $L^\bullet$.
    \item\label{item:preskprop5} $L^\bullet\cap\mathcal{M}_{\mathrm{Res}}=\varnothing$.
\end{enumerate}
We now proceed to define the partial external skeletons formally.
The construction proceeds by induction on $i\in[m]$ in decreasing order, starting with $i=m$.
We define a \emph{valid connection sequence $(x^i,y^i,\hat x^i,\hat y^i)$ for} $v_i$ as any set of distinct vertices $x^i,y^i,\hat x^i,\hat y^i\in V(\mathcal{M}_{b(i)}(v_i))$ which satisfy the following:
\begin{enumerate}[label=($\mathrm{V}\arabic*$)]
    \item\label{itm:CS1} $x^i\neqp y^i$ if $\Delta(v_i)$ is even, and $x^i\eqp y^i$ otherwise; 
    \item\label{itm:CS2} $\hat x^i\neqp x^i$, and 
    \item\label{itm:CS3} $\hat y^i\neqp y^i$.
\end{enumerate}
Given any valid connection sequence $(x^i,y^i,\hat x^i,\hat y^i)$, we will refer to $x^i$ and $\hat x^i$ as \emph{starting} vertices, and to $y^i$ and $\hat y^i$ as \emph{ending} vertices.
Throughout the construction ahead, observe that, every time we use a partial external skeleton to build a larger one, its starting and ending vertices form a valid connection sequence by construction.
The vertices $x^i$, $y^i$, etc.~will be part of $\mathcal{L}_i(x^i,y^i)$, and the vertices $\hat x^i$, $\hat y^i$, etc.~will be part of $\hat{\mathcal{L}}_i(\hat x^i,\hat y^i)$.
The vertices $x^i$, $y^i$, $\hat x^i$, $\hat y^i$ will be used by the skeleton to move from the molecule represented by $v_i$ in $\tau_0$ to  the molecule represented by its parent.
Given these vertices, the following construction provides the vertices $w^i_k$ and $\hat w^i_k$ (as well as $z^i_k$ and $\hat z^i_k$, if applicable) which are used to move to molecules represented by the children of $v_i$.
Given any vertices $(x,y,\hat x,\hat y)$ in $\cQ^n$ and any direction $f\in\mathcal{D}(\cQ^n)$, we write $f+(x,y,\hat x,\hat y)=(f+x,f+y,f+\hat x,f+\hat y)$.

Now suppose that $i\in[m]$ and that, for each $i'\in[m]\setminus[i]$, we have already constructed a partial external skeleton $L_{i'}^\bullet(x^{i'},y^{i'},\hat x^{i'},\hat y^{i'})$ for $v_{i'}$ and every valid connection sequence $(x^{i'},y^{i'},\hat x^{i'},\hat y^{i'})$ for $v_{i'}$.
We will now construct a partial external skeleton for $v_i$ and every valid connection sequence for $v_i$.
We consider several cases.

\textbf{Case 1:} $v_i\in V(\tau_0)$ is a leaf of $\tau_0$.
Assume that $(x^i,y^i,\hat x^i,\hat y^i)$ is a valid connection sequence for $v_i$.
Then, the partial external skeleton for this connection sequence is given by $L_i^\bullet(x^i,y^i,\hat x^i,\hat y^i)\coloneqq\{x^i,y^i,\hat x^i,\hat y^i\}$.

\textbf{Case 2:} $v_i\in V(\tau_0)$ is an inner tree vertex.
We construct a set of partial external skeletons for $v_i$ as follows.
\begin{enumerate}[label=\arabic*.]
    \item\label{ppreskITV1} Suppose $(x^i,y^i,\hat x^i,\hat y^i)$ is a valid connection sequence for $v_i$.
    Let $w_0^i\coloneqq x^i$, $w_{p_i}^i\coloneqq y^i$, $\hat{w}_0^i\coloneqq \hat x^i$ and $\hat{w}_{p_i}^i\coloneqq\hat y^i$.
    Let $W^i_0\coloneqq\{w_0^i, w_{p_i}^i, \hat{w}_0^i, \hat{w}_{p_i}^i\}$.
    \item For each $k\in[p_i-1]$, iteratively choose two vertices $w_k^i,\hat w_k^i\in V(\mathcal{M}_{b(i)}(v_i))\setminus W^i_{k-1}$ such that $f_k^i+(w_{k-1}^i,w_k^i,\hat{w}_{k-1}^i,\hat{w}_k^i)$ is a valid connection sequence for $u^i_k$, and let $W^i_k\coloneqq W^i_{k-1}\cup\{w_k^i,\hat{w}_k^i\}$.
    \COMMENT{Alternatively, we can think of this by following these two steps:
    \begin{enumerate}
    \item For each $k\in[p_i-1]$, iteratively let $W_k^i\coloneqq\{w_0^i,\ldots,w_{k-1}^i\}$ and choose a vertex $w_k^i\in V(\mathcal{M}_{b(i)}(v_i))\setminus (W_k^i\cup\{y^i,\hat x^i,\hat y^i\})$ with $w_k^i\neqp w_{k-1}^i$ if $\Delta(u_k^i)$ is even, or $w_k^i\eqp w_{k-1}^i$ if $\Delta(u_k^i)$ is odd.
    \item For each $k\in[p_i-1]$, iteratively let $W_k^{*i}\coloneqq\{\hat{w}_0^i,\ldots,\hat{w}_{k-1}^i\}$ and choose a vertex $\hat{w}_k^i\in V(\mathcal{M}_{b(i)}(v_i))\setminus(W_{p_i-1}^i\cup W_k^{*i}\cup\{y_i,\hat y_i,w_{p_i-1}^i\})$ with $\hat{w}_k^i\neqp w_k^i$.
    \end{enumerate}
    }
\end{enumerate}
Note that the definition of $q$ in \eqref{qandt} and the bound on $p_i$ in \eqref{equa:step9.2} ensure that we have sufficiently many vertices to choose from (similar comments apply in the other cases).
Moreover, \eqref{equa:Delta} implies that $f^i_{p_i}+(w_{p_i-1}^i,w_{p_i}^i,\hat{w}_{p_i-1}^i,\hat{w}_{p_i}^i)$ is a valid connection sequence for $u^i_{p_i}$.
The partial external skeleton for $v_i$ and connection sequence $(x^i,y^i,\hat x^i,\hat y^i)$ is defined as 
\[L_i^\bullet(x^i,y^i,\hat x^i,\hat y^i)\coloneqq\{x^i,\hat x^i\}\cup\bigcup_{k=1}^{p_i}\left(\{w_k^i,\hat{w}_k^i\}\cup L_{j_k^i}^\bullet(f_k^i+(w_{k-1}^i,w_k^i,\hat{w}_{k-1}^i,\hat{w}_k^i))\right).\]

\textbf{Case 3:} $v_i\in V(\tau_0)$ is an atomic vertex which is not a leaf.
We construct a set of partial external skeletons for $v_i$ as follows.
\begin{enumerate}[label=\arabic*.]
    \item\label{ppreskAV1} Assume $(x^i,y^i,\hat x^i,\hat y^i)$ is a valid connection sequence for $v_i$.
    Let $w^i_0\coloneqq x^i$.
    \item\label{ppreskAV3} For each $k\in[p_i]$, iteratively choose distinct vertices $z_k^i,w_k^i,\hat{z}_k^i,\hat{w}_k^i\in (V(\mathcal{M}_{b(i)+k-1}(v_i))\cap V(\cM_{\nu_{k}^i}))\setminus\{x^i,y^i,\hat x^i,\hat y^i\}$ satisfying that $z_k^i\neqp w_{k-1}^i$ and $f_k^i+(z_k^i,w_k^i,\hat{z}_k^i,\hat{w}_k^i)$ is a valid connection sequence for $u^i_k$\COMMENT{That is, we have that
    \begin{enumerate}[label=\ref{ppreskAV3}\arabic*.]
        \item $z_k^i\neqp w_{k-1}^i$;
        \item $\hat{z}_k^i\neqp z_k^i$;
        \item $w_k^i\neqp z_k^i$ if $\Delta(u_{k+1}^i)$ is even and $w_k^i\neqp w_{k-1}^i$ otherwise, and
        \item $\hat{w}_k^i\neqp w_k^i$.
    \end{enumerate}
    }.
\end{enumerate}
Then, the partial external skeleton for $v_i$ and connection sequence $(x^i,y^i,\hat x^i,\hat y^i)$ is defined as
\[
    L_i^\bullet(x^i,y^i,\hat x^i,\hat y^i)\coloneqq\{x^i,y^i,\hat x^i,\hat y^i\}
    \cup\bigcup_{k=1}^{p_i}\left(\{z_k^i,w_k^i,\hat{z}_k^i,\hat{w}_k^i\} \cup L_{j_{k}^i}^\bullet(f_k^i+(z_k^i,w_k^i,\hat{z}_k^i,\hat{w}_k^i))\right).
\]

After having constructed all these partial external skeletons for all $v_i$ with $i\in[m]$, we are now ready to construct $L^\bullet$.
\begin{enumerate}[label=\arabic*.]
    \item Choose any vertex $w_0^0\in V(\mathcal{A}_1(v_0))$.
    \item\label{presk2} For each $k\in[p_0]$, iteratively choose four distinct vertices $z_k^0,\hat{z}_k^0,w_k^0,\hat{w}_k^0\in(V(\mathcal{M}_k(v_0))\cap V(\cM_{\nu_k^0}))$ satisfying that $z_k^0\neqp w_{k-1}^0$ and $f_k^0+(z_k^0,w_k^0,\hat{z}_k^0,\hat{w}_k^0)$ is a valid connection sequence for $u_k^0$\COMMENT{That is, we have that
    \begin{enumerate}[label=\ref{presk2}\arabic*.]
        \item $z_k^0\neqp w_{k-1}^0$;
        \item $\hat{z}_k^0\neqp z_k^0$;
        \item $w_k^0\neqp z_k^0$ if $\Delta(u_k^0)$ is even, and $w_k^0\neqp w_{k-1}^0$ otherwise, and
        \item $\hat{w}_k^0\neqp w_k^0$.
    \end{enumerate}
    }.
\end{enumerate}
Then, we define 
\[L^\bullet\coloneqq\bigcup_{k=1}^{p_0}\left(\{z_k^0,w_k^0,\hat{z}_k^0,\hat{w}_k^0\}\cup L_{j_k^0}^\bullet(f_k^0+(z_k^0,w_k^0,\hat{z}_k^0,\hat{w}_k^0))\right).\]

Observe that \ref{item:preskprop1}--\ref{item:preskprop4} hold by construction.
In turn, \ref{item:preskprop5} holds because of the definition of $\tau_0$.
Indeed, observe that $V(T^*)\cap R'=\varnothing$ by \ref{itm:C1}.
Moreover, by the construction above, all vertices in $L^\bullet$ are incident to some edge in a clone of the tree $T^*$, and thus, they cannot lie in $\mathcal{M}_{\mathrm{Res}}$.\\

\textbf{Step~12: Constructing an auxiliary tree $\boldsymbol{\tau_0'}$.}
In order to extend the external skeleton into the skeleton and construct an almost spanning cycle, we first need to extend $\tau_0$ to a new auxiliary tree $\tau_0'$ which encodes information about some additional molecules.

We construct $\tau_0'$\index{tau02@$\tau_0'$} by appending some new leaves to $\tau_0$.
Note that $\tau_0$ was built by encoding all the information about $T^*$, and $\tau_0'$ will encode the information about $T''$.
In particular, by \ref{itm:C1}, each cube $C\in\mathcal{C}'$ which intersects $T''$ and does not intersect $T^*$ contains at least one vertex $u$ which is joined to $T^*$ by an edge $e'=\{u,v\}\in E(T'')$ such that $v\in V(C')$, where $C\neq C'\in\mathcal{C}''$.
Note that the construction of $\tau_0$ implies that $C'$ is represented in $\tau_0$.
For each such cube $C$, choose one such vertex $u$ and append a new vertex to the atomic vertex representing $C'$ in $\tau_0$ via an edge $e$ which originates as $e' \in E(T'')$.
We say that this newly added vertex is atomic and \emph{represents} $C$.
The resulting tree after all these leaves are appended is $\tau_0'$.
In particular, $\tau_0\subseteq\tau_0'$, and it now follows that precisely the $C\in\cC'''$ are represented in $\tau_0'$, where $\cC'''$ is as defined in Step~9.
Furthermore, it follows from \ref{itm:C2} that\COMMENT{I am not very sure how to explain this. Intuitively, we are simply reconstructing a new tree $\tau_0'$ using $T''$ instead of $T^*$. As the bounds on the degrees for $T^*$ come from $T''$, the same bounds must hold here.\\
To be a bit more careful, observe first that, by construction, we do not append any new leaves to inner tree vertices, so the bound on the degrees of inner tree vertices must satisfy \eqref{equa:step9.1}.
For the bound on the maximum degree, take a vertex of maximum degree of $\tau_0'$; we may assume that it is an atomic vertex $v$.
Observe that each new leaf that we have appended to this vertex corresponds to an edge $e$ incident to some vertex $x$ represented in $\tau_0'$ by $v$ such that $e\in E(T'')$ but $e\notin E(T^*)$.
That means that the degree of $v$ in $\tau_0$ was bounded by at most the number of such vertices $x$ ($2^\ell$) times the maximum degree of each such vertex in $T''$ ($12D$) minus the sum over all such vertices $x$ of the number of edges incident to them in $T''$ which were removed in $T^*$; these are now added back, so the overall maximum degree is at most $2^\ell\cdot12D$.}
\begin{equation}\label{equa:step11.1}
\begin{minipage}[c]{0.7\textwidth}
$d_{\tau_0'}(v)\leq12D$ for all $v\in V(\tau_0')$ which are inner tree vertices, and $\Delta(\tau_0')\leq 12\cdot2^\ell D$.
\end{minipage}\ignorespacesafterend 
\end{equation} 
For all vertices of $\tau_0'$, we will use the same notation for the vertices, cubes and molecules that they represent as we did for the vertices of $\tau_0$.
Note that, by \ref{itm:C325} and \ref{itm:C54},
\begin{enumerate}[label=$(\mathrm{CP})$]
\item\label{itm:AB5} every cube $C$ belonging to some absorbing $\ell$-cube pair in $\cC^\mathrm{sc}\cup\cC^\mathrm{\neg sc}$ is represented in $\tau_0'$.
\end{enumerate}

It will be important for us that $\tau_0'$ represents `most' vertices of the hypercube.
In particular, for each $x\in V(I)$, let $\lambda(x)$ denote the number of vertices $y\in N_I(x)$ which are represented in $\tau_0'$ by atomic vertices.
By \ref{itm:C5bis}\COMMENT{Any vertex that belongs to the intersection of the tree $T''$ and a cube $C\in\cC'$ is represented in $\tau_0'$.}, we have that
\begin{equation}\label{equa:step10}
    \lambda(x)\geq(1-2/\ell^4)n.
\end{equation}
By an averaging argument, it follows that at least $(1-2/\ell^4)2^{n-s}$ vertices $x\in V(I)$ are represented in $\tau_0'$ by atomic vertices.
We will construct an almost spanning cycle in $G'$ which contains all the clones of these vertices.

Let $m'\coloneqq|V(\tau_0')|-1$.
Label $V(\tau_0')\setminus V(\tau_0)=\{v_{m+1},\ldots,v_{m'}\}$ arbitrarily.
For each $i\in[m]$, we define $\tau_i'$ as the maximal subtree of $\tau_0'$ which contains $v_i$ and all of whose vertices have labels at least as large as $i$.
For each $i\in[m]_0$, let $p_i'\coloneqq d_{\tau_i'}(v_i)$\COMMENT{No need to define for $i>m$, as for all of these we know $d_{\tau_i'}(v_i)=0$.} and let $N_{\tau_i'}(v_i)=\{u_1^i,\ldots,u_{p_i'}^i\}$ (where the labelling is consistent with that of $N_{\tau_i}(v_i)$).
For each $i\in[m]_0$ and $k\in[p_i']\setminus[p_i]$, let $e_k^i\coloneqq\{v_i,u_k^i\}$, let $f_k^i\coloneqq D(e_k^i)$, and let $j_k^i$ be the label of $u_k^i$ in $\tau_0'$.
Furthermore, for each $i\in[m]_0$ such that $v_i$ is an atomic vertex, and for each $k\in[p_i']\setminus[p_i]$, consider the unique edge $e_{i,k}' \in E(I(G'))$ so that $e_k^i$ originates from $e_{i,k}'$ (recall that `originating' was defined in Step 10), and let \index{nuik@$\nu_k^i$}$\nu_k^i$ be the endpoint of $e_{i,k}'$ in $C(v_i)$. 
Finally, for each $i\in[m']\setminus[m]$ we set $\Delta(v_i)\coloneqq0$\COMMENT{Observe that this does not alter the value of $\Delta(v_i)$ for any $i\in[m]_0$.}.

As in Step~10, we consider the partition into slices for the new molecules arising from the newly added cubes represented by $\tau_0'$.
For each $i\in [m']\setminus[m]$, we assign an input index $b(i)\in[t]$.
To do so, for each $i\in[m]_0$ such that $v_i$ is an atomic vertex and each $k\in[p_i']\setminus[p_i]$, we set $b(j_k^i)\coloneqq b(i)+k-1\pmod{t}$.
Similarly to Step~10, \eqref{qandt} and \eqref{equa:step11.1} imply that in this case $b(j^i_k)\neq b(j^i_{k'})$ for all $k\neq k'$.
For each $i\in[m']\setminus[m]$, let $\ell_i$ be the label in $\tau_0'$ of the unique vertex adjacent to $v_i$ (i.e., the parent label of $i$), and let $m_i$ be the label of $v_i$ in $N_{\tau_{\ell_i}'}(v_{\ell_i})$.
Note that $b(i)=b(\ell_i)+m_i-1$.\\

\textbf{Step~13: Fixing absorbing $\boldsymbol{\ell}$-cube pairs for vertices that need to be absorbed.}
At this point, we can determine every vertex in $V(\cQ^n)$ that will have to be absorbed into the almost spanning cycle we are going to construct.
For every vertex $x \in V(I)$ not represented in $\tau_0'$, we will have to absorb all vertices in $\cM_x$.\COMMENT{Bondless and isolated molecules, as well as vertices not covered by either the tree or the molecules, and some vertices covered by the tree $T''$ by which are deleted in the construction of $\tau_0$ (because they were leaves not pertaining to any cube).}
Furthermore, for each $v\in V(\tau_0)$ which is an inner tree vertex, we will also need to absorb all vertices in $\mathcal{M}_v\setminus L^\bullet$.\COMMENT{We need to absorb every vertex from a bondless molecule and isolated molecule, as well of most of the vertices from an inner tree vertex molecule.}
By \ref{item:preskprop1}, this means that, in each such molecule $\cM_v$, the same number of vertices of each parity need to be absorbed.
Recall the definition of $V_{\mathrm{abs}}$ from Step~9.
This is precisely the set of vertices which are not represented in $\tau_0'$ by an atomic vertex and, therefore, it is the set of all vertices $x\in V(I)$ such that some clone of $x$ needs to be absorbed.
It follows from \eqref{equa:step10} that 
\begin{equation}\label{equa:step11}
    |V_{\mathrm{abs}}|\leq2^{n-s+1}/\ell^4.
\end{equation}
Now, for each $x\in V_{\mathrm{abs}}$, we will pair the vertices in each slice which need to be absorbed (each pair consisting of one vertex of each parity) and fix an absorbing $\ell$-cube pair for each such pair of vertices.
The absorbing $\ell$-cube pair that we fix will be the one given by \ref{itm:C3} or \ref{itm:C6} for this pair of vertices, depending on whether $x\in V_\mathrm{sc}$ or not.

For each $x\in V_{\mathrm{abs}}$ and $\mathcal{S}\in\mathcal{S}(\mathcal{M}_x)$, let $S(x,\mathcal{S})\coloneqq V(\cS)\cap L^\bullet$.
It follows by \ref{item:preskprop1}--\ref{item:preskprop4} that $|S(x,\mathcal{S})|\leq 24D$\COMMENT{Let us look at the different cases.
For each $x\in V_{\mathrm{abs}}$ which is not represented in $\tau_0'$, we have that $|S(x,\mathcal{S})|=0$ for all $\mathcal{S}\in\mathcal{S}(\mathcal{M}_x)$.
If $x$ is an inner tree vertex and $\mathcal{S}$ is not the input slice, we also have that $|S(x,\mathcal{S})|=0$.
As for the input slice, we have shown in the last steps that at most $24D$ vertices were used from an inner tree vertex; indeed, recall that $p_i\leq12D-1$.} and $S(x,\mathcal{S})$ contains the same number of vertices of each parity.
(Here we also use that $p_i\leq 12D-1$ for every inner tree vertex $v_i$ by \eqref{equa:step9.2} and \eqref{equa:step11.1}.)
Therefore, the matching $\mathfrak{M}(\mathcal{S},S(x,\mathcal{S}))$ defined in Step~3 is well defined.
Recall that each edge $e\in\mathfrak{M}(\mathcal{S},S(x,\mathcal{S}))$ gives rise to a unique index $i\in[K]$ via the relation $N(e)=A_i(x)$. 
(Here we ignore all those indices $i' \in [K]$ arising by artificially increasing the size of $\mathfrak{A}(x)$, see the beginning of Step~4.)
For each $x\in V_{\mathrm{abs}}$, let $\mathfrak{I}_{x}\subseteq[K]$ be the set of indices $i\in[K]$ which correspond to edges in $\bigcup_{\cS\in\cS(\cM_x)}\mathfrak{M}(\mathcal{S},S(x,\mathcal{S}))$\COMMENT{Basically, for any index $i\notin\mathfrak{I}_x$, we no longer care about $A_i(x)$ or $M(A_i(x))$.}.

For each $x\in V_{\mathrm{abs}}$ and $i\in\mathfrak{I}_x$, as stated in \ref{itm:C3} and \ref{itm:C6}, we have already fixed an absorbing $\ell$-cube pair for the clones of $x$ corresponding to $(x,i)$.
Let\index{Vabs2@$V^\mathrm{abs}$} 
\[V^\mathrm{abs}\coloneqq\bigcup_{x\in V_{\mathrm{abs}}}V(\mathcal{M}_x)\setminus L^\bullet.\]
As discussed above, this is the set of all vertices that need to be absorbed.
Recall that $G'$ was defined before \ref{itm:C2}--\ref{itm:C6}.
It follows from \ref{itm:C3} and \ref{itm:C6} that $((H\cup G)\setminus F) \cup G'$ contains a set $\mathcal{C}^\mathrm{abs}=\{(C^l(u),C^r(u)):u\in V^\mathrm{abs}\}$ of absorbing $\ell$-cube pairs such that
\begin{enumerate}[label=$(\mathrm{C_{\arabic*}})$]
    \item\label{itm:RL1} for all distinct $u,v\in V^\mathrm{abs}$, the absorbing $\ell$-cube pairs $(C^l(u),C^r(u))$ and $(C^l(v),C^r(v))$ for $u$ and $v$ are vertex-disjoint and $(C^l(u),C^r(u)) \cup (C^l(v),C^r(v))- \{u,v\} \subseteq G'$;
    \item\label{itm:RL22} there exists a pairing $\mathcal{U}=\{f_1, \dots, f_{K^\diamond}\}$ of $V^\mathrm{abs}$ such that
    \begin{enumerate}[label=$({\mathrm{C}_\mathrm{2.\arabic*}})$]
        \item\label{itm:RL222} for all $i\in[K^\diamond]$, if $f_i=\{u_i,u_i'\}$, then $u_i\neqp u_i'$;
        \item\label{itm:RL2} if $f_i=\{u_i,u_i'\}$, then there is a vertex $v\in V_\mathrm{abs}$ such that $u_i$ and $u_i'$ are clones of $v$ which lie in the same slice of $\mathcal{M}_v$, and $(C^l(u_i),C^r(u_i))$ and $(C^l(u_i'),C^r(u_i'))$ are clones of the same absorbing $\ell$-cube pair for $v$ in $I$ such that $(C^l(u_i),C^r(u_i))$ lies in the same layer as $u_i$ and $(C^l(u_i'),C^r(u_i'))$ lies in the same layer as $u_i'$;
        \item\label{itm:RL3} if $u,u'\in V^\mathrm{abs}$ do not form a pair $f\in\mathcal{U}$, then $(C^l(u),C^r(u))$ and $(C^l(u'),C^r(u'))$ are clones of vertex-disjoint absorbing $\ell$-cube pairs in $I$ (except in the case when $u,u'$ are clones of the same vertex $v\in V_\mathrm{abs}$, in which case $(C^l(u),C^r(u))$ and $(C^l(u'),C^r(u'))$ are clones of absorbing $\ell$-cube pairs in $I$ which intersect only in $v$);
    \end{enumerate}
    \item\label{itm:RL4} if we let $\mathcal{C}^*\coloneqq\bigcup_{(C^l(u),C^r(u))\in\mathcal{C}^\mathrm{abs}}\{C^l(u),C^r(u)\}$, then $\mathcal{C}^*$ contains either two or no clones of each cube $C\in\mathcal{C}''\cap\cC'''$, and every cube in $\mathcal{C}^*$ is a clone of some cube $C\in\mathcal{C}''\cap\cC'''$.
\end{enumerate}
The pairing described in \ref{itm:RL22} is given by the matchings $\mathfrak{M}(\cS,S(x,\cS))$.
Furthermore, it follows from \ref{itm:C31}, \ref{itm:C52} and \ref{item:preskprop5} that
\begin{enumerate}[label=$(\mathrm{C_{\arabic*}})$]\setcounter{enumi}{3}
    \item\label{itm:RL5} the set of all tips of the absorbing $\ell$-cube pairs in $\mathcal{C}^\mathrm{abs}$ is disjoint from $L^\bullet$.
\end{enumerate}

We denote by $\mathfrak{L}$, $\mathfrak{R}_1$ and $\mathfrak{R}_2$ the collections of all left absorber tips, right absorber tips, and third absorber vertices, respectively, of the absorbing $\ell$-cube pairs in $\mathcal{C}^\mathrm{abs}$.
Observe that the following properties are satisfied:
\begin{enumerate}[label=$(\mathrm{C}^*\arabic*)$]
    \item\label{item:Cstarcon1} For all $i\in[m']_0$ such that $v_i$ is an atomic vertex and all $j\in[t]$, we have that $|\mathfrak{L}\cap V(\mathcal{M}_j(v_i))|\in\{0,2\}$ and, if $|\mathfrak{L}\cap V(\mathcal{M}_j(v_i))|=2$, then these two vertices $u,u'$ lie in different atoms of the slice and satisfy that $u\neqp u'$.
    \item\label{item:Cstarcon2} For all $i\in[m']_0$ such that $v_i$ is an atomic vertex and all $j\in[t]$, we have that $|(\mathfrak{R}_1\cup\mathfrak{R}_2)\cap V(\mathcal{M}_j(v_i))|\in\{0,4\}$.
    If $|(\mathfrak{R}_1\cup\mathfrak{R}_2)\cap V(\mathcal{M}_j(v_i))|=4$, then these four vertices form two pairs such that one vertex of each pair belongs to $\mathfrak{R}_1$ and the other to $\mathfrak{R}_2$.
    Each of these pairs lies in a different atom of the slice and satisfies that its two vertices are adjacent in $G'$.
    \item\label{item:Cstarcon3} For all $i\in[m']_0$ such that $v_i$ is an atomic vertex and all $j\in[t]$, if $\mathfrak{L}\cap V(\mathcal{M}_j(v_i))\neq\varnothing$, then $(\mathfrak{R}_1\cup\mathfrak{R}_2)\cap V(\mathcal{M}_j(v_i))=\varnothing$.
    \item\label{item:Cstarcon4} The sets described in \ref{item:Cstarcon1} and \ref{item:Cstarcon2} partition $\mathfrak{L}$ and $\mathfrak{R}_1\cup\mathfrak{R}_2$, respectively\COMMENT{that is, there are no vertices of $\mathfrak{L}\cup\mathfrak{R}_1\cup\mathfrak{R}_2$ outside the cube molecules represented in $\tau_0'$.}.
\end{enumerate}
Indeed, \ref{item:Cstarcon1}--\ref{item:Cstarcon3} follow from \ref{itm:RL22} and \ref{itm:RL4}, and \ref{item:Cstarcon4} follows by \ref{itm:AB5}.

For each $u\in V^\mathrm{abs}$, we denote the edge consisting of the right absorber tip and the third absorber vertex of $(C^l(u),C^r(u))$ by $e_\mathrm{abs}(u)$, and we denote by $\mathcal{P}^\mathrm{abs}(u)$ the path of length three formed by the third absorber vertex, the left absorber tip, $u$, and the right absorber tip, visited in this order.
Note that $e_\mathrm{abs}(u)\in E(G')$ by \ref{itm:RL1}.
Moreover, recall that $\cC^\mathrm{abs}$ consists of absorbing $\ell$-cube pairs in $((H\cup G)\setminus F)\cup G')$.
Thus, $\mathcal{P}^\mathrm{abs}(u)\subseteq ((H\cup G)\setminus F)\cup G'$.\\

\textbf{Step~14: Constructing the skeleton.}
We can now define the skeleton for the almost spanning cycle.
Intuitively, this skeleton builds on the external skeleton by adding more structure that the cycle will have to follow.
In particular, the skeleton adds the edges used to traverse from each slice in a cube molecule to its neighbouring slices, and it also incorporates the cube molecules represented in $\tau_0'$ which were not represented in $\tau_0$.
(The reason why these were not incorporated earlier is the following: if we already choose the valid connection sequences for these cube molecules in Step~12, then the tips of the absorbing cubes chosen in Step~13 might have non-empty intersection with the external skeleton, which we want to avoid, see \ref{itemskelprop6} below.)
Furthermore, the skeleton gives an ordering to its vertices, and the cycle will visit the vertices of the skeleton in this order.

We will build a skeleton $\mathcal{L}=(x_1,\ldots,x_r)$, for some $r\in\mathbb{N}$, and write $\mathcal{L}^\bullet\coloneqq\{x_1,\ldots,x_r\}$.
We will construct $\mathcal{L}$ in such a way that the following properties hold:
\begin{enumerate}[label=$(\mathrm{S}\arabic*)$]
    \item\label{itemskelprop1} For all distinct $k,k'\in[r]$, we have that $x_k\neq x_{k'}$.
    \item\label{itemskelprop2} $\{x_1,x_r\}\in E(G')$.
    \item\label{itemskelprop3}  For every $k\in[r-1]$, if $x_k$ and $x_{k+1}$ do not both lie in the same slice of a cube molecule represented in $\tau_0'$, then $\{x_k,x_{k+1}\}\in E(G')$.
    Moreover, in this case, if $x_{k+1}$ lies in a cube molecule represented in $\tau_0'$, then $x_{k+2}$ lies in the same slice of this cube molecule as $x_{k+1}$.
    \item\label{itemskelprop35} For every $i\in[m']_0$ and every $j\in[t]$, no three consecutive vertices of $\mathcal{L}$ lie in $\cM_j(v_i)$ (here $\mathcal{L}$ is viewed as a cyclic sequence of vertices).
    \item\label{itemskelprop4} For every $i\in[m']$ such that $v_i$ is an atomic vertex and every $j\in[t]$, we have that $|V(\mathcal{M}_j(v_i))\cap\mathcal{L}^\bullet|$ is even and $4\leq|V(\mathcal{M}_j(v_i))\cap\mathcal{L}^\bullet|\leq12$\COMMENT{Observe that, for many slices, we only need $8$ here.
    The upper bound of $12$ only affects \emph{input slices}.
    Let $v_i$ be an atomic vertex (other than $v_0$) which is not a leaf of $\tau_0'$.
    Indeed, in this case the skeleton will contain six vertices of $\cM_{b(i)}(v_i)$ each time we traverse $\tau_0'$ (four of them given by the external skeleton and two of them added in the skeleton to `traverse' this molecule).
    (Alternatively, think of this as we follow the skeleton.
    We first see this molecule $\cM(v_i)$ when we enter its input slice.
    We keep going to the next molecule $\cM(u^i_1)$ (2 vertices for the skeleton).
    Then, we come back from this next molecule and go down to the next slice (2 vertices); then, once we finished going around the molecule, we come back into the input slice and leave back to the `parent' molecule (2 vertices).)
    Note that we could redefine the input slices and the external skeleton so that we never use an input slice to go to the `children' molecules, in which case 8 would be a general bound, but this is not the way we have written it.}. 
    In particular, $|V(\mathcal{M}_t(v_0))\cap\mathcal{L}^\bullet|=4$.
    \item\label{itemskelprop5} For all $k\in[r]$ except two values, we have that $x_k\neqp x_{k+1}$.
    The remaining two values $k_1,k_2\in[r]$ correspond to two pairs of vertices $x_{k_1},x_{k_1+1},x_{k_2},x_{k_2+1}\in V(\mathcal{M}_t(v_0))$.
    For these two values, we have that $x_{k_1}\neqp x_{k_2}$ and either 
    \begin{enumerate}[label=$(\mathrm{\roman*})$]
        \item\label{itemskelprop51} $x_{k_1}\eqp x_{k_1+1}$ and $x_{k_2}\eqp x_{k_2+1}$, or
        \item\label{itemskelprop52} $x_{k_1}\neqp x_{k_1+1}$ and $x_{k_2}\neqp x_{k_2+1}$,
    \end{enumerate}
    where $x_{k_1},x_{k_2}\in V(\mathcal{A}_{(t-1)q+1}(v_0))$ and $x_{k_1+1},x_{k_2+1}\in V(\mathcal{A}_{tq}(v_0))$.
    \item\label{itemskelprop6} $\mathcal{L}^\bullet\cap(\mathfrak{L}\cup\mathfrak{R}_1\cup V^\mathrm{abs})=\varnothing$ and $L^\bullet\subseteq\mathcal{L}^\bullet$.
\end{enumerate}

As happened with the external skeleton, the skeleton is built recursively from partial skeletons, which are defined first for the leaves.
This recursive construction means that the overall order in which the molecules are visited will be determined by a depth first search of the tree $\tau_0'$.
Moreover, as discussed in \cref{section:outline5}, for parity reasons the skeleton will actually traverse $\tau_0'$ twice.
These two traversals will be `tied together' in the final step of the construction of the skeleton.

Note that, for each $i\in[m]$, the starting and ending vertices $x^i$, $\hat{x}^i$, $y^i$, $\hat{y}^i$ for the partial skeletons for $v_i$ are determined by the external skeleton.
For each $i\in[m']\setminus[m]$, the starting and ending vertices for the partial skeletons of $v_i$ will be determined when constructing the partial skeleton for the parent vertex $v_{\ell_i}$ of $v_i$.
In particular, when constructing the partial skeleton for $v_{\ell_i}$, we will define vertices $z_{m_i}^{\ell_i},\hat{z}_{m_i}^{\ell_i},w_{m_i}^{\ell_i},\hat{w}_{m_i}^{\ell_i}\in\cM_{b(i)}(v_{\ell_i})$\COMMENT{Here we have that $b(i)=b(\ell_i)+m_i-1$, as was noted at then end of Step~12.}.
Then, the starting and ending vertices for the partial skeleton of $v_i$ will be 
\begin{equation}\label{equa:skel1}
    (x^i,y^i,\hat x^i,\hat y^i)\coloneqq f_{m_i}^{\ell_i}+(z_{m_i}^{\ell_i},w_{m_i}^{\ell_i},\hat{z}_{m_i}^{\ell_i},\hat{w}_{m_i}^{\ell_i}).
\end{equation}
(Recall that $\ell_i$, $m_i$, $b(i)$ and $f_{m_i}^{\ell_i}$ were defined at the end of Step~12.)

We are now in a position to define the partial skeletons formally.
The construction proceeds by induction on $i\in[m']$ in decreasing order, starting with $i=m'$.
Recall from the beginning of Step~11 that, for all $i\in[m]$, $x^i,y^i\in V(\mathcal{M}_{b(i)}(v_i))$ are the starting and ending vertices for the first partial skeleton $\mathcal{L}(x^i,y^i)$ for $v_i$, respectively, and $\hat x^i,\hat y^i\in V(\mathcal{M}_{b(i)}(v_i))$ are the starting and ending vertices for the second partial skeleton $\hat{\mathcal{L}}(\hat x^i,\hat y^i)$ for $v_i$, respectively.
The vertices $x^i,y^i,\hat x^i,\hat y^i$ were fixed in the construction of the external skeleton, and they form a valid connection sequence.
For each $i\in[m']\setminus[m]$, the vertices $x^i,y^i,\hat x^i,\hat y^i\in V(\mathcal{M}_{b(i)}(v_i))$ defined in \eqref{equa:skel1} will also form a valid connection sequence.

Let $\cF\coloneqq\mathfrak{L}\cup\mathfrak{R}_1\cup L^\bullet$\COMMENT{This will be a `forbidden' set. Making sure the vertices we choose are not in $\mathfrak{R}_1$ is not strictly necessary, but it will make the writing much easier.}.
For each $k\in[2^s]$, let \index{ehatk@$\hat e_k$}$\hat e_k$ be the direction of the edges in $\cQ^n$ between $L_k$ and $L_{k+1}$.
Throughout the following construction, we will often choose vertices which are used to transition between neighbouring slices, all while avoiding the set $\cF$.
Similarly to the proof of \cref{lem:slicecover}, all of these choices can be made by \ref{item:preskprop2}, \ref{item:preskprop3}, \ref{item:Cstarcon1}, \ref{item:Cstarcon2}, and because all cube molecules considered here are bonded in $G_5$ and, therefore, also in $G'$.
(The latter holds since for each atomic vertex $v\in V(\tau_0')$ the corresponding cube $C(v)$ satisfies $C(v)\in\cC'$.)
Whenever we mention a vertex that we do not define here, we refer to the vertex with the same notation defined when constructing the external skeleton in Step~11.

Suppose that $i\in[m']$ and that for every $i'\in[m']\setminus([i]\cup[m])$ and every valid connection sequence $(x^{i'},y^{i'},\hat x^{i'},\hat y^{i'})$ for $v_{i'}$ we have already defined two partial skeletons $\mathcal{L}(x^{i'},y^{i'})$, $\hat{\mathcal{L}}(\hat x^{i'},\hat y^{i'})$ for $v_{i'}$ with this connection sequence.
(As discussed above, eventually we will only use the two partial skeletons for $v_{i'}$ with connection sequence as defined in \eqref{equa:skel1}.)
Moreover, suppose that for every $i'\in[m]\setminus[i]$ we have already defined two partial skeletons $\mathcal{L}(x^{i'},y^{i'})$, $\hat{\mathcal{L}}(\hat x^{i'},\hat y^{i'})$ for $v_{i'}$ with connection sequence $(x^{i'},y^{i'},\hat x^{i'},\hat y^{i'})$ (fixed by the external skeleton).
If $i\in[m]$, let $(x^{i},y^{i},\hat x^{i},\hat y^{i})$ be the connection sequence for $v_i$ fixed by the external skeleton.
If $i\in[m']\setminus[m]$, let $(x^{i},y^{i},\hat x^{i},\hat y^{i})$ be any connection sequence for $v_i$.
We will now define the two partial skeletons for $v_i$ with connection sequence $(x^{i},y^{i},\hat x^{i},\hat y^{i})$.
We consider several cases.

\textbf{Case 1:} $v_i$ is a leaf of $\tau_0'$.
We construct the partial skeletons as follows.
Let $x_0^i\coloneqq x^i$ and $\hat{x}_0^i\coloneqq\hat{x}^i$.
For each $k\in[t-1]_0$, iteratively choose any two vertices $y_k^i,\hat{y}_k^i\in V(\mathcal{A}_{(b(i)+k)q}(v_i))\setminus (\cF\cup\{x^i,y^i,\hat x^i,\hat y^i\})$ satisfying that
    \begin{enumerate}[label=$\arabic*$.]
        \item $y_k^i\neqp x_k^i$ and $\hat{y}_k^i\neqp \hat{x}_k^i$;
        \item $x_{k+1}^i\coloneqq y_k^i+\hat{e}_{(b(i)+k)q}\notin \cF\cup\{x^i,y^i,\hat x^i,\hat y^i\}$ and $\hat{x}_{k+1}^i\coloneqq\hat{y}_k^i+\hat{e}_{(b(i)+k)q}\notin \cF\cup\{x^i,y^i,\hat x^i,\hat y^i\}$, and 
        \item $\{y_k^i,x_{k+1}^i\},\{\hat{y}_k^i,\hat{x}_{k+1}^i\}\in E(G')$.
    \end{enumerate}
Recall that we use $\bigtimes$ to denote the concatenation of sequences.
The first and second partial skeletons for $v_i$ with connection sequence $(x^i,y^i,\hat x^i,\hat y^i)$ are given by
\[\mathcal{L}_i(x^i,y^i)\coloneqq(x^i)\left(\bigtimes_{k=0}^{t-1}(y_k^i,x_{k+1}^i)\right)(y^i)\quad\text{ and }\quad\hat{\mathcal{L}}_i(\hat x^i,\hat y^i)\coloneqq(\hat x^i)\left(\bigtimes_{k=0}^{t-1}(\hat{y}_k^i,\hat{x}_{k+1}^i)\right)(\hat y^i).\]

\textbf{Case 2:} $v_i\in V(\tau_0)$ is an inner tree vertex.
Then, the first and second partial skeletons for $v_i$ with connection sequence $(x^i,y^i,\hat x^i,\hat y^i)$ are defined as
\[\mathcal{L}_i(x^i,y^i)\coloneqq(x^i)\bigtimes_{k=1}^{p_i}(\mathcal{L}_{j_k^i}(x^{j_k^i},y^{j_k^i}),w_k^i)\quad\text{ and }\quad\hat{\mathcal{L}}_i(\hat x^i,\hat y^i)\coloneqq(\hat x^i)\bigtimes_{k=1}^{p_i}(\hat{\mathcal{L}}_{j_k^i}(\hat{x}^{j_k^i},\hat{y}^{j_k^i}),\hat{w}_k^i),\]
where $j_k^i$ was defined in Step~10.

\textbf{Case 3:} $v_i\in V(\tau_0)$ is an atomic vertex which is not a leaf of $\tau_0'$.
We construct the partial skeletons for $v_i$ as follows.
(Recall that, for each $k\in[p_i']\setminus[p_i]$, the vertex $\nu^i_k$ was defined in Step~12.)
\begin{enumerate}[label=$\arabic*$.]
    \item\label{pskelatom2} For each $k\in[p_i]$, iteratively choose distinct vertices $y_k^i,\hat{y}_k^i\in V(\mathcal{A}_{(b(i)+k-1)q}(v_i))\setminus \cF$ such that
    \begin{enumerate}[label=\ref{pskelatom2}$\arabic*$.]
        \item $y_k^i\neqp w_k^i$ and $\hat{y}_k^i\neqp\hat{w}_k^i$;
        \item $x_{k+1}^i\coloneqq y_k^i+\hat{e}_{(b(i)+k-1)q}\notin \cF$ and $\hat{x}_{k+1}^i\coloneqq \hat{y}_k^i+\hat{e}_{(b(i)+k-1)q}\notin \cF$, and 
        \item $\{y_k^i,x_{k+1}^i\},\{\hat{y}_k^i,\hat{x}_{k+1}^i\}\in E(G')$.
    \end{enumerate}
    \item\label{pskelatom3} If $p_i = 0$, let $x_1^i \coloneqq x^i$ and $\hat{x}_1^i \coloneqq \hat{x}^i$.
    For each $k\in[p_i']\setminus[p_i]$, iteratively choose distinct vertices $z_k^i,w_k^i,\hat{z}_k^i,\hat{w}_k^i\in(V(\mathcal{M}_{b(i)+k-1}(v_i))\cap V(\mathcal{M}_{\nu_{k}^i}))\setminus(\cF\cup\{x^i_k,\hat x_k^i\})$ and distinct vertices $y_k^i,\hat{y}_k^i\in V(\mathcal{A}_{(b(i)+k-1)q}(v_i))\setminus(\cF\cup\{z_k^i,w_k^i,\hat{z}_k^i,\hat{w}_k^i\})$ satisfying that
    \begin{enumerate}[label=\ref{pskelatom3}$\arabic*$.]
        \item $z_k^i,\hat{w}_k^i\neqp x_k^i$ and $\hat{z}_k^i,w_k^i\eqp x_k^i$;
        \item\label{weirdnewitem} $x^{j_{k}^i},y^{j_{k}^i},\hat x^{j_{k}^i},\hat y^{j_{k}^i}\notin \cF$, where $x^{j_{k}^i}$, $y^{j_{k}^i}$, $\hat x^{j_{k}^i}$ and $\hat y^{j_{k}^i}$ are defined as in \eqref{equa:skel1}\COMMENT{Note that they form a valid connection sequence for $u_{k}^i$.};
        \item $y_k^i\neqp w_k^i$ and $\hat{y}_k^i\neqp\hat{w}_k^i$;
        \item $x_{k+1}^i\coloneqq y_k^i+\hat{e}_{(b(i)+k-1)q}\notin \cF$ and $\hat{x}_{k+1}^i\coloneqq\hat{y}_k^i+\hat{e}_{(b(i)+k-1)q}\notin \cF$, and 
        \item $\{y_k^i,x_{k+1}^i\},\{\hat{y}_k^i,\hat{x}_{k+1}^i\}\in E(G')$.
    \end{enumerate}
    As discussed earlier, observe that a choice satisfying \ref{weirdnewitem}~exists by \ref{item:Cstarcon1}, \ref{item:Cstarcon2} and \ref{item:preskprop2}.
    \item\label{pskelatom4} For each $k\in[t]\setminus[p_i']$, iteratively choose distinct vertices $y_k^i,\hat{y}_k^i\in V(\mathcal{A}_{(b(i)+k-1)q}(v_i))\setminus \cF$ satisfying that
    \begin{enumerate}[label=\ref{pskelatom4}$\arabic*$.]
        \item $y_k^i\neqp x_k^i$ and $\hat{y}_k^i\neqp\hat{x}_k^i$;
        \item $x_{k+1}^i\coloneqq y_k^i+\hat{e}_{(b(i)+k-1)q}\notin \cF$ and $\hat{x}_{k+1}^i\coloneqq\hat{y}_k^i+\hat{e}_{(b(i)+k-1)q}\notin \cF$, and 
        \item $\{y_k^i,x_{k+1}^i\},\{\hat{y}_k^i,\hat{x}_{k+1}^i\}\in E(G')$.
    \end{enumerate}
\end{enumerate}
Then, we may define the first and second partial skeletons for $v_i$ with connection sequence $(x^i,y^i,\hat x^i,\hat y^i)$ as
\[\mathcal{L}_i(x^i,y^i)\coloneqq(x^i)\left(\bigtimes_{k=1}^{p_i'}(z_k^i,\mathcal{L}_{j_{k}^i}(x^{j_{k}^i},y^{j_{k}^i}),w_k^i,y_k^i,x_{k+1}^i)\right)\Bigg(\bigtimes_{k=p_i'+1}^{t}(y_k^i,x_{k+1}^i)\Bigg)(y^i),\]
\[\hat{\mathcal{L}}_i(\hat x^i,\hat y^i)\coloneqq(\hat x^i)\left(\bigtimes_{k=1}^{p_i'}(\hat{z}_k^i,\hat{\mathcal{L}}_{j_{k}^i}(\hat x^{j_{k}^i},\hat y^{j_{k}^i}),\hat{w}_k^i,\hat{y}_k^i,\hat{x}_{k+1}^i)\right)\Bigg(\bigtimes_{k=p_i'+1}^{t}(\hat{y}_k^i,\hat{x}_{k+1}^i)\Bigg)(\hat y^i).\]

We are now ready to construct $\mathcal{L}$.
The idea is similar to that of Case 3, except that we now tie together the first and second partial skeletons in Step~1.2 below.
Recall from Step 10 that the `root vertex' $v_0$ of the tree $\tau_0$ is atomic.
This will be used in step 3 below.
\begin{enumerate}[label=$\arabic*$.]
    \item\label{skel1} Choose any two vertices $x_1^0,\hat{x}_1^0\in V(\mathcal{A}_1(v_0))\setminus \cF$ such that
    \begin{enumerate}[label=\ref{skel1}$\arabic*$.]
        \item $x_1^0\eqp w_0^0$ and $\hat{x}_1^0\neqp w_0^0$;
        \item $y_t^0\coloneqq\hat{x}_1^0+\hat{e}_{2^s}\notin \cF$ and $\hat{y}_t^0\coloneqq x_1^0+\hat{e}_{2^s}\notin \cF$, and 
        \item $\{x_1^0,\hat{y}_t^0\},\{\hat{x}_1^0,y_t^0\}\in E(G')$.
    \end{enumerate}
    \item\label{skel2} For each $k\in[p_0]$, iteratively choose two distinct vertices $y^0_k, \hat{y}^0_k \in V(\cA_{qk}(v_0))\setminus \cF$ such that
    \begin{enumerate}[label=\ref{skel2}$\arabic*$.]
        \item $y_k^0\neqp w_k^0$ and $\hat{y}_k^0\neqp\hat{w}_k^0$;
        \item $x_{k+1}^0\coloneqq y_k^0+\hat{e}_{kq}\notin \cF$ and $\hat{x}_{k+1}^0\coloneqq\hat{y}_k^0+\hat{e}_{kq}\notin \cF$, and 
        \item $\{y_k^0,x_{k+1}^0\},\{\hat{y}_k^0,\hat{x}_{k+1}^0\}\in E(G')$.
    \end{enumerate}
    \item\label{skel3} For each $k\in[p_0']\setminus[p_0]$, iteratively choose distinct vertices $z_{k}^0,w_{k}^0,\hat{z}_{k}^0,\hat{w}_{k}^0\in(V(\mathcal{M}_{k}(v_0))\cap V(\mathcal{M}_{\nu_{k}^0}))\setminus(\cF\cup\{x_{k}^0,\hat x_{k}^0\})$ and distinct vertices $y_{k}^0,\hat{y}_{k}^0\in V(\mathcal{A}_{kq}(v_0))\setminus(\cF\cup\{z_{k}^0,w_{k}^0,\hat{z}_{k}^0,\hat{w}_{k}^0\})$ satisfying that
    \begin{enumerate}[label=\ref{skel3}$\arabic*$.]
        \item $z_{k}^0,\hat{w}_{k}^0\neqp x_{k}^0$ and $\hat{z}_{k}^0,w_{k}^0\eqp x_{k}^0$;
        \item $x^{j_{k}^0},y^{j_{k}^0},\hat x^{j_{k}^0},\hat y^{j_{k}^0}\notin \cF$, where $x^{j_{k}^0}$, $y^{j_{k}^0}$, $\hat x^{j_{k}^0}$ and $\hat y^{j_{k}^0}$ are defined as in \eqref{equa:skel1}\COMMENT{Note that they form a valid connection sequence for $u_{k}^i$.};
        \item $y_{k}^0\neqp w_{k}^0$ and $\hat{y}_{k}^0\neqp\hat{w}_{k}^0$;
        \item $x_{k+1}^0\coloneqq y_{k}^0+\hat{e}_{kq}\notin \cF$ and $\hat{x}_{k+1}^0\coloneqq \hat{y}_{k}^0+\hat{e}_{kq}\notin \cF$, and 
        \item $\{y_{k}^0,x_{k+1}^0\},\{\hat{y}_{k}^0,\hat{x}_{k+1}^0\}\in E(G')$.
    \end{enumerate}
    \item\label{skel4} For each $k\in[t-1]\setminus[p_0']$, iteratively choose any two vertices $y_k^0,\hat{y}_k^0\in V(\mathcal{A}_{kq}(v_0))\setminus \cF$ satisfying that
    \begin{enumerate}[label=\ref{skel4}$\arabic*$.]
        \item $y_k^0\neqp x_k^0$ and $\hat{y}_k^0\neqp\hat{x}_k^0$;
        \item $x_{k+1}^0\coloneqq y_k^0+\hat{e}_{kq}\notin \cF$ and $\hat{x}_{k+1}^0\coloneqq \hat{y}_k^0+\hat{e}_{kq}\notin \cF$, and 
        \item $\{y_k^0,x_{k+1}^0\},\{\hat{y}_k^0,\hat{x}_{k+1}^0\}\in E(G')$.
    \end{enumerate}
\end{enumerate}
The final definition of $\mathcal{L}$ is given by
\begin{align*}
    \mathcal{L}\coloneqq(x_1^0)&\left(\bigtimes_{k=1}^{p_0'}(z_k^0,\mathcal{L}_{j_k^0}(x^{j_k^0},y^{j_k^0}),w_k^0,y_k^0,x_{k+1}^0)\right)\Bigg(\bigtimes_{k=p_0'+1}^{t-1}(y_k^0,x_{k+1}^0)\Bigg)(y_{t}^0,\hat{x}_1^0)\\
    &\left(\bigtimes_{k=1}^{p_0'}(\hat{z}_k^0,\hat{\mathcal{L}}_{j_k^0}(\hat{x}^{j_k^0},\hat{y}^{j_k^0}),\hat{w}_k^0,\hat{y}_k^0,\hat{x}_{k+1}^0)\right)\Bigg(\bigtimes_{k=p_0'+1}^{t-1}(\hat{y}_k^0,\hat{x}_{k+1}^0)\Bigg)(\hat{y}_{t}^0).
\end{align*}

Observe that \ref{itemskelprop1}--\ref{itemskelprop5} hold by construction.
In particular, \eqref{equa:Delta} together with \ref{itm:CS1} ensure that in Case~3 the final two vertices of the two partial skeletons satisfy $x_{t+1}^i\neqp y^i$ and $\hat x_{t+1}^i\neqp\hat y^i$.
Moreover, the pairs $x^0_t,y^0_t$ and $\hat x^0_t,\hat y^0_t$ will play the roles of the pairs $x_{k_1},x_{k_1+1}$ and $x_{k_2},x_{k_2+1}$ in the second part of \ref{itemskelprop5}.
Similarly, \ref{itemskelprop6} holds by combining the construction of $\mathcal{L}$, \ref{itm:RL5}, \ref{item:preskprop5} and the definition of $V^{\mathrm{abs}}$.

Recall that we write $\mathcal{L}=(x_1,\ldots,x_r)$.
For each $i\in[m']_0$ such that $v_i$ is an atomic vertex and each $j\in[t]$, let $\mathfrak{J}_{i,j}\coloneqq\{k\in[r]:x_k,x_{k+1}\in V(\mathcal{M}_j(v_i))\}$ and $S_{i,j}\coloneqq\{\{x_k,x_{k+1}\}:k\in\mathfrak{J}_{i,j}\}$.\\

\textbf{Step~15: Constructing an almost spanning cycle.}
We will now apply the connecting lemmas to obtain an almost spanning cycle in $G'$ from $\mathcal{L}=(x_1,\ldots,x_r)$.
For each $i\in[m']_0$ such that $v_i$ is an atomic vertex and each $j\in[t]$, except the pair $(0,t)$, we apply \cref{lem:slicecover} to the slice $\mathcal{M}_j(v_i)$ and the graph $G'$, with $\mathfrak{L}\cap V(\mathcal{M}_j(v_i))$, $(\mathfrak{R}_1\cup\mathfrak{R}_2)\cap V(\mathcal{M}_j(v_i))$ and $S_{i,j}$ playing the roles of $L$, $R$ and the pairs of vertices described in \cref{lem:slicecover}\ref{itm:conn1.3}, respectively.
Note that the conditions of \cref{lem:slicecover} can be verified as follows.
\ref{itm:conn1.1} and \ref{itm:conn1.2} hold by \ref{item:Cstarcon1} and \ref{item:Cstarcon2} combined with \ref{item:Cstarcon3}.
\ref{itm:conn1.3} holds by \ref{itemskelprop1} and \ref{itemskelprop3}--\ref{itemskelprop6}.
For $\mathcal{M}_t(v_0)$, we apply \cref{lem:slicecover} or \cref{lem:slicecover2} depending on whether \ref{itemskelprop52} or \ref{itemskelprop51} holds in \ref{itemskelprop5} (the conditions for \cref{lem:slicecover2} can be checked analogously).
For each $i\in[m']_0$ such that $v_i$ is an atomic vertex and each $j\in[t]$, this yields $|\mathfrak{J}_{i,j}|$ vertex-disjoint paths $(\mathcal{P}^{i,j}_k)_{k\in\mathfrak{J}_{i,j}}$ in $\mathcal{M}_j(v_i)\cup G'=G'$ such that, for each $k\in\mathfrak{J}_{i,j}$,
\begin{enumerate}[label=$(\mathrm{\roman*})$]
    \item\label{itm:AlmHCyc1} $\mathcal{P}^{i,j}_k$ is an $(x_k,x_{k+1})$-path,
    \item\label{itm:AlmHCyc2} $\bigcup_{k\in\mathfrak{J}_{i,j}}V(\mathcal{P}^{i,j}_k)=V(\mathcal{M}_j(v_i))\setminus\mathfrak{L}$, and
    \item\label{itm:AlmHCyc3} any pair of second and third absorber vertices in $\mathfrak{R}_1\cup\mathfrak{R}_2$ contained in the same atom of $\mathcal{M}_j(v_i)$ form an edge in one of the paths\COMMENT{The Hamilton cycle does not contain these edges, but it is convenient to formulate \ref{item:almHC3}.
    Furthermore, these edges do belong to the almost spanning cycle that we construct for \cref{thm:almost}.}.
\end{enumerate}

Now consider the path obtained as follows by going through $\mathcal{L}$.
Start with $x_1$. 
For each $k\in[r]$, if there exist $i\in[m']_0$ and $j\in[t]$ such that $\{x_k,x_{k+1}\}\in S_{i,j}$, add $\mathcal{P}^{i,j}_k$ to the path; otherwise, add the edge $\{x_k,x_{k+1}\}$ (this must be an edge of $G'$ by \ref{itemskelprop3}).
Finally, add the edge $\{x_r,x_1\}$ of $G'$ (this is given by \ref{itemskelprop2}) to the path to close it into a cycle $\mathfrak{H}$ in $G'$.
This cycle satisfies the following properties (recall that $e_\mathrm{abs}(u)$ was defined at the end of Step~13):
\begin{enumerate}[label=$(\mathrm{HC}\arabic*)$]
    \item\label{item:almHC1} $|V(\mathfrak{H})|\geq(1-4/\ell^4)2^n$.
    \item\label{item:almHC2} $V(\mathfrak{H})\cupdot\mathfrak{L}\cupdot V^\mathrm{abs}$ partitions $V(\cQ^n)$.
    \item\label{item:almHC3} For all $u\in V^\mathrm{abs}$, we have that $e_\mathrm{abs}(u)\in E(\mathfrak{H})$.
\end{enumerate}
Indeed, note that $\mathfrak{H}$ covers all vertices in $L^\bullet$ (since $L^\bullet\subseteq\mathcal{L}^\bullet$ by \ref{itemskelprop6}) as well as all vertices lying in cube molecules represented in $\tau_0'$ except for those in $\mathfrak{L}$ (by \ref{itm:AlmHCyc2}).
Together with the definition of $V^\mathrm{abs}$, this implies \ref{item:almHC2}.
Moreover, since $|\mathfrak{L}|=|V^\mathrm{abs}|$, \ref{item:almHC1} follows from \eqref{equa:step11}.
Finally, \ref{item:almHC3} follows by \ref{itm:AlmHCyc3}.\\

\textbf{Step~16: Absorbing vertices to form a Hamilton cycle.}
For each $u\in V^\mathrm{abs}$, replace the edge $e_\mathrm{abs}(u)$ by the path $\mathcal{P}_\mathrm{abs}(u)$ (recall from the end of Step~13 that $\mathcal{P}_\mathrm{abs}(u)$ lies in $((H\cup G)\setminus F)\cup G'$).
Clearly, this incorporates all vertices of $\mathfrak{L}\cup V^\mathrm{abs}$ into the cycle and, by \ref{item:almHC2} and \ref{item:almHC3}, the resulting cycle is Hamiltonian.
\end{proof}


\subsection{Proofs of \texorpdfstring{\cref{thm:almost,thm:main,thm:thresholdk}}{Theorems 1.2, 1.7 and 1.5}}\label{sect:thm1}

First, we show that, as a byproduct of the proof of \cref{thm:main1}, we also have a proof of \cref{thm:almost}.

\begin{proof}[Proof of \cref{thm:almost}]
Apply Steps 1, 4, 6, 7, 10, 11, 12, 14 and 15 in succession.
In general, any reference to absorbing cubes in these steps (see e.g.~the end of Step~6) should be skipped as well.\COMMENT{Should we mention that we could ignore Step~12 (defining $\tau_0'$ at Step~10) and do steps 10 and 13 simultaneously? Or would this be more confusing?}
\end{proof}

Next, we will show how \cref{thm:main1} can be used to prove \cref{thm:main}.

\begin{proof}[Proof of \cref{thm:main}]
Consider a decomposition of $H$ into $k$ edge-disjoint subgraphs $H_1 \cup \dots \cup H_k$ such that, for every $i \in [k]$, we have $\delta(H_i)\geq \alpha n/(2k)$. 
To see that this is possible, let us randomly partition the edges of $H$ so that each $e\in E(H)$ is assigned to one of the $H_i$'s uniformly at random and independently from all other edges.
Thus, for every $i \in [k]$ we have $\mathbb{P}[e \in E(H_i)]=1/k$. 
It follows by \cref{lem:Chernoff} that, for every vertex $x\in V(\cQ^n)$ and every $i\in[k]$\COMMENT{It is clear that $\alpha n/k\leq\mathbb{E}[d_{H_i}(x)]\leq n/k$.
By \cref{lem:Chernoff}, it follows that 
\[\mathbb{P}[d_{H_i}(x)\leq \alpha n/(2k)]\leq\mathbb{P}[d_{H_i}(x)\leq\mathbb{E}[d_{H_i}(x)]/2]\leq e^{-(1/2)^2\mathbb{E}[d_{H_i}(x)]/2}\leq e^{-\alpha n/(8k)}.\]
}, 
\[\mathbb{P}[d_{H_i}(x)\leq \alpha n/(2k)]\leq e^{-\alpha n/(8k)}.\]
For each $x\in V(\cQ^n)$, let $\mathcal{B}(x)$ be the event that $d_{H_i}(x)\leq\alpha n/(2k)$ for some $i\in[k]$. 
Hence, $\mathbb{P}[\mathcal{B}(x)]\leq ke^{-\alpha n/(8k)}$ for all $x\in V(\cQ^n)$.
Observe that $\mathcal{B}(x)$ is independent of the collection of events $\{\mathcal{B}(y):\dist(x,y)\geq2\}$.
A simple application of \cref{lem: LLL} shows that 
\[\mathbb{P}\left [\bigwedge_{x \in V(\cQ^n)} \overline{\mathcal{B}(x)}\right]>0\]
and, therefore, such a decomposition of $H$ exists. 

We now consider a similar decomposition of $\cQ^n_\varepsilon$.
In particular, given $\cQ_{\eps}^{n}$, we partition its edges into $k$ edge-disjoint subgraphs, $Q_1\cup \dots \cup Q_k$, in such a way that, if $e\in E(\cQ_{\eps}^{n})$, then $e$ is assigned to one of the $Q_i$ chosen uniformly at random and independently of all other edges.
Thus, for each $e\in E(\cQ_{\eps}^{n})$ we have $\mathbb{P}[e \in E(Q_i)]=1/k$ for all $i\in [k]$.
It follows that, for each $i\in [k]$, we have $Q_i\sim\cQ^{n}_{\eps/k}$.

Let $\Phi$ be a constant such that \cref{thm:main1} holds with $\varepsilon/k$, $\alpha/(2k)$ and $k+2$ playing the roles of $\varepsilon$, $\alpha$ and $c$, respectively.
For each $i\in[k]$, apply \cref{thm:main1} with $H_i$ and $Q_i$ playing the roles of $H$ and $G$, respectively.
We obtain that a.a.s.~there exists a subgraph $G_i\subseteq Q_i$ with $\Delta(G_i)\leq\Phi$ such that, for every $F_i\subseteq \cQ^{n}$ with $\Delta(F_i)\leq (k+2)\Phi$, the graph $((H_i\cup Q_i)\setminus F_i))\cup G_i$ is Hamiltonian.
Condition on the event that this holds for all $i\in[k]$ simultaneously (which holds a.a.s.~by a union bound).
Note that, since the graphs $Q_i$ are pairwise edge-disjoint, so are the different graphs $G_i$.

We are now going to find $k$ edge-disjoint Hamilton cycles $C_1,\ldots,C_k$ iteratively.
For each $i\in[k]$, we proceed as follows.
Let $F_i\coloneqq\bigcup_{j=1}^{k}G_j\cup\bigcup_{j=1}^{i-1}C_j$.
It is clear by construction that $\Delta(F_i)\leq k(\Phi+2)\leq(k+2)\Phi$.
By the conditioning above, there must be a Hamilton cycle $C_i\subseteq ((H_i\cup Q_i)\setminus F_i))\cup G_i$.
Take any such $C_i$ and proceed.

It remains to prove that $C_1,\ldots,C_k$ are pairwise edge-disjoint.
In order to see this, suppose that there exist $i,j\in[k]$ with $i<j$ such that $E(C_i)\cap E(C_j)\neq \varnothing$, and let $e\in E(C_i)\cap E(C_j)$.
In order to have $e\in E(C_i)$, since $G_j\subseteq F_i\setminus G_i$, we must have $e\notin E(G_j)$.
However, since $e\in F_{j}$ by definition, we must have $e \in E(G_j)$, a contradiction.
\end{proof}

Now, \cref{thm:thresholdk} follows as an immediate corollary. 

\begin{proof}[Proof of \cref{thm:thresholdk}]
It is well known (see e.g. \cite{Bol90}) and easy to show that $\cQ^n_{1/2 - \eps}$ a.a.s.~contains isolated vertices.
So it suffices to consider $\cQ^n_{1/2+\eps}$ for any fixed $\eps>0$ and show that a.a.s.~it contains $k$ edge-disjoint Hamilton cycles.
Let $0<\delta \ll \eps\leq1/2$. 
Let $H\sim\cQ^n_{1/2+\eps/2}$ and $G\sim\cQ^n_{\eps/2}$.
Note that $H\cup G\sim\cQ^n_{\eta}$, for some $\eta\leq1/2+\varepsilon$.
Furthermore, by \cref{lem:mindegQeps}, a.a.s.~$\delta(H)\geq\delta n$.
Applying \cref{thm:main} to $H\cup G$, we obtain the desired result.
\end{proof}


\section{Hitting time result}\label{sect:hitting}
  
In this section we prove \cref{thm: kedgehit} and \cref{thm:hitting}.
The proof of \cref{thm: kedgehit} closely follows that of \cref{thm:main}, with some necessary changes to deal with the possibility that some vertices have degree $o(n)$.
We will describe here all necessary changes.
It is worth noting that the set $\cA$ in \cref{lem: main treereshit} is only needed for this section.

\Cref{thm: kedgehit} will actually follow from a slightly more general result (see \cref{thm: important}), which we prove in \cref{section:hitmain}.
Similarly as for \cref{thm:main}, we will first prove \cref{thm: kedgehit} for the property of being Hamiltonian (that is, the case $k=2$), and then use this to prove the general result.

Let $s\in\mathbb{N}$.
Throughout this section, we consider a similar setup to that in \cref{section8}.
In particular, we have a partition of the vertex set of the hypercube into layers $L_1,\ldots,L_{2^s}$, where the labelling of the layers is given by some Hamilton cycle of $\cQ^s$.
We fix any such partition, and will write \index{ehatk@$\hat e_k$}$\hat e_i$ for the direction of the edges between the layers $L_i$ and $L_{i+1}$.
We will use the notation for atoms, molecules and slices in the same way.
We will also use part of the notation introduced at the beginning of the proof of \cref{thm:main1}; in particular, we will consider the intersection graph $I\cong\cQ^{n-s}$ of the layers, and for any $G\subseteq\cQ^n$ we define $I(G)$ analogously.

\subsection{Absorbing structures for vertices with small degree}\label{sect:absstruct}

The main reason why the proof of \cref{thm:main1} does not work for the case $k=2$ of \cref{thm: kedgehit} is the existence of vertices of very low degree (as low as degree $2$).
We cannot hope to absorb these via absorbing $\ell$-cube pairs as in the proof of \cref{thm:main1}, as they may not have any neighbours which lie in the cube factor we construct. 
Thus, we first construct alternative absorbing structures for these vertices.

To be more precise, recall that in the proof of \cref{thm:main1} we absorbed vertices in pairs, in order to compensate the parities (within a vertex molecule).
Here, we will need something similar.
To achieve this, we will define several paths.
One path will contain the vertex of low degree, while the others are used to compensate the parities of vertices in this first path.
Moreover, the paths will be constructed in such a way that they end in vertices which can be paired up so that each pair consists of clones of a given vertex of the intersection graph $I$ and these two clones lie in a (bonded) cube molecule.
These cube molecules can then be used to connect these paths.
It is worth mentioning, however, that we cannot guarantee that the pairing of the vertices can be done within a single slice, so we need to alter our approach to deal with this.

These special absorbing structures will be used in a somewhat different way to the absorbing $\ell$-cube pairs (which are still used to absorb all other vertices).
For the latter, recall that we enforce a condition on the near-spanning cycle $\mathfrak{H}$ (namely, that it contains certain edges $e_\mathrm{abs}(u)$) so that we can absorb the required vertices $u\in V^\mathrm{abs}$.
For the special absorbing structures, however, we will actually enforce that the vertices of very low degree are already incorporated into $\mathfrak{H}$.
In particular, consider the cube molecules which contain the endpoints of paths of one of these special absorbing structures.
When constructing the skeleton, we will add extra segments in these molecules which will be used to connect these paths in such a way that, when we apply the connecting lemmas, we can completely incorporate each special absorbing structure into $\mathfrak{H}$.

Suppose $x\in V(\cQ^n)$ is incident to one edge with direction $a$ and one edge with direction $b$.
We will construct three different types of special absorbing structures, to handle the cases where both, only one, or none of these two edges lie in the same layer as $x$.
Representations of these three types of special absorbing structures can be found in \cref{fig:SAS}.
Let $L$ be the layer such that $x\in V(L)$.
Given a path $P=x_1x_2\ldots x_k$ in $\cQ^n$, we define $\mathrm{end}(P)\coloneqq \{x_1,x_k\}$.

\textbf{Type~I.}
Assume that $x+a,x+b\in V(L)$.
Let $f\colon V(\cQ^n)\to V(\cQ^n)$ be defined as follows: for each $i\in[2^s]$ and each $y\in V(L_i)$, if $i$ is even, we set $f(y)\coloneqq y+\hat e_{i-1}$; otherwise, we set $f(y)\coloneqq y+\hat e_i$\COMMENT{So this just takes each vertex to its copy in the next or previous layer, in such a way that $f$ is an involution.}.
By abusing notation, for any $F\subseteq L_i$, we also consider the graph $f(F)$, where for each edge $e=\{y,z\}\in E(F)$ we define $f(e)\coloneqq\{f(y),f(z)\}$.
Let $(c,d,d_1,d_2,d_3,d_4)\in(\mathcal{D}(L)\setminus\{a,b\})^6$ be a tuple of distinct directions and define the following paths in $\cQ^n$:
\begin{itemize}
    \item $P_1\coloneqq(x+a+d_1,x+a,x,x+b,x+b+d_2)$;
    \item $P_2\coloneqq(f(x+b+d_2),f(x+b),f(x+b+c))$;
    \item $P_3\coloneqq(x+c+b,x+c,x+c+d_3)$;
    \item $P_4\coloneqq(f(x+c+d_3),f(x+c),f(x),f(x+d),f(x+d+d_4))$;
    \item $P_5\coloneqq(x+d+d_4,x+d,x+d+a)$;
    \item $P_6\coloneqq(f(x+a+d),f(x+a),f(x+a+d_1))$.
\end{itemize}
Observe that, for every $y\in V(P_1\cup\dots\cup P_6)$, we have that $f(y)=f^{-1}(y)\in V(P_1\cup\dots\cup P_6)$.
We say that $\mathit{CS}(x,a,b)\coloneqq(P_1,\ldots,P_6)$\index{CSx@$\mathit{CS}(x,a,b)$} is an $(x,a,b)$-\emph{consistent system of paths}, and let $\mathrm{end}(\mathit{CS}(x,a,b))\coloneqq \bigcup_{i=1}^{6}\mathrm{end}(P_i)$\index{endCS@$\mathrm{end}(\mathit{CS})$}.

Let $(P_1,\ldots,P_6)$ be an $(x,a,b)$-consistent system of paths.
Let $\mathcal{D}\coloneqq\{c,d,d_1,d_2,d_3,d_4\}\subseteq\mathcal{D}(L)$ be the set of directions such that $P_1,\ldots,P_6$ are as defined above.
Let $\mathbf{C}\coloneqq\{C_1,\ldots,C_{12}\}$ be a collection of vertex-disjoint $\ell$-cubes which satisfy the following:
\begin{enumerate}[label=$(\mathrm{PI}.\arabic*)$]
    \item\label{itm:specialp1} for all $i\in[12]$, we have that either $C_i\subseteq L$ or $C_i\subseteq f(L)$;
    \item\label{itm:specialp2} for all $i\in[6]$, we have $f(C_{2i})=C_{2i+1}$, where indices are taken modulo $12$\COMMENT{In particular, this property, together with the fact that the cubes are vertex-disjoint, implies that the projection of these cubes to the intersection graph gives 6 cubes which are vertex-disjoint; equivalently, we are looking at 6 vertex-disjoint cube-molecules containing each of the 12 cubes.};
    \item\label{itm:specialp3} for all $i\in[6]$, $C_{2i-1}$ contains the first vertex of $P_{i}$, and $C_{2i}$ contains the last vertex of $P_i$;
    \item\label{itm:specialp4} for all $i\in[12]$, we have $\mathcal{D}(C_i)\cap(\mathcal{D}\cup\{a,b\})=\varnothing$.
\end{enumerate}
We say that $\mathit{SA}(x,a,b)\coloneqq(P_1,\ldots,P_6,C_1,\ldots,C_{12})$\index{SAx@$\mathit{SA}(x,a,b)$} is an $(x,a,b)$-\emph{special absorbing structure}.

Finally, given any graph $G\subseteq\cQ^n$ and an $(x,a,b)$-consistent system of paths $\mathit{CS}(x,a,b)=(P_1,\ldots,P_6)$, we say that $\mathit{CS}(x,a,b)$ \emph{extends} to an $(x,a,b)$-special absorbing structure in $G$ if there is a collection $\mathbf{C}=\{C_1,\ldots,C_{12}\}$ with $C_i\subseteq G$ such that  $\mathit{SA}(x,a,b)=(P_1,\ldots,P_6,C_1,\ldots,C_{12})$ is an $(x,a,b)$-special absorbing structure.

\textbf{Type~II.}
Assume now that $x+a,x+b\notin V(L)$.
Let $(d_1,d_2)\in(\mathcal{D}(L))^2$ be a pair of distinct directions and define the following paths in $\cQ^n$:
\begin{itemize}
    \item $P_1\coloneqq(x+a+d_1,x+a,x,x+b,x+b+d_2)$;
    \item $P_2\coloneqq(x+a+b+d_2,x+a+b,x+a+b+d_1)$.
\end{itemize}
We say that $\mathit{CS}(x,a,b)\coloneqq(P_1,P_2)$ is an $(x,a,b)$-\emph{consistent system of paths}, and let $\mathrm{end}(\mathit{CS}(x,a,b))\coloneqq\mathrm{end}(P_1)\cup\mathrm{end}(P_2)$.

Let $(P_1,P_2)$ be an $(x,a,b)$-consistent system of paths.
Let $\mathcal{D}\coloneqq\{d_1,d_2\}\subseteq\mathcal{D}(L)$ be the set of directions such that $P_1$ and $P_2$ are as defined above.
Let $L_{ab}$, $L_a$ and $L_b$ be the layers such that $x+a+b\in V(L_{ab})$, $x+a\in V(L_a)$ and $x+b\in V(L_b)$, respectively.
Let $\mathbf{C}\coloneqq\{C_1,C_2,C_3,C_4\}$ be a collection of vertex-disjoint $\ell$-cubes which satisfy the following:
\begin{enumerate}[label=$(\mathrm{PII}.\arabic*)$]
    \item\label{item:type21} $C_1\subseteq L_a$, $C_2\subseteq L_b$ and $C_3,C_4\subseteq L_{ab}$;
    \item\label{item:type22} $C_1+b=C_4$ and $C_2+a=C_3$\COMMENT{As happened before, the first two properties, together with the fact that the cubes are vertex-disjoint, implies that the projection of these cubes to the intersection graph gives two vertex-disjoint cubes.};
    \item\label{item:type23} for all $i\in [2]$, $C_{2i-1}$ contains the first vertex of $P_{i}$, and $C_{2i}$ contains the last vertex of $P_{i}$;
    \item\label{item:type24} for all $i\in[4]$, we have $\mathcal{D}(C_i)\cap\mathcal{D}=\varnothing$.
\end{enumerate}
We say that $\mathit{SA}(x,a,b)\coloneqq(P_1,P_2,C_1,\ldots,C_4)$ is an $(x,a,b)$-\emph{special absorbing structure}.

Finally, given any graph $G\subseteq\cQ^n$ and an $(x,a,b)$-consistent system of paths $\mathit{CS}(x,a,b)=(P_1,P_2)$, we say that $\mathit{CS}(x,a,b)$ \emph{extends} to an $(x,a,b)$-special absorbing structure in $G$ if there is a collection $\mathbf{C}=\{C_1,\ldots,C_4\}$ with $C_i\subseteq G$ such that  $\mathit{SA}(x,a,b)=(P_1,P_2,C_1,\ldots,C_4)$ is an $(x,a,b)$-special absorbing structure.

\textbf{Type~III.}
Finally, assume that $x+a\notin V(L)$ and $x+b\in V(L)$\COMMENT{The case when $x+a\in V(L)$ and $x+b\notin V(L)$ is the same, by interchanging the roles of $a$ and $b$.}.
For each vertex $y\in V(L)$, let $f(y)\coloneqq y+a$.
By abusing notation, for any $F\subseteq L$, we also consider the graph $f(F)$, where for each edge $e=\{y,z\}\in E(F)$ we define $f(e)\coloneqq\{f(y),f(z)\}$.
Let $(d_1,d_2,d_3)\in(\mathcal{D}(L)\setminus\{b\})^3$ be a tuple of distinct directions and define the following paths in $\cQ^n$:
\begin{itemize}
    \item $P_1\coloneqq(f(x+d_1+d_2),f(x+d_1),f(x),x,x+b,x+b+d_3)$;
    \item $P_2\coloneqq(f(x+b+d_3),f(x+b),f(x+b+d_1))$;
    \item $P_3\coloneqq(x+d_1+b,x+d_1,x+d_1+d_2)$.
\end{itemize}
We say that $\mathit{CS}(x,a,b)\coloneqq(P_1,P_2,P_3)$ is an $(x,a,b)$-\emph{consistent system of paths}, and let $\mathrm{end}(\mathit{CS}(x,a,b))\coloneqq\mathrm{end}(P_1)\cup\mathrm{end}(P_2)\cup\mathrm{end}(P_3)$.

Let $(P_1,P_2,P_3)$ be an $(x,a,b)$-consistent system of paths.
Let $\mathcal{D}\coloneqq\{d_1,d_2,d_3\}\subseteq\mathcal{D}(L)$ be the set of directions such that $P_1$, $P_2$ and $P_3$ are as defined above.
Let $\mathbf{C}\coloneqq\{C_1,\ldots,C_6\}$ be a set of vertex-disjoint $\ell$-cubes which satisfy the following:
\begin{enumerate}[label=$(\mathrm{PIII}.\arabic*)$]
    \item\label{item:type31} for all $i\in[6]$, we have that $C_i\subseteq L$ or $C_i\subseteq f(L)$;
    \item\label{item:type32} for all $i\in[3]$, we have $C_{2i}=f(C_{2i+1})$, where indices are taken modulo $6$;
    \item\label{item:type33} for all $i\in[3]$, $C_{2i-1}$ contains the first vertex of $P_{i}$, and $C_{2i}$ contains the last vertex of $P_{i}$;
    \item\label{item:type34} for all $i\in[6]$, we have $\mathcal{D}(C_i)\cap(\mathcal{D}\cup\{b\})=\varnothing$.
\end{enumerate}
We say that $\mathit{SA}(x,a,b)\coloneqq(P_1,P_2,P_3,C_1,\ldots,C_6)$ is an $(x,a,b)$-\emph{special absorbing structure}.

Finally, given any graph $G\subseteq\cQ^n$ and an $(x,a,b)$-consistent system of paths $\mathit{CS}(x,a,b)=(P_1,P_2,P_3)$, we say that $\mathit{CS}(x,a,b)$ \emph{extends} to an $(x,a,b)$-special absorbing structure in $G$ if there is a collection $\mathbf{C}=\{C_1,\ldots,C_6\}$ with $C_i\subseteq G$ such that  $\mathit{SA}(x,a,b)=(P_1,P_2,P_3,C_1,\ldots,C_6)$ is an $(x,a,b)$-special absorbing structure.\\

\captionsetup[subfigure]{labelformat=empty}

\begin{figure}
    \centering
    \begin{subfigure}[b]{\textwidth}
        \centering
        \begin{tikzpicture}[every node/.style={scale=1.7}]
	\tikzset{e/.style   = {line width=1pt}}
	\tikzset{re/.style   = {line width=.5pt, red}}
	\tikzset{be/.style   = {line width=.5pt, blue}}
	\tikzset{ge/.style   = {line width=.5pt, green!60!black}}
	
\draw[gray, thin] (-18,1.5) -- (-16,4.5) -- (18,4.5) -- (16,1.5) -- (-18,1.5);
\draw[gray, thin] (-18,-3.5) -- (-16,-0.5) -- (18,-0.5) -- (16,-3.5) -- (-18,-3.5);
\node at (17.5,2.7){$f(L)$};
\node at (17.3,-2.3){$L$};

\filldraw[fill=gray!15] (-15,3) -- (-16,3) -- (-16.5,2.3) -- (-15.5,2.3) -- (-15,3);
\filldraw[fill=gray!15] (-15,-2) -- (-16,-2) -- (-16.5,-2.7) -- (-15.5,-2.7) -- (-15,-2);
\node (fxdd4) at (-15,3){} ;
\draw (-15,3) node[circle,fill,inner sep=2pt]{};
\node (xdd4) at (-15,-2){} ;
\draw (-15,-2) node[circle,fill,inner sep=2pt]{};
\draw[dashed, gray] (-15,3) -- (-15,-2);
\node[anchor = south east] at (-15,-2){$x+d+d_4$};

\node (fxd) at (-12,3){};
\draw (-12,3) node[circle,fill,inner sep=2pt]{};
\node (xd) at (-12,-2){};
\draw (-12,-2) node[circle,fill,inner sep=2pt]{};
\draw[dashed, gray] (-12,3) -- (-12,-2);
\node[anchor = north] at (-12,-2){$x+d$};

\filldraw[fill=gray!15] (-9,3) -- (-9.9,2.7) -- (-10.4,1.9) -- (-9.5,2.2) -- (-9,3);
\filldraw[fill=gray!15] (-9,-2) -- (-9.9,-2.3) -- (-10.4,-3.1) -- (-9.5,-2.8) -- (-9,-2);
\node (fxda) at (-9,3){};
\draw (-9,3) node[circle,fill,inner sep=2pt]{};
\node (xda) at (-9,-2){};
\draw (-9,-2) node[circle,fill,inner sep=2pt]{};
\draw[dashed, gray] (-9,3) -- (-9,-2);
\node[anchor = north west] at (-9.1,-2){$x+d+a$};

\filldraw[fill=gray!15] (-6,3) -- (-6.5,3.7) -- (-7.5,3.7) -- (-7,3) -- (-6,3);
\filldraw[fill=gray!15] (-6,-2) -- (-6.5,-1.3) -- (-7.5,-1.3) -- (-7,-2) -- (-6,-2);
\node (fxad1) at (-6,3){} ;
\draw (-6,3) node[circle,fill,inner sep=2pt]{};
\node (xad1) at (-6,-2){} ;
\draw (-6,-2) node[circle,fill,inner sep=2pt]{};
\draw[dashed, gray] (-6,3) -- (-6,-2);
\node[anchor = south west] at (-6.1,-2){$x+a+d_1$};

\node (fxa) at (-3,3){};
\draw (-3,3) node[circle,fill,inner sep=2pt]{};
\node (xa) at (-3,-2){};
\draw (-3,-2) node[circle,fill,inner sep=2pt]{};
\draw[dashed, gray] (-3,3) -- (-3,-2);
\node[anchor = north] at (-3,-2.05){$x+a$};

\node (fx) at (0,3){};
\draw (0,3) node[circle,fill,inner sep=2pt]{};
\node (x) at (0,-2){};
\draw (0,-2) node[circle,fill,inner sep=2pt]{};
\draw[dashed, gray] (0,3) -- (0,-2);
\node[anchor = north] at (0,-2.13){$x$};
\node[anchor = south] at (0,3.1){$f(x)$};

\node (fxb) at (3,3){};
\draw (3,3) node[circle,fill,inner sep=2pt]{};
\node (xb) at (3,-2){};
\draw (3,-2) node[circle,fill,inner sep=2pt]{};
\draw[dashed, gray] (3,3) -- (3,-2);
\node[anchor = north] at (3,-2){$x+b$};

\filldraw[fill=gray!15] (6,3) -- (6.5,3.7) -- (7.5,3.7) -- (7,3) -- (6,3);
\filldraw[fill=gray!15] (6,-2) -- (6.5,-1.3) -- (7.5,-1.3) -- (7,-2) -- (6,-2);
\node (fxbd2) at (6,3){};
\draw (6,3) node[circle,fill,inner sep=2pt]{};
\node (xbd2) at (6,-2){};
\draw (6,-2) node[circle,fill,inner sep=2pt]{};
\draw[dashed, gray] (6,3) -- (6,-2);
\node[anchor = south east] at (6.1,-2){$x+b+d_2$};

\filldraw[fill=gray!15] (9,3) -- (9.9,2.7) -- (10.4,1.9) -- (9.5,2.2) -- (9,3);
\filldraw[fill=gray!15] (9,-2) -- (9.9,-2.3) -- (10.4,-3.1) -- (9.5,-2.8) -- (9,-2);
\node (fxcb) at (9,3){};
\draw (9,3) node[circle,fill,inner sep=2pt]{};
\node (xcb) at (9,-2){};
\draw (9,-2) node[circle,fill,inner sep=2pt]{};
\draw[dashed, gray] (9,3) -- (9,-2);
\node[anchor = north east] at (9.1,-2){$x+c+b$};

\node (fxc) at (12,3){};
\draw (12,3) node[circle,fill,inner sep=2pt]{};
\node (xc) at (12,-2){};
\draw (12,-2) node[circle,fill,inner sep=2pt]{};
\draw[dashed, gray] (12,3) -- (12,-2);
\node[anchor = north] at (12,-2.05){$x+c$};

\filldraw[fill=gray!15] (15,3) -- (16,3) -- (16.5,3.7) -- (15.5,3.7) -- (15,3);
\filldraw[fill=gray!15] (15,-2) -- (16,-2) -- (16.5,-1.3) -- (15.5,-1.3) -- (15,-2);
\node (fxcd3) at (15,3){};
\draw (15,3) node[circle,fill,inner sep=2pt]{};
\node (xcd3) at (15,-2){};
\draw (15,-2) node[circle,fill,inner sep=2pt]{};
\draw[dashed, gray] (15,3) -- (15,-2);
\node[anchor = north] at (15,-2){$x+c+d_3$};

\draw[ultra thick, black] (-6,-2) -- (6,-2);
\node[anchor = south] at (1.5,-2){$P_1$};

\draw[ultra thick, black] (6,3) to (3,3) to[bend right = 25] (9,3);
\node[anchor = north] at (4,2.5){$P_2$};

\draw[ultra thick, black] (9,-2) -- (15,-2);
\node[anchor = south] at (10.5,-2){$P_3$};

\draw[ultra thick, black] (15,3) to (12,3) to[bend right = 20] (0,3) to[bend right = 20] (-12,3) to (-15,3);
\node[anchor = north] at (1.3,3.3){$P_4$};

\draw[ultra thick, black] (-15,-2) -- (-9,-2);
\node[anchor = south] at (-10.5,-2){$P_5$};

\draw[ultra thick, black] (-9,3) to[bend right = 25] (-3,3) to (-6,3);
\node[anchor = north] at (-4,2.5){$P_6$};
\end{tikzpicture}
        \caption{Type I}
        \label{fig:typeI}
    \end{subfigure}
    
    \begin{subfigure}[b]{0.29\textwidth}
        \centering
        \begin{tikzpicture}[every node/.style={scale=1.2}]
	\tikzset{e/.style   = {line width=1pt}}
	\tikzset{re/.style   = {line width=.5pt, red}}
	\tikzset{be/.style   = {line width=.5pt, blue}}
	\tikzset{ge/.style   = {line width=.5pt, green!60!black}}

\filldraw[fill=gray!15] (-3,4) -- (-3.5,3.5) -- (-4,4) -- (-3.5,4.5) -- (-3,4);
\filldraw[fill=gray!15] (-3,2) -- (-3.5,1.5) -- (-4,2) -- (-3.5,2.5) -- (-3,2);
\node (fxad1) at (-3,4){} ;
\draw (-3,4) node[circle,fill,inner sep=2pt]{};
\node (xad1) at (-3,2){} ;
\draw (-3,2) node[circle,fill,inner sep=2pt]{};
\draw[dashed, gray] (-3,4) -- (-3,2);
\node[anchor = north west] at (-3.1,4){$x+a+b+d_1$};
\node[anchor = north west] at (-3.1,2){$x+a+d_1$};

\node (x) at (0,-2){};
\draw (0,-2) node[circle,fill,inner sep=2pt]{};
\node (fx) at (0,0){};
\draw (0,0) node[circle,fill,inner sep=2pt]{};
\node (ffx) at (0,2){};
\draw (0,2) node[circle,fill,inner sep=2pt]{};
\node (fffx) at (0,4){};
\draw (0,4) node[circle,fill,inner sep=2pt]{};
\draw[dashed, gray] (0,4) -- (0,-2);
\node[anchor = west] at (0,-2){$x$};
\node[anchor = south west] at (0,0){$x+b$};
\node[anchor = south west] at (0,2){$x+a$};
\node[anchor = south] at (0,4){$x+a+b$};

\filldraw[fill=gray!15] (3,4) -- (3.5,3.5) -- (4,4) -- (3.5,4.5) -- (3,4);
\filldraw[fill=gray!15] (3,0) -- (3.5,-0.5) -- (4,0) -- (3.5,0.5) -- (3,0);
\node (fxbd2) at (3,4){};
\draw (3,4) node[circle,fill,inner sep=2pt]{};
\node (xbd2) at (3,0){};
\draw (3,0) node[circle,fill,inner sep=2pt]{};
\draw[dashed, gray] (3,4) -- (3,0);
\node[anchor = north east] at (3.1,4){$x+a+b+d_2$};
\node[anchor = north east] at (3.1,0){$x+b+d_2$};

\draw[ultra thick, black] (3,0) to (0,0) to[bend left = 25] (0,-2) to[bend left = 25] (0,2) to (-3,2);
\node[anchor = east] at (-0.45,-0.2){$P_1$};

\draw[ultra thick, black] (-3,4) to (3,4);
\node[anchor = south] at (2,4){$P_2$};
\end{tikzpicture}
        \caption{Type II}
        \label{fig:typeII}
    \end{subfigure}
    \begin{subfigure}[b]{0.69\textwidth}
        \centering
        \begin{tikzpicture}[every node/.style={scale=1.5}]
	\tikzset{e/.style   = {line width=1pt}}
	\tikzset{re/.style   = {line width=.5pt, red}}
	\tikzset{be/.style   = {line width=.5pt, blue}}
	\tikzset{ge/.style   = {line width=.5pt, green!60!black}}
	
\draw[gray, thin] (-9,1.8) -- (-7,4.2) -- (12,4.2) -- (10,1.8) -- (-9,1.8);
\draw[gray, thin] (-9,-1.2) -- (-7,1.2) -- (12,1.2) -- (10,-1.2) -- (-9,-1.2);
\node at (11.4,2.7){$f(L)$};
\node at (11.3,-0.3){$L$};

\filldraw[fill=gray!15] (-6,3) -- (-6.5,2.5) -- (-7,3) -- (-6.5,3.5) -- (-6,3);
\filldraw[fill=gray!15] (-6,0) -- (-6.5,-0.5) -- (-7,0) -- (-6.5,0.5) -- (-6,0);
\draw (-6,3) node[circle,fill,inner sep=2pt]{};
\node (xad1) at (-6,0){} ;
\draw (-6,0) node[circle,fill,inner sep=2pt]{};
\draw[dashed, gray] (-6,3) -- (-6,0);
\node[anchor = north west] at (-6.1,0){$x+d_1+d_2$};

\node (fxa) at (-3,3){};
\draw (-3,3) node[circle,fill,inner sep=2pt]{};
\node (xa) at (-3,0){};
\draw (-3,0) node[circle,fill,inner sep=2pt]{};
\draw[dashed, gray] (-3,3) -- (-3,0);
\node[anchor = south west] at (-3,-0.05){$x+d_1$};

\node (fx) at (0,3){};
\draw (0,3) node[circle,fill,inner sep=2pt]{};
\node (x) at (0,0){};
\draw (0,0) node[circle,fill,inner sep=2pt]{};
\draw[dashed, gray] (0,3) -- (0,0);
\node[anchor = east] at (-0.1,0){$x$};
\node[anchor = south] at (0,3.1){$x+a=f(x)$};

\node (fxb) at (3,3){};
\draw (3,3) node[circle,fill,inner sep=2pt]{};
\node (xb) at (3,0){};
\draw (3,0) node[circle,fill,inner sep=2pt]{};
\draw[dashed, gray] (3,3) -- (3,0);
\node[anchor = north] at (3,0){$x+b$};

\filldraw[fill=gray!15] (6,3) -- (7,3) -- (7.5,3.7) -- (6.5,3.7) -- (6,3);
\filldraw[fill=gray!15] (6,0) -- (7,0) -- (7.5,0.7) -- (6.5,0.7) -- (6,0);
\node (fxbd2) at (6,3){};
\draw (6,3) node[circle,fill,inner sep=2pt]{};
\node (xbd2) at (6,0){};
\draw (6,0) node[circle,fill,inner sep=2pt]{};
\draw[dashed, gray] (6,3) -- (6,0);
\node[anchor = south east] at (6,0){$x+b+d_3$};

\filldraw[fill=gray!15] (9,3) -- (8.5,2.3) -- (9.5,2.3) -- (10,3) -- (9,3);
\filldraw[fill=gray!15] (9,0) -- (8.5,-0.7) -- (9.5,-0.7) -- (10,0) -- (9,0);
\node (fxcb) at (9,3){};
\draw (9,3) node[circle,fill,inner sep=2pt]{};
\node (xcb) at (9,0){};
\draw (9,0) node[circle,fill,inner sep=2pt]{};
\draw[dashed, gray] (9,3) -- (9,0);
\node[anchor = south] at (9.3,0){$x+d_1+b$};

\draw[ultra thick, black] (-6,3) -- (0,3) -- (0,0) -- (6,0);
\node[anchor = west] at (0,1.5){$P_1$};

\draw[ultra thick, black] (6,3) to (3,3) to[bend right = 25] (9,3);
\node[anchor = south] at (4.5,3){$P_2$};

\draw[ultra thick, black] (9,0) to[bend left = 15] (-3,0) to (-6,0);
\node[anchor = south] at (3,-1.7){$P_3$};
\end{tikzpicture}
        \caption{Type III}
        \label{fig:typeIII}
    \end{subfigure}
    \caption{\footnotesize{A representation of the special absorbing structures.
    In each case, vertices (or cubes) represented in vertical lines are clones of the same vertex (or cube) of the intersection graph $I$, and any vertices in the same horizontal line lie in the same layer of $\cQ^n$.}}
    \label{fig:SAS}
\end{figure}

Whenever $x$, $a$ and $b$ are clear from the context, we will simply write $\mathit{CS}$ and $\mathit{SA}$ instead of $\mathit{CS}(x,a,b)$ and $\mathit{SA}(x,a,b)$.
Given any consistent system of paths $\mathit{CS}$, we let $\mathrm{endmol}(\mathit{CS})$\index{endXmol@$\mathrm{endmol}(\mathit{CS})$} be the set of vertices $v \in V(I)$ such that some clone of $v$ lies in $\mathrm{end}(\mathit{CS})$. 
We write $\cD(\mathit{CS})$\index{dirCS@$\cD(\mathit{CS})$} to denote the set of directions $\cD\cup\{a,b\}$ used to define the paths which comprise $\mathit{CS}$ as above.
If $\mathit{CS}$ extends to a special absorbing structure $\mathit{SA}$, we denote $\mathrm{end}(\mathit{SA})\coloneqq \mathrm{end}(\mathit{CS})$\index{endSA@$\mathrm{end}(\mathit{SA})$}.
Moreover, we denote by $\mathbf{C}(\mathit{SA})$\index{CRA@$\mathbf{C}(\mathit{SA})$} the collection of cubes associated with $\mathit{SA}$.
Observe that \ref{itm:specialp4}, \ref{item:type24} and \ref{item:type34} imply that 
\begin{enumerate}[label=(AS)]
    \item\label{item:AbsStruct}each cube $C \in \mathbf{C}(\mathit{SA})$ is vertex-disjoint from the paths in $\mathit{SA}$ except for the unique vertex in $\mathrm{end}(\mathit{SA})$ contained in $C$.
\end{enumerate}
We will sometimes abuse notation and treat $\mathit{CS}$ and $\mathit{SA}$ as graphs; in particular, we will write $V(\mathit{CS})$ to denote the vertices of the union of the paths which comprise $\mathit{CS}$, and $V(\mathit{SA})$ to denote the vertices of the union of the paths and cubes which comprise $\mathit{SA}$, and similarly for $E(\mathit{CS})$ and $E(\mathit{SA})$.


\subsection{Auxiliary lemmas} 

We now state and prove some auxiliary lemmas.
When taking random subgraphs of the hypercube, we will need to guarantee that, given a vertex in $I$ and a large collection of cubes in $I$ incident to this vertex, some of the cube molecules given by these cubes are bonded.

\begin{lemma}\label{lem: bondedlots}
Let $\varepsilon,\gamma\in(0,1)$ and $\ell,n\in\mathbb{N}$ with $0 < 1/n\ll1/\ell\ll\varepsilon,\gamma$, and let $s \coloneqq 10\ell$.
Let $x\in V(I)$ and let $\cC$ be a collection of $\ell$-cubes $C\subseteq I$ such that $|\cC| \ge \gamma n^\ell$ and, for all $C \in \cC$, we have $x \in V(C)$.
For each $C\in\cC$, let $\cM_C$ denote the cube molecule of $C$ in $\cQ^n$.
For any graph $G \subseteq \cQ^n$, let
\[B(G)\coloneqq \{C \in \cC: \cM_C \text{ is bonded in } G\}.\]
Then, with probability at least $1 -2^{-10n}$, we have $|B(\cQ^n_\eps)| \ge \gamma n^\ell/4$.
\end{lemma}

\begin{proof}
Let $G \sim \cQ^n_\eps$ and let $\cC' \coloneqq \{C - (N_{I}(x) \cup \{x\}): C \in \cC\}$\COMMENT{That is, this is the set of all cubes after removing $x$ and its neighbourhood.}.
Given $C' \in \cC'$, for each $i\in[2^s]$, let $C'_i$ be the $i$-th clone of $C'$.
We denote $\cM_{C'}\coloneqq C'_1 \cup \dots \cup C'_{2^s}$, and refer to it as the molecule of $C'$ in $\cQ^n$.
We say that $\cM_{C'}$ is \emph{bonded in $G$} if, for each $i \in [2^s]$, the graph $G$ contains at least $100$ edges between $C'_{i}$ and $C'_{i+1}$ whose endpoints in $C'_{i}$ have odd parity and $100$ edges whose endpoints in $C'_{i}$ have even parity, where indices are taken cyclically\COMMENT{Note that we never claim that the cubes or these $C'$ are in $G$.
In application, we will call on this lemma with a collection of cubes $\cC$ that we know exists already in one layer of probability, and use an extra layer of probability to get bondedness among them.}.
Note that, if $C' = C - (N_{I}(x) \cup \{x\})$ for some $C \in \cC$ and $C'$ is bonded in $G$, then $C$ must be bonded in~$G$.
Moreover, $|V(C')| = |V(C)|-\ell-1>9|V(C)|/10$.
Therefore, similarly to the proof of \cref{lem:moleculegood}, by \cref{lem:Chernoff} and a union bound we have that\COMMENT{For all $i\in[2^s]$, the number of edges in $\cQ^n$ between $C'_i$ and $C'_{i+1}$ with an endpoint of odd parity in $C'_i$ is at least $2^{\ell -1}-\ell$ (and the same bound holds for those edges with an endpoint of even parity in $C'_i$).
Thus, the number of such edges in $G$ stochastically dominates $\sim\mathrm{Bin}(2^{\ell-1}-\ell, \eps)$.
This stochastically dominates a variable $Y\sim\mathrm{Bin}(2^\ell/4, \eps)$.
Now apply, \cref{lem:Chernoff} in the same way as is done in \cref{lem:moleculegood} to obtain that $\mathbb{P}[Y<100]\leq2^{-\varepsilon2^{\ell-3}}$.
The claim follows by a union bound over both parities and $2^s$ (the number of layers of the molecule).}
\[\mathbb{P}[\cM_{C'} \text{ is bondless in }G] \le 2^{s - \eps 2^{\ell}/100} \le 1/2.\]
Let $X\coloneqq|\{C' \in \cC': \cM_{C'} \text{ is bonded in } G\}|$. 
It follows that \begin{equation}\label{eqn:bdlots}\mathbb{E}[|B(G)|] \ge \mathbb{E}[X] \ge \gamma n^{\ell}/2.
\end{equation}

Let $V \coloneqq \bigcup_{C' \in \cC'}V(\cM_{C'})$.
Let $e_1, \dots, e_m$ be an arbitrary ordering of the edges of $E\coloneqq \bigcup_{i \in [2^s]}E_{\cQ^n}(L_i\cap V, L_{i+1}\cap V)$.
For each $j \in [m]$, let $X_j$ be the indicator variable which takes value $1$ if $e_j \in E(G)$ and $0$ otherwise.
Consider the edge-exposure martingale $Y_j \coloneqq \mathbb{E}[X\mid X_1, \dots, X_j]$ for $j \in [m]_0$.
This is a Doob martingale with $Y_0 = \mathbb{E}[X]$ and $Y_m = X$.

We will now bound the differences $|Y_j - Y_{j-1}|$, for all $j \in [m]$.
For each $i \in [\ell]\setminus\{1\}$, let $N^i(x)\coloneqq \{y \in \bigcup_{C' \in \cC'}V(C') :\dist(x,y) = i\}$.
Let $E^{i} \subseteq E$ be the collection of edges $e = (u,v)$ where both $u$ and $v$ are clones of a vertex $z \in N^i(x)$.
Note that the sets $E^2, \dots, E^\ell$ partition~$E$.
Moreover, for each $j \in [m]$, if $e_j \in E^i$, then $|Y_j - Y_{j-1}| \le n^{\ell-i}$.\COMMENT{We have seen this argument a few times by now.
If $e_j \in E^i$ with endpoint clone $z \in V(I)$, say, and if $x$ and $z$ are in the same cube, then $i$ directions for this cube are fixed, leaving $n-s-i$ free directions.
Overall, this means there are at most $n^{\ell - i}$ choices for other directions to make up cubes.}
Furthermore, $|E^i| = 2^s|N^i(x)| \le 2^s n^i$ and, thus,
\[\sum_{j \in [m]}|Y_j - Y_{j-1}|^2 \le \sum_{i=2}^{\ell}|E^i(x)|n^{2\ell-2i} \le\sum_{i=2}^{\ell}2^s n^i n^{2\ell-2i} = \bigO(n^{2\ell-2}).\]
Therefore, we can apply \cref{lem: Azuma} and combine it with \eqref{eqn:bdlots} to obtain \COMMENT{We are applying \cref{lem: Azuma} with $\mathbb{E}[X]/2$ playing the role of $t$. 
This gives that $\mathbb{P}[|X - \mathbb{E}[X]| \ge \mathbb{E}[X]/2] \le 2e^{-\gamma^2 n^{2\ell}/8O(n^{2\ell -2})}=e^{-\Omega(n^2)}$.}
\[\mathbb{P}[|B(G)| < \gamma n^\ell/4] \le \mathbb{P}[X < \gamma n^\ell /4] \le \mathbb{P}[|X - \mathbb{E}[X]| \ge \mathbb{E}[X]/2] \le e^{-10n.}\qedhere\]
\end{proof}

The following observation will also be used repeatedly.

\begin{remark}\label{prop:cubesdirect}
Let $n,\ell \in \mathbb{N}$ and $\eta,\eta'\in(0,1)$ with $1/n\ll \eta'<\eta$.
Let $x\in V(\cQ^n)$.
Let $\cC$ be a collection of $\ell$-cubes $C \subseteq \cQ^{n}$ such that $x \in V(C)$ for all $C \in \cC$ and $|\cC| \ge \eta  n^\ell$.
Let $\cD'\subseteq \cD(\cQ^{n})$ be a set of directions with $|\cD'|\leq \eta'n$.
Then, there exists a cube $C\in \cC$ with $\cD(C)\cap \cD'=\varnothing$. 
\end{remark}

\begin{proof}
Observe that the number of $\ell$-cubes $C\subseteq\cQ^n$ with $x\in C$ and $\cD(C)\cap \cD'\neq \varnothing$ is at most $\eta'n\cdot n^{\ell-1}$. 
Since $\eta >\eta'$, we are done. 
\end{proof}

As discussed in \cref{section:outline6}, a crucial requirement for the proof of \cref{thm: kedgehit} will be that vertices of very low degree in $\cQ^n_{1/2-\varepsilon}$ are few and far apart.
Moreover, we will also require some more properties about the distribution of these vertices, and that, for all of them, we can find many candidates for special absorbing structures.
We express all this information in the following definition.

\begin{definition}\label{def:rob}
Let $n,s,\ell \in \mathbb{N}$ with $1/n\ll 1/s\leq 1/\ell$, and let $\eps_1,\eps_2,\gamma\in[0,1]$. 
Fix an ordering of the layers $L_1,\ldots,L_{2^s}$ of $\cQ^n$ induced by any Hamilton cycle in $\cQ^s$ (as defined in \cref{sect8notation}).\COMMENT{This need for the layers and the parameter $s$ is because the definition of the absorbing structures in \ref{itm:bad4} require the layer structure to be in place.}
Let $G\subseteq \cQ^{n}$ be a spanning subgraph.
For any $\varepsilon>0$, let $\cU(G, \eps) \coloneqq \{x \in V(\cQ^{n}) : d_G(x)<\eps n\}$. 
Let $\cU\subseteq V(\cQ^n)$ be a set of size $|\cU|\leq2^{\eps_2n}$.
We say that $G$ is $(s,\ell, \eps_1, \eps_2, \gamma, \cU)$\emph{-robust}\index{(sl@$(s,\ell, \eps_1, \eps_2, \gamma, \cU)$-robust} if the following properties are satisfied:
\begin{enumerate}[label=$(\mathrm{R}\arabic*)$]
  \item\label{itm:bad1} $\cU(G,\eps_1)\subseteq\cU$.
    \item\label{itm:bad2} For all $x \in \cU$ and every $y \in B_{\cQ^{n}}^{s+5\ell}(x)\setminus \{x\}$, we have $d_G(y) \geq \gamma n$.
    \item\label{itm:bad3} For all $x \in V(\cQ^n)$, we have $|\cU\cap B_{\cQ^{n}}^{\gamma n}(x)| \le 1$.\COMMENT{In particular, in every molecule there is at most one \emph{bad} vertex.}
    \item\label{itm:bad4} For all $x \in \cU$ and any distinct directions $a,b\in \cD(\cQ^{n})$, there exists a collection $\mathfrak{C}(x,a,b)$ of $(x,a,b)$-consistent systems of paths in $G \cup \{\{x,x+a\},\{x,x+b\}\}$ which satisfies the following.
    Let $L$ be the layer containing $x$.
    \begin{enumerate}[label=$(\mathrm{R}4.\mathrm{\Roman*})$]
        \item\label{itm:bad41} Suppose $x+a,x+b\in V(L)$.
        Then, there exists a collection $\cD^{(2)}(x,a,b)$ of disjoint pairs of distinct directions $c,d\in\cD(L)\setminus\{a,b\}$ such that $|\cD^{(2)}(x,a,b)|\geq \gamma n$ and, for every $(c,d)\in \cD^{(2)}(x,a,b)$, there is a collection $\cD^{(4)}(x,a,b,c,d)$ of disjoint $4$-tuples of distinct directions in $\cD(L)\setminus\{a,b,c,d\}$ with $|\cD^{(4)}(x,a,b,c,d)|\geq \gamma n$ satisfying the following property:
        for each $(c,d)\in \cD^{(2)}(x,a,b)$ and each $(d_1,d_2,d_3,d_4)\in \cD^{(4)}(x,a,b,c,d)$, the $(x,a,b)$-consistent system of paths $\mathit{CS}(c,d,d_1,d_2,d_3,d_4)=(P_1,\ldots,P_6)$ defined as in \cref{sect:absstruct} belongs to $\mathfrak{C}(x,a,b)$.
        \item\label{itm:bad42} Suppose $x+a,x+b\notin V(L)$.
        Then, there exists a collection $\cD^{(2)}(x,a,b)$ of disjoint pairs of distinct directions $d_1,d_2\in\cD(L)$ such that $|\cD^{(2)}(x,a,b)|\geq \gamma n$ and, for every $(d_1,d_2)\in \cD^{(2)}(x,a,b)$, the $(x,a,b)$-consistent system of paths $\mathit{CS}(d_1,d_2)=(P_1,P_2)$ defined as in \cref{sect:absstruct} belongs to $\mathfrak{C}(x,a,b)$.
        \item\label{itm:bad43} Suppose $x+a\notin V(L)$ and $x+b \in V(L)$.
        Then, there exists a set $\cD(x,a,b)$ of directions $d_1\in\cD(L)$ such that $|\cD(x,a,b)|\geq \gamma n$ and, for every $d_1\in \cD(x,a,b)$, there exists a collection $\cD^{(2)}(x,a,b,d_1)$ of disjoint pairs of distinct directions in $\cD(L)\setminus\{b,d_1\}$ with $|\cD^{(2)}(x,a,b,d_1)|\geq \gamma n$ satisfying the following property:
        for each $d_1\in \cD^{(2)}(x,a,b)$ and each $(d_2,d_3)\in \cD^{(2)}(x,a,b,d_1)$, the $(x,a,b)$-consistent system of paths $\mathit{CS}(d_1,d_2,d_3)=(P_1,P_2,P_3)$ defined as in \cref{sect:absstruct} belongs to $\mathfrak{C}(x,a,b)$.
    \end{enumerate}
    \item\label{itm:bad5} Let $x_1 \coloneqq \{0\}^n$, $x_2 \coloneqq \{1\}^n$, $x_3 \coloneqq \{1\}^{\lceil n/2\rceil}\{0\}^{n-\lceil n/2\rceil}$ and $x_4 \coloneqq \{0\}^{\lceil n/2\rceil}\{1\}^{n-\lceil n/2\rceil}$.
    Then, for each $i \in [4]$ we have $\cU \cap B_{\cQ^{n}}^{s+\ell}(x_i) = \varnothing$.
\end{enumerate}
\end{definition}

\begin{lemma}\label{lem: badvertices}
Let $1/n \ll 1/s\leq 1/\ell \ll \eps_1 \ll \eps \ll \eps_2 \ll \gamma \ll 1/r $ with $n, s,\ell, r \in \mathbb{N}$.
Then, 
\begin{enumerate}[label=(\roman*)]
    \item\label{lem:rob1} a.a.s.~$G\sim \cQ^{n}_{1/2-\eps}$  is $(s,\ell,\eps_1,\eps_2,\gamma,\cU(G,\varepsilon_1))$-robust, and
    \item\label{lem:rob2} given any $\cU \subseteq V(\cQ^n)$ with $|\cU|\leq2^{\eps_2n}$ and any $H\subseteq \cQ^n$ which is $(s,\ell,\eps_1,\eps_2,\gamma,\cU)$-robust, there exists an edge-decomposition  $H = H_1 \cupdot \dots \cupdot H_r$ such that for each $i \in [r]$ we have that $H_i\subseteq H$ is spanning and $(s, \ell, \eps_1/(2r), \eps_2, \gamma/r^{10},\cU)$-robust.
\end{enumerate} 
\end{lemma}

\begin{proof}
We begin with a proof of \ref{lem:rob1}.
Let $\cU\coloneqq\cU(G,\varepsilon_1)$.
For any $x \in V(\cQ^n)$, we have that\COMMENT{
Choose at most any $\eps_1 n$ of its neighbours that it will have edges too. 
We then multiply by the probability that it doesn't have edges to all the rest, making it bad.
We do not condition on whether or not there are edges with the $\eps_1 n$ neighbours chosen, and therefore this covers all cases.
\begin{align*}
\mathbb{P}[x \in \cU]  &\le \binom{n}{\eps_1 n}(1/2 +  \eps)^{n-\eps_1 n} \\
    &\le (\frac{en}{\eps_1n})^{\eps_1 n} (1/2)^{n-\eps_1 n}(1+ 2\eps)^{n-\eps_1 n}\\
    &\le 2^{\eps_1 n \log_2(e/\eps_1)} 2^{\eps_1 n - n} e^{2(1-\eps_1)\eps n}\\
    &< 2^{\eps_1 n \log_2(e/\eps_1)} 2^{\eps_1 n - n} 2^{4(1-\eps_1)\eps n}\\
    &\le 2^{-n +20\eps n},
\end{align*}
where in the last inequality we use that $\eps_1 \ll \eps$.
}
\begin{equation}\label{eqn: probvbad} 
\mathbb{P}[x \in \cU] \leq \binom{n}{\eps_1 n}(1/2 + \eps)^{n-\eps_1 n} < 2^{-n + 20\eps n}.
\end{equation}
It follows that $\mathbb{E}[|\cU|] < 2^n 2^{-n + 20\eps n} = 2^{20\eps n}$.
Therefore, by Markov's inequality we have that 
\[\mathbb{P}[|\cU| \ge 2^{\eps_2n}] < 2^{20\eps n}/2^{\eps_2n} < 2^{-\eps_2 n/2},\]
so a.a.s.~$|\cU| \leq 2^{\eps_2n}$.
\ref{itm:bad1} holds trivially by the choice of $\cU$.
Furthermore, we have that
\[\Bigg|\bigcup_{i\in[4]}B^{s+\ell}_{\cQ^n}(x_i)\Bigg|\leq 5n^{s+\ell},\]
so, by \eqref{eqn: probvbad} and a union bound, \ref{itm:bad5} also holds a.a.s.

To see that \ref{itm:bad2} holds, fix $x \in V(\cQ^{n})$ and $y \in B^{s+5\ell}_{\cQ^n}(x)$.
Then, $\mathbb{E}[d_{G-\{x\}}(y)] = (1/2 -\eps)n \pm 1$.
Thus, by \cref{lem:Chernoff},\COMMENT{By \cref{lem:Chernoff}, we have that
\[\mathbb{P}[d_{G-\{x\}}(y)\leq\gamma n-1]\leq\mathbb{P}[d_{G-\{x\}}(y)\leq\mathbb{E}[d_{G-\{x\}}(y)]/2]\leq e^{-\mathbb{E}[d_{G-\{x\}}(y)]/8}\le e^{-n/20}\le 2^{-n/20}.\] 
}
\[\mathbb{P}[d_{G-\{x\}}(y)\leq\gamma n-1]\leq 2^{-n/20}.\]
Therefore, by \eqref{eqn: probvbad}\COMMENT{Observe that the events $u\in\cU$ and $d_{G-\{x\}}(y)\leq\gamma n$ are independent.}, we have that $\mathbb{P}[x\in \cU \wedge d_{G-\{x\}}(y)\leq \gamma n-1]\leq 2^{-n + 20\eps n-n/20}\leq 2^{-31n/30}$. 
A union bound over all $x \in V(\cQ^{n})$ and over all $y \in B^{s+5\ell}_{\cQ^n}(x)$ shows that \ref{itm:bad2} holds a.a.s.\COMMENT{We have a union bound of at most $n^{s+5\ell}2^n\leq2^{61n/60}$.
Furthermore, note that $d_{G-\{x\}}(y)\leq\gamma n-1\implies d_G(y)\leq\gamma n$, so this is enough to show what we want.}

The fact that \ref{itm:bad3} holds a.a.s.~can be shown similarly. \COMMENT{
Fix $x \in V(\cQ^{n})$ and $y_1, y_2 \in \cB^{\gamma n}_{\cQ^n}(x)$.
We have that $\mathbb{E}[d_{G-\{x, y_2\}}(y_1)] = (1/2 -\eps)n \pm 2$.
By the same argument as in \eqref{eqn: probvbad}, we have that 
\[ \mathbb{P}[d_{G-\{x, y_2\}}(y_1) < \eps_1 n] < 2^{-n+20\eps_1 n}.\] 
Now, taking a union bound over all $y_1,y_2 \in \cB^{\gamma n}_{\cQ^n}(x)$ gives
\[\mathbb{P}[\text{there exist distinct }y_1,y_2\in\cB^{\gamma n}_{\cQ^n}(x)\cap\cU]<(\gamma n)^2\binom{n}{\gamma n}^22^{-2n+40\eps n}<(\gamma n)^2(e/\gamma)^{2\gamma n}2^{-2n+40\eps n}<2^{-(2-c)n},\]
for some fixed $c < 1$.
We have by a union bound over all $x \in V(\cQ^{n})$ that \ref{itm:bad3} holds a.a.s.}

Finally, consider \ref{itm:bad4}. 
Let $x \in V(\cQ^{n})$, and suppose $x\in V(L)$, for some layer $L$.
First, let $a, b \in \cD(L)$ be distinct.
We are going to show that a.a.s.~we can find the desired collection of $(x,a,b)$-consistent systems of paths.

Recall that an $(x,a,b)$-consistent system of paths $(P_1,\ldots,P_6)$, as defined in \cref{sect:absstruct}, is determined uniquely by a $6$-tuple of directions $(c,d,d_1,d_2,d_3,d_4)$.
In order to show that \ref{itm:bad41} is satisfied, we will first consider the directions $c$ and $d$, and then the rest of the tuple.
Recall that all $(x,a,b)$-consistent systems of paths contain the two edges $\{x,x+a\}$ and $\{x,x+b\}$.
Then, once $c$ and $d$ are fixed, this determines a total of $6$ more edges.
The remaining $8$ edges will be determined by the choice of $(d_1,d_2,d_3,d_4)$.

Consider a collection $\mathcal{W}$ of disjoint pairs of distinct directions $(c,d)$ with $c,d\in\mathcal{D}(L)\setminus\{a,b\}$ such that $|\cW|\geq n/4$\COMMENT{Such a collection exists trivially.}.
For each $(c,d)\in\cW$, let $E^*(c,d)\subseteq E(\cQ^n)$ be the set of six edges of an $(x,a,b)$-consistent system of paths determined by these two directions.
Observe that, since the pairs in $\cW$ are disjoint, it follows that, for any distinct $(c,d),(c',d')\in\cW$, we have $E^*(c,d)\cap E^*(c',d')=\varnothing$.
Now let $\cW_G\coloneqq\{(c,d)\in\cW:E^*(c,d)\subseteq E(G)\}$ and $X\coloneqq|\cW_G|$.
We have that $\mathbb{E}[X]\geq(1/2-\varepsilon)^6n/4$ and, by \cref{lem:Chernoff}, it follows that $\mathbb{P}[X\leq\gamma n]\leq2^{-\gamma n}$\COMMENT{We have that
\[\mathbb{P}[X\leq\gamma n]\leq\mathbb{P}[X\leq\mathbb{E}[X]/2]\leq e^{-\mathbb{E}[X]/8}\leq e^{-(1/2-\varepsilon)^6n/36}\leq e^{-\gamma n}\leq 2^{-\gamma n}.\]}.

For each $(c,d)\in\cW$, let $\cV(c,d)$ be a collection of disjoint $4$-tuples of distinct directions $(d_1,d_2,d_3,d_4)$ with $d_1,d_2,d_3,d_4\in\mathcal{D}(L)\setminus\{a,b,c,d\}$ such that $|\cV(c,d)|\geq n/5$\COMMENT{Again, such a collection exists trivially.}.
For each $(d_1,d_2,d_3,d_4)\in\cV(c,d)$, let $E^*(c,d,d_1,d_2,d_3,d_4)\subseteq E(\cQ^n)$ be the set of eight edges of an $(x,a,b)$-consistent system of paths determined by $(c,d,d_1,d_2,d_3,d_4)$ which are not in $E^*(c,d)\cup\{\{x,x+a\},\{x,x+b\}\}$.
In particular, since the tuples in $\cV(c,d)$ are disjoint, it follows that, for any distinct $(d_1,d_2,d_3,d_4),(d_1',d_2',d_3',d_4')\in\cV(c,d)$, we have $E^*(c,d,d_1,d_2,d_3,d_4)\cap E^*(c,d,d_1',d_2',d_3',d_4')=\varnothing$.
Now let $\cV_G(c,d)\coloneqq\{(d_1,d_2,d_3,d_4)\in\cV(c,d):E^*(c,d,d_1,d_2,d_3,d_4)\subseteq E(G)\}$, and let $Y(c,d)\coloneqq|\cV_G(c,d)|$.
We then have that $\mathbb{E}[Y(c,d)]\geq(1/2-\varepsilon)^8n/5$ and, again by \cref{lem:Chernoff}, it follows that $\mathbb{P}[Y(c,d)\leq\gamma n]\leq2^{-\gamma n}$\COMMENT{We have that
\[\mathbb{P}[Y(c,d)\leq\gamma n]\leq\mathbb{P}[Y(c,d)\leq\mathbb{E}[Y(c,d)]/2]\leq e^{-\mathbb{E}[Y(c,d)]/8}\leq e^{-(1/2-\varepsilon)^8n/40}\leq e^{-\gamma n}\leq 2^{-\gamma n}.\]}.
Thus, by a union bound, with probability at least $1-2^{-\gamma n/2}$, for every $(c,d)\in\cW$ we have $Y(c,d)\geq\gamma n$.

Let $\mathcal{E}(x,a,b)$ be the event that $G\cup\{\{x,x+a\},\{x,x+b\}\}$ contains a collection $\mathfrak{C}(x,a,b)$ of $(x,a,b)$-consistent systems of paths satisfying \ref{itm:bad41}.
By combining all the above, it follows that\COMMENT{We have that 
\[\mathbb{P}[\mathcal{E}(x,a,b)]\geq\mathbb{P}[\{X\geq\gamma n\}\wedge\{Y(c,d)\geq\gamma n\text{ for all }(c,d)\in\cW\}]\geq1-2^{-\gamma n/4}.\]}
\begin{equation}\label{equa:newrobustrandom}
    \mathbb{P}[\mathcal{E}(x,a,b)]\geq1-2^{-\gamma n/4}.
\end{equation}

The same bound can be proved for the cases where $a \notin \cD(L), b \in \cD(L)$ and $a, b \notin \cD(L)$.
Observe that, for any $x\in V(\cQ^n)$ and $a,b\in\cD(\cQ^n)$, the event $\mathcal{E}(x,a,b)$ is independent of the event that $x\in\cU$\COMMENT{This is because the consistent system of paths only has two edges incident with $x$, and these two edges are added deterministicly.}.
Now, by combining \eqref{equa:newrobustrandom} with \eqref{eqn: probvbad} and a union bound over all choices of $a,b\in\mathcal{D}(\cQ^n)$, and then a union bound over all $x \in V(\cQ^{n})$, we conclude that \ref{itm:bad4} holds a.a.s.\COMMENT{We have that the probability that \ref{itm:bad4} does not hold is at most $2^n(n^22^{-\gamma n/4}2^{-(1-20\varepsilon)n})=o(1)$, since $\varepsilon\ll\gamma$.}
\COMMENT{All five properties a.a.s.~hold simultaneously by one more union bound.}

The proof of \ref{lem:rob2} is similar.
Let $H\subseteq \cQ^n$ be given and consider a random partition of $E(H)$ into $r$ parts $H_1, \dots, H_r$, in such a way that each edge is assigned to one of the parts uniformly and independently of all other edges.
For each $x\in V(\cQ^n)\setminus\cU$, let $\mathcal{B}(x)$ be the event that there exists some $i\in[r]$ such that $d_{H_i}(x) < d_{H}(x)/(2r)$.
Observe that, if $\overline{\mathcal{B}(x)}$ holds for all $x\in V(\cQ^n)\setminus\cU$, then no vertex outside $\cU$ will be contained in $\cU(H_i, \eps_1/(2r))$ for any $i \in [r]$.
It would then follow that \ref{itm:bad1}, \ref{itm:bad2}, \ref{itm:bad3} and \ref{itm:bad5} all hold with the desired constants for each~$H_i$.
Fix $x \in V(\cQ^n)\setminus \cU$ and $i \in [r]$. 
Let $X\coloneqq d_{H_i}(x)$.
Then, $\mathbb{E}[X] = d_{H}(x)/r\geq\varepsilon_1n/r$.
Thus, by \cref{lem:Chernoff}, $\mathbb{P}[X \leq \mathbb{E}[X]/2] \leq e^{-\eps_1^2 n}$.\COMMENT{$\mathbb{P}[X \leq \mathbb{E}[X]/2] \leq e^{-\mathbb{E}[X]/8}\leq e^{-\varepsilon_1 n/8r}\leq e^{-\eps_1^2 n}$.}
A union bound over all $i \in [r]$ shows that $\mathbb{P}[\mathcal{B}(x)]\leq re^{-\eps_1^2 n}\leq e^{-\varepsilon_1^3n}$.

We now consider the property \ref{itm:bad4}.
For each $x\in\cU$, let $\mathcal{B}(x)$ be the event that there exist $i\in[r]$ and distinct directions $a,b \in \cD(\cQ^n)$ such that \ref{itm:bad4} does not hold for $H_i$ with $\gamma/r^{10}$ playing the role of $\gamma$.

Fix $x \in \cU$, $i \in [r]$ and distinct directions $a,b \in \cD(\cQ^n)$.
Similarly as in \ref{lem:rob1}, using \cref{lem:Chernoff} one can show that the probability that $H_i$ does not satisfy \ref{itm:bad4} for $x$ with $\gamma /r^{10}$ playing the role of $\gamma$ is at most $2^{-\gamma^2 n} + 2^{-\gamma^3 n}$.
Therefore, by a union bound over all choices of $a, b \in \cD(\cQ^n)$ and over each $i \in [r]$, we have that $\mathbb{P}[\mathcal{B}(x)]\le r n^2(2^{-\gamma^2 n} + 2^{-\gamma^3 n})\leq e^{-\varepsilon_1^3n}$.

Finally, we are interested in the event where $\cB(x)$ does not occur for any $x \in V(\cQ^n)$.
We will invoke \cref{lem: LLL}.
Note that each event $\cB(x)$ is mutually independent of all but at most $n^{10}$ other events.\COMMENT{A crude bound for consistent systems.}
We have that $\mathbb{P}[\cB(x)] \leq e^{-\eps_1^3 n}$ for every $x \in V(\cQ^n)$ and $e\cdot e^{-\eps_1^3 n}(n^{10} +1) < 1$, so by \cref{lem: LLL} there exists an edge-decomposition of $H$ with the desired properties.
\end{proof}

Finally, we need to show a result analogous to \cref{lem: scant don't clump} for robust graphs, that is, that scant molecules are not too clustered.
Recall that $\mathit{Res}(\cQ^n, \delta)$ was defined in \cref{section:tree1}.

\begin{lemma}\label{lem: scant don't clumphit}
Let $0<1/n\ll1/C\ll\varepsilon_1,\varepsilon_2\ll\gamma,\delta\leq1$ and $1/n\ll1/s\leq1/\ell$, where $n,C,s,\ell\in\mathbb{N}$.
Let $H \subseteq \cQ^n$ and $\cU\subseteq V(\cQ^n)$ be such that $H$ is $(s,\ell,\varepsilon_1,\varepsilon_2,\gamma,\cU)$-robust.
Let $\cU_I\subseteq V(I)$ be the set of vertices $u\in V(I)$ such that $\cU$ contains some clone of $u$.
For each $v \in V(I)$ and each $i\in[2^s]$, let $v_i$ be the $i$-th clone of $v$, and let $\mathcal{M}_v\coloneqq\{v_i: i\in[2^s]\}$. 
Let $R \sim \mathit{Res}(I, \delta)$ and, for each $i \in [2^s]$, let $R_i$ be the $i$-th clone of $R$.
Let 
\[B \coloneqq \{v \in  V(I): \text{there exists } i \in [2^s]\text{ with }v_i\notin\cU\text{ and } e_{H}(v_i, R_i) < \varepsilon_1\delta n/4\}.\]
Let $\mathcal{E}_1$ be the event that there exists some $u \in V(I)$ such that $|B_{I}^{10\ell}(u)\cap B|\geq C$.
Let $\mathcal{E}_2$ be the event that there exists some $u \in \cU_I$ such that $|B_{I}^{5\ell}(u)\cap B|\geq 1$.
Then, $\mathbb{P}[\mathcal{E}_1\vee\mathcal{E}_2] \le 1/n$.
\end{lemma}

\begin{proof}
Let $u \in V(I)$ and let $D\subseteq B_{I}^{10\ell}(u)\setminus\cU_I$ be a set of $C$ vertices. 
Let $D'\coloneqq\bigcup_{x,y\in D:x\neq y}N_{I}(x)\cap N_{I}(y)$. 
Since any pair of distinct vertices in $I$ share at most two neighbours, we have that $|D'| \leq 2\binom{C}{2}$.
For each $i \in [2^s]$, we denote the $i$-th clone of $D'$ by $D'_i$, and let $R'_{i}\coloneqq R_{i}\setminus D'_{i}$.

For each $x\in V(\cQ^{n})$, let $i(x)$ be the unique index $i\in[2^s]$ such that $x\in V(L_i)$.
Observe that, by \ref{itm:bad1}, we have $e_H(x, V(L_{i(x)})) >2\varepsilon_1 n/3$\COMMENT{Since $s$ is constant.} for every $x\in V(\cQ^{n})\setminus\cU$.
For each $x\in V(\cQ^{n})$, let $\mathcal{E}_x$ be the event that $e_H(x, R_{i(x)})\leq\varepsilon_1\delta n/4$, and let $\mathcal{E}'_x$ be the event that $e_H(x, R'_{i(x)})\leq\varepsilon_1\delta n/4$\COMMENT{Note that if $\mathcal{E}_x$ occurs, then $\mathcal{E}_x'$ occurs, that is, $\mathbb{P}[\mathcal{E}_x]\leq\mathbb{P}[\mathcal{E}'_x]$.}.
It follows by \cref{lem:Chernoff} that $\mathbb{P}[\mathcal{E}'_x]\leq e^{-\varepsilon_1\delta n/16}$\COMMENT{Let $Y\coloneqq e_H(x, R'_{i(x)})$.
We have that $\mathbb{E}[Y]\geq\delta(\varepsilon_1 n-s-2\binom{C}{2})\geq\varepsilon_1\delta n/2$.
It follows by \cref{lem:Chernoff} that
\[\mathbb{P}[Y\leq\varepsilon_1\delta n/4]\leq\mathbb{P}[Y\leq\mathbb{E}[Y]/2]\leq e^{-\mathbb{E}[Y]/8}\leq e^{-\varepsilon_1\delta n/16}.\]
}
for all $x\in V(\cQ^n)\setminus\cU$.
For each $v\in V(I)$, let $\mathcal{E}_v$ and $\mathcal{E}_v'$ be the events that there exists $i\in[2^s]$ with $v_i\notin\cU$ such that $\mathcal{E}_{v_i}$ and $\mathcal{E}_{v_i}'$ hold, respectively.
By a union bound, it follows that $\mathbb{P}[\mathcal{E}_v']\leq 2^se^{-\varepsilon_1\delta n/16}$ for all $v\in V(I)$.
Finally, let $\mathcal{E}_D$ and $\mathcal{E}'_D$ be the events that $\mathcal{E}_v$ and $\mathcal{E}_v'$, respectively, hold for every $v\in D$.
Note that the events in the collection $\{\mathcal{E}_v': v\in V(I)\}$ are mutually independent.
Furthermore, since the event $\mathcal{E}_x$ implies $\mathcal{E}_x'$ for all $x\in V(\cQ^n)$, we have that
\[\mathbb{P}[\mathcal{E}_D] \leq \mathbb{P}[\mathcal{E}'_D] \leq (2^se^{-\varepsilon_1\delta n/16})^{C} < e^{-5n}.\]
By a union bound over all $u\in V(I)$ and over all choices of $D$, we have $\mathbb{P}[\mathcal{E}_1]\leq e^{-n}$\COMMENT{$2^{n}n^{10C\ell}e^{-5n}<e^{-n}$.}.

Consider now any $u\in \cU_I$.
Observe that, if $v \in B_{I}^{5\ell}(u)$, then for every $i,j\in[2^s]$ we have that $\dist(u_i,v_j) \leq 5\ell +s$.
Therefore, by \ref{itm:bad2}, for all $v \in B_{I}^{5\ell}(u)$ and $i\in[2^s]$ such that $v_i\notin\cU$, we have $d_H(v_i)\geq \gamma n$.
For each $v\in B_I^{5\ell}(u)$ and each $i\in[2^s]$ with $v_i\notin\cU$, let $\mathcal{F}_{v_i}$ be the event that $e_H(v_i,R_i)\leq\varepsilon_1\delta n/4$, and let $\mathcal{F}_{v}$ be the event that there exists some $i\in[2^s]$ with $v_i\notin\cU$ such that $\mathcal{F}_{v_i}$ holds.
By \cref{lem:Chernoff} and a union bound, it follows that $\mathbb{P}[\mathcal{F}_{v}]\leq2^{-\gamma\delta n/16}$.\COMMENT{Let $Y\coloneqq e_H(v_j,R_j)$.
We have that $\mathbb{E}[Y]\geq\delta(\gamma n-s)\geq\gamma\delta n/2\ (\gg\varepsilon_1\delta n/2)$.
It follows by \cref{lem:Chernoff} that
\[\mathbb{P}[\mathcal{F}_{v_j}]=\mathbb{P}[Y\leq\varepsilon_1\delta n/4]\leq\mathbb{P}[Y\leq\gamma\delta n/4]\leq\mathbb{P}[Y\leq\mathbb{E}[Y]/2]\leq e^{-\mathbb{E}[Y]/8}\leq e^{-\gamma\delta n/16}.\]
Then, by a union bound over all $j\in[2^s]$, $\mathbb{P}[\mathcal{F}_{v}]\leq2^se^{-\gamma\delta n/16}\leq2^{-\gamma\delta n/16}$.}
Then, by a union bound over all $u\in \cU_I$ and $v \in B_{I}^{5\ell}(u)$,
\[
\mathbb{P}[\mathcal{E}_2] \leq |\cU_I|\cdot|B^{5\ell}_I(\cU_I)| \cdot \mathbb{P}[\mathcal{F}_{v}]\leq 2n^{5\ell}2^{\varepsilon_2 n}2^{-\gamma\delta n/16}.\qedhere
\]
\end{proof}


\subsection{Hamilton cycles in robust subgraphs of the cube}\label{section:hitmain}

It will be useful to prove the following result, which (together with \cref{lem: badvertices}) directly implies the case $k=2$ of \cref{thm: kedgehit} (by choosing $H$ in \cref{thm: kedgehit} to play the role of $H'$ in \cref{thm: important} and $F$ to be empty). 
As with \cref{thm:main1}, the formulation of \cref{thm: important} is designed so that the case $k>2$ can be derived easily (see \cref{sect:hit}).
To state the result, we need the following notation.

Given any integers $s\leq n$, we say that $d\in\cD(\cQ^n)$ is an \emph{$s$-direction} if its only non-zero coordinate is one of the first $s$ coordinates.
Given a graph $F \subseteq \cQ^n$, a set $\cU\subseteq V(\cQ^n)$ and $\ell,s \in \mathbb{N}$, we say that $F$ is \emph{$(\cU,\ell,s)$-good}\index{(ul@$(\cU,\ell,s)$-good} if, for each $x\in\cU$, the set $E_F(x)\coloneqq\{e\in E(F):e\cap N_{\cQ^n}(x)\neq\varnothing\}$ satisfies that, for each $d\in\mathcal{D}(\cQ^n)$ which is not an $s$-direction, we have $|\{e\in E_F(x):\mathcal{D}(e)=d\}|\leq n/\ell$.
Thus, a graph is good if locally the directions of its edges are not too correlated (ignoring $s$-directions).
The goodness of the `forbidden' graph $F$ below will be needed when finding the special absorbing structures (see Step~11).

\begin{theorem}\label{thm: important}
Let $0< 1/\ell \ll \eps_1\ll \eps_2 \ll \gamma \leq 1$ and $1/\ell\ll \eta,1/c \leq 1$, with $\ell \in \mathbb{N}$. 
Let $s\coloneqq10\ell$ and $n\in\mathbb{N}$.
Then, there exists $\Phi \in \mathbb{N}$ such that the following holds.

Let $H\subseteq \cQ^{n}$ and $\cU\subseteq V(\cQ^n)$ be such that $H$ is an $(s,\ell,\eps_1,\eps_2,\gamma,\cU)$-robust subgraph, and let $Q\sim\cQ^n_\eta$.
Then, a.a.s.~there is a $(\cU,\ell^2,s)$-good subgraph $Q'\subseteq Q$ with $\Delta(Q')\leq\Phi$ such that
\begin{itemize}
    \item for every $H'\subseteq \cQ^{n}$, where $d_{H'}(x)\geq 2$ for every $x \in \cU$, and
    \item for every $F\subseteq\cQ^n$ with $\Delta(F)\leq c\Phi$ which is $(\cU, \ell,s)$-good,
\end{itemize}
we have that $((H\cup Q)\setminus F)\cup H'\cup Q'$ contains a $(\cU, \ell^2, s)$-good Hamilton cycle $C$ such that, for all $x\in\cU$, both edges of $C$ incident to $x$ belong to $H'$.
\end{theorem}

As we have already discussed, the proof of \cref{thm: important} builds on that of \cref{thm:main1}.
Thus, we will avoid repeating all the details which are analogous, and we will often refer back the proof of \cref{thm:main1}.
A full sequential proof can be found in~\cite{Athesis}.

\begin{proof}[Proof of \cref{thm: important}]
Let $1/D, \delta' \ll 1$, let $\varepsilon_1 \ll \varepsilon_2\ll \gamma, \delta'\leq 1$, and let 
\begin{equation}\label{equa:hierarchy2}
    0 < 1/n_0 \ll \delta, \lambda \ll  1/\ell \ll 1/k^*, \alpha' \ll \beta, 1/S' \ll 1/c, 1/D, \eta, \varepsilon_1, \delta',
\end{equation}  
where $n_0,\ell,k^*,S',D\in\mathbb{N}$. 
Our proof assumes that $n$ tends to infinity; in particular, $n\geq n_0$.
Let \[\Phi\coloneqq60\ell^4\] and \[\Psi\coloneqq c\Phi.\]

We define the layers $L_1, \dots, L_{2^s}$ of $\cQ^n$, the intersection graph $I$ and, for each $G \subseteq \cQ^n$, the graphs $I(G)$ and $G_I$ as in the proof of \cref{thm:main1}. 
Similarly, for any layer $L$ and $\cG\subseteq I$, we keep the notation $\cG_L$ for the clone of $\cG$ in $L$.

Let $\cU$ be as in the  statement of \cref{thm: important}.
In particular, $\cU(H,\eps_1)\subseteq\cU$ by \ref{itm:bad1}.
Let $\cU_I\subseteq V(I)$ be the set of vertices $x\in V(I)$ such that there is some clone $x'$ of $x$ with $x'\in\cU$.
Note that, by property \ref{itm:bad2}, for each $x\in\cU_I$, there is exactly one clone $x'$ of $x$ with $x'\in\cU$.

For each $i \in [8]$, let $\eta_i \coloneqq \eta/8$ and $G_i \sim \cQ^n_{\eta_i}$, where these graphs are chosen independently.
We have that $\bigcup_{i=1}^{8} G_i \sim \cQ^n_{\eta'}$ for some $\eta' < \eta$, so it suffices to show that a.a.s.~there is a $(\cU,\ell^2,s)$-good subgraph $Q'\subseteq\bigcup_{i=1}^8G_i$ with $\Delta(Q')\leq\Phi$ and such that, for every $H'\subseteq \cQ^{n}$, where $d_{H'}(x)\geq 2$ for every $x\in\cU$, and every $F\subseteq\cQ^n$ with $\Delta(F)\leq\Psi$ which is $(\cU, \ell, s)$-good, the graph $((H\cup\bigcup_{i=1}^{8}G_i)\setminus F)\cup H'\cup Q'$ contains a Hamilton cycle of the form in the statement of the theorem.
We now split our proof into several steps. \\
  
\textbf{Step~1. Finding a tree and a reservoir.}
As in the proof of \cref{thm:main1}, consider the probability space $\Omega \coloneqq\cQ^{n-s}_{\eta_1^{2^s}} \times \mathit{Res}(\cQ^{n-s},\delta')$ and let $R\sim\mathit{Res}(I,\delta')$.
Let $\mathcal{E}_1$ be the event that there exists a tree $T \subseteq I(G_1)-(R\cup B^{5}_{I}(\cU_I))$ such that the following hold:
\begin{enumerate}[label=$(\mathrm{TR}\arabic*)$]
\item $\Delta(T) < D$, and
\item for all $x \in V(I)\setminus B^{5}_I(\cU_I)$, we have that $|N_{I}(x)\cap V(T)|\geq4(n-s)/5$.
\end{enumerate}
Note that, by \ref{itm:bad3}, for all $x,y\in \cU_I$ we have that $\dist(x,y)\geq\gamma n/2$.
Furthermore, by \ref{itm:bad5}, we have that, if we see $x_1 \coloneqq \{0\}^{n-s}$, $x_2 \coloneqq \{1\}^{n-s}$, $x_3 \coloneqq \{1\}^{\lceil (n-s)/2\rceil}\{0\}^{n-s-\lceil (n-s)/2\rceil}$ and $x_4 \coloneqq \{0\}^{\lceil (n-s)/2\rceil}\{1\}^{n-s-\lceil (n-s)/2\rceil}$ as vertices of $I$, then $\cU_I \cap B_{I}^{\ell}(x_i) = \varnothing$ for all $i \in [4]$.
Thus, it follows from \cref{lem: main treereshit}, with $n-s$, $D$, $\delta'$, $1/5$, $\eta_1^{2^{s}}$, $\gamma/2$, $5$ and $\cU_I$ playing the roles of $n$, $D$, $\delta$, $\eps'$, $\eps$, $\gamma$, $k$ and $\cA$, respectively, that $\mathbb{P}_\Omega[\mathcal{E}_1]=1-o(1)$.\\

\textbf{Step~2. Identifying scant molecules.}
For each $v \in V(I)$, let $\cM_v$ denote the vertex molecule of all clones of $v$ in $\cQ^n$.
We say $\cM_v$ is \emph{scant} if there exist some layer $L$ and some vertex $x \in V(\cM_v \cap L)\setminus\cU$ such that $e_H(x, R_L) < \varepsilon_1\delta' n/10$, where $R_L$ is the clone of $R$ in $L$.\COMMENT{Note here that scant is defined with respect to the deterministic graph $H$.}
Let $\mathcal{E}_2$ be the event that there exists some $x\in V(I)$ such that there are at least $S'$ vertices $v \in B_{I}^{10\ell}(x)$ with the property that $\mathcal{M}_v$ is scant, where $S'$ satisfies \eqref{equa:hierarchy2}.
Let $\mathcal{E}_3$ be the event that there exist $x\in \cU_I$ and $v \in B_{I}^{5\ell}(x)$ such that $\mathcal{M}_v$ is scant.\COMMENT{That is, we are about to show that vertices in $\cU$ do not have scant molecules nearby, nor are scant themselves.}
It follows from \cref{lem: scant don't clumphit} with $S'$ and $\delta'$ playing the roles of $C$ and $\delta$ that $\mathbb{P}_\Omega[\mathcal{E}_2\vee\mathcal{E}_3]=o(1)$.
Let $\mathcal{E}_1^*\coloneqq\mathcal{E}_1\wedge\overline{\mathcal{E}_2}\wedge\overline{\mathcal{E}_3}$.
Then, $\mathbb{P}_{\Omega}[\mathcal{E}_1^*]=1-o(1)$.

Condition on $\mathcal{E}_1^*$ holding.
Then, there exist a set $R\subseteq V(I)$\index{R1@$R$} and a tree $T\subseteq I(G_1)-(R\cup B_{I}^{5}(\cU_I))$\index{T1@$T$} such that the following hold:
\begin{enumerate}[label=$(\mathrm{T}\arabic*)$]
\item\label{itm:tree1hit} $\Delta(T) < D$;
\item\label{itm:tree2hit} for all $x \in V(I)\setminus B_{I}^{5}(\cU_I)$, we have that $|N_{I}(x) \cap V(T)| \ge 4(n-s)/5$;
\item\label{itm:scanthit} for every $x\in V(I)$, we have $|\{v \in B_{I}^{10\ell}(x): \cM_{v} \text{ is scant}\}| \le S'$, and
\item \label{itm: bad and scant} for every $x\in \cU_I$ and every $v\in B_{I}^{5\ell}(x)$, we have that $\mathcal{M}_v$ is not scant.
\end{enumerate}
Recall this implies clones of $T$ and $R$ satisfying \ref{itm:tree1hit}--\ref{itm: bad and scant} exist simultaneously in each layer of $G_1$.\\

\textbf{Step~3: Finding clustered robust matchings for each molecule.}
As in the proof of \cref{thm:main1}, the aim is to find auxiliary matchings which can later be used to pair up vertices which need to be absorbed.
In the proof of \cref{thm:main1}, we were able to carry out this pairing within each slice.
However, we cannot guarantee that the vertices of each special absorbing structure of Type II and III will lie within a single slice.
There will be an even number of these vertices (zero, two or four) within each vertex molecule, but there might be exactly one within a slice of this molecule, making it impossible to pair up vertices within a slice.
Thus we now consider the entire vertex molecule when finding the auxiliary matching (rather than each slice separately).
This would normally make it much more difficult to link up vertices of the skeleton in Step~18.
We are able to overcome this problem by considering matchings which are `clustered', i.e.~the endpoints of each matching edge either lie in a common slice or in two consecutive slices.

Let $q\coloneqq2^{10Dk^*}$ and $t\coloneqq2^s/q$, where $k^*$ satisfies \eqref{equa:hierarchy2}.
For each $j\in[t]$, let $S_j\coloneqq\bigcup_{i=(j-1)q+1}^{jq}L_i$.
Given any molecule $\mathcal{M}$, we define the slices $\mathcal{S}_j(\mathcal{M})\coloneqq S_j\cap\mathcal{M}$.
We denote by $\cS(\cM)$ the collection of all these slices of $\cM$.

Let $V_{\mathrm{sc}} \subseteq V(I)$\index{Vsc@$V_{\mathrm{sc}}$} be the set of all vertices $x\in V(I)$ such that $\cM_x$ is scant.
In particular, by \ref{itm: bad and scant} we have that $V_{\mathrm{sc}}\cap\cU_I=\varnothing$.
Recall $G_2\sim \cQ^{n}_{\eta_2}$.
For each $v\in V(I)\setminus (V_{\mathrm{sc}}\cup\cU_I)$, we define auxiliary bipartite graphs $H(v)\coloneqq (V(\cM_v),N_{I}(v),E_{H})$ and $G_2(v)\coloneqq (V(\cM_v),N_{I}(v),E_{G_2})$, where $E_{H}$ and $E_{G_2}$ are defined as in Step~3 of the proof of \cref{thm:main1} for vertices $v\in V(I)\setminus V_{\mathrm{sc}}$ (but now the first vertex class of $H(v)$ and $G_2(v)$ is $V(\cM_v)$ rather than $V(\cS)$ for some slice $\mathcal{S}$).
For each $v\in V_{\mathrm{sc}}$, we similarly define $H(v)$ and $G_2(v)$ as we did in Step~3 of the proof of \cref{thm:main1}, with the same modifications as above.

For each $v\in\cU_I$, we also define two such auxiliary graphs.
Let $H(v)\coloneqq (V(\cM_v),N_{I}(v),E_{H}^*)$, where $E_{H}^*$ is defined as follows.
Consider $v'\in V(\mathcal{M}_v)$ and let $L^{v'}$ be the layer which contains $v'$.
Let $w\in N_{I}(v)$, and let $w_{L^{v'}}$ be the clone of $w$ in $L^{v'}$.
Then, if $v'\in\cU$, we add $\{v',w\}$ to $E_{H}^*$ (these can be seen as purely auxiliary edges, and we will ignore their effect later).
Otherwise, $\{v',w\}\in E_{H}^*$ if and only if $w\in R$ and $\{v',w_{L^{v'}}\}\in E(H)$.
In particular, $d_{H(v)}(v') \ge \varepsilon_1\delta' n/10$ for all $v'\in V(\cM_v)$ since $\cM_v$ is a not a scant molecule.
We define $G_2(v)\coloneqq (V(\cM_v),N_{I}(v),E_{G_2})$, where $\{v',w\}\in E_{G_2}$ if and only if $\{v',w_{L^{v'}}\}\in E(G_2)$.

For every $v\in V(I)$ and every slice $\mathcal{S}\in\cS(\cM_v)$, note that the partition of $V(\mathcal{S})$ into vertices of even and odd parity is a balanced bipartition.
Define the graph $\Gamma^{\beta}_{H(v),G_2(v)}(V(\cM_v))$ as in \cref{sect5matchings}, where $\beta$ satisfies \eqref{equa:hierarchy2}.
Thus, $V(\Gamma^{\beta}_{H(v),G_2(v)}(V(\cM_v)))=V(\cM_v)$.
Furthermore, by definition, \COMMENT{It would be ok also in the following to have included edges when the neighbourhoods in $H$ are both large.
It would not be ok to include edges when the neighbourhoods in $G_2$ are both large however. 
This is because in definition of $G_2$ we included edges which do not go into the reservoir, which is the set of interest here.}
\begin{enumerate}[label=$(\mathrm{RM})$]
    \item\label{item:RMhit} given any $w_1,w_2\in V(\mathcal{M}_v)$, we have that $\{w_1,w_2\}\in E(\Gamma^{\beta}_{H(v),G_2(v)}(V(\cM_v)))$ if and only if $|N_{H(v)}(w_1)\cap N_{G_2(v)}(w_2)|\geq \beta (n-s)$ or
    $|N_{G_2(v)}(w_1)\cap N_{H(v)}(w_2)|\geq \beta (n-s)$\COMMENT{Here we are using that $|B| = n-s$ in \cref{lema:robustmatch}.}.
\end{enumerate}
For each $i \in [t]$, let $\mathfrak{A}_i(v)$ consist of all vertices of $V(\cS_i(\cM_v))$ of even parity, and let $\mathfrak{B}_i(v)$ consist of those of odd parity. 
By applying \cref{cor: robustmatch} with $d=100D$, $\alpha=\varepsilon_1\delta'/10$\COMMENT{This gives us the condition in the statement because $\mathcal{M}_v$ is not scant.}, $\eps=\eta_2$, $n= n-s$, $k=q=2^{10Dk^*}$, $\beta = \beta$, $t=t$, $G= H(\cM_v)$ and $V(\cS_1(\cM_v))\cup\ldots \cup V(\cS_t(\cM_v))$ as a partition of $V(\cM_v)$, we obtain that, with probability at least $1-2^{-9(n-s)}\geq 1-2^{-8n}$, the graph $\Gamma^{\beta}_{H(v),G_2(v)}(V(\cM_v))$ is $100D$-robust-parity-matchable with respect to $(\bigcup_{i=1}^t\mathfrak{A}_i(v), \bigcup_{i=1}^t\mathfrak{B}_i(v))$ clustered in  $(V(\cS_1(\cM_v)),\ldots, V(\cS_t(\cM_v)))$.\COMMENT{Note that here we are using the fact that $|V(\cS)|=q>k^*$, where $k^*$ is given in the starting hierarchy.}

By a union bound over all $v\in V(I)$, a.a.s.~$\Gamma^{\beta}_{H(v),G_2(v)}(V(\cM_v))$ is $100D$-robust-parity-matchable with respect to $(\bigcup_{i=1}^t\mathfrak{A}_i(v), \bigcup_{i=1}^t\mathfrak{B}_i(v))$ clustered in $(V(\cS_1(\cM_v)),\ldots, V(\cS_t(\cM_v)))$ for every $v \in V(I)$\COMMENT{This happens with probability at least $1-2^n 2^{-8n}$.}.
We condition on this event holding and call it $\cE_2^*$.
Thus, for each $v\in V(I)$ and each set $X\subseteq V(\mathcal{M}_v)$ with $|X|\leq 100D$ which contains as many odd vertices as even vertices, there exists a perfect matching $\mathfrak{M}(\mathcal{M}_v,X)$ in the bipartite graph with parts consisting of the even and odd vertices of $V(\mathcal{M}_v)\setminus X$, respectively, and edges given by $\Gamma^{\beta}_{H(v),G_2(v)}(V(\cM_v))$, with the property that, for each $e=\{w_\mathrm{e},w_\mathrm{o}\}\in\mathfrak{M}(\mathcal{M}_v,X)$, if $w_\mathrm{e}\in V(\cS_i(\cM_v))$ for some $i\in[t]$, then $w_\mathrm{o}\in V(\cS_{i-1}(\cM_v))\cup V(\cS_i(\cM_v))\cup V(\cS_{i+1}(\cM_v))$ (where indices are taken cyclically).
When we apply this in Step~15, we will have $\cU\cap V(\cM_v)\subseteq X$.

For each $v\in V(I)\setminus\cU_I$, we denote by $\mathfrak{M}(v)$\index{MS@$\mathfrak{M}(\cS)$, $\mathfrak{M}(v)$} the set of edges contained in the union (over all $X$) of the matchings $\mathfrak{M}(\mathcal{M}_v,X)$ (without multiplicity).
For each $v\in \cU_I$, we let $u(v)$ be the unique vertex in $\cM_v$ such that $u(v)\in\cU$, and we denote by $\mathfrak{M}(v)$ the set of edges contained in the union of the matchings $\mathfrak{M}(\mathcal{M}_v,X)$ over all $X$ such that $u(v)\in X$ (again, without multiplicity).\COMMENT{We can do this because we already know that we will remove the vertices in $\cU$. So all that we need is to know which vertices it will be paired with, and these will form the set $X$ which is removed.}
Furthermore, for each $v\in V(I)$ and each $e=\{w_\mathrm{e},w_\mathrm{o}\} \in \mathfrak{M}(v)$, we let $N(e)\coloneqq (N_{H(v)}(w_\mathrm{e})\cap N_{G_2(v)}(w_\mathrm{o})) \cup (N_{G_2(v)}(w_\mathrm{e})\cap N_{H(v)}(w_\mathrm{o})) $\COMMENT{Note that $N(e)\subseteq V(I)$, by definition.}.
By \ref{item:RMhit}, we have $|N(e)|\geq\beta(n-s)\geq\beta n/2$.
Let $K \coloneqq \max_{v \in V(I)}|\mathfrak{M}(v)|$.\COMMENT{In other words, $K$ is the maximum number of edges in the union of all matchings over all the robust-parity-matchable of short distance graphs given by each molecule.}
Thus, $K \leq \binom{2^s}{2}$\COMMENT{For all $v\in V(I)$, $\mathcal{M}_v$ contains $2^s$ vertices.
Since $\bigcup_{\cS \in \cS(\cM_v)} \mathfrak{M}(\mathcal{S})$ can be seen as a graph on vertex set $V(\mathcal{M}_v)$, we get the claim as a trivial upper bound.}.\\

\textbf{Step~4: Obtaining an appropriate cube factor via the nibble.}
For each $x\in V(I)$, define the sets $A_1(x),\ldots,A_K(x)$ as in Step~4 of the proof of \cref{thm:main1}.
Recall that $G_3\sim \cQ^{n}_{\eta_3}$ and $I(G_3) \sim \cQ^{n-s}_{\eta^{2^s}_3}$.
Apply \cref{thm: nibble} to the graph $I(G_3)$, so that a.a.s.~we obtain a collection $\mathcal{C}$\index{Cw1@$\mathcal{C}$} of vertex-disjoint copies of $\cQ^\ell$ in $I(G_3)$ satisfying \ref{itm:nib1}--\ref{itm:nib2} in the proof of \cref{thm:main1}.
Condition on the event that such a collection $\cC$ exists and call it $\cE_3^*$.\\

\textbf{Step~5: Absorption cubes.}
Recall $G_4 \sim \cQ^n_{\eta_4}$.
We define the event $\cE_4^*$ in the same way as in Step~5 of the proof of \cref{thm:main1}.
Upon conditioning on this event, for each $x \in V(I)$ and $i \in [K]$, we construct the matching $M'(A_i(x))$ in exactly the same way as well.

Consider $A_i(x)$, for some $x\in V(I)$ and $i\in[K]$, and let $x_1$, $x_2$ be the clones of $x$ which correspond to $(x,i)$.
For each $j \in [2]$, let $L^j$ be the layer containing $x_j$.
Similarly as in the proof of \cref{thm:main1}, the following holds. 
\begin{enumerate}[label=$(\mathrm{AB}\arabic*)$]
    \item \label{itm:abss2hit} For each edge $(\hat e, \hat e')\in M'(A_i(x))$, there is an absorbing $\ell$-cube pair $(C^l,C^r)$ for $x$ in $I$ such that, for each $j \in [2]$, the clone $(C^l_j, C^r_j)$ of $(C^l,C^r)$ in $L^j$ is an absorbing $\ell$-cube pair for $x_j$ in $H \cup G_2 \cup G_3 \cup G_4$.
    In particular, the edge joining the left absorber tip to the third absorber vertex lies in $G_4$.
    Moreover, $C^l, C^r \in \cC_x(A_i(x))\subseteq \cC$ and $(C^l,C^r)$ has left and right absorber tips $x+\hat e'$ and $x+\hat e$, respectively.
    Furthermore, for each $x\in V(I)\setminus V_\mathrm{sc}$, these tips lie in $R$.
    We refer to $(C_1^l,C_1^r)$ and $(C_2^l,C_2^r)$ as the absorbing $\ell$-cube pairs for $x_1$ and $x_2$ \emph{associated} with $(\hat e,\hat e')$.
\end{enumerate}

Recall by the construction in Step~3 that, for each $x\in V(I)$ and $i\in[K]$, we have $x_1, x_2 \notin \cU$.
In particular, this means that we do not choose absorbing $\ell$-cube pairs for the vertices in $\cU$ and, thus, the auxiliary edges at the vertices in $\cU$ introduced in the definition of the sets $E^*_H$ in Step~3 will never be used.
As discussed before, the vertices in $\cU$ will instead be incorporated into the Hamilton cycle using the special absorbing structures introduced in \cref{sect:absstruct}.\\

\textbf{Step~6: Removing bondless molecules.}
Recall $G_5\sim \cQ^{n}_{\eta_5}$.
We define the collections $\cC'$, $\cC''$\index{Cw2@$\cC'$}\index{Cw3@$\cC''$} and $\cC_{\mathrm{bs}}$ in the same way as in Step~6 of the proof of \cref{thm:main1}, and we also define the event $\cE^*_5$ in the same way.
We condition on this event and, as before, for each $x\in V(I)$ and each $i\in[K]$, we modify the matching $M'(A_i(x))$ into a matching $M(A_i(x))$ described in \ref{itm:abss1}.\\

\textbf{Step~7: Extending the tree $\boldsymbol{T}$.}
For each $x \in V(I)\setminus B^{5}_I(\cU_I)$, let $Z(x)\coloneqq N_{I}(x)\cap V(T) \cap (\bigcup_{C\in\cC''}V(C))$.
As in the proof of \cref{thm:main1}, we have that
\begin{equation}\label{equa:Zsize}
    |Z(x)| \geq 3n/4.
\end{equation}

Recall $G_6\sim \cQ^{n}_{\eta_6}$. 
We apply \cref{thm: maintreeres} in the same way as in the proof of \cref{thm:main1}, but this time with $2\ell$ and $B^{5}_{I}(\cU_I)$ playing the roles of $\ell$ and $W$, respectively.
Combining this with \ref{itm:tree1hit}, we conclude that\COMMENT{with probability at least $1 - 2e^{-\eps^{6}n}$} a.a.s.~there exists a tree $T'$\index{T2@$T'$} such that $T\subseteq T'\subseteq (I(G_6)\cup T)-B^{5}_I(\cU_I)$ and the following hold:
\begin{enumerate}[label=$(\mathrm{ET}\arabic*)$]
\item\label{itm: ET1hit} $\Delta(T') < D+1$;
\item\label{itm: ET2hit} for all $x \in V(I)$, we have that $|B^{2\ell}_{I}(x)\setminus (V(T')\cup B^{5}_I(\cU_I))| \le n^{3/4}$;
\item\label{itm:ET3hit} for each $x \in V(T')\cap R$, we have that $d_{T'}(x)=1$ and the unique neighbour $x'$ of $x$ in $T'$ is such that $x' \in Z(x)$.
\end{enumerate}
We condition on the above event holding and call it $\cE^*_6$.

As in the proof of \cref{thm:main1}, for each $x \in V(I)$ and each $i \in [K]$, we now redefine the set $M(A_i(x))$ so that \ref{itm:AB3} holds.
It again follows that
\begin{equation}\label{eqn:maix1hit}
    |M(A_i(x))| \ge  n/\ell^2 - n^{3/4} \geq 4n/\ell^3.
\end{equation}

\textbf{Step~8. Consistent systems of paths and cubes.}
Recall $G_{7}\sim \cQ^{n}_{\eta_{7}}$.
For each $v\in V(I)$, let $\cC(v)\coloneqq\{C\subseteq I(G_7):C\cong\cQ^\ell,v\in V(C)\}$.
Let $\cP\coloneqq \{v\in V(I): |\cC(v)|\geq  \lambda n^{\ell}\}$, where $\lambda$ satisfies \eqref{equa:hierarchy2}.
By \cref{remark:RNdegree} (applied with $r=10$ and $\eta_7^{2^s}$ playing the role of $\varepsilon$\COMMENT{And $n-s$ playing the role of $n$, and any value of $a$, e.g. $a=3/4$.}), the following property holds a.a.s.
\begin{enumerate}[label=$(\mathrm{D}\arabic*)$]
    \item\label{itm:degreehit} For every $v\in V(I)$ we have $|B_{I}^{10}(v)\setminus\cP|\leq n^{7/8}$.
\end{enumerate}
For each $v \in \cP$, a straightforward application of \cref{lem: bondedlots} with $\lambda$ and $\cC(v)$ playing the roles of $\gamma$ and $\cC$, respectively, shows that the following holds with probability at least $1-2^{-10n}$:
there exists a subcollection $\cC'(v)\subseteq \cC(v)$, with $|\cC'(v)|\geq \lambda n^{\ell}/4$, with the property that, for every $C\in \cC'(v)$, the molecule $\cM_{C}$ is bonded in $G_7$\COMMENT{Observe that, in order to obtain \ref{itm:degreehit}, we only reveal the edges inside the layers, whereas here we only reveal edges between layers.
So the events are completely independent.}.
By a simple union bound over all vertices in $\cP$, we obtain that the following holds a.a.s.
\begin{enumerate}[label=$(\mathrm{D}\arabic*)$]\setcounter{enumi}{1}
    \item\label{itm:degreehitbonded} For every $v\in \cP$, there exists a collection $\cC'(v)\subseteq\cC(v)$ with $|\cC'(v)|\geq \lambda n^{\ell}/4$ such that, for every $C\in \cC'(v)$, we have that $\cM_C$ is bonded in $G_7$.
\end{enumerate}
Condition on the event that \ref{itm:degreehit} and \ref{itm:degreehitbonded} hold and call it $\mathcal{E}_7^{*}$.

We will show that we may extend many of the consistent systems of paths given by \ref{itm:bad4} into special absorbing structures.
Recall that, since $H$ is $(s,\ell,\varepsilon_1,\varepsilon_2,\gamma,\cU)$-robust, for every $x\in \cU$ and each pair of directions $a,b\in \cD(\cQ^{n})$, there exists a collection $\mathfrak{C}(x,a,b)$ of $(x,a,b)$-consistent systems of paths in $H\cup\{\{x,x+a\}\{x,x+b\}\}$ satisfying \ref{itm:bad4}. 
By \ref{itm:degreehit}, for every $x\in \cU$ and $a,b\in \cD(\cQ^{n})$
\begin{enumerate}[label=(CS)]
    \item\label{itm:conssys} there exists a subcollection $\mathfrak{C}'(x,a,b)\subseteq \mathfrak{C}(x,a,b)$ which satisfies \ref{itm:bad4} with $\gamma/2$ playing the role of $\gamma$ and such that, for every $\mathit{CS}\in \mathfrak{C}'(x,a,b)$, we have $\mathrm{endmol}(\mathit{CS})\subseteq \cP$.
\end{enumerate}

Let $L$ be any layer of $\cQ^n$.
For each $x\in \cU$, let $\mathrm{End}(x)\coloneqq \bigcup_{a,b\in \cD(\cQ^{n})}\bigcup_{\mathit{CS}\in \mathfrak{C}'(x,a,b)} \mathrm{endmol}(\mathit{CS})$, let $\mathrm{End\,II}(x)\coloneqq \bigcup_{a,b\in \cD(\cQ^{n})\setminus\cD(L)}\bigcup_{\mathit{CS}\in \mathfrak{C}'(x,a,b)} \mathrm{endmol}(\mathit{CS})$, and let $x_I\in V(I)$ be such that $x$ is a clone of $x_I$. 
Given any tree $T^\bullet\subseteq I$ and any cube $C\subseteq V(I)$, we say that $C$ \emph{meets $T^\bullet$} if $V(C)\cap V(T^\bullet) \neq \varnothing$.
Recall that, for any $y,z\in V(I)$, we denote the set of directions in a shortest path in $I$ between $y$ and $z$ by $\cD(y,z)$. 

\begin{claim}\label{claim:cubesinendpointsofCSP}
For each $x\in\cU$ and each $z\in \mathrm{End}(x)$, there exists a collection of cubes $\cC''(z)\subseteq \cC'(z)$ with $|\cC''(z)|= 20\ell$ which satisfies the following properties.
\begin{enumerate}[label=$(\mathrm{\roman*})$]
    \item\label{itm:cubesprimepecial3} Every $C\in \cC''(z)$ meets $T'$.
    \item\label{itm:cubesprimepecial4} For every $C\in \cC''(z)$, we have $\cD(C)\cap \cD(x_I,z)=\varnothing$.
    \item\label{itm:cubesprimepecial5} For every $C_1,C_2\in \cC''(z)$, we have $V(C_1)\cap V(C_2)=\{z\}$.
    \item\label{itm:cubesprimepecial6} For each $x\in\cU$, let $\boldsymbol{\mathcal{C}''}(x)\coloneqq\bigcup_{z\in\mathrm{End\,II}(x)}\cC''(z)$.
    For each $d\in\cD(L)$ and each $x\in\cU$, let $\alpha(d,x)\coloneqq|\{C\in\boldsymbol{\mathcal{C}''}(x):d\in\cD(C)\}|$.
    For all $d\in\cD(L)$ and $x\in\cU$, we have that $\alpha(d,x)\leq n/\ell^3$.
\end{enumerate}
\end{claim}

\begin{claimproof}
Fix any vertex $x\in\cU$.
Let $\mathfrak{m}\coloneqq|\mathrm{End}(x)|$ and $m\coloneqq|\mathrm{End\,II}(x)|$, and consider an arbitrary labelling $z_1,\ldots,z_m,z_{m+1},\ldots,z_\mathfrak{m}$ of the vertices in $\mathrm{End}(x)$ such that all vertices in $\mathrm{End\,II}(x)$ come first.
Observe that, by the definition of Type~II consistent systems of paths (see \cref{sect:absstruct}), we have that $m\leq n$\COMMENT{In fact, $m\leq n-s$.}.
We will now iteratively define the sets $\cC''(z_i)$ for each $i\in[\mathfrak{m}]$.

Fix first any $i\in[m]$, and suppose that a set $\cC''(z_j)$ satisfying the claim is already defined for all $j\in[i-1]$.
Let $\boldsymbol{\mathcal{C}''}(x,i)\coloneqq\bigcup_{j=1}^{i-1}\cC''(z_j)$.
For each $d\in\cD(L)$, let $\alpha(d,x,i)\coloneqq|\{C\in\boldsymbol{\mathcal{C}''}(x,i):d\in\cD(C)\}|$.
Let $\cD^*\coloneqq\{d\in\cD(L):\alpha(d,x,i)\geq n/\ell^3\}$\COMMENT{In particular, we must have equality, since each of the cubes we had so far satisfied the statement of the claim.}.
Observe that $|\cD^*|\leq 20\ell^5$\COMMENT{Indeed, if we let $x\coloneqq|\cD^*|$, we have that $xn/\ell^3\leq (i-1)20\ell^2\leq 20\ell^2n$.}.
Let $\cC'''(z_i)$ be the set of all cubes $C\in\cC'(z_i)$ such that $C$ meets $T'$ and $\cD(C)\cap (\cD(x_I,z_i)\cup\cD^*)=\varnothing$ (i.e., they satisfy \ref{itm:cubesprimepecial3} and \ref{itm:cubesprimepecial4} and, if added to the collection, would not violate \ref{itm:cubesprimepecial6}).
We claim that $|\cC'''(z_i)|\geq\lambda n^\ell/5$.
Indeed, observe that $\dist(x_I,z_i)\leq 2$ and, thus, the number of cubes $C\in\cC'(z_i)$ such that $\cD(C)\cap(\cD(x_I,z_i)\cup\cD^*)\neq\varnothing$ is at most $(20\ell^5+2)n^{\ell-1}$. 
(To see this, note that, since $\dist(x_I,z_i)\leq2$, we have that $|\cD(x_I,z_i)|\leq2$.
Now, for each direction $\hat e\in\cD(x_I,z_i)\cup\cD^*$, the number of $\ell$-cubes in $\cQ^n$ containing $z$ and some edge in direction $\hat e$ is $\binom{n}{\ell-1}<n^{\ell-1}$.)
Moreover, by \ref{itm: ET2hit} and \ref{itm:bad3}, the number of such cubes which do not meet $T'$ is at most $n^{\ell-2}$\COMMENT{By \ref{itm: ET2hit} (here we also use \ref{itm:bad3} to guarantee that there are no other vertices of $\cU$ close to $x$), the number of vertices $y\in B^\ell_I(x)\setminus B^5_I(x)$ which do not belong to $T'$ is at most $n^{3/4}$.
Each cube $C\in\cC'(z_i)$ contains many vertices in this set.
For each such vertex $y$, there are at most $\binom{n^{3/4}}{\ell}n^5 < n^{\ell-5}$ cubes $C\in\cC'(z_i)$ such that $y\in V(C)$.
The claim above is a very crude upper bound.}.
The bound then follows by \ref{itm:degreehitbonded}.

We can now construct $\cC''(z_i)$ by obtaining cubes $C_1(z_i),\ldots, C_{20\ell}(z_i)$ iteratively.
Note that, for any pair of cubes $C_1,C_2\in\cC'(z_i)$, we have that $z\in V(C_1)\cap V(C_2)$.
Then, \ref{itm:cubesprimepecial5} is equivalent to having that $\cD(C_1)\cap\cD(C_2)=\varnothing$.
For each $k\in[20\ell]$, we proceed as follows.
Let $\cD'_k\coloneqq \bigcup_{j=1}^{k-1}\cD(C_j(z_i))$. 
Note that $|\cD'_k|\leq 20\ell^2\leq \lambda n/8$. 
Now, applying \cref{prop:cubesdirect} with $n-s$, $z$, $\lambda/5$, $\lambda/8$ and $\cD'_k$ playing the roles of $n$, $x$, $\eta$, $\eta'$ and $\cD'$, respectively, we deduce that there is a cube $C_k(z_i)\in\cC'''(z_i)$ with $\cD(C_k(z_i))\cap\cD'_k=\varnothing$.
By enforcing that \ref{itm:cubesprimepecial5} holds, it follows that each direction is used at most once in the cubes that were added in this step. 
It then follows that \ref{itm:cubesprimepecial6} holds as well.

Consider now any $i\in[\mathfrak{m}]\setminus[m]$, and suppose that a set $\cC''(z_j)$ satisfying the claim is already defined for all $j\in[i-1]$.
Let $\cC'''(z_i)$ be the set of all cubes $C\in\cC'(z_i)$ such that $C$ meets $T'$ and $\cD(C)\cap \cD(x_I,z_i)=\varnothing$ (i.e., they satisfy \ref{itm:cubesprimepecial3} and \ref{itm:cubesprimepecial4}).
As above, we claim that $|\cC'''(z_i)|\geq\lambda n^\ell/5$.
Indeed, the number of cubes $C\in\cC'(z_i)$ such that $\cD(C)\cap \cD(x_I,z_i)\neq\varnothing$ is at most $2n^{\ell-1}$\COMMENT{Note that $\dist(x_I,z_i)\leq2$. 
It follows that $|\cD(x_I,z_i)|\leq2$.
Now we must remove from $\cC'(z_i)$ all cubes that contain one of these directions, of which, for each direction, there are at most $\binom{n}{\ell-1}<n^{\ell-1}$ in $\cQ^n$.}, and, again, the number of such cubes which do not meet $T'$ is at most $n^{\ell-2}$\COMMENT{By \ref{itm: ET2hit} (here we also use \ref{itm:bad3} to guarantee that there are no other vertices of $\cU$ close to $x$), the number of vertices $y\in B^\ell_I(x)\setminus B^5_I(x)$ which do not belong to $T'$ is at most $n^{3/4}$.
Each cube $C\in\cC'(z)$ contains many vertices in this set.
For each such vertex $y$, there are at most $\binom{n^{3/4}}{\ell}n^5 < n^{\ell-5}$ cubes $C\in\cC'(z)$ such that $y\in V(C)$.
The claim above is a very crude upper bound.}.
The bound then follows by \ref{itm:degreehitbonded}.

We can now construct $\cC''(z_i)$ as above.
For each $k\in[20\ell]$, we proceed as follows.
Let $\cD'_k\coloneqq \bigcup_{j=1}^{k-1}\cD(C_j(z_i))$. 
Note that $|\cD'_k|\leq 20\ell^2\leq \lambda n/8$. 
Now, applying \cref{prop:cubesdirect} with $n-s$, $z$, $\lambda/5$, $\lambda/8$ and $\cD'_k$ playing the roles of $n$, $x$, $\eta$, $\eta'$ and $\cD'$, respectively, we deduce that there is a cube $C_k(z_i)\in\cC'''(z_i)$ with $\cD(C_k(z_i))\cap\cD'_k=\varnothing$.
\end{claimproof}

Let $J^{1}\coloneqq \bigcup_{x \in \cU}\bigcup_{z\in \mathrm{End}(x)}\bigcup_{C\in \cC''(z)} \cM_C$, where $\cC''(z)$ are the sets given by \cref{claim:cubesinendpointsofCSP}, and let $G^{*}_7\subseteq G_7$ consist of all edges of $G_7$ which have endpoints in different layers.

\begin{claim}\label{itm: boundeddegreerand}
$J^{1}\cup G^{*}_7$ is $(\cU,\ell^3,s)$-good and $\Delta(J^{1}\cup G^{*}_7)\leq 50\ell^{4}$.
\end{claim}

\begin{claimproof}
In order to see that $J^{1}\cup G^{*}_7$ is $(\cU,\ell^3,s)$-good, observe first that the edges of $G^*_7$ do not affect this definition, so it suffices to see that $J^{1}$ is $(\cU,\ell^3,s)$-good.
By \cref{claim:cubesinendpointsofCSP}\ref{itm:cubesprimepecial4}, for all $x\in\cU$, $z\in\mathrm{End}(x)$ and $C\in\cC''(z)$ we have that $\dist(x_I,C)=\dist(x_I,z)$.
In particular, by the definition of the different consistent systems of paths (see \cref{sect:absstruct}), it follows that the only cubes which affect whether $J^{1}$ is $(\cU,\ell^3,s)$-good or not are those of the collection $\boldsymbol{\mathcal{C}''}(x)$ described in \cref{claim:cubesinendpointsofCSP}\ref{itm:cubesprimepecial6}.
But then, by \cref{claim:cubesinendpointsofCSP}\ref{itm:cubesprimepecial6}, we have that no direction $d\in\cD(L)$ is used more than $n/\ell^3$ times, as required.

Note that $\Delta(G_7^{*})\leq s=10\ell$, by construction.
We will now show that $\Delta(J^1)\leq 45\ell^{4}$.
Observe that $J^1$ does not contain any edges with endpoints in different layers.
In particular, $J^1$ consists of clones of the same subgraph of $I$, that is $I(J^1)=\bigcup_{x \in \cU}\bigcup_{z\in \mathrm{End}(x)}\bigcup_{C\in \cC''(z)}C$.
By this observation, it is enough to show that $\Delta(I(J^1))\leq45\ell^{4}$.

Recall that, by \ref{itm:bad3}, given any distinct $x,y\in\cU_I$, we have that $\dist(x,y)\geq\gamma n/2$.
In particular, by this observation and \cref{claim:cubesinendpointsofCSP}\ref{itm:cubesprimepecial4}, it follows that, for all $x\in\cU_I$, we have $d_{I(J^1)}(x)=0$.
We also note that, for every $z\in V(I)$ for which $\dist(z,\cU_I)\geq \ell+3$, we have $d_{I(J^1)}(z)=0$.\COMMENT{Using that in $I$ any vertex in a consistent system for $x$ is at most distance 2 from $x$.}
Now, fix any $x\in \cU_I$ and $z\in V(I)$ such that $\dist(z,x)=t$, for some $t\in[\ell+2]$.
We claim that $d_{I(J^1)}(z)\leq 2t^{2}\cdot 20\ell^{2}\leq 45\ell^4$.

Suppose first that $t=1$.
Then, by \cref{claim:cubesinendpointsofCSP}\ref{itm:cubesprimepecial4}, for every edge $e\in E(I(J^1))$ incident with $z$, we have $e\in E(C)$ for some $C\in \cC''(z)$, and hence $d_{I(J^1)}(z)\leq \ell |\cC''(z)|=20\ell^{2}$, as we wanted to show. 
Suppose now that $t\geq 2$ and let $\cD(z,x)=\{d_1,\ldots, d_t\}$.
Every edge $e$ incident with $z$ must come from the edges of a cube $C\in \cC''(w)$, for some $w\in \mathrm{End}(x)$.
Moreover, by \cref{claim:cubesinendpointsofCSP}\ref{itm:cubesprimepecial4}, we must have $\cD(x,w)\subseteq \cD(x,z)$.
As there are at most $t+t^{2}\leq 2t^{2}$ vertices $w\in V(I)$ such that $\cD(x,w)\subseteq \cD(x,z)$ and $\dist(x,w)\in[2]$, we have have that $d_{I(J^1)}(z)\leq 2t^{2} 20\ell^{2}$, which concludes the proof of the claim. 
\end{claimproof}

\textbf{Step~9: Fixing a collection of absorbing $\boldsymbol{\ell}$-cube pairs for the vertices in scant molecules.}
In this step, we use $G_8 \sim \cQ^n_{\eta_8}$ to alter $T'$ so that any tips of absorbing $\ell$-cube pairs for vertices $x\in V_\mathrm{sc}$ which do not lie in $R$ are relocated from the tree $T'$ to the reservoir.
We follow the same approach as in Step~8 of the proof of \cref{thm:main1}.
In particular, we define an event $\cE_8^*$ (which is analogous to $\cE_7^*$ in the proof of \cref{thm:main1}) and condition that it holds.
This then gives a set of absorbing $\ell$-cube pairs $\cC^\mathrm{sc}_1=\{(C^l(x,j,k),C^r(x,j,k))\subseteq I:x\in V_\mathrm{sc}, j\in[K], k\in[2^{s+1}\Psi]\}$, where each $(C^l(x,j,k),C^r(x,j,k))$ is an absorbing $\ell$-cube pair for $x$, which satisfies that 
\begin{enumerate}[label=$(\mathrm{CD})$]
    \item\label{itm:DisjAbsCubeshit} for all distinct $(x,j,k),(x',j',k')\in V_\mathrm{sc}\times[K]\times[2^{s+1}\Psi]$, $C^l(x,j,k)$ and $C^r(x,j,k)$ are both vertex-disjoint from $C^l(x',j',k')$ and $C^r(x',j',k')$.
\end{enumerate}

We define $P'$ and $P$ as in Step~8 of the proof of \cref{thm:main1}.
Observe that \ref{itm: bad and scant}\COMMENT{and the assumption on the size of the graphs $P(x_i,j,k,M'_{i,j,k},B_{i,j,k})$, which are no longer defined here...} implies that $(P'\cup V(P))\cap B^{2\ell}(\cU_I)=\varnothing$.\COMMENT{This holds because all graphs $P(x_i,j,k,M'_{i,j,k},B_{i,j,k})$ are at distance at most $\ell$ from $V_\mathrm{sc}$ and all vertices in $P'$ are at distance $1$ from $V_\mathrm{sc}$, and by assumption $\dist(V_\mathrm{sc}, \cU_I)\geq 5\ell$.}
Let $T^\textsc{iv}\coloneqq T'[V(T')\setminus P'] \cup P$, which is connected, and let $T''$ be a spanning tree of $T^\textsc{iv}$.
In particular, it follows from the above and the definitions of $T$ and $T'$ in Steps~1 and 7 that\index{T3@$T''$}
\begin{equation}\label{equa:T''containedhit}
    T''\subseteq I(G_1\cup G_6\cup G_8)-B_I^5(\cU_I)\subseteq I-B_I^5(\cU_I).
\end{equation}
Furthermore, as in Step~8 of the proof of \cref{thm:main1}, we have that
\begin{equation}\label{equa:T''boundhit}
    \Delta(T'')\leq 12D.
\end{equation}
Define the (new) reservoir $R' \coloneqq (R \cup P')\setminus V(P)$\index{R2@$R'$}.

For each $x \in V(I)\setminus B^{5}_I(\cU_I)$, let $Z'(x)\coloneqq Z(x)\cap V(T'')$ (where $Z(x)$ is as defined in Step~7).
It follows by \eqref{equa:Zsize} and \ref{itm:scanthit} that
\[|Z'(x)| \geq 3n/4 - 4\cdot2^s\Psi KS' \geq n/2.\]
Choose any vertex $x \in V(I)$ with $\dist(x,\cU_I)\geq3\ell$\COMMENT{Note that such a vertex exists by \ref{itm:bad3}.}.
Again by \ref{itm:scanthit}, there are at most $4\cdot2^{s+\ell}\Psi KS'$ vertices in $Z'(x)$ which lie in cubes of absorbing $\ell$-cube pairs of $\cC^\mathrm{sc}_1$.
Choose any vertex $y\in Z'(x)$ which does not lie in any of those cubes.
Denote the cube $C\in\cC''$ which contains $y$ by $C^\bullet$.

For each $x \in V(I)\setminus V_\mathrm{sc}$ and each $i \in [K]$, we now redefine the set $M(A_i(x))$ as follows.
\begin{enumerate}[label=$(\mathrm{AB}\arabic*)$]\setcounter{enumi}{3}
\item\label{itm:AB4hit} Let $M(A_i(x))$ retain only those edges whose associated absorbing $\ell$-cube pair $(C^l, C^r)$ satisfies that both $C^l$ and $C^r$ are different from $C^\bullet$ and vertex-disjoint from both cubes of all absorbing $\ell$-cube pairs of $\cC^\mathrm{sc}_1$, and both tips $x^l$ and $x^r$ satisfy that $x^l,x^r\in R\setminus V(P)\subseteq R'$.
\end{enumerate}
Note that, by \ref{itm:scanthit}, we have $|B^{\ell+1}_I(x)\cap V(P)|\leq21\cdot2^s\Psi DKS'$ and $|B^{\ell+1}_I(x)\cap V(\bigcup_{(C^l,C^r)\in\cC^\mathrm{sc}_1}(C^l\cup C^r))|\leq4\cdot2^{\ell+s}\Psi KS'$.
Combining this with \eqref{eqn:maix1hit} and \ref{itm:abss2hit}, it follows that 
\begin{equation}\label{eqn:maix2hit}
    |M(A_i(x))| \geq 4n/\ell^3-(21D+4\cdot2^\ell)2^s\Psi KS' - 1 \geq 2n/\ell^3.
\end{equation}

\textbf{Step~10: Fixing a collection of absorbing $\boldsymbol{\ell}$-cube pairs for vertices in non-scant molecules and vertices near $\boldsymbol{\cU}$.}
At this point, it is not yet clear which vertices will need to eventually be absorbed into the long cycle we construct.
For vertices in $I$ which are `far' from $\cU_I$, we can already determine those which will have clones that will need to be absorbed (though we cannot yet determine the precise clones).
However, for vertices which are `near' $\cU_I$, we still cannot say which of them will have clones that need to be absorbed (this depends on the special absorbing structure which is fixed once the edges of $H'$ are revealed). 
As a result, we proceed as if all of the clones of vertices in $I$ near $\cU_I$ will need to be absorbed.
Recall that $\cC'$ and $\cC''$ were introduced in Step~6.
Let\index{Cw4@$\cC'''$}\index{Vabs1@$V_\mathrm{abs}'$} 
\[\cC'''\coloneqq \{C\in \cC': V(C)\cap V(T'')\neq \varnothing\}\quad\text{ and }\quad V_\mathrm{abs}'\coloneqq (V(I)\setminus \bigcup_{C\in \cC'''} V(C))\cup B_I^{3\ell}(\cU_I).\]
We will now fix a collection of absorbing $\ell$-cube pairs for all vertices in each vertex molecule $\mathcal{M}_x$ with $x\in V_\mathrm{abs}'\setminus V_\mathrm{sc}$, except for the vertices of $\cU$.

Similarly as in \eqref{eqn:moveC}, we have
\begin{equation}\label{eqn:moveChit}
|B^{2\ell}_{I}(x)\setminus (V(T'')\cup B_I^5(\cU_I))| \leq 2n^{3/4}.
\end{equation}
For all $x \in \bigcup_{C\in\cC''}V(C)$, by combining \ref{itm:nib1}, \eqref{eqn:moveChit} and \ref{itm:bad3} with the definition of bondlessly surrounded molecules, we claim that
\begin{equation}\label{eqn:movedCihit}
|\cC_{x}\cap \cC'''|\geq(1-2^{-\ell-5s+1})n,
\end{equation}
where $\cC_x$ is the collection of all those $C\in\cC$ such that $x\notin V(C)$ and $N_I(x)\cap V(C)\neq\varnothing$.
Indeed, combining \ref{itm:nib1} with the definition of bondlessly surrounded, it follows that $|\cC_x\cap\cC'|\geq(1-\delta)n-n/2^{\ell+5s}$. 
Now, since $|B_I^{2\ell}(x)\cap\cU_I|\leq1$ by \ref{itm:bad3}, every cube in $\cC_x$ has at least half of its vertices outside $B_{I}^5(\cU_I)$.
For any cube in $\cC_x\cap\cC'$ to not belong in $\cC'''$, these at least $2^{\ell-1}$ vertices which do not lie in $B_{I}^5(\cU_I)$ must also avoid $V(T'')$.
Therefore, by \eqref{eqn:moveChit}, at most $2n^{3/4}$ of these cubes in $\cC_x\cap\cC'$ to not lie in $\cC'''$.\COMMENT{This is a trivial bound, we could in fact divide this by $2^{\ell-1}$.}
Thus, in total, we have $|\cC_{x}\cap \cC'''|\geq(1-\delta)n-n/2^{\ell+5s}-2n^{3/4}\geq(1-2\cdot2^{-\ell-5s})n$.

Recall that, for any $x\in V(I)$, each index $i\in[K]$ is given by a unique edge $e \in \mathfrak{M}(x)$ via the relation $N(e) = A_i(x)$.
Recall also the definition of $\mathfrak{M}(x)$ from Step~3.
We now prove the following claim, which is similar to \cref{claim:rainbowapplnew} (apart from the new property \ref{item:rainbowapplnew2hit}).

\begin{claim}\label{claim:rainbowapplnewhit}
For each $x\in V_\mathrm{abs}'\setminus V_\mathrm{sc}$ and each $e\in \mathfrak{M}(x)$, there exists a set $\mathcal{C}_1^\mathrm{abs}(e)$ of $2^{s+1}\Psi$ absorbing $\ell$-cube pairs $(C_{k}^l(e),C_{k}^r(e))\subseteq I$, one for each $k\in[2^{s+1}\Psi]$, which satisfies the following:
\begin{enumerate}[label=$(\mathrm{\roman*})$]
    \item\label{item:rainbowapplnew1hit} for all $x\in V_\mathrm{abs}'\setminus V_\mathrm{sc}$, $e\in \mathfrak{M}(x)$ and $k\in[2^{s+1}\Psi]$, the absorbing $\ell$-cube pair $(C_{k}^l(e),C_{k}^r(e))$ is associated with some edge in $M(A_j(x))$, for some $j \in [K]$;
    \item\label{item:rainbowapplnew2hit} for all $x\in B_I^{3\ell}(\cU_I)$\COMMENT{Note that, by \ref{itm: bad and scant}, there is no need to remove scant molecules from here.}, $e\in \mathfrak{M}(x)$ and $k\in[2^{s+1}\Psi]$, the absorbing $\ell$-cube pair $(C_{k}^l(e),C_{k}^r(e))$ has tips $x_{k}^l(e)$ and $x_{k}^r(e)$ which satisfy that $\dist(x,\cU_I)<\dist(x_{k}^l(e),\cU_I)$ and $\dist(x,\cU_I)<\dist(x_{k}^r(e),\cU_I)$, and
    \item\label{item:rainbowapplnew3hit} for all $x,x'\in V_\mathrm{abs}'\setminus V_\mathrm{sc}$, all $e\in \mathfrak{M}(x)$ and $e'\in \mathfrak{M}(x')$, and all $k,k'\in [2^{s+1}\Psi]$ with $(x,e,k)\neq (x',e',k')$, the absorbing $\ell$-cube pairs $(C_{k}^l(e),C_{k}^r(e))$ and $(C_{k'}^l(e'),C_{k'}^r(e'))$ satisfy that $(V(C_{k}^l(e))\cup V(C_{k}^r(e)))\cap(V(C_{k'}^l(e'))\cup V(C_{k'}^r(e')))=\varnothing$.
\end{enumerate}
\end{claim}

\begin{claimproof}
Let $\cV \coloneqq\bigcup_{x\in V_\mathrm{abs}'\setminus V_\mathrm{sc}}\mathfrak{M}(x)$.
Let $K'\coloneqq|\cV|$, and let $f_1, \dots, f_{K'}$ be an ordering of the edges in $\cV$.
Given any $i \in [K']$, the edge $f_{i}$ corresponds to a pair $(x,j(i))$, where $x\in V_\mathrm{abs}'\setminus V_\mathrm{sc}$ and $j(i)\in [K]$.
If $x\notin B_I^{3\ell}(\cU_I)$, let $\mathfrak{C}_i$ be the collection of absorbing $\ell$-cube pairs for $x$ in $I$ associated with some edge of $M(A_{j(i)}(x))$.
Otherwise, let $\mathfrak{C}_i$ be the same collection, after removing all those absorbing $\ell$-cube pairs for which \ref{item:rainbowapplnew2hit} does not hold.
Since each $x\in B_I^{3\ell}(\cU_I)$ has at most $3\ell$ neighbours $y\in N_I(x)$ such that $\dist(x,\cU_I)\geq\dist(y,\cU_I)$, it follows by \eqref{eqn:maix2hit} that $|\mathfrak{C}_i|\geq n/\ell^3$ for all $i\in[K']$.
In particular, by \ref{itm:abss1}, each of the absorbing $\ell$-cube pairs $(C^l, C^r)$ in any of the collections $\mathfrak{C}_i$ satisfies that $C^l, C^r \in \cC''$.

Let $\cH$ be the $2^{s+1}\Psi K'$-edge-coloured auxiliary multigraph with $V(\mathcal{H})\coloneqq\mathcal{C}''$, which contains one edge of colour $(i,k)\in [K']\times [2^{s+1}\Psi]$ between $C$ and $C'$ whenever $(C,C') \in \mathfrak{C}_i$ or $(C',C)\in \mathfrak{C}_i$.
In particular, $\cH$ contains at least $n/\ell^3$ edges of each colour.
We now bound $\Delta(\cH)$.
Consider any $C \in V(\cH)$.
Note that, for each edge $e$ of $\cH$ incident to $C$, there exists some $x=x(e) \in V_\mathrm{abs}'\setminus V_\mathrm{sc}$ such that $C$ together with some other cube $C' \in V(\cH)$ forms an absorbing $\ell$-cube pair for $x$.
In particular, $x$ must be adjacent to $C$ in $I$.
Let $\overline\partial(C)\coloneqq\{x(e):e\in E(\cH) \text{ is incident to }C\}$.
Moreover, if $e$ has colour $(i,z)$, then $f_i\in \mathfrak{M}(x)$ (and $f_i$ has corresponding pair $(x,j(i))$ for some $j(i)\in[K]$).
Since $f_i \in \mathfrak{M}(x)$ and $|\mathfrak{M}(x)|\leq \binom{2^{s}}{2}$, it follows that each vertex $y$ which is adjacent to $C$ in $I$ can play the role of $x$ for at most $2^{s+1}\Psi\cdot 2^{2s}$ edges of $\cH$ incident to $C$.
Thus, $d_{\cH}(C)\leq2^{3s+1}\Psi|\overline\partial(C)|$.

Fix a cube $C\in V(\cH)$.
In order to bound $|\overline\partial(C)|$, consider first $|\overline\partial(C)\cap B_I^{3\ell}(\cU_I)|$.
Recall that, by \ref{itm:bad3}, there is at most one vertex $z\in\cU_I\cap B_I^{3\ell}(V(C))$.
Furthermore, since the property described in \ref{item:rainbowapplnew2hit} holds for all absorbing $\ell$-cube pairs for vertices in $B_I^{3\ell}(\cU_I)$ represented in $\cH$, it follows that each vertex $x\in V(C)$ has at most $3\ell+1$ neighbours in $B_I^{3\ell}(\cU_I)\cap\overline\partial(C)$.
Thus, in total, $|\overline\partial(C)\cap B_I^{3\ell}(\cU_I)|\leq(3\ell+1)2^\ell$.
Consider now $|\overline\partial(C)\setminus B_I^{3\ell}(\cU_I)|$.
Note that $\overline\partial(C)\setminus B_I^{3\ell}(\cU_I)\subseteq (V_{\mathrm{abs}}'\cap N_I(V(C)))\setminus B_I^{3\ell}(\cU_I)\subseteq N_I(V(C)) \setminus \bigcup_{C'\in \cC'''} V(C')$.
By \eqref{eqn:movedCihit}, the number of vertices in $V_{\mathrm{abs}}'\setminus B_I^{3\ell}(\cU_I)$ which are adjacent to $C$ is at most $2|C|n/{2^{\ell +5s}}$, that is, $|\overline\partial(C)\setminus B_I^{3\ell}(\cU_I)|\leq2n/2^{5s}$.
We conclude that $|\overline\partial(C)|\leq 3n/2^{5s}$\COMMENT{$|\overline\partial(C)|=|\overline\partial(C)\cap B_I^{3\ell}(\cU_I)|+|\overline\partial(C)\setminus B_I^{3\ell}(\cU_I)|\leq(3\ell+1)2^\ell+2n/2^{5s}\leq 3n/2^{5s}$.} and, thus, $d_{\cH}(C) \le 2^{3s+1}\Psi 3n/{2^{5s}}\le n/\ell^4$.\COMMENT{Here we used that $\Phi=60\ell^4$ and $\Psi=c\Phi\leq 60\ell^5$.}

Since each colour class has size at least $n/\ell^3$ and $\Delta(\mathcal{H})\leq n/\ell^4$, by \cref{lem: rainbow}, $\mathcal{H}$ contains a rainbow matching of size $2^{s+1}\Psi K'$.\COMMENT{We apply \cref{lem: rainbow} with $r= 2$, $m=n/\ell^3$, and note that $\Delta(\cH) < n/\ell^{4} < m/(6r)$.}
For each $(i,z)\in [K']\times [2^{s+1}\Psi]$, let $(C_{z}^l(f_i),C_{z}^r(f_i))\in\mathfrak{C}_i$ be the absorbing $\ell$-cube pair of colour $(i,z)$ in this rainbow matching.
This ensures that \ref{item:rainbowapplnew3hit} holds, while \ref{item:rainbowapplnew1hit} and \ref{item:rainbowapplnew2hit} follow by construction.
\end{claimproof}

For each $x\in V_\mathrm{abs}'\setminus V_\mathrm{sc}$ and each $i\in[K]$, let $\cC^\mathrm{abs}_1(x,i)\coloneqq\cC^\mathrm{abs}_1(e)$ be the set of absorbing $\ell$-cube pairs guaranteed by \cref{claim:rainbowapplnewhit}, where $e\in\mathfrak{M}(x)$ is the unique edge such that $A_i(x)=N(e)$.
Similarly, for each $k\in[2^{s+1}\Psi]$, let $(C^l(x,i,k),C^r(x,i,k))\coloneqq(C_{k}^l(e),C_{k}^r(e))$.
Let $\cC^\mathrm{abs}_1\coloneqq\bigcup_{x\in V_\mathrm{abs}'\setminus V_\mathrm{sc}}\bigcup_{i\in[K]}\cC^\mathrm{abs}_1(x,i)$.\COMMENT{This change of notation is not needed here or in the next paragraph, as taking union over all $i$ is the same as taking unions over all $e$. However, this notation is more convenient later for the list of properties \ref{itm:C6} and their proof, as we want to have consistent notation there.}\\

Let $G\coloneqq\bigcup_{i=1}^{8} G_i$.
Recall that $G^*_7$ and $J^1$ were defined in Step~8.
We let $Q'\subseteq G$ be the spanning subgraph with edge set 
\[E(Q')\coloneqq E(J^1)\cup E(G_4^{*})\cup E(G_5^*)\cup E(G_7^{*})\cup\bigcup_{C\in\cC'}E(\cM_C)\cup\bigcup_{i=1}^{2^s}E(T''_{L_i}),\]
where $G_4^*$ and $G_5^*$ are as defined in Step~9 of the proof of \cref{thm:main1} (but with $V_\mathrm{abs}'$ playing the role of $V_\mathrm{abs}$).
Note that $\Delta(G_4^{*})=1$, $\Delta(G_5^*)\leq s=10\ell$ and, since the molecules $\cM_C$ with $C\in\cC'$ are vertex disjoint, each vertex in $\bigcup_{C\in\cC'}E(\cM_C)$ has degree at most $\ell$.
Combining these bounds with \cref{itm: boundeddegreerand} and \eqref{equa:T''boundhit}, we have that $\Delta(Q')\leq 50\ell^4+1+s+\ell+12D\leq60\ell^4=\Phi$.

\begin{claim}\label{claim:Q'good}
$Q'$ is $(\cU,2\ell^2,s)$-good.
\end{claim}

\begin{claimproof}
Indeed, observe that this fact only depends on those edges contained within a layer which are incident to a neighbour of $x$ in $\cQ^n$, for some $x\in\cU$.
Therefore, the graph $G^*_5$ has no effect here.
By property \ref{itm:C4hit} below, the graph $\bigcup_{i=1}^{2^s}T''_{L_i}$ also has no effect.
Now, by \cref{itm: boundeddegreerand} we have that $J^1\cup G_7^*$ is $(\cU,\ell^3,s)$-good, and $\bigcup_{C\in\cC'}\cM_C$ is $(\cU,\ell^3,s)$-good by \ref{itm:nib3} combined with \ref{itm:bad3}.
Finally, consider $G^*_4$.
For each $x\in V_\mathrm{abs}'\cup V_\mathrm{sc}$, $i\in[K]$ and $k\in[2^{s+1}\Psi]$, let $e(x,i,k)$ be the edge between the left absorber tip and the third absorber vertex of $(C^l(x,i,k),C^r(x,i,k))\in\cC^\mathrm{sc}_1\cup\cC^\mathrm{abs}_1$.
Observe that, for all $x\in V_\mathrm{abs}'\cup V_\mathrm{sc}$ such that $\dist(x,\cU_I)\geq5$, all $i\in[K]$ and all $k\in[2^{s+1}\Psi]$, we have that $e(x,i,k)$ does not affect whether $Q'$ is $(\cU,2\ell^2,s)$-good or not.
In particular, by \ref{itm: bad and scant}, this is true for all $x\in V_\mathrm{sc}$.
Now consider each $x\in V_\mathrm{abs}'$ such that $\dist(x,\cU_I)<5$.
By \cref{claim:rainbowapplnewhit}\ref{item:rainbowapplnew2hit}, it follows that $e(x,i,k)$ only affects our claim when $x\in\cU_I$.
Observe that, for each $i\in[K]$ and $k\in[2^{s+1}\Psi]$, the direction of $e(x,i,k)$ is the same as that of the edge $e'(x,i,k)$ joining $x$ to the right absorber tip of $(C^l(x,i,k),C^r(x,i,k))$.
By \cref{claim:rainbowapplnewhit}\ref{item:rainbowapplnew3hit}, all cubes of absorbing $\ell$-cube pairs in $\bigcup_{i\in[K]}\cC^\mathrm{abs}_1(x,i)$ are vertex disjoint, which implies that each edge $e'(x,i,k)$ with $i\in[K]$ and $k\in[2^{s+1}\Psi]$ uses a different direction.
Hence, $G^*_4$ is $(\cU,n,s)$-good, and the claim follows.
\end{claimproof}

Note that $T'' \subseteq I(Q')$, $R'\subseteq V(I)$, and $C\subseteq I(Q')$ for all $C\in\cC'$.
Recall the definitions of $\cC''$ from Step~6 and $\cC'''$ from Step~10.
For any $u\in\cU$, recall the definitions of $\mathrm{End}(x)$ given in Step~8.
Recall also the definitions of $P$, $P'$ and $C^\bullet$ from Step~9.
Combining all the previous steps, we claim that the following hold (conditioned on the events $\cE_1^*,\ldots,\cE_8^*$, which occur a.a.s.). 
\begin{enumerate}[label=$(\mathrm{C}\arabic*)$]
    \item\label{itm:C2hit} $\Delta(T'')\leq 12D$.
    \item\label{itm:C1hit} Any vertex $x\in R'\cap V(T'')$ is a leaf of $T''$.
    Furthermore, if $x\in R'\cap V(T'')$, then $x\notin V(T)$ and its unique neighbour $x'$ in $T''$ satisfies that $x'\in Z(x)$ (where $Z(x)$ is as defined in Step~7).
    \item\label{itm:C4hit} $B_{I}^{5}(\cU_{I})\cap V(T'')=\varnothing$.
    \item\label{itm:C52hit} For all $x \in V(I)$ we have that $|\cC_x \cap \cC'''| \geq (1-3/2\ell^4)n$.
    \item\label{itm:C3hit} For each $x \in V_{\mathrm{sc}}$ and $i \in [K]$, there is a collection $\cC^\mathrm{sc}_1(x,i)$ of $2^{s+1}\Psi$ absorbing $\ell$-cube pairs $(C^l(x,i,k),C^r(x,i,k))$ for $x$ in $I$ (defined in Step~9), each of which is associated with some edge $e \in M(A_i(x))$.
    In particular, $(C^l(x,i,k),C^r(x,i,k))$ is as described in \ref{itm:abss2hit} (recall also \ref{itm:abss1}), that is, there are two absorbing $\ell$-cube pairs $(C^l_1(x,i,k),C^r_1(x,i,k))$ and $(C^l_2(x,i,k),C^r_2(x,i,k))$ in $H\cup G$, associated with $e\in M(A_i(x))$, for the clones $x_1$ and $x_2$ of $x$ which correspond to $(x,i)$.
    Moreover, each of these absorbing $\ell$-cube pairs $(C^l(x,i,k),C^r(x,i,k))$ satisfies the following:
    \begin{enumerate}[label=$(\mathrm{C}5.\arabic*)$]
        \item\label{itm:C30hit} $(C^l_1(x,i,k), C^r_1(x,i,k)) \cup (C^l_2(x,i,k), C^r_2(x,i,k)) - V(\cM_x) \subseteq Q'$;
        \item\label{itm:C31hit} the tips of $C^l(x,i,k)$ and $C^r(x,i,k)$ lie in $R'\setminus V(T'')$;
        \item\label{itm:C325hit} $C^l(x,i,k),C^r(x,i,k)\in\mathcal{C}''\cap\cC'''$, and
        \item\label{itm:C33hit} for any $x' \in V_{\mathrm{sc}}$, $i' \in [K]$ and $k'\in[2^{s+1}\Psi]$ with $(x',i',k')\neq(x,i,k)$, we have that $C^l(x,i,k)$, $C^r(x,i,k)$, $C^l(x',i',k')$ and $C^r(x',i',k')$ are vertex-disjoint.
    \end{enumerate}
    \item\label{itm:C6hit} For each $x \in V_\mathrm{abs}'\setminus V_{\mathrm{sc}}$ and $i \in [K]$, there is a collection $\cC^\mathrm{abs}_1(x,i)$ of $2^{s+1}\Psi$ absorbing $\ell$-cube pairs $(C^l(x,i,k),C^r(x,i,k))$ for $x$ in $I$ (defined in Step~10), each of which is associated with an edge $e\in M(A_i(x))$.
    In particular, $(C^l(x,i,k),C^r(x,i,k))$ is as described in \ref{itm:abss2hit} (recall also \ref{itm:abss1}), that is, there are two absorbing $\ell$-cube pairs $(C^l_1(x,i,k),C^r_1(x,i,k))$ and $(C^l_2(x,i,k),C^r_2(x,i,k))$ in $H\cup G$, associated with $e\in M(A_i(x))$, for the clones $x_1$ and $x_2$ of $x$ which correspond to $(x,i)$.
    Moreover, each of these absorbing $\ell$-cube pairs $(C^l(x,i,k),C^r(x,i,k))$ satisfies the following:
    \begin{enumerate}[label=$(\mathrm{C}6.\arabic*)$]
        \item\label{itm:C60hit} $(C^l_1(x,i,k), C^r_1(x,i,k)) \cup (C^l_2(x,i,k), C^r_2(x,i,k)) - V(\cM_x) \subseteq Q'$;
        \item\label{itm:C62hit} the tips of $C^l(x,i,k)$ and $C^r(x,i,k)$ lie in $R'$;
        \item\label{itm:C64hit} $C^l(x,i,k),C^r(x,i,k)\in\mathcal{C}''\cap\cC'''$;
        \item\label{itm:C65hit} for any $x' \in V_\mathrm{abs}'\setminus V_{\mathrm{sc}}$, $i' \in [K]$ and $k'\in[2^{s+1}\Psi]$ with $(x',i',k')\neq(x,i,k)$, we have that $C^l(x,i,k)$, $C^r(x,i,k)$, $C^l(x',i',k')$ and $C^r(x',i',k')$ are vertex-disjoint, and
        \item\label{itm:C61hit} both $C^l(x,i,k)$ and $C^r(x,i,k)$ are vertex-disjoint from all cubes of absorbing $\ell$-cube pairs in $\cC^\mathrm{sc}_1$.
    \end{enumerate}
    \item \label{itm:C7hit} For every $x\in \cU$ and every $z\in \mathrm{End}(x)$, there exists a collection of cubes $\cC''(z)$ in $I(Q')$ with $|\cC''(z)|= 20\ell$ which satisfies the following properties:
    \begin{enumerate}[label=$(\mathrm{C}7.\arabic*)$]
        \item\label{itm:C70hit} for every $C\in\cC''(z)$, we have that $z\in V(C)$;
        \item\label{itm:C72hit} for every $C\in \cC''(z)$, the molecule $\cM_{C}$ is bonded in $Q'$;
        \item\label{itm:C73hit} every $C\in \cC''(z)$ meets $T''$;
        \item\label{itm:C74hit} for every $C_1,C_2\in \cC''(z)$, we have $V(C_1)\cap V(C_2)=\{z\}$, and
        \item\label{itm:C71hit} for every $C^{*}\in \cC^\mathrm{sc}_1$ and $C\in \cC''(z)$, we have $V(C)\cap V(C^{*})= \varnothing$.
    \end{enumerate}
    \item\label{itm:C8hit} $C^\bullet$ intersects $V(T)\cap V(T'')$ and is different from all cubes described in \ref{itm:C3hit}, \ref{itm:C6hit} and \ref{itm:C7hit}.
\end{enumerate}
Indeed, \ref{itm:C2hit} is given by \eqref{equa:T''boundhit}.
\ref{itm:C1hit} holds by \ref{itm:ET3hit} and the fact that $P'\cap V(T'')=\varnothing$.\COMMENT{Note here it is crucial that all vertices $z$ in $T''\cap R'$ are crucially still in fact joined to a vertex in $Z(z)$.
If a vertex $z \in T' \cap R$ was affected by the repatching step, then recall that the connecting paths $P$ are now part of the new tree $T''$ so we do not need to consider these vertices in this sense anymore.}
\ref{itm:C4hit} follows directly by \eqref{equa:T''containedhit}.
\ref{itm:C52hit} follows by combining \ref{itm:nib1}, the conditioning on $\mathcal{E}_5^*$, \eqref{eqn:moveChit} and \ref{itm:bad3}.
(Indeed, for every $x \in V(I)$ and every $C \in \cC_x$ we have by \ref{itm:bad3} that $C$ contains a vertex that lies in $B^{2\ell}_I(x) \setminus B^{5}_I(\cU_I)$.
By \eqref{eqn:moveChit}, at most $2n^{3/4}$ vertices lying in $B^{2\ell}_I(x) \setminus B^{5}_I(\cU_I)$ do not lie in $T''$.
Therefore, by combining \ref{itm:nib1}, the conditioning on $\mathcal{E}_5^*$, and \eqref{eqn:moveChit} we see that at least $(1-\delta)n-n/\ell^4-2n^{3/4}\geq(1-3/2\ell^4)n$ of the neighbours of $x$ in $I$ lie in some bonded cube $C\in\cC_x\cap\cC'$ which intersects $T''$.)
\ref{itm:C3hit} follows from the construction of $P$ and $T''$ in Step~9.
Indeed, \ref{itm:C30hit} follows from the definition of $Q'$ combined with \ref{itm:abss2hit}, and \ref{itm:C31hit} holds by the definition of $R'$ and $T''$ combined with \ref{itm:abss2hit}, while \ref{itm:C325hit} follows because of the definition of the set $M(A_i(x))$ in \ref{itm:abss1} and \ref{itm:AB3}, and \ref{itm:C33hit} holds by \ref{itm:DisjAbsCubeshit}.
Consider now \ref{itm:C6hit}.
For each $x \in V_\mathrm{abs}'\setminus V_{\mathrm{sc}}$ and $i \in [K]$, consider $\cC^\mathrm{abs}_1(x,i)$.
All absorbing $\ell$-cube pairs of $\cC^\mathrm{abs}_1(x,i)$ satisfy \ref{itm:C60hit} and \ref{itm:C62hit} by the definition of $Q'$, \ref{itm:abss2} and \ref{itm:AB4hit}.
Similarly, they satisfy \ref{itm:C64hit} by \ref{itm:abss1}, \ref{itm:AB3} and the fact that, by \ref{itm:AB4hit}, their intersection with $T''$ contains their intersection with $T'$.
Moreover, \ref{itm:C65hit} holds by \cref{claim:rainbowapplnewhit}, and \ref{itm:C61hit} holds because of \ref{itm:AB4hit}.
Now, \ref{itm:C7hit} holds by \cref{claim:cubesinendpointsofCSP} and \ref{itm: bad and scant}.
Indeed, let $\cC''(z)$ be the collection of cubes given by \cref{claim:cubesinendpointsofCSP}, so \ref{itm:C70hit}, \ref{itm:C72hit} and \ref{itm:C74hit} follow directly.
\ref{itm:C73hit} follows by using again the observation that, by \ref{itm: bad and scant} and the construction of $P$, for any $x\in\cU_I$, we have that $T'$ and $T''$ coincide in $B_I^{2\ell}(x)$.
Now recall that, by \ref{itm: bad and scant}, all vertices $x\in\cU_I$ are at distance at least $5\ell$ from $V_{\mathrm{sc}}$, so \ref{itm:C71hit} follows by construction.
Finally, consider \ref{itm:C8hit}.
The fact that $C^\bullet$ intersects $V(T)\cap V(T'')$ follows by its definition in Step~9, as does the fact that it is different from all cubes described in \ref{itm:C3hit}.
The fact that it is different from all cubes in \ref{itm:C6hit} follows by \ref{itm:AB4hit}.
Finally, the fact that it is different from the cubes in \ref{itm:C7hit} follows since $\dist(V(C^\bullet),\cU_I)\geq\ell$ by the definition of $C^\bullet$ in Step~9.\\

\textbf{Step~11: Fixing special absorbing structures.}
From this point onward, every step will be deterministic.
Let $F\subseteq\cQ^n$ be any graph with $\Delta(F)\leq\Psi$ which is $(\cU, \ell, s)$-good, that is, for each $x\in\cU$, the set $E_F(x)\coloneqq\{e\in E(F):e\cap N_{\cQ^n}(x)\neq\varnothing\}$ satisfies the following:
\begin{enumerate}[label=($F^*$)]
    \item\label{itm:Fdirec} for each layer $L$ of $\cQ^n$ and all $d\in\mathcal{D}(L)$, we have $|\{e\in E_F(x):\mathcal{D}(e)=d\}|\leq n/\ell$.\COMMENT{Alberto: Actually, we could also restrict $E_F$ to layers, but I think we don't need to (meaning that this condition follows for free for these vertices, I think, since all of them are absorbed, and at most two of them are absorbed using the same directions).}\COMMENT{Note that this is the same as the condition written in the statement of the theorem we are proving.}
\end{enumerate}
Let $H'\subseteq \cQ^n$ be any graph such that, for every $x\in\cU$, we have $d_{H'}(x)\geq2$.
For each $x \in \cU$, let $\{x,x+a(x)\},\{x, x+b(x)\}\in E(H')$, where $a(x),b(x)\in \cD(\cQ^{n})$.  
Our goal is to find a $(\cU, \ell^2, s)$-good Hamilton cycle in $((H\cup G)\setminus F)\cup H'\cup Q'$ which, for each $x\in\cU$, contains the edges $\{x,x+a(x)\}$ and $\{x, x+b(x)\}$.
Recall that $\mathfrak{C}'(x,a(x),b(x))$ was defined in Step~8.

\begin{claim}\label{claim:10.4}
For every $x\in \cU$, there exists an $(x,a(x),b(x))$-consistent system of paths $\mathit{CS}(x)\in \mathfrak{C}'(x,a(x),b(x))$\index{CS@$\mathit{CS}(x)$} such that $(E(\mathit{CS}(x))\setminus\{\{x,x+a(x)\},\{x,x+b(x)\}\})\cap E(F)=\varnothing$. 
\end{claim}

\begin{claimproof}
Suppose $x \in V(L)$, for some layer $L$.
Suppose first that $x+a(x), x+b(x) \in V(L)$.
Thus, we must show the existence of an $(x,a(x),b(x))$-consistent system of paths of Type~I in $\mathfrak{C}'(x,a(x),b(x))$ with the desired property.
Recall all the notation for consistent systems of paths introduced in \cref{sect:absstruct}, as well as \cref{def:rob}.
By \ref{itm:conssys}, there is a collection $\cD^{(2)}(x, a(x), b(x))$ of at least $\gamma n/2$ disjoint pairs of distinct directions $c,d \in \cD(L)\setminus\{a(x),b(x)\}$ such that, for each $(c,d)\in\cD^{(2)}(x, a(x), b(x))$, there is a collection $\cD^{(4)}(x,a(x),b(x),c,d)$ of at least $\gamma n/2$ disjoint $4$-tuples of distinct directions in $\mathcal{D}(L)\setminus\{a(x),b(x),c,d\}$ satisfying the following: for each $(c,d)\in\cD^{(2)}(x, a(x), b(x))$ and each $(d_1,d_2,d_3,d_4)\in\cD^{(4)}(x,a(x),b(x),c,d)$, $\mathfrak{C}'(x,a(x),b(x))$ contains the $(x,a(x),b(x))$-consistent system of paths $\mathit{CS}(c,d,d_1,d_2,d_3,d_4)$ defined as in \cref{sect:absstruct}.
We will now show that there are many such consistent systems of paths which avoid $F$.

The choice of $(c,d)\in\cD^{(2)}(x,a(x),b(x))$ determines six edges of the consistent system of paths: $e_1\coloneqq\{f(x+a(x)), f(x+a(x)+d)\}$, $e_2\coloneqq \{f(x+b(x)), f(x+b(x)+c)\}$, $e_3\coloneqq \{f(x+c), f(x)\}$, $e_4\coloneqq \{f(x), f(x+d)\}$, $e_5\coloneqq\{x+c, x+c+b(x)\}$ and $e_6\coloneqq\{x+d, x+d+a(x)\}$.
Since $f(x)$, $f(x+a(x))$ and $f(x+b(x))$ are fixed and $\Delta(F)\le \Psi$, for each $i \in [4]$ there are at most $\Psi$ choices of $(c,d)\in\cD^{(2)}(x,a(x),b(x))$ such that $e_i \in E(F)$.
Furthermore, by \ref{itm:Fdirec}, for each $i \in \{5,6\}$ there are at most $n/\ell$ choices $(c,d)\in \cD^{(2)}(x, a(x), b(x))$ such that $e_i \in E(F)$.
Thus, there exist at least $\gamma n/2 - 4 \Psi - 2n/\ell \geq \gamma n/4$ choices $(c,d)\in \cD^{(2)}(x, a(x), b(x))$ such that $e_i \in E(H)\setminus E(F)$ for all $i \in [6]$.
For any such choice of $(c,d)$, the choice of $(d_1,d_2,d_3,d_4)\in\cD^{(4)}(x,a(x),b(x),c,d)$ now determines the remaining eight edges of an $(x,a(x),b(x))$-consistent system of paths, each with a unique endpoint in $\{x+ a(x), x+b(x), f(x+b(x)), x+c, f(x+c), f(x+d), x+d, f(x+a(x))\}$.
It now follows by the fact that $\Delta(F) \le \Psi$ that there are at most $8\Psi$ choices of $(d_1,d_2,d_3,d_4)\in\cD^{(4)}(x,a(x),b(x),c,d)$ such that some of these eight edges lies in $E(F)$.
In particular, we may fix a consistent system of paths $\mathit{CS}(x)\in \mathfrak{C}'(x,a(x),b(x))$ which satisfies the statement of the claim.

The cases where $x+a(x) \notin L$, $x+b(x) \in L$ and where $x+a(x), x+b(x) \notin L$ can be shown similarly.
\end{claimproof}

Note that $\mathit{CS}(x)\subseteq(H\setminus F)\cup H'$ for each $x\in\cU$.

\begin{claim}\label{claim:lastabs}
For every $x\in \cU$, we can extend $\mathit{CS}(x)$ into an $(x,a(x),b(x))$-special absorbing structure $\mathit{SA}(x)$\index{AS@$\mathit{AS}(x)$} such that the following hold:
\begin{enumerate}[label=$(\mathrm{SA}_{\mathrm{\roman*}})$]
    \item\label{itm:finalabs1} for every $C\in \mathbf{C}(\mathit{SA}(x))$, we have that $\cM_C\subseteq Q'$ and $\cM_C$ is bonded in $Q'$, and
    \item\label{itm:finalabs3} every $C\in \mathbf{C}(\mathit{SA}(x))$ meets $T''$.
\end{enumerate}
\end{claim}

\begin{claimproof}
For each $x \in \cU$, we iterate through each $z\in \mathrm{endmol}(\mathit{CS}(x))$ fixing a cube $C(z) \in \cC''(z)$.
This will then determine $\mathbf{C}(\mathit{SA}(x))$, by taking the appropriate clones of $C(z)$.
To see that this can be done, note that $|\mathrm{endmol}(\mathit{CS}(x))|\le 6$.
For each $z\in \mathrm{endmol}(\mathit{CS}(x))$, by \ref{itm:C7hit}, there exist at least $20\ell$ choices of $C(z) \in \cC''(z)$ for which \ref{itm:finalabs1} and \ref{itm:finalabs3} hold.
Finally, by \ref{itm:C74hit} it follows that we can fix $C(z) \in \cC''(z)$ such that $\cD(C(z))\cap\cD(\mathit{CS}(x))=\varnothing$\COMMENT{Since $|\cD(\mathit{CS})|\leq8$ and, by \ref{itm:C74hit}, all cubes in $\cC''(z)$ use disjoint sets of directions, this means we may remove at most eight cubes from $\cC''(z)$ and guarantee that all remaining cubes satisfy this property.}.
In particular, this implies that $C(z)$ is vertex-disjoint from all $C(z')$ already fixed with $z\neq z' \in \mathrm{endmol}(\mathit{CS}(x))$\COMMENT{This follows because of the geometry of the cube. 
Indeed, observe that any two vertices of $\mathrm{endmol}(\mathit{CS}(x))$ are at distance at most $4$ from each other, and their differing directions are some of the directions of $\mathit{CS}$.
Therefore, if the cubes we put on each vertex do not contain any of these directions, they cannot intersect each other.} and, therefore, this process forms a valid extension of $\mathit{CS}(x)$ into an $(x,a(x),b(x))$-special absorbing structure.
\end{claimproof}

For each $x\in\cU$, let $\mathit{SA}(x)$ be an $(x,a(x),b(x))$-special absorbing structure which extends $\mathit{CS}(x)$, as determined by \cref{claim:lastabs}.
Note that, by \ref{itm:bad3} and the fact that $V(\mathit{SA}(x))\subseteq B_{\cQ^{n}}^{2\ell}(x)$, the special absorbing structures in the collection $\{\mathit{SA}(x):x\in\cU\}$ are pairwise vertex-disjoint.
Denote by $\mathit{SA}^{v}\coloneqq \bigcup_{x\in \cU} V(\mathit{SA}(x))$\index{SAv@$\mathit{SA}^{v}$}.
Recall that, for any $C\in\mathbf{C}(\mathit{SA}(x))$, $C_I\subseteq I$ denotes the cube which $C$ is a clone of.
Given any tree $T^\bullet\subseteq I$ and any $x\in\cU$, we say that $\mathit{SA}(x)$ \emph{meets $T^\bullet$} if, for all $C\in \mathbf{C}(\mathit{SA}(x))$, we have $V(C_I)\cap V(T^\bullet) \neq \varnothing$.

Recall that $\cC'$ and $\cC''$ were defined in Step~6.
Let $\cC_1^{*}\coloneqq\{C\in \cC': V(\cM_{C})\cap \mathit{SA}^{v}\neq \varnothing\}$. 
Note that, by \ref{itm:bad3},
\begin{enumerate}[label=$(\mathrm{CB})$]
\item\label{itm:CB} for each $x\in\cU$, there are at most $2^{2\ell}$ $\ell$-cubes $C \in \cC^*_1$ such that $V(C) \cap B^{10\ell}(x) \neq \varnothing$.\COMMENT{This follows because special absorbing structures are very far apart and each such special absorbing structure contains at most $2^{2\ell}$ vertices.}
\end{enumerate}
Let $\cC'_{1}\coloneqq \cC' \setminus \cC_1^{*}$\index{Cw7@$\cC'_{1}$} and   $\cC''_1=\cC''\setminus \cC^*_1$\index{Cw8@$\cC''_{1}$}.
We now define a tree $T'''\subseteq T''$ in the following way.
Consider each $x\in R'\cap V(T'')$ such that $x\in V(C)$ for some $C\in\cC'_1$\COMMENT{This extra condition that it lies in some cube that has not been deleted is here to guarantee that the cubes for the special absorbing structures meet $T'''$. 
Indeed, observe that for all vertices which do not satisfy this condition we have that they either lie in one of the special absorbing structures' cubes, or do not lie in any cube (and therefore we don't mind if they are retained as edges of the tree).}.
By \ref{itm:C1hit}, we have that $x$ has a unique neighbour $x'$ in $T''$, and $x' \in Z(x)$.
By the definition of $Z(x)$ (see Step~7), it follows that $x' \in V(C')$ for some $C' \in \cC''$.
If $C' \notin \cC''_1$, then we remove $x$ from $T''$.
We denote the resulting tree by $T'''$\index{T4@$T'''$}.
Let $\cC'''_{1}\coloneqq \{C\in\cC'_1:V(C)\cap V(T''')\neq\varnothing\}$\index{Cw9@$\cC'''_{1}$}.
By using \ref{itm:C2hit}--\ref{itm:C6hit}, the definition of $\cC'_1$, $\cC''_1$ and $\cC'''_1$, the construction of $T'''$, and the maximum degree of $F$, we claim that the following now hold.
\begin{enumerate}[label=$(\mathrm{C'}\arabic*)$]
    \item\label{itm:C2primehit} $\Delta(T''')\leq 12D$.
    \item\label{itm:C1primehit} Any vertex $x\in R'\cap V(T''')$ is a leaf of $T'''$.
    Furthermore, if $x\in R'\cap V(T''')$, then $x\notin V(T)$ and its unique neighbour $x'$ in $T'''$ satisfies that $x'\in N_{I}(x)\cap V(T) \cap (\bigcup_{C\in\cC''_1}V(C))\subseteq Z(x)$.
    \item\label{itm:C4primehit} $B_{I}^{5}(\cU_{I})\cap V(T''')=\varnothing$.
    \item\label{itm:C5primebis2hit} For all $x\in V(I)$ we have that $|\cC_x \cap \cC'''_1| \geq (1-2/\ell^4)n$.
    \item \label{itm:C3primehit} For each $x \in V_{\mathrm{sc}}$ and $i \in [K]$, there is an absorbing $\ell$-cube pair $(C^l(x,i),C^r(x,i))$ for $x$ in $I$, which is associated with some edge $e \in M(A_i(x))$.
    In particular, $(C^l(x,i),C^r(x,i))$ is such that there are two absorbing $\ell$-cube pairs $(C^l_1(x,i),C^r_1(x,i))$ and $(C^l_2(x,i),C^r_2(x,i))$ in $H\cup G$, associated with $e\in M(A_i(x))$, for the clones $x_1$ and $x_2$ of $x$ which correspond to $(x,i)$.
    Additionally, each of these absorbing $\ell$-cube pairs $(C^l(x,i),C^r(x,i))$ satisfies the following:
    \begin{enumerate}[label=$(\mathrm{C'}5.\arabic*)$]
        \item \label{itm:C35primehit} $(C^l_1(x,i), C^r_1(x,i)) \cup (C^l_2(x,i), C^r_2(x,i)) - V(\cM_x) \subseteq Q'$;
        \item\label{itm:C31primehit} the tips $x^{l}$ of $C^l(x,i)$ and $x^{r}$ of $C^r(x,i)$ lie in $R'\setminus V(T''')$, and $\{x,x^{l}\},\{x,x^{r}\}\notin E(F_I)$; in particular, the tips $x_1^l,x_1^r$ of $(C^l_1(x,i),C^r_1(x,i))$ and $x_2^l,x_2^r$ of $(C^l_2(x,i),C^r_2(x,i))$ satisfy that $\{x_1,x_1^l\},\{x_1,x_1^r\},\{x_2,x_2^l\},\{x_2,x_2^r\}\in E((H\cup G)\setminus F)$;
        \item\label{itm:C32primehit} $C^l(x,i),C^r(x,i)\in\mathcal{C}''_1\cap \cC'''_1$, and
        \item\label{itm:C33primehit} for any $x' \in V_{\mathrm{sc}}$ and $i' \in [K]$ with $(x',i')\neq(x,i)$, we have that $C^l(x,i)$, $C^r(x,i)$, $C^l(x',i')$ and $C^r(x',i')$ are vertex-disjoint.
    \end{enumerate}
    Let $\cC^\mathrm{sc}$ denote the collection of these absorbing $\ell$-cube pairs.
    \item\label{itm:C6primehit} For each $x \in V_\mathrm{abs}'\setminus V_{\mathrm{sc}}$ and $i \in [K]$, there is an absorbing $\ell$-cube pair $(C^l(x,i),C^r(x,i))$ for $x$ in $I$, which is associated with an edge $e\in M(A_i(x))$.
    In particular, $(C^l(x,i),C^r(x,i))$ is such that there are two absorbing $\ell$-cube pairs $(C^l_1(x,i),C^r_1(x,i))$ and $(C^l_2(x,i),C^r_2(x,i))$ in $H\cup G$,  associated with $e\in M(A_i(x))$, for the clones $x_1$ and $x_2$ of $x$ which correspond to $(x,i)$.
    Moreover, each of these absorbing $\ell$-cube pairs $(C^l(x,i),C^r(x,i))$ satisfies the following:
    \begin{enumerate}[label=$(\mathrm{C'}6.\arabic*)$]
        \item\label{itm:C66primehit} $(C^l_1(x,i), C^r_1(x,i)) \cup (C^l_2(x,i), C^r_2(x,i)) - V(\cM_x) \subseteq Q'$;
        \item\label{itm:C62primehit} the tips $x^{l}$ of $C^l(x,i)$ and $x^{r}$ of $C^r(x,i)$ lie in $R'$, and $\{x,x^{l}\},\{x,x^{r}\}\notin E(F_I)$; in particular, the tips $x_1^l,x_1^r$ of $(C^l_1(x,i),C^r_1(x,i))$ and $x_2^l,x_2^r$ of $(C^l_2(x,i),C^r_2(x,i))$ satisfy that $\{x_1,x_1^l\},\{x_1,x_1^r\},\{x_2,x_2^l\},\{x_2,x_2^r\}\in E((H\cup G)\setminus F)$;
        \item\label{itm:C64primehit} $C^l(x,i),C^r(x,i)\in\mathcal{C}''_1\cap\cC'''_1$;
        \item\label{itm:C65primehit} for any $x' \in V_\mathrm{abs}'\setminus V_{\mathrm{sc}}$ and $i' \in [K]$ with $(x',i')\neq(x,i)$, we have that $C^l(x,i)$, $C^r(x,i)$, $C^l(x',i')$ and $C^r(x',i')$ are vertex-disjoint, and
        \item\label{itm:C61primehit} both $C^l(x,i)$ and $C^r(x,i)$ are vertex-disjoint from all cubes of absorbing $\ell$-cube pairs in $\cC^\mathrm{sc}$.
    \end{enumerate}
    Let $\cC^{\neg\mathrm{sc}}$ denote the set of these absorbing $\ell$-cube pairs.
    \item \label{itm:C7primehit} For every $x\in \cU$, there is an $(x,a(x),b(x))$-consistent system of paths $\mathit{CS}(x)$ in $(H\setminus F)\cup H'$ which extends into an $(x,a(x),b(x))$-special absorbing structure $\mathit{SA}(x)$ which meets $T'''$ and with the property that, for every $C\in \mathbf{C}(SA(x,a,b))$, we have that $\cM_{C_I}\subseteq Q'$ and $\cM_{C_I}$ is bonded in $Q'$.
    Moreover, $\{x, x+a(x)\}, \{x, x+b(x)\} \in E(H')$.
    \item\label{itm:C8primehit} $C^\bullet$ intersects $V(T)\cap V(T''')$ (so, in particular, $C^\bullet \in \cC'''_1$) and is different from all cubes described in \ref{itm:C3primehit}, \ref{itm:C6primehit} and \ref{itm:C7primehit}.
    \item\label{itm:C9primehit} $\bigcup_{C\in\cC'''\setminus\cC'''_1}V(C)\subseteq B_I^{3\ell}(\cU_I)$.
\end{enumerate}
Indeed, since $T'''\subseteq T''$, \ref{itm:C2primehit}--\ref{itm:C4primehit} follow immediately by \ref{itm:C2hit}--\ref{itm:C4hit}, respectively.
\ref{itm:C5primebis2hit} follows from \ref{itm:C52hit}, \ref{itm:bad3} and \ref{itm:CB}.
Now, for each $x\in V_\mathrm{sc}$ and $i\in[K]$, consider the set $\cC^\mathrm{sc}_1(x,i)$ described in \ref{itm:C3hit}.
We first remove from this set all absorbing $\ell$-cube pairs any of whose cubes do not belong to $\cC'_1$.
Then, we remove all absorbing $\ell$-cube pairs any of whose cubes do not intersect $T'''$.
Finally, we remove all absorbing $\ell$-cube pairs such that any of the edges joining its tips to $x$ belong to $F_I$.
Observe that by \ref{itm:CB} and \ref{itm:C2primehit} it follows that, for any $y \in V(I)$, we have $|B_I^{10\ell}(y) \cap (V(T'') \setminus V(T'''))| \leq 12D \cdot 2^{3\ell}$.\COMMENT{The reason is as follows. 
By \ref{itm:CB} we removed at most $2^{2\ell}$ cubes around a special absorbing structure.
Then in the definition of $T'''$ we removed edges from $T''$ which were leaves into the reservoir joined to these cubes.
By construction, every vertex in $T''$ could have at most one of these leaves (but since we didn't write this explicitly, we use the upper bound for the degree of $T'''$ here, which is $12D$) and so we removed at most $12D 2^\ell 2^{2\ell}= 12D2^{3\ell}$ leaves.}
Using this fact, \ref{itm:CB}, and the fact that $\Delta(F) \le \Psi$ (and, thus, $\Delta(F_I)\leq2^s\Psi$), it follows that there is at least one absorbing $\ell$-cube pair remaining in the collection\COMMENT{By \ref{itm:CB} and \ref{itm:bad3}, we remove at most $2^{2\ell}$ absorbing $\ell$-cube pairs in the first deletion.
The second deletion, using the fact about leaves mentioned above, removes at most $12D 2^{3\ell}$ more pairs (indeed, the worst case scenario is given when each of those leaves belongs to a cube of a different absorbing $\ell$-cube pair and these cubes are not connected to the tree $T''$ by any other edge).
Finally, at most $2^s\Psi$ absorbing $\ell$-cube pairs are removed because of $F$.
Therefore, at least $2^{s+1}\Psi-2^s\Psi-12D2^{3\ell}-2^{2\ell}=2^{10\ell}\Psi-12D2^{3\ell}-2^{2\ell}\geq1$ absorbing $\ell$-cube pairs remain.}.
Let $(C^l(x,i),C^r(x,i))$ be such an absorbing $\ell$-cube pair.
Then, \ref{itm:C31primehit} and \ref{itm:C32primehit} hold by the choice above, and \ref{itm:C35primehit} and \ref{itm:C33primehit} hold by \ref{itm:C30hit} and \ref{itm:C33hit}, respectively.
For each $x\in V_\mathrm{abs}'\setminus V_\mathrm{sc}$ and $i\in[K]$, we proceed similarly from the set $\cC^\mathrm{abs}_1(x,i)$ to fix an absorbing $\ell$-cube pair $(C^l(x,i),C^r(x,i))\in\cC^\mathrm{abs}_1(x,i)$ which satisfies \ref{itm:C62primehit} and \ref{itm:C64primehit}.
Then, \ref{itm:C66primehit}, \ref{itm:C65primehit} and \ref{itm:C61primehit} hold by \ref{itm:C60hit}, \ref{itm:C65hit} and \ref{itm:C61hit}, respectively.
Furthermore, \ref{itm:C7primehit} holds by \cref{claim:10.4}, \cref{claim:lastabs} and the construction of $T'''$ above\COMMENT{Here is where we use the condition that we only removed edges of $T''$ for vertices $x\in R'\cap V(T'')$ such that $x\in V(C)$ for some $C\in\cC'_1$, which guarantees that all cubes of the special absorbing structures meet $T'''$.}.
(Indeed, to see that each $SA(x)$ still meets $T'''$, note that
$(V(T'') \setminus V(T''')) \cap \mathit{SA}^v = \varnothing$.)
For \ref{itm:C8primehit}, by construction $C^\bullet \in \cC'$ and $V(C^\bullet) \cap B^{3\ell/2}_I(\cU_I) = \varnothing$.
Therefore, $C^\bullet \in \cC'_1$.
The fact that $C^\bullet$ intersects $V(T)\cap V(T''')$ follows by \ref{itm:C8hit} and the fact that, in constructing $T'''$, none of the leaves which are removed are vertices of $T$.
Thus, in particular, $C^\bullet\in\cC'''_1$.
The rest of  \ref{itm:C8primehit} follows immediately from \ref{itm:C8hit}.
Finally, \ref{itm:C9primehit} follows by the definition of $\cC^*_1$.
Indeed, consider the set $\mathit{SA}^v_I\subseteq V(I)$ of vertices such that each vertex in $\mathit{SA}^v$ is a clone of some vertex in $\mathit{SA}^v_I$.
It follows by construction (see \cref{sect:absstruct}) that for any $x\in\mathit{SA}^v_I$ we have $\dist(x,\cU_I)\leq\ell+2$. 
The claim follows since any $\ell$-cube $C\in\cC^*_1$ must intersect $\mathit{SA}^v_I$ and any two vertices in $C$ are at distance at most $\ell$.

Let $\cC'_2\coloneqq \bigcup_{x\in\cU}\{C_{I}: C\in \mathbf{C}(\mathit{SA}(x))\}$\index{Cw99@$\cC'_{2}$} and $R''\coloneqq R'\setminus\bigcup_{C\in\cC'_2}V(C)$\index{R3@$R''$}.
Finally, let $\cC'''_3\coloneqq \cC'''_1\cup \cC'_2$\index{Cw999@$\cC'''_3$}.\COMMENT{This collection contains precisely the cubes needed for absorbing the low degree vertices.}
Note that, by construction, any two cubes in $\cC'''_3$ are vertex-disjoint.\COMMENT{$\cC'''_1$ was constructed before the above list of properties to contain the remaining nibble cubes which avoided all vertices of all fixed special absorbing structures. 
$\cC'_2$ contained precisely those cubes which were part of the fixed special absorbing structures (and were disjoint from each other by construction).
Therefore by definition $\cC'''$ consists of vertex-disjoint cubes.}\\

\textbf{Step~12: Constructing auxiliary trees $\boldsymbol{T^*}$ and $\boldsymbol{\tau_0}$.}
Let $T^*$ be obtained from $T'''$ by removing all leaves of $T'''$ which lie in $R''$\COMMENT{By doing so, we retain all cubes for special absorbing structures.}.
In particular, by \ref{itm:C1primehit} and \ref{itm:C8primehit}, we have that $C^\bullet$ intersects $T^*$.

We now construct an auxiliary tree $\tau_0$, which will be used in the construction of an almost spanning cycle.
The construction of $\tau_0$\index{tau0@$\tau_0$} is identical to that in Step~10 of the proof of \cref{thm:main1}, except that $\cC'''_3$ plays the role of $\cC'$ in the definition of $\Gamma_1\coloneqq T^*\cup\bigcup_{C\in\cC'''_3}C$ and the subsequent steps, and that, for the depth-first search on $\Gamma'$, the root vertex $v_0\in V(\Gamma')$ is chosen to be the vertex which resulted from contracting $C^\bullet$ (in particular, as in the proof of \cref{thm:main1}, $v_0$ is an atomic vertex).

Let $m\coloneqq|V(\tau_0)|-1$.
We define $v_1,\ldots,v_m$, \index{Mofv@$\mathcal{M}(v)$}$\cM(v_i)$, $\cA_j(v_i)$, $\tau_i$, $p_i$, $u_1^i,\ldots,u_{p_i}^i$, $e_k^i$, $f_k^i$, $j_k^i$, $\nu_k^i$, $\Delta(v_i)$ and $b(i)$ analogously to Step~10 of the proof of \cref{thm:main1}.
In particular, we again have that
\begin{equation}\label{equa:step9.2hit}
\begin{minipage}[c]{0.7\textwidth}
$p_i \leq 12D-1$ if $v_i$ is an inner tree vertex\COMMENT{We need the number of layers in each slice to be larger than this (say, at least $10$ times).}, and $\Delta(\tau_0)\leq12\cdot2^\ell D$\COMMENT{We need the number of slices in each molecule to be larger than this (say, at least $10$ times).}.
\end{minipage}\ignorespacesafterend 
\end{equation}

\textbf{Step~13: Finding an external skeleton for $\boldsymbol{T^*}$.}
We now generate an external skeleton, following Step~11 of the proof of \cref{thm:main1}.
Using this external skeleton, we will construct a first skeleton in Step~16 and then extend it in Step~17 by incorporating the special absorbing structures for the vertices in $\cU$.

Let $\mathcal{M}_{\mathrm{Res}}\subseteq V(\cQ^n)$\index{MRes@$\mathcal{M}_{\mathrm{Res}}$} be the union of all the clones of $R''$.
For each $x\in\cU$, consider the graph $\mathit{CS}(x)_I\subseteq I$, and let $\mathcal{M}_{\mathit{CS}}\subseteq \cQ^n$ be the union of all the clones of $\bigcup_{x\in\cU}\mathit{CS}(x)_I$.
We construct an external skeleton $L^\bullet$ which satisfies properties \ref{item:preskprop1}--\ref{item:preskprop4} as in the proof of \cref{thm:main1} and the following variant of \ref{item:preskprop5}:
\begin{enumerate}[label=$(\mathrm{ES}\arabic*)$,start=5]
    \item\label{item:preskprop5hit} $L^\bullet\cap(\mathcal{M}_{\mathrm{Res}} \cup V(\cM_\mathit{CS}))=\varnothing$.
\end{enumerate}

The construction of $L^\bullet$ is identical to Step~11 of the proof of \cref{thm:main1}.
The new \ref{item:preskprop5hit} holds because of the definition of $\tau_0$.
Indeed, by \ref{itm:C1primehit} and \ref{itm:C4primehit} together with the definition of $R''$ and $T^*$, observe that $V(T^*)\cap (R''\cup \bigcup_{x \in \cU}V(\mathit{CS}(x)_I))=\varnothing$.
Moreover, by construction, all vertices in $L^\bullet$ are incident to some edge in a clone of the tree $T^*$, and thus they cannot lie in $\mathcal{M}_{\mathrm{Res}} \cup V(\cM_\mathit{CS})$.\\

\textbf{Step~14: Constructing an auxiliary tree $\boldsymbol{\tau_0'}$.}
We now extend $\tau_0$ to a new auxiliary tree $\tau_0'$\index{tau02@$\tau_0'$} which encodes information about all cube molecules which intersect $T'''$.
This is done as in Step~12 of the proof of \cref{thm:main1}, except that, again, $\cC_3'''$ plays the role of $\cC'''$ and $T'''$ plays the role of $T''$.
Then, the cubes represented in $\tau_0'$ are precisely all those in $\cC'''_3$.

Analogously to the proof of \cref{thm:main1}, it follows from \ref{itm:C2primehit} that
\begin{equation}\label{equa:step11.1hit}
\begin{minipage}[c]{0.7\textwidth}
$d_{\tau_0'}(v)\leq12D$ for all $v\in V(\tau_0')$ which are inner tree vertices, and $\Delta(\tau_0')\leq 12\cdot2^\ell D$.
\end{minipage}\ignorespacesafterend 
\end{equation} 
By \ref{itm:C32primehit} and \ref{itm:C64primehit}, we have that
\begin{enumerate}[label=$(\mathrm{CP})$]
\item\label{itm:AB5hit} every cube $C$ belonging to some absorbing $\ell$-cube pair in $\cC^\mathrm{sc}\cup\cC^{\neg\mathrm{sc}}$ is represented in $\tau_0'$.
\end{enumerate}
Finally, for each $x\in V(I)$, let $\zeta(x)$ denote the number of vertices $y\in N_I(x)$ which are represented in $\tau_0'$ by atomic vertices.
By \ref{itm:C5primebis2hit}\COMMENT{Any vertex that belongs to the intersection of the tree $T'''$ and a cube $C\in\cC'_1$ is represented in $\tau_0'$.}, we have that
\begin{equation}\label{equa:step10hit}
    \zeta(x)\geq(1-2/\ell^4)n.
\end{equation}

Let $m'\coloneqq|V(\tau_0')|-1$ and label $V(\tau_0')\setminus V(\tau_0)=\{v_{m+1},\ldots,v_{m'}\}$ arbitrarily.
We define $\tau_i'$, $p_i'$, $u_1^i,\ldots,u_{p_i'}^i$, $e_k^i$, $f_k^i$, $j_k^i$, $\nu_k^i$, $\Delta(v_i)$, $b(i)$, $\ell_i$ and $m_i$ as in Step~12 of the proof of \cref{thm:main1}.\\

\textbf{Step~15: Fixing absorbing $\boldsymbol{\ell}$-cube pairs for the vertices that need to be absorbed.}
We can now determine every vertex in $V(\cQ^n)$ that will have to be absorbed via absorbing $\ell$-cube pairs into the almost spanning cycle we are going to construct.
Recall from Step~11 that $\mathit{SA}^v=\bigcup_{y\in\cU}V(\mathit{SA}(y))$.
For every vertex $x \in V(I)$ not represented in $\tau_0'$, we will have to absorb all vertices in $\cM_x \setminus\mathit{SA}^v$.\COMMENT{Bondless and isolated molecules, as well as vertices not covered by either the tree or the molecules, and some vertices covered by the tree $T'''$ by which are deleted in the construction of $\tau_0$ (because they were leaves not pertaining to any cube).
Note here that all vertices in cubes of absorbing structures do not have to be absorbed, as they are represented in $\tau_0'$.}
Furthermore, for each $v\in V(\tau_0)$ which is an inner tree vertex, we will also need to absorb all vertices in $\mathcal{M}_v\setminus L^\bullet=\mathcal{M}_v\setminus (L^\bullet\cup\mathit{SA}^v)$.\COMMENT{We need to absorb every vertex from a bondless molecule and isolated molecule, as well of most of the vertices from an inner tree vertex molecule.}\COMMENT{Note that the tree avoided balls around $\cU$ so that no special absorbing vertices intersect with these molecules.}
(The fact that $\mathcal{M}_v\cap\mathit{SA}^v=\varnothing$ follows by \ref{itm:C4primehit}.)
Recall the definition of $V_{\mathrm{abs}}'$ from Step~10.
Let $V_{\mathrm{abs}}\subseteq V(I)$\index{Vabs@$V_\mathrm{abs}$} be the set of all vertices which are not represented in $\tau_0'$ by an atomic vertex.
Therefore, $V_{\mathrm{abs}}$  is the set of all vertices $x\in V(I)$ such that some clone of $x$ needs to be absorbed.
Moreover, $V_\mathrm{abs}=V(I)\setminus\bigcup_{C\in\cC'''_3}V(C)$ and, thus, \ref{itm:C9primehit} and the definition of $\cC_3'''$ at the end of Step~11 imply that $V_{\mathrm{abs}} \subseteq V_{\mathrm{abs}}'$.\COMMENT{As well as some extras. See definition of $\cN$ where $V_{\mathrm{abs}}'$ was defined. These vertices which lie near vertices in $\cU$ might be represented in $\tau_0'$ and therefore do not need absorbing.}
It follows from \eqref{equa:step10hit} that 
\begin{equation}\label{equa:step11hit}
    |V_{\mathrm{abs}}|\leq2^{n-s+1}/\ell^4.
\end{equation}
Now, for each $x\in V_{\mathrm{abs}}$, we will pair all those vertices in $\cM_x$ which need to be absorbed (each pair consisting of one vertex of each parity) and fix an absorbing $\ell$-cube pair for each such pair of vertices.
Recall that a difference between this pairing and the pairing in Step~13 of the proof of \cref{thm:main1} is that, in \cref{thm:main1}, we could guarantee that each pair was contained in one of the slices defined in Step~3.
Since a special absorbing structure might not lie in a single slice, we now cannot guarantee this anymore.
Instead, we can impose that each pair either lies in a slice or in two consecutive slices (with respect to their labelling).
The absorbing $\ell$-cube pair that we fix for each pair will be the one given by \ref{itm:C3primehit} or \ref{itm:C6primehit}, depending on whether $x\in V_\mathrm{sc}$ or not.

For each $x\in V_{\mathrm{abs}}$, let $S(x)\coloneqq V(\cM_x)\cap (L^\bullet \cup \mathit{SA}^v)=V(\cM_x)\cap (L^\bullet \cup V(\cM_\mathit{CS}))$.
It follows by \ref{item:preskprop1}--\ref{item:preskprop5hit}, \ref{itm:bad3} and the definition of our special absorbing structures that $|S(x)|\leq 25D$\COMMENT{$24D$ is enough. Indeed, this bound follows directly from the bound given by the external skeleton and the fact that the tree and the special absorbing structures are disjoint.} and $S(x)$ contains the same number of vertices of each parity\COMMENT{By definition of special structures they compensate in each molecule to ensure this. \ref{itm:bad3} is not really needed here, just emphasise that at most one special structure intersects any given molecule.}.
(Here we also use that $p_i\leq 12D-1$ for every inner tree vertex $v_i$ by \eqref{equa:step9.2hit} and \eqref{equa:step11.1hit}.)
Therefore, the matching $\mathfrak{M}(\cM_x,S(x))$ defined in Step~3 is well defined, and we can use this matching to define our pairing of the vertices in $V(\cM_x)\setminus S(x)$.
Recall that each edge $e\in\mathfrak{M}(\cM_x,S(x))$ gives rise to a unique index $i\in[K]$ via the relation $N(e)=A_i(x)$. 
(Here we ignore all those indices $i' \in [K]$ arising by artificially increasing the size of $\mathfrak{A}(x)$, see of Step~4 as well as Step~4 in the proof of \cref{thm:main1}.)
For each $x\in V_{\mathrm{abs}}$, let $\mathfrak{I}_{x}\subseteq[K]$ be the set of indices $i\in[K]$ which correspond to edges in $\mathfrak{M}(\cM_x,S(x))$\COMMENT{Basically, for any index $i\notin\mathfrak{I}_x$, we no longer care about $A_i(x)$ or $M(A_i(x))$.}.

For each $x\in V_{\mathrm{abs}}$ and $i\in\mathfrak{I}_x$, as stated in \ref{itm:C3primehit} and \ref{itm:C6primehit}, we have already fixed an absorbing $\ell$-cube pair for the clones of $x$ corresponding to $(x,i)$.
Let\index{Vabs2@$V^\mathrm{abs}$} 
\[V^\mathrm{abs}\coloneqq\bigcup_{x\in V_{\mathrm{abs}}}V(\mathcal{M}_x)\setminus (L^\bullet\cup\mathit{SA}^v).\]
(Thus, in particular, $V^\mathrm{abs}\cap\cU=\varnothing$.)
As discussed above, this is the set of all vertices that need to be absorbed via absorbing $\ell$-cube pairs.
Recall that $Q'$ was defined before \ref{itm:C2hit}--\ref{itm:C8hit}.
It follows from \ref{itm:C3primehit} and \ref{itm:C6primehit} that $((H\cup G)\setminus F) \cup Q'$ contains a set $\mathcal{C}^\mathrm{abs}=\{(C^l(u),C^r(u)):u\in V^\mathrm{abs}\}$ of absorbing $\ell$-cube pairs such that \ref{itm:RL1}--\ref{itm:RL5} in the proof of \cref{thm:main1} hold with $Q'$, $\cC''_1$ and $\cC_1'''$ playing the roles of $G'$, $\cC''$ and $\cC'''$, except that \ref{itm:RL2} is now replaced by the following:
\begin{enumerate}[label=$(\mathrm{C_{2.\arabic*}})$,start=2]
    \item\label{itm:RL2hit} if $f_i=\{u_i,u_i'\}$, then there is a vertex $v\in V_\mathrm{abs}$ such that $u_i$ and $u_i'$ are clones of $v$ which lie in either the same or consecutive slices of $\mathcal{M}_v$, and $(C^l(u_i),C^r(u_i))$ and $(C^l(u_i'),C^r(u_i'))$ are clones of the same absorbing $\ell$-cube pair for $v$ in $I$ such that $(C^l(u_i),C^r(u_i))$ lies in the same layer as $u_i$ and $(C^l(u_i'),C^r(u_i'))$ lies in the same layer as $u_i'$.
\end{enumerate}

We denote by $\mathfrak{L}$, $\mathfrak{R}_1$ and $\mathfrak{R}_2$ the collections of all left absorber tips, right absorber tips, and third absorber vertices, respectively, of the absorbing $\ell$-cube pairs in $\mathcal{C}^\mathrm{abs}$.
Observe that by \ref{itm:RL1}--\ref{itm:RL4} the following properties are satisfied:
\begin{enumerate}[label=$(\mathrm{C}^*\arabic*)$]
    \item\label{item:Cstarcon1hit} For all $i\in[m']_0$ such that $v_i$ is an atomic vertex, we have that $|\mathfrak{L}\cap V(\mathcal{M}(v_i))|\in\{0,2\}$. 
    If $|\mathfrak{L}\cap V(\mathcal{M}(v_i))|=2$, then these two vertices $u,u'$ lie in different atoms of either the same or consecutive slices of $\cM(v_i)$, and satisfy that $u\neqp u'$.
    \item\label{item:Cstarcon2hit} For all $i\in[m']_0$ such that $v_i$ is an atomic vertex, we have that  $|(\mathfrak{R}_1\cup\mathfrak{R}_2)\cap V(\mathcal{M}(v_i))|\in \{0,4\}$.
    If $|(\mathfrak{R}_1\cup\mathfrak{R}_2)\cap V(\mathcal{M}(v_i))|=4$, then these four vertices form two pairs such that one vertex of each pair belongs to $\mathfrak{R}_1$ and the other to $\mathfrak{R}_2$.
    Each of these pairs lies in a different atom of the same or consecutive slices of $\cM(v_i)$ and satisfies that its two vertices are adjacent in $Q'$.
    \item\label{item:Cstarcon3hit} For all $i\in[m']_0$ such that $v_i$ is an atomic vertex, if $\mathfrak{L}\cap V(\mathcal{M}(v_i))\neq\varnothing$, then $(\mathfrak{R}_1\cup\mathfrak{R}_2)\cap V(\mathcal{M}(v_i))=\varnothing$.
    \item\label{item:Cstarcon4hit} The sets described in \ref{item:Cstarcon1hit} and \ref{item:Cstarcon2hit} partition $\mathfrak{L}$ and $\mathfrak{R}_1\cup\mathfrak{R}_2$, respectively\COMMENT{that is, there are no vertices of $\mathfrak{L}\cup\mathfrak{R}_1\cup\mathfrak{R}_2$ outside the cube molecules represented in $\tau_0'$.}.
\end{enumerate}
Indeed, \ref{item:Cstarcon1hit}--\ref{item:Cstarcon3hit} follow from \ref{itm:RL1}--\ref{itm:RL4}, and \ref{item:Cstarcon4hit} follows by \ref{itm:AB5hit}. 

For each $u\in V^\mathrm{abs}$, we denote the edge consisting of the right absorber tip and the third absorber vertex of $(C^l(u),C^r(u))$ by $e_\mathrm{abs}(u)$, and we denote by $\mathcal{P}^\mathrm{abs}(u)$ the path of length three formed by the third absorber vertex, the left absorber tip, $u$, and the right absorber tip, visited in this order.
Note that $e_\mathrm{abs}(u)\in E(Q')$ by \ref{itm:RL1}.
Moreover, recall that $\cC^\mathrm{abs}$ consists of absorbing $\ell$-cube pairs in $((H\cup G)\setminus F)\cup Q')$.
Thus, $\mathcal{P}^\mathrm{abs}(u)\subseteq ((H\cup G)\setminus F)\cup Q'$.\\

\textbf{Step~16: Constructing the skeleton.}
As in Step~14 of the proof of \cref{thm:main1}, we now define the skeleton $\mathcal{L}=(x_1,\ldots,x_r)$ for the almost spanning cycle.
We again let $\mathcal{L}^\bullet\coloneqq\{x_1,\ldots,x_r\}$.
The skeleton $\mathcal{L}$ again satisfies \ref{itemskelprop1}--\ref{itemskelprop6} in the proof of \cref{thm:main1}, except that $Q'$ now plays the role of $G'$ and, \ref{itemskelprop3} and \ref{itemskelprop6} are now replaced by the following.
\begin{enumerate}[label=$(\mathrm{S}\arabic*)$]\setcounter{enumi}{2}
    \item\label{itemskelprop3hit} For every $k\in[r-1]$, if $x_k$ and $x_{k+1}$ do not both lie in the same slice of a cube molecule represented in $\tau_0'$, then $\{x_k,x_{k+1}\}\in E(Q')$.
    Moreover, in this case, if $x_{k+1}$ lies in a cube molecule represented in $\tau_0'$, then $x_{k+2}$ lies in the same slice of this cube molecule as $x_{k+1}$.\setcounter{enumi}{6}
    \item\label{itemskelprop6hit} $\mathcal{L}^\bullet\cap(\mathfrak{L}\cup\mathfrak{R}_1\cup V^\mathrm{abs} \cup V(\mathcal{M}_{\mathit{CS}})) =\varnothing$ and $L^\bullet\subseteq\mathcal{L}^\bullet$.
\end{enumerate}
The construction of $\cL$ is identical to that in \cref{thm:main1}. 
The only difference is that the `forbidden set $\cF$' which the skeleton has to avoid is replaced by  $\cF\coloneqq\mathfrak{L}\cup\mathfrak{R}_1\cup L^\bullet \cup V(\mathcal{M}_{\mathit{CS}})$ (this is required to ensure that \ref{itemskelprop6hit} holds).
For each $i\in[m']_0$ such that $v_i$ is an atomic vertex and each $j\in[t]$, let $\mathfrak{J}_{i,j}\coloneqq\{k\in[r]:x_k,x_{k+1}\in V(\mathcal{M}_j(v_i))\}$ and $S_{i,j}\coloneqq\{\{x_k,x_{k+1}\}:k\in\mathfrak{J}_{i,j}\}$.\\

\textbf{Step~17. Incorporating special absorbing structures into the skeleton.}
In this step, we are going to incorporate all special absorbing structures fixed in Step~11 into the skeleton we just constructed.
Note that all cube molecules referred to are represented by atomic vertices in $\tau'_0$, and all slices referred to are one of the $t$ slices of each of these molecules defined in Step~3.
For each $x \in \cU$, consider the consistent system of paths $\mathit{CS}(x)$ and the special absorbing structure $\mathit{SA}(x)$ given by \ref{itm:C7primehit}.
By \ref{itemskelprop6hit}, we have that $\cL$ avoids $\mathit{CS}(x)$.
Moreover, by the definition of $\cC'''_3$ at the end of Step~11, each $ C' \in \mathbf{C}(\mathit{SA}(x))$ is a clone of some $C \in \cC'''_3$.
Thus, by \ref{itemskelprop4} we have that $\cL$ has positive intersection with each slice which contains a vertex of $\mathrm{end}(\mathit{CS}(x))$.

Recall that a special absorbing structure is a tuple of paths and cubes (see \cref{sect:absstruct}).
For each $z \in \cU$, let $P^z_1$ denote the first path of $\mathit{SA}(z)$. 
Let $x(z) \in \cL^\bullet$ be the first vertex in $\cL$ that is contained in the slice which contains the first vertex of $P^z_1$.
Let $x'(z)$ be the successor of $x(z)$ in $\cL$ (in particular, by \ref{itemskelprop3hit}, both $x(z)$ and $x'(z)$ lie in the same slice).
Now, for each $z\in\cU$, depending on the type of the special absorbing structure $\mathit{SA}(z)$, we will update $\cL$ in different ways.
\begin{enumerate}[label=(\Roman*)]
    \item\label{caseI} If $\mathit{SA}(z)$ is a special absorbing structure of Type~I, proceed as follows.
    Let $P^z_1,\dots,P^z_6$ be the six paths of $\mathit{SA}(z)$.
    Let $S\coloneqq\bigtimes_{i=1}^6P^z_i$ and let $S^{-1}$ denote the same sequence of vertices in reverse order.
    If $x(z)$ has opposite parity to the initial vertex of $P^z_1$, then we replace the segment $(x(z), x'(z))$ of $\cL$ by $(x(z),S,x'(z))$; otherwise, we we replace the segment $(x(z), x'(z))$ by $(x(z),S^{-1},x'(z))$.\COMMENT{We use that the consistent systems start and end on vertices of opposite parity here and in what follows.
    This can be observed in their definition.
    For a more detailed proof of this see the comment in \ref{itemskelprop5hitbis} below.
    By imposing that the vertex we look at is the first one in $P^z_1$, we are guaranteed that all paths are ordered consistently for us, with the order given in \cref{sect:absstruct}.}
    \item  If $\mathit{SA}(z)$ is a special absorbing structure of Type~II, we proceed as follows.
    Let $P^z_1$ and $P^z_2$ be the two paths of $\mathit{SA}(z)$.
    Let $y^1$ and $x^1_0$ be the first and last vertices of $P^z_1$, and let $y^2$ and $x^2_0$ be the first and last vertices of $P^z_2$, respectively.
    Let $v,v'\in V(\tau_0')$ and $t_1,t_2,t_3\in[t]$ be such that $y^1 \in  V(\cM_{t_1}(v))$, $x^1_0 \in  V(\cM_{t_2}(v'))$, and $y^2 \in  V(\cM_{t_3}(v'))$ (this implies $x_0^2\in V(\cM_{t_3}(v))$).
    
    We now define two sequences of vertices $S^z_1$ and $S^z_2$ following similar ideas to Step~16.
    Recall that, for each $i\in[2^s]$, we use $\hat{e}_i$ to denote the direction of the edges between $L_i$ and $L_{i+1}$.
    If $t_3 \geq t_2$, let $m_1 \coloneqq t_3 - t_2$; otherwise, let $m_1 \coloneqq t - (t_2 - t_3)$.
    For each $k\in[m_1-1]_0$, iteratively choose a vertex $y^1_k\in V(\mathcal{A}_{(t_2+k)q}(v'))\setminus\cL^\bullet$ satisfying that\COMMENT{We have the definition of bondedness to thank for being able to do this again. 
    Note that in Step~16 there were at most $10$ vertices in any slice incorporated into the skeleton, and bondedness gives us $100$ edges of each parity between all atoms, so a suitable choice for such a sequence does exist.}
    \begin{enumerate}[label=$\arabic*$.]
        \item $y_k^1\neqp x_k^1$;
        \item $x_{k+1}^1\coloneqq y_k^1+\hat{e}_{(t_2+k)q}\notin \cL^\bullet$\COMMENT{We do not need to forbid the endpoints of the paths because we are guaranteed that they will be avoided by parity.}, and 
        \item $\{y_k^1,x_{k+1}^1\}\in E(Q')$.
    \end{enumerate}
    We set $S^z_1\coloneqq\bigtimes_{k=0}^{m_1-1}(y_k^1,x_{k+1}^1)$. 
    
    In order to construct $S_2^z$, we proceed similarly.
    If $t_3 \geq t_1$, let $m_2 \coloneqq t_3 - t_1$; otherwise, let $m_2 \coloneqq t - (t_1 - t_3)$.
    For each $k\in[m_2-1]_0$, iteratively choose a vertex $y^2_k\in V(\mathcal{A}_{(t_3-k-1)q+1}(v))\setminus\cL^\bullet$ satisfying that
    \begin{enumerate}[label=$\arabic*$.]
        \item $y_k^2\neqp x_k^2$;
        \item $x_{k+1}^2\coloneqq y_k^2+\hat{e}_{(t_3-k-1)q}\notin \cL^\bullet$, and
        \item $\{y_k^2,x_{k+1}^2\}\in E(Q')$.
    \end{enumerate}
    Now, let $S^z_2\coloneqq\bigtimes_{k=0}^{m_2-1}(y_k^2,x_{k+1}^2)$.
    
    Let $S\coloneqq P^z_1\times S^z_1\times P^z_2\times S^z_2$, and let $S^{-1}$ denote the same sequence of vertices in reverse order.
    Finally, we replace the segment $(x(z),x'(z))$ of $\cL$ by $(x(z),S,x'(z))$ if $x(z)$ has parity opposite to the initial vertex of $P^z_1$; otherwise, we replace $(x(z),x'(z))$ by $(x(z),S^{-1},x'(z))$.
    \item If $\mathit{SA}(z)$ is a special absorbing structure of Type~III, we proceed as follows.
    Let $P^z_1$, $P^z_2$ and $P^z_3$ be the three paths of $\mathit{SA}(z)$. 
    For each $i \in [3]$, let $y^i$ and $x_0^i$ be the first and last vertices of $P^z_i$, respectively.
    Let $v_1,v_2,v_3\in V(\tau_0')$ and $t_1,t_2\in[t]$ be such that $y^1\in V(\cM_{t_1}(v_1))$, $y^2\in V(\cM_{t_1}(v_2))$ and $y^3\in V(\cM_{t_2}(v_3))$ (note this implies that $x^1_0\in V(\cM_{t_2}(v_2))$, $x^2_0\in V(\cM_{t_1}(v_3))$ and $x^3_0\in V(\cM_{t_2}(v_1))$).
    
    For each $i\in[3]$, we define a sequence $S^z_i$ as follows.
    If $t_2 \geq t_1$, let $m^* \coloneqq t_2 - t_1$; otherwise, let $m^* \coloneqq t - (t_1 - t_2)$.
    For each $k\in[m^*-1]_0$, iteratively choose three vertices $y^1_k\in V(\mathcal{A}_{(t_2-k-1)q+1}(v_2))\setminus\cL^\bullet$, $y^2_k\in V(\mathcal{A}_{(t_1+k)q}(v_3))\setminus\cL^\bullet$ and $y^3_k\in V(\mathcal{A}_{(t_2-k-1)q+1}(v_1))\setminus\cL^\bullet$ satisfying that
    \begin{enumerate}[label=$\arabic*$.]
        \item $y_k^i\neqp x_k^i$ for all $i\in[3]$;
        \item $x_{k+1}^1\coloneqq y_k^1+\hat{e}_{(t_2-k-1)q}\notin \cL^\bullet$, $x_{k+1}^2\coloneqq y_k^2+\hat{e}_{(t_1+k)q}\notin \cL^\bullet$ and $x_{k+1}^3\coloneqq y_k^3+\hat{e}_{(t_2-k-1)q}\notin \cL^\bullet$, and 
        \item $\{y_k^i,x_{k+1}^i\}\in E(Q')$ for all $i\in[3]$.
    \end{enumerate}
    Then, for each $i\in[3]$, we define $S^z_i\coloneqq\bigtimes_{k=0}^{m^*-1}(y_k^i,x_{k+1}^i)$.
    
    Let $S\coloneqq\bigtimes_{i=1}^3(P^z_i\times S^z_i)$, and let $S^{-1}$ denote the same sequence in reverse order.
    Finally, we replace the segment $(x(z), x'(z))$ of $\cL$ by $(x(z),S,x'(z))$ if $x(z)$ has parity opposite to the initial vertex of $P^z_1$; otherwise, we replace $(x(z), x'(z))$ by $(x(z),S^{-1},x'(z))$.
\end{enumerate}

Write $\mathcal{L}=(x_1,\ldots,x_r)$, for some $r\in\mathbb{N}$, for the extended skeleton into which all the special absorbing structures $\mathit{SA}(z)$ for $z \in \cU$ have been incorporated, and let $\mathcal{L}^\bullet\coloneqq\{x_1,\ldots,x_r\}$. 
It follows from \ref{itemskelprop1}--\ref{itemskelprop6hit} and the construction above (together with the choice of $v_0$ in Step~12) that the following properties hold:
\begin{enumerate}[label=$(\mathrm{S'}\arabic*)$]
    \item\label{itemskelprop1hitbis} For all distinct $k,k'\in[r]$, we have that $x_k\neq x_{k'}$.
    \item\label{itemskelprop2hitbis} $\{x_1,x_r\}\in E(Q')$.
    \item\label{itemskelprop3hitbis} For every $k\in[r-1]$, if $x_k$ and $x_{k+1}$ do not both lie in the same slice of a cube molecule represented in $\tau_0'$, then $\{x_k,x_{k+1}\}\in E((H\setminus F)\cup H'\cup Q')$.
    Moreover, in this case $\{x_k,x_{k+1}\} \in E(Q')$ unless both $x_k$ and $x_{k+1}$ lie in $\mathit{SA}^v$.
    \item\label{itemskelprop35hitbis} For every $i\in[m']_0$ and every $j\in[t]$, no three consecutive vertices of $\mathcal{L}$ lie in $\cM_j(v_i)$ (here $\mathcal{L}$ is viewed as a cyclic sequence of vertices).
    \item\label{itemskelprop4hitbis} For every $i\in[m']$ such that $v_i$ is an atomic vertex and every $j\in[t]$, we have that $|V(\mathcal{M}_j(v_i))\cap\mathcal{L}^\bullet|$ is even and $4\leq|V(\mathcal{M}_j(v_i))\cap\mathcal{L}^\bullet|\leq14$\COMMENT{It follows from the construction above that we add either none or two vertices inside each slice.}.
    Moreover, $|V(\mathcal{M}_t(v_0))\cap\mathcal{L}^\bullet|=4$.\COMMENT{It is for this particular case that we need the choice of $v_0$.}
    \item\label{itemskelprop5hitbis} For all $k\in[r]$ except two values, we have that $x_k\neqp x_{k+1}$.\COMMENT{We took care to swap the consistent systems to ensure the first vertex of the consistent system we add had parity opposite to that of $x(z)$.
    We now observe that each consistent system is has start and end vertices of opposite parities.
    The reason why this is preserved has to be checked in each case.\\
    For Type~I special absorbing structures, we have added six paths in between two consecutive vertices, each of which has an odd number of vertices.
    This gives a total of an even number of vertices, and they alternate parities, so the overall parity imbalance is $0$.\\
    For Type~II special absorbing structures, we have added two paths and two extra segments.
    For the segments, we are guaranteed that they alternate parities by construction, and they have even length.
    For the paths, again both have an odd number of vertices, so they overall compensate.\\
    Finally, consider Type~III special absorbing structures.
    We have added three paths and three extra segments.
    For the segments, again, we are guaranteed that they alternate parities by construction, and they have even length.
    For the paths, now two of them have an odd number of vertices, and the other one has an even number of vertices.}
    The remaining two values $k_1,k_2\in[r]$ correspond to two pairs of vertices $x_{k_1},x_{k_1+1},x_{k_2},x_{k_2+1}\in V(\mathcal{M}_t(v_0))$.
    For these two values, we have that $x_{k_1}\neqp x_{k_2}$ and either 
    \begin{enumerate}[label=$(\mathrm{\roman*})$]
        \item\label{itemskelprop51hitbis} $x_{k_1}\eqp x_{k_1+1}$ and $x_{k_2}\eqp x_{k_2+1}$, or
        \item\label{itemskelprop52hitbis} $x_{k_1}\neqp x_{k_1+1}$ and $x_{k_2}\neqp x_{k_2+1}$,
    \end{enumerate}
    where $x_{k_1},x_{k_2}\in V(\mathcal{A}_{(t-1)q+1}(v_0))$ and $x_{k_1+1},x_{k_2+1}\in V(\mathcal{A}_{tq}(v_0))$.
    \item\label{itemskelprop6hitbis} $\mathcal{L}^\bullet\cap(\mathfrak{L}\cup\mathfrak{R}_1\cup V^\mathrm{abs})=\varnothing$ and $L^\bullet\subseteq\mathcal{L}^\bullet$.
\end{enumerate}
Indeed, for \ref{itemskelprop3hitbis}--\ref{itemskelprop5hitbis} we make use of the properties of the paths $P_j^z$ defined in \cref{sect:absstruct} as well as \ref{item:AbsStruct} (see \cref{sect:absstruct}) and  \ref{itm:C7primehit}\COMMENT{This is needed to see that the new paths lie in $(H\setminus F) \cup H'$.}.
We also use that the set of all cube molecules represented in $\tau'_0$ is precisely $\cC'''_3 = C'''_1\cupdot\cC'_2$.
Recall that, by the definition of $\cC'''_1$, for each $C \in \cC'''_1$ we have $V(\cM_C) \cap \mathit{SA}^v = \varnothing$.
Moreover, by \ref{itm:bad3} and \ref{item:AbsStruct}, for each $C \in \cC'_2$ we have $|V(\cM_C) \cap \mathit{SA}^v|=2$, and these two vertices $x_C^1, x_C^2\in V(\cM_C) \cap \mathit{SA}^v$ satisfy the following properties:
\begin{enumerate}[label=(\roman*)]
    \item there is some $z\in\cU$ and two consecutive paths $P^z_i$, $P^z_{i+1}$ in $\mathit{SA}(z)$ (with indices taken cyclically) such that $x_C^1$ is the final vertex of $P^z_i$ and $x_C^2$ is the first vertex of $P^z_{i+1}$, and
    \item $\mathbf{C}(\mathit{SA}(z))$ contains two clones $C_1, C_2$ of $C$, where $x^1_C \in V(C_1)$ and $x_C^2 \in V(C_2)$.
\end{enumerate}
Moreover, to check \ref{itemskelprop4hitbis} for case \ref{caseI}, note that the definition of $f$ in a consistent system of paths of Type~I in \cref{sect:absstruct} implies that, for each $i\in[6]$, the final vertex of $P^z_i$ and the first vertex of $P^z_{i+1}$ lie in the same slice.\\

\textbf{Step~18: Constructing an almost spanning cycle.}
Similarly to Step~15 in the proof of \cref{thm:main1}, we will now apply the connecting lemmas to obtain an almost spanning cycle in $(H\setminus F)\cup H'\cup Q'$ from $\mathcal{L}=(x_1,\ldots,x_r)$.
For each $i\in[m']$ such that $v_i$ is an atomic vertex, by \ref{item:Cstarcon1hit} there is at most one value $k(i) \in [t]$ such that $|\mathfrak{L}\cap V(\mathcal{M}_{k(i)}(v_i))|=1$ and $|\mathfrak{L}\cap V(\mathcal{M}_{k(i)+1}(v_i))|=1$.
If such $k(i)$ exists, then we denote by $k^*(i)$ an additional index not in $[t]$ and let $\mathfrak{J}_{i, k^*(i)} \coloneqq \mathfrak{J}_{i, k(i)}\cup \mathfrak{J}_{i, k(i)+1}$, $S_{i, k^*(i)}\coloneqq S_{i, k(i)}\cup S_{i, k(i)+1}$ and $\mathcal{M}_{k^*(i)}(v_i)\coloneqq\mathcal{M}_{k(i)}(v_i) \cup \mathcal{M}_{k(i)+1}(v_i)$.
Let
\[\mathfrak{T}(i)\coloneqq
\begin{cases}
[t] & \text{if there is no }k(i)\text{ as above,}\\
([t]\cup\{k^*(i)\})\setminus\{k(i),k(i)+1\} & \text{otherwise.}
\end{cases}\]
Observe that the definition of $v_0$ in Step~12 together with \ref{itm:C8primehit} ensures that $\mathfrak{T}(0)=[t]$.

For each $i\in[m']_0$ such that $v_i$ is an atomic vertex and for each $j\in\mathfrak{T}(i)$, except the pair $(0,t)$, we apply \cref{lem:slicecover} to  $\mathcal{M}_j(v_i)$ and the graph $Q'$, with $\mathfrak{L}\cap V(\mathcal{M}_j(v_i))$, $(\mathfrak{R}_1\cup\mathfrak{R}_2)\cap V(\mathcal{M}_j(v_i))$ and $S_{i,j}$ playing the roles of $L$, $R$ and the pairs of vertices described in \cref{lem:slicecover}\ref{itm:conn1.3}, respectively.
For $\mathcal{M}_t(v_0)$, we apply \cref{lem:slicecover} or \cref{lem:slicecover2} depending on whether \ref{itemskelprop52hitbis} or \ref{itemskelprop51hitbis} holds in \ref{itemskelprop5hitbis}.
For each $i\in[m']_0$ such that $v_i$ is an atomic vertex and each $j\in\mathfrak{T}(i)$, this yields $|\mathfrak{J}_{i,j}|$ vertex-disjoint paths $(\mathcal{P}^{i,j}_k)_{k\in\mathfrak{J}_{i,j}}$ in $\mathcal{M}_j(v_i)\cup Q'=Q'$ such that, for each $k\in\mathfrak{J}_{i,j}$, properties \ref{itm:AlmHCyc1}--\ref{itm:AlmHCyc3} in Step~15 of the proof of \cref{thm:main1} hold.

Now consider the path obtained as follows by going through $\mathcal{L}$.
Start with $x_1$. 
For each $k\in[r]$, if there exist $i\in[m']_0$ and $j\in\mathfrak{T}(i)$ such that $\{x_k,x_{k+1}\}\in S_{i,j}$, add $\mathcal{P}^{i,j}_k$ to the path; otherwise, add the edge $\{x_k,x_{k+1}\}$ (this must be an edge of $(H\setminus F)\cup H'\cup Q'$ by \ref{itemskelprop3hitbis}).
Finally, add the edge $\{x_r,x_1\}$ of $Q'$ (this is given by \ref{itemskelprop2hitbis}) to the path to close it into a cycle $\mathfrak{H}$ in $(H\setminus F)\cup H'\cup Q'$.
This cycle $\mathfrak{H}$ satisfies \ref{item:almHC1}--\ref{item:almHC3} as in the proof of \cref{thm:main1} as well as the following:
\begin{enumerate}[label=$(\mathrm{HC}\arabic*)$]\setcounter{enumi}{3}
    \item\label{item:almHC4hit} For all $x\in\cU$, we have that $\{x,x+a(x)\},\{x,x+b(x)\}\in E(\mathfrak{H})$.
    \item\label{item:almHC5hit} For all $x \in V(\mathfrak{H})\setminus\mathit{SA}^v$, each of the two edges of $\mathfrak{H}$ incident to $x$ lies in $Q'$.
\end{enumerate}
Indeed, \ref{item:almHC4hit} follows immediately by the definition of $P_1$ in each of the three types of special absorbing structures defined in \cref{sect:absstruct}, and \ref{item:almHC5hit} follows from \ref{itemskelprop3hitbis}.\\

\textbf{Step~19: Absorbing vertices to form a Hamilton cycle.}
Similarly as in Step~16 of the proof of \cref{thm:main1}, for each $u\in V^\mathrm{abs}$ we now replace the edge $e_\mathrm{abs}(u)$ by the path $\mathcal{P}_\mathrm{abs}(u)$ (recall from the end of Step~15 that $\mathcal{P}_\mathrm{abs}(u)$ lies in $((H\cup G)\setminus F)\cup Q'$).
Clearly, this incorporates all vertices of $\mathfrak{L}\cup V^\mathrm{abs}$ into the cycle and, by \ref{item:almHC2} and \ref{item:almHC3}, the resulting cycle $\mathfrak{H}'$ is Hamiltonian.
Moreover, since by \ref{itm:RL4} the endvertices of each edge $e_{\mathrm{abs}}(u)$ lie in cubes belonging to $\cC'''_1$, all these endvertices avoid $\cU$. 
Thus, by \ref{item:almHC4hit}, for each $x \in \cU$ the edges at $x$ in $\mathfrak{H}'$ are still $\{x, x+a(x)\}$ and $\{x, x+b(x)\}$, and so, in particular, by \ref{itm:C7primehit} these edges belong to $H'$.

It now remains to show that $\mathfrak{H}'$ is $(\cU,\ell^2,s)$-good.
Fix any vertex $x\in\cU$.
Let $Y_x\coloneqq N_{\cQ^n}(x)\setminus (V(\mathit{SA}(x))\cup V^\mathrm{abs})$ (that is, by \ref{itm:C4primehit}, $Y_x$ is the set of all vertices in $N_{\cQ^n}(x)\setminus\mathit{SA}^v = N_{\cQ^n}(x)\setminus V(\mathit{SA}(x))$ which lie in clones of cubes which are represented in $\tau_0'$ by atomic vertices).
By \eqref{equa:step10hit}, we have that $|Y_x|\geq(1-2/\ell^4)n-|V(\mathit{SA}(x))|\geq(1-1/\ell^3)n$.
\Cref{claim:rainbowapplnewhit}\ref{item:rainbowapplnew2hit}\COMMENT{\Cref{claim:rainbowapplnewhit}\ref{item:rainbowapplnew2hit} implies that no vertex $y\in Y_x$ can play the role of the left or right absorber tip or of the third absorber vertex for any vertex other than $x$, and we know they do not play this role for $x$ since $x\in\cU$ (so it does not partake in any of the absorbing $\ell$-cube pairs that have been defined at any point throughout the proof).} implies that $Y_x\cap (\mathfrak{L} \cup \mathfrak{R}_1 \cup \mathfrak{R}_2) = \varnothing$, so by definition we have that $Y_x\cap \bigcup_{u \in V^{\mathrm{abs}}}V(\cP_{\mathrm{abs}}(u)) = \varnothing$.
It then follows by \ref{item:almHC5hit} that, for each $y \in Y_x$, each of the two edges of $\mathfrak{H}'$ incident to $y$ lies in $Q'$.
But $Q'$ is $(\cU, 2\ell^2, s)$-good by \cref{claim:Q'good}.
Now, even if all the edges incident to the remaining vertices $y\in N_{\cQ^n}(x)\setminus Y_x$ used the same pair of directions, it follows that the edges of $\mathfrak{H}'$ incident to the vertices in $N_{\cQ^n}(x)$ use each direction of $\cQ^n$ which is not an $s$-direction at most $n/\ell^3 + n/(2\ell^2) \leq n/\ell^2$ times.
\end{proof}


\subsection{Proofs of \texorpdfstring{\cref{thm: kedgehit}}{Theorem~\ref{thm: kedgehit}} and \texorpdfstring{\cref{thm:hitting}}{Theorem~\ref{thm:hitting}}}\label{sect:hit}

We now deduce \cref{thm: kedgehit,thm:hitting} from \cref{thm: important}. 

\begin{proof}[Proof of \cref{thm: kedgehit}]
Let $0<1/n\ll 1/\ell \ll \eps_1\ll \eps\ll \eps_2 \ll \gamma\ll 1/k \leq 1$\COMMENT{We do not need $\gamma\ll 1/k$, but we need $\gamma\ll 1$, so this gives that condition.}.
Let $s\coloneqq 10\ell$.
Let $H^{*}\sim\cQ^{n}_{1/2-2\eps}$ and $Q\sim\cQ^n_\varepsilon$.
Observe that $H^*\cup Q\sim\cQ^n_{1/2-\varepsilon'}$ for some $\varepsilon'\geq\varepsilon$, so it suffices to prove that $H\cup H^*\cup Q$ contains the desired Hamilton cycles and perfect matchings.

By \cref{lem: badvertices} with $2\varepsilon$ playing the role of $\varepsilon$, we have that a.a.s.~$H^{*}$ is $(s,\ell,\eps_1,\eps_2,\gamma,\cU(H^{*},\eps_1))$-robust. 
Condition on this event and let $\cU\coloneqq\cU(H^{*},\eps_1)$.
By an application of \cref{lem: badvertices}\ref{lem:rob2}, it follows that there exists a decomposition of $H^{*}$ into $r\coloneqq\lceil k/2 \rceil $ edge-disjoint spanning subgraphs $H^{*}_1,\ldots, H^{*}_r$ such that, for every $i\in [r]$, we have that $H^{*}_i$ is $(s, \ell, \eps_1/(2r), \eps_2, \gamma/r^{10},\cU)$-robust.

Consider a random decomposition of $Q$ into $r$ edge-disjoint spanning subgraphs $Q_1,\ldots,Q_r$ in such a way that, if $e\in Q$, then $e$ is assigned to one of the $Q_i$ chosen uniformly at random and independently of all other edges. 
It follows that, for all $i\in [r]$, we have $Q_i\sim \cQ^{n}_{\eps/r}$.

Let $\Phi$ be a constant such that \cref{thm: important} holds with $\varepsilon_1/(2r)$, $\gamma/r^{10}$, $\varepsilon/r$ and $r+2$ playing the roles of $\varepsilon_1$, $\gamma$, $\eta$ and $c$, respectively.
(In particular, $\Phi\geq r$.)
For each $i\in [r]$, apply \cref{thm: important} with $H_i^*$, $Q_i$, $\varepsilon_1/(2r)$, $\gamma/r^{10}$, $\varepsilon/r$ and $r+2$ playing the roles of $H$, $Q$, $\varepsilon_1$, $\gamma$, $\eta$ and $c$, respectively, to conclude that a.a.s.~there is a $(\cU,\ell^2,s)$-good subgraph $Q_i'\subseteq Q_i$ with $\Delta(Q_i')\leq\Phi$ such that, for every $H'\subseteq \cQ^{n}$ such that $d_{H'}(x)\geq 2$ for every $x\in \cU$, and every $F\subseteq\cQ^n$ with $\Delta(F)\leq (r+2)\Phi$ which is $(\cU,\ell,s)$-good, we have that $((H^{*}_i\cup Q_i)\setminus F)\cup H'\cup Q_i'$ contains a $(\cU,\ell^2,s)$-good Hamilton cycle $C$ such that, for all $x\in \cU$, both edges of $C$ incident to $x$ belong to $H'$.
Condition on the event that this holds for all $i\in[r]$ (which holds a.a.s.~by a union bound).
Since the $Q_i$ are pairwise edge-disjoint, so are the $Q_i'$.

Now consider the graph $H$ from the statement of \cref{thm: kedgehit}.
By \ref{itm:bad3} in \cref{def:rob}, we can greedily find $r$ edge-disjoint subgraphs $H_1,\ldots ,H_r\subseteq H$ such that
\begin{enumerate}[label=(\roman*)]
    \item\label{itemlastproof1} for each $i\in[\lfloor k/2 \rfloor]$, we have that $|E(H_i)|=2|\cU|$ and $d_{H_i}(x)=2$ for every $x\in \cU$\COMMENT{That is, $H_i$ is a union of cherries, each centered at each $x\in\cU$.}, and
    \item\label{itemlastproof2} if $2r=k+1$, then $H_r$ is a matching of size $|\cU|$ such that $d_{H_r}(x)=1$ for all $x\in\cU$.
\end{enumerate}
\COMMENT{By \ref{itm:bad3} we have that the vertices in $\cU$ are far apart, so we may choose two edges greedily for each part (and at least one final edge for the last part in the case of odd parity).}

Suppose first that $2r=k$.
We are going to find $r$ edge-disjoint $(\cU,\ell^2,s)$-good Hamilton cycles $C_1,\ldots,C_r$ with $H_i\subseteq C_i$ iteratively.
Suppose that for some $i\in[r]$ we have already found $C_1,\ldots,C_{i-1}$.
Let $F_i\coloneqq\bigcup_{j=1}^rQ_j'\cup\bigcup_{j=1}^{i-1}C_j$.
It follows by construction that $F_i$ is $(\cU,\ell,s)$-good\COMMENT{It is the union of at most $2r\ll\ell$ $(\cU,\ell^2,s)$-good graphs.} and $\Delta(F_i)\leq r(\Phi+2)\leq(r+2)\Phi$. 
Then, by the conditioning above, the graph $((H_i^*\cup Q_i)\setminus F_i)\cup H_i\cup Q_i'$ must contain a $(\cU,\ell^2,s)$-good Hamilton cycle $C_i$ such that, for each $u\in\cU$, both edges of $C_i$ incident to $x$ belong to $H_i$.
In particular, $H_i\subseteq C_i$.
Take one such cycle and proceed.

In order to see that these $r$ cycles are pairwise edge-disjoint, suppose that there exist $i,j\in[r]$ with $i<j$ such that $E(C_i)\cap E(C_j)\neq\varnothing$, and let $e\in E(C_i)\cap E(C_j)$.
Observe that $e\notin E(H_i)\cup E(H_j)$ because, otherwise, we would have $e$ incident to some vertex $x\in\cU$, and we know that both edges incident to $x$ in $C_i$ and $C_j$ belong to $H_i$ and $H_j$, respectively, which are edge-disjoint.
Therefore, since $e\in E(C_i)$ and $Q_j'\subseteq F_i\setminus Q_i'$, we must have that $e\notin E(Q_j')$.
However, since $e\in E(C_j)$ and $e\in E(F_j)$ by definition, we must have $e\in E(Q_j')$, a contradiction.

Suppose now that $2r=k+1$.
Let $F_1\coloneqq\bigcup_{j=1}^rQ_j'$, so $\Delta(F_1)\leq r\Phi$ and it is $(\cU,\ell,s)$-good.
By the conditioning above, $((H^{*}_1\cup Q_1)\setminus F_1)\cup H_1\cup Q_1'$ contains a $(\cU,\ell^2,s)$-good Hamilton cycle $C$ with $H_1\subseteq C$. 
We split $C$ into two perfect matchings $M_1\cup M_2$ (observe that both of them are $(\cU,\ell^2,s)$-good) and redefine $H_r\coloneqq H_r\cup\{e\in M_2:\cU\cap e\neq\varnothing\}$, so that $H_r$ now satisfies~\ref{itemlastproof1}. 
Now, for each $i\in \{2,\ldots, r\}$, we proceed as follows.
Let $F_i\coloneqq M_1\cup\bigcup_{j=1}^rQ_j'\cup\bigcup_{j=2}^{i-1}C_j$.
It follows by construction that $F_i$ is $(\cU,\ell,s)$-good and $\Delta(F_i)\leq r(\Phi+2)\leq (r+2)\Phi$.
Then, by the conditioning above, the graph $((H_i^*\cup Q_i)\setminus F_i)\cup H_i\cup Q_i'$ must contain a $(\cU,\ell^2,s)$-good Hamilton cycle $C_i$ with $H_i\subseteq C_i$.
Take one such cycle and proceed.
The fact that the graphs $M_1,C_2,\ldots,C_r$ are pairwise edge-disjoint can be proved as in the previous case.\COMMENT{The proof that the cycles are edge-disjoint is identical to the above.
To see the same for the matching, assume that there is an edge $e\in E(M_1)\cap E(C_i)$, for some $i\in[r]\setminus\{1\}$.
By the same argument as above, we can't have that $e\in E(H_1)\cup E(H_r)$ (where here we abuse notation and use $H_1$ to refer to $H_1\setminus M_2$).
Then, the rest of the argument also follows as above.}
\end{proof}

We now prove \cref{thm:hitting}. 
Recall from \cref{introduction3} that, for any $k\in\mathbb{N}$ and any graph $G\subseteq\cQ^n$, we say that $G\in\boldsymbol{\delta}k$ if $\delta(G)\geq k$, and $G\in\mathcal{HM}k$ if it contains $\lfloor k/2\rfloor$ edge-disjoint Hamilton cycles and $k-2\lfloor k/2\rfloor$ perfect matchings which are edge-disjoint from these cycles.
We say that $G\in\cP k$ if, for every spanning subgraph $H\subseteq \cQ^{n}$ with $H\in\boldsymbol{\delta}k$, we have $G\cup H\in\cH\cM k$.

\begin{proof}[Proof of \cref{thm:hitting}]
The case $k=1$ of the statement was proved by \citet{Bol90}, so we may assume $k\geq2$.
Let $0<\varepsilon\ll1/k$ and $G\sim\cQ^{n}_{1/2-\varepsilon}$.
By \cref{thm: kedgehit}, we have $\mathbb{P}[G\in\cP k]=1-o(1)$.
Also note that, by \cref{lem:Chernoff}, we have that $\mathbb{P}[e(G)\geq(1/2-\varepsilon/2)n2^{n-1}]=o(1)$.
Hence,
\[\mathbb{P}[\{G\in\cP k\} \wedge \{e(G)<(1/2-\varepsilon/2)n2^{n-1}\}]=1-o(1).\]
Thus, by a simple conditioning argument, there exists a positive integer $m<(1/2-\varepsilon/2)n2^{n-1}$ such that 
\begin{equation}\label{equa:hitproof}
    \mathbb{P}[G\in\cP k \mid e(G)=m]=1-o(1).
\end{equation}

Let $G_m\subseteq\cQ^{n}$ be a uniformly random subgraph of $\cQ^{n}$ with exactly $m$ edges.
Since $\mathbb{P}[G\in\cP k \mid e(G)=m]=\mathbb{P}[G_m\in\cP k]$, by \eqref{equa:hitproof} we have 
$\mathbb{P}[G_m\in\cP k]=1-o(1)$. 
Now, condition on the event that $G_m\in\cP k$ and $\tau_{\boldsymbol{\delta}k}(\tilde{\cQ^{n}}(\sigma))\geq (1/2-\varepsilon/4)n2^{n-1}$, which holds a.a.s.\COMMENT{Of course, the hitting time for minimum degree $k$ is larger than the hitting time for having no isolated vertices.
To bound the latter, let $m'\coloneqq(1/2-\varepsilon/4)n2^{n-1}$.
Then, $\mathbb{P}[x \text{ is isolated in } G_{m'}]=\binom{n2^{n-1}-n}{m'}/\binom{n2^{n-1}}{m'}= (1/2+\eta)^{n}$, for some $0<\eta$.
Let $X$ be the number of isolated vertices in $G_m$.
By the above, $\mathbb{E}[X]=2^{n}(1/2+\eta)^{n}$. 
Finally, we note that a simple calculation gives that $\Var[X]=o(\mathbb{E}[X])^2$, and the claim follows by Chebyshev's inequality.
Indeed, $\Var[X]=\sum_{v,v'\in V(\cQ^{n})} \Cov[X_v,X_{v'}]\leq \mathbb{E}[X]+\sum_{v\neq v'\in V(\cQ^{n})} \Cov[X_v,X_{v'}]=\mathbb{E}[X]+2^{n}n\mathbb{P}[v \text{ and } v' \text{ are isolated, where }v,v'\text{ are fixed and adjacent in }\cQ^n]\leq \mathbb{E}[X]+2^{n}n(1/2+\eta)^{n}=o(2^{2n}(1/2+\eta)^{2n})$.}
Then, let $H\coloneqq \cQ^n_{\tau_{\boldsymbol{\delta}k}(\tilde{\cQ^{n}}(\sigma))}(\sigma)$.
It follows that $\delta(H)=k$ and, by the definition of property $\cP k$, we have that $H=G_m\cup H\in\cH\cM k$, so $\tau_{\cH\cM k}(\tilde{\cQ^{n}}(\sigma))\leq\tau_{\boldsymbol{\delta}k}(\tilde{\cQ^{n}}(\sigma))$, as required.
\end{proof}

\section*{Acknowledgement}

We are grateful to the anonymous referees for their helpful comments.

\printindex

\bibliographystyle{afstyle}
\bibliography{cubes_bib}

\begin{thebibliography}{53}
\providecommand{\natexlab}[1]{#1}
\providecommand{\url}[1]{\texttt{#1}}
\providecommand{\urlprefix}{URL }
\providecommand{\selectlanguage}[1]{\relax}
\providecommand{\bibAnnoteFile}[1]{%
  \IfFileExists{#1}{\begin{quotation}\noindent\textsc{Key:} #1\\
  \textsc{Annotation:}\ \input{#1}\end{quotation}}{}}
\providecommand{\bibAnnote}[2]{%
  \begin{quotation}\noindent\textsc{Key:} #1\\
  \textsc{Annotation:}\ #2\end{quotation}}
\providecommand{\eprint}[2][]{\url{#2}}

\bibitem[Ajtai, Koml\'os and Szemer\'edi(1985)]{AKS85}
M.~Ajtai, J.~Koml\'os and E.~Szemer\'edi, First occurrence of {H}amilton cycles
  in random graphs, \emph{Cycles in graphs ({B}urnaby, {B}.{C}., 1982)}, vol.
  115 of \emph{North-Holland Math. Stud.}, 173--178, North-Holland, Amsterdam
  (1985).
\bibAnnoteFile{AKS85}

\bibitem[Alon, Kim and Spencer(1997)]{AKS97}
N.~Alon, J.-H. Kim and J.~Spencer, Nearly perfect matchings in regular simple
  hypergraphs, \emph{Israel J. Math.} \textbf{100} (1997), 171--187.
\bibAnnoteFile{AKS97}

\bibitem[Alon and Spencer(2016)]{AS16}
N.~Alon and J.~H. Spencer, \emph{The probabilistic method}, Wiley Series in
  Discrete Mathematics and Optimization, John Wiley \& Sons, Inc., Hoboken, NJ,
  fourth ed. (2016).
\bibAnnoteFile{AS16}

\bibitem[Alon and Krivelevich(2022)]{AK19}
Y.~Alon and M.~Krivelevich, Hitting time of edge disjoint {Hamilton} cycles in
  random subgraph processes on dense base graphs, \emph{SIAM J. Discrete Math.}
  \textbf{36} (2022), 728--754.
\bibAnnoteFile{AK19}

\bibitem[Azuma(1967)]{Azu67}
K.~Azuma, Weighted sums of certain dependent random variables, \emph{T\^{o}hoku
  Math. J. (2)} \textbf{19} (1967), 357--367.
\bibAnnoteFile{Azu67}

\bibitem[Balogh, Treglown and Wagner(2019)]{BTW19}
J.~Balogh, A.~Treglown and A.~Z. Wagner, Tilings in randomly perturbed dense
  graphs, \emph{Combin. Probab. Comput.} \textbf{28} (2019), 159--176.
\bibAnnoteFile{BTW19}

\bibitem[Bhatt and Cai(1988)]{BC}
S.~Bhatt and J.-Y. Cai, Take a walk, grow a tree, \emph{29th Annual Symposium
  on Foundations of Computer Science (STOC)}, 469--478 (1988).
\bibAnnoteFile{BC}

\bibitem[Bhatt, Chung, Leighton and Rosenberg(1992)]{BCLR92}
S.~N. Bhatt, F.~R.~K. Chung, F.~T. Leighton and A.~L. Rosenberg, Efficient
  embeddings of trees in hypercubes, \emph{SIAM J. Comput.} \textbf{21} (1992),
  151--162.
\bibAnnoteFile{BCLR92}

\bibitem[Bohman, Frieze and Martin(2003)]{BFM03}
T.~Bohman, A.~Frieze and R.~Martin, How many random edges make a dense graph
  {H}amiltonian?, \emph{Random Structures Algorithms} \textbf{22} (2003),
  33--42.
\bibAnnoteFile{BFM03}

\bibitem[Bollob\'{a}s(1983)]{Bol83}
B.~Bollob\'{a}s, The evolution of the cube, \emph{Combinatorial mathematics
  ({M}arseille-{L}uminy, 1981)}, vol.~75 of \emph{North-Holland Math. Stud.},
  91--97, North-Holland, Amsterdam (1983).
\bibAnnoteFile{Bol83}

\bibitem[Bollob\'as(1984)]{Bol84}
B.~Bollob\'as, The evolution of sparse graphs, \emph{Graph theory and
  combinatorics ({C}ambridge, 1983)}, 35--57, Academic Press, London (1984).
\bibAnnoteFile{Bol84}

\bibitem[Bollob\'{a}s(1990)]{Bol90}
B.~Bollob\'{a}s, Complete matchings in random subgraphs of the cube,
  \emph{Random Structures Algorithms} \textbf{1} (1990), 95--104.
\bibAnnoteFile{Bol90}

\bibitem[Bollob\'{a}s(2020)]{BBPC}
B.~Bollob\'{a}s, personal communication (2020).
\bibAnnoteFile{BBPC}

\bibitem[Bollob\'{a}s and Frieze(1985)]{BF85}
B.~Bollob\'{a}s and A.~M. Frieze, On matchings and {H}amiltonian cycles in
  random graphs, \emph{Random graphs '83 ({P}ozna\'{n}, 1983)}, vol. 118 of
  \emph{North-Holland Math. Stud.}, 23--46, North-Holland, Amsterdam (1985).
\bibAnnoteFile{BF85}

\bibitem[Bollob\'{a}s, Kohayakawa and {\L}uczak(1992)]{BKL}
B.~Bollob\'{a}s, Y.~Kohayakawa and T.~{\L}uczak, The evolution of random
  subgraphs of the cube, \emph{Random Structures Algorithms} \textbf{3} (1992),
  55--90.
\bibAnnoteFile{BKL}

\bibitem[Bollob\'{a}s and Thomason(1985)]{BT85}
B.~Bollob\'{a}s and A.~Thomason, Random graphs of small order, \emph{Random
  graphs '83 ({P}ozna\'{n}, 1983)}, vol. 118 of \emph{North-Holland Math.
  Stud.}, 47--97, North-Holland, Amsterdam (1985).
\bibAnnoteFile{BT85}

\bibitem[Borgs, Chayes, van~der Hofstad, Slade and Spencer(2006)]{BCVSS}
C.~Borgs, J.~T. Chayes, R.~van~der Hofstad, G.~Slade and J.~Spencer, Random
  subgraphs of finite graphs. {III}. {T}he phase transition for the {$n$}-cube,
  \emph{Combinatorica} \textbf{26} (2006), 395--410.
\bibAnnoteFile{BCVSS}

\bibitem[B\"{o}ttcher, Han, Kohayakawa, Montgomery, Parczyk and
  Person(2019)]{BHKMPP19}
J.~B\"{o}ttcher, J.~Han, Y.~Kohayakawa, R.~Montgomery, O.~Parczyk and
  Y.~Person, Universality for bounded degree spanning trees in randomly
  perturbed graphs, \emph{Random Structures Algorithms} \textbf{55} (2019),
  854--864.
\bibAnnoteFile{BHKMPP19}

\bibitem[B\"{o}ttcher, Montgomery, Parczyk and Person(2020)]{BMPP20}
J.~B\"{o}ttcher, R.~Montgomery, O.~Parczyk and Y.~Person, Embedding spanning
  bounded degree graphs in randomly perturbed graphs, \emph{Mathematika}
  \textbf{66} (2020), 422--447.
\bibAnnoteFile{BMPP20}

\bibitem[Burtin(1977)]{Burt77}
J.~D. Burtin, The probability of connectedness of a random subgraph of an
  {$n$}-dimensional cube, \emph{Problemy Pereda\v{c}i Informacii} \textbf{13}
  (1977), 90--95.
\bibAnnoteFile{Burt77}

\bibitem[Caha and Koubek(2007)]{CK07}
R.~Caha and V.~Koubek, Spanning multi-paths in hypercubes, \emph{Discrete
  Math.} \textbf{307} (2007), 2053--2066.
\bibAnnoteFile{CK07}

\bibitem[Chan and Lee(1991)]{CL91}
M.~Y. Chan and S.-J. Lee, On the existence of {H}amiltonian circuits in faulty
  hypercubes, \emph{SIAM J. Discrete Math.} \textbf{4} (1991), 511--527.
\bibAnnoteFile{CL91}

\bibitem[Chen(2013)]{Chen13}
X.-B. Chen, Paired many-to-many disjoint path covers of the hypercubes,
  \emph{Inform. Sci.} \textbf{236} (2013), 218--223.
\bibAnnoteFile{Chen13}

\bibitem[Dvo\v{r}\'{a}k and Gregor(2008)]{DG08}
T.~Dvo\v{r}\'{a}k and P.~Gregor, Partitions of faulty hypercubes into paths
  with prescribed endvertices, \emph{SIAM J. Discrete Math.} \textbf{22}
  (2008), 1448--1461.
\bibAnnoteFile{DG08}

\bibitem[Dvo\v{r}\'{a}k, Gregor and Koubek(2017)]{DGK17}
T.~Dvo\v{r}\'{a}k, P.~Gregor and V.~Koubek, Generalized {G}ray codes with
  prescribed ends, \emph{Theoret. Comput. Sci.} \textbf{668} (2017), 70--94.
\bibAnnoteFile{DGK17}

\bibitem[Dyer, Frieze and Foulds(1987)]{Frieze87}
M.~E. Dyer, A.~M. Frieze and L.~R. Foulds, On the strength of connectivity of
  random subgraphs of the {$n$}-cube, \emph{Random graphs '85 ({P}ozna\'{n},
  1985)}, vol. 144 of \emph{North-Holland Math. Stud.}, 17--40, North-Holland,
  Amsterdam (1987).
\bibAnnoteFile{Frieze87}

\bibitem[{Erde}, {Kang} and {Krivelevich}(2021)]{EKK21}
J.~{Erde}, M.~{Kang} and M.~{Krivelevich}, {Expansion in supercritical random
  subgraphs of the hypercube and its consequences}, \emph{arXiv e-prints}
  (2021), \eprint{2111.06752}.
\bibAnnoteFile{EKK21}

\bibitem[Erd\H{o}s and Spencer(1979)]{ES79}
P.~Erd\H{o}s and J.~Spencer, Evolution of the {$n$}-cube, \emph{Comput. Math.
  Appl.} \textbf{5} (1979), 33--39.
\bibAnnoteFile{ES79}

\bibitem[Espuny~D\'iaz(2020)]{Athesis}
A.~Espuny~D\'iaz, \emph{Hamiltonicity problems in random graphs}, Ph.D. thesis,
  University of Birmingham (2020).
\bibAnnoteFile{Athesis}

\bibitem[Fill and Pemantle(1993)]{FP93}
J.~A. Fill and R.~Pemantle, Percolation, first-passage percolation and covering
  times for {R}ichardson's model on the {$n$}-cube, \emph{Ann. Appl. Probab.}
  \textbf{3} (1993), 593--629.
\bibAnnoteFile{FP93}

\bibitem[Frieze(2014)]{FriezeICM14}
A.~Frieze, Random structures and algorithms, \emph{Proceedings of the
  {I}nternational {C}ongress of {M}athematicians---{S}eoul 2014. {V}ol. 1},
  311--340, Kyung Moon Sa, Seoul (2014).
\bibAnnoteFile{FriezeICM14}

\bibitem[{Frieze}(2019)]{Frieze19}
A.~{Frieze}, {Hamilton Cycles in Random Graphs: a bibliography}, \emph{arXiv
  e-prints}  (2019), \eprint{1901.07139}.
\bibAnnoteFile{Frieze19}

\bibitem[Gregor and Dvo\v{r}\'{a}k(2008)]{GD08}
P.~Gregor and T.~Dvo\v{r}\'{a}k, Path partitions of hypercubes, \emph{Inform.
  Process. Lett.} \textbf{108} (2008), 402--406.
\bibAnnoteFile{GD08}

\bibitem[Hahn-Klimroth, Maesaka, Mogge, Mohr and Parczyk(2021)]{HMMMP20}
M.~Hahn-Klimroth, G.~S. Maesaka, Y.~Mogge, S.~Mohr and O.~Parczyk, Random
  perturbation of sparse graphs, \emph{Electron. J. Comb.} \textbf{28} (2021),
  research paper p2.26, 12 pages.
\bibAnnoteFile{HMMMP20}

\bibitem[Hoeffding(1963)]{Hoef63}
W.~Hoeffding, Probability inequalities for sums of bounded random variables,
  \emph{J. Amer. Statist. Assoc.} \textbf{58} (1963), 13--30.
\bibAnnoteFile{Hoef63}

\bibitem[van~der Hofstad and Nachmias(2014)]{HN14}
R.~van~der Hofstad and A.~Nachmias, Unlacing hypercube percolation: a survey,
  \emph{Metrika} \textbf{77} (2014), 23--50.
\bibAnnoteFile{HN14}

\bibitem[van~der Hofstad and Nachmias(2017)]{HN17}
R.~van~der Hofstad and A.~Nachmias, Hypercube percolation, \emph{J. Eur. Math.
  Soc. (JEMS)} \textbf{19} (2017), 725--814.
\bibAnnoteFile{HN17}

\bibitem[Janson, \L{}uczak and Ruci\'{n}ski(2000)]{JLR}
S.~Janson, T.~\L{}uczak and A.~Ruci\'{n}ski, \emph{Random graphs},
  Wiley-Interscience Series in Discrete Mathematics and Optimization,
  Wiley-Interscience, New York (2000).
\bibAnnoteFile{JLR}

\bibitem[Johansson(2020)]{Joha18}
T.~Johansson, On {Hamilton} cycles in {Erd{\H{o}}s}-{R{\'e}nyi} subgraphs of
  large graphs, \emph{Random Struct. Algorithms} \textbf{57} (2020), 132--149.
\bibAnnoteFile{Joha18}

\bibitem[Knox, K\"{u}hn and Osthus(2015)]{KKO15}
F.~Knox, D.~K\"{u}hn and D.~Osthus, Edge-disjoint {H}amilton cycles in random
  graphs, \emph{Random Structures Algorithms} \textbf{46} (2015), 397--445.
\bibAnnoteFile{KKO15}

\bibitem[Knuth(2005)]{Knuth05}
D.~E. Knuth, \emph{The Art of Computer Programming, Volume 4, Fascicle 2:
  Generating All Tuples and Permutations (Art of Computer Programming)},
  Addison-Wesley Professional (2005).
\bibAnnoteFile{Knuth05}

\bibitem[Kohayakawa, Kreuter and Osthus(2000)]{KKO}
Y.~Kohayakawa, B.~Kreuter and D.~Osthus, The length of random subsets of
  {B}oolean lattices, \emph{Random Structures Algorithms} \textbf{16} (2000),
  177--194.
\bibAnnoteFile{KKO}

\bibitem[Koml\'{o}s and Szemer\'{e}di(1983)]{KS83}
J.~Koml\'{o}s and E.~Szemer\'{e}di, Limit distribution for the existence of
  {H}amiltonian cycles in a random graph, \emph{Discrete Math.} \textbf{43}
  (1983), 55--63.
\bibAnnoteFile{KS83}

\bibitem[Kor\v{s}unov(1977)]{Kor77}
A.~D. Kor\v{s}unov, Solution of a problem of {P}. {E}rd\H{o}s and {A}.
  {R}\'{e}nyi on {H}amiltonian cycles in undirected graphs, \emph{Metody
  Diskretn. Anal.} \textbf{31} (1977), 17--56.
\bibAnnoteFile{Kor77}

\bibitem[Krivelevich, Kwan and Sudakov(2017)]{KKS17}
M.~Krivelevich, M.~Kwan and B.~Sudakov, Bounded-degree spanning trees in
  randomly perturbed graphs, \emph{SIAM J. Discrete Math.} \textbf{31} (2017),
  155--171.
\bibAnnoteFile{KKS17}

\bibitem[Krivelevich, Lee and Sudakov(2014)]{KLS14}
M.~Krivelevich, C.~Lee and B.~Sudakov, Robust {H}amiltonicity of {D}irac
  graphs, \emph{Trans. Amer. Math. Soc.} \textbf{366} (2014), 3095--3130.
\bibAnnoteFile{KLS14}

\bibitem[Krivelevich and Samotij(2012)]{KS12}
M.~Krivelevich and W.~Samotij, Optimal packings of {H}amilton cycles in sparse
  random graphs, \emph{SIAM J. Discrete Math.} \textbf{26} (2012), 964--982.
\bibAnnoteFile{KS12}

\bibitem[K\"{u}hn and Osthus(2014{\natexlab{a}})]{KOICM14}
D.~K\"{u}hn and D.~Osthus, Hamilton cycles in graphs and hypergraphs: an
  extremal perspective, \emph{Proceedings of the {I}nternational {C}ongress of
  {M}athematicians---{S}eoul 2014. {V}ol. {IV}}, 381--406, Kyung Moon Sa, Seoul
  (2014{\natexlab{a}}).
\bibAnnoteFile{KOICM14}

\bibitem[K\"{u}hn and Osthus(2014{\natexlab{b}})]{KO14}
D.~K\"{u}hn and D.~Osthus, Hamilton decompositions of regular expanders:
  applications, \emph{J. Combin. Theory Ser. B} \textbf{104}
  (2014{\natexlab{b}}), 1--27.
\bibAnnoteFile{KO14}

\bibitem[Leighton(1992)]{Leigh92}
F.~T. Leighton, \emph{Introduction to parallel algorithms and architectures},
  Morgan Kaufmann, San Mateo, CA (1992), arrays, trees, hypercubes.
\bibAnnoteFile{Leigh92}

\bibitem[{McDiarmid}, {Scott} and {Withers}(2021)]{MSW18}
C.~{McDiarmid}, A.~{Scott} and P.~{Withers}, The component structure of dense
  random subgraphs of the hypercube, \emph{Random Struct. Algorithms}
  \textbf{59} (2021), 3--24.
\bibAnnoteFile{MSW18}

\bibitem[P\'osa(1976)]{Posa76}
L.~P\'osa, Hamiltonian circuits in random graphs, \emph{Discrete Math.}
  \textbf{14} (1976), 359--364.
\bibAnnoteFile{Posa76}

\bibitem[Savage(1997)]{Sava97}
C.~Savage, A survey of combinatorial {G}ray codes, \emph{SIAM Rev.} \textbf{39}
  (1997), 605--629.
\bibAnnoteFile{Sava97}

\end{thebibliography}

\appendix

\section{Proof of \texorpdfstring{\cref{lem:slicecover2}}{Lemma 8.9}}\label{app:connect}

\begin{proof}[Proof of \cref{lem:slicecover2}]
The proof is similar (but easier) to that of \cref{lem:slicecover}.
By relabelling the atoms, we may assume that $\mathcal{M^*}=\mathcal{A}_1\cup\dots\cup\mathcal{A}_t$.
Without loss of generality, we may assume that, for each $r \in [2]$, if $u_r \in R$, then $u_r = z_1$, and if $v_r \in R$, then $v_r=w_t$.
Moreover, we may assume that $x\neqp u_1$.
Let $S\coloneqq\{u_1,v_1,u_2,v_2\}$.
Let $I_R\coloneqq\{k\in[t]: R\cap V(\mathcal{A}_k)\cap S\neq\varnothing\}$, $R^*\coloneqq R \setminus \bigcup_{k \in I_R}V(\cA_{k})$ and $I_{R^*}\coloneqq \{k\in[t]: R^*\cap V(\mathcal{A}_k)\neq\varnothing\}$.
For each $r\in[2]$, let $I_R^r\subseteq\{1,t\}$ be such that $1\in I_R^r$ if and only if $u_r\in R$ and $t\in I_R^r$ if and only if $v_r\in R$.
Note that $I_R=I_R^1\cup I_R^2$. 
Fix an index $t^*\in[t-1]\setminus(I_{R^*}\cup\{1\})$.
If $|L|=2$, let $I_L^1\coloneqq\{i\}$, $I_L^2\coloneqq\{j\}$ and $I_L\coloneqq\{i,j\}$; otherwise, let $I_L^1\coloneqq I_L^2\coloneqq I_L\coloneqq\{t^*\}$\COMMENT{Note that $I_L \cap (I_R\cup I_{R^*}) = \varnothing$.}.

For each $r\in[2]$, we create an ordered list $\mathcal{L}_r$ of vertices, which will be used to construct the vertex-disjoint paths $\mathcal{P}_r$.
Given any list of vertices $\mathcal{L}_r$, we write $L^*_r$ to denote the (unordered) set of vertices in $\mathcal{L}_r$, and whenever $\mathcal{L}_r$ is updated, we implicitly update $L_r^*$.
In the end, for each $r\in[2]$ we will have a list of vertices $\mathcal{L}_r = (x_1^r, \dots, x_{\ell_r}^r)$.
For each $r\in[2]$ and $k \in [t]$, let $I_r(k) \coloneqq \{h\in[\ell_r-1]:2\nmid h\text{ and }x^r_h, x^r_{h+1} \in V(\cA_{k})\}$.
We will require $\cL_1$ and $\cL_2$ to be vertex-disjoint and to satisfy the following properties:
\begin{enumerate}[label=$(\mathcal{L}'\arabic*)$]
    \item\label{itm:2conn02} $\ell_1$ and $\ell_2$ are even.
    \item\label{itm:2conn12} For each $r\in[2]$, for all $h \in[\ell_r-1]$, if $h$ is odd, then $x^r_h,x^r_{h+1}\in V(\cA_{k})$, for some $k \in [t]$; if $h$ is even, then $x^r_h x^r_{h+1} \in E(G\cup\mathcal{M}^*)$.
    \item\label{itm:2conn22} For all $k\in[t]$ we have that $|I_1(k)|,|I_2(k)|\geq1$ and $2 \le |I_1(k)| + |I_2(k)| \le 3$.
    \item\label{itm:2conn31} For each $r\in[2]$, the following holds.
    For each $k\in[t]\setminus(I_L^r\cup I_R^r)$ and each $h\in I_r(k)$, we have $x_h^r \neqp x_{h+1}^r$.
    For each $k\in I_L^r\cup I_R^r$, we have that $|I_r(k)|=1$ and for the unique index $h\in I_r(k)$ we have $x_h^r \eqp x_{h+1}^r$, with the same parity as $u_r$ in the case when $k\in I_L^r$, and with parity opposite to that of the unique vertex in $\{w_k, z_k\} \cap \{u_r, v_r\}$ in the case when $k\in I_R^r$.
    \item\label{itm:2conn4} For each $r\in[2]$, we have the following.
    If $u_r \notin R$, then $u_r = x_1^r$.
    If $v_r \notin R$, then $v_r = x^r_{\ell_r}$.
    If $u_r \in R$ (and thus $u_r = z_{1}$), then $w_{1} = x_1^r$ and $u_r \notin L^*_1 \cup L^*_2$.
    If $v_r \in R$ (and thus $v_r = w_{t}$), then $z_{t} = x_{\ell_r}^r$ and $v_r \notin L^*_1 \cup L^*_2$.
    \item\label{itm:2conn5} Every pair $(w_k, z_k)$ with $\{w_k, z_k\} \subseteq R^*$ is contained in $\cL_1$ and $z_k$ directly succeeds $w_k$ or vice versa.
\end{enumerate}

If $R^*\cap V(\cA_1) = \{w_1, z_1\}$, then let $\mathcal{L}_1 \coloneqq (u_1, w_1, z_1)$, where we assume that $w_1\neqp u_1$; otherwise, let $\mathcal{L}_1\coloneqq(u_1)$.
Observe once more that, in what follows, the existence of each alternating parity sequence follows from the bondedness of $\mathcal{M}$.\COMMENT{It is a lot easier than in the previous lemma.
Here we just form two sequences, each going straight from the bottom to the top, and so that existence of $100$ edges of each parity between consecutive atoms is more than enough at all stages to build our alternating parity sequences.}

Let $F_1 \coloneqq L \cup R \cup S$ and let $t^\bullet_1\in I_L^1$. 
Let $\cS_1$ be a $(u_1, t^\bullet_1, F_1, R)$-alternating parity sequence.
If $u_1\in R$, update $\mathcal{L}_1\coloneqq\mathcal{S}_1$; otherwise, update $\mathcal{L}_1\coloneqq\mathcal{L}_1\cS_1^{-}$\COMMENT{Here there is no need to update $F_1$ (nor $F_2$ in the next paragraph) because we never `reuse' the atoms.}.
Choose any vertex $u_{t^\bullet_1} \in V(\mathcal{A}_{t^\bullet_1})$ with $u_{t^\bullet_1}\neqp u_1$, and let $\cS_2$ be a $(u_{t^\bullet_1}, t, F_1, R^*)$-alternating parity sequence.
Update $\mathcal{L}_1\coloneqq\mathcal{L}_1\cS_2^{-}$.
If $v_1 \in R$, update $\cL_1\coloneqq \cL_1(z_t)$. 
Otherwise, update $\cL_1\coloneqq \cL_1(v_1)$.

Next, let $F_2 \coloneqq F_1 \cup L_1^*$ and let $t^\bullet_2\in I_L^2$. 
Let $\cS_3$ be a $(u_2, t^\bullet_2, F_2, R\cap V(\cA_{1}))$-alternating parity sequence, and let $\mathcal{L}_2\coloneqq\cS_3$.
Choose any vertex $u'_{t^\bullet_2} \in V(\mathcal{A}_{t^\bullet_2})$ with $u'_{t^\bullet_2}\neqp u_2$, and let $\cS_4$ be a $(u'_{t^\bullet_2}, t, F_2, \varnothing)$-alternating parity sequence.
Update $\mathcal{L}_2\coloneqq\mathcal{L}_2\cS_4^{-}$.
Finally, if $v_2 \in R$, update $\cL_2\coloneqq \cL_2(z_t)$. 
Otherwise, update $\cL_2\coloneqq \cL_2(v_2)$.

Observe that $\cL_1$ and $\cL_2$ satisfy \ref{itm:2conn02}--\ref{itm:2conn5}.\COMMENT{Each property follows similarly to how it did in the last lemma. 
See explanation there if needs be.}
We are now in a position to apply \cref{lem:connectcubes}.
For each $k \in [t]$, let $t_k\coloneqq|I_1(k)|+|I_2(k)|$.
Again, for any $r\in[2]$ and $k\in[t]$, for each $h\in I_r(k)$, we refer to the pair $x^r_h, x^r_{h+1}$ as a $\emph{matchable pair}$.
By \ref{itm:2conn22}, \ref{itm:2conn31} and \cref{lem:connectcubes}\ref{lem:connectcubesnormal}, each $\mathcal{A}_k$ with $k\in[t]\setminus(I_L\cup I_{R})$ can be covered by $t_k$ vertex-disjoint paths, each of whose endpoints are a matchable pair contained in $\mathcal{A}_k$.
Similarly, by \ref{itm:2conn22}, \ref{itm:2conn31} and \cref{lem:connectcubes}\ref{lem:connectcubesavoid}, each $\mathcal{A}_k$ with $k \in I_R$ contains $t_k$ vertex-disjoint paths, each of whose endpoints are a matchable pair in $\mathcal{A}_k$, such that the union of these $t_k$ paths covers precisely $V(\mathcal{A}_k)\setminus (S\cap R)$.
Similarly, if $L\neq\varnothing$ and $k\in I_L$, then $\mathcal{A}_k$ contains $t_k$ paths, each of whose endpoints are a matchable pair in $\mathcal{A}_k$, such that the union of these $t_k$ paths covers precisely $V(\mathcal{A}_k)\setminus L$.
Finally, by \ref{itm:2conn22}, \ref{itm:2conn31} and \cref{lem:connectcubes}\ref{lem:connectcubesRareParities}, if $L=\varnothing$ and $k\in I_L$ (that is, $k=t^*$), then $\mathcal{A}_k$ can be covered by $t_k$ paths, each of whose endpoints are a matchable pair in $\mathcal{A}_k$.
For each matchable pair $x^r_h, x^r_{h+1}$ in $\cA_k$, let us denote the corresponding path by $\cP_{x^r_h, x^r_{h+1}}$.

The paths $\mathcal{P}_r$ required for \cref{lem:slicecover2} can now be constructed as follows.
For each $r\in[2]$, let $\mathcal{P}_r$ be the path obtained from the concatenation of the paths $\mathcal{P}_{x^r_{h},x^r_{h+1}}$, for each odd $h\in[\ell_r]$, via the edges $x^r_{h}x^r_{h+1}$ for $h\in[\ell_r-1]$ even.
By \ref{itm:2conn4}, if $\mathcal{P}_r$ does not contain $u_r$, then $\mathcal{P}_r$ starts in $w_1$, and $u_r$ does not lie in any other path; therefore, we can update $\mathcal{P}_r$ as $\mathcal{P}_r\coloneqq u_r\mathcal{P}_r$.
Similarly, if $\mathcal{P}_r$ does not contain  $v_r$, then $\mathcal{P}_r$ ends in $z_t$ and $v_r$ does not lie in any other path, hence we can update $\mathcal{P}_r$ as $\mathcal{P}_r\coloneqq\mathcal{P}_rv_r$.
It follows that $V(\mathcal{P}_1\cup\mathcal{P}_2) = V(\cM^*)\setminus L$, and thus the paths $\mathcal{P}_r$ are as required for \cref{lem:slicecover2}.
\end{proof}

\end{document}